\documentclass[12pt, twoside]{article}
\usepackage{enumitem}
\usepackage{amsmath}

\usepackage{tocloft}
\renewcommand{\cftsecleader}{\cftdotfill{\cftdotsep}}
\usepackage{titletoc}
\usepackage[linguistics]{forest}

\usepackage{bold-extra}
\usepackage{pdfpages}
\usepackage{amsmath,amssymb}
\usepackage{amsthm}
\usepackage{epsfig}
\usepackage{graphics}
\usepackage{graphicx}
\usepackage{colordvi}
\usepackage{mathrsfs} 
\usepackage{stmaryrd}
\usepackage{geometry}
\usepackage{tikz}
\usetikzlibrary{automata}
\usetikzlibrary{snakes}
\usetikzlibrary{decorations.pathmorphing}
\usepackage{leftidx}
\usepackage{comment}
\excludecomment{details}
\usepackage[normalem]{ulem}

\usepackage{dsfont}
\usetikzlibrary{automata}
\usetikzlibrary{decorations.markings}
\usetikzlibrary{chains,shadows.blur}
\usepackage{url}
\usepackage[hidelinks]
{hyperref}
\hypersetup{
    colorlinks=false,
   linkcolor=blue!75!black,
    filecolor=blue!75!black,      
   urlcolor=blue!75!black, 
   citecolor=blue!75!black,
}
\usepackage{enumitem}
\usepackage[textsize=small]{todonotes}
\usepackage[margin=1cm, size=small]{caption}
\usepackage{multirow}
\allowdisplaybreaks[4]
\oddsidemargin
-.25in \evensidemargin -.25in\marginparwidth 0.4in
\usepackage{color,soul}
\usepackage{xcolor}

\usepackage{etoolbox}
\usepackage{fancyvrb}
\usepackage{fvextra}
\usepackage{fancyhdr}
\newcommand{\ps}[1]

\usepackage{empheq}
\usetikzlibrary{automata}
\definecolor{blue2}{cmyk}{.94,.11,0,0}
\usepackage{enumitem}
\setitemize{itemsep=0pt}
\definecolor{myblue}{rgb}{.8, .8, 1}
\newlength\mytemplen
\newsavebox\mytempbox

\DeclareFontFamily{U}{mathx}{}
\DeclareFontShape{U}{mathx}{m}{n}{<-> mathx10}{}
\DeclareSymbolFont{mathx}{U}{mathx}{m}{n}
\DeclareMathAccent{\widecheck}{0}{mathx}{"71}

\makeatletter
\renewenvironment{thebibliography}[1]
     {\section*{\refname}
      \@mkboth{\MakeUppercase\refname}{\MakeUppercase\refname}
      \begin{enumerate}[label={[\arabic{enumi}]},itemindent=*,leftmargin=2.5em]
      \@openbib@code
      \sloppy
      \clubpenalty4000
      \@clubpenalty \clubpenalty
      \widowpenalty4000
      \sfcode`\.\@m}
     {\def\@noitemerr
       {\@latex@warning{Empty `thebibliography' environment}}
      \end{enumerate}}
\makeatother

\makeatletter
\newsavebox{\@brx}
\newcommand{\lla}[1][]{\savebox{\@brx}{\(\m@th{#1\langle}\)}
  \mathopen{\copy\@brx\kern-0.5\wd\@brx\usebox{\@brx}}}
\newcommand{\rra}[1][]{\savebox{\@brx}{\(\m@th{#1\rangle}\)}
  \mathclose{\copy\@brx\kern-0.5\wd\@brx\usebox{\@brx}}}
\makeatother

\makeatletter
\newcommand\cb{
    \@ifnextchar[
       {\@cb}
       {\@cb[5pt]}}

\def\@cb[#1]{
    \@ifnextchar[
       {\@@cb[#1]}
       {\@@cb[#1][5pt]}}

\def\@@cb[#1][#2]#3{
    \sbox\mytempbox{#3}
    \mytemplen\ht\mytempbox
    \advance\mytemplen #1\relax
    \ht\mytempbox\mytemplen
    \mytemplen\dp\mytempbox
    \advance\mytemplen #2\relax
    \dp\mytempbox\mytemplen
    \colorbox{myblue}{\hspace{1em}\usebox{\mytempbox}\hspace{1em}}}
\makeatother
\patchcmd{\thebibliography}{\section*}{\section}{}{}

\newcommand{\cedge}{{\,\leftrightsquigarrow\,}}
\newcommand{\edge}{{\;\leftarrow\!\!\!\!\rightarrow\;}}
\newcommand{\spin}{{\bs \sigma}}

\newcommand{\wt}{\widetilde}
\newcommand{\wh}{\widehat}
\newcommand{\ol}{\overline}

\newcommand{\less}{\lesssim}
\newcommand{\more}{\gtrsim}
\newcommand{\e}{{\mathrm e}}
\newcommand{\DBG}{{\rm DBG}}
\newcommand{\rDBG}{\mathring{\DBG}\mbox{}^{\beta}}

\newcommand{\1}{\mathds 1}

\newcommand{\C}{{\mathscr C}}
\newcommand{\F}{\mathscr F}
\newcommand{\B}{\mathscr B}

\newcommand{\vep}{{\varepsilon}}

\newcommand{\la}{{\langle}}
\newcommand{\ra}{{\rangle}}

\newcommand{\supp}{{\rm supp}}

\newcommand{\bs}{\boldsymbol}

\renewcommand{\P}{{\mathbb P}}
\newcommand{\E}{{\mathbb E}}

\newcommand{\mc}{\mathcal}
\renewcommand{\d}{{\mathrm d}}

\newcommand{\R}{{\Bbb R}}

\newcommand{\lv}{\Lambda_\vep}

\newcommand{\s}{{\mathfrak s}}

\renewcommand{\i}{{\mathtt i}}
\newcommand{\defeq}{{\stackrel{\rm def}{=}}}
\newcommand{\M}{{M}}
\newcommand{\m}{\mathfrak m}

\renewenvironment{proof}[1][\proofname]{\noindent {\bfseries #1.}\;}{\hfill\ensuremath{\blacksquare}\\}

\newcommand{\EM}{\gamma_{\mathsf E\mathsf M}}
\newcommand{\cc}{{^\circ}}

\newcommand{\eqspace}{{\quad\;}}

\newcommand{\bi}{\mathbf i}
\newcommand{\bj}{\mathbf j}

\newcommand{\two}{{\sqrt{2}}}

\newcommand{\SHF}{{\rm SHF}}

\newcommand{\vertiii}[1]{{\left\vert\kern-0.25ex\left\vert\kern-0.25ex\left\vert #1 
    \right\vert\kern-0.25ex\right\vert\kern-0.25ex\right\vert}}
    \newcommand{\nvertiii}[1]{{\vert\kern-0.25ex\vert\kern-0.25ex\vert #1 
    \vert\kern-0.25ex\vert\kern-0.25ex\vert}}
    
\newcommand{\f}{f}

\newcommand{\ovl}{{\vep_{\ref{def:ovl}}(\lambda)}}

\newcommand\ringring[1]{
  {
   \mathop{\kern0pt #1}\limits^{
     \vbox to-1.85ex{
       \kern-2ex 
       \hbox to 0pt{\hss\normalfont\kern.1em \r{}\kern-.45em \r{}\hss}
       \vss 
     }
   }
  }
}

\newtheoremstyle{slantthm}{10pt}{10pt}{\slshape}{}{\bfseries}{}{.5em}{\thmname{#1}\thmnumber{ #2}\thmnote{ (#3)}.}
\newtheoremstyle{slantrmk}{10pt}{10pt}{\rmfamily}{}{\bfseries}{}{.5em}{\thmname{#1}\thmnumber{ #2}\thmnote{ (#3)}.}

\begin{document}
\theoremstyle{slantthm}
\newtheorem{thm}{Theorem}[section]
\newtheorem{prop}[thm]{Proposition}
\newtheorem{lem}[thm]{Lemma}
\newtheorem{cor}[thm]{Corollary}
\newtheorem{disc}[thm]{Discussion}
\newtheorem{conj}[thm]{Conjecture}
\newtheorem{heuristic}[thm]{Heuristic}
\newtheorem{prob}[thm]{Problem}

\theoremstyle{slantrmk}
\newtheorem{ass}[thm]{Assumption}
\newtheorem{rmk}[thm]{Remark}
\newtheorem{defi}[thm]{Definition}
\newtheorem{eg}[thm]{Example}
\newtheorem{que}[thm]{Question}
\newtheorem{exe}[thm]{Exercise}
\numberwithin{equation}{section}
\newtheorem{quest}[thm]{Quest}
\newtheorem{iprob}[thm]{Informal Problem}
\newtheorem{discussion}[thm]{Discussion}
\newtheorem{nota}[thm]{Notation}
\newtheorem{convent}[thm]{Convention}
\newcommand{\thetitle}{
Semimartingale characteristics of the two-dimensional stochastic heat equation at criticality}

\title{\vspace{-2cm}\bf \thetitle\footnote{Support from an NSERC Discovery grant is gratefully acknowledged.}}

\author{Yu-Ting Chen\,\footnote{Department of Mathematics and Statistics, University of Victoria, British Columbia, Canada.}\,\,\footnote{Email: \url{chenyuting@uvic.ca}}\vspace{0cm}}

\date{\today  \vspace{-1cm}}

\maketitle
\abstract{We characterize the sample-path structure of the two-dimensional stochastic heat equation at criticality. The main theorems provide a complete, pathwise description of its semimartingale characteristics and establish explicit, exact formulas for the covariation measure. The proof derives an asymptotic expansion for the covariation measures of regularized approximate solutions by integrating techniques from quantum-solvable operators into stochastic integration theory. \smallskip  
 }


\vspace{-.5cm}

\setlength{\cftbeforesecskip}{-1pt}
\setlength\cftaftertoctitleskip{-1pt}
\renewcommand{\cftsecleader}{\cftdotfill{\cftdotsep}}
\setcounter{tocdepth}{2}
{\small \tableofcontents}
\setlength{\parskip}{0pt}

\section{Introduction}\label{sec:intro}
Following the construction of the critical two-dimensional stochastic heat flow (SHF) as a continuous random field by Caravenna, Sun, and Zygouras~\cite{CSZ:21}, and its subsequent axiomatic characterization via moment matching by Tsai~\cite{Tsai:24}, we consider the problem of characterizing the underlying space-time semimartingale dynamics. Our main results fully determine the semimartingale characteristics of the scaling limits in the stochastic heat equation (SHE) formulation, extending the classical Jacod--Shiryaev framework~\cite{JS:LT} to the infinite-dimensional setting. Explicit, exact representations for the space-time quadratic variations are established. These representations capture the continuous-time dynamics of the limiting SHF while embedding an essential moment formula to determine the underlying coefficients. 

In general dimensions $d\geq 1$, the SHE is described by the following SPDE:
\begin{align}\label{def:SHE}
\frac{\partial X}{\partial t}(x,t)=\frac{\Delta X}{2}(x,t)+\Lambda^{1/2} X(x,t)\xi(x,t), \quad x\in \R^d.
\end{align}
Here, $\Delta$ denotes the $d$-dimensional Laplacian, $\Lambda>0$ is a coupling constant, and $\xi$ represents space-time white noise, a distribution-valued Gaussian field formally characterized by the correlation function
\begin{linenomath*}\begin{align}\label{xi:covar}
\E[\xi(x,t)\xi(y,s)]=\delta_0(x-y)\delta_0(t-s).
\end{align}\end{linenomath*}
In the physics literature, the SHE arises from the formal Cole--Hopf transformation of the Kardar--Parisi--Zhang (KPZ) equation \cite{KPZ}, which describes interface growth, and as a continuum model of directed polymers in random environments (DPRE) \cite{HH:85,KPZ}. A rigorous, fully probabilistic treatment of DPRE originated in the discrete setting with Imbrie and Spencer \cite{IS:88}. This paved the way for the theory of continuum directed random polymers, first introduced within the SHE framework by Alberts, Khanin, and Quastel \cite{AKQ:14a,AKQ:14b} for $d=1$ and also regarded as central to higher-dimensional regimes. For further background and results on DPRE in statistical physics, particularly regarding phase transitions and localization, we refer to Comets's survey \cite{Comets:17} and the references therein.

Owing to its close resemblance to the heat equation with a potential, the SHE described above may appear to be “readily solved” from a heuristic standpoint using Duhamel's formula. In one dimension, this intuition is supported by a rigorous framework: replacing the heuristic field $\xi(x,t)$ with a space-time white noise random measure $\xi(\d x,\d t)$, formally related via
\begin{align}\label{def:xi}
\xi(\d x,\d t)=\xi(x,t)\d x\d t,
\end{align}
renders the SHE \eqref{def:SHE} for $d=1$ classically well-posed via Walsh's extension of It\^{o}'s stochastic integration theory~\cite{Walsh}; see also \cite{Shiga:94} for a survey. The equation thus serves as a canonical model within the one-dimensional KPZ universality class \cite{Corwin:12, QS:15}. The situation changes dramatically for $d \geq 2$, however. In this context, the previously discussed formal solvability encounters significant challenges. This is particularly evident when preserving the multiplication in $\Lambda^{1/2} X\xi$ from \eqref{def:SHE}, which now involves the product of distributions. Furthermore, standard solution theories developed for subcritical singular SPDEs (such as \cite{Hairer:14, GIP:15,Kup:16}) are inherently restricted to subcritical regimes and do not apply here. 

We study the two-dimensional SHE in the critical regime introduced by 
Bertini and Cancrini~\cite{BC:2D}, focusing on bounded, nonnegative initial conditions (see Section~\ref{sec:smt-se} for extensions). Following \cite{GQT}, we refer to the resulting continuous limit as the ``two-dimensional SHE at criticality." Within this framework, solutions are constructed via regularization and noise-normalized approximations. Specifically, for any constant $\lambda\in \R$ and mollifier $\varrho \in \C_c(\R^2)$, we consider 
\begin{linenomath*}\begin{align}\label{def:SHE:vep}
\frac{\partial X_\vep}{\partial t}(x,t)=\frac{\Delta X_\vep}{2}(x,t)+\lv^{1/2} X_\vep(x,t)\xi_\vep(x,t),\quad x\in \R^2,
\end{align}\end{linenomath*}
with fixed function $X_0(\cdot)$ independnet of $\vep$. Here,  the coupling constant  and mollified noise  are defined by
\begin{align}
\lv=\lv(\lambda)&\,\defeq\, \frac{2\pi }{\log \vep^{-1}}+\frac{2\pi \lambda}{\log^2 \vep^{-1}},\label{def:lv}\\
\xi_\vep(x,t)&\,\defeq\, \int_{\R^2}\varrho_\vep(x-y)\xi(y,t)\d y,\label{def:xivep}
\end{align}
where 
\begin{align}
\varrho_\vep(x)&\;\defeq\, \vep^{-2}\varrho(\vep^{-1}x).\label{def:phivep}
\end{align}
Modern SPDE theory ensures the existence of function-valued solutions $X_\vep(x,t)$, whose weak formulation, using
\begin{align}\label{def:Xfintro}
X_\vep(f,t)\,\defeq\, \int_{\R^2}f(x)X_\vep(\d x,t)\quad \mbox{($X_\vep(\d x,t)\,\defeq\, X_\vep(x,t)\d x$)},
\end{align}
is detailed in \eqref{def:weak:vep}. 

The distinctive feature of the approximate equations above is the presence of special coupling constants in \eqref{def:lv}. In \cite{BC:2D}, these constants serve as the crucial condition for the tightness of finite-dimensional marginals of the approximate solutions. This also yields a characterization of the limiting covariance functions that explicitly depend on $(\lambda,\varrho)$. Furthermore, developments over the last decade or so show that, by analyzing higher moments, the limiting solutions are neither Gaussian nor Gaussian multiplicative chaos (GMC). In contrast, markedly different behaviour of the limiting solutions has been established when the definition of $\lv$ in \eqref{def:lv} is replaced by $2\pi \lambda_0/\log \vep^{-1}$ for any fixed $\lambda_0 \in (0,1)$.  In that context, the approximate equations, while converging to the classical heat equation, exhibit Gaussian fluctuations. Specifically, after normalization of the central limit theorem type, the approximate solutions have finite-dimensional marginals that converge to those of the following \emph{additive} SHEs, which fall again within the scope of modern SPDE solution theory:
\begin{align}\label{def:ASHEd=2}
\frac{\partial X}{\partial t}(x,t)=\frac{\Delta X}{2}(x,t)+\left(\frac{1}{1 - \lambda_0}\right)^{1/2}P_tX_0(x)\xi(x,t), \quad x\in \R^2.
\end{align}
Here, $X_0(\cdot)$ refers to the initial condition as above, and $\{P_t\}$ denotes the probability semigroup of two-dimensional standard Brownian motion, with kernels given by the following formula: 
\begin{linenomath*}\begin{align}\label{def:Pt}
P_t(x,y)=P_t(x-y)\,\defeq\, \frac{1}{2\pi t}\exp\left(-\frac{\lvert x-y\rvert^2}{2t}\right),\quad x,y\in \R^2.
\end{align}\end{linenomath*}
Further references and context regarding the effects of these special coupling constants are expanded in Section~\ref{sec:intro-context}.

At criticality, the two-dimensional SHE nevertheless presents a fundamental barrier to standard stochastic integration theory. When a sufficiently regular test function $f(x)$ is applied, the weak formulation of the SHE in \eqref{def:SHE} with $d=2$ naturally yields a formal semimartingale whose quadratic variation process is given by
\begin{align}\label{SHEqv}
\int_0^T\int_{\R^2}\Lambda X(x,t)^2f(x)^2\d x \d t.
\end{align}
The main difficulty with directly interpreting \eqref{SHEqv} is that, at criticality in two dimensions, the SHE admits only measure-valued solutions (see, e.g., \cite{CSZ:25}), rendering the pointwise term $X(x,t)^2$ ill-defined and infinite. Thus, while the ill-defined multiplication in the noise term of \eqref{def:SHE} already complicates its interpretation as a stochastic integral, the difficulty persists even when considering quadratic variations. 

We address this semimartingale structure problem, particularly the characterization of quadratic variations, from a new angle focused on moments. Precise limiting formulas for these moments are well-established via duality with the two-dimensional many-body delta-Bose gas. However, a major obstacle to using them directly is that the $N$-th moments of the two-dimensional SHE at criticality grow too rapidly with $N$ (cf. \cite[Theorem~1.5]{CSZ:23}). Therefore, they fail to uniquely determine the probability distribution via standard moment problems. Nevertheless, these explicit continuum formulas reveal that the \emph{second moment} serves as a fundamental building block for all higher-order moments and imparts their unique character. 

Considering the classical principle that the noise coefficient in a stochastic differential equation is determined by its instantaneous variance, we revisit our earlier analysis of second moments for the two-dimensional SHE at criticality~\cite{C:DBGmix,C:Additive}. This perspective reinforces the fundamental role of second moments in determining all higher-order moments. Furthermore, comparing the formal dynamics in \eqref{SHEqv} with the exact limiting second moments reveals that, even at the level of expectations, the heuristic formulation conceals several subtle, genuinely \emph{second-order} effects. These subtleties extend beyond the specific choice of coupling constants in \eqref{def:lv} or the classical interpretation of quadratic variations as second-order corrections in stochastic calculus, yet evade all heuristic approaches known to us. Section~\ref{sec:smt-mom} provides further context and motivates the key formal expansions of covariation measures introduced in Section~\ref{sec:smt-expansion}. These expansions are the starting point of our methods.

\subsection{Overview of main results}\label{sec:intro-overview}
The main theorems of this paper provide a full pathwise characterization of the semimartingale characteristics of the two-dimensional SHE at criticality. We present explicit formulas for the quadratic variations in terms of the solution, obtained via a pathwise calculus approach. Our methods complement the constructive approach of Caravenna, Sun, and Zygouras~\cite{CSZ:21} and the law characterization framework of Tsai~\cite{Tsai:24}. Also, while our derivations start with the approximate solutions from \eqref{def:SHE:vep}, the resulting formulas are established without a limiting procedure, making them directly applicable to stochastic calculus.

By establishing key moment identities and bounds for subsequent applications, our first main result  (Theorem~\ref{thm:main1}) proves the joint weak convergence of the approximate solution $X_\vep(\d x,t)$ (defined in \eqref{def:Xfintro}) and the covariation measure $\la M_\vep,M_\vep\ra$ of its driving martingale measure $M_\vep(\d x,\d t)$ under a suitable topology. In the limit, this yields the pair consisting of the solution $X$ and the covariation measure $\langle M, M \rangle$ of the limiting martingale measure $M$. Specifically, the martingale measure $M_\vep$ takes the following explicit form:
\begin{align}\label{def:M:vep}
\M_\vep(\d x,\d t)&\,\defeq\,\lv^{1/2}X_\vep(x,t)\xi_\vep(\d x,\d t),\quad (x,t)\in \R^2\times \R_+,
\end{align}
where 
\[
\xi_\vep(\d x,\d t)\,\defeq\int_{\R^2}\varrho_\vep(x-y)\xi(\d y,\d t)\d x
\]
is the spatially mollified space-time white noise. Consequently, the standard weak formulation can be equivalently expressed as
\begin{align}
X_\vep(f,T)
&=X_\vep(f,0)+\int_0^T X_\vep\left(\frac{\Delta f}{2},t\right)\d t+\int_{0}^T \int_{\R^2} f(x)M_\vep(\d x,\d t),\quad f\in \C^2_{pd}(\R^2),\label{def:weak:vep}
\end{align}
where $X_\vep(f,t)$ is defined as in \eqref{def:Xfintro}, and $\C^k_{pd}(\R^2)$ denotes the space of $\C^k(\R^2)$-functions whose partial derivatives up to order $k$ exhibit arbitrary polynomial decay. Furthermore, spatial mollification of the white noise $\xi$ introduces non-local interactions that regularize the system. This results in a covariation measure for $M_\vep$ that requires two spatial variables, denoted by $\la M_\vep, M_\vep \ra(\d x',\d x,\d t)$. In contrast, the limiting martingale measure $M(\d x,\d t)$ for the process $X$, defined as the $\vep\to 0$ analogue of \eqref{def:weak:vep}, is {\bf orthogonal}: its covariation measure is supported on the spatial diagonal $\{(x',x)\in \R^2\times \R^2;x'=x\}$ for all $t\geq 0$, and is accordingly denoted by $\la M,M\ra (\d x,\d t)$. Section~\ref{sec:smt-se} provides further details on the setup, the statement of Theorem~\ref{thm:main1}, and the role of Tsai's law characterization~\cite{Tsai:24} in our main theorems. For background and foundational results on martingale and covariation measures, we refer to \cite{Dalang, EKM:90, Walsh}.

While the weak formulation for $X$ applies to many SPDEs in general and induces semimartingales, it is now distinguished by the following explicit formula for the covariation measure, which crucially reflects the critical nature of the approximate solutions under consideration: for $f\in \C_{pd}^0(\R^2\times \R_+)$ (defined analogously to $\C^0_{pd}(\R^2)$ above),
\begin{align}
& \int_{t=0}^{T}\int_{\R^2} f(x,t)\langle M, M\rangle (\d x,\d t)\notag\\
&\quad =  \int_{t=0}^{T}\int_{\R^2}
 f(x,t)\rDBG_{t}X_0^{\otimes 2}(x)\d x\d t\notag\\
&\quad\quad  +2\int_{t=0}^T \int_{t'=0}^t\int_{\R^2}
 \int_{\R^2}\int_{\R^2}
f(x,t)
\rDBG_{t-t'}(x;z_1,z_2)
X(\d z_2,t') \d x M(\d z_1,\d t')\d t.\label{main:intro}
\end{align}
Here, $X_0^{\otimes 2}(z_1,z_2)\,\defeq\, X_0(z_1)X_0(z_2)$, and $\{P_t\}$ denotes the semigroup corresponding to two-dimensional standard Brownian motion, with kernels specified in \eqref{def:Pt}. The formula in \eqref{main:intro} explicitly involves the kernel $\rDBG_{t}(x;z_1,z_2)$ and the associated linear functional $\rDBG_{t}F(x)$, which generate the semigroup for the two-dimensional two-body delta-Bose gas with parameter $\beta > 0$: 
\vspace{-\abovedisplayskip}
\begin{subequations}\label{def:rDBG}
\begin{align}
\rDBG_{t}(x;z_1,z_2)&\,\defeq\,\int_{s=0}^{t}\int_{\R^{2}}
P_{(t-s)/2}(x,y)\s^\beta(t-s)
\prod_{j=1}^2
P_{s}(y,z_j)\d y\d s,\\
\rDBG_{t}F(x)&\,\defeq\,\int_{\R^{2}}
\int_{\R^2}\rDBG_{t}(x;z_1,z_2)
F(z_1,z_2)\d z_1\d z_2.
\end{align}
\end{subequations}
Here, $\s^\beta(\tau)$ is defined via Euler's gamma function $\Gamma(\cdot)$ as
\begin{gather}
\s^\beta(\tau)\,\defeq\,\int_{0}^\infty \frac{4\pi\beta^u\tau^{u-1}}{\Gamma(u)}\d u.\label{def:sbeta}
\end{gather}
Furthermore, the parameter $\beta$ depends explicitly on the noise strength and mollifier $(\lambda,\varrho)$ from the approximate SHE in \eqref{def:SHE:vep} via
\begin{gather}
\frac{\log \beta}{2}=\log 2+\lambda-\EM-\int_{\R^2}\int_{\R^2}\phi(y')\phi(y)\log \lvert y'-y\rvert \d y'\d y,\label{def:betaintro}
\end{gather}
where $\EM$ is the Euler--Mascheroni constant, and
\begin{gather}\label{def:phi}
\phi(y)\,\defeq\,\int_{\R^2} \varrho(y')\varrho(y'-y)\d y'.
\end{gather}

The formula for the covariation measure in \eqref{main:intro}, along with our technical framework, is fundamentally based on the system's second-order characteristics, aligning with our original motivation (recall the discussion preceding Section~\ref{sec:intro-overview}). The delta-Bose gas kernel $\rDBG_t(x;z_1,z_2)$ not only serves as the primary building block encoding these second-order characteristics, but also has the following direct consequences:
\begin{itemize}
\item Whereas the subcritical additive SHE in \eqref{def:ASHEd=2} is governed by the classical heat kernel, \eqref{main:intro} reveals that the critical two-dimensional regime is instead controlled by the highly nonlinear kernel $\rDBG_t(x;z_1,z_2)$. 
\item Furthermore, this kernel emerges directly from \eqref{main:intro} since taking the expectation on both sides eliminates the second term on the right-hand side (Corollary~\ref{cor:2ndmom}). As a result, \eqref{main:intro} generalizes the classical second-moment formula for the two-dimensional SHE at criticality to the pathwise semimartingale setting. 
\end{itemize}
More broadly, our analysis shows that the precise sample-path behavior of the two-dimensional SHE at criticality is shaped predominantly by second-order characteristics. Quadratic variations and second-moment formulas are essential, while higher moments are invoked exclusively to establish key bounds. Further features of \eqref{main:intro}, alongside several structural variants, are discussed in Section~\ref{sec:smt-se}.

Our semimartingale characterization of the two-dimensional SHE at criticality naturally raises the question of whether the weak formulation and our explicit formulas for the covariation measures uniquely determine the probability distributions. This uniqueness is particularly relevant given that these equations constitute a standard martingale problem formulation. In particular, establishing uniqueness in law via classical methods, as in \cite{EK:MP,SV:MDP}, appears to require techniques distinct from those of the moment-based axiomatic framework of Tsai~\cite{Tsai:24}. Due to the length and scope of this paper, we leave this line of inquiry to separate work, which is currently under our investigation.

\subsection{Context and literature review}\label{sec:intro-context}
To expand on the context surrounding the special coupling constant regimes and the associated limiting behavior, let us highlight several key lines of work in the literature. Regarding evidence of non-Gaussianity, early work by Feng~\cite{Fen:16} provided the first analysis; explicit higher-moment formulas were subsequently developed in \cite{CSZ:19a, C:DBG, GQT} (see Section~\ref{sec:smt-mom} for further discussion); and Caravenna, Sun, and Zygouras~\cite{CSZ:23} established that the limiting solutions do not fall within the class of GMC. For results on the Gaussian fluctuations and the closely related KPZ formulations in the subcritical regime, we refer the reader to \cite{CSZ:19b, CSZ:20, CD:20, Gu:19, NN:23}, among others. In particular, the identification of \eqref{def:ASHEd=2} across the entire subcritical regime for general initial conditions is due to Nakajima and Nakashima~\cite{NN:23}.

Caravenna, Sun, and Zygouras~\cite{CSZ:21} constructed the critical two-dimensional stochastic heat flow (SHF) as a unique scaling limit. This provides a rigorous framework for canonical solutions to the two-dimensional stochastic heat equation at criticality, with the SHF serving as the fundamental solution. The constructions in \cite{CSZ:21} are centred on passing the limit of a specific coarse-grained model that is capable of handling two regularizations widely applied in the earlier literature simultaneously. The first regularization, also the main focus of study in \cite{CSZ:21}, is formulated using a discrete polymer model, in which the underlying random-walk generator and i.i.d. disorder serve as the discrete two-dimensional Laplacian and space-time white noise, respectively. Subject only to mild moment conditions on the disorder, the corresponding limit further exhibits disorder universality. The second regularization considered in \cite{CSZ:21} uses the mollified model described by \eqref{def:SHE:vep}, and is thus directly related to the framework employed in this paper. 

The stochastic heat flow associated with the approximate model in \eqref{def:SHE:vep}  captures key properties of the critical two-dimensional SHF at the pre-limit level. In this context, it is precisely represented by a family of random bi-measures
\[
{\rm SHE}_{\vep,(\lambda,\varrho)}(s,t)={\rm SHE}_{\vep,(\lambda,\varrho)}(\d x,s;\d y,t), \quad x,y\in \R^2,\;0\leq s\leq  t<\infty, 
\]
whose kernels are defined by the following rigorous, though approximate, mathematical formulation for the continuum directed random polymer previously discussed:
\begin{align}
&\hspace{-.2cm}\frac{\SHF_{\vep,(\lambda,\varrho)}(\d x,s;\d y,t)}{\d x\d y}\notag\\
&\hspace{-.2cm}\quad \,\defeq\,\E^{W}_{x,s}\left[\exp\left\{\lv^{1/2}\int_s^t\int_{\R^2}\varrho_\vep(W_{r}-z)\xi(\d z,\d r)-\frac{(t-s)\lv }{2}\int_{\R^2}\varrho_\vep(z)^2\d z\right\}\delta_y(W_{t})\right].\label{SHF:vep}
\end{align}
Here, $W$ denotes a two-dimensional standard Brownian motion, independent of $\xi$, starting from $x$ at time $s$ under $\E^W_{x,s}$, and $\E^W_{x,s}$ leaves $\xi$ unintegrated. For further discussion of the connection between \eqref{SHF:vep} and the continuum directed random polymer, see \cite[Section~1]{C:DBG}, which utilizes a time-reversed Brownian motion approach instead. In the manner of \eqref{SHF:vep}, the approximate stochastic heat flow acts as the fundamental solution to \eqref{def:SHE:vep}, relating to $X_\vep$ via the following equation for appropriate test functions $f$: 
\begin{align}\label{SHEtoSHF}
\int_{\R^2}f(y)X_\vep(y,t)\d y=\int_{\R^2}\int_{\R^2} X_\vep(x,s)
f(y)\SHF_{\vep,(\lambda,\varrho)}(\d x,s;\d y,t),\quad 0\leq s\leq t<\infty.
\end{align}
Further properties of $\SHF_{\vep,(\lambda,\varrho)}$ include the following independence and flow properties: For any integer $n \geq 1$ and $0 \leq s_1 \leq s_2 \leq \cdots \leq s_{n+1}$, the bi-measures over the associated non-overlapping intervals, given by ${\rm SHE}_{\vep,(\lambda,\varrho)}(s_1,s_2), \ldots, {\rm SHE}_{\vep,(\lambda,\varrho)}(s_n,s_{n+1})$, are independent, reflecting the correlation structure of space-time white noise. Additionally, the flow property is directly exhibited by the following identity, generalizing \eqref{SHEtoSHF}:
\begin{align*}
&{\rm SHE}_{\vep,(\lambda,\varrho)}(\d x,r;\d y,t)\\
&\quad =\int_{\R^2}{\rm SHE}_{\vep,(\lambda,\varrho)}(\d x,r;\d z,s)\frac{{\rm SHE}_{\vep,(\lambda,\varrho)}(\d z,s;\d y,t)}{\d z},\quad 0\leq r\leq s\leq t.
\end{align*}
For the existence of the limiting analogues of these independence and flow properties, see Clark and Mian~\cite{CM:26}.

Tsai’s axiomatic characterization \cite{Tsai:24} stands in contrast to the constructive approach by Caravenna, Sun, and Zygouras~\cite{CSZ:21}. The approach in \cite{Tsai:24} focuses directly on the key role of known moment formulas for the two-dimensional SHE at criticality, together with independence and flow properties analogous to those discussed earlier for $X_\vep$. In essence, the conditions isolated in \cite{Tsai:24} allow an extension of Lindeberg’s method \cite{Lind:22} to the SHF context, replacing probability distributions over all small time intervals via moments while maintaining control over cumulative errors. It is notable that Tsai's theorem~\cite{Tsai:24} applies to the mollified SHE in \eqref{def:SHE:vep}, proving the law convergence of its solutions. Specifically, \cite[Section~1.2]{Tsai:24} considers an extension of the SHF driven by space-time white noise defined across the entire time domain $(-\infty,\infty)$, and the SHE is also formulated this way.  

Nakashima~\cite{Nakashima:25} studied the critical two-dimensional SHF by examining the continuum limit of a discretized model, which closely follows the discretization in Caravenna, Sun, and Zygouras \cite{CSZ:21}. By testing the discrete stochastic heat flow against sufficiently regular functions for both initial and terminal conditions (as in \eqref{SHEtoSHF} with $s=0$), Nakashima demonstrated that the weak formulation of the distributional limit produces continuous semimartingales. This leads to a limiting bi-measure-valued process that shares the same finite-dimensional marginals as the critical two-dimensional SHF constructed in \cite{CSZ:21}. Additionally, Nakashima's proof establishes $C$-tightness for the martingale components and their predictable quadratic variations. In Section~\ref{sec:smt-se}, we examine the description of the quadratic variations for these limiting semimartingales given in \cite{Nakashima:25} and discuss how our formulation provides a complementary perspective.\medskip 

\noindent {\bf Frequently used notation.} We write $C(T)\in(0,\infty)$ to denote a constant depending only on $T$, which may change from inequality to inequality unless explicitly indexed. Other constants are defined analogously, and we write $A=\mathcal O_\chi(B)$ if $\lvert A\rvert \leq C(\chi)\lvert B\rvert $. Furthermore, $A\less B$ (or $B\more A$) indicates that $\lvert A\rvert \leq C\lvert B\rvert $ for a universal constant $C\in (0,\infty)$, and $A\asymp B$ means that both $A\less B$ and $B\more A$. Other frequently used notations are listed below:  
\begin{itemize}
\item For a stochastic process $Y$, we use $\E^Y_y$ to denote expectation and $\P^Y_y$ to denote probability, both conditioned on $Y$ starting at $y$.
\item $\log $ denotes the natural logarithm, and $\log^ba\,\defeq\, (\log a)^b$ for all $b\in \R$.
\item For any subset $\Theta$ of $\Bbb R^d$ with $d \in \Bbb N$, $\C_{pd}^k(\Theta)$ denotest the space of functions $f$ on $\Theta$ that extend to a $\C^k$-function on some open set containing $\Theta$, whose all partial derivatives up to order $k$ exhibit arbitrary polynomial decay. 
\end{itemize}

\section{Statement of main theorems}\label{sec:smt}
\subsection{From second moments to covariation measures}\label{sec:smt-mom}
In this paper, we establish our main theorems by integrating asymptotic expansion techniques for quantum-solvable operators into stochastic integration theory. The following expectation formula provides the foundational link to these operators:
\begin{align}\label{Mom:dual}
\E\biggl[\prod_{j=1}^2X_\vep(x_j,t)\biggr]=\E^{B^1,B^2}_{(x_1,x_2)}\biggl[
\exp\biggl\{\lv \int_0^t \phi_\vep(B^2_s-B^1_s)\d s\biggr\}
  \prod_{j=1}^2 X_0(B^j_t)\biggr].
\end{align}
Here, $\Lambda_\vep$ is as defined in \eqref{def:lv}, $\phi_\vep(x)\,\defeq\,\vep^{-2}\phi(\vep^{-1}x)$ is an approximation to the identity based on the mollifier $\phi$ from \eqref{def:phi}, and $B^1,B^2$ are independent two-dimensional standard Brownian motions, starting from $x_1,x_2$ under the expectation $\E^{B^1,B^2}_{(x_1,x_2)}$. Quantum-solvable operators give the purely analytic counterpart of the right-hand side of \eqref{Mom:dual}, because the Brownian exponential moment generates the Feynman--Kac semigroup associated with the following two-particle Schr\"odinger operator:
\[
-\frac{1}{2}\sum_{j=1}^2\Delta_{x_j}-\lv \phi_\vep(x_2-x_1),\quad x_1,x_2\in \R^2,
\]
where $\Delta_x$ denotes the two-dimensional Laplacian with respect to $x$.
This operator rigorously approximates the Hamiltonian of the {\bf two-dimensional two-body delta-Bose gas}:
\[
-\frac{1}{2}\sum_{j=1}^2\Delta_{x_j}-\Lambda' \delta_0(x_2-x_1),\quad x_1,x_2\in \R^2.
\] 

A central feature of the Brownian representation in \eqref{Mom:dual} is the applicability of the following explicit limiting formula, which holds for all sufficiently large $q>0$:
\begin{align}
&\lim_{\vep\to 0}\int_0^\infty \e^{-qt}\E^B_{x}\biggl[\exp\biggl\{\lv\int_0^t\phi_\vep(\two B_s)\d s\biggr\}g(B_t)\biggr]\d t\notag\\
&\quad=\int_0^\infty \e^{-qt}P_t g(x)+\int_0^\infty \e^{-qt}P_t(x)\d t\times \frac{2\pi}{\log (q/\beta)}\int_0^\infty \e^{-qt}P_tg(0)\d t.\label{SingLap}
\end{align}
Here, $B$ is a two-dimensional standard Brownian motion, starting from $x$ under $\E^B_x$, and the parameter $\beta$ is defined as in \eqref{def:betaintro} via $(\lambda,\phi)$. Equation~\eqref{SingLap} highlights a subtle spatial behavior: for any nonzero, nonnegative, and bounded function $g$, the limit diverges at $x=0$ but converges for all $x\neq 0$. Furthermore, this formula enables the exact representation of second moments by utilizing the independence of $B^2 + B^1$ and $B^2 - B^1$, where $B^2-B^1$ has the same distribution as $\two B$. Applying appropriate inverse Laplace transforms then leads to the ultimate untransformed representation (see \cite[Section~3.2]{ABD:95} and \cite[Proposition~5.1]{C:DBG}). Analytically, this probabilistic decoupling corresponds to the classical change of variables $x_2+x_1$ and $x_2-x_1$ for the two-dimensional two-body delta-Bose gas. The reduced operator in the relative coordinate $x_2-x_1$ is the {\bf two-dimensional Laplacian singularly perturbed at the origin}, a crucial relation that forms the foundation of the approach by Bertini and Cancrini~\cite{BC:2D}. 

The study of \eqref{SingLap} originated in the work of Albeverio, Gesztesy, H\o{}egh-Krohn, and Holden~\cite{AGHH:2D}. They introduced functional-analytic techniques for resolvent decompositions to study the two-dimensional Laplacian singularly perturbed at the origin; see also \cite{AGHH:Solvable} for a review and the general framework. Functional-analytic approaches extending this operator-theoretic setting to the many-body delta-Bose gas were subsequently developed in \cite{DFT:94, DR:04}. This connection provides a purely analytic route to determining limiting moments, expressing the higher moments of the two-dimensional SHE at criticality via the many-body semigroup. By utilizing this connection, Gu, Quastel, and Tsai~\cite{GQT} introduced resolvent expansions that generalize the Neumann series to construct the many-body delta-Bose gas semigroup in two dimensions, which directly applies to higher moments of the critical SHE. The motivation for considering the Brownian representation in \eqref{Mom:dual} is detailed in Remark~\ref{rmk:interactingpath}.

Our approach in \cite[Section~4.4]{C:Additive} to proving \eqref{SingLap} differs from earlier methods by recasting the problem in terms of asymptotic recursions, bypassing the heavy reliance on constructing preliminary exact identities. (This viewpoint is particularly useful for developing an asymptotic expansion in the more subtle scenario where the associated mollifier has zero integral. We do not employ such results in this paper, though.) In the current setting, where $\phi$ is a probability density, the approach of \cite{C:Additive} focuses on the following Green functions, evaluated with \emph{a small initial condition} $\vep x$ and an arbitrary nonzero bounded nonnegative function $g$:
\begin{align}\label{RM:Green}
G_\vep(x)\stackrel{\rm def}{=} \int_0^\infty \e^{-qt}\E^B_{\vep x}\biggl[\exp\biggl\{\lv\int_0^t\phi_\vep(\two B_s)\d s\biggr\}g(B_t)\biggr]\d t.
\end{align}
The use of small conditions allows capturing the primary mechanism underlying the limiting formula in \eqref{SingLap}, as we have the following limit: 
\begin{align}\label{SingLap1}
\lim_{\vep\to 0} \lv G_\vep(x)=\frac{2\pi}{\log (q/\beta)}\int_0^\infty \e^{-qt}P_tg(0)\d t,
\end{align}
which holds uniformly for $x$ in compact subsets of $\R^2$ and all sufficiently large $q>0$.

The proof of \eqref{SingLap1} in \cite[Section~4.4]{C:Additive} centers on constructing preliminary asymptotic expansions for $G_\vep(x)$, using only the standard asymptotic expansion of the Green’s function of two-dimensional Brownian motion (see \eqref{rmk:EM}) and a priori bounds for $G_\vep(x)$, which are independently obtained via a convolution-type Gr\"onwall inequality in \cite{C:DBG}. This leads to the conclusion that $G_\vep(x)$ admits the following asymptotic expansion, which can also be interpreted as an asymptotic recursion:
\begin{gather}\label{ARE:G}
G_\vep(x)=C_1+(1-s_\vep C_2) G_\vep(x)+o(1),\quad \vep\to 0.
\end{gather}
Moreover, we have the following specifications of $s_\vep$, $C_1$, and $C_2$:
\[
s_\vep\equiv \lv, \quad C_1\equiv \int_0^\infty \e^{-qt}P_tg(0)\d t,\quad C_2\equiv \frac{\log (q/\beta)}{2\pi}. 
\]
A direct consequence of \eqref{ARE:G} is $\lim_{\vep\to 0}s_\vep G_\vep(x)=C_1/C_2$. Here, the seemingly simple prefactor $1-s_\vep C_2$ plays a central role. The constant $1$ identifies the system’s criticality: it is replaced by $\lambda_0$ in the subcritical regime, where $\lv$ from \eqref{def:lv} is replaced by $2\pi \lambda_0/\log \vep^{-1}$ for $\lambda_0 \in (0,1)$. Furthermore, the term $s_\vep C_2$ highlights a crucial second-order contribution that underscores the system’s critical nature. In the subcritical limit, such a term plays no role, nor does introducing a second-order correction to the renormalizing constant.

We propose that the asymptotic recursion method from \cite{C:Additive}, discussed above, is deeply connected to the microscopic sample-path dynamics of the two-dimensional SHE at criticality. An immediate challenge in this pursuit is determining whether the relationship between the second moments and the critical SHE is sufficiently robust. Nonetheless, as with SDEs (recall the discussion before Section~\ref{sec:intro-overview}), the noise coefficient of an SPDE can be interpreted as an instantaneous variance. This perspective is particularly prominent in the derivation of SPDEs and martingale problems for super-Brownian motion, and more generally, the Dawson--Watanabe superprocesses (see, e.g., \cite{Dawson:75, KS:88, MP:92, Perkins:DW, R:89,Watanabe:68}), which serve as natural analogues of the SHE in \eqref{def:SHE}. In the context of super-Brownian motion, the noise coefficient $\Lambda^{1/2}X$ in \eqref{def:SHE} is replaced by $\Lambda^{1/2}\sqrt{X}$. Moreover, an analogue of the formal variance of $\Lambda^{1/2} X(x,t)$ for the regularized SHE in two dimensions satisfies the following equation (see \cite[Section~2.3.1]{C:DBG} for a detailed derivation):
\begin{align*}
&{\rm Var}\big(\lv^{1/2}X_\vep(x,t)\big)+ \lv [P_t X_0(x)]^2\\
&\quad  =\lv \E^B_0\biggl[\exp\left\{\lv \int_0^t \phi_\vep(\two B_s)\d s\right\}
\int_{\R^2} P_t(\two x ,y)\prod_{\varsigma\in \{\pm\}}
X_0\left(\frac{y+\varsigma B_t}{\two}\right)
\d y
\biggr].
\end{align*}
This illustrates a direct connection to the Green function $G_\vep(0)$ as defined in \eqref{RM:Green}.

The following problem translates the above heuristics into a concrete formulation that reflects our initial idea for the method of this paper.

\begin{prob}\label{prob:ARE}
Rigorously formulate the following pathwise stochastic asymptotic recursion extending \eqref{ARE:G}, and derive the associated limiting equation:
\[
X_\vep(x,t)^2=F_0(x,t)+\left[1-s_\vep F_1(x,t)\right] X_\vep (x,t)^2+F_{2,\vep}(x,t)+o(1),\quad \vep\to 0,
\]
where $s_\vep\to 0$, $F_0, F_1$ are $\vep$-independent random fields, and $F_{2,\vep}$ is a mean-zero random field that converges as $\vep\to 0$, all with explicit representations. 
\end{prob}

We stress that the expansion considered in Problem~\ref{prob:ARE} is for $X_\vep(x,t)^2$, not $\lv X_\vep(x,t)^2$. Therefore, applications to quadratic variations require $s_\vep=\Theta(\log^{-1}\vep^{-1})$ as $\vep\to 0$.
 
An immediate technical challenge posed by Problem~\ref{prob:ARE} is the lack of resolvents, tools that have proven crucial in existing proofs for deriving limiting second moments and analyzing related quantum-solvable models. Nevertheless, by directly targeting the sample-path dynamics, our formulation in Problem~\ref{prob:ARE} bypasses iterations of the mild formulation (see \eqref{def:mild:vep}). Specifically, this avoids relying on chaos series expansions. While effective throughout the subcritical regime, these expansions do not seem to provide a clear path for determining precise sample-path descriptions, such as the semimartingale characteristics considered here.

\begin{rmk}[Interacting paths in moments]\label{rmk:interactingpath}
In contrast to the operator-theoretic techniques discussed above, our emphasis on the path-measure structure in \eqref{Mom:dual} continues the line of development initiated in \cite{C:DBG}, where we established an alternative, probabilistic-analytic construction of higher moments. This path-measure structure enriches the SPDE and semimartingale analysis in the present work. 

The physical formulation of this approach goes back to Kardar and Zhang~\cite{KZ:87}, who studied $N$-th moments in the un-mollified regime using formal, real-valued Feynman path integrals. The associated $N$ polymers were considered as paths that interact via the correlation function of the space-time white noise. A comprehensive review is provided in \cite{Kardar:SPF}. Bertini and Cancrini~\cite{BC:2D} provided the first mathematical formulation of this idea, viewing the second moments for $\varepsilon>0$ as exponential moments of approximate intersection local times of two Brownian paths. In the continuum limit, these paths were expected to exhibit the attractive contact interactions predicted by Kardar and Zhang. However, a direct proof was obstructed by the absence of Brownian local times in two dimensions. A rigorous construction confirming that the limiting second moments indeed support such interacting paths via precise Feynman--Kac-type formulas was proved in \cite{C:BES}. Subsequently, Feynman--Kac-type formulas for higher-order moments were established across a series of works \cite{C:SDBG1, C:SDBG2, C:SDBG3, C:SDBG4}, where the underlying interacting paths are constructed as multiple Brownian motions ``conditioned for indefinite pairwise contact interactions.'' Notably, this higher-order construction relies heavily on the second-moment diffusions established in \cite{C:BES}, underscoring the central role that second moments play throughout the perspective introduced in Section~\ref{sec:intro} and elaborated on above. \hfill $\blacklozenge$
\end{rmk}

\subsection{Formal expansions of covariation measures}\label{sec:smt-expansion}
In this subsection, we present a formal solution to Problem~\ref{prob:ARE} in Heuristic~\ref{heu:formaleq} below, leaving the second-order term unspecified. This approach elevates moment analysis to the pathwise setting and constructs a fully independent dynamical asymptotic expansion, positioning It\^o's formula as the natural tool. Remark~\ref{rmk:subcritical} addresses how this formal solution extends to the subcritical case. This provides a heuristic explanation for the emergence of the known additive SHE limit in \eqref{def:ASHEd=2} and unifies these approaches up to criticality. Finally, Remark~\ref{rmk:supercritical} discusses the fundamental obstacles faced when extending our approach to the supercritical regime.

We adopt the following notation for the formal solution. Since our analysis examines the limiting behaviour of $X_\vep(x,t)^2$ as $\vep \to 0$, we suppress explicit dependence on $\vep$ and let $X(x,t)$ denote the formal limit of $X_{\vep} ( x, t)$ as $\vep  \to 0$. In this context, rather than taking $\lim_{\vep\to 0}\Lambda_\vep=0$, we preserve the renormalizing effect of $\Lambda_\vep$ by using the placeholder $\Lambda$, which is interpreted as $0+$ when we evaluate limiting expressions. {\bf This notation is therefore fully consistent with the formal SHE (\ref{def:SHE}) for $d=2$ and applies throughout Section~\ref{sec:smt-expansion} unless stated otherwise.} Additionally, we work with the following purely formal martingale measure: 
\[
M(\d x,\d t)=\Lambda^{1/2}X(x,t)\xi(\d x,\d t),
\]
which in turn has the following purely formal covariation measure:
\begin{align}\label{Mqv:formal}
\la M,M\ra(\d x,\d t)=\Lambda X(x,t)^2 \d x\d t;
\end{align}
see \eqref{def:M:vep} for analogues at the approximate level. Note that while we specify the components of $M$ for its description, we will avoid the ill-posed product $\Lambda^{1/2}X\xi$ as in \eqref{def:SHE} by treating $M$ as an indivisible object.

\begin{heuristic}\label{heu:formaleq}
Under the notation introduced above, we have
\begin{align}
[P_tX_0(x)]^2
&=X(x,t)^2-\int_{s=0}^t \int_{\R^2} P_{t-s}(x,y)^2 \Lambda X(y,s)^2 \d y\d s\notag\\
&\quad -2\int_{t'=0}^t\int_{\R^2} \int_{\R^2} 
\biggl(\prod_{j=1}^2P_{t-t'}(x,z_j)\biggr)
X(z_2,t')\d z_2\M(\d z_1,\d t')
,\label{asymp:2}\\
[P_tX_0(x)]^2&=X(x,t)^2-X(x,t)^2\notag\\
&\quad -2\int_{t'=0}^t\int_{\R^2} \int_{\R^2} 
\biggl(\prod_{j=1}^2P_{t-t'}(x,z_j)\biggr)
X(z_2,t')\d z_2\M(\d z_1,\d t'),\label{formaleq}
\end{align}
where $P_t(x, y)$ denotes the transition density of two-dimensional standard Brownian motion as defined in \eqref{def:Pt}, and $\{P_t\}$ is the associated heat semigroup.
\end{heuristic}

Equation~\eqref{formaleq} plays a fundamental role in this work and refines the expansion in \eqref{asymp:2}. We highlight two key features of \eqref{formaleq}, deferring further discussion to Remark~\ref{rmk:expansion}:
\begin{itemize}
\item \emph{Heuristic interpretation of the formal components:} The terms on the right-hand side of \eqref{formaleq} reflect two distinct roles:
\begin{itemize}
\item We emphasize that the difference $X(x,t)^2 - X(x,t)^2$ in \eqref{formaleq} is a purely symbolic representation rather than a vanishing real-valued difference. This notation is chosen both to signify a formal $\infty - \infty$ indeterminate limit resulting from singular cancellation (because the limiting process is measure-valued, as discussed around \eqref{SHEqv}) and to encapsulate the second-order term considered in Section~\ref{sec:smt-mom}. Specifically, the derivation below reveals the second $X(x,t)^2$ as the leading-order term in the divergent integral $\int_{s=0}^t \int_{\R^2} P_{t-s}(x,y)^2 \Lambda X(y,s)^2 \d y\d s$ in \eqref{asymp:2}. In turn, this integral arises naturally as the $\vep\to 0$ limit of a regularized integral whose precise form must be identified through precise $\vep$-level calculations. 

\item The formal equation \eqref{formaleq} identifies a mean-zero random field given by the stochastic integral
\begin{align}\label{RF:mean0}
-2\int_{t'=0}^t\int_{\R^2} \int_{\R^2} 
\biggl(\prod_{j=1}^2P_{t-t'}(x,z_j)\biggr)
X(z_2,t')\d z_2\M(\d z_1,\d t').
\end{align}
\end{itemize}
\item \emph{Analytical challenges in rigorous justification:} The formal derivation of \eqref{formaleq} proceeds via heuristic manipulations that fully bypass potential singularities. However, rigorously establishing the complete asymptotic expansion and extracting the second-order term pose two major analytical challenges:

\begin{itemize}
\item The first major challenge arises from omitting perturbative terms in these calculations, an act that conceals the crucial but latent second-order contribution emphasized in Section~\ref{sec:smt-mom}, and from managing any resulting singularities in the limiting process. To address this, our forthcoming proofs rely on the weak formulation to absorb pointwise singularities. We then systematically track all $\vep$-perturbations, establish uniform-in-$\vep$ moment bounds, and identify the exact limiting second-order term. This analysis draws on baseline results for SPDEs with colored noise due to Dalang~\cite{Dalang} and Mytnik, Perkins, and Sturm~\cite{MPS:06}, alongside a general theorem on $C$-tightness of laws of measure-valued processes due to Perkins~\cite{Perkins:DW}. Our earlier methods for bounding moments of the two-dimensional SHE at criticality are also adapted and extended.

\item Another major technical challenge is to provide rigorous justifications for the existence of the stochastic integral in \eqref{RF:mean0}, whose pre-limiting processes prove to be mean-zero random fields sought in Problem~\ref{prob:ARE}, as well as for the convergence of the approximating sequences. To address the convergence, the Kurtz--Protter method~\cite{KP:91} offers crucial tools for establishing the convergence of ordinary stochastic integrals $\int Y_n \d Z_n \to \int Y \d Z$ as $(Y_n, Z_n) \to (Y, Z)$.  In this context, we also make use of Cho's extension~\cite{Cho:95} to stochastic integrals with respect to martingale measures.
\end{itemize}
\end{itemize}
	
\begin{proof}[Derivation of Heuristic~\ref{heu:formaleq}]
We start by expanding $X(x,t)^2$, which determines the covariation measure 
in \eqref{Mqv:formal}. To proceed, we assume and utilize the mild form satisfied by the formal limit $X(x,t)$ mentioned above:
\begin{align}\label{mild:formal}
X(x,t)=P_tX_0(x)+\int_{s=0}^t \int_{\R^2}P_{t-s}(x,y)M(\d y, \d s).
\end{align}
We then apply \eqref{mild:formal} to expand $X(x,t)^2$ and rearrange the resulting equation via the algebraic identity 
\[
a^2=c^2-2ab-b^2\quad \mbox{given }c=a+b.
\]
This yields: 
\begin{align}
[P_tX_0(x)]^2
&=X(x,t)^2-2P_tX_0(x)\int_{s=0}^t \int_{\R^2} P_{t-s}(x,y)\M(\d y,\d s)\notag\\
&\quad -\left(\int_{s=0}^t\int_{\R^2} P_{t-s}(x,y)\M(\d y,\d s)\right)^2\notag\\
&=X(x,t)^2-2P_tX_0(x)\int_{s=0}^t \int_{\R^2} P_{t-s}(x,y)\M(\d y,\d s)\notag\\
&\quad -2\int_{s'=0}^t \int_{\R^2}\int_{s=0}^{s'}\int_{\R^2} P_{t-s}(x,y)\M(\d y,\d s)  P_{t-s'}(x,y')\M(\d y',\d s')\notag\\
&\quad -
\int_{s=0}^t \int_{\R^2} P_{t-s}(x,y)^2 \Lambda X(y,s)^2 \d y\d s.\label{asymp:1}
\end{align}
Here, we apply It\^{o}'s formula to the martingale $\int_{s=0}^{s'}\int_{\R^2} P_{t-s}(x,y)\M(\d y,\d s)$, $0\leq s'\leq t$, which produces the last two terms in \eqref{asymp:1}.

Next, the following calculation suffices to demonstrate that \eqref{asymp:1} simplifies to \eqref{asymp:2}:
\begin{align}
&-2\int_{s=0}^t \int_{\R^2}\int_{s'=0}^s\int_{\R^2} P_{t-s'}(x,y')\M(\d y',\d s')  P_{t-s}(x,y)\M(\d y,\d s)\notag\\
&\quad =-2\int_{s=0}^t\int_{\R^2} \int_{s'=0}^s\int_{\R^2} P_{t-s}(x,z)P_{s-s'}(z,y')\d z\M(\d y',\d s')P_{t-s}(x,y)\M(\d y,\d s)\notag\\
&\quad =-2\int_{s=0}^t \int_{\R^2} \int_{\R^2}  P_{t-s}(x,z)[-P_sX_0(z)+X(z,s)]\d zP_{t-s}(x,y)\M(\d y,\d s)\notag\\
&\quad =2P_tX_0(x)\int_{s=0}^t\int_{\R^2}P_{t-s}(x,y) \M(\d y,\d s)\notag\\
&\quad \quad -2\int_{s=0}^t\int_{\R^2} \int_{\R^2} 
P_{t-s}(x,y)
P_{t-s}(x,z)
X(z,s)\d z\M(\d y,\d s),\notag
\end{align}
where the first and last equalities rely on the Chapman--Kolmogorov equation, and the second equality follows from the mild form \eqref{mild:formal}. 

The crucial step in this derivation, leading to \eqref{formaleq}, is to apply the following observation. Because we are effectively in the regime $\vep = 0+$, a heat kernel expansion [see Lemma~\ref{lem:heat1} (1$\cc$)] suggests that, by neglecting subleading terms, 
\begin{align}\label{special:ati}
\Lambda P_{t-s}(x,y)^2\d s\d y= \delta_t(\d s)\delta_x(\d y),\quad 0\leq s\leq t,\;x,y\in \R^2,
\end{align}
where $\delta_t(\d s)$ and $\delta_x(\d y)$ represent Dirac measures. Thus, we expect the following equation:
\[
\int_{s=0}^t \int_{\R^2} P_{t-s}(x,y)^2 \Lambda X(y,s)^2 \d y\d s= X(x,t)^2.
\]
Applying this equation to \eqref{asymp:2} yields the formal equation \eqref{formaleq}.
\end{proof}
\vspace{-\abovedisplayskip}

\begin{rmk}\label{rmk:expansion}
(1$\cc$) In \eqref{asymp:1}, a double integrator $M(\d y, \d s)M(\d y', \d s')$ is introduced. We reformulate this in \eqref{asymp:2} using $X(z_2, t')\d z_2 M(\d z_1, \d t')$, as Lebesgue integration with respect to $X(z_2, t')\d z_2$ offers significant technical advantages, such as straightforward interchange of the order of integration. The moment bound in \eqref{X4BMbdd}, which serves as a central tool in the proof of Theorem~\ref{thm:main3}, provides another example of this convenience.\medskip 

\noindent (2$\cc$)   Proposition~\ref{prop:spdelta} provides a rigorous formulation of \eqref{special:ati}. \medskip

\noindent (3$\cc$) A review of the derivation for \eqref{asymp:2} reveals that its counterpart for the one-dimensional SHE is mathematically rigorous. Specifically, for any function-valued solution $Y$ to the one-dimensional SHE, given by \eqref{def:SHE} with $d=1$ using an arbitrary order-$1$ constant $\Lambda>0$, the following equation holds: 
\begin{align}
&Y(x,t)^2-\int_{s=0}^t \int_\R Q_{t-s}(x,y)^2\Lambda Y(y,s)^2\d y \d s\notag\\
&\quad =[Q_tY_0(x)]^2+2\int_{t'=0}^t\int_{\R} \int_{\R} 
\biggl(\prod_{j=1}^2Q_{t-t'}(x,z_j)\biggr)
Y(z_2,t')\d z_2M_Y(\d z_1,\d t').\label{eq:Y3}
\end{align}
Here, $Q_t(x, y)$, $x, y \in \R$, are the transition densities of one-dimensional standard Brownian motion, and $M_Y(\d x, \d t)$ denotes the martingale measure $\Lambda^{1/2}Y(x, t)\xi(\d x, \d t)$.
Furthermore, with appropriate adjustments to the contributions from covariation measures, \eqref{asymp:2}  and \eqref{asymp:1} apply rigorously to any well-defined mild form of a general stochastic heat equation, including the approximate solutions $X_\vep$ to \eqref{def:SHE:vep}. Proposition~\ref{prop:qvdec} examines the details for the approximate solutions, focusing specifically on the formulation in \eqref{asymp:1}.
\hfill $\blacklozenge$  
\end{rmk}

\begin{rmk}[Subcritical regime]\label{rmk:subcritical}
The above argument can be adapted to formally derive Gaussian fluctuations in the subcritical case, where $\lv$ is replaced by $ 2\pi \lambda_0 / \log \vep^{-1}$ for $\lambda_0 \in (0,1)$. This derivation can also be made rigorous using the methods presented in this paper, and the corresponding proof is expected to be comparatively simpler. 

In this subcritical regime, the formal equation in \eqref{formaleq} becomes
\begin{align}
[P_tX_0(x)]^2&=(1-\lambda_0)X(x,t)^2\notag\\
&\quad -2\int_{t'=0}^t\int_{\R^2} \int_{\R^2} 
\biggl(\prod_{j=1}^2P_{t-t'}(x,z_j)\biggr)
X(z_2,t')\d z_2\M(\d z_1,\d t').\label{formaleqv0}
\end{align}
Hence, multiplying both sides by $\Lambda=0+$ yields
\[
0=(1-\lambda_0)\Lambda X(x,t)^2+0\Longrightarrow \la M,M\ra(\d x,\d t)=0\Longrightarrow \M(\d x,\d t)=0,
\]
and \eqref{formaleqv0} simplifies to
\[
[P_tX_0(x)]^2=(1-\lambda_0)X(x,t)^2\Longrightarrow 
X(x,t)^2\d x\d t=\frac{1}{1-\lambda_0}[P_tX_0(x)]^2\d x\d t.
\]
The final equation indicates that the covariation measure of $\Lambda^{-1/2}M(\d x,\d t)$ is deterministic. This martingale measure is thus expected to be Gaussian; see \cite[Proposition~2.10, p.302]{Walsh} for a general theorem of this property. 
\hfill $\blacklozenge$
\end{rmk}

\begin{rmk}[Supercritical regime]\label{rmk:supercritical}
The mathematical justification for the heuristic argument in Heuristic~\ref{heu:formaleq} and Remark~\ref{rmk:subcritical} relies fundamentally on the tightness of the approximate solutions and their associated covariation measures. Consequently, extending this approach to the supercritical regime, where $\lv$ is replaced by $ 2\pi \lambda_0 / \log \vep^{-1}$ for $\lambda_0 \in (1, \infty)$, raises immediate questions about its mathematical validity. Moreover, although the same formal calculation yields \eqref{formaleqv0} in this setting, taking formal expectations to eliminate the stochastic integral forces the condition $\E[X(x,t)^2]=0$, as the alternative leads to the contradictory inequality $[P_tX_0(x)]^2<0$. Rather than indicating an absolute barrier, this degenerate limiting behavior may simply reflect an intrinsic limitation of applying these solutions to the supercritical case. For precise asymptotic expansions of second moments of approximate solutions approached via time-changed additive functionals, see \cite{C:Additive}; for a separate study establishing local extinction and collapse of mass via discrete polymers and stochastic heat flow, see \cite{BCT:25}. \hfill $\blacklozenge$
\end{rmk}

\subsection{Precise equations for the limiting covariation measures}\label{sec:smt-se}
In this subsection, we present the precise statements of our main theorems (Theorems~\ref{thm:main1}, \ref{thm:main2}, and \ref{thm:main3}). The proof organization is provided at the end.

The first step in making the heuristic argument in Section~\ref{sec:smt-expansion} rigorous is to establish process-level weak convergence for the approximate solutions $X_\vep$ defined in \eqref{def:SHE:vep}, reformulated as measure-valued processes, along with their covariation measures. To achieve this goal, our formulation below highlights three key features:
\begin{itemize}
\item We temper these measures at infinity so that a suitable state space exists on which the Skorokhod representation theorem applies, thereby greatly streamlining the subsequent asymptotic expansions. 

\item Alongside weak convergence, we derive uniform-in-$\varepsilon$ moment bounds; these extend naturally to the limiting process, which is the version stated in the theorem below for brevity. While similar weak convergence results and moment bounds appear in Nakashima \cite{Nakashima:25}, our approach differs by utilizing the continuous approximation model \eqref{def:SHE:vep} rather than a fully discrete one. Related results also appear in Tsai \cite{Tsai:24}. Furthermore, the moment bounds established here involve new formulations essential to our main theorems and yield improved H\"older exponents.

\item We assume that the initial condition $X_0(\cdot)$ is bounded and nonnegative. Nonnegativity is imposed primarily for analytical convenience. The main theorems presented in this subsection remain fully valid without it, with the sole exception of \eqref{MMbdd} and \eqref{X4BMbdd}. Because \eqref{MMbdd} and \eqref{X4BMbdd} serve only as auxiliary bounds, their absence does not affect the extensions of our main formulas. We elaborate on this in Remark~\ref{rmk:IC}.
\end{itemize}

Theorem~\ref{thm:main1} offers a minimal summary of our main results from the first step. In this theorem, $C_{E}[0,\infty)$ denotes the space of continuous functions from $[0,\infty)$ into a topological space $E$, equipped with the topology of uniform convergence on compact sets. Additionally, $\mathcal M_f(E)$ refers to the space of (nonnegative) finite measures on $E$, endowed with the topology of weak convergence. The pre-limit processes are then taken to be the following $\mathcal M_f(\R^2)$- and $\mathcal M_f(\R^4)$-valued processes, constructed from the approximate solutions $X_\vep$ to \eqref{def:SHE:vep} and the associated martingale measure $M_\vep$ defined in \eqref{def:M:vep}:
\begin{gather}
X^\sharp_\vep(\d x,t)\,\defeq\, \psi(x)X_\vep(\d x,t)=\psi(x)X_\vep(x,t)\d x,\quad x\in \R^2,\label{def:Xsharp1}\\
\la M_\vep,M_\vep\ra^\sharp(\d x',\d x,\d t)\,\defeq\, \psi^{\otimes 2}(x',x)\la M_\vep,M_\vep\ra (\d x',\d x,\d t),\quad (x',x,t)\in \R^2\times \R^2\times \R_+.\label{def:Mqvsharp1}
\end{gather}
Here, $\psi$ is an auxiliary function introduced for convenience. Its explicit form, along with the associated tensor product $\psi^{\otimes 2}$, is given below:
\begin{align}\label{def:psi}
\psi(x)\,\defeq\,(1+\lvert x\rvert^{10})^{-1}, \quad \psi^{\otimes 2}(x',x)\,\defeq\,\psi(x')\psi(x).
\end{align} 
Consequently, $X_\vep^\sharp(\d x,0)$ is a finite measure, even though we only assume that $X_0(\cdot)$ is bounded. Additionally, the covariation measure of $M_\vep$ refers to the unique measure $\la M_\vep,M_\vep\ra (\d x',\d x,\d t)$ on $\B(\R^2\times \R^2\times \R_+)$ such that, for any $ f\in \C^0_{pd}(\R^2)$, 
\[
\biggl(\int_0^T\int_{\R^2}f(x)M_\vep(\d x,\d t)\biggr)^2-\int_0^T\int_{\R^2}\int_{\R^2}f(x')f(x)\la M_\vep,M_\vep\ra (\d x',\d x,\d t)
\]
is a continuous martingale. Finally, beyond the weak formulation \eqref{def:weak:vep} that expresses $M_\vep$ in terms of $X_\vep$, the following theorem also considers the mild formulation (specifically, the Duhamel form) as well as the mild formulation integrated against a test function, as captured by the next two equations:
\begin{subequations}\label{def:mild:vep}
\begin{align}
X_\vep(x',T)
&=P_TX_0(x')+\int_{0}^T\int_{\R^2} P_{T-t}(x',x)M_\vep(\d x,\d t),\quad x'\in \R^2,\\
\begin{split}
X_\vep(f,T)
&=\int_{\R^2} f(x)P_TX_0(x)\d x+\int_{0}^T \int_{\R^2} P_{T-t}f(x)M_\vep(\d x,\d t),\quad f\in \C_{pd}^0(\R^2),
\end{split}
\end{align} 
\end{subequations}
Here, $\C_{pd}^k(\R^2)$ denotes the subspace of $\C^k(\R^2)$ specified at the end of Section~\ref{sec:intro}. 

\begin{thm}\label{thm:main1}
Fix $\lambda\in \R$, and let $\varrho \in \C_c(\R^2)$ be a mollifier. Recall that $\phi\in \C_c(\R^2)$ is the probability density defined in \eqref{def:phi}. Given a bounded nonnegative initial condition $X_0(\cdot)$ independent of $\vep$, the solutions $X_\vep$ to \eqref{def:SHE:vep} satisfy the following properties:  \medskip 

\noindent {\rm (1$\cc$)} The family of laws of $\{(X^\sharp_\vep(\d x,t), \langle M_\vep,M_\vep\rangle^\sharp(\d x',\d x,(0,t])\}_{\vep>0}$ converges in law as probability measures on the Polish space $C_{\mathcal M_f(\R^2)}[0,\infty)\times C_{\mathcal M_f(\R^2\times \R^2)}[0,\infty)$ as $\vep\to 0$. 

Let $(X^\sharp(\d x,t),\mu^\sharp(\d x',\d x,(0,t]))$ denote the limiting process, and define
\begin{gather}
X(\d x,t)\,\defeq\, \psi(x)^{-1}X^\sharp(\d x,t),\quad x\in \R^2,\label{def:Xns}\\
X(f,t)\,\defeq\,\int_{\R^2} f(x)X(\d x,t),\quad f\in \C^0_{pd}(\R^2),\notag\\
 \mu(\d x',\d x,\d t)\,\defeq\, \psi^{\otimes 2}(x',x)^{-1}\mu^\sharp(\d x',\d x,\d t),\quad (x',x,t)\in \R^2\times \R^2\times \R_+.\label{def:muns}
\end{gather}
Then $M$ is uniquely determined as an orthogonal martingale measure by the relation
\begin{align}
X(f,T)&=X(f,0)+\int_0^T X\left(\frac{\Delta f}{2},t\right)\d t+\int_{0}^T \int_{\R^2} f(x)M(\d x,\d t),\quad f\in \C^2_{pd}(\R^2),\label{def:weak}
\end{align}
where $\mu$ and the covariation measure $\la M,M\ra(\d x,\d t)$ satisfy
\begin{align}\label{Mqv:det}
\int_0^T\int_{\R^2}g(x,x) \la M,M\ra(\d x,\d t)
=\int_0^T \int_{\R^2}\int_{\R^2}g(x',x)\mu(\d x',\d x,\d t),\;\; g\in \C^0_{pd}(\R^2\times \R^2).
\end{align}
Moreover, the following weak formulation of the mild form holds:
\begin{align}
X(f,T)&=\int_{\R^2}f(x)P_TX_0(x)\d x+\int_{0}^T \int_{\R^2} P_{T-t}f(x)M(\d x,\d t),\quad f\in \C^0_{pd}(\R^2).\label{def:mild}
\end{align}

\noindent {\rm (2$\cc$)} The covariation measure $\la M,M\ra(\d x,\d t)$ satisfies the following measure identity:
\begin{align}
\E[\la M,M\ra(\d x,\d t)] = \rDBG_t X_0^{\otimes 2}(x)\d x\d t,\quad (x,t)\in \R^2\times \R_+,\label{m1:apriori}
\end{align}
where $X_0^{\otimes 2}(z_1,z_2) \defeq X_0(z_1)X_0(z_2)$, and $\rDBG_t F(x)$ is defined in \eqref{def:rDBG}. Furthermore, the following measure inequalities hold for all $0 \le T_1 \le T_2 \le T < \infty$, $p \in (0,1)$, $a, a' \in (0,T]$, and $y, y' \in \R^2$:
\begin{align}
&\int_{t=T_1}^{T_2} \int_{t'=t}^{T_2} \E[\la M,M\ra(\d x',\d t')\la M,M\ra(\d x,\d t) ] \d t' \d t \notag\\
&\qquad \le C(\lambda,\phi,\lVert X_0\rVert_\infty,T,p) (T_2-T_1)^{2-p}\d x'\d x,\quad (x',x)\in \R^2\times \R^2,\label{MMbdd}\\[1ex]
&\int_{\R^2}\int_{\R^2} P_{a'}(y',z') P_{a}(y,z)\E[X_\infty(\d z' ,t)X_\infty(\d z,t)\la M_\infty,M_\infty\ra(\d x,\d t)] \notag\\
&\qquad \le \frac{C(\lambda,\phi,\lVert X_0\rVert_\infty,T,p)}{\min\{a',a\}^p}\d x\d t,\quad (x,t)\in \R^2\times (0,T]. \label{X4BMbdd}
\end{align}
\end{thm}

\begin{rmk}
As in the case of $M_\vep$, the following martingale property provides a characterization of $\la M,M\ra(\d x,\d t)$: for any $ f\in \C^0_{pd}(\R^2)$,
\[
\biggl(\int_0^T\int_{\R^2}f(x)M(\d x,\d t)\biggr)^2-\int_0^T\int_{\R^2}f(x)^2\la M,M\ra (\d x,\d t)
\]
is a continuous martingale.
\hfill $\blacklozenge$
\end{rmk}

The following immediate corollary of Theorem~\ref{thm:main1} reproduces the known second-moment formulas, which were first established in \cite{BC:2D} in a different yet equivalent form.
 
\begin{cor}\label{cor:2ndmom}
For all $f_1,f_2\in \C_{pd}^0(\R^2)$, we have
\begin{align*}
\E\biggl[\prod_{j=1}^2 X_T(f_j,T)\biggr]
&=\prod_{j=1}^2\int_{\R^2}
P_TX_0(w_j)f_j(w_j)\d w_j +\int_0^T\int_{\R^2}\prod_{j=1}^2
P_{T-t}f(x)
\rDBG_{t}X_0^{\otimes 2}(x)
\d x\d t.
\end{align*}
\end{cor}
\begin{proof}
The required identity is obtained more directly by applying the moment duality \eqref{Mom:dual} (see the proof of \eqref{m1:apriori}), although it can also be derived a posteriori from Theorem~\ref{thm:main1} as follows. Expressing $X(f_j, T)$ as the sum of a deterministic component and a stochastic integral via the mild formulation \eqref{def:mild}, we expand $\prod_{j=1}^2 X(f_j, T)$ and take expectations. By \eqref{MMbdd}, the cross term vanishes in expectation. Consequently, since the covariation measure determines the quadratic variation of products of stochastic integrals with respect to $M(\mathrm{d}x, \mathrm{d}t)$, we get
\begin{align*}
\E\biggl[\prod_{j=1}^2 X(f_j,T)\biggr]
=\prod_{j=1}^2\int_{\R^2}P_TX_0(w_j)f_j(w_j)\d w_j+
\E\biggl[\int_{0}^T\int_{\R^2}\prod_{j=1}^2P_{T-t}f_j(x)\la M,M\ra(\d x,\d t)\biggr].
\end{align*}
The required equality follows upon applying \eqref{m1:apriori} to the last expectation.
\end{proof}

Our proof of weak convergence in Theorem~\ref{thm:main1} (1$\cc$) aims specifically to establish tightness. The unique determination of $\mu$ in \eqref{Mqv:det} rests on the fact that the quadratic variation of a continuous local martingale $N$ is the unique continuous, finite variation process $A$ with zero initial condition for which $N^2-A$ remains a continuous local martingale. Therefore, establishing uniqueness for the limit in Theorem~\ref{thm:main1} reduces to identifying the limit law of $X_\vep$ as $\vep \to 0$, which follows from Tsai’s axiomatic characterization~\cite{Tsai:24}. Independently of this full convergence, all proofs in this paper applies along any convergent subsequence $\{(X^\sharp_{\vep_n}(\d x, t), \la M_{\vep_n}, M_{\vep_n} \ra^\sharp(\d x', \d x, (0,t]))\}_{n\in\mathbb{N}}$ as $\vep_n \to 0$. Thus, our results apply to the corresponding limiting process, denoted by $\{(X^\sharp_{\infty}(\d x, t), \la M_{\infty}, M_{\infty} \ra^\sharp(\d x', \d x,(0,t]))\}$, which we use throughout the proofs and for future reference. Analogues of the processes $(X, \mu)$ defined in \eqref{def:Xns} and \eqref{def:muns} can also be identified. In particular, the limiting measure $\mu_\infty$ associated with $\mu$ is uniquely determined as the covariation measure $\la M_\infty, M_\infty \ra$ of the orthogonal martingale measure $M_\infty$ of the limiting solution $X_\infty$. 

The following theorem offers a rigorous counterpart to the formal asymptotic expansions discussed in Section~\ref{sec:smt-expansion}, using a weak formulation to absorb pointwise singularities. 

\begin{thm}\label{thm:main2}
The covariation measure $\langle M, M \rangle $ of the limiting martingale measure $M$ in Theorem~\ref{thm:main1} satisfies the following identity for all $T \in (0, \infty)$ and $f(x, t) \in \C_{pd}^1(\R^2 \times [0, T])$:
\begin{align}
&\int_{t=0}^T\int_{\R^2}\mathfrak Lf(x,t,T)\langle M,M\ra(\d x,\d t)\notag\\
&\quad =\int_{t=0}^T\int_{\R^2}f(x,t)[P_tX_0(x)]^2\d x\d t\notag\\
&\quad \quad +2\int_{t'=0}^T\int_{\R^2} \int_{t=t'}^T \int_{\R^2} f(x,t)\int_{\R^2}\biggl(
\prod_{j=1}^2P_{t-t'}(x,z_j)\biggr) X(\d z_2,t')\d x\d t \M(\d z_1,\d t').\label{mild:goal1}
\end{align}
Here, $\mathfrak L$ refers to the following integro-multiplication operator:
\begin{align}
\mathfrak Lf(x,t,T)&\,\defeq-\frac{1}{2\pi}
\biggl(\frac{\log (T-t)}{2}+\log 2+\lambda -\frac{\gamma_{\sf EM}}{2}-\int_{\R^2}\int_{\R^2}\phi(y')\phi(y)\log \rvert y'-y\rvert \d y'\d y\biggr)f(x,t)\notag\\
&\quad\,-\int_{\R^2}\int_{t'=0}^{T-t}[f(x-x',t+t')-f(x,t)]P_{t'}(x')P_{t'}(x')\d t'\d x',\label{def:L}
\end{align}
where $(\lambda,\phi)$ is as given in \eqref{def:lv} and \eqref{def:phi}, and $\EM$ denotes the Euler--Mascheroni constant. All the stochastic integrals and Lebesgue integrals in \eqref{mild:goal1} and \eqref{def:L} are well-defined. 
\end{thm}

By formulating the integro-multiplication operator $\mathfrak L$, Theorem~\ref{thm:main2} identifies the crucial second-order term absent from formal expansions. However, the presence of $\mathfrak L$ also accounts for why Theorem~\ref{thm:main2} does not explicitly capture the covariation measure $\langle M, M \rangle(\d x, \d t)$. Furthermore, the analytical properties of $\mathfrak L$, discussed in more detail below, make it difficult to invert directly when deriving explicit formulas for $\langle M, M \rangle(\d x, \d t)$, so the next theorem manages to bypass this technical issue. Ultimately, owing to its direct link to the expansions in Section~\ref{sec:smt-expansion}, Theorem~\ref{thm:main2} is essential for deriving our full formulas and provides a general formulation that may prove useful for studying related SPDEs.

The individual components of $\mathfrak L$ reveal both classical and new characteristics of the two-dimensional SHE at criticality. On the classical side, the derivation of $\mathfrak L$ (Sections~\ref{sec:I1} and \ref{sec:I2}) closely parallels one step in deriving the asymptotic Green function expansion recalled in \eqref{RM:Green} (see especially the final display in \cite[Section~4.4]{C:Additive}). In particular, the multiplication symbol in $\mathfrak L$ nearly reproduces the constant $(\log \beta)/2$ from \eqref{def:betaintro}, where $\beta$ parametrizes the second moments at criticality. The new characteristics stem from the logarithmic multiplier $-(4\pi)^{-1}\log (T-t)$ and the {\bf integral part} (the iterated integral term in \eqref{def:L} with respect to $\d t' \d x'$). While this logarithmic multiplier may account for the remaining $-\EM/2$ in \eqref{def:betaintro}, using the formula $\EM = -\int_0^\infty \e^{-t}\log t\d t$ \cite[(1.3.19), p.~8]{Lebedev} (see also Remark~\ref{rmk:Green}), it is inherently time-reversed (appearing with $T-t$). This time reversal complicates the picture and highlights the essential role of the integral part. Further delicate properties of $\mathfrak L$ are discussed in Section~\ref{sec:lcd2-ime}.

The following theorem gives a full account of the covariation measure $\la M, M \ra(\d x, \d t)$, thereby completely determining the semimartingale characteristics of the two-dimensional SHE at criticality via \eqref{def:weak} and \eqref{def:mild}. 

\begin{thm}\label{thm:main3}
Let $\rDBG_t(x;z_1,z_2)$ and $\s^\beta(\cdot)$ be as defined in \eqref{def:rDBG} and \eqref{def:sbeta}, respectively. The covariation measure $\la M,M\ra(\d x,\d t)$ of the martingale measure $M$ described in Theorem~\ref{thm:main1} admits the following formulas: for all $T\in (0,\infty)$ and
 $f(x,t)\in\C^0_{pd}(\R^2\times [0,T])$, 
\begin{align}
& \int_{t=0}^{T}\int_{\R^2} f(x,t)\langle M,M\rangle (\d x,\d t)\notag\\
&\quad =  \int_{t=0}^{T}\int_{\R^2}
 f(x,t)\rDBG_{t}X^{\otimes 2}_0(x)\d x\d t\notag\\
&\quad \quad +2\int_{t'=0}^T\int_{\R^2} \int_{t=t'}^T
 \int_{\R^2}
 \int_{\R^2}
f(x,t)
\rDBG_{t-t'}(x;z_1,z_2) X(\d z_2,t') \d x \d t M(\d z_1,\d t'), \label{eq:main3-1}
\end{align}
and
\begin{align}
&\int_{t=0}^{T}\int_{\R^2} f(x,t)\langle M,M\rangle (\d x,\d t)\notag\\
&\quad  =  \int_{t=0}^{T}\int_{\R^2}
 f(x,t)\rDBG_{t}X^{\otimes 2}_0(x)\d x\d t\notag\\
&\quad \quad +2\int_{t=0}^T\int_{t'=0}^t\int_{\R^2} 
 \int_{\R^2}\int_{\R^2}
f(x,t)\rDBG_{t-t'}(x;z_1,z_2) X(\d z_2,t') \d x M(\d z_1,\d t')\d t, \label{eq:main3-2}
\end{align}
where 
\[
X_0^{\otimes }(z_1,z_2)\,\defeq\, X_0(z_1)X_0(z_2). 
\]
Additionally, for all $f(x)\in \C^0_{pd}(\R^2)$, we have
\begin{align}
&\int_{t=0}^{T}\int_{\R^2} f(x)\langle M,M\rangle (\d x,\d t)\notag\\
&\quad =  \int_{t=0}^{T}\int_{\R^2}
 f(x)\rDBG_{t}X_0^{\otimes 2}(x)\d x\d t\notag\\
&\quad \quad  +2 \int_{s=0}^{T}\int_{\R^2}P_{s/2}f(y)
\s^\beta(s) \int_{t'=0}^{T-s}\int_{\R^2}
\int_{\R^2}
\int_{t=s}^{T-t'} \prod_{j=1}^2
P_{t-s}(y,z_j) \notag\\
&\quad \quad \times \d t
X(\d z_2,t')  M(\d z_1,\d t')\d y\d s.\label{eq:main3-3}
\end{align}
All the stochastic and Lebesgue integrals in \eqref{eq:main3-1}--\eqref{eq:main3-3} are well-defined.
\end{thm}

\begin{rmk}\label{rm:rDBG}
(1$\cc$) Feynman--Kac-type representations involving $\rDBG_t(x;z_1,z_2)$, along with related analytic properties, have been developed in \cite{C:SDBG1,C:BES}.\medskip 

\noindent (2$\cc$) The bound in \eqref{X4BMbdd} provides a key tool for extending \eqref{eq:main3-1}--\eqref{eq:main3-3} to more general test functions. However, \eqref{eq:main3-3}  should remain restricted to time-independent test functions.
\hfill $\blacklozenge$
\end{rmk}

Theorem~\ref{thm:main3} highlights both differences and connections with several established results, some of which were discussed in Section~\ref{sec:intro}. We also emphasize the following observations:
\begin{itemize}
\item Equation~\eqref{main:intro} reflects the quadratic structure of the covariation measure in the solution $X$, analogous to \eqref{SHEqv}, via $X_0^{\otimes 2}$ and the integrator $X(\d z_2, t') M(\d z_1, \d t')$. 

\item Combined with the weak formulation (see \eqref{def:weak}), our formulas for the covariation measures as in \eqref{eq:main3-2} provide precise descriptions of the martingale problem for the two-dimensional SHE at criticality. 
However, in standard settings, such as finite-dimensional SDEs or Dawson--Watanabe superprocesses, the semimartingale characteristics are specified by non-random functions or functionals of the state. In contrast, specifying the coefficients in our setting inherently requires It\^o integration with respect to the martingale measure $M(\d z_1, \d t')$, reflecting a covariation structure that is fundamentally stochastic and path-dependent.
\end{itemize}

A further aspect of Theorem~\ref{thm:main3} concerns the singularities of the kernel $\rDBG_t(x;z_1,z_2)$, which render changes in the order of integration delicate. Consequently, although formulas \eqref{eq:main3-1}--\eqref{eq:main3-3} differ only in the order of integration of their final terms, our analysis shows that any such rearrangement requires thorough justification. In particular, this motivates the restriction on \eqref{eq:main3-3} highlighted in Remark~\ref{rmk:rDBG} (2$\cc$). If, instead, one assumes full interchangeability of the order of integration without rigorous justification, the covariation measure $\langle M, M \rangle(\d x, \d t)$ appears to have the density
\begin{align}
\rDBG_tX_0^{\otimes 2}(x)+2\int_{t'=0}^t\int_{\R^2}\int_{\R^2} \rDBG_{t-t'}(x;z_1,z_2)X(\d z_2,t')M(\d z_1,\d t')\label{eq:qvM1}
\end{align} 
with respect to the Lebesgue measure $\d x\d t$. Standard arguments would then imply that the martingale $\int_{t=0}^T\int_{\R^2}f(x)M(\d x,\d t)$ in the weak formulation \eqref{def:weak} can be represented as a stochastic integral with respect to a space-time white noise. However, beyond the technical obstacles that the current analytical bounds encounter in justifying this interchange, the presence of space-time white noise as the driving noise seems overly idealized given the singular nature of the two-dimensional SHE at criticality, making its validity implausible. See~\cite{GT:25} for a rigorous proof of this obstruction. 

Although formal, the density in \eqref{eq:qvM1} provides an analogue of the right-hand side of \eqref{eq:Y3} for the one-dimensional SHE, replacing, in particular, the product of heat kernels with the special kernel $\rDBG_t(x;z_1,z_2)$. This transition from the Gaussian regime to the delta-Bose gas, along with that from the subcritical additive SHE discussed in Section~\ref{sec:intro-overview}, further illustrates the heuristics introduced in \cite{C:SDBG2,C:BES}. These heuristics reflect the fundamental interplay between the Laplacian and noise at criticality in two dimensions, motivating the Feynman--Kac formulas for the two-body delta-Bose gas and the related analytical properties established in those works.

Some additional connections between the formulas in Theorem~\ref{thm:main3} and existing literature are outlined below:
\begin{itemize}
\item The explicit semimartingale characteristics of the two-dimensional SHE at criticality also bear a close structural resemblance to those found in the superprocess within a random environment constructed by Mytnik~\cite{Mytnik:96}, despite the significant differences in frameworks and methods between \cite{Mytnik:96} and the present work. One of the main theorems from \cite{Mytnik:96} establishes the existence of $Z\in C_{\mathcal M_f(\R^d)}[0,\infty)$ for any dimension $d \geq 1$, such that for every $f \in \C_c^2(\R^d)$, the process
\begin{align}\label{eq:Mytnik}
M_Z(f,T)\,\defeq\, Z(f,T)-Z(f,0)-\int_0^T Z\left(\frac{\Delta f}{2},t\right)\d t
\end{align}
is an $L^2$-martingale, where $Z(f,T)\,\defeq\,\int_{\R^d}f(x)Z(\d x,T)$, and $\Delta$ now refers to the $d$-dimesnsional Laplacian. Equation~\eqref{eq:Mytnik}, in turn, induces a martingale measure given by $M_Z(\d x, \d t)$, so that $M_Z(f, T)=\int_{t=0}^T\int_{\R^d}f(x)M_Z(\d x,\d t)$. Furthermore, the covariation measure of $M_Z(\d x, \d t)$ is non-orthogonal and given by the following formula: for all nonnegative functions $g\in \C_b(\R^d \times \R^d)$,
 \begin{align}
&\int_{t=0}^T\int_{\R^d}\int_{\R^d}g(x',x)\langle M_Z,M_Z\mathcal \ra(\d x',\d x,\d t)\notag\\
&\quad =\int_{t=0}^T\int_{\R^d}g(x,x)Z(\d x, t)\d t  +\int_{t=0}^T\int_{\R^d}\int_{\R^d}g(x',x) \kappa(x',x)
Z(\d x', t)
Z(\d x,t)
\d t.\label{covar:Mytnik}
\end{align}
Here, $\kappa(x',x)$ is a deterministic, bounded, and continuous function characterizing this process $Z$. Notably, the original results in \cite{Mytnik:96} are formulated in far greater generality, extending well beyond the specific triple $(\frac{\Delta}{2},\R^d,\C_c^2(\R^d))$ presented here.

The close structural resemblance mentioned above emerges when comparing \eqref{eq:main3-2} and \eqref{covar:Mytnik}, though the martingale measure $M_Z(\d x,\d t)$ is not orthogonal. For \eqref{covar:Mytnik}, taking $\kappa \equiv 0$ reduces $Z$ to the standard superprocess, while the random environment generates an additional quadratic term involving $Z(\d x',t)Z(\d x,t)$. This mirrors the right-hand side of \eqref{eq:main3-2}: the first term alone defines a Gaussian covariation measure, as in the subcritical regime, but at criticality, noise additionally produces a quadratic term involving $X(\d z_2,t')M(\d z_1,\d t')$. 

\item Nakashima~\cite{Nakashima:25} recently provided an alternative description of the covariation measure for the martingale measure in the critical two-dimensional stochastic heat flow, along with additional results discussed in Section~\ref{sec:intro}. To restate this description from \cite{Nakashima:25}, let $\SHF'(\d x,0;\d y,t)$ represent the limiting stochastic heat flow from \cite{Nakashima:25}, and define
\[
X'(f,T)\,\defeq\, \int_{\R^2}\int_{\R^2}X_0(x)f(y)\SHF'(\d x,0;\d y,T),
\]
with $M'(f)=M'(f,T)$ denoting the martingale part. (We do not directly identify $X'$ with the limiting process $X$ from Theorem~\ref{thm:main1}, as parameter identification remains unresolved in the corresponding limiting schemes; see \cite[Remark~1.10]{Nakashima:25}.) Specifically, Nakashima \cite{Nakashima:25} proves that for all nonnegative $X_0(\cdot) \in \C_c(\R^2)$, $f \in \C^2_b(\R^2)$, and $0 < T_0 < \infty$, the following convergence in probability holds: 
\begin{align}\label{SHEqv:approx}
&\sup_{0\leq T\leq T_0}\left|\frac{4\pi}{\log\eta^{-1}}\int_0^T\int_{\R^2}\big[X'(P_\eta(x,\cdot),t\big) f(x)\big]^2\d x\d t -\langle M'(f),M'(f)\rangle_T\right|\xrightarrow[\eta\searrow 0]{\P}0,
\end{align}
where $P_\eta(x,y)$, as defined in \eqref{def:Pt}, denotes the transition probability of two-dimensional standard Brownian motion.

The descriptions in Theorem~\ref{thm:main3} and \eqref{SHEqv:approx} differ in key analytical aspects. In particular, \eqref{SHEqv:approx} provides an approximation scheme for the covariation measure, whereas Theorem~\ref{thm:main3} yields an exact expression that explicitly incorporates the dependence on the parameter $\beta$ from \eqref{def:betaintro}.
\end{itemize}

The following theorem presents an analogue of Nakashima's description discussed above. We provide a brief sketch of the proof, as the technical argument requires only straightforward modifications of the methods developed later.

\begin{thm}\label{thm:qvsq}
Assume that $X_0(\cdot)\in \C_b(\R^2)$, and let $X$ denote the associated limiting process from Theorem~\ref{thm:main1}. Then for all $f\in \C_{pd}^1(\R^2)$ and $0<T_0<\infty$, 
\begin{align*}
&\sup_{0\leq T\leq T_0}\biggl|\frac{4\pi}{\log\eta^{-1}}\!\int_{0}^T\int_{\R^2}\big[X\big(P_\eta(x,\cdot),t\big) f(x)\big]^2\d x\d t-\int_{0}^T\int_{\R^2} f(x)^2\langle M,M\rangle(\d x,\d t)\bigg|\xrightarrow[\eta\searrow 0]{\P}0,
\end{align*}
where the limit denotes convergence in probability. 
\end{thm}
\begin{proof}[Sketch of the proof of Theorem~\ref{thm:qvsq}] 
It is enough to show that the convergence in probability under consideration holds for each fixed $0\leq T\leq T_0$, rather than uniformly for all $0\leq T\leq T_0$. To justify this weaker requirement, we use the following fact along with a Helly-selection-type argument: if a sequence of non-decreasing functions $f_n:[a,b]\to \R$ converges pointwise to a continuous function $f$, then the convergence is actually uniform.

The central step in this proof is to demonstrate that, for any fixed $T > 0$, the left-hand side of the following equation converges in probability to the required limit as $\eta \to 0$:
\begin{align}\label{sqqv:1}
\frac{4\pi}{\log\eta^{-1}}\!\int_{t=0}^T\int_{\R^2}X\big(P_\eta(x,\cdot),t\big)^2 f(x)^2\d x\d t
=\sum_{k=1}^3 I^k_{\ref{appsq}}(T),
\end{align}
where we define
\begin{subequations}
\label{appsq}
\begin{align}
I^1_{\ref{appsq}}(T)&\,\defeq\,
\frac{4\pi}{\log\eta^{-1}}\int_{t=0}^T \int_{\R^2}
\langle P_\eta(x,\cdot),P_tX_0\ra^2 f(x)^2\d x\d t,\\
I^2_{\ref{appsq}}(T)&\,\defeq\,\frac{8\pi}{\log \eta^{-1}}\int_{t=0}^T\int_{\R^2}
\int_{t'=0}^t\int_{\R^2}X\big(P_{t-t'+\eta}(x,\cdot),t'\big) P_{t-t'+\eta}(x,x')\notag\\
&\quad \times M(\d x',\d t')f(x)^2\d x\d t,\\
I^3_{\ref{appsq}}(T)&\,\defeq\, \frac{4\pi}{\log \eta^{-1}}\int_{t=0}^T\int_{\R^2}\int_{t'=0}^t \int_{\R^2}P_{t-t'+\eta}(x,x')^2\la M,M\ra (\d x',\d t') f(x)^2\d x\d t.
\end{align}
\end{subequations}
Here, the decomposition in \eqref{sqqv:1} is derived by first applying the smooth test function $P_{\eta}(x,\cdot)$ to \eqref{def:mild}, which produces the following analogue of \eqref{asymp:2}: 
\begin{align*}
X\big(P_\eta(x,\cdot),t\big)^2&=\langle P_\eta(x,\cdot),P_tX_0\ra^2 +2\int_{t'=0}^t\int_{\R^2}X(P_{t-t'+\eta}(x,\cdot),t')
P_{t-t'+\eta}(x,x')M(\d x',\d t')\notag\\
&\quad  +\int_{t'=0}^t \int_{\R^2}P_{t-t'+\eta}(x,x')^2\la M,M\ra (\d x',\d t').
\end{align*}
Then, integrating both sides of this equation against $(4\pi/\log \eta^{-1})f(x)^2 \d x\d t$ for $(x,t)\in \R^2\times [0,T]$ yields \eqref{sqqv:1}.

To complete the proof, it suffices to show that for all $1 \leq j \leq 2$, $I^j_{\ref{appsq}}(T)$ converges to zero in probability as $\eta \searrow 0$, ensuring that the required limit is determined by the behavior of $I^3_{\ref{appsq}}(T)$. The argument for $ j =2$ involves first justifying the use of the stochastic Fubini theorem, so that $\la M, M \ra (\d x', \d t') \d x \d t$ in $I^2_{\ref{appsq}}(T)$ can be interchanged with $\d x \d t\la M, M \ra (\d x', \d t')$ for any \emph{fixed} $\eta \in (0,1)$. Here, the verification uses Walsh's $L^2$-condition \cite[Theorem~2.6, pp.296--297]{Walsh} and the moment bounds from Theorem~\ref{thm:main1}. These moment bounds are also essential for proving the required zero limits. To identify the limit in probability of $I^3_{\ref{appsq}}(T)$, we consider the following calculation:
\begin{align*}
\lim_{\eta\searrow 0}I^3_{\ref{appsq}}(T)
&=\lim_{\eta\searrow 0}\frac{4\pi}{\log \eta^{-1}}\int_{t'=0}^T \int_{\R^2}\biggl(\int_{t=t'}^T \int_{\R^2} P_{t-t'+\eta}(x,x')^2\d x\d t\biggr)f(x')^2 \la M,M\ra (\d x',\d t')\\
&=\lim_{\eta\searrow 0}\int_{t'=0}^T \int_{\R^2}\biggl(\frac{4\pi}{\log \eta^{-1}}\frac{\log (T+\eta)-\log \eta}{4\pi}\biggr)f(x')^2 \la M,M\ra (\d x',\d t')\\
&=\int_{t'=0}^T\int_{\R^2}f(x')^2 \la M,M\ra(\d x',\d t'),
\end{align*}
where the first equality replaces $f(x)^2$ by $f(x')^2$, as justified by \eqref{ineq:heat2-1}. This limit suffices to establish the theorem. 
\end{proof}

We conclude this discussion of the main theorems with a remark on the use of signed initial conditions and their connection to the critical two-dimensional stochastic heat flow:

\begin{rmk}[Signed initial conditions and SHF]\label{rmk:IC}
The main results presented above, with the sole exception of the auxiliary bounds \eqref{MMbdd} and \eqref{X4BMbdd}, naturally extend to general signed bounded initial conditions (see (1$\cc$) below). Further extensions apply to the stochastic heat flow setting discussed in Section~\ref{sec:intro-context}; see (2$\cc$) for the precise setting.\medskip

\noindent (1$\cc$) To extend the exact formulas of covariation measures to a general signed, bounded initial condition $X_0$, consider a polar decomposition approach using three solutions: $X^+$, $X^-$, and $\ol{X}\,\defeq\,X^++X^-$. They correspond to the bounded nonnegative initial conditions $X_0^+$, $X_0^-$, and $\ol{X}_0 =\lvert X_0\rvert  =X_0^+ + X_0^-$, respectively, where $X_0^\pm$ denote the positive and negative parts of $X_0$. The associated martingale measures are denoted by $M^+$, $M^-$ for $X^+$ and $X^-$, and the martingale measure $\ol{M}$ of $\ol{X}$ satisfies $\ol{M}=M^++M^-$. By linearity, the solution and martingale measure for the two-dimensional SHE at criticality with initial condition $X_0$ are recovered in the following way:
\[
X\;\defeq\, X^+-X^-,\quad M\;\defeq\, M^+-M^-. 
\]
Then $M$ satisfies \eqref{m1:apriori} since by bilinearity, 
 \begin{align}\label{MM:sign}
\langle M,M\ra(\d x,\d t)
 &= 2\sum_{\varsigma\in \{\pm\}} \la M^\varsigma,M^\varsigma\ra(\d x,\d t) -\la \ol{M},\ol{M}\ra(\d x,\d t).
\end{align}
We also have extensions of Theorems~\ref{thm:main2} and \ref{thm:main3} to the pair $(X,M)$ under consideration since
\begin{align}\label{XM:sign}
X(\d z_2,t')M(\d z_1,\d t')
&=2\sum_{\varsigma\in \{\pm\}}X^\varsigma(\d z_2,t') M^\varsigma(\d z_1,\d t')-\ol{X}(\d z_2,t')\ol{M}(\d z_1,\d t').
\end{align}
By virtue of \eqref{MM:sign} and \eqref{XM:sign}, directly extending \eqref{MMbdd} and \eqref{X4BMbdd} to $(X,M)$ is unnecessary for applications to the exact formulas for covariation measures, thereby avoiding the need to treat the resulting cross-terms directly. 

To make the above argument rigorous, note that the algebraic identities in \eqref{MM:sign} and \eqref{XM:sign} hold analogously for the approximate solutions $(X_\vep, M_\vep)$, whose tightness is justified by \cite[Proposition~2.4, p.~107]{EK:MP}, as well as the tightness behind Theorem~\ref{thm:main1} (1$\cc$). Consequently, tightness of the approximations alone is sufficient to pass to the exact formulas.\medskip

\noindent (2$\cc$) The present framework applies to any countable set of initial times $\mathcal S \subset [0, \infty)$ rather than restricting to $s=0$. Given any collection of bounded initial conditions $\{X_s^k\}_{k \in \mathbb{N}, s \in \mathcal S}$, the tightness result in Theorem~\ref{thm:main1} extends to the corresponding families of approximate solutions $\{X^k_{s,t,\varepsilon}\}_{k \in \mathbb{N}, t \ge s}$ constructed on a shared probability space, again justified by \cite[Proposition~2.4, p.~107]{EK:MP}. Applying a density argument and the Kakutani--Riesz--Markov representation theorem (the latter as in Tsai's axiomatic characterization \cite{Tsai:24}) then connects this construction to the stochastic heat flow for all $s \in \mathcal S$ and $t \ge s$.
 \hfill $\blacklozenge$
\end{rmk}

\begin{figure}[h!]
\centering
\begin{tikzpicture}[
  node distance=0.8cm and 1.2cm,
  box/.style={draw, rounded corners, minimum height=1.3cm, minimum width=2.6cm, align=center, font=\small},
  sec/.style={box},
  thm/.style={box},
  arrow/.style={->, >=stealth, thick}
]

  \node[sec] (sec2) {\underline{Section~\ref{sec:smt-se}}\\ \footnotesize Statements of Theorems~\ref{thm:main1},\\ \footnotesize \ref{thm:main2}, and \ref{thm:main3}};
  \node[sec, right=of sec2] (sec3) {\underline{Section~\ref{sec:soln+t}}\\ \footnotesize Baseline stochastic analytic framework\\ \footnotesize and statements of Theorems~\ref{thm:mombdd} and \ref{thm:tight}};

  \node[thm, below=of sec2] (sec6) {\underline{Section~\ref{sec:mom}}\\ \footnotesize Proofs of Theorems~\ref{thm:mombdd},\\ \footnotesize \ref{thm:tight}, and \ref{thm:main1}~(2)};
  \node[thm, below=of sec3] (sec4) {\underline{Section~\ref{sec:lcd1}}\\ \footnotesize Proofs of Theorems~\ref{thm:main1}~(1)\\ \footnotesize and \ref{thm:main2}};
\node[thm, below=.851cm of sec6] (sec5) {\underline{Section~\ref{sec:lcd2}}\\ \footnotesize Proof of Theorem~\ref{thm:main3}};

  \node[thm, below=of sec4] (sec7) {\underline{Section~\ref{sec:heat}}\\ \footnotesize Heat kernel expansions\\ \footnotesize and integral bounds};

  \draw[arrow] (sec2) -- (sec3);
  \draw[arrow] (sec2) -- (sec6);
  \draw[arrow] (sec3) -- (sec4);
  \draw[arrow] (sec3) -- (sec6);
  \draw[arrow] (sec2) -- (sec4);
  \draw[arrow] (sec7) -- (sec4);
  \draw[arrow] (sec7) -- (sec5);
  \draw[arrow] (sec2.west) -- ++(-0.6cm,0) |- (sec5.west);

\end{tikzpicture}
\caption{Dependencies among Section~\ref{sec:smt-se} and Sections~\ref{sec:soln+t}--\ref{sec:heat}}
\end{figure}

\noindent {\bf Organization of the proofs of Theorems~\ref{thm:main1}, \ref{thm:main2}, and \ref{thm:main3}.} The proof of Theorem~\ref{thm:main1}~(1$\cc$), excluding the uniqueness of limits, \eqref{Mqv:det}, and \eqref{def:mild}, relies on the more refined moment bounds and tightness estimates in Theorems~\ref{thm:mombdd} and \ref{thm:tight}. These finer theorems are stated in Section~\ref{sec:soln+t} alongside the baseline stochastic analytic framework. To prioritize our main quantitative results (Theorems~\ref{thm:main2} and \ref{thm:main3}), the proofs of Theorems~\ref{thm:mombdd}, \ref{thm:tight}, and Theorem~\ref{thm:main1}~(2$\cc$) are deferred to Section~\ref{sec:mom}. The arguments for Theorems~\ref{thm:main2} and \ref{thm:main3} build on Corollary~\ref{cor:Sk}, which introduces distinct approximate solutions that need not share the same space-time white noise (see Convention~\ref{convent:Sk} for detailed conventions). Section~\ref{sec:lcd1} then concludes the proof of Theorem~\ref{thm:main1}(1$\cc$) and presents the proof of Theorem\ref{thm:main2}. Susequently, Section~\ref{sec:lcd2} is dedicated to the proof of Theorem~\ref{thm:main3}. Finally, Section~\ref{sec:heat} provides the necessary heat kernel expansions and integral bounds underlying the proofs in Sections~\ref{sec:lcd1} and \ref{sec:lcd2}.\medskip

\noindent {\bf Setup for proofs.} From this point forward, we employ the following specialized notation:
\begin{itemize}
\item We focus exclusively on $\lv=\lv(\lambda)$, as defined in \eqref{def:lv} for any $\lambda\in \R$, under the following condition:
\begin{align}\label{def:ovl}
\ovl\in (0,1/4)\quad \mbox{and}\quad 
\lv\geq \pi/\log \vep^{-1},\; \forall\;\vep\in (0,\ovl].
\end{align}
\item For higher-order iterated integrals, we often use the shorthand $\int_x f(x)\d x$ in place of $\int_{\R^2} f(x)\d x$. Similar notation applies for integrals over other spatial variables in $\R^2$.
\item In the context of two-body dynamics, the following notation is commonly used to represent products of quantities and measures, typically indicating that the factors operate in the same time interval:
\begin{align}
\label{def:column} 
\begin{bmatrix} a_{2}\\
a_{1} \end{bmatrix}_{\times } \stackrel{\mathrm{def}} {=}
a_{1} \times a_{2},\quad \begin{bmatrix} \d x\\
\d y\end{bmatrix}_{\otimes } \stackrel{\mathrm{def}} {=}
\d x\d y.
\end{align}
\end{itemize}

\section{Expansions and tightness of approximate solutions}\label{sec:soln+t}
In this section, we provide a baseline stochastic analytic framework for solutions to \eqref{def:SHE}, formulate a precise version of the associated decomposition in \eqref{asymp:1}, and state Theorems~\ref{thm:mombdd} and \ref{thm:tight} mentioned in the above organization of proofs. 

\subsection{Pathwise properties and moment bounds}\label{sec:soln+t-ppb}
The following proposition provides explicit expressions for the covariation measures, specifies nonnegativity, and presents essential moment bounds for the approximate solutions. Most results in part (2$\cc$) follow from general theorems in \cite{Dalang, MPS:06}.

\begin{prop}\label{prop:XM} Let $\vep\in (0,\ovl]$. \medskip 

\noindent {\rm (1$\cc$)} The mollified noise $\xi_\vep$, as defined in \eqref{def:xivep}, admits an explicit covariation measure given by the following formula:
\begin{align}
\la \xi_\vep,\xi_\vep\ra(\d x',\d x,\d s)=\phi_\vep(x'-x)\d x'\d x\d s,\quad (x',x,s)\in \R^2\times \R^2\times \R_+.\label{eq:covarxi}
\end{align}
Here, $\{\phi_\vep\}$ is an approximation to the identity defined as follows:
\begin{align}\label{def:phivep+phi}
\phi_\vep(x)\,\defeq\,\vep^{-2}\phi(\vep^{-1}x),\quad 
\phi(x)\,\defeq\,\int_{y}\varrho(y)\varrho(y-x)\d y,
\end{align}
and $\phi$ is an even probability density with compact support. 
In particular, the martingale measure $M_\vep$, defined in \eqref{def:M:vep}, satisfies the following formulas: for all nonnegative $g(x',x,s)$, 
\begin{subequations}
\label{eq:covarM}
\begin{align}
&\int_{s=0}^\infty \int_{x}\int_{x'}{g}(x',x,s) \la M_\vep,M_\vep\ra(\d x',\d x,\d s)\notag\\
&\quad =\Lambda_\vep \int_{s=0}^\infty\int_x\int_{x'}{g}(x',x,s) X_\vep(x',s)X_\vep(x,s)\phi_\vep(x'-x)\d x'\d x \d s\label{eq:covarM-1}\\
&\quad =\Lambda_\vep \int_{s=0}^\infty\int_x\int_{\wt{x}}{g}(\vep \wt{x}+x,x,s) X_\vep(\vep \wt{x}+x,s)X_\vep(x,s)\phi(\wt{x})\d \wt{x}\d x \d s.\label{eq:covarM-2}
\end{align}
\end{subequations}

\noindent {\rm (2$\cc$)} Given any bounded, nonnegative initial condition $X_0(\cdot)$, the mild formulation \eqref{def:mild:vep} admits a unique pathwise strong solution $X_\vep$ taking values in $C_{\C_{\rm tem}(\R^2)}[0,\infty)$; moreover, $X_\vep(x,t)$  remains nonnegative for all $t \geq 0$ and $x\in \R^2$, and becomes strictly positive for all $t > 0$ and $x\in \R^2$ whenever $\int_{\R^2} X_0(x)\d x > 0$. Here, $\C_{\rm tem}(\R^2)$ denotes a topological space defined by  the following norms:
\begin{align}\label{def:Ctem}
\lVert f\rVert_{\infty,\gamma}\,\defeq\,\sup_{x\in\R^2}\lvert f(x)\rvert \e^{-\gamma x},\quad \gamma\in (0,\infty).
\end{align}
Specifically, this space consists of all functions $f\in \C(\R^2)$ satisfying $\lVert f\rVert_{\infty,\gamma}<\infty$ for all $\gamma>0$, with the topology induced by the norms in \eqref{def:Ctem}. Moreover, the following moment bounds hold: 
\begin{gather}
\sup_{0\leq t\leq T}\sup_{x\in \R^2}\E[X_\vep(x,t)^p]<\infty,\quad \forall\;T,p\in (0,\infty),\label{Xvep:mom}\\
\E\biggl[\sup_{0\leq t\leq T}\sup_{x\in \R^2}\lvert X_\vep(x,t)\rvert^p\e^{-\gamma \lvert x\rvert }\biggr]<\infty,\quad \forall\;T,p,\gamma,\in (0,\infty).
\end{gather}
\end{prop}
\begin{proof}
{\rm (1$\cc$)} We first show that $\phi$ is even and has compact support. First, the even parity follows upon applying the change of variables $y'=y-x$ to the definition of $\phi$ in \eqref{def:phivep+phi}. To see the compact support property of $\phi$, choose $M>0$ such that $\varrho(\wt{y})=0$ for all $\lvert \wt{y}\rvert >M$. Then $\varrho(y)\varrho(y-x)=0$ for all $y\in \R^2$ and $\lvert x\rvert >2M$ since $\lvert y \rvert\leq M$ implies $\rvert y-x\rvert \geq \lvert x\rvert -\lvert y \rvert>M$.

To prove \eqref{eq:covarxi}, let $f\in \C_c(\R^2)$, and consider
\begin{align}\label{nrep}
\int_{s=0}^t \int_xf(x)\xi_\vep(\d x,\d s)=\int_{s=0}^t \int_{\wt{x}}\left(\int_x f(x)\varrho_\vep(x-\wt{x})\d x\right) \xi(\d \wt{x},\d s),
\end{align}
which follows from the definition \eqref{def:xivep} of $\xi_\vep$ and the stochastic Fubini theorem for stochastic integrals with respect to space-time white noise (cf. \cite[Theorem~2.4.1, pp.82--83]{DS:SPDE} or \cite[Theorem~2.6, pp.296--297]{Walsh}). By \eqref{nrep}, we have
\begin{align}
&\left\langle\int_{s=0}^\cdot \int_x f(x)\xi_\vep(\d x,\d s),\int_{s=0}^\cdot \int_x f(x)\xi_\vep(\d x,\d s)\right\rangle_t\notag
\\
&\quad=\int_{s=0}^t\int_{\wt{x}}\left(\int_{x'} f(x')\varrho_\vep(x'-\wt{x})\d x'\right)\left(\int_x f(x)\varrho_\vep(x-\wt{x})\d x\right)\d \wt{x}\d s\notag\\
&\quad =\int_{s=0}^t \int_{x}\int_{x'}f(x')f(x) \left(\int_{\wt{x}}\varrho_\vep(x'-\wt{x})\varrho_\vep(x-\wt{x})\d \wt{x}\right) \d x'\d x \d s\notag\\
&\quad =\int_{s=0}^t \int_{x}\int_{x'}f(x')f(x)\left(\int_{y}\varrho_\vep(y)\varrho_\vep(x+y-x')\d y\right) \d x'\d x \d s\quad\mbox{($\because y=x'-\wt{x}$)}\notag\\
&\quad =\int_{s=0}^t \int_{x}\int_{x'}f(x')f(x)\phi_\vep(x'-x) \d x'\d x \d s,\label{mu:proof1}
\end{align}
where \eqref{mu:proof1} uses the definition of $\phi_\vep$ in \eqref{def:phivep+phi}. The last equality is enough to prove \eqref{eq:covarxi}. \medskip 

\noindent {\rm (2$\cc$)} All properties stated in (2$\cc$), except for nonnegativity, follow by applying general theorems from \cite{Dalang, MPS:06}. For convenience, we rely on \cite[Theorem 1.2, Proposition 1.8 (a), and Theorem~A.1]{MPS:06} throughout. The required conditions are met by taking $k(x,y)\equiv \phi_\vep(x-y)$. In particular, condition $(A)_1$ in  \cite[Theorem~A.1]{MPS:06} is satisfied because, using the Fourier transform $\mathcal Ff(z)\,\defeq \int \exp \{-2\pi \i z\cdot x\}f(x)\d x$, we have
\begin{gather}
\int_{\R^2} f(x)\phi_\vep(x)\d x=\int_{\R^2}\mathcal Ff(z) \lvert\mathcal F\varrho_\vep(z)\rvert^2\d z,\label{WellX:1}\\
\int_{\R^2}\frac{\lvert\mathcal F\varrho_\vep(z)\rvert^2}{1+\lvert z\rvert ^2}\d z<\infty,\label{WellX:2}\\
\sum_{j,k=1}^nc_j\overline{c}_k\phi_\vep(y_j-y_k)=\int_{\R^2}\sum_{j,k=1}^nc_j\overline{c}_k\varrho_\vep(y+y_j)\varrho_\vep(y+y_k)\d y\geq 0,\label{WellX:3}
\end{gather}
for any rapidly decreasing test function $f$ and any $n\in \Bbb N$, $y_1,\cdots,y_n\in \R^2$, and $c_1,\cdots,c_n\in \Bbb C$. Here, Plancherel's theorem and \eqref{def:phivep+phi}
 provide \eqref{WellX:1}. Moreover, \eqref{WellX:2} holds because $\varrho\in \C_c(\R^2)$ by assumption, and \eqref{WellX:3} uses the definition of $\phi_\vep$ from \eqref{def:phivep+phi}. 

Finally, the nonnegativity of $X_\vep$ and its strict positivity under the condition $\int_{\R^2}X_0(x)\d x>0$ follow directly from pathwise uniqueness. Specifically, in this framework, the Feynman--Kac representation (see, e.g., \cite[p.129]{C:DBG}) produces an explicit solution that inherits nonnegativity from $X_0(\cdot)$ and becomes strictly positive for all $t > 0$ whenever the initial mass is non-zero.
\end{proof}

\subsection{Quasi-It\^o expansion of covariation measures}\label{sec:soln+t-dec}
The main result of this subsection, stated as Proposition~\ref{prop:qvdec}, formulates a precise version of the decomposition in \eqref{asymp:1} at the level of $\vep>0$. We use a weak formulation by integrating over space and time with respect to a test function $f(x,t)$ for $x \in \R^2$ and $0 \leq t \leq T$. Thus, the decomposition in \eqref{asymp:1} is made precise through \eqref{mild:goal} and the following identification, utilizing rigorously defined functionals $I_\vep^j(f,T)$, $1\leq j\leq 4$, introduced below in \eqref{def:I1}--\eqref{def:I4}:
\begin{align*}
 I_\vep^1(f,T)+I_\vep^2(f,T)&\longleftrightarrow 
X(x,t)^2- 
\int_{s=0}^t \int_y P_{t-s}(x,y)^2 \Lambda X(y,s)^2 \d y\d s ,\\
-2I^3_\vep(f,T) &\longleftrightarrow  -2P_tX_0(x)\int_{s=0}^t \int_{y} P_{t-s}(x,y)\M(\d y,\d s),\\
-2I^4_\vep(f,T)&\longleftrightarrow -2\int_{s'=0}^t \int_{y'}\int_{s=0}^{s'}\int_{y}P_{t-s'}(x,y') P_{t-s}(x,y)\M(\d y,\d s) \M(\d y',\d s') .
\end{align*}
In particular, the functionals $I_\vep^j(f, T)$ employ either the standard or stochastic Fubini theorem, both of which allow the test function $f(x, t)$ to regularize singularities in the limit.

\begin{prop}\label{prop:qvdec}
For all $\vep\in (0,\ovl]$, $T>0$, and $f(x,t)\in \C_{pd}^0(\R^2\times [0,T])$, it holds that
\begin{align}\label{mild:goal}
\int_{t=0}^T\int_xf(x,t)[P_tX_0(x)]^2\d x\d t=I_\vep^1(f,T)+I_\vep^2(f,T)-2I_\vep^3(f,T)-2I_\vep^4(f,T),
\end{align}
where
\begin{align}
\begin{split}\label{def:I1}
I^{1}_{\vep}(f,T)&\,\defeq\,\int_{t=0}^T\int_x f(x,t)X_\vep(x,t)^2\d x\d t-\int_{s=0}^T \int_{y} \int_{\wt{y}}\mathfrak L_\vep^1f (\wt{y},y,s,T) X_\vep(y,s)^2\d \wt{y}\d y\d s,
\end{split}\\
\begin{split}\label{def:I2}
I^2_\vep(f,T)&\,\defeq\,
-\int_{s=0}^T\int_{y} \int_{\wt{y}} \mathfrak L_\vep^1f (\wt{y},y,s,T) [X_\vep(\vep \wt{y}+y,s)-X_\vep(y,s)]X_\vep(y,s)\d \wt{y}\d y\d s,
\end{split}\\
\begin{split}\label{def:I3}
I^3_\vep(f,T)&\,\defeq\,\int_{s=0}^T \int_{y} \mathfrak L_0^3f(y,s,T)\M_\vep(\d y,\d s),
\end{split}\\
\begin{split}\label{def:I4}
I^4_\vep(f,T)&\,\defeq\,\int_{s'=0}^T \int_{y'}  \int_{s=0}^{s'}\int_{y} \mathfrak L_0^4f(y',y,s',s,T) \M_\vep(\d y,\d s)\M_\vep(\d y',\d s'),
\end{split}
\end{align}
and the integrands use the following integral transformations of $f(x,t)$:
\begin{align}
\mathfrak L^1_\vep f(\wt{y},y,s,T)&\,\defeq\;\mathring{\mathfrak L}^1_\vep f(\wt{y},y,s,T) \phi(\wt{y}),\label{def:L1}\\
\mathring{\mathfrak L}^1_\vep f(\wt{y},y,s,T)&
\,\defeq\,\Lambda_\vep\int_{t=s}^T\int_{x}f(x,t)P_{t-s}(\vep \wt{y}+y,x) P_{t-s}(x,y)\d x \d t,\label{def:L1r}
\\
\mathfrak L_0^3f(y,s,T)&\,\defeq\,\int_{t=s}^T\int_{x}f(x,t)P_tX_0(x) P_{t-s}(x,y)\d x \d t,\label{def:L3}\\
\mathfrak L_0^4f(y',y,s',s,T)&\,\defeq\,\int_{t=s'}^T\int_{x}f(x,t)P_{t-s'}(x,y')P_{t-s}(x,y) \d x\d t .\label{def:L4}
\end{align}
\end{prop}
\begin{proof}
We first derive the main equation of this proof.
By squaring both sides of \eqref{def:mild:vep} and rearranging, we have
\begin{align*}
[P_tX_0(x)]^2&=X_\vep(x,t)^2-2P_tX_0(x)\int_{s=0}^t \int_{y} P_{t-s}(x,y)M_\vep(\d y,\d s)\\
&\quad -\left(\int_{s=0}^t \int_{y} P_{t-s}(x,y)M_\vep(\d y,\d s)\right)^2.
\end{align*}
Integrating both sides against $f(x,t)$ over $x\in\R^2$ and $0\leq t\leq T$ yields
\begin{align}
&\int_{t=0}^T\int_{x}f(x,t)[P_tX_0(x)]^2\d x\d t \notag\\
&\quad =\int_{t=0}^T\int_{x} f(x,t) X_\vep(x,t)^2\d x\d t  -2\int_{t=0}^T\int_{x}f(x,t) P_tX_0(x)\int_{s=0}^t \int_{y} P_{t-s}(x,y)M_\vep(\d y,\d s)\d x\d t \notag\\
&\quad\quad -\int_{t=0}^T\int_{x} f(x,t)\left(\int_{s=0}^t\int_y P_{t-s}(x,y)M_\vep(\d y,\d s)\right)^2\d x\d t \notag\\
&\quad =\int_{t=0}^T\int_{x}f(x,t) X_\vep(x,t)^2\d x\d t  -2\int_{t=0}^T\int_{x}f(x,t) P_tX_0(x)\int_{s=0}^t \int_{y} P_{t-s}(x,y)\M_\vep(\d y,\d s)\d x\d t \notag\\
&\quad\quad  -2\int_{t=0}^T\int_{x} f(x,t)\int_{s'=0}^t  \int_{y'} \int_{s=0}^{s'}\int_{y} P_{t-s}(x,y)P_{t-s'}(x,y')\M_\vep(\d y,\d s)\M_\vep(\d y',\d s')\d x\d t \notag\\
&\quad \quad -\Lambda_\vep \int_{t=0}^T\int_{x}f(x,t) \int_{s=0}^t \int_{y}\int_{\wt{y}} P_{t-s}(x,\vep \wt{y}+y)P_{t-s}(x,y)\notag\\
&\quad \quad \times X_\vep(\vep \wt{y}+y,s)X_\vep(y,s)\phi(\wt{y})  \d \wt{y}\d y\d s \d x\d t  .\label{SMF}
\end{align}
Here, \eqref{SMF} is derived by applying It\^{o}'s formula and \eqref{eq:covarM-2} to expand the squared process $(\int_{s=0}^{t} \int_y P_{t-s}(x, y) \M_\vep(\d y, \d s))^2$, as $\int_{s=0}^{s'} \int_y P_{t-s}(x, y) \M_\vep(\d y, \d s)$, for $0 \leq s' \leq t$, forms an $L^2$-bounded martingale. This martingale property and $L^2$-integrability are ensured by \eqref{eq:covarM}, \eqref{Xvep:mom}, and the estimate below, which employs the change of variables $\wt{y} = y' - y$ and the Chapman--Kolmogorov equation: 
\begin{align}
&\sup_{0\leq t\leq T}\sup_{x\in \R^2}\int_{s=0}^{t}\int_{y}\int_{y'}P_{t-s}(x, y')P_{t-s}(x,y)\phi_\vep(y'-y)\d y'\d y\d s\notag\\
&\quad =\sup_{0\leq t\leq T}\int_{s=0}^{t}\int_{\wt{y}}P_{2(t-s)}(\wt{y})\phi_\vep(\wt{y})\d \wt{y}\d s
\leq \sup_{0\leq t\leq T} t\lVert\phi_\vep\rVert_\infty<\infty.\label{SFT:main}
\end{align}
Equation~\eqref{SMF} is the main equation of this proof, from which we will derive \eqref{mild:goal}.

To proceed, we demonstrate that the right-hand sides of \eqref{mild:goal} and \eqref{SMF} are equal. To derive $I^1_{\vep}(f,T)$ and $I^2_{\vep}(f,T)$ from \eqref{SMF}, we utilize the standard Fubini theorem to switch the order of integration in the last term of \eqref{SMF} from $\d \wt{y}\d y\d s \d x\d t$ to $\d x\d t\d \wt{y}\d y\d s $:
\begin{align}
&\Lambda_\vep \int_{t=0}^T\int_{x}f(x,t) \int_{s=0}^t \int_{y}\int_{\wt{y}} P_{t-s}(x,\vep \wt{y}+y)P_{t-s}(x,y)\notag\\
&\quad \times X_\vep(\vep \wt{y}+y,s)X_\vep(y,s)\phi(\wt{y}) \d \wt{y}\d y\d s \d x\d t\notag\\
&\quad =\int_{s=0}^T \int_{y}\int_{\wt{y}}\left(\Lambda_\vep\int_{t=s}^T\int_{x}f(x,t) P_{t-s}(\vep \wt{y}+y,x)P_{t-s}(x,y)\d x \d t\phi(\wt{y})\right) \notag\\
&\quad\quad  \times X_\vep(\vep \wt{y}+y,s)X_\vep(y,s)\d \wt{y}\d y\d s \notag\\
&\quad =\int_{s=0}^T \int_{y}\int_{\wt{y}} \mathfrak L_\vep^1f (\wt{y},y,s,T) X_\vep(\vep \wt{y}+y,s)X_\vep(y,s)\d \wt{y}\d y\d s,\label{SMF:1}
\end{align}
where the last equality uses the definition \eqref{def:L1} of $\mathfrak L_\vep^1f (\wt{y},y,s,T)$. Hence, the sum of the first and last terms on the right-hand side of \eqref{SMF} equals $I^1_{\vep}(f,T)+I^2_{\vep}(f,T)$.

Next, to derive $I^3_{\vep}(f,T)$ and $I^4_{\vep}(f,T)$ from \eqref{SMF}, we use the stochastic Fubini theorem (see \cite[Theorem~2.4.1, p.82]{DS:SPDE} or \cite[Theorem~2.6, pp.296--297]{Walsh}). This theorem allows us to repeatedly interchange the order of integration in the second and third terms on the right-hand side of \eqref{SMF}. The following two calculations give the details of these interchanges: 
\begin{align}
&-2\int_{t=0}^T\int_{x}f(x,t) P_tX_0(x)\int_{s=0}^t \int_{y} P_{t-s}(x,y)\M_\vep(\d y,\d s)\d x\d t\notag\\
&\quad =-2\int_{s=0}^T \int_{y}\left(\int_{t=s}^T\int_{x}f(x,t)P_tX_0(x) P_{t-s}(x,y)\d x\d t\right)\M_\vep(\d y,\d s)\label{SFT:000-1}\\
&\quad =-2I^3_\vep(f,T),\label{SMF:2}
\end{align}
and
\begin{align}
&-2\int_{t=0}^T\int_{x} f(x,t)\int_{s'=0}^t  \int_{y'} \int_{s=0}^{s'}\int_{y} P_{t-s}(x,y)P_{t-s'}(x,y')\M_\vep(\d y,\d s)\M_\vep(\d y',\d s')\d x\d t \notag\\
&\quad =-2\int_{s'=0}^T  \int_{y'} \int_{t=s'}^T\int_{x}f(x,t)\int_{s=0}^{s'}\int_{y} P_{t-s}(x,y)P_{t-s'}(x,y')\M_\vep(\d y,\d s)\d x\d t\notag\\
&\quad \quad \times \M_\vep(\d y',\d s')\label{SFT:000-2}\\
&\quad =-2\int_{s'=0}^T  \int_{y'} \int_{s=0}^{s'}\int_{y}\left(\int_{t=s'}^T\int_{x}f(x,t) P_{t-s}(x,y)P_{t-s'}(x,y')\d x\d t\right)\M_\vep(\d y,\d s)\notag\\
&\quad \quad \times \M_\vep(\d y',\d s')\label{SFT:000-3}\\
&\quad=-2I^4_\vep(f,T).\label{SMF:3}
\end{align}
The verification of the integrability conditions associated with \eqref{SFT:000-1}, \eqref{SFT:000-2}, and \eqref{SFT:000-3} is deferred until the end of the proof. By \eqref{SMF:1}, \eqref{SMF:2}, and \eqref{SMF:3}, we have proved that the right-hand sides of  \eqref{mild:goal} and \eqref{SMF} are equal. 

To complete the proof, it remains to justify the previous uses of the stochastic Fubini theorem. We verify the relevant $L^2$-condition as stated in \cite[Theorem~2.6, pp.296--297]{Walsh}
(see also \cite[Corollary~2.4.3, p.90]{DS:SPDE}). By It\^{o}'s isometry, \eqref{eq:covarM-2}, and \eqref{Xvep:mom}, this verification amounts to establishing the following three bounds, which justify \eqref{SFT:000-1}, \eqref{SFT:000-2}, and \eqref{SFT:000-3}, respectively: 
\begin{align}
&\int_{t=0}^T\int_{x} \lvert f(x,t)\rvert P_tX_0(x)\E\biggl[\,\int_{s=0}^t\int_{y_1} \int_{y_2}P_{t-s}(x,y_2)P_{t-s}(x,y_1)\la \M_\vep,M_\vep\ra(\d y_2,\d y_1,\d s)\biggr]\d x\d t<\infty,\notag\\
&\int_{t=0}^T\int_{x}\lvert f(x,t)\rvert \E\biggl[\,\int_{s'=0}^t\int_{y_1'}  \int_{y_2'} P_{t-s'}(x,y_2') P_{t-s'}(x,y_1')
\biggl(\int_{s=0}^{s'}\int_{y} P_{t-s}(x,y)\M_\vep(\d y,\d s)\biggr)^2\notag\\
&\quad \times \la \M_\vep,\M_\vep\ra(\d y_2',\d y_1',\d s')\biggr]\d x\d t<\infty,\label{SFT:000-2-1}\\
&  \int_{t=s'}^T\int_{x}\lvert f(x,t)\rvert P_{t-s'}(x,y')\E\biggl[\,\int_{s=0}^{s'}\int_{y_1}\int_{y_2} P_{t-s}(x,y_2)P_{t-s}(x,y_1)\la \M_\vep,\M_\vep\ra(\d y_2,\d y_1,\d s)\biggr]\notag\\
&\quad\times \d x\d t<\infty,\notag
\end{align}
where $(s',y')$ in the third bound ranges over $0<s'<T$ and $y'\in \R^2$. To obtain the first and third of these bounds, we apply \eqref{eq:covarM-1}, \eqref{Xvep:mom}, and \eqref{SFT:main}. Additionally, to obtain \eqref{SFT:000-2-1}, we apply \eqref{SFT:main} again, now together with the following bound:
\begin{align*}
&\sup_{0\leq t\leq T}\sup_{0\leq s'\leq t}\sup_{x,y_2',y_1'\in \R^2}\E\biggl[\biggl(\int_{s=0}^{s'}\int_{y} P_{t-s}(x,y)\M_\vep(\d y,\d s)\biggr)^2X_\vep(y_2',s')X_\vep(y_1',s')\biggr]\\
&\quad \less\sup_{0\leq t\leq T} \sup_{x\in \R^2}\E\biggl[\biggl(\int_{s=0}^{t}\int_{y_1}\int_{y_2} P_{t-s}(x,y_2) P_{t-s}(x,y_1)\la \M_\vep,\M_\infty\ra(\d y_2,\d y_1,\d s)\biggr)^2\,\biggr]^{1/2}\\
&\quad \quad \times \sup_{0\leq t\leq T}\sup_{y_2',y_1'\in \R^2}\E\left[X_\vep(y_2',t)^2X_\vep(y_1',t)^2\right]^{1/2}<\infty,
\end{align*}
where the first inequality follows from the Cauchy--Schwarz inequality and the Burkholder--Davis--Gundy inequality \cite[(4.1) Theorem, p.160]{RY}, and the last inequality can be deduced from \eqref{Xvep:mom} and \eqref{SFT:main}. We have justified all of the uses of the stochastic Fubini theorem in \eqref{SFT:000-1}, \eqref{SFT:000-2}, and \eqref{SFT:000-3}. The proof is complete.
\end{proof}

The bound in \eqref{SFT:main} will be used frequently in the upcoming applications of the stochastic Fubini theorem when $\vep>0$. For clarity, we state it separately in the following lemma.

\begin{lem}\label{lem:SFTint}
For all $\vep\in (0,\ovl]$ and $T>0$,
\[
\sup_{0\leq t\leq T}\sup_{x\in \R^2}\int_{s=0}^{t}\int_{y}\int_{y'}P_{t-s}(x, y')P_{t-s}(x,y)\phi_\vep(y'-y)\d y'\d y\d s<\infty.
\]
\end{lem}

\subsection{From tightness to Skorokhod's representation}\label{sec:lcd1-sk}
In this section, we refine Theorem~\ref{thm:main1} and introduce the framework for the Skorokhod representation, setting the stage for the proof of Theorem~\ref{thm:main2} alongside. The essential results focus on specific moment bounds detailed in the next theorem: parts (1$\cc$) and (2$\cc$) parallel \eqref{MMbdd} and \eqref{X4BMbdd}, while part (3$\cc$) corresponds to \eqref{m1:apriori}, though it is stated here in a form tailored for the forthcoming applications. 

\begin{thm}\label{thm:mombdd}
Let $\vep\in (0,\ovl]$, $p\in (0,1)$, $T\in (0,\infty)$ and $x_0^1,x_0^0,\wt{x}_0^1,\wt{x}_0^0\in \R^2$.
\begin{itemize}
\item [\quad \rm (1$\cc$)] For all $0\leq T_1\leq T_2\leq T$, 
\begin{align}
&\int_{t=T_1}^{T_2}\int_{t'=t}^{T_2}\lv^2\E[X_\vep(x_0^1,t)X_\vep(x_0^1+\vep x_0^0,t)X_\vep(\wt{x}_0^1,t')X_\vep(\wt{x}_0^1+\vep \wt{x}_0^0,t')]\d t'\d t\notag\\
 &\quad \leq C(\lambda,\phi,\lVert X_0\rVert_\infty,T,p)\times
\biggl(\lv\int_{s=0}^{T}  \int_{\ol{z}}\phi(\ol{z})P_{2s}(\vep x_0^0-\vep\ol{z})\d\ol{z}\d s+\lv\biggr)\notag\\
&\quad \quad \times \biggl(\lv\int_{s=0}^{T}  \int_{\ol{z}}\phi(\ol{z})P_{2s}(\vep\wt{x}_0^0-\vep\ol{z})\d\ol{z}\d s+\lv\biggr) 
  \times (T_2-T_1)^{2-p}.\label{modulus:bdd4}
\end{align}
\item [\rm (2$\cc$)] For all $t,a',a\in (0,T]$,
\begin{align}
&\int_{z_2}\int_{z_1}P_{a'}(\wt{x}_0^1+\vep \wt{x}_0^0,z_2) P_{a}(\wt{x}_0^1,z_1)\lv\E[X_\vep ( z_1,t)X_\vep(z_2,t)X_\vep(x_0^1,t)X_\vep(x_0^1+\vep x_0^0,t)]\d z_1\d z_2\notag\\
&\quad \leq \frac{ C(\lambda,\phi,\lVert X_0\rVert_\infty,T,p)}{\min\{a',a\}^p}\times\biggl( \lv\int_{s=0}^{T} \int_{\ol{z}}\phi(\ol{z})P_{2s}(\vep x_0^0-\vep\ol{z})\d\ol{z}\d s+\lv\biggr) .
\label{modulus:bdd2}
\end{align}
\item [\rm (3$\cc$)] For all $t\in (0,T]$,
\begin{align}
&\lv \E[X_\vep(x_0^1,t)X_\vep(x_0^1+\vep x_0^0,t)]\notag\\
&\quad  \leq C(\lambda,\phi,\lVert X_0\rVert_\infty,T)\times\biggl( \lv\int_{s=0}^{T}  \int_{\ol{z}}\phi(\ol{z})P_{2s}(\vep x_0^0-\vep\ol{z})\d\ol{z}\d s+\lv\biggr) .\label{mod:final1}
\end{align}
\end{itemize}
Moreover, we have
\begin{align}\label{thmmombdd:aux}
\lv\int_{s=0}^{T} \int_{\ol{z}}\phi(\ol{z})P_{2s}(\vep x_0^0-\vep\ol{z})\d\ol{z}\d s
&\leq C(\lambda,T)\left(\int_{\ol{z}}\phi(\ol{z})\lvert \log \lvert x_0^0-\ol{z}\rvert \rvert \d \ol{z}+1\right),
\end{align}
which leads to improved bounds in (1$\cc$)--(3$\cc$) when $x_0^0$ and $\wt{x}_0^0$ are both of order one as $\vep \to 0$.
\end{thm}

The proof of Theorem~\ref{thm:mombdd}, presented in Section~\ref{sec:dual}, builds upon the diagrammatic expansion method of \cite{C:DBG} (see also \cite{CSZ:19a, GQT}). The argument begins by extending the uniform-in-$\vep$ bound for $\E\big[\prod_{j=1}^N X_\vep(x_0^j, t)\big]$ from \cite[Theorem~2.3]{C:DBG}, which holds for all integers $N \geq 3$ given fixed initial conditions $x_0 = (x_0^1, \ldots, x_0^N) \in \mathbb{R}^{2N}$ with distinct components. However, new challenges emerge when initial points coalesce, such as in \eqref{modulus:bdd4} and \eqref{modulus:bdd2} where $x_0^0$ and $\wt{x}_0^0$ remain of order one as $\vep \to 0$. This proximity induces new singularities that lead to divergent moments. Resolving these singularities requires both explicit renormalization via direct multiplication by $\lv$ and substantially more delicate analysis. By contrast, while the moments in \eqref{mod:final1} exhibit a comparable formulation, their second-moment structure renders them far more tractable.

Theorem~\ref{thm:tight} below offers a partial refinement of the weak convergence statement in Theorem~\ref{thm:main1} (1$\cc$) by introducing two additional measure-valued processes. These processes perturb the covariation measures $\la M_\vep, M_\vep \ra(\d x', \d x, \d s)$ and arise from the terms $I^j_\vep(f, T)$ in Proposition~\ref{prop:qvdec}. Specifically, Theorem~\ref{thm:tight} addresses the joint tightness of four measure-valued processes defined by the following equations using $f_0\in \C_b(\R^2)$ and $g_0\in \C_b(\R^2\times \R^2)$:
 \begin{subequations}
\label{def:sharp}
\begin{align}
X^\sharp_\vep(f_0,t)&\,\defeq\, X_\vep(\psi f_0,t),\quad \;\;\;
\nu^\sharp_\vep(f_0,t)\,\defeq\, \nu_\vep(\psi f_0,t),\label{def:sharp1}\\
\mu^\sharp_\vep(g,t)&\,\defeq\,\mu_\vep( \psi^{\otimes 2} g_0,t),\quad \mathring{\mu}^\sharp_\vep(g_0,t)\,\defeq\,\mathring{\mu}_\vep(\psi^{\otimes 2} g_0,t).\label{def:sharp2}
\end{align}
\end{subequations}
Here, we use the auxiliary function $\psi$ from \eqref{def:psi}, as well as the following definitions, where we permit time-dependent test functions for subsequent applications:
\begin{subequations}
\label{def:processes}
 \begin{align}
 X_\vep(f,t)&\,\defeq\,\int_x f(x,t) X_\vep(x,t)\d x,\label{def:Xf}\\
  \nu_\vep(f,t)&\,\defeq\,\int_{s=0}^t\int_{x}f(x,s)\lv X_\vep(x,s)^2 \d x\d s,\label{def:covarprocessnu}\\
   \mu_\vep(g,t)&\,\defeq\,\int_{s=0}^t\int_{x}\int_{x'} g(x',x,s)\la M_\vep,M_\vep\ra(\d x',\d x,\d s)\notag\\
 &=\int_{s=0}^t\int_{x}\int_{\wt{x}}g(\vep \wt{x}+x,x,s)\lv X_\vep(\vep \wt{x}+x,s)X_\vep(x,s)\phi(\wt{x})\d \wt{x}\d x\d s,\label{def:covarprocessmu}\\
 \mathring{\mu}_\vep(g,t)&\,\defeq\,\int_{s=0}^t\int_{x}\int_{\wt{x}}g(\wt{x},x,s)\lv X_\vep(\vep \wt{x}+x,s)X_\vep(x,s)\phi(\wt{x})\d \wt{x}\d x\d s.\label{def:covarprocessmur}
\end{align}
\end{subequations}
Owing to the choice of $\psi$ and \eqref{Xvep:mom}, all the processes in \eqref{def:sharp} are finite-valued. 

\begin{thm}\label{thm:tight}
The family of joint laws of $\{(X^\sharp_\vep,\nu^\sharp_\vep,\mu^\sharp_\vep,\mathring{\mu}^\sharp_\vep)\}_{\vep\in (0,\ovl]}$ is tight as a family of probability measures on 
the Polish space $C_{\mathcal M_f(\R^2)}[0,\infty)\times  C_{\mathcal M_f(\R^2)}[0,\infty)\times C_{\mathcal M_f(\R^4)}[0,\infty)\times C_{\mathcal M_f(\R^4)}[0,\infty)$, where $\mathcal M_f(\R^d)$ denotes the space of (nonnegative) finite measures equipped with the topology of weak convergence.
\end{thm}

The proof of this theorem appears in Section~\ref{sec:jointC}. As an immediate consequence, we obtain the following corollary, which is justified by Theorem~\ref{thm:tight} together with Prohorov's theorem~\cite[2.2 Theorem, p.104]{EK:MP} and Skorokhod's representation theorem~\cite[1.8 Theorem, p.102]{EK:MP}.

\begin{cor}\label{cor:Sk}
Given any sequence $\{\vep_n\}_{n\in\mathbb{N}}$ with $\vep_n \in (0,\ovl]$ and $\vep_n \searrow 0$, there exist a subsequence, still denoted by $\{\vep_n\}_{n\in \Bbb N}$, and a change of the underlying probability space (making the space-time white noise $\xi$ to possibly depend on $n \in \mathbb{N}$) such that 
\begin{align}\label{Skseq}
(X^\sharp_{n},\nu^\sharp_{n},\mu^\sharp_{n},\mathring{\mu}^\sharp_{n})
\,\defeq\,
(X^\sharp_{\vep_n},\nu^\sharp_{\vep_n},\mu^\sharp_{\vep_n},\mathring{\mu}^\sharp_{\vep_n}),\quad n\in \Bbb N,
\end{align}
forms an almost-sure convergent sequence in $C_{\mathcal M_f(\R^2)}[0,\infty)\times  C_{\mathcal M_f(\R^2)}[0,\infty)\times C_{\mathcal M_f(\R^4)}[0,\infty)\times C_{\mathcal M_f(\R^4)}[0,\infty)$. Hence, with probability one,
\begin{align}\label{eq:Sko}
\begin{split}
&\lim_{n\to\infty}(X^\sharp_{n},\nu^\sharp_{n},\mu^\sharp_{n},\mathring{\mu}^\sharp_{n})=(X^\sharp_{\infty},\nu^\sharp_{\infty},\mu^\sharp_{\infty},\mathring{\mu}^\sharp_{\infty})
\end{split}
\end{align}
for some $(X^\sharp_{\infty},\nu^\sharp_{\infty},\mu^\sharp_{\infty},\mathring{\mu}^\sharp_{\infty})\in C_{\mathcal M_f(\R^2)}[0,\infty)\times  C_{\mathcal M_f(\R^2)}[0,\infty)\times C_{\mathcal M_f(\R^4)}[0,\infty)\times C_{\mathcal M_f(\R^4)}[0,\infty)$. 
\end{cor}

{\bf Unless stated otherwise, the rest of this paper will use a fixed sequence as given in Corollary~\ref{cor:Sk}.} Furthermore, for simplicity, we will often write ``$n$'' instead of the subscript ``$\vep_n$'' from this point onward.

Although the quadruplets in \eqref{Skseq} may no longer share the same space-time white noise, they still uniquely determine the corresponding processes without the ``$\sharp$'' notation. These processes, together with the extension to $n=\infty$, are explicitly defined via the following equations:
\begin{align}
&\big(X_{n}(f,t),\nu_{n}(f,t),\mu_{n}(g,t),\mathring{\mu}_{n}(g,t)\big)\notag\\
&\quad\,\defeq\,\big(X^\sharp_{n}(f/\psi,t),\nu^\sharp_{n}(f/\psi,t),\mu^\sharp_{n} (g/\psi\otimes \psi,t),\mathring{\mu}^\sharp_{n}(g/\psi\otimes \psi,t)\big),\quad n\in \Bbb N\cup \{\infty\},\label{def:Xinfty}
\end{align}
where $f\in \C^0_{pd}(\R^2)$ and $g\in \C^0_{pd}(\R^2\times\R^2)$. Additionally, we define
\begin{align}\label{XM:relation}
M_n(f,t)\,\defeq\,X_n(f,t)-X_n(f,0)-\int_0^t X_n\left(\frac{\Delta f}{2},s\right)\d s,\quad f\in \C^2_{pd}(\R^2),\; n\in \Bbb N\cup \{\infty\}.
\end{align}
Proposition~\ref{prop:thm1-1-prelim} (to be established below) demonstrates that for $n=\infty$, the triple $(\nu_n, \mu_n, \mathring{\mu}_n)$ and the martingale measure $M_n$, as defined by \eqref{XM:relation}, are uniquely determined by $X_n$. For finite $n$, this property follows immediately. 

\begin{convent}\label{convent:Sk}
Under the notation introduced after Corollary~\ref{cor:Sk} up to \eqref{XM:relation}, the equations in \eqref{def:processes} specialized to $\vep=\vep_n$ hold as equalities rather than definitions throughout the remainder of this work. This same interpretation applies to all preceding notation, definitions, and equations, provided they do not explicitly reference the underlying space-time white noise.\mbox{\quad} \hfill $\blacklozenge$
\end{convent}

\begin{cor}\label{cor:muinfty}
We can choose an event $\Omega_0$ with $\P(\Omega_0)=1$ on which for all $f_0\in \C_b(\R^2) $ and $g_0\in \C_b(\R^2\times \R^2)$, the following limit holds in $C_{\R}[0,\infty)\times C_{\R}[0,\infty)\times C_{\R}[0,\infty)\times C_{\R}[0,\infty)$:
\begin{align}\label{eq:muinfty1}
\begin{split}
&\lim_{n\to\infty}\big(X_n^\sharp (f_0,t), \nu_n^\sharp (f_0,t),\mathring{\mu}_{n}^\sharp(g_0,t), \mu_n^\sharp (g_0,t)\big)\\
&\quad =\big(X_\infty^\sharp (f_0,t), \nu_\infty^\sharp (f_0,t),\mathring{\mu}_{\infty}^\sharp(g_0,t), \mu_\infty^\sharp (g_0,t)\big).
\end{split}
\end{align}
\end{cor}
\begin{proof}
The proof readily follows from \eqref{eq:Sko} and the observation that for any integer $d\geq 1$ and $f\in \C_b(\R^d)$, the following mappings are continuous:
\begin{gather}
\mu\mapsto \mu(f): \mathcal M_f(\R^d)\to \R,\label{cont:Mf1}\\
\{\mu_t;t\geq 0\}\mapsto \{\mu(f,t);t\geq 0\}:C_{\mathcal M_f(\R^d)}[0,\infty)\to C_{\R}[0,\infty).\label{cont:Mf2}
\end{gather}
Here, $\mu(f,t)\,\defeq\,\int f\d \mu_t$ as before, and 
 the continuity of \eqref{cont:Mf1} follows directly from the topology of $\mathcal M_f(\R^d)$. Additionally, the continuity of \eqref{cont:Mf2} is established via a subsequence argument: $\mu^n(f, t_n) \to \mu^\infty(f, t_\infty)$ whenever $t_n \to t_\infty \in \R_+$ and $\{\mu^n_t; t \geq 0\} \to \{\mu^\infty_t; t \geq 0\}$ in $C_{\mathcal M_f(\R^d)}[0, \infty)$.
\end{proof}

\section{Limiting expansions of covariation measures}\label{sec:lcd1}
This section is dedicated to proving Theorem~\ref{thm:main2}. The main arguments are organized across Sections~\ref{sec:I1}, \ref{sec:I2}, and \ref{sec:I34}, where we establish the precise subsequential limits of the four terms $I^j_{n}(f,T)=I^j_{\vep_n}(f,T)$ for $1\leq j\leq 4$ defined in \eqref{mild:goal}. Section~\ref{sec:lcd1-thm2} completes the section by combining these results to finalize the proof of Theorem~\ref{thm:main1}(1$\cc$). Recall that from Corollary~\ref{cor:Sk} onward, we work with the fixed sequence chosen there and adhere to Convention~\ref{convent:Sk}.

\addtocontents{toc}{\protect\setcounter{tocdepth}{1}}
\subsection{Limit theorem for $\bs I^{\bs 1}_{\bs n}\bs (\bs f\bs ,\bs T\bs )$}\label{sec:I1}
\addtocontents{toc}{\protect\setcounter{tocdepth}{2}}
\addcontentsline{toc}{subsection}{\numberline{\thesubsection}{Limit theorem for $I^1_n(f,T)$}}

The following proposition initiates the identification of subsequential limits of the functionals $I_{\vep}^j(f,T)$ defined in Proposition~\ref{prop:qvdec}, specifically focusing on the case $j=1$. 
\begin{prop}\label{prop:I1}
For all $0<T<\infty$ and $f(x,t)\in \C_{pd}^1(\R^2\times[0,T])$, it holds that 
\begin{align*}
&I^{1}_{n}(f,T)\,\defeq\,\int_{s=0}^T\int_x f(x,t)X_n(x,t)^2\d x\d t-\int_{s=0}^T \int_{x}\left(\int_{y} \mathfrak L^1_nf (y,x,s,T)\d y\right) X_n(x,s)^2\d x\d s\\
&\quad \xrightarrow[n\to\infty]{\P} I^{1}_{\infty}(f,T)\,\defeq\,\int_{t=0}^T\int_x \mathfrak L_0f(x,t,T)\nu_\infty(\d x,\d t),
\end{align*}
and the limit defines an absolutely convergent integral. Here, $\mathfrak L^1_n\,\defeq\,\mathfrak L^1_{\vep_n}$ for $\mathfrak L^1_\vep $ in \eqref{def:L1}, and
\begin{align}
\begin{split}
\mathfrak L_0f(x,t,T)&\,\defeq -\frac{1}{2\pi}
\biggl(\frac{\log [4(T-t)]}{2}+\lambda -\frac{\gamma_{\sf EM}}{2}-\int_{y}\phi(y)\log \lvert y \rvert\d y\biggr)f(x,t)\\
&\quad\,-\int_{x'}\int_{t'=0}^{T-t}[f(x-x',t+t')-f(x,t)]P_{t'}(x')P_{t'}(x')\d t'\d x'.\label{def:L0f}
\end{split}
\end{align}
\end{prop}

Lemma~\ref{lem:I1} below provides the asymptotic expansion of $\mathring{\mathfrak L}^1_\vep f$ as $\vep\to0$, accompanied by relevant estimates. This expansion is a key ingredient in the proof of Proposition~\ref{prop:I1}, as it allows the cancellation of the divergent terms ($\infty-\infty$) originating from $X_n(x,s)^2\d x\d s$ in the two integrals defining $I^1_n(f,T)$, achieved through the use of test functions. Lemma~\ref{lem:I1} also provides results that are not used here but will be needed in Section~\ref{sec:I2}, as $\mathfrak L^1_\vep$ appears in the definition \eqref{def:I2} of $I^2_\vep(f,T)$.

\begin{lem}\label{lem:I1}
Let $0<T<\infty$, $\vep\in (0,\ovl]$, $f(x,s)\in \C_{pd}^1(\R^2\times[0,T])$, and $m\in \Bbb Z_+$.\medskip 

\noindent {\rm (1$\cc$)} 
For all $M\in (0,\infty)$, $(y,x)\in \R^2\times \R^2$ with $0<\lvert y \rvert\leq M$, and $0<s< T$, 
\begin{align}\label{key}
\begin{split}
\mathring{\mathfrak L}^1_\vep f(y,x,s,T)
&=f(x,s)+\frac{1}{\log \vep^{-1}}
\biggl(\frac{\log [4(T-s)]}{2}+\lambda -\frac{\gamma_{\sf EM}}{2}-\log \lvert y \rvert\biggr)f(x,s)\\
&\quad+\frac{2\pi  }{\log \vep^{-1}}\int_{x'}\int_{t'=0}^{T-s}[f(x-x',s+t')-f(x,s)]P_{t'}(x')P_{t'}(x')\d t'\d x'\\
&\quad+R_{\ref{key}}(y,x,s,T).
\end{split}
\end{align}
Here, $\lambda$ defines $\lv$ through \eqref{def:lv}, $\EM$ is the Euler--Mascheroni constant, the double integral is absolutely convergent
with 
\begin{align}
 \int_{x'}\int_{t'=0}^{T-s}\lvert f(x-x',s+t')-f(x,s)\rvert P_{t'}(x')P_{t'}(x')\d t'\d x'
\leq \frac{ C(T,f,m)
}{1+\lvert x\rvert ^m},\label{I22:bound}
\end{align}
and $R_{\ref{key}}(y,x,s,T)=R_{\ref{key},\vep}(y,x,s,T)$, independent of $m$, satisfies the following estimate:
\begin{align}
\begin{split}\label{R:bdd}
\lvert R_{\ref{key}}(y,x,s,T)\rvert &\leq \frac{C(\lambda,T,f,m,M)
}{(\log^2\vep^{-1})(1+\lvert x\rvert ^m)}\left(\lvert \log (T-s)\rvert +1+\lvert \log \lvert y \rvert\rvert\right)\\
&\quad+\frac{C(\lambda,T,f,m,M)
}{(\log \vep^{-1})(1+\lvert x\rvert ^m)}   \vep_{\ref{limT:glue}}
\left(\frac{4(T-s)\vep^{-2}}{\lvert y \rvert^2}\right),
\end{split}
\end{align}
where $0\leq \vep_{\ref{limT:glue}}(a)\less a^{-1}\1(a\geq 1)+(1-\log a)\1 (a<1)$. Moreover, 
\begin{align}
\begin{split}
\lvert \mathring{\mathfrak L}^1_\vep f(y,x,s,T)-f(x,s)\rvert &\leq \frac{C(\lambda,T,f,m,M)
}{(\log\vep^{-1})(1+\lvert x\rvert ^m)}\\
&\quad \times \left[\lvert \log (T-s)\rvert+1+\lvert\log \lvert y \rvert\rvert+\vep_{\ref{limT:glue}}
\left(\frac{4(T-s)\vep^{-2}}{\lvert y \rvert^2}\right)\right].\label{Phi1diff:bdd}
\end{split}
\end{align}

\noindent {\rm (2$\cc$)} For all $(y,x)\in \R^2\times \R^2$, we have
\begin{align*}
\nabla_{x}\mathring{\mathfrak L}^1_\vep f(y,x,s,T)&=\lv\int_{x'}\left(\int_{t'=0}^{T-s}\nabla_{x}f(x-x',s+t')P_{t'}(\vep y+x')P_{t'}(x')\d t'\right)\d x',\;0\leq s\leq T.
\end{align*}
In particular, for all $y,x\in \R^2$, $0<s<T$, and $m\in \Bbb N$, 
\begin{align}
\lvert\nabla_{x}\mathring{\mathfrak L}^1_\vep f(y,x,s,T)\rvert \leq \frac{C(f,\lambda,m)}{1+\lvert x\rvert ^m}\left[\lvert\log (T-s)\rvert+1+\lvert\log \lvert y \rvert\rvert+\vep_{\ref{limT:glue}}
\left(\frac{4(T-s)\vep^{-2}}{\lvert y \rvert^2}\right)\right].\label{I1:derbdd}
\end{align}

\noindent {\rm (3$\cc$)} Suppose that $m\geq 5$. For $0<s<T$, it holds that
\begin{align}\label{EMremainder:bdd}
\begin{aligned}
\int_{y}\frac{1}{1+\lvert y \rvert^m}\left|\vep_{\ref{limT:glue}}
\left(\frac{4(T-s)\vep^{-2}}{\lvert y \rvert^2}\right)\right|\d y &\less  C(m)\min\left\{\frac{1}{4(T-s)\vep^{-2}},1\right\}\\
&\quad +C(m)\lvert\log [4(T-s)\vep^{-2}]\rvert\1_{\{4(T-s)\vep^{-2}<1\}}. 
\end{aligned}
\end{align}
\end{lem}

\begin{proof}[Proof of Lemma~\ref{lem:I1}]
{(1$\cc$)} By \eqref{def:L1r}, we can write
\begin{align}
\mathring{\mathfrak L}^1_\vep f(y,x,s,T)&=\Lambda_\vep\int_{t=s}^T\int_{x_0}f(x_0,t)P_{t-s}(\vep y+x,x_0) P_{t-s}(x_0,x)\d x_0 \d t\notag\\
&=\Lambda_\vep \int_{t'=0}^{T-s}\int_{x'} f(x-x',s+t')P_{t'}(\vep y+x')P_{t'}(x')\d x'\d t'\label{asymp11prep}\\
&=I^1_{\ref{asymp11}} f(x,s)+\lv I^2_{\ref{asymp11}}+\lv  I^3_{\ref{asymp11}},\label{asymp11}
\end{align}
where
\begin{align}
I^1_{\ref{asymp11}}&\,\defeq\,\Lambda_\vep \int_{t'=0}^{T-s}\int_{x'}P_{t'}(\vep y+x')P_{t'}(x')\d x'\d t',\label{asymp11prep-1}\\
I^2_{\ref{asymp11}}&\,\defeq\, \int_{x'}\int_{t'=0}^{T-s}[f(x-x',s+t')-f(x,s)]P_{t'}(x')P_{t'}(x')\d t'\d x',\notag\\
I^3_{\ref{asymp11}}&\,\defeq\, \int_{x'}\int_{t'=0}^{T-s}[f(x-x',s+t')-f(x,s)]  [P_{t'}(\vep y+x')-P_{t'}(x')]P_{t'}(x')\d t'\d x'.\notag
\end{align}
The remainder of the proof of (1$\cc$) proceeds as follows. 
Steps~1--3 below obtain, respectively: a second-order expansion of $I^1_{\ref{asymp11}}$ in terms of $1/\log\vep^{-1}$; the absolute convergence of $I^2_{\ref{asymp11}}$ as an integral; and a bound on $\lv I^3_{\ref{asymp11}}$ of order $\lv \vep^{1/2}$. Additionally, these results proven for $I^2_{\ref{asymp11}}$ and $\lv I^3_{\ref{asymp11}}$ include bounds that decay polynomially in $\lvert x\rvert $ as $\lvert x\rvert \to\infty$. Finally, Step~4 concludes the proof. \medskip 

\noindent {\bf Step 1.}
To obtain the required property of $I^1_{\ref{asymp11}}$, we first recall the definition \eqref{def:lv} of $\lv$. We then apply Lemma~\ref{lem:heat1} (1$\cc$) after integrating out $x'$, applying the Chapman-Kolmogorov equation, and making the substitution $t'' = \vep^{-2} t'$:
\begin{align}
I^1_{\ref{asymp11}}&=\lv  \int_{t''=0}^{(T-s)\vep^{-2}}P_{2t''}( y)\d t''\notag\\
&=\biggl(\frac{2\pi}{\log \vep^{-1}}+\frac{2\pi\lambda}{\log^2 \vep^{-1}}\biggr) \biggl[\frac{\log [4(T-s)\vep^{-2}]}{4\pi} -\frac{\gamma_{\sf EM}}{4\pi}-\frac{\log \lvert y \rvert}{2\pi}+ \vep_{\ref{limT:glue}}
\left(\frac{4(T-s)\vep^{-2}}{\lvert y \rvert^2}\right)\biggr]\notag\\
\begin{split}\label{asymp11-1}
&=1+\frac{1}{\log \vep^{-1}}
\biggl(\frac{\log[4 (T-s)]}{2}+\lambda -\frac{\gamma_{\sf EM}}{2}-\log \lvert y \rvert\biggr)\\
&\quad +\frac{\lambda}{\log^2\vep^{-1}}\biggl(\frac{\log [4(T-s)]}{2} -\frac{\gamma_{\sf EM}}{2}-\log \lvert y \rvert\biggr) +\lv   \vep_{\ref{limT:glue}}
\left(\frac{4(T-s)\vep^{-2}}{\lvert y \rvert^2}\right).
\end{split}
\end{align}

\noindent {\bf Step 2.}
To obtain the absolute convergence of $I^2_{\ref{asymp11}}$ as an integral, we apply the assumption $f\in \C^1_{pd}(\R^2\times [0,T])$ and the mean-value theorem
 to get the following bound for any $m\in \Bbb Z_+$: 
\begin{align}
&\sup_{\scriptstyle 0\leq s\leq T\atop \scriptstyle s+t'\leq T}\lvert f(x-x',s+t')-f(x,s)\rvert \leq (\lvert x'\rvert+t')\left(\frac{C(T,f,m
)}{1+\lvert x-x'\rvert^m}+\frac{C(T,f,m
)}{1+\lvert x\rvert ^m}\right).\label{varphi:diff}
\end{align}
When applied to $I^2_{\ref{asymp11}}$,
these bounds allow \eqref{ineq:heat2-1} with $(p_1,p_2,p)=(1,1,m)$ and $(p_1,p_2,p)=(1,1,0)$. Therefore, \eqref{I22:bound} holds, and
$I^2_{\ref{asymp11}}$ defines an absolutely convergent integral with
\begin{align}
\lvert I^2_{\ref{asymp11}}\rvert\leq \frac{C(T,f,m)}{1+\lvert x\rvert^m}.\label{I22:boundproof}
\end{align}
Moreover, by the definition \eqref{def:lv} of $\lv$ along with \eqref{I22:bound}, we get 
\begin{align}\label{asymp11-2}
\lv I^2_{\ref{asymp11}}=\frac{2\pi}{\log \vep^{-1}}I^2_{\ref{asymp11}}+\frac{2\pi}{\log^2\vep^{-1}}\mathcal O_{\lambda,T,f,m}\left(\frac{1}{1+\lvert x\rvert ^m}\right).
\end{align}

\noindent {\bf Step~3.}
We bound $I^3_{\ref{asymp11}}$ in this step. 
For fixed $\delta_\vep\in (0,\infty)$ such that $\vep/\delta_\vep^{1/2}\leq 1$, we write the integral in $I^3_{\ref{asymp11}}$ 
as the sum of two integrals according to $1=\1_{(0,\delta_\vep)}(t')+\1_{[\delta_\vep,\infty)}(t')$. The bound in \eqref{varphi:diff} is also applied again. 
Hence, for any $m\in \Bbb Z_+$ and $0<\lvert y \rvert\leq M$,
\begin{align*}
\lvert I^3_{\ref{asymp11}}\rvert &\leq C(T,f,m)\int_{x'}\int_{t'=0}^{T-s}\1_{(0,\delta_\vep)}(t')(\lvert x'\rvert+t')\left(\frac{1}{1+\lvert x-x'\rvert^m}+\frac{1}{1+\lvert x\rvert ^m}\right)
\\
&\quad \times [P_{t'}(\vep y+x')+P_t(x')]P_{t'}(x')\d t'\d x'\\
&\quad +C(T,f,m,M)\frac{\vep}{\delta_\vep^{1/2}}\int_{x'}\int_{t'=0}^{T-s}\1_{[\delta_\vep,\infty)}(t')(\lvert x'\rvert+t')\\
&\quad \times \left(\frac{1}{1+\lvert x-x'\rvert^m}+\frac{1}{1+\lvert x\rvert ^m}\right)P_{2t'}(x')P_{t'}(x')\d t'\d x',
\end{align*}
where the last term uses \eqref{gauss:vep}. We take $\delta_\vep=\vep$. Then by \eqref{ineq:heat2-1} of Lemma~\ref{lem:heat2},
\begin{align}\label{asymp11-3}
\lv I^3_{\ref{asymp11}}=\lv \vep^{1/2}\mathcal O_{T,f,m,M}\biggl(\frac{1}{1+\lvert x\rvert ^m}\biggr).
\end{align}

\noindent {\bf Step 4.}
Finally, we get \eqref{key} by applying 
\eqref{asymp11-1}, \eqref{asymp11-2} and \eqref{asymp11-3} to
\eqref{asymp11} with 
\begin{align}
R_{\ref{key}}(y,x,s,T)&\,\defeq\,\frac{\lambda }{\log^2\vep^{-1}}\biggl(\frac{\log [4(T-s)]}{2} -\frac{\gamma_{\sf EM}}{2}-\log \lvert y \rvert\biggr)f(x,s)\notag\\
&\quad +\lv   \vep_{\ref{limT:glue}}
\left(\frac{4(T-s)\vep^{-2}}{\lvert y \rvert^2}\right)f(x,s) +\frac{2\pi}{\log^2\vep^{-1}}\mathcal O_{\lambda,T,f,m
}\left(\frac{1}{1+\lvert x\rvert ^m}\right)\notag\\
&\quad +\lv \vep^{1/2}\mathcal O_{T,f,m,M
}\biggl(\frac{1}{1+\lvert x\rvert ^m}\biggr).\notag
\end{align}
The remaining properties in (1$\cc$) can be obtained as follows: \eqref{I22:bound} has been established as \eqref{I22:boundproof}; the function $R_{\ref{key}}(y,x,s,T)$ satisfies \eqref{R:bdd} by direct application of the previous display and the assumption $f\in \C^1_{pd}(\R^2\times [0,T])$; and the bound in \eqref{Phi1diff:bdd} is obtained from \eqref{key} using the assumption $f\in \C^1_{pd}(\R^2\times [0,T])$, together with \eqref{I22:bound} and \eqref{R:bdd}. \medskip 

\noindent {(2$\cc$)} Both derivative formulas follow from the integral representation given in \eqref{asymp11prep}. Specifically, the time derivative formula is obtained by applying the Leibniz integral rule. Furthermore, the bound in \eqref{I1:derbdd} holds since, by the first derivative formula and \eqref{ratiot}, we have
\begin{align*}
\lvert\nabla_{x}\mathring{\mathfrak L}^1_\vep f(y,x,s,T)\rvert &\leq \lv\int_{x'}\left(\int_{t'=0}^{T-s}\frac{C(f,m,T)(1+\lvert x'\rvert^m)}{1+\lvert x\rvert ^m}
P_{t'}(\vep y+x')P_{t'}(x')\d t'\right)\d x',
\end{align*}
so that the bound in \eqref{asymp11-1} ($I^1_{\ref{asymp11}}$ is defined by \eqref{asymp11prep-1}) and \eqref{ineq:heat2-1} apply. 
\medskip 

\noindent {(3$\cc$)}
As stated below \eqref{R:bdd}, $0\leq \vep_{\ref{limT:glue}}(a)\less a^{-1}\1(a\geq 1)+(1-\log a)\1 (a<1)$. Hence, 
\begin{align}
&\int_{y}\frac{1}{1+\lvert y \rvert^m}\left|\vep_{\ref{limT:glue}}
\left(\frac{4(T-s)\vep^{-2}}{\lvert y \rvert^2}\right)\right|\d y\notag\\
&\quad \leq \int_y\frac{1}{1+\lvert y \rvert^m}\left(\frac{4(T-s)\vep^{-2}}{\lvert y \rvert^2}\right)^{-1}\1\left(\frac{4(T-s)\vep^{-2}}{\lvert y \rvert^2}\geq 1\right)\d y\notag\\
&\quad\quad +\int_y \frac{1}{1+\lvert y \rvert^m}\1\left( \frac{4(T-s)\vep^{-2}}{\lvert y \rvert^2}< 1\right)\d y\notag\\
&\quad\quad + \int_y \frac{1}{1+\lvert y \rvert^m}\left(-\log \frac{4(T-s)\vep^{-2}}{\lvert y \rvert^2}\right)\1\left(\frac{4(T-s)\vep^{-2}}{\lvert y \rvert^2}< 1\right)\d y.\label{EMremainder:bdd1}
\end{align} 
Then observe that for $m\geq 5$,  the first and second integrals on the right-hand side are bounded by $C(m)\min\{1/[4(T-s)\vep^{-2}],1\}$. In particular, the bound $C(m)/[4(T-s)\vep^{-2}]$ applies to the second integral since for $A>0$,
\[
\int_y\frac{1}{1+\lvert y \rvert^m}\1\left(\frac{A}{\lvert y \rvert^2}<1\right)\d y\less \int_{r=\sqrt{A}}^\infty \frac{r}{1+r^m}\d r\less \int_{r=\sqrt{A}}^\infty \frac{1}{r^3}\d r\less \frac{1}{A}. 
\]
To bound the last integral in \eqref{EMremainder:bdd1}, note that for any $A\geq 1/2$, 
\begin{align}
\int_{\frac{A}{\lvert y \rvert^2}<1} \frac{1}{1+\lvert y \rvert^m}\left\lvert\log \frac{A}{\lvert y \rvert^2}\right|\d y\leq \int_{\frac{\lvert y \rvert^2}{A}>1}\frac{1}{1+\lvert y \rvert^m}\left|\frac{\lvert y \rvert^2}{A}\right|\d y\leq C(m)\min\left\{\frac{1}{A},1\right\},\notag
\end{align}
where the first inequality follows since $\log x\less x $ for all $x\geq 1$, and the last inequality holds as $m\geq 5$ by assumption. Additionally, for any $0<A<1/2$,
\begin{align}
\int_{\frac{A}{\lvert y \rvert^2}<1} \frac{1}{1+\lvert y \rvert^m}\left\lvert\log \frac{A}{\lvert y \rvert^2}\right|\d y&\leq \int_y \frac{\lvert \log A\rvert+\lvert\log \lvert y \rvert^2\rvert}{1+\lvert y \rvert^m}\d y\notag\\
&\leq C(m)(\lvert \log A\rvert +1)\leq C(m)\lvert \log A\rvert .\notag 
\end{align}
Hence, \eqref{EMremainder:bdd1} is sufficient to obtain \eqref{EMremainder:bdd}.
\end{proof}

\begin{proof}[Proof of Proposition~\ref{prop:I1}]
We begin by proving the following two limits: 
\begin{align}
&\lim_{n\to\infty}\E\left[\int_{s=0}^T\int_{x}\left(\int_y \lvert R_{\ref{key},n}(y,x,s,T)\rvert \phi(y)\d y\right)X_n(x,s)^2\d x\d s\right]=0,\label{I1:Rconv1}\\
&\lim_{n\to\infty}\E\biggl[\int_{s=0}^T\int_{x}\biggl|-\frac{1}{2\pi}
\biggl(\frac{\log [4(T-s)]}{2}+\lambda -\frac{\gamma_{\sf EM}}{2}-\int_{y}\phi(y)\log \lvert y \rvert\d y\biggr)f(x,s)\notag\\
&\quad -\int_{x'}\int_{t'=0}^{T-s}[f(x-x',s+t')-f(x,s)]P_{t'}(x')P_{t'}(x')\d t'\d x'\biggr|
\frac{X_n(x,s)^2}{\log^2\vep_n^{-1}} \d x\d s\biggr]=0.\label{I1:Rconv2}
\end{align}
Here, \eqref{I1:Rconv1} follows from Theorem~\ref{thm:mombdd} (3$\cc$), \eqref{R:bdd}, and Lemma~\ref{lem:I1} (3$\cc$), and \eqref{I1:Rconv2} follows from Theorem~\ref{thm:mombdd} (3$\cc$), \eqref{I22:bound} and the assumption $f(x,s)\in \C_{pd}^1(\R^2\times [0,T])$.  

We are ready to prove the required limit of $I^1_n(f,T)$. Recall that we have $\mathfrak L^1_nf(y,x,s,T)=\mathring{\mathfrak L}^1_nf(y,x,s,T)\phi(y)$ by \eqref{def:L1}, where $\phi$ is a probability density by Proposition~\ref{prop:XM} (1$\cc$). Hence, by \eqref{def:I1}, we can rewrite $I^1_n(f,T)$ as
\[
I^{1}_{n}(f,T)=\int_{s=0}^T \int_{x}\left( \int_{y}[f(x,s)-\mathring{\mathfrak L}_n^1f (y,x,s,T)]\phi(y)\d y\right) X_n(x,s)^2\d x\d s.
\]
With this expression,
\eqref{key} and \eqref{I1:Rconv1} imply the first limit in probability below:
\begin{align*}
&\P\mbox{-}\lim_{n\to\infty}I^1_n(f,T)\\
&\quad =\P\mbox{-}\lim_{n\to\infty}
\int_{s=0}^T\int_{x}\biggl[-\frac{1}{\log \vep_n^{-1}}
\biggl(\frac{\log [4(T-s)]}{2}+\lambda -\frac{\gamma_{\sf EM}}{2}-\int_y\phi(y) \log \lvert y \rvert\d y\biggr)f(x,s)\\
&\quad\quad -\frac{2\pi  }{\log \vep_n^{-1}}\int_{x'}\int_{t'=0}^{T-s}[f(x-x',s+t')-f(x,s)]P_{t'}(x')P_{t'}(x')\d t'\d x'\biggr]X_n(x,s)^2\d x\d s\\
&\quad=\P\mbox{-}\lim_{n\to\infty}
\frac{2\pi}{\log \vep_n^{-1}}\int_{s=0}^T\int_{x}\biggl[-\frac{1}{2\pi}
\biggl(\frac{\log [4(T-s)]}{2}+\lambda -\frac{\gamma_{\sf EM}}{2}-\int_y\phi(y)\log \lvert y \rvert\d y\biggr)f(x,s)\\
&\quad\quad -\int_{x'}\int_{t'=0}^{T-s}[f(x-x',s+t')-f(x,s)]P_{t'}(x')P_{t'}(x')\d t'\d x'\biggr] X_n(x,s)^2\d x\d s\\
&\quad =\P\mbox{-}\lim_{n\to\infty}\int_{s=0}^T\int_{x}\biggl[-\frac{1}{2\pi}
\biggl(\frac{\log [4(T-s)]}{2}+\lambda -\frac{\gamma_{\sf EM}}{2}-\int_{y}\phi(y)\log \lvert y \rvert\d y\biggr)f(x,s)\\
&\quad\quad -\int_{x'}\int_{t'=0}^{T-s}[f(x-x',s+t')-f(x,s)]P_{t'}(x')P_{t'}(x')\d t'\d x'\biggr] \nu_n(\d x,\d s)\\
&\quad =\int_{s=0}^T\int_{x}\mathfrak L_0f(x,s,T) \nu_\infty(\d x,\d s).
\end{align*}
Here, the third equality follows from  \eqref{I1:Rconv2}, the definition \eqref{def:lv} of $\Lambda_{n}=\Lambda_n$, and the relationship in \eqref{def:covarprocessnu} between $X_n=X_{\vep_n}$ and $\nu_{n}=\nu_{\vep_n}$. Additionally, 
the last equality can be deduced from the definition \eqref{def:L0f} of $\mathfrak L_0f(x,s,T)$ and Corollary~\ref{cor:muinfty}. In more detail, Corollary~\ref{cor:muinfty} is applicable in this context for two reasons. First, the singularity of $\log (T-s)$ at $s=T$ does not pose an obstruction, as Theorem~\ref{thm:mombdd} (3$\cc$) implies 
\[
\limsup_{\eta\searrow 0+}\limsup_{n\to\infty}\E\left[\int_{s=T-\eta}^T\int_x \frac{\lvert \log (T-s)\rvert}{1+\lvert x\rvert ^m}\Lambda_n X_n(x,s)^2\d x\d s\right]=0.
\]
Second,
\[
(x,s)\in \R^2\times [0,T]\mapsto \int_{x'}\int_{t'=0}^{T-s}[f(x-x',s+t')-f(x,s)]P_{t'}(x')P_{t'}(x')\d t'\d x'\in \C^0_{pd}(\R^2\times [0,T]),
\]
where the continuity can be deduced from dominated convergence, thanks to \eqref{varphi:diff} and \eqref{ineq:heat2-1}, and the polynomial decay of any order follows from \eqref{I22:bound}. The proof of Proposition~\ref{prop:I1} is complete. 
\end{proof}

\addtocontents{toc}{\protect\setcounter{tocdepth}{1}}
\subsection{Limit theorem for $\bs I^{\bs 2}_{\bs n}\bs (\bs f\bs ,\bs T\bs)$}\label{sec:I2}
\addtocontents{toc}{\protect\setcounter{tocdepth}{2}}
\addcontentsline{toc}{subsection}{\numberline{\thesubsection}{Limit theorem for $I^2_n(f,T)$}}

The following proposition specifies the limit and
also establishes relations among the three random measures $\nu_\infty,\mu_\infty,\mathring{\mu}_\infty$ from Corollary~\ref{cor:Sk}.

\begin{prop}\label{prop:I2}
Fix $0<T<\infty$. For all nonnegative $ g\in \C(\R^2\times \R^2\times [0,T])$, 
\begin{align}
\int_{s=0}^T\int_{x}\int_{x'}g(x',x,s)\mu_\infty(\d x',\d x,\d s)
&=\int_{s=0}^T\int_{x}g(x,x,s)\nu_\infty(\d x,\d s),\label{id:mu=nu}\\
\int_{s=0}^T\int_{x}\int_{x'}g(x',x,s)\mathring{\mu}_\infty(\d x',\d x,\d s)
&=\int_{s=0}^T\int_x \left(\int_{x'} g(x',x,s)\phi(x')\d x'\right)\nu_\infty(\d x,\d s).\label{id:mur=nu}
\end{align}
Additionally, it holds that for all $f\in \C_{pd}^1(\R^2\times[0,T])$, 
\begin{align*}
I^2_n(f,T)&\,\defeq\,\int_{s=0}^T \int_{x}\int_{y} \mathfrak L^1_nf (y,x,s,T) X_n(x,s)^2\d y\d x\d s\\
&\quad  -\int_{s=0}^T \int_{x}\int_{y} \mathfrak L^1_nf (y,x,s,T) X_n(\vep_n y+x,s)X_n(x,s)\d y\d x\d s\\
&\quad \xrightarrow[n\to\infty]{\P}  I^2_\infty(f,T)\,\defeq\,\int_{r'=0}^T \int_{y'}\biggl(\frac{1}{2\pi }\int_{y''}\int_{y}\phi(y'')(\log\rvert y''-y\rvert )\phi(y)\d y \d y''\\
&\quad \quad   -\frac{1}{2\pi }\int_{y''}\int_{y}\phi(y'')(\log\lvert y'' \rvert)\phi(y)\d y \d y''\biggr) f(y',r')\nu_\infty(\d y',\d r'),
\end{align*}
where $\mathfrak L^1_n\,\defeq\,\mathfrak L^1_{\vep_n}$ for $\mathfrak L^1_\vep $ defined in \eqref{def:L1}.
\end{prop}

The terms $I^1_n(f,T)$ and $I^2_n(f,T)$ share a similar structure, as each is defined as a difference of integrals. While $I^1_n(f,T)$ involves only the measure $X_n(x,s)^2\d x\d s$, $I^2_n(f,T)$ incorporates both $X_n(x,s)^2\d x\d s$ and $X_n(\vep_n y + x,s)X_n(x,s)\d x\d s$. Accordingly, our approach for analyzing $I^2_n(f,T)$ adapts the methods from Proposition~\ref{prop:I1}. In particular, we employ the Fubini theorem and the stochastic Fubini theorem to bring Lebesgue integrals inside stochastic integrals, thereby absorbing pointwise singularities in the limit. Nevertheless, determining the precise limit of $I^{2}_{n}(f,T)$ still necessitates a far more refined analysis. This increased complexity stems from perturbing the spatial variable in $X_n$, specifically by using $X_n(\vep_n y+x,s)X_n(x,s)\d x\d s$.

Lemmas~\ref{lem:I12} and \ref{lem:I2ap} below form the core technical machinery for Proposition~\ref{prop:I2}, formulated under the general parameter $\vep$ rather than $\vep_n$, as the analysis depends only on the setting established before Corollary~\ref{cor:Sk}. Lemma~\ref{lem:I12} analyzes the difference of integrals defining $I^{2}_{\vep}(f,T)$ by replacing $\mathfrak L^1_\vep f(y,x,s,T)$ with a general function $g(y,x,s,T)$ to simplify notation. Using the mild formulation \eqref{def:mild:vep} to expand $X_\vep(\vep y+x,s)-X_\vep(x,s)$, together with the Fubini arguments considered above, we obtain the decomposition \eqref{dec:I12} involving the refined functionals $I^{2,k}_{\vep}(g,T)$. Lemma~\ref{lem:I2ap} then investigates the limits of these functionals and related terms. At the conclusion of the proof of Proposition~\ref{prop:I2}, we set $g = \mathfrak L^1_\vep f$. Ultimately, starting from the original expression for $I^{2}_{\vep}(f,T)$, the goal of this proposition is to transfer the spatial perturbation in $X_\vep(\vep y+x,s)$ to the spatial variables of $\mathfrak L^1_\vep f(y,x,s,T)$. The decomposition using $I^{2,k}_{\vep}(g,T)$ facilitates this transfer, allowing Lemma~\ref{lem:I1} to be applied where appropriate to analyze the limiting behavior.

\begin{lem}\label{lem:I12}
For all $T\in (0,\infty)$ and ${g}={g}(y,x,s,T)$ with $(y,x,s)\mapsto {g}(y,x,s,T)\in L^1(\R^2\times \R^2\times [0,T])$, 
\begin{align}
\begin{split}\label{dec:I12}
&\int_{s=0}^T\int_{x}\int_{y}{g}(y,x,s,T) [X_\vep(\vep y+x,s)-X_\vep(x,s)]X_\vep(x,s)\d y\d x \d s\\
&\quad =I^{2,1}_{\vep}({g},T)+I^{2,2}_{\vep}({g},T)+I^{2,3}_{\vep}({g},T)+I^{2,4}_{\vep}({g},T).
\end{split}
\end{align}
Here, the left-hand side of \eqref{dec:I12} defines an absolutely convergent integral, as justified by \eqref{Xvep:mom}, and
\begin{align*}
I^{2,1}_{\vep}({g},T)&\,\defeq\,\int_{s=0}^T \int_{x}\int_{y}{g}(y,x,s,T) 
[P_sX_0(\vep y+x)-P_sX_0(x)] P_sX_0(x)\d y\d x\d s,\\
I^{2,2}_{\vep}({g},T)&\,\defeq\,
\int_{r'=0}^T\int_{y'}  \int_{z''} \int_{s=r'}^T \int_{x}\int_{y}{g}(y,x,s,T) [P_{s-r'}(\vep y+x,y')-P_{s-r'}(x,y')] \\
&\quad \times P_{s-r'}(x,z'') \d y\d x \d s X_\vep(z'',r')
\d z''   \M_\vep(\d y',\d r')
,\\
I^{2,3}_{\vep}({g},T)&\,\defeq\,  \int_{r''=0}^T \int_{y''} \int_{z'} \int_{s=r''}^T \int_{x}\int_{y} {g}(y,x,s,T)
[P_{s-r''}(\vep y+x,z')-P_{s-r''}(x,z')]\\
 &\quad \times P_{s-r''}(x,y'')\d y \d x\d sX_\vep(z',r'')\d z' \M_\vep(\d y'',\d r''),\\
I^{2,4}_{\vep}({g},T)&\,\defeq\,
 \int_{r'=0}^T \int_{y'}\int_{y''} \int_{s=r'}^T \int_{x}\int_{y}{g} (y,x,s,T) 
 [P_{s-r'}(\vep y+x,\vep y''+y')-P_{s-r'}(x,\vep y''+y')]\\
 &\quad \times P_{s-r'}(x,y')\d y \d x\d s \lv  X_\vep(\vep y''+y',r') X_\vep(y',r') \phi(y'')\d y''\d y'
\d r',
\end{align*}
where all Lebesgue integrals are absolutely convergent, and all stochastic integrals are well-defined in It\^{o}'s sense. Moreover, \eqref{dec:I12} applies with the choice of ${g}(y,x,s,T)\equiv \mathfrak L_\vep^1f(y,x,s,T)$ for $f(x,s)\in \C_{pd}^0(\R^2\times [0,T])$, in which case
 $I^{2}_{\vep}(f,T)$ defined by \eqref{def:I2} satisfies the following decomposition: 
 \begin{align}\label{dec0:I2}
 I^{2}_{\vep}(f,T)=-I^{2,1}_{\vep}(\mathfrak L_\vep^1f,T)-I^{2,2}_{\vep}(\mathfrak L_\vep^1f,T)-I^{2,3}_{\vep}(\mathfrak L_\vep^1f,T)-I^{2,4}_{\vep}(\mathfrak L_\vep^1f,T).
 \end{align}
\end{lem}

\begin{proof} 
We begin by proving \eqref{dec:I12}. To do so, we first expand $[X_\vep(\vep y + x,s) - X_\vep(x,s)] X_\vep(x,s)$ using the mild form \eqref{def:mild:vep}, which yields the following equations:
\begin{align*}
X_\vep(\vep y+x,s)-X_\vep(x,s)
&=[P_sX_0(\vep y+x)-P_sX_0(x)]\\
&\quad+\int_{r'=0}^s \int_{y'}[P_{s-r'}(\vep y+x,y')-P_{s-r'}(x,y')]\M_\vep(\d y',\d r'),\\
X_\vep(x,s)&=P_sX_0(x)+\int_{r'=0}^s \int_{y'} P_{s-r'}(x,y')\M_\vep(\d y',\d r').
\end{align*}
We regard the right-hand sides as $a+b$ and $c+d$ and use a direct multiplication: $(a+b)(c+d)=a c+bc+a d+b d$. The product of the two stochastic integrals is further expanded by using integration by parts. Hence, we obtain the following expansion: 
\begin{align} 
& [X_\vep(\vep y+x,s)-X_\vep(x,s)]X_\vep(x,s)\notag\\
&\quad =[P_sX_0(\vep y+x)-P_sX_0(x)] P_sX_0(x)\notag\\
&\quad \quad +\int_{r'=0}^s \int_{y'}[P_{s-r'}(\vep y+x,y')-P_{s-r'}(x,y')]\M_\vep(\d y',\d r')P_sX_0(x)\notag\\
&\quad\quad  +[P_sX_0(\vep y+x)-P_sX_0(x)] \int_{r'=0}^s \int_{y'} P_{s-r'}(x,y')\M_\vep(\d y',\d r')\notag\\
&\quad\quad + \int_{r'=0}^s\int_{y'} [P_{s-r'}(\vep y+x,y')-P_{s-r'}(x,y')] \notag\\
&\quad\quad  \times \biggl(\int_{r''=0}^{r'} \int_{y''} P_{s-r''}(x,y'') \M_\vep(\d y'',\d r'')\biggr)\M_\vep(\d y',\d r')\notag\\
&\quad\quad  +  \int_{r''=0}^s \int_{y''}\biggl(\int_{r'=0}^{r''} \int_{y'} [P_{s-r'}(\vep y+x,y')-P_{s-r'}(x,y')] \M_\vep(\d y',\d r')\biggr)\notag \\
&\quad\quad  \times P_{s-r''}(x,y'')\M_\vep(\d y'',\d r'')\notag\\
&\quad \quad +  \int_{r'=0}^s \int_{y'}\int_{y''} [P_{s-r'}(\vep y+x,\vep y''+y')-P_{s-r'}(x,\vep y''+y')] P_{s-r'}(x,y')\notag\\
&\quad\quad \times  \lv  X_\vep(\vep y''+y',r')X_\vep(y',r')\phi(y'')\d y''\d y'\d r'.\label{DeltaXX}
\end{align}
To elaborate, the final three terms result from applying integration by parts to stochastic integrals, 
with the last term specifically using the formula in \eqref{eq:covarM-2}.

We proceed to the weak formulations of both sides of \eqref{DeltaXX}. 
With ${g}={g}(y,x,s,T)$ to lighten notation, integration against ${g} \d y\d x \d s$ over $0\leq s\leq T$ and $y,x\in \R^2$ for both sides of \eqref{DeltaXX} yields
\begin{align*}
&\int_{s=0}^T \int_{x}\int_{y}{g}[X_\vep(\vep y+x,s)-X_\vep(x,s)]X_\vep(x,s)
\d y\d x\d s\\
&\quad =I^{2,1}_{\vep}({g},T)+\tilde{I}^{2,2}_{\vep}({g},T)+\tilde{I}^{2,3}_{\vep}({g},T)+\tilde{I}^{2,4}_{\vep}({g},T)+\tilde{I}^{2,5}_{\vep}({g},T)+\tilde{I}^{2,6}_{\vep}({g},T),
\end{align*}
where $I^{2,1}_{\vep}({g},T)$ has been defined below \eqref{dec:I12}, and the remaining terms are defined as folllows: 
\begin{align*}
\tilde{I}^{2,2}_{\vep}({g},T)&\,\defeq\,\int_{s=0}^T \int_{x}\int_{y}\int_{r'=0}^s \int_{y'}{g} [P_{s-r'}(\vep y+x,y')-P_{s-r'}(x,y')] P_sX_0(x)\M_\vep(\d y',\d r')\d y\d x\d s,\\
\tilde{I}^{2,3}_{\vep}({g},T)&\,\defeq\,\int_{s=0}^T \int_{x}\int_{y}\int_{r'=0}^s \int_{y'} {g}[P_sX_0(\vep y+x)-P_sX_0(x)] P_{s-r'}(x,y')\M_\vep(\d y',\d r') \d y\d x\d s,\\
\tilde{I}^{2,4}_{\vep}({g},T)&\,\defeq\,
\int_{s=0}^T \int_{x}\int_{y}
\int_{r'=0}^s\int_{y'}{g} [P_{s-r'}(\vep y+x,y')-P_{s-r'}(x,y')] \\
&\quad\times \biggl(\int_{r''=0}^{r'} \int_{y''}P_{s-r''}(x,y'')\M_\vep(\d y'',\d r'')\biggr)   \M_\vep(\d y',\d r')
\d y\d x \d s,\\
\tilde{I}^{2,5}_{\vep}({g},T)&\,\defeq\, \int_{s=0}^T \int_{x}\int_{y} \int_{r''=0}^s \int_{y''} {g} \biggl( \int_{r'=0}^{r''} \int_{y'}[P_{s-r'}(\vep y+x,y')-P_{s-r'}(x,y')]  \M_\vep(\d y',\d r')\biggr)\\
 &\quad \times P_{s-r''}(x,y'') \M_\vep(\d y'',\d r'')\d y \d x\d s,\\
\tilde{I}^{2,6}_{\vep}({g},T)&\,\defeq\, \int_{s=0}^T \int_{x}\int_{y}
\int_{r'=0}^s \int_{y'}\int_{y''}{g}  [P_{s-r'}(\vep y+x,\vep y''+y')-P_{s-r'}(x,\vep y''+y')] P_{s-r'}(x,y')\\
&\quad \times \lv X_\vep(\vep y''+y',r') X_\vep(y',r')\phi(y'')\d y''\d y'\d r'
\d y \d x\d s.
\end{align*}
By the Fubini theorem, 
$\tilde{I}^{2,6}_\vep({g},T)=I^{2,4}_\vep({g},T)$.
To complete the proof of  \eqref{dec:I12}, we show below that
\begin{align}
\tilde{I}^{2,2}_{\vep}({g},T)+\tilde{I}^{2,4}_{\vep}({g},T)&=I^{2,2}_{\vep}({g},T),\label{eq:I2-1}\\
\tilde{I}^{2,3}_{\vep}({g},T)+\tilde{I}^{2,5}_{\vep}({g},T)&=I^{2,3}_{\vep}({g},T).\label{eq:I2-2}
\end{align}

To prove \eqref{eq:I2-1}, we first consider the following equality, which is obtained by
applying the Chapman--Kolmorgorov equation to $P_{s-r''}(x,y'')$:
\begin{align*}
\tilde{I}^{2,4}_{\vep}({g},T)
&=
\int_{s=0}^T \int_{x}\int_{y}
\int_{r'=0}^s\int_{y'}{g} [P_{s-r'}(\vep y+x,y')-P_{s-r'}(x,y')] \\
&\quad \biggl(\int_{r''=0}^{r'} \int_{y''}\int_{z''}P_{s-r'}(x,z'')P_{r'-r''}(z'',y'')\d z''\M_\vep(\d y'',\d r'')\biggr)   \M_\vep(\d y',\d r')
\d y\d x \d s\\
&=
\int_{s=0}^T \int_{x}\int_{y}
\int_{r'=0}^s\int_{y'}{g} [P_{s-r'}(\vep y+x,y')-P_{s-r'}(x,y')] \\
&\quad \biggl(\int_{z''}P_{s-r'}(x,z'')\int_{r''=0}^{r'} \int_{y''}P_{r'-r''}(z'',y'')\M_\vep(\d y'',\d r'')\d z''\biggr)   \M_\vep(\d y',\d r')
\d y\d x \d s.
\end{align*}
Here, the final equality employs the version of the stochastic Fubini theorem found in \cite[Theorem~2.6, pp.296--297]{Walsh}. The required $L^2$ condition for that theorem in \cite{Walsh} is satisfied because, for any fixed $0 < r' < s$ and $x \in \R^2$, $P_{s-r'}(x, z'')\, d z''$ over $z'' \in \R^2$ defines a finite measure, and
\begin{align*}
&\int_{z''}\int_{r''=0}^{r'}\int_{y''}\int_{y'''} P_{s-r'}(x,z'')P_{r'-r''}(z'',y''')P_{r'-r''}(z'',y'')\langle M_\vep,M_\vep\ra (\d y''',\d y'',\d r'')\d z''\\
&\quad =\int_{z''}\int_{r''=0}^{r'}\int_{y''}\int_{\wt{y}''} P_{s-r'}(x,z'')P_{r'-r''}(z'',\vep\wt{y}''+y'')P_{r'-r''}(z'',y'')\phi(\wt{y}'')\d \wt{y}''\d y'' \d r''\d z''<\infty,
\end{align*}
where the last inequality follows upon applying Lemma~\ref{lem:SFTint}.

To proceed, we use the mild form \eqref{def:mild:vep} for a further representation of $\tilde{I}^{2,4}_{\vep}({g},T)$:
\begin{align*}
\tilde{I}^{2,4}_{\vep}({g},T)&=
\int_{s=0}^T \int_{x}\int_{y}
\int_{r'=0}^s\int_{y'}{g} [P_{s-r'}(\vep y+x,y')-P_{s-r'}(x,y')] \\
&\quad \left(\int_{z''}P_{s-r'}(x,z'')[X_\vep(z'',r')-P_{r'}X_0(z'')]
\d z''\right)   \M_\vep(\d y',\d r')
\d y\d x \d s\\
&=
\int_{s=0}^T \int_{x}\int_{y}
\int_{r'=0}^s\int_{y'}{g} [P_{s-r'}(\vep y+x,y')-P_{s-r'}(x,y')] \\
&\quad\times \int_{z''}P_{s-r'}(x,z'')X_\vep(z'',r')
\d z''   \M_\vep(\d y',\d r')
\d y\d x \d s\\
&\quad-\int_{s=0}^T \int_{x}\int_{y}
\int_{r'=0}^s\int_{y'}{g} [P_{s-r'}(\vep y+x,y')-P_{s-r'}(x,y')] \\
&\quad\times P_{s}X_0(x)   \M_\vep(\d y',\d r')
\d y\d x \d s,
\end{align*}
where the last equality holds since $\int_{z''}P_{s-r'}(x,z'')P_{r'}X_0(z'') \d z''=P_sX_0(x)$ for $0< r'< s$. Since the last term equals $\tilde{I}^{2,2}_{\vep}({g},T)$, we get 
\begin{align*}
\tilde{I}^{2,2}_{\vep}({g},T)+\tilde{I}^{2,4}_{\vep}({g},T)
&=\int_{s=0}^T \int_{x}\int_{y}
\int_{r'=0}^s\int_{y'}{g} [P_{s-r'}(\vep y+x,y')-P_{s-r'}(x,y')] \\
&\quad\times \left(\int_{z''}P_{s-r'}(x,z'')X_\vep(z'',r')
\d z''\right)   \M_\vep(\d y',\d r')
\d y\d x \d s\\
&=
\int_{r'=0}^T\int_{y'}\int_{s=r'}^T \int_{x}\int_{y}{g} [P_{s-r'}(\vep y+x,y')-P_{s-r'}(x,y')] \\
&\quad \times\int_{z''}P_{s-r'}(x,z'')X_\vep(z'',r')
\d z'' \d y\d x \d s  \M_\vep(\d y',\d r').
\end{align*}
Here, the final equality is justified by the stochastic Fubini theorem, as in the preceding paragraph, now also invoking Proposition~\ref{prop:XM} (2$\cc$). 
 In particular, we interchange the order of integration so that 
\[
M_\vep(\d y',\d r')g(y,x,s,T)^{\pm}\d x\d y\d s = g(y,x,s,T)^{\pm} \d x \d y \d s M_\vep(\d y',\d r'),
\]
where $g^\pm$ are the positive and negative parts of $g$. 
Equation~\eqref{eq:I2-1} then follows by applying the Fubini theorem to the integrand of $M_\vep(\d y',\d r')$, using $\d z'' \d y \d x \d s = \d y \d x \d s \d z''$. 

The proof of \eqref{eq:I2-2} follows a similar argument. To rewrite $\tilde{I}^{2,5}_{\vep}({g},T)$, we again use the Chapman-Kolmogorov equation to decompose $P_{s-r'}(\vep y + x, y')$ and $P_{s-r'}(x, y')$, followed by applications of the stochastic Fubini theorem. This produces the first equality below:
\begin{align*}
\tilde{I}^{2,5}_{\vep}({g},T)
&= \int_{s=0}^T \int_{x}\int_{y} \int_{r''=0}^s \int_{y''}  {g}\biggl(
\int_{z'}P_{s-r''}(\vep y+x,z') \int_{r'=0}^{r''} \int_{y'} P_{r''-r'}(z',y')\M_\vep(\d y',\d r')\d z'\\
&\quad-\int_{z'}P_{s-r''}(x,z') \int_{r'=0}^{r''} \int_{y'} P_{r''-r'}(z',y')\M_\vep(\d y',\d r')\d z'\biggr)\\
 &\quad \times P_{s-r''}(x,y'') \M_\vep(\d y'',\d r'')\d y \d x\d s\\
 &= \int_{s=0}^T \int_{x}\int_{y} \int_{r''=0}^s \int_{y''} {g}  \biggl(
\int_{z'}P_{s-r''}(\vep y+x,z')[X_\vep(z',r'')-P_{r''}X_0(z')]\d z'\\
&\quad-\int_{z'}P_{s-r''}(x,z') [X_\vep(z',r'')-P_{r''}X_0(z')]\d z'\biggr) P_{s-r''}(x,y'') \M_\vep(\d y'',\d r'')\d y \d x\d s,
\end{align*}
where the last equality follows from the mild form \eqref{def:mild:vep} of $X_\vep$. To simplify the right-hand side, we apply the Chapman--Kolmogorov equation to get $\int_{z'}P_{s-r''}(\wt{y},z')P_{r''}X_0(z')\d z'=P_sX_0(\wt{y})$ for $\wt{y}\in\{\vep y+x,x\}$ and $0< r''< s$.
Hence, after rearrangement, the last equality gives 
 \begin{align*}
\tilde{I}^{2,5}_{\vep}({g},T) 
   &= \int_{s=0}^T \int_{x}\int_{y} \int_{r''=0}^s \int_{y''} {g} \biggl(
\int_{z'}[P_{s-r''}(\vep y+x,z')-P_{s-r''}(x,z')]X_\vep(z',r'')\d z'\biggr)\\
 &\quad \times P_{s-r''}(x,y'') \M_\vep(\d y'',\d r'')\d y \d x\d s\\
    &\quad- \int_{s=0}^T \int_{x}\int_{y} \int_{r''=0}^s \int_{y''}{g}[P_{s}X_0(\vep y+x)-P_{s}X_0(x)] P_{s-r''}(x,y'') \M_\vep(\d y'',\d r'')\d y \d x\d s.
\end{align*}
Since the last integral equals $\tilde{I}^{2,3}_{\vep}({g},T)$, we get
\begin{align*}
&\tilde{I}^{2,3}_{\vep}({g},T)+\tilde{I}^{2,5}_{\vep}({g},T)\\
&\quad = \int_{s=0}^T \int_{x}\int_{y} \int_{r''=0}^s \int_{y''} {g} \left(
\int_{z'}[P_{s-r''}(\vep y+x,z')-P_{s-r''}(x,z')]X_\vep(z',r'')\d z'\right)\\
 &\quad\quad  \times P_{s-r''}(x,y'') \M_\vep(\d y'',\d r'')\d y \d x\d s\\
 &\quad =  \int_{r''=0}^T \int_{y''}\int_{s=r''}^T \int_{x}\int_{y} {g} \left(
\int_{z'}[P_{s-r''}(\vep y+x,z')-P_{s-r''}(x,z')]X_\vep(z',r'')\d z'\right)\\
 &\quad \quad \times P_{s-r''}(x,y'')\d y \d x\d s \M_\vep(\d y'',\d r''),
\end{align*}
where the last equality applies the stochastic Fubini theorem. By applying the Fubini theorem to the last equality, we rearrange the order of integration as $\d z'\d y \d x\d s = \d y \d x\d s\d z'$. This proves \eqref{eq:I2-2}.

Finally, to justify \eqref{dec0:I2}, observe that the mapping $(y,x,s)\mapsto \mathfrak L_\vep^1 f(y,x,s,T)$, as defined in \eqref{def:L1}, belongs to $L^1(\R^2\times \R^2\times [0,T])$. This integrability is ensured by the assumption $f(x,s)\in \C_{pd}^0(\R^2\times [0,T])$ and the inequality below, which follows from Lemma~\ref{lem:SFTint}:
\[
\int_{s=0}^T\int_x\int_y\left(\int_{t=s}^T\int_{\wt{x}}\frac{1}{1+\lvert\wt{x}\rvert^{10}}P_{t-s}(\vep y+x,\wt{x}) P_{t-s}(\wt{x},x)\d \wt{x} \d t\phi(y)\right)\d y\d x\d s<\infty.
\]
The proof is complete.
\end{proof}

The following lemma further develops the functionals $g \mapsto I^{2,k}_\vep(g,T)$ introduced in Lemma~\ref{lem:I12}. It also involves the remainder terms $R_{\ref{key}}(y, x, s, T)$ from Lemma~\ref{lem:I1}, which depend on $f$ via \eqref{key}.

\begin{lem}\label{lem:I2ap}
\noindent {\rm (1$\cc$)} Let $0<T<\infty$. Suppose that $G=\{{g}_\vep(y,x,s,T)\}_{\vep\in (0,\ovl]}$ is a family of functions such that 
$x\mapsto {g}_\vep(y,x,s,T)\in \C^1(\R^2)$ and,
for some constant $M\in (0,\infty)$, the following bound holds for all $\vep\in (0,\ovl]$:
\begin{align*}
&\max\{\lvert {g}_\vep(y,x,s,T)\rvert,\lvert\nabla_{x}{g}_\vep(y,x,s,T)\rvert \}\leq \frac{C({G},m)\1_{\{\lvert y \rvert\leq M\}}}{1+\lvert x\rvert ^m}\\
&\quad \times \left[\lvert \log (T-s)\rvert+1+\lvert\log \lvert y \rvert\rvert+\vep_{\ref{limT:glue}}
\left(\frac{4(T-s)\vep^{-2}}{\lvert y \rvert^2}\right)\right],\;y,x\in \R^2,\;0< s< T,\;m\in \Bbb N.
\end{align*}
Then Lemma~\ref{lem:I12} applies to the functions $g_\vep$, and we have the following estimates:
\begin{align}
\begin{split}\label{I2ap-1}
&\lim_{\vep\to 0}\lvert I^{2,1}_{\vep}({g}_\vep,T)\rvert=0,\quad\quad \;\;  \lim_{\vep\to 0}\E[\lvert I^{2,2}_{\vep}({g}_\vep,T)\rvert^2]=0,\\
&\lim_{\vep\to 0}\E[\lvert I^{2,3}_{\vep}({g}_\vep,T)\rvert^2]=0,\quad 
\limsup_{\vep\to 0}\E[\lvert I^{2,4}_{\vep}({g}_\vep,T)\rvert]<\infty.
\end{split}
\end{align}

\noindent {\rm (2$\cc$)} For all $0<T<\infty$ and $f\in \C^1_{pd}(\R^2\times [0,T])$, the following estimate holds:
\begin{align}
&\lim_{\vep\to 0}\E\biggl[\int_{r'=0}^T \int_{y'}\int_{y''}\int_{y}\biggl(\lv^{-1}\lvert R_{\ref{key}}(y''-y,y',r',T)\rvert+\lv^{-1}\lvert R_{\ref{key}}(y'',y',r',T)\rvert\biggr)\phi(y)\d y\notag\\
&\quad \times \Lambda_{\vep}  X_\vep(\vep y''+y',r') X_\vep(y',r') \phi(y'')\d y''\d y'\d r'\biggr]=0.\label{I2ap-2}
\end{align}
\end{lem}

\begin{proof}
{(1$\cc$)} First, observe that Lemma~\ref{lem:I12} applies to each $g_\vep$, since the assumed bound on $g_\vep$ together with \eqref{EMremainder:bdd} ensures $(y, x, s) \mapsto g_\vep(y, x, s, T)\in L^1(\mathbb{R}^2 \times \mathbb{R}^2 \times [0, T])$ for all $\vep \in (0, \ovl]$. In the remainder of the proof of (1$\cc$), for $j = 1, 2, 3, 4$, we establish the necessary property of $I^{2, j}_{\vep}(g_\vep, T)$, as stated in \eqref{I2ap-1}, in Step~$j$ below. Henceforth, we abbreviate $ {g}_\vep(y, x, s, T)$ as $g_\vep$ whenever the meaning is clear from context.
 \medskip 

\noindent {\bf Step 1.} We utilize the heat kernel estimate for suitably perturbed initial conditions; see \eqref{gauss:vep} with $\delta_\vep = \vep$. Consequently,
\begin{align*}
\lvert I^{2,1}_{\vep}({g}_\vep,T)\rvert&\leq \int_{s=0}^{\vep\wedge T}\int_{x}\int_{y}\lvert {g}_\vep\rvert
[P_sX_0(\vep y+x)+P_sX_0(x)] P_sX_0(x)\d y\d x\d s\\
&\quad +C(M)\vep^{1/2}\int_{s=0}^T \int_{x}\int_{y}\lvert {g}_\vep\rvert
P_{2s}X_0(x)\d y\d x\d s.
\end{align*}
Observe that, by invoking the boundedness of $X_0(\cdot)$, the assumed bound on $\lvert{g}_\vep\rvert$, and \eqref{EMremainder:bdd}, the right-hand side approaches zero as $\vep\to 0$. In particular, this convergence holds as we apply the following limits:
\begin{align*}
&\lim_{\vep\to 0}\int_{s=0}^{\vep\wedge T}\lvert\log [4(T-s)\vep^{-2}]\rvert\1_{\{4(T-s)\vep^{-2}<1\}}\d s=0,\\
&\lim_{\vep\to 0} \vep^{1/2}\int_{s=0}^T\lvert\log [4(T-s)\vep^{-2}]\rvert\1_{\{4(T-s)\vep^{-2}<1\}}\d s=0.
\end{align*}
We have proved the property of $I^{2,1}_{\vep}({g}_\vep,T)$ required in \eqref{I2ap-1}.  \medskip

\noindent {\bf Step~2.} To show that $\E[\lvert I^{2,2}_{\vep}({g}_\vep,T)\rvert^2]$ vanishes, as required in \eqref{I2ap-1}, we again apply the heat kernel estimate in \eqref{gauss:vep} with $\delta_\vep = \vep$. To carry this out, we first apply It\^{o}'s isometry and \eqref{eq:covarM-2} and write
\begin{align}\label{dec:I22}
\E[\lvert I^{2,2}_{\vep}({g}_\vep,T)\rvert^2]=I^{\1,\1}_{\vep,\ref{I22:error}}=\sum\left\{I^{\kappa,\wt{\kappa}}_{\vep,\ref{I22:error}};\kappa,\wt{\kappa}\in \{\1_{[0,\vep)},\1_{[\vep,\infty)}\}\right\},
\end{align}
where $I^{\kappa,\wt{\kappa}}_{\vep,\ref{I22:error}}$ is bilinear in $(\kappa,\wt{\kappa})$, defined by
\begin{align}
I^{\kappa,\wt{\kappa}}_{\vep,\ref{I22:error}}&\,\defeq\,\E\biggl[\int_{r'=0}^T\int_{y'}\int_{y''} \notag\\
&\quad \times \biggl(\int_{z''} \int_{s=r'}^T \kappa(s-r')\int_{x}\int_{y}{g}_\vep(y,x,s,T) [P_{s-r'}(\vep y+x,\vep y''+y')-P_{s-r'}(x,\vep y''+y')]\notag\\
&\quad \times  P_{s-r'}(x,z'') \d y\d x \d s X_\vep(z'',r')
\d z''\biggr)\notag\\
&\quad \times  \biggl(\int_{\wt{z}''}\int_{\wt{s}=r'}^T\wt{\kappa}(\wt{s}-r') \int_{\wt{x}}\int_{\wt{y}}{g}_\vep(\wt{y},\wt{x},\wt{s},T) [P_{\wt{s}-r'}(\vep \wt{y}+\wt{x},y')-P_{\wt{s}-r'}(\wt{x},y')]\notag\\
&\quad \times  P_{\wt{s}-r'}(\wt{x},\wt{z}'') \d \wt{y}\d \wt{x} \d \wt{s} X_\vep(\wt{z}'',r')
\d \wt{z}''\biggr)\notag\\
&\quad \times \lv X_\vep(\vep y''+ y',r')X_\vep( y',r')\phi(y'')\d y''\d y'\d r'
\biggr].\label{I22:error}
\end{align}
Observe that by substituting $\hat{y} = \vep y'' + y'$ for $y'$ and utilizing the even parity of $\phi$ from Proposition~\ref{prop:XM} (1$\cc$), we have:
\begin{align}\label{I22:abelian}
I^{\1_{[\vep,\infty)},\1_{[0,\vep)}}_{\vep,\ref{I22:error}}=I^{\1_{[0,\vep)},\1_{[\vep,\infty)}}_{\vep,\ref{I22:error}}.
\end{align}
(Alternatively, this identity follows directly from the symmetry of quadratic covariations.) Hence, to obtain the limit of $\E[\lvert I^{2,2}_{\vep}({g}_\vep,T)\rvert^2]$ required in \eqref{I2ap-1}, it suffices to verify the following limits for $(\kappa,\wt{\kappa})=(\1_{[\vep,\infty)},\1_{[\vep,\infty)})$, $(\1_{[\vep,\infty)},\1_{[0,\vep)})$, and $(\1_{[0,\vep)},\1_{[0,\vep)})$:
\begin{align}\label{I22:abeliankappa}
\lim_{\vep\to 0}I^{\kappa,\wt{\kappa}}_{\vep,\ref{I22:error}}=0.
\end{align}
The verifications proceed through the next three sub-steps, with the last sub-step giving a conclusion. \medskip 

\noindent {\bf Step 2-1.} We verify \eqref{I22:abeliankappa} in the case of $(\kappa,\wt{\kappa})=(\1_{[\vep,\infty)},\1_{[\vep,\infty)})$. 
The assumed bound on $g_\vep$ ensures that ${g}_\vep(y, x, s, T) = 0$ whenever $\lvert y \rvert > M$. Hence, by  \eqref{gauss:vep} with $\delta_\vep= \vep$,
\begin{align}
& \lvert I^{\1_{[\vep,\infty)},\1_{[\vep,\infty)}}_{\vep,\ref{I22:error}}\rvert\notag\\
&\quad \leq \int_{r'=0}^T \int_{y'}\int_{y''}\int_{z''}\int_{\wt{z}''}\int_{s=r'}^T \int_{x}\int_{y}
\int_{\wt{s}=r'}^T \int_{\wt{x}}\int_{\wt{y}}
\lvert{g}_\vep\rvert(y,x,s,T)\lvert{g}_\vep\rvert(\wt{y},\wt{x},\wt{s},T)\notag\\
&\quad\quad \times C(M)\vep P_{2(s-r')}(x,\vep y''+y')P_{2(\wt{s}-r')}(\wt{x},y')\notag\\
&\quad \quad \times P_{s-r'}(x,z'')P_{\wt{s}-r'}(\wt{x},\wt{z}'') \lv\E[X_\vep(z'',r')X_\vep(\wt{z}'',r')X_\vep(\vep y''+y',r')X_\vep(y',r')]\phi(y'')\notag\\
&\quad\quad  \times \d \wt{y}\d \wt{x}\d \wt{s} \d y\d x\d s\d \wt{z}'' \d z''\d y''\d y'\d r'\notag\\
&\quad \leq \int_{r'=0}^T \int_{y'}\int_{y''}\int_{s=r'}^T \int_{x}\int_{y}
\int_{\wt{s}=r'}^T \int_{\wt{x}}\int_{\wt{y}}\lvert {g}_\vep\rvert(y,x,s,T)\lvert{g}_\vep\rvert(\wt{y},\wt{x},\wt{s},T)\notag\\
&\quad\quad \times  C(M) \vep  P_{2(s-r')}(x,\vep y''+y')P_{2(\wt{s}-r')}(\wt{x},y')\notag\\
&\quad\quad  \times  \frac{ C(\lambda,\phi,\lVert X_0\rVert_\infty,T)}{\sqrt{\min\{s,\wt{s}\}-r'}}\left(\int_{\ol{z}}\phi(\ol{z})\lvert \log \lvert y''-\ol{z}\rvert \rvert \d\ol{z}+1\right)\phi(y'')\notag\\
&\quad \quad \times  \d \wt{y}\d \wt{x}\d \wt{s} \d y\d x\d s\d y''\d y'\d r'.\label{I21:11}
\end{align}
Here, \eqref{I21:11} can be justified by using Theorem~\ref{thm:mombdd} (2$\cc$) with $p=1/2$ after integrating out $\wt{z}''$ and $z''$ since we can write
\begin{align*}
&\int_{z''}\int_{\wt{z}''}P_{s-r'}(x,z'')P_{\wt{s}-r'}(\wt{x},\wt{z}'') \lv\E[X_\vep(z'',r')X_\vep(\wt{z}'',r')X_\vep(\vep y''+y',r')X_\vep(y',r')]\d \wt{z}''\d z''\\
&\quad = \lv\E^{(B^1,B^2)}_{(x,\wt{x})}[X_\vep(B^1_{s-r'},r')X_\vep(B^2_{\wt{s}-r'},r')X_\vep(\vep y''+y',r')X_\vep(y',r')],
\end{align*}
where $B^1$ and $B^2$ are independent two-dimensional standard Brownian motions, also independent of $X_\vep$.  

To use \eqref{I21:11}, we consider the effect of integrating out $\wt{y}$ and $y$ therein. After integrating out $y$, we use
 \eqref{EMremainder:bdd} and the assumption on ${g}_\vep$ to get the following bound:
\begin{align}
&\int_{y}\lvert {g}_\vep\rvert(y,x,s,T)\d y\notag\\
&\quad \leq \frac{C(G,m,M)}{1+\lvert x\rvert ^m}\left[\lvert\log (T-s)\rvert+1+\lvert\log [4(T-s)\vep^{-2}]\rvert\1_{\{4(T-s)\vep^{-2}<1\}}\right].\label{I21:psiaux}
\end{align}
An analogous inequality holds with $(y,x,s)$ replaced by $(\wt{y},\wt{x},\wt{s})$ as we integrate out $\wt{y}$.
Hence, by \eqref{I21:11} with $m\geq 10$,
\begin{align}
\lvert I^{\1_{[\vep,\infty)},\1_{[\vep,\infty)}}_{\vep,\ref{I22:error}}\rvert &\leq C(\lambda,\phi,\lVert X_0\rVert_\infty,T,G,m,M)\notag\\
&\quad \times \int_{r'=0}^T \int_{y'}\int_{y''}\int_{s=r'}^T \int_{x}
\int_{\wt{s}=r'}^T \int_{\wt{x}}
\frac{1}{(1+\lvert x\rvert ^m)(1+\lvert\wt{x}\rvert^m)}
\notag\\
&\quad \times \left[\lvert\log (T-s)\rvert+1+\lvert\log [4(T-s)\vep^{-2}]\rvert\1_{\{4(T-s)\vep^{-2}<1\}}\right]\notag\\
&\quad \times \left[\lvert\log (T-\wt{s})\rvert+1+\lvert\log [4(T-\wt{s})\vep^{-2}]\rvert\1_{\{4(T-\wt{s})\vep^{-2}<1\}}\right]\notag\\
&\quad\times \vep P_{2(s-r')}(x,\vep y''+y')P_{2(\wt{s}-r')}(\wt{x},y')\notag\\
&\quad \times  \frac{ 1}{\sqrt{\min\{s,\wt{s}\}-r'}}\left(\int_{\ol{z}}\phi(\ol{z})\lvert \log \lvert y''-\ol{z}\rvert \rvert \d \ol{z}+1\right)\phi(y'')\notag\\
&\quad \times \d \wt{x}\d \wt{s}\d x\d s \d y''\d y'\d r'.\label{I21:psiaux-1}
\end{align}
To bound the right-hand side, observe that the integral with respect to $\wt{x},x,y'',y'$ can be handled as follows, where the first inequality uses \eqref{ineq:heat2-2}:
\begin{align}
&\int_{y'}\int_{y''}\int_{x}\int_{\wt{x}}\frac{ \vep P_{2(s-r')}(x,\vep y''+y')P_{2(\wt{s}-r')}(\wt{x},y')}{(1+\lvert x\rvert ^m)(1+\lvert\wt{x}\rvert^m)}\notag\\
&\quad\times\left(\int_{\ol{z}}\phi(\ol{z})\lvert \log \lvert y''-\ol{z}\rvert \rvert\d\ol{z}+1\right)\phi(y'')\d \wt{x}\d x\d y''\d y'\notag\\
&\quad \leq C(T,m)  \vep \int_{y'}\int_{y''}\frac{\phi(y'')}{(1+\lvert \vep y''+y'\rvert^m)(1+\lvert y' \rvert^m)}\left(\int_{\ol{z}}\phi(\ol{z})\lvert \log \lvert y''-\ol{z}\rvert \rvert\d \ol{z}+1\right)\d y''\d y'\notag\\
&\quad \leq  C(\phi,T,m)\vep.\label{I21:psiaux1}
\end{align}
Note that the last inequality can be seen by using \eqref{ratiot}, which yields 
\[
(1+\lvert \vep y''+y'\rvert^m)^{-1}\leq C(m) (1+\lvert\vep y''\rvert^m)(1+\lvert y' \rvert^m)^{-1}.
\]
Therefore, the bound in \eqref{I21:psiaux-1} implies the following one:
\begin{align}
\lvert I^{\1_{[\vep,\infty)},\1_{[\vep,\infty)}}_{\vep,\ref{I22:error}}\rvert &\leq  C(\lambda,\phi,\lVert X_0\rVert_\infty,T,G,m,M) \vep 
\int_{r'=0}^T \int_{s=r'}^T\int_{\wt{s}=r'}^T\frac{1}{\sqrt{\min\{s,\wt{s}\}-r'}}\notag\\
&\quad \times \left[\lvert\log (T-s)\rvert+1+\lvert\log [4(T-s)\vep^{-2}]\rvert\1_{\{4(T-s)\vep^{-2}<1\}}\right]\notag\\
&\quad \times \left[\lvert\log (T-\wt{s})\rvert+1+\lvert\log [4(T-\wt{s})\vep^{-2}]\rvert\1_{\{4(T-\wt{s})\vep^{-2}<1\}}\right]
\d \wt{s}\d s\d r'.
\end{align}
We deduce from the last inequality that
\begin{align}\label{I2psi:1}
\lim_{\vep\to 0}I^{\1_{[\vep,\infty)},\1_{[\vep,\infty)}}_{\vep,\ref{I22:error}}=0,
\end{align}
which is \eqref{I22:abeliankappa} in the case of $(\kappa,\wt{\kappa})=(\1_{[\vep,\infty)},\1_{[\vep,\infty)})$.\medskip 

\noindent {\bf Step~2-2.} We verify \eqref{I22:abeliankappa} in the case of $(\kappa,\wt{\kappa})=(\1_{[\vep,\infty)},\1_{[0,\vep)})$. By \eqref{gauss:vep} with the choice of $\delta_\vep= \vep$, we have the following bound:
\begin{align}
&\lvert I^{\1_{[\vep,\infty)},\1_{(0,\vep]}}_{\vep,\ref{I22:error}}\rvert \notag\\
&\quad \leq \int_{r'=0}^T \int_{y'}\int_{y''}\int_{z''}\int_{\wt{z}''}\int_{s=r'}^T \int_{x}\int_{y}
\int_{\wt{s}=r'}^{T\wedge (r'+\vep)} \int_{\wt{x}}\int_{\wt{y}} \lvert {g}_\vep\rvert(y,x,s,T)\lvert {g}_\vep\rvert(\wt{y},\wt{x},\wt{s},T)\notag\\
&\quad\quad \times C(M)\vep^{1/2} P_{2(s-r')}(x,\vep y''+y')[P_{\wt{s}-r'}(\vep \wt{y}+\wt{x},y')+P_{\wt{s}-r'}(\wt{x},y')]\notag\\
&\quad\quad  \times P_{s-r'}(x,z'')P_{\wt{s}-r'}(\wt{x},\wt{z}'') \lv\E[X_\vep(z'',r')X_\vep(\wt{z}'',r')X_\vep(\vep y''+y',r')X_\vep(y',r')]\phi(y'')\notag\\
&\quad\quad  \times  \d \wt{y}\d \wt{x}\d \wt{s} \d y\d x\d s\d \wt{z}'' \d z''\d y''\d y'\d r'\notag\\
&\quad \leq \int_{r'=0}^T \int_{y'}\int_{y''}\int_{s=r'}^T \int_{x}\int_{y}
\int_{\wt{s}=r'}^{T\wedge (r'+\vep)} \int_{\wt{x}}\int_{\wt{y}} \lvert {g}_\vep\rvert(y,x,s,T)\lvert {g}_\vep\rvert(\wt{y},\wt{x},\wt{s},T)\notag\\
&\quad\quad \times C(M)\vep^{1/2} P_{2(s-r')}(x,\vep y''+y')[P_{\wt{s}-r'}(\vep \wt{y}+\wt{x},y')+P_{\wt{s}-r'}(\wt{x},y')]
\notag\\
&\quad\quad  \times  \frac{ C(\lambda,\phi,\lVert X_0\rVert_\infty,T)}{\sqrt{\min\{s,\wt{s}\}-r'}}\left(\int_{\ol{z}}\phi(\ol{z})\lvert \log \lvert y''-\ol{z}\rvert \rvert\d \ol{z}+1\right)\phi(y'')\notag\\
&\quad\quad  \times \d \wt{y}\d \wt{x}\d \wt{s} \d y\d x\d s\d y''\d y'\d r',  \label{I22:Step2}
\end{align}
where the last inequality follows by integrating out $\wt{z}''$ and $z''$ and using Theorem~\ref{thm:mombdd} (2$\cc$) with $p=1/2$.

To estimate the right-hand side of \eqref{I22:Step2}, we first note that the iterated integral over all spatial variables satisfies the first inequality below, which is justified by using the assumed bound for ${g}_\vep$ with $m \geq 10$:
\begin{align*}
& \int_{y'}\int_{y''} \int_{x}\int_{y} \int_{\wt{x}}\int_{\wt{y}} \lvert {g}_\vep\rvert(y,x,s,T)\lvert {g}_\vep\rvert(\wt{y},\wt{x},\wt{s},T)\\
&\quad\times C(M)\vep^{1/2} P_{2(s-r')}(x,\vep y''+y')[P_{\wt{s}-r'}(\vep \wt{y}+\wt{x},y')+P_{\wt{s}-r'}(\wt{x},y')]
\\
&\quad \times\left(\int_{\ol{z}}\phi(\ol{z})\lvert \log \lvert y''-\ol{z}\rvert \rvert\d \ol{z}+1\right)\phi(y'')\d \wt{y}\d \wt{x} \d y\d x \d y''\d y'\\
&\quad \leq C(G,m,M)\int_{y'}\int_{y''}\int_{x}\int_{y}\int_{\wt{x}}\int_{\wt{y}} \\
&\quad\quad  \times \frac{\1_{\{\lvert y \rvert\leq M\}}}{1+\lvert x\rvert ^m} \left[\lvert\log (T-s)\rvert+1+\lvert\log \lvert y \rvert\rvert+\vep_{\ref{limT:glue}}
\left(\frac{4(T-s)\vep^{-2}}{\lvert y \rvert^2}\right)\right]\\
&\quad \quad \times\frac{\1_{\{\lvert\wt{y}\rvert\leq M\}}}{1+\lvert\wt{x}\rvert^m} \left[\lvert\log (T-\wt{s})\rvert+1+\lvert \log \lvert \wt{y}\rvert \rvert+\vep_{\ref{limT:glue}}
\left(\frac{4(T-\wt{s})\vep^{-2}}{\lvert\wt{y}\rvert^2}\right)\right]\\
&\quad\quad  \times \vep^{1/2} P_{2(s-r')}(x,\vep y''+y')[P_{\wt{s}-r'}(\vep \wt{y}+\wt{x},y')+P_{\wt{s}-r'}(\wt{x},y')]
\\
&\quad\quad   \times\left(\int_{\ol{z}}\phi(\ol{z})\lvert \log \lvert y''-\ol{z}\rvert \rvert\d \ol{z}+1\right)\phi(y'') \d \wt{y}\d \wt{x} \d y\d x \d y''\d y'\\
&\quad \leq C(T,G,m,M)\int_{y'}\int_{y''}\int_{y}\int_{\wt{y}} \\
&\quad \quad  \times \1_{\{\lvert y \rvert\leq M\}} \left[\lvert\log (T-s)\rvert+1+\lvert\log \lvert y \rvert\rvert+\vep_{\ref{limT:glue}}
\left(\frac{4(T-s)\vep^{-2}}{\lvert y \rvert^2}\right)\right]\\
&\quad \quad \times \1_{\{\lvert\wt{y}\rvert\leq M\}}\left[\lvert\log (T-\wt{s})\rvert+1+\lvert\log \lvert\wt{y}\rvert\rvert+\vep_{\ref{limT:glue}}
\left(\frac{4(T-\wt{s})\vep^{-2}}{\lvert\wt{y}\rvert^2}\right)\right]\\
&\quad  \quad \times \vep^{1/2} \left(\frac{1+\lvert\vep y''\rvert^m}{1+\lvert y' \rvert^m}\frac{1+\lvert\vep \wt{y}\rvert^m}{1+\lvert y' \rvert^m}\right)
\\
&\quad \quad \times \left(\int_{\ol{z}}\phi(\ol{z})\lvert\log \lvert y''-\ol{z}\rvert\rvert\d \ol{z}+1\right)\phi(y'') \d \wt{y} \d y \d y''\d y',
\end{align*}
where the last inequality follows by integrating out $\wt{x}$ and $x$ and using \eqref{ratiot} and \eqref{ineq:heat2-2}.  
Also, by integrating out $y'$, $y''$, $y$ and $\wt{y}$ in the same order and using \eqref{EMremainder:bdd}, the last inequality implies
\begin{align*}
& \int_{y'}\int_{y''} \int_{x}\int_{y} \int_{\wt{x}}\int_{\wt{y}} \lvert {g}_\vep\rvert(y,x,s,T)\lvert {g}_\vep\rvert(\wt{y},\wt{x},\wt{s},T)\\
&\quad\times C(M)\vep^{1/2} P_{2(s-r')}(x,\vep y''+y')[P_{\wt{s}-r'}(\vep \wt{y}+\wt{x},y')+P_{\wt{s}-r'}(\wt{x},y')]
\\
&\quad \times \left(\int_{\ol{z}}\phi(\ol{z})\lvert \log \lvert y''-\ol{z}\rvert \rvert \d \ol{z}+1\right)\phi(y'')\d \wt{y}\d \wt{x} \d y\d x \d y''\d y'\\
&\quad \leq  C(\phi,T,G,m,M)\vep^{1/2}\left[\lvert\log (T-s)\rvert+1+\lvert\log [4(T-s)\vep^{-2}]\rvert\1_{\{4(T-s)\vep^{-2}<1\}}\right]\\
&\quad\quad  \times \left[\lvert\log (T-\wt{s})\rvert+1+\lvert\log [4(T-\wt{s})\vep^{-2}]\rvert\1_{\{4(T-\wt{s})\vep^{-2}<1\}}\right].
\end{align*}
Applying the right-hand side to the right-hand side of \eqref{I22:Step2}, we get
\begin{align*}
\lvert I^{\1_{[\vep,\infty)},\1_{(0,\vep]}}_{\ref{I22:error}}\rvert &\leq  C(\lambda,\phi,\lVert X_0\rVert_\infty,T,G,m,M) \vep^{1/2} \int_{r'=0}^T\int_{s=r'}^T
\int_{\wt{s}=r'}^{T\wedge (r'+\vep)}\frac{1}{\sqrt{\min\{s,\wt{s}\}-r'}}\\
&\quad \times  \left[\lvert\log (T-s)\rvert+1+\lvert\log [4(T-s)\vep^{-2}]\rvert\1_{\{4(T-s)\vep^{-2}<1\}}\right]\\
&\quad \times \left[\lvert\log (T-\wt{s})\rvert+1+\lvert\log [4(T-\wt{s})\vep^{-2}]\rvert\1_{\{4(T-\wt{s})\vep^{-2}<1\}}\right]\d \wt{s}\d s\d r'.
\end{align*}
It follows that 
\begin{align}
\lim_{\vep\to 0}I^{\1_{[\vep,\infty)},\1_{(0,\vep]}}_{\ref{I22:error}}=0,\label{I2psi:2}
\end{align}
which is \eqref{I22:abeliankappa} in the case of $(\kappa,\wt{\kappa})=(\1_{[\vep,\infty)},\1_{[0,\vep)})$.\medskip  

\noindent {\bf Step~2-3.} We now address \eqref{I22:abeliankappa} in the final case, where $(\kappa, \wt{\kappa}) = (\1_{[0,\vep)}, \1_{[0,\vep)})$. In this situation, we have only the following straightforward bound available: 
\begin{align}
& \lvert I^{\1_{(0,\vep]},\1_{(0,\vep]}}_{\ref{I22:error}}\rvert \notag\\
&\quad \leq \int_{r'=0}^T \int_{y'}\int_{y''}\int_{z''}\int_{\wt{z}''}\int_{s=r'}^{T\wedge(r'+\vep)} \int_{x}\int_{y}
\int_{\wt{s}=r'}^{T\wedge (r'+\vep)} \int_{\wt{x}}\int_{\wt{y}} \lvert {g}_\vep\rvert(y,x,s,T)\lvert {g}_\vep\rvert(\wt{y},\wt{x},\wt{s},T)\notag\\
&\quad \quad\times      [P_{s-r'}(\vep y+x,\vep y''+y')+P_{s-r'}(x,\vep y''+y')] [P_{\wt{s}-r'}(\vep \wt{y}+\wt{x},y')+P_{\wt{s}-r'}(\wt{x},y')]\notag\\
&\quad\quad  \times P_{s-r'}(x,z'')P_{\wt{s}-r'}(\wt{x},\wt{z}'') \lv\E[X_\vep(z'',r')X_\vep(\wt{z}'',r')X_\vep(\vep y''+y',r')X_\vep(y',r')]\phi(y'')\notag\\
&\quad \quad \times \d \wt{y}\d \wt{x}\d \wt{s} \d y\d x\d s\d \wt{z}'' \d z''\d y''\d y'\d r'\notag\\
&\quad \leq \int_{r'=0}^T \int_{y'}\int_{y''}\int_{s=r'}^{T\wedge (r'+\vep)} \int_{x}\int_{y}
\int_{\wt{s}=r'}^{T\wedge (r'+\vep)} \int_{\wt{x}}\int_{\wt{y}}\lvert {g}_\vep\rvert(y,x,s,T)\lvert {g}_\vep\rvert(\wt{y},\wt{x},\wt{s},T)\notag\\
&\quad\quad \times      [P_{s-r'}(\vep y+x,\vep y''+y')+P_{s-r'}(x,\vep y''+y')] [P_{\wt{s}-r'}(\vep \wt{y}+\wt{x},y')+P_{\wt{s}-r'}(\wt{x},y')]\notag\\
&\quad\quad  \times  \frac{ C(\lambda,\phi,\lVert X_0\rVert_\infty,T)}{\sqrt{\min\{s,\wt{s}\}-r'}}\left(\int_{\ol{z}}\phi(\ol{z})\lvert \log \lvert y''-\ol{z}\rvert \rvert \d \ol{z}+1\right)\phi(y'')\notag\\
&\quad\quad  \times  \d \wt{y}\d \wt{x}\d \wt{s} \d y\d x\d s\d y''\d y'\d r'.\label{I22:Step2-2}
\end{align}
Here, the last inequality is justified by integrating out $\wt{z}''$ and $z''$ and then using
Theorem~\ref{thm:mombdd} (2$\cc$) with $p=1/2$. 

To integrate out the remaining spatial variables on the right-hand side of \eqref{I22:Step2-2}, we isolate the relevant multiple integral and use the assumed bound on ${g}_\vep$ for $m \geq 10$. This yields
\begin{align*}
& \int_{y'}\int_{y''} \int_{x}\int_{y} \int_{\wt{x}}\int_{\wt{y}} \lvert {g}_\vep\rvert(y,x,s,T)\lvert {g}_\vep\rvert(\wt{y},\wt{x},\wt{s},T)\\
&\quad\times  [P_{s-r'}(\vep y+x,\vep y''+y')+P_{s-r'}(x,\vep y''+y')] [P_{\wt{s}-r'}(\vep \wt{y}+\wt{x},y')+P_{\wt{s}-r'}(\wt{x},y')]\notag\\
&\quad \times \left(\int_{\ol{z}}\phi(\ol{z})\lvert \log \lvert y''-\ol{z}\rvert \rvert \d \ol{z}+1\right)\phi(y'')\d \wt{y}\d \wt{x} \d y\d x \d y''\d y'\\
&\quad \leq C(G,m)\int_{y'}\int_{y''}\int_{x}\int_{y}\int_{\wt{x}}\int_{\wt{y}} \\
&\quad \quad \times \frac{\1_{\{\lvert y \rvert\leq M\}}}{1+\lvert x\rvert ^m} \left[\lvert\log (T-s)\rvert+1+\lvert\log \lvert y \rvert\rvert+\vep_{\ref{limT:glue}}
\left(\frac{4(T-s)\vep^{-2}}{\lvert y \rvert^2}\right)\right]\\
&\quad\quad  \times\frac{\1_{\{\lvert\wt{y}\rvert\leq M\}}}{1+\lvert\wt{x}\rvert^m} \left[\lvert\log (T-\wt{s})\rvert+1+\lvert \log \lvert \wt{y}\rvert \rvert +\vep_{\ref{limT:glue}}
\left(\frac{4(T-\wt{s})\vep^{-2}}{\lvert\wt{y}\rvert^2}\right)\right]\\
&\quad\quad \times  [P_{s-r'}(\vep y+x,\vep y''+y')+P_{s-r'}(x,\vep y''+y')] [P_{\wt{s}-r'}(\vep \wt{y}+\wt{x},y')+P_{\wt{s}-r'}(\wt{x},y')]\notag\\
&\quad\quad  \times \left(\int_{\ol{z}}\phi(\ol{z})\lvert \log \lvert y''-\ol{z}\rvert \rvert \d \ol{z}+1\right)\phi(y'') \d \wt{y}\d \wt{x} \d y\d x \d y''\d y'\\
&\quad \leq C(T,G,m)\int_{y'}\int_{y''}\int_{y}\int_{\wt{y}} \\
&\quad\quad \times \1_{\{\lvert y \rvert\leq M\}} \left[\lvert\log (T-s)\rvert+1+\lvert\log \lvert y \rvert\rvert+\vep_{\ref{limT:glue}}
\left(\frac{4(T-s)\vep^{-2}}{\lvert y \rvert^2}\right)\right]\\
&\quad\quad  \times \1_{\{\lvert\wt{y}\rvert\leq M\}}\left[\lvert\log (T-\wt{s})\rvert+1+\lvert \log \lvert \wt{y}\rvert \rvert +\vep_{\ref{limT:glue}}
\left(\frac{4(T-\wt{s})\vep^{-2}}{\lvert\wt{y}\rvert^2}\right)\right]\\
&\quad \quad \times  \left(\frac{1+\lvert\vep y\rvert^m+\lvert \vep y''\rvert^m}{1+\lvert y' \rvert^m}\frac{1+\lvert\vep \wt{y}\rvert^m}{1+\lvert y' \rvert^m}\right)
 \left(\int_{\ol{z}}\phi(\ol{z})\lvert \log \lvert y''-\ol{z}\rvert \rvert \d \ol{z}+1\right)\phi(y'') \d \wt{y} \d y \d y''\d y'.
\end{align*}
Here, the last inequality is obtained by integrating out $x$ and $\wt{x}$, and applying the following: the inequalities in \eqref{ratiot} and \eqref{ineq:heat2-2}, along with $(a+b)^m \leq C(m)(a^m+b^m)$ for all $a,b \geq 0$. Furthermore, by integrating out $y''$, $y'$, $y$ and $\wt{y}$ in the same order and using \eqref{EMremainder:bdd}, the last inequality gives
\begin{align*}
&\int_{y'}\int_{y''} \int_{x}\int_{y} \int_{\wt{x}}\int_{\wt{y}} \lvert {g}_\vep\rvert(y,x,s,T)\lvert {g}_\vep\rvert(\wt{y},\wt{x},\wt{s},T)\\
&\quad\times  [P_{s-r'}(\vep y+x,\vep y''+y')+P_{s-r'}(x,\vep y''+y')] [P_{\wt{s}-r'}(\vep \wt{y}+\wt{x},y')+P_{\wt{s}-r'}(\wt{x},y')]\notag\\
&\quad \times \left(\int_{\ol{z}}\phi(\ol{z})\lvert \log \lvert y''-\ol{z}\rvert \rvert \d \ol{z}+1\right)\phi(y'')\d \wt{y}\d \wt{x} \d y\d x \d y''\d y'\\
&\quad \leq  C(\phi,T,G,m,M)\left[\lvert\log (T-s)\rvert+1+\lvert\log [4(T-s)\vep^{-2}]\rvert\1_{\{4(T-s)\vep^{-2}<1\}}\right]\\
&\quad\quad  \times \left[\lvert\log (T-\wt{s})\rvert+1+\lvert\log [4(T-\wt{s})\vep^{-2}]\rvert\1_{\{4(T-\wt{s})\vep^{-2}<1\}}\right].
\end{align*}

Finally, applying the right-hand side of the foregoing inequality to the right-hand side of \eqref{I22:Step2-2}, we get
\begin{align*}
\lvert I^{\1_{(0,\vep]},\1_{(0,\vep]}}_{\vep,\ref{I22:error}}\rvert &\leq  C(\lambda,\phi,\lVert X_0\rVert_\infty,T,G,m,M) \int_{r'=0}^{T}\int_{s=r'}^{T\wedge (r'+\vep)}
\int_{\wt{s}=r'}^{T\wedge (r'+\vep)}\frac{1}{\sqrt{\min\{s,\wt{s}\}-r'}}\\
&\quad  \times \left[\lvert\log (T-s)\rvert+1+\lvert\log [4(T-s)\vep^{-2}]\rvert\1_{\{4(T-s)\vep^{-2}<1\}}\right]\\
&\quad \times \left[\lvert\log (T-\wt{s})\rvert+1+\lvert\log [4(T-\wt{s})\vep^{-2}]\rvert\1_{\{4(T-\wt{s})\vep^{-2}<1\}}\right]\d \wt{s}\d s\d r'.
\end{align*}
It follows that 
\begin{align}
\lim_{\vep\to 0}I^{\1_{(0,\vep]},\1_{(0,\vep]}}_{\vep,\ref{I22:error}}=0,\label{I2psi:3}
\end{align}
which is \eqref{I22:abeliankappa} in the case of $(\kappa,\wt{\kappa})=(\1_{[0,\vep)},\1_{[0,\vep)})$.\medskip 

\noindent {\bf Step~2-4.} We now complete Step~2. 
Recall \eqref{dec:I22} and \eqref{I22:abelian}. Hence, combining \eqref{I2psi:1}, \eqref{I2psi:2} and \eqref{I2psi:3} yields the following limit:
\[
\lim_{\vep\to 0}\E[\lvert I^{2,2}_{\vep}({g}_\vep,T)\rvert^2]=0,
\]
as required in \eqref{I2ap-1}.\medskip 

\noindent {\bf Step~3.} Observe that $I^{2,2}_{\vep}({g}_\vep,T)$ and $I^{2,3}_{\vep}({g}_\vep,T)$ are quite similar; the differences become apparent by reducing $''$ in the variables to $'$, and then swapping $z'$ with $y'$ in the terminal conditions of the heat kernels used. Hence, the argument here largely resembles that in Step~2, aside from a few minor changes.
Again, we begin by applying It\^{o}'s isometry and \eqref{eq:covarM-2} to compute $\E[\lvert I^{2,3}_{\vep}({g}_\vep,T)\rvert^2]$. This yields the next equality:
\begin{align}\label{dec:I23}
\E[\lvert I^{2,3}_{\vep}({g}_\vep,T)\rvert^2]=I^{\1,\1}_{\vep,\ref{I23:error}}=\sum\left\{I^{\kappa,\wt{\kappa}}_{\vep,\ref{I23:error}};\kappa,\wt{\kappa}\in \{\1_{[0,\vep)},\1_{[\vep,\infty)}\}\right\},
\end{align}
where $I^{\kappa,\wt{\kappa}}_{\vep,\ref{I23:error}}$ is blinear in $(\kappa,\wt{\kappa})$, defined by
\begin{align}
I^{\kappa,\wt{\kappa}}_{\vep,\ref{I23:error}}&\,\defeq\,\E\bigg[\int_{r''=0}^T\int_{y''}\int_{y'''} \notag\\
&\quad \times \bigg(\int_{z'} \int_{s=r''}^T \kappa(s-r'')\int_{x}\int_{y}{g}_\vep(y,x,s,T) [P_{s-r''}(\vep y+x,z')-P_{s-r''}(x,z')]\notag\\
&\quad \times  P_{s-r''}(x,\vep y'''+y'') \d y\d x \d s X_\vep(z',r'')
\d z'\bigg)\notag\\
&\quad \times  \bigg(\int_{\wt{z}'}\int_{\wt{s}=r''}^T\wt{\kappa}(\wt{s}-r'') \int_{\wt{x}}\int_{\wt{y}}{g}_\vep(\wt{y},\wt{x},\wt{s},T) [P_{\wt{s}-r''}(\vep \wt{y}+\wt{x},\wt{z}')-P_{\wt{s}-r''}(\wt{x},\wt{z}')]\notag\\
&\quad \times  P_{\wt{s}-r''}(\wt{x},y'') \d \wt{y}\d \wt{x} \d \wt{s}X_\vep(\wt{z}',r'')
\d \wt{z}'\bigg)\notag\\
&\quad \times \lv X_\vep(\vep y'''+ y'',r'')X_\vep( y'',r'')\phi(y''')\d y'''\d y''\d r''
\bigg].\label{I23:error}
\end{align}
In this case, we still have $I^{\kappa,\wt{\kappa}}_{\vep,\ref{I23:error}}=I^{\wt{\kappa},\kappa}_{\vep,\ref{I23:error}}$, which can be seen again by using the substitution $\hat{y}=\vep y'''+y''$ for $y''$ and the even parity of $\phi$ from Proposition~\ref{prop:XM} (1$\cc$). 
Hence, to prove the required limit, given by $\E[\lvert I^{2,3}_{\vep}({g}_\vep,T)\rvert^2]\to 0$,
it suffices to show that $I^{\kappa,\wt{\kappa}}_{\vep,\ref{I23:error}}\to 0$ as $\vep\to 0$, for $(\kappa,\wt{\kappa})=(\1_{[\vep,\infty)},\1_{[\vep,\infty)})$, $(\1_{[\vep,\infty)},\1_{[0,\vep)})$ and $(\1_{[0,\vep)},\1_{[0,\vep)})$. The verifications are shown in the next two sub-steps, with the last sub-step concluding Step~3. 
\medskip 

\noindent {\bf Step 3-1.} In the case of $(\kappa,\wt{\kappa})=(\1_{[\vep,\infty)},\1_{[\vep,\infty)})$, we have the following bound, which is obtained by applying  \eqref{gauss:vep} with the choice of $\delta_\vep=\vep$:
\begin{align}
& \lvert I^{\1_{[\vep,\infty)},\1_{[\vep,\infty)}}_{\vep,\ref{I23:error}}\rvert\notag\\
&\quad \leq \int_{r''=0}^T \int_{y''}\int_{y'''}\int_{z'}\int_{\wt{z}'}\int_{s=r''}^T \int_{x}\int_{y}
\int_{\wt{s}=r''}^T \int_{\wt{x}}\int_{\wt{y}}
\lvert {g}_\vep\rvert(y,x,s,T)\lvert {g}_\vep\rvert(\wt{y},\wt{x},\wt{s},T)\notag\\
&\quad\quad \times C(M) \vep P_{2(s-r'')}(x,z')P_{2(\wt{s}-r'')}(\wt{x},\wt{z}')\notag\\
&\quad\quad  \times P_{s-r''}(x,\vep y'''+y'')P_{\wt{s}-r''}(\wt{x},y'') \lv\E[X_\vep(z',r'')X_\vep(\wt{z}',r'')X_\vep(\vep y'''+y'',r'')X_\vep(y'',r'')]\phi(y''')\notag\\
&\quad\quad  \times \d \wt{y}\d \wt{x}\d \wt{s} \d y\d x\d s \d \wt{z}' \d z'\d y'''\d y''\d r''\notag\\
&\quad \leq \int_{r''=0}^T \int_{y''}\int_{y'''}\int_{s=r''}^T \int_{x}\int_{y}
\int_{\wt{s}=r''}^T \int_{\wt{x}}\int_{\wt{y}}\lvert {g}_\vep\rvert(y,x,s,T)\lvert {g}_\vep\rvert(\wt{y},\wt{x},\wt{s},T)\notag\\
&\quad\quad \times C(M) \vep P_{s-r''}(x,\vep y'''+y'')P_{\wt{s}-r''}(\wt{x},y'')\notag\\
&\quad\quad  \times  \frac{ C(\lambda,\phi,\lVert X_0\rVert_\infty,T)}{\sqrt{\min\{s,\wt{s}\}-r''}}\left(\int_{\ol{z}}\phi(\ol{z})\lvert \log \lvert y'''-\ol{z}\rvert \rvert \d \ol{z}+1\right)\phi(y''')\notag\\
&\quad\quad  \times \d \wt{y}\d \wt{x}\d \wt{s} \d y\d x\d s \d y'''\d y''\d r'',\label{I21:22}
\end{align}
where the last inequality uses Theorem~\ref{thm:mombdd} (2$\cc$) by integrating out $z'$ and $\wt{z}'$. Observe that by replacing $y'',y''',r''$ with $y',y'',r'$, respectively, and using the inequality $P_{a}(z)\leq 2 P_{2a}(z)$, the right-hand side of \eqref{I21:22} can be bounded by the right-hand side of \eqref{I21:11}, possibly larger multiplicative constants $C(M)$ and $C(\lambda,\phi,\lVert X_0\rVert_\infty,T)$. Hence, the calculation below \eqref{I21:11} applies, and we have a direct analogue of \eqref{I2psi:1}, given by the following limit: 
\begin{align}\label{I23psi:1}
\lim_{\vep\to 0}I^{\1_{[\vep,\infty)},\1_{[\vep,\infty)}}_{\vep,\ref{I23:error}}=0.
\end{align}

\noindent {\bf Step~3-2.} The required limits for $(\kappa,\wt{\kappa}) = (\1_{[\vep, \infty)}, \1_{(0, \vep]})$ and $(\kappa, \wt{\kappa}) = (\1_{(0, \vep]}, \1_{(0, \vep]})$ follow analogously by applying the arguments from Steps~2-2 and 2-3; therefore, we present them together in this sub-step.

We first show the bounds that link to those arguments. 
In the case of $(\kappa,\wt{\kappa})=(\1_{[\vep,\infty)},\1_{(0,\vep]})$
applying \eqref{gauss:vep} with the choice of $\delta_\vep=\vep$ leads to 
\begin{align}
& \lvert I^{\1_{[\vep,\infty)},\1_{(0,\vep]}}_{\vep,\ref{I23:error}}\rvert\notag\\
&\quad \leq \int_{r''=0}^T \int_{y''}\int_{y'''}\int_{z'}\int_{\wt{z}'}\int_{s=r''}^T \int_{x}\int_{y}
\int_{\wt{s}=r''}^{T\wedge (r''+\vep)} \int_{\wt{x}}\int_{\wt{y}} \lvert {g}_\vep\rvert(y,x,s,T)\lvert {g}_\vep\rvert(\wt{y},\wt{x},\wt{s},T)\notag\\
&\quad\quad \times C(M) \vep^{1/2} P_{2(s-r'')}(x,z')[P_{\wt{s}-r''}(\vep \wt{y}+\wt{x},\wt{z}')+P_{\wt{s}-r''}(\wt{x},\wt{z}')]\notag\\
&\quad\quad  \times P_{s-r''}(x,\vep y'''+y'')P_{\wt{s}-r''}(\wt{x},y'') \lv\E[X_\vep(z',r'')X_\vep(\wt{z}',r'')X_\vep(\vep y'''+y'',r'')X_\vep(y'',r'')]\phi(y''')\notag\\
&\quad\quad  \times \d \wt{y}\d \wt{x}\d \wt{s}\d y\d x\d s\d  \wt{z}'\d z'\d y'''\d y''\d r''\notag\\
&\quad \leq \int_{r''=0}^T \int_{y''}\int_{y'''}\int_{s=r''}^T \int_{x}\int_{y}
\int_{\wt{s}=r''}^{T\wedge (r''+\vep)} \int_{\wt{x}}\int_{\wt{y}} \lvert {g}_\vep\rvert(y,x,s,T)\lvert {g}_\vep\rvert(\wt{y},\wt{x},\wt{s},T)\notag\\
&\quad\quad \times C(M)\vep^{1/2} P_{s-r''}(x,\vep y'''+y'')P_{\wt{s}-r''}(\wt{x},y'')
\notag\\
&\quad\quad  \times  \frac{ C(\lambda,\phi,\lVert X_0\rVert_\infty,T)}{\sqrt{\min\{s,\wt{s}\}-r''}}\left(\int_{\ol{z}}\phi(\ol{z})\lvert \log \lvert y'''-\ol{z}\rvert \rvert \d \ol{z}+1\right)\phi(y''')\notag\\
&\quad\quad  \times \d \wt{y}\d \wt{x}\d \wt{s} \d y\d x\d s \d y'''\d y''\d r'',\label{I22:Step22}
\end{align}
where the last inequality is followed by integrating out $\wt{z}'$ and $z'$ and
using Theorem~\ref{thm:mombdd} (2$\cc$) with $p=1/2$, as well as the bound $\min\{2A,B\}\geq \min\{A,B\}$ for all $A,B\geq 0$. Additionally, in the case of $(\kappa,\wt{\kappa})=(\1_{(0,\vep]},\1_{(0,\vep]})$, we have
\begin{align}
& \lvert I^{\1_{(0,\vep]},\1_{(0,\vep]}}_{\vep,\ref{I23:error}}\rvert\notag\\
&\quad \leq \int_{r''=0}^T \int_{y''}\int_{y'''}\int_{z'}\int_{\wt{z}'}\int_{s=r''}^{T\wedge(r''+\vep)} \int_{x}\int_{y}
\int_{\wt{s}=r''}^{T\wedge (r''+\vep)} \int_{\wt{x}}\int_{\wt{y}} \lvert {g}_\vep\rvert(y,x,s,T)\lvert {g}_\vep\rvert(\wt{y},\wt{x},\wt{s},T)\notag\\
&\quad\quad \times      [P_{s-r''}(\vep y+x,z')+P_{s-r''}(x,z')] [P_{\wt{s}-r''}(\vep \wt{y}+\wt{x},\wt{z}')+P_{\wt{s}-r''}(\wt{x},\wt{z}')]\notag\\
&\quad\quad  \times P_{s-r''}(x,\vep y'''+y'')P_{\wt{s}-r''}(\wt{x},y'') \lv\E[X_\vep(z',r'')X_\vep(\wt{z}',r'')X_\vep(\vep y'''+y'',r'')X_\vep(y'',r'')]\phi(y''')\notag\\
&\quad\quad  \times \d \wt{y}\d \wt{x}\d \wt{s}\d y\d x\d s\d  \wt{z}'\d z'\d y'''\d y''\d r''\notag\\
&\quad \leq \int_{r''=0}^T \int_{y''}\int_{y'''}\int_{s=r''}^{T\wedge (r''+\vep)}  \int_{x}\int_{y}
\int_{\wt{s}=r''}^{T\wedge (r''+\vep)} \int_{\wt{x}}\int_{\wt{y}} \lvert {g}_\vep\rvert(y,x,s,T)\lvert {g}_\vep\rvert(\wt{y},\wt{x},\wt{s},T)\notag\\
&\quad\quad \times    P_{s-r''}(x,\vep y'''+y'')P_{\wt{s}-r''}(\wt{x},y'') \notag\\
&\quad\quad  \times  \frac{ C(\lambda,\phi,\lVert X_0\rVert_\infty,T)}{\sqrt{\min\{s,\wt{s}\}-r''}}\left(\int_{\ol{z}}\phi(\ol{z})\lvert \log \lvert y'''-\ol{z}\rvert \rvert \d \ol{z}+1\right)\phi(y''')\notag\\
&\quad\quad  \times \d \wt{y}\d \wt{x}\d \wt{s} \d y\d x\d s \d y'''\d y''\d r'',\label{I22:Step2-22}
\end{align}
where the last inequality is justified by integrating out $\wt{z}'$ and $z'$ and then using Theorem~\ref{thm:mombdd} (2$\cc$) with $p=1/2$. 

We now invoke the arguments from Steps~2-2 and 2-3. 
By enlarging the constants if necessary, the right-hand side of \eqref{I22:Step2} bounds the right-hand side of  \eqref{I22:Step22}. This implies
\begin{align}
\lim_{\vep\to 0}I^{\1_{[\vep,\infty)},\1_{(0,\vep]}}_{\vep,\ref{I23:error}}=0.\label{I23psi:2}
\end{align}
Additionally, as the right-hand side of \eqref{I22:Step2-2} bounds the right-hand side of  \eqref{I22:Step2-22}, we deduce 
\begin{align}
\lim_{\vep\to 0}I^{\1_{(0,\vep]},\1_{(0,\vep]}}_{\vep,\ref{I23:error}}=0.\label{I23psi:3}
\end{align}

\noindent {\bf Step~3-3.} Recall \eqref{dec:I23} and, as explained below \eqref{I23:error}, $I^{\1_{[\vep,\infty)},\1_{(0,\vep]}}_{\vep,\ref{I23:error}}=I^{\1_{(0,\vep]},\1_{[\vep,\infty)}}_{\vep,\ref{I23:error}}$. Hence, by \eqref{I23psi:1}, \eqref{I23psi:2} and \eqref{I23psi:3}, we obtain
\[
\lim_{\vep\to 0}\E[\lvert I^{2,3}_{\vep}({g}_\vep,T)\rvert^2]=0,
\]
as required in \eqref{I2ap-1}. \medskip 

\noindent {\bf Step~4.} To prove that the limit superior of  $\E[\lvert I^{2,4}_{\vep}({g}_\vep,T)\rvert]$ as $\vep \to 0$ is finite, as asserted in \eqref{I2ap-1}, we introduce two decompositions to work with. The first one is similar to those in \eqref{dec:I22} and \eqref{dec:I23}. We write
\begin{align}
I^{2,4}_{\vep}({g}_\vep,T)=I^{2,4,1}_{\vep}({g}_\vep,T)+I^{2,4,2}_{\vep}({g}_\vep,T),\label{dec:I24}
\end{align}
where 
\begin{align*}
I^{2,4,1}_{\vep}({g}_\vep,T)&\,\defeq\,
 \int_{r'=0}^T \int_{y'}\int_{y''} \int_{s=r'}^{ T} \int_{x}\int_{y}\1_{[0,\vep)}(s-r'){g}_\vep (y,x,s,T)\\
 &\quad \times 
 [P_{s-r'}(\vep y+x,\vep y''+y')-P_{s-r'}(x,\vep y''+y')] P_{s-r'}(x,y')\d y \d x\d s \\
 &\quad \times \lv  X_\vep(\vep y''+y',r') X_\vep(y',r') \phi(y'')\d y''\d y'
\d r',\\
I^{2,4,2}_{\vep}({g}_\vep,T)&\,\defeq\,
 \int_{r'=0}^T \int_{y'}\int_{y''} \int_{s=r'}^{ T} \int_{x}\int_{y}\1_{[\vep,\infty)}(s-r'){g}_\vep (y,x,s,T)\\
 &\quad \times
 [P_{s-r'}(\vep y+x,\vep y''+y')-P_{s-r'}(x,\vep y''+y')] P_{s-r'}(x,y')\d y \d x\d s\\
 &\quad \times  \lv  X_\vep(\vep y''+y',r') X_\vep(y',r') \phi(y'')\d y''\d y'
\d r'.
\end{align*}

The second decomposition we need is stated as follows: 
\begin{align}
I^{2,4,1}_{\vep}({g}_\vep,T)=I^{2,4,1,1}_{\vep}({g}_\vep,T)+I^{2,4,1,2}_{\vep}({g}_\vep,T),\label{dec:I241}
\end{align}
where
\begin{align}
I^{2,4,1,1}_{\vep}({g}_\vep,T)&\,\defeq\,
\int_{r'=0}^T \int_{y'}\int_{y''} \int_{s'=0}^{\vep^{-2} (T-r')} \int_{\wt{x}}\int_{y}\1_{[0,\vep^{-1})}(s')\notag\\
&\quad \times{g}_\vep (y,y',r'+\vep^2 s',T)\notag\\
&\quad \times [P_{s'}( y+ \wt{x}, y'')-P_{s'}( \wt{x}, y'')]  P_{s'}( \wt{x})\d y \d \wt{x}\d s'  \notag\\
&\quad \times\lv  X_\vep(\vep y''+y',r') X_\vep(y',r') \phi(y'')\d y''\d y'
\d r',\notag\\
I^{2,4,1,2}_{\vep}({g}_\vep,T)&\,\defeq\,
\int_{r'=0}^T \int_{y'}\int_{y''} \int_{s'=0}^{\vep^{-2} (T-r')} \int_{\wt{x}}\int_{y}\1_{[0,\vep^{-1})}(s')\notag\\
&\quad \times [{g}_\vep (y,y'+\vep\wt{x},r'+\vep^2 s',T)-{g}_\vep (y,y',r'+\vep^2 s',T) ] \notag\\
&\quad \times [P_{s'}( y+ \wt{x}, y'')-P_{s'}( \wt{x}, y'')]P_{s'}( \wt{x})\d y \d \wt{x}\d s' \notag\\
&\quad \times \lv  X_\vep(\vep y''+y',r') X_\vep(y',r') \phi(y'')\d y''\d y'
\d r'.\notag
\end{align}
To see \eqref{dec:I241}, we apply the substitions $s'=\vep^{-2}(s-r')$ and $\wt{x}=\vep^{-1} (x-y')$ for $s$ and $x$ to  the definition of $I^{2,4,1}_{\vep}({g}_\vep,T)$. This gives
\begin{align*}
I^{2,4,1}_{\vep}({g}_\vep,T)&= \int_{r'=0}^T \int_{y'}\int_{y''}\int_{s'=0}^{\vep^{-2} (T-r')} \int_{\wt{x}}\int_{y}\1_{[0,\vep^{-1})}(s'){g}_\vep (y,y'+\vep\wt{x},r'+\vep^2 s',T)\\
&\quad \times [P_{\vep^2s'}(\vep y+\vep \wt{x},\vep y'')-P_{\vep^2s'}(\vep \wt{x},\vep y'')] P_{\vep^2s'}(\vep \wt{x})\d y\vep^{2} \d \wt{x} \vep^2\d s' \\
&\quad \times \lv  X_\vep(\vep y''+y',r') X_\vep(y',r') \phi(y'')\d y''\d y'
\d r'\\
&= \int_{r'=0}^T \int_{y'}\int_{y''} \int_{s'=0}^{\vep^{-2} (T-r')} \int_{\wt{x}}\int_{y}\1_{[0,\vep^{-1})}(s'){g}_\vep (y,y'+\vep\wt{x},r'+\vep^2 s',T) \\
&\quad \times  [P_{s'}( y+ \wt{x}, y'')-P_{s'}( \wt{x}, y'')] P_{s'}( \wt{x})\d y \d \wt{x}\d s' \\
&\quad \times \lv  X_\vep(\vep y''+y',r') X_\vep(y',r') \phi(y'')\d y''\d y'
\d r',
\end{align*}
where the last equality is obtained by making some cancellations of $\vep$ on its left-hand side. The decomposition in \eqref{dec:I241} now follows at once by adding the integrands of $I^{2,4,1,1}_{\vep}({g}_\vep,T)$ and $I^{2,4,1,2}_{\vep}({g}_\vep,T)$.

In the following sub-steps, we prove, in the same order, that the limit superiors of the two terms $\E[\lvert I^{2,4,1,1}_{\vep}({g}_\vep,T)\rvert]$ and $\E[\lvert I^{2,4,1,2}_{\vep}({g}_\vep,T)\rvert]$ as $\vep\to 0$ are both finite, whereas the limit $\E[\lvert I^{2,4,2}_{\vep}({g}_\vep,T)\rvert]$ is zero. A subsequent sub-step concludes Step~4. \medskip 

\noindent {\bf Step 4-1.} To prove the required property of $I^{2,4,1,1}_{\vep}({g}_\vep,T)$, first, we integrate out $\wt{x}$ in its integral representation and use the Chapman--Kolmogorov equation. This gives the following expression:
\begin{align*}
I^{2,4,1,1}_{\vep}({g}_\vep,T)&=\int_{r'=0}^T \int_{y'}\int_{y''}\int_{s'=0}^{\vep^{-2} (T-r')}\int_{y}\1_{[0,\vep^{-1})}(s')\notag\\
&\quad \times{g}_\vep (y,y',r'+\vep^2 s',T)  [P_{2s'}(  y''-y)-P_{2s'}( y'')]\d y\d s' \\
&\quad\times \lv  X_\vep(\vep y''+y',r') X_\vep(y',r') \phi(y'')\d y''\d y'
\d r'.
\end{align*}
Since $\phi(\cdot)$ has a compact support, it follows from Theorem~\ref{thm:mombdd} (3$\cc$) that
\begin{align*}
\E[\lvert I^{2,4,1,1}_{\vep}({g}_\vep,T)\rvert]&\leq C(\lambda,\phi,\lVert X_0\rVert_\infty,T)\int_{r'=0}^T \int_{y'}\int_{y''} \int_{s'=0}^{\vep^{-2} (T-r')}\int_{y}\1_{[0,\vep^{-1})}(s')\notag\\
&\quad \times\lvert{g}_\vep \rvert(y,y',r'+\vep^2 s',T) \lvert P_{2s'}(  y''-y)-P_{2s'}( y'')\rvert\d y\d s' \phi(y'')\d y''\d y'
\d r'.
\end{align*}

To improve the bound in the preceding inequality, observe that
\begin{align*}
\forall\;\lvert y \rvert\leq M,\quad 
&\int_{y''}\lvert P_{2s'}(  y''-y)-P_{2s'}( y'')\rvert\phi(y'')\d y''\\
&\quad \leq \1_{\{s'\leq 1\}}C(\phi)+\1_{\{s'>1\}}\int_{y''}\left|P_{2s'}(  y''-y)-\frac{1}{4\pi s'}+\frac{1}{4\pi s'}-P_{2s'}( y'')\right|\phi(y'')\d y''\\
&\quad \leq \1_{\{s'\leq 1\}}C(\phi)+\1_{\{s'>1\}}\int_{y''}\biggl(\frac{1}{4\pi s'}\cdot \frac{\rvert y''-y\rvert ^2}{4s'}+\frac{1}{4\pi s'}\cdot \frac{\lvert y'' \rvert^2}{4s'}\biggr)\phi(y'')\d y''\\
&\quad \leq  \frac{C(\phi,M)}{(s')^2\vee 1},
\end{align*}
where the second inequality follows since $0\leq 1-\e^{-x}\leq x$ for all $x\geq 0$, and the last inequality follows since $\phi$ has a compact support and we assume $\lvert y \rvert\leq M$. Since 
${g}_\vep(y,x,s,T)=0$ for all $\lvert y \rvert>M$ by the assumed bound, the foregoing bound applies, and we get
\begin{align}
&\E[\lvert I^{2,4,1,1}_{\vep}({g}_\vep,T)\rvert]\notag\\
&\quad \leq C(\lambda,\phi,\lVert X_0\rVert_\infty,T,M)\int_{r'=0}^T \int_{y'} \int_{s'=0}^{\vep^{-2} (T-r')} \int_{y}\frac{\1_{[0,\vep^{-1})}(s')}{(s')^2\vee 1}\lvert {g}_\vep\rvert (y,y',r'+\vep^2 s',T)\d y\d s'\d y'
\d r'\notag\\
\begin{split}
&\quad \leq C(\lambda,\phi,\lVert X_0\rVert_\infty,T,M)\int_{r'=0}^T \int_{y'}\int_{s=0}^{T-r'} \int_{y}\frac{\vep^{-2}\1_{[0,\vep)}(s)}{(\vep^{-2}s)^2\vee 1}\lvert {g}_\vep\rvert (y,y',r'+ s,T)\d y\d s\\
&\quad\quad \times \d y'
\d r',\label{I2411-original}
\end{split}
\end{align}
where the last inequality uses the substitution $s=\vep^2 s'$. 

To proceed, we use the finer details of the assumed bound on ${g}_\vep$, with $m\geq 10$, to obtain a further bound on $\E[\lvert I^{2,4,1,1}_{\vep}({g}_\vep,T)\rvert]$. Hence, by \eqref{I2411-original},
\begin{align}
\E[\lvert I^{2,4,1,1}_{\vep}({g}_\vep,T)\rvert]&\leq C(\lambda,\phi,\lVert X_0\rVert_\infty,T,G,m,M)\int_{r'=0}^T \int_{y'}\int_{s=0}^{T-r'} \int_{y}\frac{\vep^{-2}\1_{[0,\vep)}(s)}{(\vep^{-2}s)^2\vee 1}\times\frac{\1_{\{\lvert y \rvert\leq M\}}}{1+\lvert y' \rvert^m}\notag\\
&\quad \times  \left[\lvert\log (T-r'-s)\rvert+1+\lvert\log \lvert y \rvert\rvert+\vep_{\ref{limT:glue}}
\left(\frac{4(T-r'-s)\vep^{-2}}{\lvert y \rvert^2}\right)\right]
\d y\d s \d y'
\d r'\notag\\
&\leq C(\lambda,\phi,\lVert X_0\rVert_\infty,T,G,m,M)\int_{r''=0}^T\int_{s=0}^{r''} \frac{\vep^{-2}\1_{[0,\vep)}(s)}{(\vep^{-2}s)^2\vee 1}\notag\\
&\quad\times \biggl[\lvert\log (r''-s)\rvert+1+\lvert\log [4(r''-s)\vep^{-2}]\rvert \1_{\{4(r''-s)\vep^{-2}<1\}}\biggr]\d s\d r'',\label{I2411-aux}
\end{align}
where the last inequality follows by applying \eqref{EMremainder:bdd} and substituting $r'' = T - r'$, then integrating out $y$ and $y'$. 

Finally, we use \eqref{I2411-aux} by writing $\1_{[0,\vep)}(s)=\1_{[\vep^2,\vep)}(s)+\1_{[0,\vep^2)}(s)$, noting that $\vep^{-2}s\geq 1$ if and only if $s\geq \vep^2$, and interchanging the order of integration from $\d s\d r''$ to $\d r''\d s$. Hence, 
\begin{align*}
\E[\lvert I^{2,4,1,1}_{\vep}({g}_\vep,T)\rvert]&\leq 
 C(\lambda,\phi,\lVert X_0\rVert_\infty,T,G,m,M)\int_{s=0}^{T} \frac{\vep^{-2}\1_{[\vep^2,\vep)}(s)}{(\vep^{-2}s)^2}\\
&\quad\times \int_{r''=s}^T \biggl[\lvert\log (r''-s)\rvert+1+\lvert\log [4(r''-s)\vep^{-2}]\rvert\1_{\{4(r''-s)\vep^{-2}<1\}}\biggr]\d r''\d s\\
&\quad +C(\lambda,\phi,\lVert X_0\rVert_\infty,T,G,m,M)\int_{s=0}^{T} \vep^{-2}\1_{[0,\vep^2)}(s)\\
&\quad\times \int_{r''=s}^T\biggl[\lvert\log (r''-s)\rvert+1+\lvert\log [4(r''-s)\vep^{-2}]\rvert\1_{\{4(r''-s)\vep^{-2}<1\}}\biggr]\d r''\d s.
\end{align*}
To bound the limit superior of the right-hand side, it suffices to observe that
\begin{align}
&\int_{s=0}^{T} \frac{\vep^{-2}\1_{[\vep^2,\vep)}(s)}{(\vep^{-2}s)^2}\d s=\int_{s=0}^{\vep^{-2}T}\frac{\1_{[1,\vep^{-1})(s)}}{s^2}\d s,\label{I2411-aux-conclusion1}\\
&\int_{s=0}^T\vep^{-2}\1_{[0,\vep^2)}(s)\d s=\int_{s=0}^{\vep^{-2}T}\1_{[0,1)}(s)\d s,\label{I2411-aux-conclusion2}\\
&\int_{r''=s}^T\lvert\log [4(r''-s)\vep^{-2}]\rvert\1_{\{4(r''-s)\vep^{-2}<1\}}\d r''\notag\\
&\quad \leq \int_{r''=s}^{s+\frac{\vep^2}{4}}\lvert\log [4(r''-s)\vep^{-2}]\rvert\d r''
=\int_{r''=0}^{\frac{\vep^2}{4}}\lvert\log (4r''\vep^{-2})\rvert\d r''=\frac{\vep^2}{4}\int_{r'''=0}^{1}\lvert\log r'''\rvert\d r''',\label{I2411-aux-conclusion3}
\end{align}
where the last equality uses the substitution $r'''=4r''\vep^{-2}$. Hence, we have
\begin{align}\label{I2411:result}
\lim\sup_{\vep\to 0}\E[\lvert I^{2,4,1,1}_{\vep}({g}_\vep,T)\rvert]<\infty,
\end{align}
which achieves the goal of this sub-step. \medskip 

\noindent {\bf Step~4-2.} Recall that $I^{2,4,1,2}_{\vep}({g}_\vep,T)$ is defined below \eqref{dec:I241}. To prove that the $L^1(\P)$-norm of $I^{2,4,1,2}_{\vep}({g}_\vep,T)$ remains bounded as $\vep \to 0$, we start by showing an estimate for the $\d \wt{x}$-integral in $I^{2,4,1,2}_\vep({g}_\vep,T)$. Recall that, by assumption, $x\mapsto g_\vep(y,x,s,T)\in \C^1(\R^2)$. Hence, by the mean value theorem and the assumed bound on $\nabla{g}_\vep$, we get
\begin{align}
&\lvert{g}_\vep (y,y'+\vep\wt{x},r'+\vep^2 s',T)-{g}_\vep (y,y',r'+\vep^2 s',T)\rvert\notag\\
&\quad \less \lvert \vep \wt{x}\rvert
\sup_{x:\lvert x-y'\rvert \leq \lvert \vep\wt{x}\rvert }\lvert\nabla_x{g}_\vep (y,x,r'+\vep^2 s',T)\rvert\notag\\
&\quad \less C(G,m)\lvert\vep \wt{x}\rvert\sup_{x:\lvert x-y'\rvert\leq  \lvert\vep\wt{x}\rvert}\frac{\1_{\{\lvert y \rvert\leq M\}}}{1+\lvert x\rvert ^m}\notag\\
&\quad\quad \times  \left[\lvert\log (T-r'-\vep^2 s')\rvert+1+\lvert\log \lvert y \rvert\rvert+\vep_{\ref{limT:glue}}
\left(\frac{4(T-r'-\vep^2 s')\vep^{-2}}{\lvert y \rvert^2}\right)\right]\notag\\
&\quad \less C(G,m)(\lvert\vep \wt{x}\rvert+\lvert\vep \wt{x}\rvert^{m+1})\frac{\1_{\{\lvert y \rvert\leq M\}}}{1+\lvert y' \rvert^m}\notag\\
&\quad \quad \times  \left[\lvert\log (T-r'-\vep^2 s')\rvert+1+\lvert\log \lvert y \rvert\rvert+\vep_{\ref{limT:glue}}
\left(\frac{4(T-r'-\vep^2 s')\vep^{-2}}{\lvert y \rvert^2}\right)\right],\label{bdd:I2412-1} 
\end{align}
where the last inequality holds since, whenever $\lvert x-y'\rvert\leq \lvert\vep\wt{x}\rvert$, by \eqref{ratiot}, 
\[
1/(1+\lvert x\rvert ^m)\leq C(m)(1+\lvert x-y'\rvert^m)/(1+\lvert y' \rvert^m)\leq C(m)(1+\lvert\vep \wt{x}\rvert^m)/(1+\lvert y' \rvert^m).
\] 
Additionally, for any integer $k\geq 1$ and any $\hat{y}\in \R^2$,
\begin{align}
&\int_{\wt{x}}\lvert\wt{x}\rvert^k P_{s'}(\hat{y},\wt{x})P_{s'}(\wt{x})\d \wt{x}= (\sqrt{2s'})^{k}\int_{\wt{x}}\left|\frac{\wt{x}}{\sqrt{2s'}}\right|^k P_{s'}(\hat{y},\wt{x})P_{s'}(\wt{x})\d \wt{x}\notag\\
&\quad \leq C(k)(\sqrt{2s'})^k\int_{\wt{x}} P_{s'}(\hat{y},\wt{x})P_{2s'}(\wt{x})\d \wt{x}
=C(k) (\sqrt{s'})^k P_{3s'}(\hat{y}),\label{bdd:I2412-2} 
\end{align}
where the inequality follows from the fact that $x^k e^{-x} \leq C(k) \e^{-x/2}$ for all $x \geq 0$, and the final equality is a consequence of the Chapman--Kolmogorov equation.
By \eqref{bdd:I2412-1} and \eqref{bdd:I2412-2}, the $\d \wt{x}$-integral in $I^{2,4,1,2}_\vep({g}_\vep,T)$ satisfies the following bound:
\begin{align}
&\Bigg|\int_{\wt{x}} [{g}_\vep (y,y'+\vep\wt{x},r'+\vep^2 s',T)-{g}_\vep (y,y',r'+\vep^2 s',T) ] \notag\\
&\quad \times [P_{s'}( y+ \wt{x}, y'')-P_{s'}( \wt{x}, y'')] P_{s'}( \wt{x}) \d \wt{x}\Bigg|\notag\\
&\quad \leq C({G},m)(\lvert\vep \sqrt{s'}\rvert+\lvert\vep \sqrt{s'}\rvert^{m+1})\frac{\1_{\{\lvert y \rvert\leq M\}}}{1+\lvert y' \rvert^m}\notag\\
&\quad \quad \times  \left[\lvert\log (T-r'-\vep^2 s')\rvert+1+\lvert\log \lvert y \rvert\rvert+\vep_{\ref{limT:glue}}
\left(\frac{4(T-r'-\vep^2 s')\vep^{-2}}{\lvert y \rvert^2}\right)\right]\notag\\
&\quad \quad\times [P_{3s'}(y''-y)+P_{3s'}(y'')].\notag
\end{align}

We proceed to the first bound of $\E[\lvert I^{2,4,1,2}_{\vep}({g}_\vep,T)\rvert]$. 
By the last inequality and Theorem~\ref{thm:mombdd} (3$\cc$), we have
\begin{align}
\E[\lvert I^{2,4,1,2}_{\vep}({g}_\vep,T)\rvert]
&\leq C(\lambda,\phi,\lVert X_0\rVert_\infty,T,{G},m,M)\notag\\
&\quad \times\int_{r'=0}^T \int_{y'}\int_{y''}\int_{s'=0}^{\vep^{-2} (T-r')}\int_{y}\1_{[0,\vep^{-1})}(s') (\lvert \vep \sqrt{s'}\rvert+\lvert\vep \sqrt{s'}\rvert^{m+1})\frac{\1_{\{\lvert y \rvert\leq M\}}}{1+\lvert y' \rvert^m}\notag\\
&\quad \times  \left[\lvert\log (T-r'-\vep^2 s')\rvert+1+\lvert\log \lvert y \rvert\rvert+\vep_{\ref{limT:glue}}
\left(\frac{4(T-r'-\vep^2 s')\vep^{-2}}{\lvert y \rvert^2}\right)\right]\notag\\
&\quad \times [P_{3s'}(y''-y)+P_{3s'}(y'')] \d y\d s' \phi(y'')\d y''\d y'
\d r'.\label{I2412-original1}
\end{align}
Here, the integral with respect to $\d y\d y''\d y'$ satisfies the following bound:
\begin{align*}
&\int_{y'}\int_{y''}\int_{y} \frac{\1_{\{\lvert y \rvert\leq M\}}}{1+\lvert y' \rvert^m} \left[\lvert\log (T-r'-\vep^2 s')\rvert+1+\lvert\log \lvert y \rvert\rvert+\vep_{\ref{limT:glue}}
\left(\frac{4(T-r'-\vep^2 s')\vep^{-2}}{\lvert y \rvert^2}\right)\right]\\
&\quad \times [P_{3s'}(y''-y)+P_{3s'}(y'')]\phi(y'')\d y\d y''\d y'\\
&\quad \leq  \frac{C(\phi,m,M)}{(s')\vee 1}\biggl[\lvert\log (T-r'-\vep^2s')\rvert+1+\lvert\log [4(T-r'-\vep^2s')\vep^{-2}]\rvert\1_{\{4(T-r'-\vep^2s')\vep^{-2}<1\}}\biggr],
\end{align*}
where the last inequality follows by integrating out $y''$, $y'$ and $y$ in the same order and using  the following bound when integrating out $y''$:
\[
\int_{y}P_{3s'}(\hat{y},y)\phi(y)\d y\leq\frac{C(\phi)}{(s')\vee 1},\quad\forall\;\hat{y}\in \R^2,
\]
as well as \eqref{EMremainder:bdd} when integrating out $y$. By using the bound just obtained for the integral with respect to $\d y\d y''\d y'$, as well as
 the substitutions $r''=T-r'$ and $s=\vep^2 s'$, \eqref{I2412-original1} implies the first inequality below:
\begin{align}
&\E[\lvert I^{2,4,1,2}_{\vep}({g}_\vep,T)\rvert]\notag\\
&\quad \leq C(\lambda,\phi,\lVert X_0\rVert_\infty,T,{G},m,M)
\int_{r''=0}^T \int_{s=0}^{ r''}\1_{[0,\vep)}(s)\frac{(\sqrt{s}+\sqrt{s}^{m+1})}{(\vep^{-2}s)\vee 1}\notag\\
&\quad\quad \times   \biggl[\lvert\log (r''-s)\rvert+1+\lvert\log [4(r''-s)\vep^{-2}]\rvert\1_{\{4(r''-s)\vep^{-2}<1\}}\biggr]\vep^{-2}\d s 
\d r''\label{I2412-original2}\\
&\quad \leq C(\lambda,\phi,\lVert X_0\rVert_\infty,T,{G},m,M)
\int_{r''=0}^T \int_{s=0}^{ r''}\frac{\vep^{-2}\sqrt{s}\1_{[0,\vep)}(s)}{(\vep^{-2}s)\vee 1} \notag\\
&\quad\quad \times  \biggl[\lvert\log (r''-s)\rvert+1+\lvert\log [4(r''-s)\vep^{-2}]\rvert\1_{\{4(r''-s)\vep^{-2}<1\}}\biggr]\d s 
\d r''.\label{I2412-original}
\end{align}
Observe that the right-hand side is similar to the right-hand side of \eqref{I2411-aux}, though not the same. Hence, we rewrite the foregoing inequality as follows: 
\begin{align}
&\E[\lvert I^{2,4,1,2}_{\vep}({g}_\vep,T)\rvert]\notag\\
&\quad \leq  C(\lambda,\phi,\lVert X_0\rVert_\infty,T,G,m,M)\int_{s=0}^{T} \frac{\vep^{-2}\sqrt{s}\1_{[\vep^2,\vep)}(s)}{\vep^{-2}s}\notag\\
&\quad\times \int_{r''=s}^T \biggl[\lvert\log (r''-s)\rvert+1+\lvert\log [4(r''-s)\vep^{-2}]\rvert\1_{\{4(r''-s)\vep^{-2}<1\}}\biggr]\d r''\d s\notag\\
&\quad +C(\lambda,\phi,\lVert X_0\rVert_\infty,T,G,m,M)\int_{s=0}^{T} \vep^{-2}\sqrt{s}\1_{[0,\vep^2)}(s)\notag\\
&\quad\times \int_{r''=s}^T\biggl[\lvert\log (r''-s)\rvert+1+\lvert\log [4(r''-s)\vep^{-2}]\rvert\1_{\{4(r''-s)\vep^{-2}<1\}}\biggr]\d r''\d s.\notag
\end{align}
To show that the right-hand side remains bounded as $\vep \to 0$, we apply \eqref{I2411-aux-conclusion3} once more, but this time we replace \eqref{I2411-aux-conclusion1} and \eqref{I2411-aux-conclusion2} with the following estimates, using the substitution $s' = \vep^{-2} s$:
\begin{align*}
\int_{s=0}^{T} \frac{\vep^{-2}\sqrt{s}\1_{[\vep^2,\vep)}(s)}{\vep^{-2}s}\d s&=\int_{s'=0}^{\vep^{-2}T}\frac{\sqrt{\vep^2s'}\1_{[1,\vep^{-1})}(s')}{s'}\d s'
\leq \vep \cdot 2[\sqrt{\vep^{-2}(T\vee 1)}-1],\\
\int_{s=0}^{T} \vep^{-2}\sqrt{s}\1_{[0,\vep^2)}(s)\d s&=\int_{s'=0}^{\vep^{-2}T}\sqrt{\vep^2s'}\1_{[0,1)}(s')\d s'.
\end{align*}
 Hence, we have the following analogue of \eqref{I2411:result}:
\begin{align}\label{I2412:result}
\limsup_{\vep\to 0}\E[\lvert I^{2,4,1,2}_{\vep}({g}_\vep,T)\rvert]<\infty. 
\end{align}

\noindent {\bf Step~4-3.} Lastly, we show that the $L^1(\P)$-norm of $I^{2,4,2}_{\vep}({g}_\vep,T)$ tends to zero, where $I^{2,4,2}_{\vep}({g}_\vep,T)$ is defined below \eqref{dec:I24}. Since ${g}_\vep(y,x,s,T)=0$ whenever $\lvert y \rvert>M$ by assumption, \eqref{gauss:vep} with $\delta=\vep$ gives the following bound:
\begin{align}
\lvert I^{2,4,2}_{\vep}({g}_\vep,T)\rvert
&\leq \int_{r'=0}^T \int_{y'}\int_{y''}\int_{s=r'}^{ T} \int_{x}\int_{y}
\lvert {g}_\vep\rvert(y,x,s,T)\notag\\ 
&\quad \times  C(M)\vep^{1/2} P_{2(s-r')}(x,\vep y''+y') P_{s-r'}(x,y')\d y \d x\d s\notag\\
 &\quad\times  \lv  X_\vep(\vep y''+y',r') X_\vep(y',r') \phi(y'')\d y''\d y'
\d r'\notag\\
&\leq C(G,m,M,T)\vep^{1/2} \int_{r'=0}^T \int_{y'}\int_{y''} \int_{s=r'}^{ T}\int_{y}\1_{\{\lvert y \rvert\leq M\}}
\frac{P_{4(s-r')}(\vep y'')}{1+\lvert y' \rvert^m}\notag\\
&\quad \times \left[\lvert\log (T-s)\rvert+1+\lvert\log \lvert y \rvert\rvert+\vep_{\ref{limT:glue}}
\left(\frac{4(T-s)\vep^{-2}}{\lvert y \rvert^2}\right)\right]\d y\d s\notag\\
 &\quad\times  \lv  X_\vep(\vep y''+y',r') X_\vep(y',r') \phi(y'')\d y''\d y'
\d r',\label{I242-original1}
\end{align}
where the last inequality follows from applying the assumed bound on ${g}_\vep$ with $m \geq 10$, together with \eqref{ineq:heat2-4} when integrating out $x$. By integrating out $y$ and using \eqref{EMremainder:bdd}, \eqref{I242-original1} implies
\begin{align*}
\lvert I^{2,4,2}_{\vep}({g}_\vep,T)\rvert
&\leq C({G},m,M,T)\vep^{1/2} \int_{r'=0}^T \int_{y'}\int_{y''} \int_{s=r'}^{ T}
\frac{P_{4(s-r')}(\vep y'')}{1+\lvert y' \rvert^m}\\
&\quad \times \biggl[\lvert\log (T-s)\rvert+1+\lvert\log [4(T-s)\vep^{-2}]\rvert\1_{\{4(T-s)\vep^{-2}<1\}} \biggr]\d s \\
 &\quad\times  \lv  X_\vep(\vep y''+y',r') X_\vep(y',r') \phi(y'')\d y''\d y'
\d r'.
\end{align*} 
By Theorem~\ref{thm:mombdd} (3$\cc$), the foregoing bound implies
\begin{align}
&\E[\lvert I^{2,4,2}_{\vep}({g}_\vep,T)\rvert]\notag\\
&\quad \leq C(\lambda,\phi,\lVert X_0\rVert_\infty,T,{G},m,M)\vep^{1/2} \int_{r'=0}^T \int_{y'}\int_{y''} \int_{s=r'}^{ T}
\frac{P_{4(s-r')}(\vep y'')}{1+\lvert y' \rvert^m}\notag\\
&\quad \quad \times \biggl[\lvert\log (T-s)\rvert+1+\lvert\log [4(T-s)\vep^{-2}]\rvert\1_{\{4(T-s)\vep^{-2}<1\}} \biggr]\d s   \phi(y'')\d y''\d y'
\d r'\label{I242-original}\\
&\quad \leq C(\lambda,\phi,\lVert X_0\rVert_\infty,T,{G},m,M)\vep^{1/2} \int_{s=0}^T \int_{y''}\int_{r'=0}^{s}
P_{4(s-r')}(\vep y'')\d r'\notag\\
&\quad\quad  \times \biggl[\lvert\log (T-s)\rvert+1+\lvert\log [4(T-s)\vep^{-2}]\rvert\1_{\{4(T-s)\vep^{-2}<1\}} \biggr]  \phi(y'')\d y''
\d s,\notag
\end{align} 
where the right-hand side follows by integrating out $y'$ and exchanging the order of integration so that $\d s\d r'=\d r'\d s$. By using Lemma~\ref{lem:heat1} (1$\cc$), the right-hand side of the last inequality tends to zero as $\vep\to 0$. This proves
\begin{align}\label{I242:result}
\lim_{\vep\to 0}\E[\lvert I^{2,4,2}_{\vep}({g}_\vep,T)\rvert]=0.
\end{align}

\noindent {\bf Step~4-4.} Recall \eqref{dec:I24} and \eqref{dec:I241}. Hence, by  \eqref{I2411:result}, \eqref{I2412:result}, \eqref{I242:result}, we have proved that 
\[
\limsup_{\vep\to 0}\E[\lvert I^{2,4}_{\vep}({g}_\vep,T)\rvert]<\infty.
\]
All required properties of $I^{2,j}_{\vep}({g}_\vep,T)$ for $j = 1,2,3,4$, as stated in \eqref{I2ap-1}, have now been proved. This completes the proof of (1$\cc$). \medskip 

\noindent {(2$\cc$)} By the definition \eqref{def:ovl} of $\ovl$, $\lv\geq \pi/ \log \vep^{-1}$ for all $\vep\in (0,\ovl]$. In this context, we can bound part of the integrand in \eqref{I2ap-2} as follows by using \eqref{R:bdd}: for $\varsigma\in \{0,1\}$,
\begin{align*}
\lv^{-1}\lvert R_{\ref{key}}(y''-\varsigma y,y',r',T)\rvert&\leq \frac{C(\lambda,T,f,m,M)
}{(\log\vep^{-1})(1+\lvert y' \rvert^m)}\left(\lvert\log (T-r')\rvert+1+\lvert \log \lvert y''-\varsigma y\rvert \rvert \right)\\
&\quad+\frac{C(\lambda,T,f,m,M)
}{(1+\lvert y' \rvert^m)}   \vep_{\ref{limT:glue}}
\left(\frac{4(T-r')\vep^{-2}}{\lvert y''-\varsigma y\rvert^2}\right).
\end{align*}
Additionally, since $\phi$ has a compact support, 
 we can use Theorem~\ref{thm:mombdd} (3$\cc$) to bound the expectation $\E[\Lambda_{\vep}  X_\vep(\vep y''+y',r') X_\vep(y',r')]$, present in \eqref{I2ap-2}, simply by a constant $C(\lambda,\phi,\lVert X_0\rVert_\infty,T)$. Hence,  \eqref{I2ap-2}
follows as soon as we show the following limits for all $m\geq 10$ and $\varsigma\in \{0,1\}$:
\begin{gather*}
\lim_{\vep\to 0}\int_{r'=0}^T\int_{y'}\int_{y''} \int_y
\frac{1}{(1+\lvert y' \rvert^m)}   \left|\vep_{\ref{limT:glue}}
\left(\frac{4(T-r')\vep^{-2}}{\lvert y''-\varsigma y\rvert^2}\right)\right|\phi(y)\phi(y'')\d y\d y''\d y'\d r'=0.
\end{gather*}
To observe these limits, recall once again that $\phi$ has compact support. Therefore, by \eqref{EMremainder:bdd}, the two limits discussed above can be justified directly using the following formulas:
\begin{align}\label{vep:limtool}
\begin{aligned}
\int_{r''=0}^T \min\left\{\frac{1}{4r''\vep^{-2}},1\right\}\d r''&=\frac{\vep^2[\log T-\log (T\wedge \frac{\vep^2}{4})]}{4}+T\wedge \frac{\vep^2}{4},\\
\int_{r''=0}^T\lvert\log (4r''\vep^{-2})\rvert\1_{\{4r''\vep^{-2}<1\}}\d r''&=\vep^2\int_0^{(\vep^{-2}T)\wedge \frac{1}{4}}\lvert\log (4r''')\rvert\d r'''.
\end{aligned}
\end{align}
We have proved the required limit in \eqref{I2ap-2}. The proof of Lemma~\ref{lem:I2ap} is complete.
\end{proof}
\mbox{}

\begin{proof}[End of the proof of Proposition~\ref{prop:I2}]
First, to prove \eqref{id:mu=nu} and \eqref{id:mur=nu}, it suffices to show the following particular case, using test functions $g\in \C^1_c(\R^2\times \R^2)$:
\begin{align}
\int_{s=0}^T\int_x\int_{\wt{x}}g(\wt{x},x)\mu^\sharp_\infty(\d \wt{x}, \d x,\d s)&=\int_{s=0}^T\int_{x} g(x,x)\psi(x)\nu^\sharp_\infty(\d x,\d s),\label{id1:mu}\\
\int_{s=0}^T\int_{x}\int_{\wt{x}}g(\wt{x},x)\mathring{\mu}^\sharp_\infty(\d \wt{x}, \d x,\d s)&=\int_{s=0}^T\int_{x}\left(\int_{y}g(y,x)\psi(y)\phi(y)\d y\right)\nu^\sharp_\infty(\d x,\d s).\label{id2:mu}
\end{align}
First, to obtain \eqref{id1:mu}, we consider the following equalities, where the limits use convergence in probability: 
\begin{align}
&\int_{s=0}^T\int_{x}\int_{\wt{x}}g(\wt{x},x)\mu^\sharp_\infty(\d \wt{x}, \d x,\d s)\notag\\
&\quad =\P\mbox{-}\lim_{n\to\infty}\int_{s=0}^T\int_{x}\int_{\wt{x}}g(\wt{x},x)\mu^\sharp_n(\d \wt{x}, \d x,\d s)\label{id2:mu-1-1}\\
&\quad =\P\mbox{-}\lim_{n\to\infty}\int_{s=0}^T\int_{x}\int_{y}g(\vep_n y+x,x)\psi(\vep_ny+x)\psi(x)\phi(y)\Lambda_n X_{n}(\vep_n y+x,s)X_n(x,s)\d y\d x\d s\notag\\
&\quad =\P\mbox{-}\lim_{n\to\infty}\int_{s=0}^T\int_{x}\int_{y}g(\vep_n y+x,x)\psi(\vep_ny+x)\psi(x)\phi(y)\Lambda_n  X_{n}(x,s)^2\d y\d x\d s\label{id2:mu-1-3}\\
&\quad =\P\mbox{-}\lim_{n\to\infty}\int_{s=0}^T\int_{x}\left(\int_{y}g(\vep_n y+x,x)\psi(\vep_ny+x)\phi(y)\d y\right)\nu^\sharp_n(\d x,\d s)\notag\\
&\quad =\int_{s=0}^T\int_{x} g(x,x)\psi(x)\nu^\sharp_\infty(\d x,\d s).\label{id2:mu-1-5}
\end{align}
The justifications for \eqref{id2:mu-1-1}, \eqref{id2:mu-1-3}, and \eqref{id2:mu-1-5} are as follows. Equation~\eqref{id2:mu-1-1} follows from the almost-sure convergence of $\mu^\sharp_n$ to $\mu^\sharp_\infty$ in $C_{\mathcal M_f(\R^4)}[0,\infty)$ (see Corollary~\ref{cor:Sk}). Equation~\eqref{id2:mu-1-3} is based on Lemma~\ref{lem:I12} and Lemma~\ref{lem:I2ap} (1$\cc$), using the following choice of $(g,g_\vep)$: 
\[
g(y,x,s,T)=g_\vep(y,x,s,T)\equiv g(\vep y+x,x)\psi(\vep y+x)\psi(x)\phi(y).
\]
Finally, \eqref{id2:mu-1-5} is justified by the uniform bound $\sup_{n}\E[\int_{s=0}^{T}\int_x (1+\lvert x\rvert ^{10})^{-1}\nu_n(\d x,\d s)]<\infty$ (from Theorem~\ref{thm:mombdd} (3$\cc$)) and the almost-sure convergence of $\nu^\sharp_n$ to $\nu^\sharp_\infty$ in $C_{\mathcal M_f(\R^2)}[0,\infty)$. This completes the proof of \eqref{id1:mu}.
The proof of \eqref{id2:mu} is almost the same:
\begin{align}
&\int_{s=0}^T\int_{x}\int_{\wt{x}}g(\wt{x},x)\mathring{\mu}^\sharp_\infty(\d \wt{x}, \d x,\d s)\notag\\
&\quad =\P\mbox{-}\lim_{n\to\infty}\int_{s=0}^T\int_{x}\int_{\wt{x}}g(\wt{x},x)\mathring{\mu}^\sharp_n(\d \wt{x}, \d x,\d s)\label{id2:mu-2-1}\\
&\quad =\P\mbox{-}\lim_{n\to\infty}\int_{s=0}^T\int_{x}\int_{y}g( y,x)\psi(y)\psi(x)\phi(y)\Lambda_n X_{n}(\vep_n y+x,s)X_n(x,s)\d y\d x\d s\notag\\
&\quad =\P\mbox{-}\lim_{n\to\infty}\int_{s=0}^T\int_{x}\int_{y}g( y,x)\psi(y)\psi(x)\phi(y)\Lambda_n  X_{n}(x,s)^2\d y\d x\d s\label{id2:mu-2-3}\\
&\quad =\P\mbox{-}\lim_{n\to\infty}\int_{s=0}^T\int_{x}\left(\int_{y}g(y,x)\psi(y)\phi(y)\d y\right)\nu^\sharp_n(\d x,\d s )\notag\\
&\quad =\int_{s=0}^T\int_{x}\left(\int_{y}g(y,x)\psi(y)\phi(y)\d y\right)\nu^\sharp_\infty(\d x,\d s).\label{id2:mu-2-5}
\end{align}
In more detail, \eqref{id2:mu-2-1} holds by using the almost-sure convergence of $\mathring{\mu}^\sharp_n$ to $\mathring{\mu}^\sharp_\infty$ in $C_{\mathcal M_f(\R^4)}[0,\infty)$. Equation~\eqref{id2:mu-2-3} is justified by Lemma~\ref{lem:I12} and Lemma~\ref{lem:I2ap} (1$\cc$), using the following choice of $(g,g_\vep)$: 
\[
g(y,x,s,T)=g_\vep(y,x,s,T)\equiv g( y,x)\psi(y)\psi(x)\phi(y).
\]
By \eqref{id2:mu-1-5} and \eqref{id2:mu-2-5}, we have verified both  \eqref{id1:mu} and \eqref{id2:mu}. 

The remainder of this proof is to establish the required limit of $I^2_n(f,T)$ for $f\in \C^1_{pd}(\R^2\times [0,T])$, where $I^2_n(f,T)$ has been defined in the statement of Proposition~\ref{prop:I2}. As the first step toward this goal, observe that $g_\vep(y,x,s,T)\equiv \mathfrak L_\vep^1f(y,x,s,T)$ satisfies the assumption in Lemma~\ref{lem:I2ap} (1$\cc$):
\begin{itemize}
\item The $\C^1(\R^2)$ regularity in the variable $x$ follows from Lemma~\ref{lem:I1} (2$\cc$) and the dominated convergence theorem. 
\item To see that $\lvert g_\vep(y,x,s,T)\rvert=\lvert \mathfrak L_\vep^1f(y,x,s,T)\rvert$ can be bounded as in Lemma~\ref{lem:I2ap} (1$\cc$), recall that $\phi$ has a compact support, and $\mathfrak L_\vep^1f(y,x,s,T)=\mathring{\mathfrak L}_\vep^1f(y,x,s,T)\phi(y)$ by \eqref{def:L1}. Hence, by \eqref{Phi1diff:bdd}, the required bound on $\lvert \mathfrak L_\vep^1f(y,x,s,T)\rvert$ holds. 
\item The required bound on $\lvert \nabla_x g_\vep(y, x, s, T)\rvert = \lvert \nabla_x \mathfrak{L}_\vep^1 f(y, x, s, T)\rvert$ follows directly from \eqref{I1:derbdd} and the compact support of $\phi$.
\end{itemize}
Given that $I^{2}_n(f, T)$ admits the decomposition in \eqref{dec0:I2} and Lemma~\ref{lem:I2ap} (1$\cc$) applies to our choice of $g_\vep(y, x, s, T)$, we obtain the first equality below, assuming the existence of the limit: 
\begin{align}
\P\mbox{-}\lim_{n\to\infty}I^{2}_n(f,T)&=\P\mbox{-}\lim_{n\to\infty}-I^{2,4}_{n}(\mathfrak L^1_nf,T)\notag\\
&=\P\mbox{-}\lim_{n\to\infty}-I^{2,4}_{n}\big((y,x,s,T)\mapsto f(x,s)\phi(y),T\big),\label{I2:final}
\end{align}
where the second equality uses \eqref{Phi1diff:bdd} and again Lemma~\ref{lem:I2ap} (1$\cc$). 

The next step is to verify the existence of the limit in probability of $-I^{2,4}_{n}\big((y, x, s, T) \mapsto f(x, s) \phi(y), T\big)$ and to compute its value. For this purpose, we examine alternative limits in probability, as expressed by the following equations: 
\begin{align}
&\P\mbox{-}\lim_{n\to\infty}-I^{2,4}_{n}\big((y,x,s,T)\mapsto f(x,s)\phi(y),T\big)\notag\\
&\quad =\P\mbox{-}\lim_{n\to\infty}-\int_{r'=0}^T \int_{y'}\int_{y''} \int_{s=r'}^T \int_{x}\int_{y}f(x,s)\phi(y)\notag\\
&\quad\quad \times [P_{s-r'}(\vep_n y+x,\vep_n y''+y')-P_{s-r'}(x,\vep_n y''+y')] P_{s-r'}(x,y')\d y \d x\d s \notag\\
&\quad\quad  \times \Lambda_n  X_n(\vep_n y''+y',r') X_n(y',r') \phi(y'')\d y''\d y'
\d r'\notag\\
 &\quad =\P\mbox{-}\lim_{n\to\infty}-\int_{r'=0}^T \int_{y'}\int_{y''} \biggl( \int_{y}\int_{s=r'}^T \int_{x}f(x,s)
 P_{s-r'}(\vep_n y''-\vep_n y+y',x)P_{s-r'}(x,y')\d x\d s \phi(y)\d y \notag \\
&\quad\quad  - \int_{y}\int_{s=r'}^T \int_{x}f(x,s) P_{s-r'}(\vep_n y''+y',x) P_{s-r'}(x,y')
\d x\d s \phi(y)\d y\biggr)\notag\\
&\quad \quad \times \Lambda_n  X_n(\vep_n y''+y',r') X_n(y',r') \phi(y'')\d y''\d y'\d r'
\notag\\
&\quad =\P\mbox{-}\lim_{n\to\infty}-\int_{r'=0}^T \int_{y'}\int_{y''}\int_{y} \biggl(\Lambda_n^{-1}\mathring{\mathfrak L}^1_nf(y''-y,y',r',T) -\Lambda_n^{-1}\mathring{\mathfrak L}^1_nf(y'',y',r',T)\biggr)\phi(y)\d y\notag\\
&\quad\quad  \times \Lambda_n  X_n(\vep_n y''+y',r') X_n(y',r') \phi(y'')\d y''\d y'\d r',\label{I2:proof-istep}
\end{align}
where the first equality applies directly the formula of $I^{2,4}_\vep(g,T)$, as defined in 
Lemma~\ref{lem:I12}, and the last equality follows by using the definition \eqref{def:L1r} of $\mathring{\mathfrak L}^1_nf=\mathring{\mathfrak L}^1_{\vep_n}f$. To continue, we use \eqref{key} to get the first equality below:
\begin{align*}
&\int_{y} \biggl(\Lambda_n^{-1}\mathring{\mathfrak L}^1_nf(y''-y,y',r',T)-\Lambda_n^{-1}\mathring{\mathfrak L}^1_nf(y'',y',r',T)\biggr)\phi(y)\d y\\
&\quad =\left[-\frac{1}{2\pi}\int_{y}(\log \rvert y''-y\rvert )\phi(y)\d y+\frac{1}{2\pi}\int_{y}(\log \lvert y'' \rvert)\phi(y)\d y\right]f(y',r')\\
&\quad \quad +\left(\frac{\Lambda_n^{-1}}{\log \vep_n^{-1}}-\frac{1}{2\pi}\right)\left[-\int_{y}(\log \rvert y''-y\rvert )\phi(y)\d y+\int_{y}(\log \lvert y'' \rvert)\phi(y)\d y\right]f(y',r')\\
&\quad \quad +\int_{y} \biggl(\Lambda_n^{-1}R_{\ref{key}}(y''-y,y',r',T)-\Lambda_n^{-1}R_{\ref{key}}(y'',y',r',T)\biggr)\phi(y)\d y,
\end{align*}
where $R_{\ref{key}}(y,x,s,T)$ satisfies \eqref{R:bdd}. By Lemma~\ref{lem:I2ap}~(2$\cc$) and Theorem~\ref{thm:mombdd} (3$\cc$), applying the last equality to the right-hand side of \eqref{I2:proof-istep} gives
\begin{align}
&\P\mbox{-}\lim_{n\to\infty}-I^{2,4}_{n}\big((y,x,s,T)\mapsto f(x,s)\phi(y),T\big)\notag\\
&\quad =\P\mbox{-}\lim_{n\to\infty}-\int_{r'=0}^T \int_{y'}\int_{y''}\left[-\frac{1}{2\pi}\int_{y}(\log \rvert y''-y\rvert )\phi(y)\d y+\frac{1}{2\pi}\int_{y}(\log \lvert y'' \rvert)\phi(y)\d y\right]f(y',r')\notag\\
&\quad\quad  \times \Lambda_n  X_n(\vep_n y''+y',r) X_n(y',r') \phi(y'')\d y''\d y'\d r'\notag\\
&\quad =\int_{r'=0}^T \int_{y'}\left(\frac{1}{2\pi }\int_{y''}\int_{y}\phi(y'')(\log\rvert y''-y\rvert )\phi(y)\d y \d y''\right) f(y',r')\nu_\infty(\d y',\d r')\notag\\
&\quad\quad  -\int_{r'=0}^T \int_{y'}\left(\frac{1}{2\pi }\int_{y''}\int_{y}\phi(y'')(\log\lvert y'' \rvert)\phi(y)\d y \d y''\right) f(y',r')\nu_\infty(\d y',\d r')\label{I2:finalfinal}
 \end{align}
where the last inequality follows by using \eqref{id:mur=nu}. We have proved
the existence of the limit in probability of $-I^{2,4}_{n}\big((y, x, s, T) \mapsto f(x, s) \phi(y), T\big)$ and found its value. 

Finally, combining \eqref{I2:final} and \eqref{I2:finalfinal} proves the required limit of $I^2_n(f,T)$. The proof of Proposition~\ref{prop:I2} is complete.
\end{proof}

\addtocontents{toc}{\protect\setcounter{tocdepth}{1}}
\subsection{Limit theorem for $\bs I^{\bs 3}_{\bs n}\bs (\bs f\bs,\bs T )$ and $\bs I^{\bs 3}_{\bs n}\bs (\bs f\bs ,\bs T\bs )+\bs I^{\bs 4}_{\bs n}\bs (\bs f\bs ,\bs T\bs )$}\label{sec:I34}
\addtocontents{toc}{\protect\setcounter{tocdepth}{2}}
\addcontentsline{toc}{subsection}{\numberline{\thesubsection}{Limit theorems for $I^3_n(f,T)$ and $I^3_n(f,T)+I^4_n(f,T)$}}

We now consider the limiting formulas for the stochastic integral terms $I^3_n(f, T)$ and $I^3_n(f,T)+I^4_n(f, T)$ as $n \to \infty$. The main result is stated in the following proposition.

\begin{prop}\label{prop:I3I4}
Let $T\in (0,\infty)$ and $f(x, t) \in \C_{pd}^0(\mathbb{R}^2 \times [0, T])$. Recall that $\mathfrak{L}^3_0 f$ and $\mathfrak{L}^4_0 f$ are defined as in \eqref{def:L3} and \eqref{def:L4}. For all $n\in \Bbb N$, define the following stochastic integrals: 
\begin{align}
I^3_n(f,T)&\,\defeq\, \int_{s=0}^T \int_{y} \mathfrak L_0^3f(y,s,T)\M_n(\d y,\d s),\label{def:I3n}\\
I^{4}_n(f,T)&\,\defeq\, \int_{s'=0}^T \int_{y'}  \int_{s=0}^{s'}\int_{y} \mathfrak L_0^4f(y',y,s',s,T) \M_n(\d y,\d s)\M_n(\d y',\d s').\label{def:I4n}
\end{align} 
Then $I^{3,4}_{n}(f,T)\,\defeq\, I^3_n(f,T)+I^4_{n}(f,T)$ satisfies the following formula:
\begin{align}
I^{3,4}_{n}(f,T)&=\int_{s'=0}^T\int_{y'} \biggl( \int_{t=s'}^T \int_x f(x,t)P_{t-s'}(x,y')\notag\\
&\quad \times \int_z P_{t-s'}(x,z)X_n(\d z,s')\d x\d t\biggr)\M_n(\d y',\d s'),\quad \forall\;n\in\Bbb N,\label{eq:I34n}
\end{align}
and the following limits in probability hold: 
\begin{align}
I^3_n(f,T)
& \xrightarrow[n\to\infty]{\P}I^3_\infty(f,T)\,\defeq\,\int_{s=0}^T \int_{y} \mathfrak L_0^3f(y,s,T)\M_\infty(\d y,\d s),\label{I3:conv}\\
I^{3,4}_n(f,T)
&\xrightarrow[n\to\infty]{\P}I^{3,4}_\infty(f,T)\,\defeq\,\int_{s'=0}^T\int_{y'} \biggl( \int_{t=s'}^T \int_x f(x,t)P_{t-s'}(x,y')\notag\\
&\quad \quad \quad \times\int_z P_{t-s'}(x,z)X_\infty(\d z,s')\d x\d t\biggr)\M_\infty(\d y',\d s').\label{I4:conv}
\end{align}
Here, the stochastic integrals in \eqref{def:I3n}--\eqref{I4:conv} are well-defined, and the following bounds hold:
\begin{gather}
\sup_{n\in \Bbb N\cup\{\infty\}}\E[I^3_n(f,T)^2]<\infty,\label{I3n:sup}\\
\sup_{n\in \Bbb N\cup\{\infty\}}\E[I^{3,4}_n(f,T)^2]<\infty.\label{I4n:sup}
\end{gather}
\end{prop}

The proof of Proposition~\ref{prop:I3I4} relies on Theorem~\ref{thm:main1} (2$\cc$), Lemma~\ref{lem:Mfconv}, and Lemma~\ref{lem:KP}. The main goal of Lemma~\ref{lem:Mfconv} is to specify a class of functions integrable against the martingale measure $M_\infty$ and establish the limit behavior of $M_n$ as $n \to \infty$. Although $M_\infty$ is 
defined via \eqref{XM:relation} using only $\C^2_{pd}(\R^2)$-test functions, it uniquely determines a measure $\langle M_\infty, M_\infty\ra (\d x, \d s)$ on $\B(\R^2)\otimes \B(\R_+)$ such that the first and second processes in
\begin{align}\label{Minfty:mgproperty}
\{M_\infty(f,t)\}_{t\geq 0},\quad \left\{M_\infty(f,t)^2-\int_0^t\int_{\R} \la M_\infty,M_\infty\ra (\d x,\d s)\right\}_{t\geq 0}
\end{align}
are $L^4$- and $L^2$-martingales, respectively, for any $f\in \C^2_{pd}(\R^2)$. Part (1$\cc$) of Lemma~\ref{lem:Mfconv} establishes the existence of $\la M_\infty,M_\infty\ra (\d x,\d s)$, while introducing an identity that links this measure to $\mu_\infty(\d x', \d x, \d s)$, and subsequently to $\nu_\infty(\d x,\d s)$ and $\mathring{\mu}_\infty(\d x',\d x,\d s)$ via \eqref{id:mu=nu} and \eqref{id:mur=nu}. Part (2$\cc$) 
uses $\la M_\infty, M_\infty\ra(\d x, \d s)$ to define an auxiliary norm, thereby enabling the standard $L^2$-construction of stochastic integrals against $M_\infty$ under conditions analogous to those in \eqref{I3n:sup} and \eqref{I4n:sup} for $n = \infty$. For our purposes, part (2$\cc$) restricts focus to $\C^0_{pd}(\R^2)$-test functions, for which the $L^4$- and $L^2$-martingale properties of the processes in \eqref{Minfty:mgproperty}  remain valid. Finally, parts (3$\cc$) and (4$\cc$) address the limiting behavior of the projected martingale measures $M_n$ as $n\to\infty$ under different test function regularities, yielding distinct modes of convergence that apply directly to our proof of Proposition~\ref{prop:I3I4}.

\begin{lem}\label{lem:Mfconv}
 {\rm (1$\cc$)}  The martingale measure $M_\infty$ is orthogonal, meaning that its covariation measure $\la M_\infty, M_\infty \ra(\d x', \d x, \d s)$ is supported on $\{(x', x) \in \R^2 \times \R^2; x' = x\} \times \R_+$. Accordingly, it is written as $\la M_\infty, M_\infty \ra(\d x, \d s)$ henceforth. Additionally, for all $T\in (0,\infty)$ and nonnegative $g\in \C(\R^2\times \R^2\times[0,T])$, 
\begin{align}\label{Mncovar:derivation:eq}
\int_{s=0}^T\int_x g(x,x,s)\la M_\infty,M_\infty\ra(\d x,\d s)=\int_{s=0}^T\int_x \int_{x'} g(x',x,s)\mu_\infty(\d x',\d x,\d s).
\end{align}

\noindent {\rm (2$\cc$)} Let $f\in \C_{pd}^0(\R^2)$.
Whenever $f_n\in \C_{pd}^2(\R^2)$, $n\in \Bbb N$, approximate $f$ in the following norm:
\begin{align}\label{eq:Mfconv2}
\lim_{n\to\infty}\E\left[\int_0^T\int_x [f_n(x)-f(x)]^2 \la M_\infty,M_\infty\ra(\d x,\d s)\right]=0,\quad \forall\; T\in (0,\infty),
\end{align}
the sequence $\{M_\infty(f_n,t);t\geq 0\}_{n\in \Bbb N}$ converges uniformly on compacts in probability as $n\to\infty$. The limiting process does not depend on the choice of $\{f_n\}_{n\in \Bbb N}$ and defines a martingale denoted by $\{M_\infty(f,t);t\geq 0\}$. Moreover, the corresponding processes defined in \eqref{Minfty:mgproperty} are $L^4$- and $L^2$-martingales. \medskip 

\noindent {\rm (3$\cc$)} With probability one, the following limit holds for all $f\in \C_{pd}^2(\R^2)$:
\begin{align}\label{M:conv}
\{M_n(f,t)\}_{t\geq 0}\xrightarrow[n\to \infty]{}\{M_\infty(f,t)\}_{t\geq 0}\mbox{ in }C_{\R}[0,\infty),
\end{align}
where $C_{\R}[0,\infty)$ is equipped with the topology of uniform convergence on compacts. \medskip

\noindent {\rm (4$\cc$)} For any $T\in (0,\infty)$ and $f\in \C^0_{pd}(\R^2)$, the following limit in probability holds:
\begin{align}\label{Mf:convunif}
\sup_{0\leq t\leq T}\lvert M_n(f,t)-M_\infty(f,t)\rvert\xrightarrow[n\to\infty]{\P}0.
\end{align}
\end{lem}

\begin{proof}
(1$\cc$) The aim of this proof is to relate the covariation measure of $M_\infty$ to $\mu_\infty$. For this purpose, observe that for any $n\in \Bbb N$, $ T\in (0,\infty)$ and $f \in \C_{pd}^2(\R^2)$, 
the two processes below are martingales: 
\begin{align}\label{eq1:Mfconv}
\{M_n(f,t)\}_{0\leq t\leq T}, \quad \left\{M_n(f, t)^2 - \int_{s=0}^t \int_{x} \int_{x'} f(x') f(x) \mu_n(\d x', \d x, \d s)\right\}_{0\leq t\leq T}.
\end{align}
This is a direct consequence of the definitions of $M_n$ and $\mu_n$ for finite $n$. Additionally, the martingales in \eqref{eq1:Mfconv} for all $n\in \Bbb N$ are uniformly integrable since
\begin{gather}
\sup_{n\in \Bbb N} \sup_{0\leq t\leq T}\E[M_n(f, t)^4]<\infty,\label{Mncovar:derivation1}\\
\sup_{n\in \Bbb N}\sup_{0\leq t\leq T}\E\left[\left(M_n(f, t)^2 - \int_{s=0}^t \int_{x} \int_{x'} f(x') f(x) \mu_n(\d x', \d x, \d s)\right)^2\right]<\infty.\label{Mncovar:derivation2}
\end{gather}
Both bounds hold since by the Burkholder--Davis--Gundy inequality \cite[(4.1) Theorem, p.160]{RY}, 
\begin{align}
\sup_{0\leq t\leq T}\E[M_n(f, t)^4]&\less \E[\langle M_n(f, \cdot),M_n(f, \cdot)\ra_T^2]\notag\\
&=\E\left[\left(\int_{s=0}^T \int_{x} \int_{x'} f(x') f(x) \mu_n(\d x', \d x, \d s)\right)^2\right]\notag\\
&\leq C(\lambda,\phi,\lVert X_0\rVert_\infty,T,p)\lVert f\rVert_{L^1(\R^2)}^2,
\label{eq:Mfconv1-1}
\end{align}
where the last inequality follows from \eqref{def:covarprocessmu} and Theorem~\ref{thm:mombdd} (1$\cc$). 

We now prove the required properties. By \eqref{Mncovar:derivation1} and \eqref{Mncovar:derivation2}, a standard result of uniform integrability implies that the martingale property of the two processes in \eqref{eq1:Mfconv} extends to $n = \infty$. Therefore, $\int_{s = 0}^t \int_{x} \int_{x'} f(x') f(x) \mu_\infty(\d x', \d x, \d s)$, for $0 \leq t \leq T$, must coincide with the quadratic variation of $M_\infty(f, t)$. Consequently, we deduce that the measure $\mu_\infty(\d x', \d x, \d s)$, which is diagonal by \eqref{def:covarprocessmu}, represents the covariation measure of $M_\infty$ via \eqref{Mncovar:derivation:eq}. This is enough to establish all the required properties in (1$\cc$). \medskip 
 
\noindent (2$\cc$) All required properties, except the $L^4$- and $L^2$-boundeness of
 the corresponding processes defined in \eqref{Minfty:mgproperty}, follow directly from Doob's weak $L^2$-inequality \cite[(1.7) Theorem, p.54]{RY}. The $L^4$- and $L^2$-boundedness can be deduced from \eqref{eq:Mfconv1-1}.  \medskip

\noindent {(3$\cc$)} By \eqref{XM:relation}, the required property follows directly from Corollary~\ref{cor:muinfty}.
\medskip 

\noindent {(4$\cc$)} Any $f \in \C^0_{pd}(\R^2)$ vanishes at infinity. Therefore, by mollification, there exists a sequence $\{f_k\} \subset \C^2_{pd}(\R^2)$ such that $f_k \to f$ uniformly on $\R^2$ as $k \to \infty$, and $\sup_{k\in \Bbb Z_+}\lvert f_k(x)\rvert\leq C(f,m)(1+\lvert x\rvert ^m)^{-1}$ for all $m\in \Bbb Z_+$. (This bound is justified directly by \eqref{ratiot}.) Moreover, for any $T, \eta, \vep > 0$, there exists an integer $k_0 \geq 1$ such that
\begin{align}\label{part4:Mfconv}
\sup_{n\in \Bbb N\cup\{\infty\}}\E\biggl[\sup_{0\leq t\leq T}\Big(M_n(f,t)-M_n(f_{k_0},t)\Big)^2\biggr]<\vep\eta^2.
\end{align}
This follows by applying Doob's strong $L^2$-inequality, in combination with \eqref{eq:covarM-2} and Theorem~\ref{thm:mombdd} (3$\cc$) for $n \in \mathbb{N}$, and with an extension to \eqref{eq:Mfconv1-1} to $n=\infty$ for $n = \infty$. For the same $T,\eta,\vep>0$,
\begin{align*}
&\P\biggl(\sup_{0\leq t\leq T}\lvert M_n(f,t)-M_\infty(f,t)\rvert>3\eta\biggr)\\
&\quad\leq \frac{1}{\eta^2}\E\biggl[\sup_{0\leq t\leq T}\Big(M_n(f,t)-M_n(f_{k_0},t)\Big)^2\biggr] +\P\biggl(\sup_{0\leq t\leq T}\lvert M_n(f_{k_0},t)-M_\infty(f_{k_0},t)\rvert>\eta\biggr)+\\
&\quad\quad +\frac{1}{\eta^2}\E\biggl[\sup_{0\leq t\leq T}\Big(M_\infty(f_{k_0},t)-M_\infty(f,t)\Big)^2\biggr]\\
&\quad \leq 2\vep +\P\biggl(\sup_{0\leq t\leq T}\lvert M_n(f_{k_0},t)-M_\infty(f_{k_0},t)\rvert>\eta\biggr)\xrightarrow[n\to\infty]{}2\vep,
\end{align*}
where the first inequality uses Markov's inequality, the second inequality uses \eqref{part4:Mfconv},
and the limit follows from (3$\cc$) since $f_{k_0}\in \C_{pd}^2(\R^2)$. Since $\vep>0$ is arbitrary, the foregoing limit is enough to get \eqref{Mf:convunif}.
The proof is complete.
\end{proof}

The next lemma collects key results on the uniform approximation of c\`adl\`ag functions with step functions, as established in \cite[pp.1066--1067]{KP:91}. The general framework and main properties can be restated as follows. Let $(E,r)$ be a metric space, and consider a sequence of i.i.d. random variables $\{U_k(\omega')\}_{k\in \Bbb Z_+}$ defined on a probability space $(\Omega',\F',\P')$, each uniformly distributed over the interval $[\frac{1}{2},1]$. For convenience, assume that $U_k(\omega')$ takes values in the same interval for all $\omega'$. Given any $\theta > 0$ and $\mathfrak{z} \in D_E[0, \infty)$, inductively define the following sequence of times for all $k \in \mathbb{Z}_+$:
\begin{subequations}\label{def:tautheta}
\begin{align}
\tau^\theta_0(\mathfrak z)&=\tau^{\theta,\omega'}_0(\mathfrak z)\,\defeq\,0,\\
\tau^\theta_{k+1}(\mathfrak z)&=\tau^{\theta,\omega'}_{k+1}(\mathfrak z)\,\defeq\,\inf\left\{t>\tau^\theta_k(\mathfrak z);r\big(\mathfrak z(t),\mathfrak z(\tau^\theta_k(\mathfrak z))\big)\vee r\big(\mathfrak z(t-),\mathfrak z(\tau^\theta_k(\mathfrak z))\big)\geq \theta U_k\right\}.
\end{align}
\end{subequations}
The associated step function that approximates $\mathfrak{z}$ is defined as follows:
\begin{align}\label{def:Ithetaz}
I^\theta(\mathfrak z)(t)=I^{\theta,\omega'}(\mathfrak z)(t)\,\defeq\,\sum_{k=0}^\infty\mathfrak z(\tau^\theta_k(\mathfrak z))\1_{ [\tau^\theta_k(\mathfrak z), \tau^\theta_{k+1}(\mathfrak z))}(t).
\end{align}

\begin{lem}[Kurtz and Protter~\cite{KP:91}]\label{lem:KP}
Fix $\theta>0$.\smallskip 

\noindent {\rm (1$\cc$)} For all $\mathfrak z\in D_E[0,\infty)$, $\tau^{\theta,\omega'}_k(\mathfrak z)\nearrow\infty$ as $k\to\infty$ for all $\omega'\in \Omega'$.\medskip 

\noindent {\rm (2$\cc$)} For all $\mathfrak z\in D_E[0,\infty)$,
 $r\big(\mathfrak z(t),I^{\theta,\omega'}(\mathfrak z)(t)\big)\leq \theta$ for all $t\geq 0$ and $\omega'\in \Omega'$.\medskip 

\noindent {\rm (3$\cc$)} For every $\mathfrak z\in D_E[0,\infty)$ and every sequence $\{\mathfrak z_n\} \subset D_E[0, \infty)$ with $\mathfrak z_n \to \mathfrak z$ in the Skorokhod topology, the following properties hold with $\P'$-probability one:
\begin{itemize}
\item For every $k\geq 1$, $\tau^{\theta,\omega'}_k(\mathfrak z_n)\to \tau^{\theta,\omega'}_k(\mathfrak z)$ as $n\to\infty$. 
\item For every $k\geq 1$, $\mathfrak z_n(\tau^{\theta,\omega'}_k(\mathfrak z_n))\to \mathfrak z(\tau^{\theta,\omega'}_k(\mathfrak z))$ as $n\to\infty$. 
\item As $n\to\infty$, $(\mathfrak z_n,I^{\theta,\omega'}(\mathfrak z_n))\to (\mathfrak z,I^{\theta,\omega'}(\mathfrak z))$ in the Skorokhod topology on $D_{E\times E}[0,\infty)$. 
\end{itemize}
\end{lem}

To use Lemma~\ref{lem:KP} in the proof of Proposition~\ref{prop:I3I4}, we take $(E, r)$ to be the metric space associated with the normed space $(E_m, \lVert\cdot\rVert_m)$, where $m \in \mathbb{N}$ is chosen sufficiently large. Here, $\C_0(\mathbb{R}^2)$ denotes the set of continuous functions $f(x)$ on $\mathbb{R}^2$ that vanish at infinity, i.e., $\lim_{\lvert x\rvert  \to \infty} f(x) = 0$, and we define
\begin{align}\label{def:Em}
E_m\,\defeq\,\{f\in \C_0(\R^2);\lVert f\rVert_m<\infty\},
\end{align}
where
\[
\lVert f\rVert_m\,\defeq\,\sup_{x\in \R^2}\lvert f(x)\rvert(1+\lvert x\rvert ^m).
\]
Observe that for any $m \in \mathbb{N}$, the space $E_m$, endowed with the metric induced by $\lVert \cdot\rVert_m$, is a Polish space. In particular, its separability follows from the Stone--Weierstrass theorem.

We now turn to the proof of Proposition~\ref{prop:I3I4}, verifying \eqref{eq:I34n}, \eqref{I3n:sup}, \eqref{I3:conv}, \eqref{I4n:sup}, and \eqref{I4:conv} individually and in that order. \medskip  

\begin{proof}[Proof of (\ref{eq:I34n}) in Proposition~\ref{prop:I3I4}]
It suffices to show that for any $n\in \Bbb N$, 
\begin{align}\label{eq:4=-3+34}
I^4_n(f,T)=-I^{3}_n(f,T)+\wt{I}^{\,3,4}_n(f,T),
\end{align}
where $\wt{I}^{\,3,4}_n(f,T)$ denotes the term on the right-hand side of \eqref{eq:I34n}. 

We begin by deriving an equivalent expression of $I^4_n(f,T)$.
By the definition of $I^4_n(f,T)$ and the stochastic Fubini theorem, 
\begin{align*}
I^4_n(f,T)&=\int_{t=0}^T \int_x f(x,t)\int_{s'=0}^t\int_{y'}\int_{s=0}^{s'}\int_y P_{t-s'}(x,y')P_{t-s}(x,y)\M_n(\d y,\d s)\M_n(\d y',\d s')\d x\d t,
\end{align*}
which is also equivalent to \eqref{SMF:3} if we set $\vep$ there to be $\vep_n$. By the Chapman-Kolmogorov equation, the last equality gives
\begin{align*}
I^4_n(f,T)
&= \int_{t=0}^T \int_x f(x,t)\int_{s'=0}^t \int_{y'}P_{t-s'}(x,y')\\
&\quad \times \bigg(\int_{s=0}^{s'}\int_y \int_z P_{t-s'}(x,z)P_{s'-s}(z,y)\d z\M_n(\d y,\d s)\biggr)\M_n(\d y',\d s')\d x\d t\\
&= \int_{t=0}^T \int_x f(x,t)\int_{s'=0}^t \int_{y'}P_{t-s'}(x,y')\\
&\quad \times \biggl(\int_z P_{t-s'}(x,z)[-P_{s'}X_0(z)+X_n(z,s')]\d z\biggr)\M_n(\d y',\d s')\d x\d t\\
&=- \int_{s'=0}^T \int_{y'} \biggl(\int_{t=s'}^T \int_x f(x,t)P_{t-s'}(x,y')\int_z P_{t-s'}(x,z)P_{s'}X_0(z)\d z \d x\d t\biggr)\M_n(\d y',\d s')\\
&\quad+\int_{s'=0}^T\int_{y'} \biggl( \int_{t=s'}^T \int_x f(x,t)P_{t-s'}(x,y')\int_z P_{t-s'}(x,z)X_n(z,s')\d z\d x\d t\biggr)\M_n(\d y',\d s').
\end{align*}
Here, the second equality applies the mild form equation \eqref{def:mild:vep} tested against the function $z \mapsto P_{t-s'}(x, z)$ for $s' < t$. The final equality uses the stochastic Fubini theorem to interchange the order of integration,
\[
\M_n(\d y', \d s')\, f(x,t)\d x\, \d t =f(x,t) \d x\, \d t\, \M_n(\d y', \d s'). 
\]
This application can be justified by Proposition~\ref{prop:XM} (2$\cc$) and \eqref{SFT:main}.

Equation~\eqref{eq:4=-3+34} follows from the previous equality because the measure $X_n(\d z, s')$ is defined to have density $X_n(z, s')$ with respect to Lebesgue measure, and we have
\begin{align*}
 \int_{t=s'}^T \int_x f(x,t)P_{t-s'}(x,y')\int_z P_{t-s'}(x,z)P_{s'}X_0(z)\d z \d x\d t=\mathfrak L_0^3f(y',s',T).
 \end{align*}
This uses the Chapman--Kolmogorov equation and the definition of $\mathfrak L_0^3f$ in \eqref{def:L3}. The proof is complete.
\end{proof}

\begin{proof}[Proof of (\ref{I3n:sup}) in Proposition~\ref{prop:I3I4}]
The definition \eqref{def:L3} of $\mathfrak L_0^3f$ gives the following equation: 
\begin{align}\label{Phi3:alt}
\mathfrak L_0^3f(y,s,T)=\int_{t=0}^{T-s} \int_x f(x,t+s)P_{t+s}X_0(x)P_{t}(x,y) \d x\d t,\quad 0\leq s\leq T.
\end{align}
Since we assume $f\in \C_{pd}^0(\R^2\times [0,T])$ and bounded $X_0(\cdot)$, by \eqref{ineq:heat2-2}, \eqref{Phi3:alt} implies
\begin{align}\label{I3:ito-iso1}
\lvert \mathfrak L_0^3f(y,s,T)\rvert\leq \frac{ C(\lVert X_0\rVert_\infty,T,f,m)}{1+\lvert y \rvert^m},\quad 0\leq s\leq T,\;y\in \R^2,\;m\in \Bbb Z_+.
\end{align}
Additionally, applying Itô's isometry yields
\begin{gather*}
\E[I^3_n(f,T)^2]
=\E\biggl[\int_{s=0}^T \int_{x}\int_{x'}
\mathfrak L_0^3f(x',s,T)\mathfrak L_0^3f(x,s,T)\la M_n,M_n\ra(\d x',\d x,\d s)\biggr],\quad n\in \Bbb N,\\
\E[I^3_\infty(f,T)^2]
=\E\biggl[\int_{s=0}^T \int_{x}
\mathfrak L_0^3f(x,s,T)^2\la M_\infty,M_\infty\ra(\d x,\d s)\biggr].
\end{gather*}
Consequently, \eqref{I3n:sup} follows upon applying Theorem~\ref{thm:mombdd} (3$\cc$) for $n\in \Bbb N$ and Theorem~\ref{thm:main1} (2$\cc$) for $n=\infty$. 
\end{proof}

\begin{proof}[Proof of (\ref{I3:conv}) in Proposition~\ref{prop:I3I4}]
First, we show that $s \mapsto \mathfrak L_0^3f(\cdot, s \wedge T, T)$ takes values in $C_{E_m}[0, \infty)$, where $m \geq 10$ is an arbitrary integer, and $E_m$ is defined in \eqref{def:Em} and equipped with the metric induced by $\lVert \cdot\rVert_m$. For any fixed $0 \leq s \leq T$, by substituting $x' = x - y$ in \eqref{Phi3:alt} and applying both \eqref{ineq:heat2-2} and the dominated convergence theorem, we see that the map $y \mapsto \mathfrak L_0^3f(y, s, T)$ belongs to $E_m$. Thus, it remains to show that $s\mapsto \mathfrak L_0^3f(\cdot,s\wedge T,T)$ is continuous with respect to the norm $\lVert \cdot\rVert_m$. To establish this, observe that by \eqref{ineq:heat2-2},
\begin{align}
\sup_{y\in \R^2}(1+\lvert y \rvert^m)\int_{t=0}^{\delta} \int_x  \frac{P_{t}(x,y)}{1+\lvert x\rvert ^m} \d x\d t\searrow 0,\quad \delta\searrow 0.\label{Phi3:small}
\end{align}
Therefore, by \eqref{Phi3:alt}, verifying the required continuity reduces to showing 
 that for any $0<\delta<T$,
\begin{align}\label{Phi3:cont}
\lim_{s_2-s_1\searrow  0\atop 0\leq s_1\leq s_2\leq T-\delta}\lVert \mathfrak  L_0^3f(\cdot,s_2,T)-\mathfrak L_0^3f(\cdot,s_1,T)\rVert_m=0.
\end{align}
To this end, we write $\mathfrak L_0^3f(y,s_2,T)-\mathfrak L_0^3f(y,s_1,T)$ as the sum of two differences using the original definition \eqref{def:L3} of $\mathfrak L_0^3f$:
\begin{align*}
&\mathfrak L_0^3f(y,s_2,T)-\mathfrak L_0^3f(y,s_1,T)\\
&\quad =\biggl(\int_{t=s_2}^T\int_{x}f(x,t)P_tX_0(x) P_{t-s_2}(x,y)\d x \d t-\int_{t=s_2}^T\int_{x}f(x,t)P_tX_0(x) P_{t-s_1}(x,y)\d x \d t\biggr)\\
&\quad \quad +\biggl(\int_{t=s_2}^T\int_{x}f(x,t)P_tX_0(x) P_{t-s_1}(x,y)\d x \d t-\int_{t=s_1}^T\int_{x}f(x,t)P_tX_0(x) P_{t-s_1}(x,y)\d x \d t\biggr).
\end{align*}
Here, as $s_2-s_1\searrow 0$, the $\lVert\cdot\rVert_m$-norm of the first difference vanishes by \eqref{Phi3:small}, the second bound in \eqref{kernel:tmod}, and \eqref{ineq:heat2-2}, while the $\lVert\cdot\rVert_m$-norm of the second difference vanishes by \eqref{Phi3:small}. Therefore, $s \mapsto \mathfrak L_0^3f(\cdot, s \wedge T, T)$ indeed maps into $C_{E_m}[0, \infty)$.

To prove \eqref{I3:conv}, we proceed as follows. For $\theta>0$ and $m\geq 10$ as above, define an $E_m$-valued step function by 
\begin{align}\label{def:Phi3theta}
\mathfrak L_0^3f_\theta(y,s,T)\,\defeq\, \sum_{k=0}^\infty \mathfrak L_0^3f(y,\tau^\theta_{k},T)\1_{[\tau^\theta_{k},\tau^{\theta}_{k+1})}(s),
\end{align}
where $\tau^\theta_{k}= \tau_k^\theta(s\mapsto \mathfrak L_0^3f(\cdot,s,T))$ is as defined before Lemma~\ref{lem:KP} with respect to the metric space $E_m$. By using Markov's inequality and Itô's isometry, we obtain:
\begin{align*}
&\P\left(\left|\int_{s=0}^T\int_y\mathfrak L_0^3f_\theta(y,s,T)\M_n(\d y,\d s)-\int_{s=0}^T\int_y\mathfrak L_0^3f(y,s,T)\M_n(\d y,\d s)\right|>\delta\right)\\
&\quad\leq \frac{1}{\delta^2}\E\left[\int_{s=0}^T\int_{x}\int_{x'} (\mathfrak L_0^3f_\theta-\mathfrak L_0^3f)(x',s,T)(\mathfrak L_0^3f_\theta-\mathfrak L_0^3f)(x,s,T)\mu_n(\d x', \d x,\d s)\right]\\
&\quad\leq \frac{1}{\delta^2}\E\left[\int_{s=0}^T\int_{x}\int_{x'}\frac{\theta^2}{(1+\lvert x'\rvert^m)(1+\lvert x\rvert ^m)}\mu_n(\d x', \d x,\d s)\right],
\end{align*}
where the last inequality uses Lemma~\ref{lem:KP} (2$\cc$) and the definition of $\lVert\cdot\rVert_m$.
Hence, by Theorem~\ref{thm:mombdd} (3$\cc$) for $n\in \Bbb N$ and Theorem~\ref{thm:main1} (2$\cc$) for $n=\infty$,
\begin{equation}\label{I3:unifPbdd}
\begin{gathered}
\forall\;\delta,\eta>0,\;\exists\;\theta_0>0\mbox{ such that }\forall\;n\in \Bbb N\cup\{\infty\},\\
\P\left(\left|\int_{s=0}^T\int_y\mathfrak L_0^3f_{\theta_0}(y,s,T)\M_n(\d y,\d s)-\int_{s=0}^T\int_y\mathfrak L_0^3f(y,s,T)\M_n(\d y,\d s)\right|>\delta\right)\leq \eta.
\end{gathered}
\end{equation}
Furthermore, $y\mapsto \mathfrak L_0^3f(y,\tau^\theta_{k},T)\in \C_{pd}^0(\R^2)$ by \eqref{I3:ito-iso1} for any $\theta>0$ and $k\in \Bbb Z_+$, and there is an integer $N'\geq 1$ independent of $(y,s)$ such that 
\[
\mathfrak L_0^3f_\theta(y,s,T)=\sum_{k=0}^{N'} \mathfrak L_0^3f(y,\tau^\theta_{k},T)\1_{[\tau^\theta_{k},\tau^{\theta}_{k+1})}(s),\quad \forall\;0\leq s\leq T.
\]
Hence, by Lemma~\ref{lem:Mfconv} (4$\cc$) with the function $x\mapsto \mathfrak L_0^3f(x,\tau^\theta_{k},T)$ for every $1\leq k\leq N'$,
\begin{align}\label{I3:unifPbdd1}
\int_{s=0}^T\int_y\mathfrak L_0^3f_\theta(y,s,T)\M_n(\d y,\d s)\xrightarrow[n\to\infty]{\P}\int_{s=0}^T\int_y\mathfrak L_0^3f_\theta(y,s,T)\M_\infty(\d y,\d s),\quad \forall\;\theta>0. 
\end{align}
Combining \eqref{I3:unifPbdd} and \eqref{I3:unifPbdd1} proves \eqref{I3:conv}.
\end{proof}

The proofs of \eqref{I4n:sup} and \eqref{I4:conv} are based on Lemma~\ref{lem:I4alt} below, which provides an a priori bound.

\begin{lem}\label{lem:I4alt}
Let $T\in (0,\infty)$, $p\in (0,1)$, and $ f(x,t)\in \B(\R^2\times [0,T])$ with $0\leq f(x,t)\leq C(f,m,T)/(1+\lvert x\rvert ^m)$ for some integer $m\geq 10$.
Then for any $n\in \Bbb N\cup \{\infty\}$, any $0\leq T_1\leq T_2\leq T$, any function $T_3(s)$ such that $0\leq s\leq T_3(s)\leq T$, the following bound holds:
\begin{align}
&\E\biggl[\int_{s={T_1}}^{T_2} \int_{x}\int_{x'} \biggl(\int_{t'''=s}^{T_3(s)} \int_{x'''} f(x''',t''')P_{t'''-s}(x''',x')\int_{z'} P_{t'''-s}(x''',z')X_n(\d z',s)\d x'''\d t'''\biggr) \notag\\
&\quad \times \biggl(\int_{t''=s}^{T_3(s)} \int_{x''} f(x'',t'')P_{t''-s}(x'',x)\int_{z} P_{t''-s}(x'',z)X_n(\d z,s)\d x''\d t''\biggr)\mu_n(\d x',\d x,\d s)\biggr]\notag\\
&\quad \leq \int_{s=T_1}^{T_2}\int_{t''=0}^{T_3(s)-s}\int_{t'''=0}^{T_3(s)-s}\frac{ C(\lambda,\phi,\lVert X_0\rVert_\infty,T,p,f,m)}{\min\{t''',t''\}^p}\d t'''\d t''\d s.\label{I4:approx2}
\end{align}
\end{lem}
\begin{proof}
To prove \eqref{I4:approx2}, we first consider the case where $n \in \mathbb{N}$. In this situation,
we use \eqref{eq:covarM-2} to get the following equality: 
\begin{align}
&\E\biggl[\int_{s={T_1}}^{T_2} \int_{x}\int_{x'}\biggl( \int_{t'''=s}^{T_3(s)} \int_{x'''} f(x''',t''')P_{t'''-s}(x''',x')\int_{z'} P_{t'''-s}(x''',z')X_n(\d z',s)\d x'''\d t'''\biggr)\notag\\
&\quad \times\biggl( \int_{t''=s}^{T_3(s)} \int_{x''} f(x'',t'')P_{t''-s}(x'',x)\int_{z} P_{t''-s}(x'',z)X_n(\d z,s)\d x''\d t'' \biggr) \mu_n(\d x',\d x,\d s)\biggr]\notag\\
&\quad =\int_{s=T_1}^{T_2} \int_{x}\int_{\wt{x}}\int_{t''=s}^{T_3(s)}\int_{t'''=s}^{T_3(s)}\int_{x''}\int_{x'''} f(x''',t''')f(x'',t'')P_{t'''-s}(x''',\vep_n \wt{x}+x) P_{t''-s}(x'',x)\notag\\
&\quad\quad  \times \biggl(\int_{z}\int_{z'}P_{t'''-s}(x''',z') P_{t''-s}(x'',z)\E[X_n( z',s) X_n (z,s)\notag\\
&\quad\quad  \times \Lambda_{n} X_n(\vep_n \wt{x}+x,s)X_n(x,s)]\d z'\d z\biggr)\phi(\wt{x})\d x'''\d x''\d t'''\d t''\d \wt{x}\d x\d s\notag\\
&\quad \leq \int_{s=T_1}^{T_2}\int_{x}\int_{\wt{x}}\int_{t''=s}^{T_3(s)}\int_{t'''=s}^{T_3(s)}\int_{x''}\int_{x'''} \frac{C(f,m,T)P_{t'''-s}(x''',\vep_n \wt{x}+x) P_{t''-s}(x'',x)}{(1+\lvert x'''\rvert^m)(1+\lvert x''\rvert^m)} \notag\\
&\quad\quad  \times \frac{ C(\lambda,\phi,\lVert X_0\rVert_\infty,T,p)}{(\min\{t''',t''\}-s)^p}\left(\int_{\ol{z}}\phi(\ol{z})\lvert \log \lvert \wt{x}-\ol{z}\rvert \rvert \d\ol{z}+1\right)\phi(\wt{x})\d x'''\d x''\d t'''\d t''\d \wt{x}\d x\d s\notag\\
&\quad \leq \int_{s=T_1}^{T_2}\int_{x}\int_{\wt{x}}\int_{t''=0}^{T_3(s)-s}\int_{t'''=0}^{T_3(s)-s}\frac{C(f,m,T)(1+\lvert\vep_n\wt{x}\rvert)}{(1+\lvert x\rvert ^m)(1+\lvert x\rvert ^m)}\notag
 \\
&\quad \quad \times \frac{ C(\lambda,\phi,\lVert X_0\rVert_\infty,T,p)}{\min\{t''',t''\}^p}\left(\int_{\ol{z}}\phi(\ol{z})\lvert\log \lvert\wt{x}-\ol{z}\rvert\rvert\d\ol{z}+1\right)\phi(\wt{x})\d t'''\d t''\d \wt{x}\d x\d s.\label{I4:approx2proof}
\end{align}
Here, the penultimate inequality follows from Theorem~\ref{thm:mombdd} (2$\cc$), while \eqref{I4:approx2proof} relies on \eqref{ratiot} and \eqref{ineq:heat2-2}. Thus, \eqref{I4:approx2} holds for all $n \in \mathbb{N}$ by \eqref{I4:approx2proof}. For $n=\infty$, the proof is straightforward since by Theorem~\ref{thm:main1} (2$\cc$) and \eqref{Mncovar:derivation:eq}, the left-hand of \eqref{I4:approx2} with $n=\infty$ is bounded by
\begin{align*}
 &\int_{s=T_1}^{T_2}\int_{x}\int_{t''=s}^{T_3(s)}\int_{t'''=s}^{T_3(s)}\int_{x''}\int_{x'''} \frac{C(f,m,T)P_{t'''-s}(x''',x) P_{t''-s}(x'',x)}{(1+\lvert x'''\rvert^m)(1+\lvert x''\rvert^m)} \notag\\
&\quad  \times \frac{ C(\lambda,\phi,\lVert X_0\rVert_\infty,T,p)}{(\min\{t''',t''\}-s)^p}\d x'''\d x''\d t'''\d t''\d x\d s,
\end{align*}
which can also be bounded by the right-hand side of \eqref{I4:approx2} using an argument analogous to the one for deriving \eqref{I4:approx2proof}. Therefore, \eqref{I4:approx2} holds for all $n \in \mathbb{N} \cup \{\infty\}$. 
\end{proof}

\begin{proof}[Proof of (\ref{I4n:sup}) in Proposition~\ref{prop:I3I4}]
Recall that $I^{3,4}_n(f,T)$ is defined in Proposition~\ref{prop:I3I4} for all $n \in \mathbb{N} \cup \{\infty\}$. By It\^{o}'s isometry, $\mathbb{E}[I^{3,4}_n(f,T)^2]$ for $n \in \mathbb{N}$ equals the left-hand side of \eqref{I4:approx2} with $T_1 = 0$, $T_2 = T$, and $T_3(s) \equiv T$. The same result holds for $n = \infty$ upon incorporating \eqref{Mncovar:derivation:eq}. Since the right-hand side of \eqref{I4:approx2}  is independent of $n$, these identities establish the required uniform bound. 
\end{proof}

In the proof below, Steps 2 and 3 adapt Cho's approach~\cite[Section~2]{Cho:95} to our framework. In Cho's original work, the weak convergence of stochastic integrals driven by martingale measures is studied as processes valued in spaces of tempered distributions. \medskip 

\begin{proof}[Proof of (\ref{I4:conv}) in Proposition~\ref{prop:I3I4}]
We divide the proof into four steps. Fix $T\in (0,\infty)$, let $\delta, \eta > 0$ denote error thresholds for probability bounds of the form $\mathbb{P}(\lvert X-Y\rvert > c\delta) \leq c\eta$, and recall that $E_m$ is defined in \eqref{def:Em} and equipped with the metric induced by $\lVert\cdot\rVert_m$.
 \medskip

\noindent {\bf Step~1.} Our main goal in this step is to show two preliminary properties: 
\begin{itemize}
\item [\rm (1)] For any $n \in \mathbb{N} \cup \{\infty\}$ and $0 < \kappa < T$, the mapping $s' \mapsto H_{n,\kappa}(\cdot, s')$ belongs to $C_{E_m}[0, \infty)$ for every $m \in \mathbb{Z}_+$, where
\[
H_{n,\kappa}(y',s')=H_{n,\kappa}(y',s'\wedge (T-\kappa)), \quad s'\in \R_+, 
\]
is defined by the following equation for $0\leq s'\leq T-\kappa$:
\begin{align}
H_{n,\kappa}(y',s')&\,\defeq\, \int_{t=s'+\kappa}^T \int_x f(x,t)P_{t-s'}(x,y')\int_z P_{t-s'}(x,z)X_n(\d z,s')\d x\d t\notag\\
&=\int_{z}\biggl(\int_{t=0}^{T-s'-\kappa} \int_x f(x,t+s'+\kappa)P_{t+\kappa}(x,y') P_{t+\kappa}(x,z)\d x\d t\biggr)\notag\\
&\quad \times X_n(\d z,s').
\label{def:Hnprocess}
\end{align}
\item [\rm (2)] There exists some $0<\kappa_0<T/3$ such that 
\begin{align}
&\sup_{n\in \Bbb N\cup \{\infty\}}\P\biggl(\biggl|I^{3,4}_n(f,T)
-\int_{s'=0}^{T-\kappa_0}\int_{y'}H_{n,\kappa_0}(y',s')\M_n(\d y',\d s')
\biggr|>\delta\biggr)
\leq \eta.\label{I4:probg0}
\end{align}
\end{itemize}
We establish properties (1) and (2) in Steps~1-1 and 1-2 below. Step~1-3 provides a setup to connect Step~1-2 with Steps~2 and 3.\medskip 

\noindent {\bf Step 1-1.} We prove property (1). The continuity of $y' \mapsto H_{n,\kappa}(y',s')$ for each fixed $0 \leq s' \leq T-\kappa$ follows from Proposition~\ref{prop:XM} (2$\cc$), the condition $f(x, t) \in \C_{pd}^0(\R^2 \times [0, T])$, and the dominated convergence theorem. The rest of property (1) will be established once we show that 
\begin{align}
\lim_{0\leq s_1'\leq s_2'\leq T-\kappa\atop s_2'-s_1'\searrow 0}\sup_{y'}\lvert H_{n,\kappa}(y',s_1')-H_{n,\kappa}(y',s_2')\rvert(1+\lvert y' \rvert^m)=0.\label{I4:probg0-1}
\end{align}

We begin by establishing three bounds. By using \eqref{def:Hnprocess}, it is immediate that for any integer $m \geq 10$, 
\begin{align}
&\sup_{y'}\lvert H_{n,\kappa}(y',s_1')-H_{n,\kappa}(y',s_2')\rvert(1+\lvert y' \rvert^m)\notag\\
&\quad \leq \sup_{y'} \int_{z}\biggl( \int_{t=T-s_2'-\kappa}^{T-s'_1-\kappa} \int_x\lvert f(x,t+s_1'+\kappa)\rvert P_{t+\kappa}(x,y') P_{t+\kappa}(x,z)\d x\d t\biggr) X_n(\d z, s_1')(1+\lvert y' \rvert^m)\notag\\
&\quad \quad +\sup_{y'} \int_{z} \int_{t=0}^{T-s'_2-\kappa} \int_x \lvert f(x,t+s_1'+\kappa)-f(x,t+s_2'+\kappa)\rvert P_{t+\kappa}(x,y')P_{t+\kappa}(x,z)\d x\d t\notag\\
&\quad\quad  \times  
X_n(\d z, s_1')(1+\lvert y' \rvert^m)\notag\\
&\quad\quad  +\sup_{y'} \biggl|\int_{z} \int_{t=0}^{T-s'_2-\kappa} \int_x f(x,t+s_2'+\kappa)P_{t+\kappa}(x,y')P_{t+\kappa}(x,z)\d x\d t\notag\\
&\quad \quad \times  [X_n(\d z, s_1')-X_n(\d z, s_2')](1+\lvert y' \rvert^m)\biggr|.\label{Hn:regularity}
\end{align}
 The remaining two bounds are provided by \eqref{Ptym:bdd} and \eqref{XnP:bdd} as follows: \begin{align}
\sup_{y'}P_{t+\kappa}(x,y')(1+\lvert y' \rvert^m)
&\leq C(m)\sup_{y'}P_{t+\kappa}(y')(1+\lvert x\rvert ^m+\lvert y' \rvert^m)\notag\\
&\leq C(m,T,\kappa)(1+\lvert x\rvert ^m),\quad \forall\;0\leq t\leq T, \label{Ptym:bdd}\\
\int_z P_{t+\kappa}(x,z)X_n(\d z,s_1')&\leq \int_{z} \frac{1}{1+\lvert z\rvert ^{m}}\biggl(\sup_{z'}(1+\lvert z'\rvert^{m})P_{t+\kappa}(x,z')\biggr)X_n(\d z,s_1')\notag\\
&\leq C(m,T,\kappa)\int_{z} \frac{1}{1+\lvert z\rvert ^{m}}X_n(\d z,s_1')(1+\lvert x\rvert ^{m}),
\label{XnP:bdd}
\end{align}
where the last inequality follows by \eqref{Ptym:bdd}. For the rest of Step~1-1, we show that the three terms on the right-hand side of \eqref{Hn:regularity} converge uniformly to zero with respect to the mode defined in \eqref{I4:probg0-1}. The proofs proceed in the same order as these terms appear.

The first term on the right-hand side of \eqref{Hn:regularity} admits the bound below, which follows by applying \eqref{Ptym:bdd} and \eqref{XnP:bdd}:
\begin{align}
&\sup_{y'} \int_{z}\biggl( \int_{t=T-s_2'-\kappa}^{T-s'_1-\kappa} \int_x\lvert f(x,t+s_1'+\kappa)\rvert P_{t+\kappa}(x,y') P_{t+\kappa}(x,z)\d x\d t\biggr) X_n(\d z, s_1')(1+\lvert y' \rvert^m)\notag\\
&\quad \leq C(m,T,\kappa)  \int_{t=T-s_2'-\kappa}^{T-s'_1-\kappa} \int_x\lvert f(x,t+s_1'+\kappa)\rvert(1+\lvert x\rvert ^{m})\d x\d t\notag\\
&\quad \quad \times \sup_{0\leq s\leq T}\int_{z} \frac{1}{1+\lvert z\rvert ^{m}}X_n(\d z,s).\label{XnP:bdd-3-1}
\end{align}
Observe that by Corollary~\ref{cor:muinfty}, 
\begin{align}\label{I4-cond}
\sup_{n\in \Bbb N\cup\{\infty\}}\sup_{0\leq s\leq T_0}\int_z\frac{1}{1+\lvert z\rvert ^m} X_{n}(\d z,s)<\infty\quad \forall\;m\geq 10,\;m\in \Bbb N,\;T_0\in(0,\infty).
\end{align}
Hence, the assumption that $f(x, t) \in \C_{pd}^0(\R^2 \times [0, T])$, the right-hand side tends to zero uniformly as $s_2' - s_1' \searrow 0$ for $0 \leq s_1' \leq s_2' \leq T - \kappa$.

The second term on the right-hand side of \eqref{Hn:regularity} can be estimated in a manner similar to the previous display:
\begin{align}
&\sup_{y'} \int_{z} \int_{t=0}^{T-s'_2-\kappa} \int_x \lvert f(x,t+s_1'+\kappa)-f(x,t+s_2'+\kappa)\rvert P_{t+\kappa}(x,y')P_{t+\kappa}(x,z)\d x\d t\notag\\
&\quad\quad  \times  
X_n(\d z, s_1')(1+\lvert y' \rvert^m)\notag\\
&\quad \leq C(m,T,\kappa)  \int_{t=0}^{T-s'_2-\kappa} \int_x\lvert f(x,t+s_1'+\kappa)- f(x,t+s_s'+\kappa)\rvert(1+\lvert x\rvert ^{m})\d x\d t\notag\\
&\quad \quad \times \sup_{0\leq s\leq T}\int_{z} \frac{1}{1+\lvert z\rvert ^{m}}X_n(\d z,s).\label{XnP:bdd-3-2}
\end{align}
Since $(1+\lvert x\rvert ^{m})f(x,t)$ belongs to $\C_{pd}^0(\R^2\times [0,T])$, \eqref{I4-cond} and the assumption that $f(x,t) \in \C_{pd}^0(\R^2 \times [0,T])$ apply again and imply that the right-hand side also converges to zero uniformly as $s_2'-s_1'\searrow 0$ for $0\leq s_1'\leq s_2'\leq T-\kappa$.

To complete the proof, we need to show that the last term in \eqref{Hn:regularity} approaches zero uniformly. We use the following bound for this term, with $M > 0$ to be sent to infinity at the end:
\begin{align}
&\sup_{y'}\biggl|\int_{z} \int_{t=0}^{T-s'_2-\kappa} \int_x f(x,t+s_2'+\kappa)P_{t+\kappa}(x,y')P_{t+\kappa}(x,z)\d x\d t\notag\\
&\quad\times   [X_n(\d z, s_1')-X_n(\d z, s_2')](1+\lvert y' \rvert^m)\biggr|\notag\\
&\quad \leq \int_{t=0}^{T} \int_{\lvert x\rvert \leq M}\frac{C(f,m,T,\kappa)}{1+\lvert x\rvert ^m} \int_z P_{t+\kappa}(x,z)\lvert X_n(z, s_1')-X_n( z, s_2')\rvert\d z\d x\d t\notag\\
&\quad \quad +\int_{t=0}^{T} \int_{\lvert x\rvert > M}\frac{C(f,m,T,\kappa)}{1+\lvert x\rvert ^{2m}} \int_z P_{t+\kappa}(x,z)\lvert X_n(z, s_1')-X_n( z, s_2')\rvert\d z\d x\d t, \label{XnP:bdd3-1}
\end{align}
which follows from \eqref{Ptym:bdd} and the assumption that $f(x, t) \in \C_{pd}^0(\R^2 \times [0, T])$.
The first term on the right-hand side satisfies 
\begin{align}
&\int_{t=0}^{T} \int_{\lvert x\rvert \leq M}\frac{C(f,m,T,\kappa)}{1+\lvert x\rvert ^m} \int_z P_{t+\kappa}(x,z)\lvert X_n(z, s_1')-X_n( z, s_2')\rvert\d z\d x\d t\notag\\
&\quad \leq C(f,m,T,\kappa,M)\sup_{z}\e^{-\lvert z\rvert /[2(T+\kappa)]}\lvert X_n(z, s_1')-X_n( z, s_2')\rvert,\label{XnP:bdd3-2}
\end{align}
which uses the inequality $\e^{-\lvert x-y\rvert^2/a}\leq C(a)\e^{-\lvert x-y\rvert/a}\leq C(a)\e^{\lvert x\rvert /a-\lvert y \rvert/a}$ for all $a>0$ and $x,y\in \R^2$. The other term satisfies 
\begin{align}
&\int_{t=0}^{T} \int_{\lvert x\rvert > M}\frac{C(f,m,T,\kappa)}{1+\lvert x\rvert ^{2m}} \int_z P_{t+\kappa}(x,z)\lvert X_n(z, s_1')-X_n( z, s_2')\rvert\d z\d x\d t\notag\\
&\quad =\int_{t=0}^{T}\int_z  \int_{\lvert x\rvert > M}\frac{C(f,m,T,\kappa)}{1+\lvert x\rvert ^{2m}}  P_{t+\kappa}(x,z)\d x\lvert X_n(z, s_1')-X_n( z, s_2')\rvert\d z\d t\notag\\
&\quad \leq C(f,m,T,\kappa)\int_{t=0}^T\int_z \biggl(\int_{\lvert x\rvert >M}\frac{\d x}{1+\lvert x\rvert ^{2m}}\biggr)^{1/2}\biggl(\int_{x}\frac{P_{t+\kappa}(x,z)^2}{1+\lvert x\rvert ^{2m}}\d x\biggr)^{1/2}\notag\\
&\quad \quad \times [X_n(z, s_1')+X_n( z, s_2')]\d z\d t\notag\\
&\quad\leq  C(f,m,T,\kappa) \biggl(\int_{\lvert x\rvert >M}\frac{\d x}{1+\lvert x\rvert ^{2m}}\biggr)^{1/2}\biggl(\sup_{0\leq s\leq T}\int_z \frac{1}{1+\lvert z\rvert ^m}X_n(\d z,s)\biggr),\label{XnP:bdd3-3}
\end{align}
where the first inequality follows from the Cauchy--Schwarz inequality and the nonnegativity of $X_n$ as stated in Proposition~\ref{prop:XM} (2$\cc$), while the last inequality is a consequence of \eqref{ineq:heat2-2}. By invoking Proposition~\ref{prop:XM} (2$\cc$) and \eqref{I4-cond} and applying \eqref{XnP:bdd3-2} and \eqref{XnP:bdd3-3} to \eqref{XnP:bdd3-1}, we get
\begin{align*}
&\limsup_{M\to\infty}\limsup_{0\leq s_1'\leq s_2'\leq T-\kappa\atop s_2'-s_1'\searrow 0}\sup_{y'}\biggl|\int_{z} \int_{t=0}^{T-s'_2-\kappa} \int_x f(x,t+s_2'+\kappa)P_{t+\kappa}(x,y')P_{t+\kappa}(x,z)\d x\d t\notag\\
&\quad\times   [X_n(\d z, s_1')-X_n(\d z, s_2')](1+\lvert y' \rvert^m)\biggr|=0.
\end{align*}
Therefore, the third term in \eqref{Hn:regularity} vanishes uniformly as required. This completes the proof of property (1). \medskip 

\noindent {\bf Step 1-2.} We verify property (2) by establishing that \eqref{I4:probg0} is satisfied for some $0<\kappa_0<T/3$. By Markov's inequality, it suffices to show
\begin{align*}
&\lim_{\kappa\searrow 0}\sup_{n\in \Bbb N\cup \{\infty\}}\E\biggl[\biggl(I^{3,4}_n(f,T)-\int_{s'=0}^{T-\kappa}\int_{y'} \\
&\quad \times \biggl( \int_{t=s'+\kappa}^T \int_x f(x,t)P_{t-s'}(x,y')\int_z P_{t-s'}(x,z)X_n(\d z,s')\d x\d t\biggr)\M_n(\d y',\d s')
\biggr)^2\biggr]=0.
\end{align*}
To justify this, observe that $\int_{s'=0}^T\int_{t=s'}^{T}-\int_{s'=0}^{T-\kappa}\int_{t=s'+\kappa}^T=\int_{s'=T-\kappa}^T \int_{t=s'}^T+\int_{s'=0}^{T-\kappa}\int_{t=s'}^{s'+\kappa}$. Hence, applying It\^{o}'s isometry, Lemma~\ref{lem:I4alt} with $(T_1,T_2,T_3(s'))=(T-\kappa,T, T)$ and $(T_1,T_2,T_3(s'))=(0,T-\kappa,s'+\kappa)$, and \eqref{Mncovar:derivation:eq} for $n=\infty$
 yield the result. Thus, \eqref{I4:probg0} is verified.
\medskip 

\noindent {\bf Step~1-3.}
With the choice of $\kappa_0$ ensuring \eqref{I4:probg0}, we demonstrate in Steps~2 and 3 below
 that there exists an integer $N_0\geq 1$ such that for all $n\geq N_0$, $H_{n,\kappa_0}(y',s')$ satisfies 
 \begin{align}
&\P\biggl(\biggl|\int_{s'=0}^{T-\kappa_0}\int_{y'} H_{n,\kappa_0}(y',s')\M_n(\d y',\d s')
-\int_{s'=0}^{T-\kappa_0}\int_{y'}H_{\infty,\kappa_0}(y',s')\M_\infty(\d y',\d s')\biggr|>3\delta\biggr)\notag\\
&\quad \leq 3\eta,\label{I4:probg2}
\end{align}
Since we focus exclusively on $H_{n,\kappa_0}(y',s')$ for the rest of this proof, we will simply write
\[
H_n(y',s')\,\defeq\, H_{n,\kappa_0}(y',s'),\quad n\in\Bbb N\cup\{\infty\}.
\]

\noindent {\bf Step~2. } In this step, we show that for small enough $\theta_0>0$, 
\begin{align}
&\sup_{n\in \Bbb N\cup\{\infty\}}\P\left(\left|\int_{s'=0}^{T-\kappa_0}\int_{y'}H_n(y',s')M_n(\d y',\d s')-\int_{s'=0}^{T-\kappa_0}\int_{y'}H^{\theta_0}_n(y',s'-)M_n(\d y',\d s')\right|> \delta\right)\notag\\
&\quad \leq \eta.\label{I4:probg3}
\end{align}
In \eqref{I4:probg3}, $H_n^\theta(y', s')$ denotes a random, function-valued step function, defined as follows. We work with $E = E_m$, where $(E_m, \lvert \cdot\rvert_m)$ is given below Lemma~\ref{lem:KP}, and $m \geq 10$ is fixed.
\begin{align}
H_n^\theta(y',s')&\,\defeq\,I^\theta(s''\mapsto H_n(\cdot,s''))(s')(y')\notag\\
&=\sum_{k=0}^\infty H_n(y',\tau^\theta_{k,n})\1_{[\tau^\theta_{k,n},\tau^{\theta}_{k+1,n})}(s'),\quad n\in \Bbb N\cup\{\infty\}.\label{def:Hntheta}
\end{align}
Note that property (1) in Step~1 is used to justify applying Lemma~\ref{lem:KP}. The notation
$I^\theta(s''\mapsto H_n(\cdot,s''))(s')(y')$ refers to the value at $y'$ of the $E_m$-valued function $I^\theta(s''\mapsto H_n(\cdot,s''))(s')$. Additionally, 
\begin{align}\label{def:tauthetakn}
\tau^\theta_{k,n}\,\defeq\, \tau_k^\theta(s''\mapsto H_n(\cdot,s'')) 
\end{align}
is independent of $y'$ since its definition only uses the norm $\lVert\cdot\rVert_m$.

To prove \eqref{I4:probg3}, we consider the following inequalities for any $\theta>0$:
\begin{align*}
&\P\left(\left|\int_{s'=0}^{T-\kappa_0}\int_{y'}H_n(y',s')M_n(\d y',\d s')-\int_{s'=0}^{T-\kappa_0}\int_{y'}H^\theta_n(y',s'-)M_n(\d y',\d s')\right|> \delta\right)\\
&\quad \leq\frac{1}{\delta^2}\E\left[\int_{s=0}^{T-\kappa_0}\int_{x}\int_{x'}[H_n(x',s)-H_n^\theta(x',s-)][H_n(x,s)-H_n^{\theta}(x,s-)]\mu_n(\d x',\d x,\d s)\right] \\
&\quad \leq \frac{1}{\delta^2}\E\left[\int_{s=0}^{T-\kappa_0}\int_{x}\int_{x'}\frac{\theta^2}{(1+\lvert x'\rvert^m)(1+\lvert x\rvert ^m)}
\mu_n(\d x',\d x,\d s)\right],
\end{align*}
where the last inequality holds by using Lemma~\ref{lem:KP} (2$\cc$) and the continuity of $s'\mapsto H_n(\cdot ,s')$ in $C_{E_m}[0,\infty)$. Furthermore, Theorem~\ref{thm:mombdd} (3$\cc$) for $n \in \mathbb{N}$, together with Theorem~\ref{thm:main1} (2$\cc$) and \eqref{Mncovar:derivation:eq} for $n = \infty$, yields 
\[
\sup_{n\in \Bbb N\cup\{\infty\}}\E\left[\int_{s=0}^{T-\kappa_0}\int_{x}\int_{x'}\frac{1}{(1+\lvert x'\rvert^m)(1+\lvert x\rvert ^m)}
\mu_n(\d x',\d x,\d s)\right]<\infty.
\]
Hence, we can choose $\theta_0>0$ small enough to attain \eqref{I4:probg3}. \medskip 

\noindent {\bf Step~3.}
This step gives the last approximation. Our goal is to show that for a.s.-$\omega$,
\begin{align}\label{I4:probg40}
\int_{s'=0}^{T-\kappa_0}\int_{y'}H^{\theta_0}_n(y',s'-)M_n(\d y',\d s')\xrightarrow[n\to\infty]{} \int_{s'=0}^{T-\kappa_0}\int_{y'}H^{\theta_0}_\infty(y',s'-)M_\infty(\d y',\d s').
\end{align}
Consequently, we can find an integer $N(\delta,\eta, \theta_0,\kappa_0)\geq 1$ large enough such that  
\begin{align}
&\P\left(\left|\int_{s'=0}^{T-\kappa_0}\int_{y'}H^{\theta_0}_n(y',s'-)M_n(\d y',\d s')- \int_{s'=0}^{T-\kappa_0}\int_{y'}H^{\theta_0}_\infty(y',s'-)M_\infty(\d y',\d s')\right|> \delta\right)\notag\\
&\quad \leq \eta,\quad \forall\;n\geq N(\delta,\eta, \theta_0,\kappa_0).\label{I4:probg4}
\end{align}

\noindent {\bf Step 3-1.} To prove \eqref{I4:probg40}, we need two preliminary convergence results. The first one iterates the following limit from Corollary~\ref{cor:muinfty}:
\begin{align}
\P\Bigl(\,\lim_{n\to\infty}X^\sharp_n=X^\sharp_\infty\mbox{ in $C_{\mathcal M_f(\R^2)}[0,\infty)$}\Bigr)=1;\label{I4-00} 
\end{align}
see Corollary~\ref{cor:muinfty}.  Additionally, we claim that
\begin{align}\label{I4-0} 
\P\Bigl(\,\lim_{n\to\infty}H_n=H_\infty\mbox{ in $C_{E_m}[0,\infty)$}\Bigr)=1,\quad \forall\;m\geq 10,\;m\in \Bbb N,
\end{align}
which is equivalent to  
the condition that
\begin{align*}
&\forall\;\mbox{ subsequence } \{n_k\},\;\forall\;\{s_{n_k}\}\subset [0,T-\kappa_0],\\
&\lim_{k\to\infty}\sup_{y'}\lvert H_{n_k}(y',s_{n_k})-H_\infty(y',s_{n_k})\rvert(1+\lvert y' \rvert^m)=0.
\end{align*}
To confirm this, property (1) from Step~1-1 guarantees that the following condition is sufficient:\begin{align}
\begin{split}
&\forall\;\mbox{ subsequence } \{n_k\},\;\forall\;s_{n_k}\subset [0,T-\kappa_0] \mbox{ with }\lim_{k\to\infty}s_{n_k}=s_\infty,\\
&\lim_{k\to\infty}\sup_{y'}\lvert H_{n_k}(y',s_{n_k})-H_\infty(y',s_{\infty})\rvert(1+\lvert y' \rvert^m)=0.\label{I4-1}
\end{split}
\end{align}

To prove \eqref{I4-1}, we use \eqref{def:Hnprocess} and \eqref{I4-cond}. Hence, an inspection of the proof of \eqref{Hn:regularity} reveals that
\begin{align*}
&\lim_{k\to\infty}\sup_{y'}\biggl|
\int_{z}\biggl(\int_{t=0}^{T-s_{n_k}-\kappa_0} \int_x f(x,t+s_{n_k}+\kappa_0)P_{t+\kappa_0}(x,y') P_{t+\kappa_0}(x,z)\d x\d t\biggr) X_{n_k}(\d z,s_{n_k})\\
&\quad -\int_{z}\biggl(\int_{t=0}^{T-s_\infty-\kappa_0} \int_x f(x,t+s_\infty+\kappa_0)P_{t+\kappa_0}(x,y') P_{t+\kappa_0}(x,z)\d x\d t\biggr) X_{n_k}(\d z,s_{n_k})
\biggr|(1+\lvert y' \rvert^m)\\
&\quad =0.
\end{align*}
(More specifically, only the analogues of the first two terms on the right-hand side of \eqref{Hn:regularity} are relevant, which means that only the analogues of \eqref{XnP:bdd-3-1} and \eqref{XnP:bdd-3-2} arise.) Consequently, \eqref{I4-1} holds as soon as we can show that 
\begin{align*}
&\lim_{k\to\infty}\sup_{y'}\biggl|
\int_{z}\biggl(\int_{t=0}^{T-s_\infty-\kappa_0} \int_x f(x,t+s_\infty+\kappa_0)P_{t+\kappa_0}(x,y') P_{t+\kappa_0}(x,z)\d x\d t\biggr) X_{n_k}(\d z,s_{n_k})\\
&\quad -\int_{z}\biggl(\int_{t=0}^{T-s_\infty-\kappa_0} \int_x f(x,t+s_\infty+\kappa_0)P_{t+\kappa_0}(x,y') P_{t+\kappa_0}(x,z)\d x\d t\biggr) X_\infty(\d z,s_{\infty})
\biggr|(1+\lvert y' \rvert^m)\\
&\quad =0.
\end{align*}
This, in turn, can be verified by showing that for any further subsequence, still denoted by $\{n_k\}$,
we have the following limit whenever $\lvert y_{n_k}'\rvert\to\infty$ or $y_{n_k}'\to y'_\infty$ for some $y'_\infty\in \R^2$:
\begin{align}
&\lim_{k\to\infty}\biggl|
\int_{z}\biggl(\int_{t=0}^{T-s_\infty-\kappa_0} \int_x f(x,t+s_\infty+\kappa_0)P_{t+\kappa_0}(x,y_{n_k}') P_{t+\kappa_0}(x,z)\d x\d t\biggr) X_{n_k}(\d z,s_{n_k})\notag\\
&\quad -\int_{z}\biggl(\int_{t=0}^{T-s_\infty-\kappa_0} \int_x f(x,t+s_\infty+\kappa_0)P_{t+\kappa_0}(x,y_{n_k}') P_{t+\kappa_0}(x,z)\d x\d t\biggr) X_\infty(\d z,s_{\infty}) \biggr|\notag\\
&\quad \times(1+\lvert y_{n_k}'\rvert^m)=0.\label{I4-2}
\end{align}

First, we show that \eqref{I4-2} holds whenever $\lvert y_{n_k}'\rvert\to\infty$. Let $n'\in \{n_k,\infty\}$, and consider $M>0$, which will ultimately be passed to infinity. Then
\begin{align*}
&\biggl|\int_{z}\biggl(\int_{t=0}^{T-s_\infty-\kappa_0} \int_x f(x,t+s_\infty+\kappa_0)P_{t+\kappa_0}(x,y_{n_k}') P_{t+\kappa_0}(x,z)\d x\d t\biggr) X_{n'}(\d z,s_{n'})\biggr|(1+\lvert y_{n_k}'\rvert^m)\\
&\quad \leq C(f,m,T)\int_{z}\int_{t=0}^T \int_{\lvert x\rvert \leq M}\frac{1}{1+\lvert x\rvert ^m}P_{t+\kappa_0}(x,y_{n_k}') P_{t+\kappa_0}(x,z)\d x\d tX_{n'}(\d z,s_{n'})(1+\lvert y_{n_k}'\rvert^m)\\
&\quad \quad +C(f,m,T)\int_{z}\int_{t=0}^T \int_{\lvert x\rvert > M}\frac{1}{1+\lvert x\rvert ^{3m}}P_{t+\kappa_0}(x,y_{n_k}') P_{t+\kappa_0}(x,z)\d x\d tX_{n'}(\d z,s_{n'})(1+\lvert y_{n_k}'\rvert^m)\\
&\quad \leq  C(f,m,T,\kappa_0,M)\e^{-\lvert y'_{n_k}\rvert /[2(T+\kappa_0)]} \sup_{0\leq s\leq T} \int_z \frac{1}{1+\lvert z\rvert ^m}X_{n'}(\d z,s)(1+\lvert y_{n_k}'\rvert^m)\\
&\quad \quad +C(f,m,T,\kappa_0)\biggl(\int_{\lvert x\rvert \geq M}\frac{1}{1+\lvert x\rvert ^{2m}}\d x\biggr)^{1/2} \sup_{0\leq s\leq T}\int_z \frac{1}{1+\lvert z\rvert ^m}X_{n'}(\d z,s),
\end{align*}
where the first term arises from reasoning analogous to that used for \eqref{XnP:bdd3-2}, while the second term follows from \eqref{Ptym:bdd} and a similar argument as in \eqref{XnP:bdd3-3}. Applying \eqref{I4-cond} to the right-hand side of the last inequality, we conclude that \eqref{I4-2} holds by taking $k\to\infty$ first, followed by $M\to\infty$. 

Next, we show that \eqref{I4-2} holds whenever $y_{n_k}'\to y'_\infty\in \R^2$ . Since $\sup_{k}\lvert y'_{n_k}\rvert <\infty$, it suffices to prove the following three limits, where $n'$ belongs to $ \{n_k,\infty\}$ for the first limit: 
\begin{align}
&\lim_{k\to\infty}\int_{z}\biggl(\int_{t=0}^{T-s_\infty-\kappa_0} \int_x \lvert f(x,t+s_\infty+\kappa_0)\rvert \lvert P_{t+\kappa_0}(x,y_{n_k}')-
P_{t+\kappa_0}(x,y_{\infty}')\rvert 
 P_{t+\kappa_0}(x,z)\d x\d t\biggr)\notag\\
 &\quad \times X_{n'}(\d z,s_{n'})=0,\label{I4-3}\\
&\lim_{k\to\infty}
\int_{z}\biggl(\int_{t=0}^{T-s_\infty-\kappa_0} \int_x f(x,t+s_\infty+\kappa_0)P_{t+\kappa_0}(x,y_{\infty}') P_{t+\kappa_0}(x,z)\d x\d t\biggr) X_{n_k}(\d z,s_{n_k})\notag\\
&\quad -\int_{z}\biggl(\int_{t=0}^{T-s_\infty-\kappa_0} \int_x f(x,t+s_\infty+\kappa_0)P_{t+\kappa_0}(x,y_{\infty}') P_{t+\kappa_0}(x,z)\d x\d t\biggr)X_\infty(\d z,s_{\infty})\notag\\
&\quad  =0.\label{I4-4}
\end{align}
To get \eqref{I4-3}, we use \eqref{gauss:vep} with the chocie of $\delta\equiv \kappa_0$, $\vep\equiv \lvert y_{n_k}'-y_{\infty}'\rvert\kappa_0^{1/2}/(2M')$, and $M\equiv (2M')/\kappa_0^{1/2}$, 
where $M'>0$ is such that $\lvert y_{n_k}'\rvert\leq M'$ for all $k$. Hence,
\begin{align*}
&\int_{z}\biggl(\int_{t=0}^{T-s_\infty-\kappa_0} \int_x \lvert f(x,t+s_\infty+\kappa_0)\rvert \lvert P_{t+\kappa_0}(x,y_{n_k}')-
P_{t+\kappa_0}(x,y_{\infty}')\rvert
 P_{t+\kappa_0}(x,z)\d x\d t\biggr) X_{n'}(\d z,s_{n'})\\
 &\quad \leq C(\kappa_0, M')\lvert y_{n_k}'-y_{\infty}'\rvert\int_{z}\biggl(\int_{t=0}^{T-s_\infty-\kappa_0} \int_x \lvert f(x,t+s_\infty+\kappa_0) \rvert P_{2(t+\kappa_0)}(x,y_{\infty}')
 P_{t+\kappa_0}(x,z)\d x\d t\biggr)\\
 &\quad \quad \times X_{n'}(\d z,s_{n'})\\
 &\quad \leq C(f,m,T,\kappa_0,M')\lvert y_{n_k}'-y_{\infty}'\rvert \sup_{0\leq s\leq T}\int_{z}\frac{1}{1+\lvert z\rvert ^m}X_{n'}(\d z,s)\xrightarrow[k\to\infty]{}0,
\end{align*}
where the last inequality follows by using the condition $f(x,t)\in \C^0_{pd}(\R^2\times [0,T])$, as well as \eqref{ineq:heat2-1}, and
the limit follows by using \eqref{I4-cond}. Additionally, \eqref{I4-4} follows directly from Corollary~\ref{cor:muinfty}, since the integrands with respect to the measures $X_{n_k}(\d z,s_{n_k})$ and $X_{\infty}(\d z,s_{\infty})$ are identical and independent of $k$, and the function of $z$ defining the integrands is in $\C^0_{pd}(\R^2\times [0,T])$ due to the condition $f(x,t)\in \C^0_{pd}(\R^2\times [0,T])$. This establishes \eqref{I4-2}, and therefore \eqref{I4-0}. \medskip 

\noindent {\bf Step~3-2.} We present further preliminary results, listed as (A)--(D) below, to prepare the proof of  \eqref{I4:probg40} in this step. \medskip 

\noindent {\bf (A)} We can choose an event $\Omega_1$ of probability one on which \eqref{XM:relation}, \eqref{I4-cond} and the following two limits hold:
\begin{align}\label{I4-5}
\begin{split}
\lim_{n\to\infty}X^\sharp_n&= X^\sharp_\infty\mbox{ in $C_{\mathcal M_f(\R^2)}[0,\infty)$},\\
\lim_{n\to\infty}H^{\theta_0}_n&= H^{\theta_0}_\infty\mbox{ in }D_{E_m}[0,\infty),
\end{split}
\end{align}
where $H_n^\theta=H_n^\theta(y',s')$ is the random step function defined in \eqref{def:Hntheta}. 

The existence of $\Omega_1$ follows from the two preliminary convergence results from Step~3-1. More precisely, the second limit in \eqref{I4-5} is ensured by \eqref{I4-0} and Lemma~\ref{lem:KP} (3$\cc$). {\bf Fix $\omega\in \Omega_1$ and, for notational simplicity, omit explicit reference to $\omega$ in the rest of Step~3.} \medskip 

\noindent {\bf (B).} There exist strictly increasing, Lipschitz-continuous functions $\lambda_n:\R_+\to \R_+$, $n\in \Bbb N$, such that all of the following limits hold: 
\begin{gather}
\sup_{0\leq t\leq T_0}\lvert \lambda_n(t)- t\rvert \xrightarrow[n\to\infty]{}0,\quad \forall\; T_0\in (0,\infty),\label{Ln:conv}\\
\sup_{0\leq t\leq T_0}\lVert H^{\theta_0}_n(\cdot,\lambda_n(t))- H^{\theta_0}_\infty(\cdot,t) \rVert_m\xrightarrow[n\to\infty]{} 0,\quad \forall\; T_0\in (0,\infty).\label{Hn:conv}
\end{gather}
This follows from \eqref{I4-5} and a standard property of the Skorokhod space \cite[5.3 Proposition, p.119]{EK:MP}.\medskip 

\noindent {\bf (C).} We have
\begin{gather}
X_n^\sharp(\cdot,\lambda_n(t))\xrightarrow[n\to\infty]{} X^\sharp_\infty(\cdot,t)\quad  \mbox{ in }C_{\mathcal M_f(\R^2)}[0,\infty),\label{Xn:conv}\\
\label{IX_n:conv}
\int_0^{\lambda_n(t)}X^\sharp_n(\cdot,s)\d s
\xrightarrow[n\to\infty]{}\int_0^{t}X^\sharp_\infty(\cdot,s)\d s\quad\mbox{ in } C_{\mathcal M_f(\R^2)}[0,\infty).
\end{gather}

To establish \eqref{Xn:conv}, we apply \eqref{Ln:conv} to the first limit in \eqref{I4-5}. Furthermore, the limit in \eqref{IX_n:conv} follows from \eqref{Ln:conv} and \eqref{Xn:conv}.\medskip 

\noindent{\bf (D).} We have the following two properties of the random step functions $H_n^{\theta_0}$, defined in \eqref{def:Hntheta},
 for all $n\in \Bbb N\cup \{\infty\}$:
\begin{itemize} 
\item []{\bf \bf (D-1).} For each $n \in \mathbb{N} \cup \{\infty\}$, let $\{t_k^n\} = \{t_0^n < t_1^n \leq t_2^n\leq \cdots\}$, where $t_0^n=0$, and each $t_k^n$ for $k\geq 1$ denotes either a unique discontinuity of $H^{\theta_0}_n$ or equals $\infty$. Then $\{t^n_k\}= \{\tau^{\theta_0}_{k,n}\}$, where $\tau^{\theta_0}_{k,n}$ is defined in \eqref{def:tauthetakn}. 
\item []{\bf (D-2).} There exists an integer $K\geq 1$, independent of $n$, such that the number of discontinuities of $H^{\theta_0}_n$ in $[0,T-\kappa_0]$ is bounded by $K$, for any $n\in \Bbb N\cup \{\infty\}$.
\end{itemize}
To justify the identity $\{t^n_k\} = \{\tau^{\theta_0}_{k,n}\}$, observe that for $I^\theta(\mathfrak{z})$ (as defined in \eqref{def:Ithetaz}), a genuine jump necessarily occurs at every finite $\tau^\theta_k$, as such a jump must have size at least $\theta/2$ with respect to the associated metric. Additionally, this observation, combined with \eqref{Hn:conv} and the continuity $H_\infty \in C_{E_m}[0, \infty)$, establishes the second property listed above. \medskip 

\noindent {\bf Step~3-3.}
We show the required convergence in \eqref{I4:probg40} with respect to the $\omega$ under consideration in this step. Set
\begin{align}
\begin{split}\label{def:sknrn}
s_k^n\,\defeq\,\lambda_n^{-1}(t_k^n),\quad
r_n\,\defeq\,\lambda_n^{-1}(T-\kappa_0).
\end{split} 
\end{align}
Then the pre-limit integral in \eqref{I4:probg40} satisfies: 
 \begin{align}
&\int_{s'=0}^{T-\kappa_0}\int_{y'} H^{\theta_0}_n(y',s'-)M_n(\d y',\d s')\notag\\
&\quad =\int_{s'=0}^{\lambda_n(r_n)}\int_{y'} H^{\theta_0}_n(y',s'-)M_n(\d y',\d s')\notag\\
&\quad =\sum_{k=0}^{K}\int_{s'=\lambda_n(r_n\wedge s_k^n)}^{\lambda_n(r_n\wedge s_{k+1}^n)}\int_{y'} H_n(y',\lambda_n(r_n\wedge s_k^n))M_n(\d y',\d s')\notag\\
&\quad =\sum_{k=0}^{K}\int_{s'=\lambda_n(r_n\wedge s_k^n)}^{\lambda_n(r_n\wedge s_{k+1}^n)}\int_{y'} [H_n(y',\lambda_n(r_n\wedge s_k^n))-H_\infty(y',r_n\wedge s_k^n)]M_n(\d y',\d s')\notag\\
&\quad\quad +\sum_{k=0}^{K}\int_{s'=\lambda_n(r_n\wedge s_k^n)}^{\lambda_n(r_n\wedge s_{k+1}^n)}\int_{y'} H_\infty(y',r_n\wedge s_k^n)(M_n-M_\infty)(\d y',\d s')\notag\\
&\quad\quad  +\sum_{k=0}^{K}\biggl(\int_{s'=\lambda_n(r_n\wedge s_k^n)}^{\lambda_n(r_n\wedge s_{k+1}^n)}-\int_{s'=r_n\wedge s_k^n}^{r_n\wedge s_{k+1}^n}\biggr)\int_{y'} H_\infty(y',r_n\wedge s_k^n)M_\infty(\d y',\d s')\notag\\
&\quad\quad  +\sum_{k=0}^{K}\int_{s'=r_n\wedge s_k^n}^{r_n\wedge s_{k+1}^n}\int_{y'} H_\infty(y',r_n\wedge s_k^n)M_\infty(\d y',\d s'),\label{weakintegral:final}
\end{align}
where the second equality follows from the identity $\{t^n_k\}= \{\tau^{\theta_0}_{k,n}\}$ in (D-1) of Step~3-2. The following list evaluates the limits of the four sums on the right-hand side as $n\to\infty$, demonstrating that the first three sums converge to zero and the last sum approaches the limit stated in \eqref{I4:probg40}:
\begin{itemize}
\item For the first sum, we observe that 
\begin{align}\label{def:Psink}
\Psi^n_k(y')\,\defeq\,H_n(y',\lambda_n(r_n\wedge s_k^n))-H_\infty(y',r_n\wedge s_k^n)\in \C_{pd}^2(\R^2)
\end{align}
(though it is a random $\C^2_{pd}(\R^2)$-function). The $\C_{pd}^2(\R^2)$-property follows from the integral representation in \eqref{def:Hnprocess}, which establishes the $\C^2$-regularity. The decay property is then obtained by applying \eqref{XnP:bdd}, \eqref{I4-cond}, the condition $f(x,t)\in \C_{pd}^0(\R^2\times [0,T])$, and \eqref{ineq:heat2-2}. 
 Hence, \eqref{XM:relation} applies
\begin{align}
&\int_{s'=\lambda_n(r_n\wedge s_k^n)}^{\lambda_n(r_n\wedge s_{k+1}^n)}\int_{y'} [H_n(y',\lambda_n(r_n\wedge s_k^n))-H_\infty(y',r_n\wedge s_k^n)]M_n(\d y',\d s')\notag\\
&\quad=X_n(\Psi^n_k,\lambda_n(r_n\wedge s_{k+1}^n))-X_n(\Psi^n_k,\lambda_n(r_n\wedge s_k^n))\notag \\
&\quad \quad+\int_{s'=\lambda_n(r_n\wedge s_k^n)}^{\lambda_n(r_n\wedge s_{k+1}^n)}X_n\biggl(\frac{\Delta \Psi^n_k}{2},s\biggr)\d s.\label{I4:3-1}
\end{align}
We demonstrate that, for each $0\leq k\leq K$, the right-hand side of the preceding equality converges to zero as $n\to\infty$. Recall that the spatial partial derivatives of the two-dimensional heat kernel yield the Hermite polynomials. For example, the primary case of interest here is
\[
\partial_{y_j}^2P_t(x,y)=\frac{(x_j-y_j)^2-t}{t^2}P_t(x,y),
\]
provided that $x=(x_1,x_2)$ and $y=(y_1,y_2)$. Thus, employing the integral representation in \eqref{def:Hnprocess}, together with the condition $f(x,t) \in \C^0_{pd}(\R^2 \times [0,T])$, property (1) of Step 1 for $n=\infty$, and \eqref{I4-0}, we obtain that
\[
\lim_{n\to\infty}\sup_{0\leq t\leq T}\lVert\partial^\alpha_{y'}H_n(y',\lambda_n(t))-\partial^\alpha_{y'} H_\infty(y',t)\rVert_m=0
\] 
for any multi-index $\alpha$ and integer $m\geq 10$.
In particular, for all $k\in \{0,\cdots,K\}$, 
\[
\lVert\partial^\alpha_{y'} \Psi^n_k(y')\rVert_m=\lVert\partial^\alpha_{y'}H_n(y',\lambda_n(r_n\wedge s_k^n))-\partial^\alpha_{y'} H_\infty(y',r_n\wedge s_k^n)\rVert_m\xrightarrow[n\to\infty]{} 0.
\] 
Therefore, by \eqref{I4-cond}, each term on the right-hand side of \eqref{I4:3-1} converges to zero as $n\to\infty$ for every $0 \leq k \leq K$. Since $K$ is independent of $n$, the first sum on the right-hand side of \eqref{weakintegral:final} also converges to zero as $n \to \infty$.

\item For the second sum, note that by \eqref{Ln:conv} and \eqref{def:sknrn}, for each $0 \leq k \leq K$, the sequence $\{r_n \wedge s^n_k\}_{n \in \Bbb N}$ is bounded. Consequently, every subsequence admits a further subsequence, which we denote by $
\{r_{n_\ell}\wedge s^{n_\ell}_k\}_{\ell\in \Bbb N}$, such that 
\begin{align}\label{cond:snell}
s^\infty_k\,\defeq\,\lim_{\ell\to\infty}r_{n_\ell}\wedge s^{n_\ell}_k\quad \mbox{ exists in }[0,T-\kappa_0].
\end{align}
Set
\[
\Phi^n_k(y')\,\defeq\,H_\infty(y',r_n\wedge s_k^n),\quad \Phi^\infty_k(y')\,\defeq\,
H_\infty(y', s_k^{\infty}),
\]
all of which belong to $\C^2_{pd}(\R^2)$, 

We now present the main argument, demonstrating that each term in the second sum on the right-hand side of \eqref{weakintegral:final} converges to zero as $n \to \infty$. By \eqref{XM:relation},
\begin{align}
&\int_{s'=\lambda_n(r_n\wedge s_k^n)}^{\lambda_n(r_n\wedge s_{k+1}^n)}\int_{y'} H_\infty(y',r_n\wedge s_k^n)(M_n-M_\infty)(\d y',\d s')\notag\\
&\quad =(X_n-X_\infty)(\Phi^n_k-\Phi^\infty_k,\lambda_n(r_n\wedge s_{k+1}^n))-(X_n-X_\infty)(\Phi^n_k-\Phi^\infty_k, \lambda_n(r_n\wedge s_k^n))\notag\\
&\quad \quad-\int_{s'=\lambda_n(r_n\wedge s_k^n)}^{\lambda_n(r_n\wedge s_{k+1}^n)}(X_n-X_\infty)\biggl(\frac{\Delta (\Phi^n_k-\Phi^\infty_k)}{2},s\biggr)\d s\notag\\
&\quad \quad +(X_n-X_\infty)(\Phi^\infty_k,\lambda_n(r_n\wedge s_{k+1}^n))-(X_n-X_\infty)(\Phi^\infty_k, \lambda_n(r_n\wedge s_k^n))\notag\\
&\quad \quad-\int_{s'=\lambda_n(r_n\wedge s_k^n)}^{\lambda_n(r_n\wedge s_{k+1}^n)}(X_n-X_\infty)\biggl(\frac{\Delta \Phi^\infty_k}{2},s\biggr)\d s.
\label{I4:3-2}
\end{align}
Similarly to the argument for the first sum, but without applying \eqref{I4-0},
it follows that for all $k \in \{0, \ldots, K\}$, 
\[
\lVert\partial^\alpha_{y'} \Phi^{n_\ell}_k(y')-\partial_{y'}\Phi^\infty_k(y')\rVert_m\xrightarrow[n\to\infty]{}0
\]
for any multi-index $\alpha$ and any integer $m \geq 10$. Thus, for each $0 \leq k \leq K$, the first three terms on the right-hand side of \eqref{I4:3-1} converge to zero by \eqref{I4-cond}, while the last three terms do so by the first limit in \eqref{I4-5}. Because $\{r_{n_\ell} \wedge s^{n_\ell}_k\}_{\ell \in \Bbb N}$ is an arbitrary subsequence for which \eqref{cond:snell} holds, it follows that the second sum on the right-hand side of \eqref{weakintegral:final} also converges to zero as $n \to \infty$.
\item For the third sum on the right-hand side of  \eqref{weakintegral:final}, we use \eqref{XM:relation} to rewrite it as follows:
\begin{align*}
&\biggl(\int_{s'=\lambda_n(r_n\wedge s_k^n)}^{\lambda_n(r_n\wedge s_{k+1}^n)}-\int_{s'=r_n\wedge s_k^n}^{r_n\wedge s_{k+1}^n}\biggr)\int_{y'} H_\infty(y',r_n\wedge s_k^n)M_\infty(\d y',\d s')\\
&\quad =X_\infty(\Phi^n_k,\lambda_n(r_n\wedge s_{k+1}^n) )-X_\infty(\Phi^n_k,\lambda_n(r_n\wedge s_k^n))\\
&\quad \quad -X_\infty(\Phi^n_k,r_n\wedge s_{k+1}^n )+X_\infty(\Phi^n_k,r_n\wedge s_k^n)
\\
&\quad \quad +\biggl(\int_{s'=\lambda_n(r_n\wedge s_k^n)}^{\lambda_n(r_n\wedge s_{k+1}^n)}-\int_{s'=r_n\wedge s_k^n}^{r_n\wedge s_{k+1}^n}\biggr)X_\infty\biggl(\frac{\Delta \Phi^n_k}{2},s\biggr)\d s.
\end{align*}
Adapting the subsequence argument previously used to establish the convergence to zero of the second sum on the right-hand side of \eqref{weakintegral:final} and substituting $\Phi^{n_\ell}_k$ with $\Phi^\infty_k$, we find, using \eqref{Ln:conv} as the key condition, that the third sum on the right-hand side of \eqref{weakintegral:final} also converges to zero as $n \to \infty$.

\item For the fourth sum in \eqref{weakintegral:final}, observe that 
$\{t_k^\infty\wedge T\}\subset \{s_k^n\wedge T\}$ for all large enough $n$, since \eqref{Hn:conv} implies 
\[
\lim_{n\to\infty}\sup_{0\leq t\leq T}\bigl\lVert[H^{\theta_0}_n(\cdot,\lambda_n(t))-H^{\theta_0}_n(\cdot,\lambda_n(t-))]-[ H^{\theta_0}_\infty(\cdot,t)-H^{\theta_0}_\infty(\cdot,t-)]\bigr\rVert_m=0.
\]
It follows that  the last sum on the right-hand side of  \eqref{weakintegral:final} satisfies 
\begin{align*}
&\sum_{k=0}^{K}\int_{s'=r_n\wedge s_k^n}^{r_n\wedge s_{k+1}^n}\int_{y'} H_\infty(y',r_n\wedge s_k^n)M_\infty(\d y',\d s')\\
&\quad =\int_{s'=0}^{r_n}\int_{y'}H^{\theta_0}_\infty(y', s'-)\M_\infty(\d y',\d s')\xrightarrow [n\to\infty]{}\int_{s'=0}^{T-\kappa_0}\int_{y'}H^{\theta_0}_\infty(y', s'-)\M_\infty(\d y',\d s'),
\end{align*}
where the limit holds by \eqref{Ln:conv} and the choice of $r_n$: $\lambda_n(r_n)=T-\kappa_0$.\medskip 
\end{itemize}
By using all of the limits we have obtained,  \eqref{I4:probg40} follows. We have proved \eqref{I4:probg40}, and thus \eqref{I4:probg4}. \medskip 

\noindent {\bf Step~4.} From \eqref{I4:probg3} and \eqref{I4:probg4}, we deduce \eqref{I4:probg2} for all sufficiently large integers $n$. Combining \eqref{I4:probg0} and \eqref{I4:probg2}, we conclude that for all large $n \geq 1$,
\[
\P\left(\lvert I^{3,4}_n(f,T)-I^{3,4}_\infty(f,T)\rvert >5\delta\right)\leq 5\eta.
\]
This bound proves the required convergence in probability in \eqref{I4:conv}. 
\end{proof}

Let us close Section~\ref{sec:I34} with the following lemma.

\begin{lem}\label{lem:mildinfty}
For all $n\in \Bbb N\cup\{\infty\}$, $f\in \C_{pd}^0(\R^2)$, and $T\in (0,\infty)$,
\begin{align}\label{id:mildinfty}
X_n(f,T)=\int_{\R^2} f(x) P_TX_0(x)\d x+\int_{s=0}^T\int_y P_{T-s}f(y)M_n(\d y,\d s).
\end{align}
\end{lem}
\begin{proof}
For $n \in \mathbb{N}$, \eqref{id:mildinfty} follows directly from the weak formulation \eqref{def:mild:vep} and the stochastic Fubini theorem: the term $\int_x f(x)\int_{s=0}^T P_{T-s}(x,y)M_n(\d y,\d s)\d x$ exactly matches the stochastic integral in \eqref{id:mildinfty}. Here, the justification for applying the stochastic Fubini theorem relies on the $L^2$ condition in \cite[pp.296--297]{Walsh} and Proposition~\ref{prop:XM} (1$\cc$). For $n = \infty$, the result is obtained by taking the limit as $n \to \infty$ in \eqref{id:mildinfty}, similar to the proof of \eqref{I3:conv} but with a simpler argument; details are omitted. 
\end{proof}

\subsection{End of the proofs of Theorem~\ref{thm:main1} (1$\cc$) and Theorem~\ref{thm:main2}}\label{sec:lcd1-thm2}
To conclude the proof of Theorem~\ref{thm:main1} (1$\cc$), we apply Theorem~\ref{thm:tight} for tightness and Lemma~\ref{lem:mildinfty} to verify \eqref{def:mild}. Additionally, Tsai's theorem~\cite{Tsai:24} ensures that $X_\vep$ admits a unique distributional limit as $\vep \to 0$. The next proposition summarizes the unique characterizations of the remaining limits in Theorem~\ref{thm:tight}, as well as that of the martingale measure $M_\infty$ initially defined in \eqref{XM:relation}. Specifically, it restates \eqref{id:mu=nu}, \eqref{id:mur=nu}, and \eqref{Mncovar:derivation:eq}, and refers to Lemma~\ref{lem:Mfconv} (2$\cc$)--(4$\cc$) for extensions of $M_\infty$ to more general test functions. For future reference, these characterizations are presented in the framework of Section~\ref{sec:lcd1-sk}, independently of Tsai's theorem.

\begin{prop}\label{prop:thm1-1-prelim}
The martingale measure $M_\infty$ and the measures $\nu_\infty$, $\mu_{\infty}$, and $\mathring{\mu}_{\infty}$ are each uniquely determined by $X_\infty$. Moreover, 
\begin{align*}
\int_{t=0}^T\int_x \int_{x'} g(x',x,t)\mu_\infty(\d x',\d x,\d t)&=
\int_{t=0}^T\int_x g(x,x,t)\la M_\infty,M_\infty\ra(\d x,\d t),\\
\int_{t=0}^T\int_{x}g(x,x,t)\nu_\infty(\d x,\d t)
&=\int_{s=0}^T\int_{x}\int_{x'}g(x',x,t)\mu_\infty(\d x',\d x,\d t),\\
\int_{t=0}^T\int_{x}\int_{x'}g(x',x,t)\mathring{\mu}_\infty(\d x',\d x,\d t)
&=\int_{t=0}^T\int_x \left(\int_{x'} g(x',x,t)\phi(x')\d x'\right)\nu_\infty(\d x,\d t),
\end{align*}
which holds for all $g(x',x,t)\in \C^0_{pd}(\R^2\times \R^2\times \R_+)$.
\end{prop}

\begin{proof}[Proof of Theorem~\ref{thm:main2}]
Let $f \in \C_{pd}^1(\R^2 \times [0,T])$. By Propositions~\ref{prop:I1}, \ref{prop:I2}, and \ref{prop:I3I4}, the limits of the three terms $I^1_n(f,T)$,  $I^2_n(f,T)$, and $I^{3}_n(f,T)+I^{4}_n(f,T)$ as $n \to \infty$ are determined. By applying these limits to \eqref{mild:goal}, and using $\mathfrak{L}_0 f(x,t,T)$ as defined in \eqref{def:L0f} and the notation $[\cdot]_{\times}$ from \eqref{def:column}, we obtain the following equation:
\begin{align}
&\int_{t=0}^T\int_xf(x,t)[P_tX_0(x)]^2\d x\d t\notag\\
&\quad =\int_{t=0}^T\int_x \mathfrak L_0f(x,t,T)\nu_\infty(\d x,\d t) +\int_{t=0}^T \int_{x}\biggl(\frac{1}{2\pi }\int_{y''}\int_{y}\phi(y'')(\log\rvert y''-y\rvert )\phi(y)\d y \d y''\notag\\
&\quad \quad -\frac{1}{2\pi }\int_{y''}\int_{y}\phi(y'')(\log\lvert y'' \rvert)\phi(y)\d y \d y''\biggr) f(x,t)\nu_\infty(\d x,\d t)\notag\\
&\quad \quad -2\int_{t'=0}^T\int_{z_1} \biggl( \int_{t=t'}^T \int_x f(x,t)\int_{z_2} 
\begin{bmatrix}
P_{t-t'}(x,z_2)\\
P_{t-t'}(x,z_1)
\end{bmatrix}_\times
X_\infty(\d z_2,t')\d x\d t\biggr)\M_\infty(\d z_1,\d t').\label{thm:main2-final}
\end{align}
Because $\phi$ is a probability density by \eqref{eq:covarM-2}, the first two integrals on the right-hand side combine to give $\int_{t=0}^T \int_x \mathfrak{L} f(x,t,T) \, \nu_\infty(\mathrm{d}x, \mathrm{d}t)$, where $\mathfrak{L}$ is defined in \eqref{def:L}. Since $\langle M_\infty, M_\infty\ra = \nu_\infty$, rearranging \eqref{thm:main2-final} completes the proof of Theorem~\ref{thm:main2}, where the notation $(X, M)$ is now used for the unique limit, as ensured by Theorem~\ref{thm:main1} (1$\cc$).
\end{proof}

\section{Formulas for the limiting covariation measures}\label{sec:lcd2}
\subsection{Solutions to the integro-multiplication equations}\label{sec:lcd2-ime}
In this subsection, we establish Lemma~\ref{lem:main3-prelim}, with its proof deferred to the end of the section. At this stage, \eqref{eq:main3-1} can only be justified by restricting to special test functions. We subsequently demonstrate in Section~\ref{sec:thmmain3end} how to extend this result to all formulas in Theorem~\ref{thm:main2}.

\begin{lem}\label{lem:main3-prelim}
For any nonnegative $g\in \C_{pd}^1(\R^2)$ and $T\in (0,\infty)$, the following formula holds:
\begin{align}
\forall\;0\leq S\leq T,\quad & \int_{t=0}^{S}\int_{\R^2} [P_{T-t}g(x)]^2\langle M_\infty,M_\infty\rangle (\d x,\d t)\notag\\
&\quad =  \int_{t=0}^{S}\int_{\R^2}
[P_{T-t}g(x)]^2\rDBG_{t}X^{\otimes 2}_0(x)\d x\d t\notag\\
&\quad\quad  +2\int_{t'=0}^{S}\int_{z_1} \int_{t=t'}^S\int_x
 \int_{z_2}
[P_{T-t}g(x)]^2
\rDBG_{t-t'}(x;z_1,z_2) X_\infty(\d z_2,t') \d x \d t\notag\\
&\quad \quad \times
 M_\infty(\d z_1,\d t'),\label{mild:goal2}
\end{align}
where $\rDBG_t(x;z_1,z_2)$ and $\rDBG_t F(x)$ are defined as in \eqref{def:rDBG}.
\end{lem}

The following theorem is a crucial step in proving Lemma~\ref{lem:main3-prelim} and may also be of independent interest. In light of Theorem~\ref{thm:main2}, we focus on solving the equation $\mathfrak{L}f(x, t, T) = [P_{T-t}g(x)]^2$, where $\mathfrak{L}f(x, t, T)$ is defined as in \eqref{def:L} and is independent of the initial condition $X_0(\cdot)$. Note that the integral term in $\mathfrak{L}f(x, t, T)$ is well-defined only under sufficient regularity on $f$, such as the $\C^1$-regularity required in Theorem~\ref{thm:main1}; recall Lemma~\ref{lem:I1} (1$\cc$) for details.

\begin{thm}\label{thm:qv}
For any fixed $T\in(0,\infty)$ and any nonnegative $g\in \C_{pd}^1(\R^2)$, there exists a solution $f(x,t)$ to the following problem:
\begin{gather}\label{problem:qv}
\begin{cases}
\mbox{$\mathfrak Lf(x,t,T)$ is well-defined}&\mbox{$\forall\; x\in\R^2$ and $0\leq t<T$},\\
\mathfrak Lf(x,t,T)=[P_{T-t}g(x)]^2& \forall\; x\in \R^2\mbox{ and }0\leq t<T,
\end{cases}
\end{gather}
where 
\begin{align}
\mathfrak Lf(x,t,T)&\,\defeq-\frac{1}{2\pi}
\biggl(\frac{\log (T-t)}{2}+\log 2+\lambda -\frac{\gamma_{\sf EM}}{2}-\int_{y'}\int_{y}\phi(y')\phi(y)\log \rvert y'-y\rvert \d y\d y'\biggr)f(x,t)\notag\\
&\quad\,-\int_{x'}\int_{s'=0}^{T-t}[f(x-x',t+t')-f(x,t)]P_{t'}(x')P_{t'}(x')\d t'\d x'.\label{def:Lv2}
\end{align}
Specifically, the problem stated in \eqref{problem:qv} admits a solution given by 
\[
f(x,t)\equiv \mathfrak m^T_{g}(x,t)\in \C_{pd}^0(\R^2\times [0,T]), 
\]
where, with $g_0^{\otimes 2}(x',x)\,\defeq\,g_0(x')g_0(x)$, $\m_{g_0}^{T}$ is defined as follows:  
\begin{align}
\mathfrak m_{g_0}(x,t)
&\,\defeq\;  \rDBG_t g_0^{\otimes 2}(x),\label{def:mg0}\\
\m^T_{g_0}(x,t)&\,\defeq\, \m_{g_0}(x,T-t),\quad 0\leq t\leq T.\label{def:mTg0}
\end{align}
Moreover, there exists a sequence $\{f_n(x,t)\}_{n\in \Bbb N}\subset \C_{pd}^1(\R^2\times [0,T])$ such that  for all $m\in \Bbb Z_+$,
\begin{align}
\begin{aligned}
f_n(x,t)&\xrightarrow[n\to\infty]{}  \mathfrak m_{g}(x,T-t),&\lvert f_n(x,t)\rvert &\leq\frac{C(m,g,T,\beta)}{1+\lvert x\rvert ^m}, \\
\mathfrak Lf_n(x,t,T)&\xrightarrow[n\to\infty]{}\mathfrak [P_{T-t}g(x)]^2,
&
\lvert \mathfrak Lf_n(x,t,T)\rvert&\leq \frac{C(m,g,T,\beta)}{1+\lvert x\rvert ^m},\label{bound:Lfn}
\end{aligned}
\end{align}
all of which hold pointwise in $(x,t)\in \R^2\times [0,T)$.
\end{thm}

\begin{rmk}\label{rmk:rDBG}
(1$\cc$) For any bounded $F$, the function $\rDBG_t F(\cdot)$ defined in \eqref{def:rDBG} obeys the following uniform bound:
\begin{align}\label{rDBG:bdd}
\sup_{x\in \R^2}\lvert \rDBG_t F(x)\rvert\leq \frac{C(F,T,\beta)}{\log(t\wedge \tfrac{1}{2})^{-1}},\quad \forall\;0<t\leq T<\infty.
\end{align}
This result follows because 
\begin{align}\label{ineq:sbeta}
\s^\beta(t)\leq \frac{C(\beta,T)}{t\log^2(t\wedge \tfrac{1}{2})},\quad 0<t\leq T<\infty;
\end{align} 
see \cite[Proposition~4.5]{C:SDBG1}. We also use the identity $\int \d t/(t\log^2t) = -1/(\log t) + C$ for all $0<t<1$.\medskip 

\noindent (2$\cc$) The second condition in \eqref{problem:qv} indicates that $\mathfrak L\m^T_{g}(x,t,T)$ converges to $g(x)^2$ as $t\nearrow T$. Alternatively, this result follows from the following calculation:
\begin{align*}
\lim_{s\nearrow T}\mathfrak L\m^T_{g}(x,t,T)&=\lim_{t\nearrow T}-\frac{\log (T-t)}{4\pi }\m^T_{g}(x,t)=\lim_{t\searrow 0}-\frac{\log t}{4\pi}\int_0^t \s^\beta(\tau)\d \tau g(x)^2\\
&=\lim_{t\searrow 0}-\log t\int_0^\infty \frac{\beta^u t^u}{\Gamma(u+1)}\d u g(x)^2=g(x)^2.
\end{align*}
Here, the first equality results from combining the observation in (1$\cc$) with Lemma~\ref{lem:m} (2$\cc$), the third equality applies the formula \eqref{def:sbeta} for $\s^\beta(\cdot)$, and the final equality uses the identity $\int t^u\d u = t^u/\log t+ C$ for all $0 < t \neq 1$. Additionally, examining the derivation of the first and second equalities shows that these limits hold uniformly for $x\in \R^2$.
\hfill $\blacklozenge$
\end{rmk}

Our proof of Theorem~\ref{thm:qv} is grounded in Proposition~\ref{prop:m} below, combined with essential regularity properties from Lemma~\ref{lem:m}. Roughly speaking, 
rather than treating \eqref{problem:qv} as a purely analytic problem, our argument utilizes established formulas for the limiting moments of the approximate SHE \eqref{def:SHE:vep}, evaluated with the initial condition $g$ specified in Theorem~\ref{thm:qv}. We denote the corresponding approximate solutions by $\wt{X}_\vep$, with the tilde notation extending analogously to all related objects. As in Corollary~\ref{cor:Sk}, we may assume that there exists a sequence $\{\vep_n\} \subset (0, \ovl]$ with $\vep_n \searrow 0$ such that
\begin{align}\label{lim:tildeX}
(\wt{X}_{n},\wt{\nu}_{n})\,\defeq\,(\wt{X}_{\vep_n},\wt{\nu}_{\vep_n})\xrightarrow[n\to\infty]{} (\wt{X}_\infty,\wt{\nu}_\infty)\mbox{ in $C_{\mathcal M_f(\R^2)}[0,\infty)\times C_{ \mathcal M_f(\R^2)}[0,\infty)$ a.s.},
\end{align}
where $\wt{\nu}_\vep$ is defined as in \eqref{def:covarprocessnu}.

\begin{prop}\label{prop:m}
Let $g\in \C_{pd}^1(\R^2)$ be nonnegative, and let $T\in (0,\infty)$.\medskip 

\noindent {\rm (1$\cc$)} For all $\wt{f}(x,t)\in \C_{pd}^1(\R^2\times [0,T])$,
\begin{align}\label{L:aux}
\int_{t=0}^T \int_x \mathfrak L\wt{f}(x,t,T)\mathfrak m_g(x,t)\d x\d t&=\int_{t=0}^T\int_x \wt{f}(x,t)[P_tg(x)]^2\d x\d s.
\end{align}

\noindent {\rm (2$\cc$)} For all $\wt{f}(x,t)\in \C_{pd}^1(\R^2\times [0,T])$,
\begin{align}\label{m3:goal}
\int_{t=0}^T\int_x \mathfrak L\wt{f}(x,t,T)\mathfrak m_g(x,t)\d x\d t
&=\int_{t=0}^T\int_x \wt{f}(x,t)\mathfrak L\mathfrak m^T_g(x,T-t,T)\d x\d t,
\end{align}
where $\m^T_g$ is as defined in \eqref{def:mTg0}. \medskip 

\noindent {\rm (3$\cc$)} For every $x \in \R^2$ and $0 \leq t < T$, the following identity holds: 
\begin{align}\label{problem:qv0}
\mathfrak L\mathfrak m^T_g(x,t,T)=[P_{T-t}g(x)]^2,
\end{align} 
where $\mathfrak L\mathfrak m^T_g(x,t, T)$ is well defined. 
\end{prop}

Proposition~\ref{prop:m} (3$\cc$) demonstrates that $\m_g^T$ solves the problem in \eqref{problem:qv}, establishing the first claim of Theorem~\ref{thm:qv}. The weak formulation of \eqref{problem:qv0} follows directly by comparing the right-hand sides of \eqref{L:aux} and \eqref{m3:goal} from Proposition~\ref{prop:m} (1$\cc$) and (2$\cc$). Notably, \eqref{m3:goal} links $\mathfrak m_g$ and $\m_g^T$ through $\mathfrak L$, reflecting a time reversal perspective. A similar property arises in the definition of $\rDBG F(x)$ in \eqref{def:rDBG}, with further clarification provided in \eqref{rDBG-time1} and \eqref{rDBG-time2} below.
\medskip 

\begin{rmk}
When $\wt{f}(x,t)\equiv \m_h^T(x,t)$ for nonnegative $h\in \C_{pd}^1(\R^2)$, \eqref{m3:goal} follows directly from \eqref{L:aux}. Specifically, the definitions in \eqref{def:rDBG}, \eqref{def:mg0}, and \eqref{def:mTg0} imply
\begin{gather}
\int_{t=0}^T\int_x \m^T_h(x,t)[P_tg(x)]^2\d x\d t=\int_{t=0}^T\int_x \m_h(x,T-t)[P_tg(x)]^2\d x\d t,\label{rDBG-time1}\\
\int_{t=0}^T\int_x \m^T_g(x,t)[P_th(x)]^2\d x\d t=\int_{t=0}^T\int_x \m_g(x,T-t)[P_th(x)]^2\d x\d t.\label{rDBG-time2}
\end{gather}
Hence, by \eqref{L:aux}, 
\begin{align}\label{L:aux1}
\int_{t=0}^T\int_x \mathfrak L\m^T_h(x,t,T)\m_g(x,t)\d x\d t=\int_{t=0}^T\int_x \m_h(x,t)\mathfrak L\m^T_g(x,t,T)\d x\d t,
\end{align}
which in turn is equivalent to 
\[
\int_{t=0}^T\int_x \mathfrak L\m^T_h(x,t,T)\m_g(x,t)\d x\d t=\int_{t=0}^T\int_x \m^T_h(x,t)\mathfrak L\m^T_g(x,T-t,T)\d x\d t,
\]
as claimed above. In more detail, although we do not directly verify that $(x,t) \mapsto \m_g^T(x,t)$ and $(x,t) \mapsto \m_h^T(x,t)$ belong to $\C^1_{pd}(\R^2 \times [0,T])$, Lemma~\ref{lem:m} below ensures that \eqref{L:aux} still holds when $\wt{f}(x,t)$ is taken as either of these functions. Consequently, \eqref{L:aux1} is valid. \hfill $\blacklozenge$
\end{rmk}

\begin{proof}[Proof of Proposition~\ref{prop:m} (1$\cc$)]
We apply the analogue of \eqref{mild:goal1} for $\wt{X}_\vep$. Specifically, one side of this equation yields the following expectation: 
\begin{align}
&\E\biggl[\int_{t=0}^T\int_{x}\wt{f}(x,t)P_tX_0(x)^2\d x\d t\notag\\
&\quad +
2\int_{t=0}^T\int_{y} \biggl( \int_{s=t}^T \int_{x} f(x,s)\int_{z}\begin{bmatrix}
P_{s-t}(x,y)\\
P_{s-t}(x,z)
\end{bmatrix}_\times X_\infty(\d z,t)\d x\d s\biggr) \M_\infty(\d y,\d t)
\biggr]\notag\\
&\quad=\int_{t=0}^T\int_{x}\wt{f}(x,t)[P_tX_0(x)]^2\d x\d t,\label{eq:mqvrep1}
\end{align}
where the equality holds because \eqref{I4n:sup} ensures that the following process is a martingale:
\[
\int_{t=0}^{t'}\int_{y} \biggl( \int_{s=t}^T \int_{x} f(x,s)\int_{z}\begin{bmatrix}
P_{s-t}(x,y)\\
P_{s-t}(x,z)
\end{bmatrix}_\times X_\infty(\d z,t)\d x\d s\biggr) \M_\infty(\d y,\d t),\quad 0\leq t'\leq T.
\]
By using the analogue of \eqref{mild:goal1} for $\wt{X}_\vep$, the expectation on the left-hand side of \eqref{eq:mqvrep1} also equals the next expectation:
\begin{align}
\E\left[\int_{t=0}^T\int_{x}\mathfrak L\wt{f}(x,t,T)\la \wt{M}_\infty,\wt{M}_\infty\ra(\d x,\d t)\right]=\int_{t=0}^T\int_x \mathfrak L\wt{f}(x,t,T)\m_g(x,t)\d x\d t,\label{eq:mqvrep2}
\end{align}
where the equality holds by Theorem~\ref{thm:main1} (2$\cc$) and \eqref{I22:bound}. Combining \eqref{eq:mqvrep1} and \eqref{eq:mqvrep2} proves Proposition~\ref{prop:m} (1$\cc$). 
\end{proof}

To prove the remainder of Proposition~\ref{prop:m}, parts (2$\cc$) and (3$\cc$), we require the following lemma. This lemma establishes regularity properties of $\m_g$ and ensures that $\mathfrak L\mathfrak m^T_g$, particularly its integral part, is well defined.

\begin{lem}\label{lem:m}
Let $g\in \C_{pd}^1 (\R^2)$ be noonnegative, and let $T\in (0,\infty)$. \medskip 

\noindent {\rm (1$\cc$)} For all $x,x'\in \R^2$, all $t',t\in [0,T]$, and all $m\in \Bbb Z_+$, we have
\begin{align}
\lvert\m_g(x-x',t)-\m_g(x,t)\rvert&\leq C(m,g,T)\left(\frac{\lvert x'\rvert+\lvert x'\rvert^{m+1}}{1+\lvert x\rvert ^m}\right)\int_{\tau=0}^t \s^\beta(\tau)\d \tau, \label{m:spatial}\\
\lvert\m_g(x,t+t')-\m_g(x,t)\rvert&\leq  \frac{C(m,g,T)}{1+\lvert x\rvert ^m}\int_{\tau=t}^{t+t'}\s^\beta(\tau)\d \tau+ \frac{C(m,g,T)}{1+\lvert x\rvert ^m}\int_{\tau=0}^t \s^\beta(\tau)\d \tau\sqrt{t'}.\label{m:temporal}
\end{align}

\noindent {\rm (2$\cc$)} The integral part of $\mathfrak L\mathfrak m^T_g(x,t,T)$, with $x\in \R^2$ and $0\leq t<T$, is absolutely convergent and satisfies the following bound for every $m\in \Bbb Z_+$:
\begin{align}\label{mg:abint}
\int_{x'}\int_{t'=0}^{T-t}\lvert\mathfrak m^T_g(x-x',t+t')-\mathfrak m_g^T(x,t)\rvert P_{t'}(x')P_{t'}(x')\d t'\d x'\leq \frac{C(m,g,T,\beta)}{1+\lvert x\rvert ^m}.
\end{align}
Hence, $\mathfrak L\mathfrak m^T_g(x,t,T)$ is well-defined for any $x\in \R^2$ and $0\leq t<T$. \medskip 

\noindent {\rm (3$\cc$)} The function $(x,t)\mapsto \mathfrak L\mathfrak m^T_g(x,t,T)$ is continuous on $\R^2\times [0,T)$ and admits a unique continuous extension to $\R^2\times \{T\}$.
\end{lem}

\begin{proof}[Proof of Lemma~\ref{lem:m} (1$\cc$)]
To prove \eqref{m:spatial}, recall the definition \eqref{def:mg0} of $\mathfrak m_g$, and consider the following expression, which follows from substitutions:
\begin{align}
&\m_g(x-x',t)-\m_g(x,t)\notag\\
&\quad =\int_{\tau=0}^{t} \int_{y}\begin{bmatrix}
P_{\tau/2}(y)\\
\s^\beta(\tau)
\end{bmatrix}_\times
\int_{\R^4}
\begin{bmatrix}
P_{t-\tau}(x-x'+y, z_2)\\
P_{t-\tau}(x-x'+y, z_1)
\end{bmatrix}_\times 
\begin{bmatrix}
g(z_2)\\
g(z_1)
\end{bmatrix}_\times \begin{bmatrix}
\d z_2\\
\d z_1
\end{bmatrix}_\otimes\d y\d \tau \notag\\
&\quad\quad -\int_{\tau=0}^{t} \int_{y}\begin{bmatrix}
P_{\tau/2}(y)\\
\s^\beta(\tau)
\end{bmatrix}_\times
\int_{\R^4}
\begin{bmatrix}
P_{t-\tau}(x+y, z_2)\\
P_{t-\tau}(x+y, z_1)
\end{bmatrix}_\times 
\begin{bmatrix}
g(z_2)\\
g(z_1)
\end{bmatrix}_\times \begin{bmatrix}
\d z_2\\
\d z_1
\end{bmatrix}_\otimes\d y\d \tau.\label{mg:xdiff}
\end{align}
To bound the right-hand side, let $\{\mathcal W_s\}$ denote a four-dimensional Brownian motion with zero initial condition. Observe that for any $G\in \C_{pd}^1(\R^4)$ and constants $A_0,A_1\in \R^4$, the mean-value theorem gives, for all $0<s\leq T$, 
\begin{align}
&\E[\lvert G(A_0+\mathcal W_s)-G(A_1+\mathcal W_s)\rvert]\notag\\
&\quad \leq  \E\biggl[\lvert A_0-A_1\rvert\sup_{0\leq a\leq 1}\lvert\nabla G(aA_0+(1-a)A_1+\mathcal W_s)\rvert\biggr]\notag\\
&\quad \leq C(m,G)\lvert A_0-A_1\rvert\E\biggl[\,\sup_{0\leq a\leq 1} \frac{1}{1+\lvert aA_0+(1-a)A_1+\mathcal W_s\rvert^m}\biggr]\notag\\
&\quad \leq C(m,G) \lvert A_0-A_1\rvert\E\biggl[\,\sup_{0\leq a\leq 1}  \frac{1+\lvert\mathcal W_s\rvert^m }{1+\lvert aA_0+(1-a)A_1\rvert^m}\,\biggr]\notag\\
&\quad \leq C(m,G,T) \lvert A_0-A_1\rvert\sup_{0\leq a\leq 1}\left(\frac{1}{1+\lvert aA_0+(1-a)A_1\rvert^m}\right),\label{4BM:GTA}
\end{align}
where the third equality uses \eqref{ratiot}, and the fourth equality follows since, by Doob's inequality, $\E[\max_{0\leq r\leq 1}\lvert W_r\rvert^p]<\infty$ for all $p\in (1,\infty)$ for a one-dimensional Brownian motion $W$. By applying \eqref{4BM:GTA} with the choice of $s\equiv t-\tau$, $A_0\equiv (x-x'+y,x-x'+y)$, $A_1\equiv (x+y,x+y)$, and $G(z_1,z_2)\equiv g(z_1)g(z_2)$ to \eqref{mg:xdiff}, we obtain
\begin{align*}
\forall\;0\leq t\leq T,\quad & \lvert\m_g(x-x',t)-\m_g(x,t)\rvert\\
&\quad \leq C(m,g,T)\int_{\tau=0}^t  \int_{y}\begin{bmatrix}
P_{\tau/2}(y)\\
\s^\beta(\tau)
\end{bmatrix}_\times \lvert x'\rvert\left(\sup_{0\leq a\leq 1}\frac{1}{1+\lvert x-ax'+y\rvert^m}\right)\d y\d \tau\\
&\quad \leq C(m,g,T)\int_{\tau=0}^t  \int_{y}\begin{bmatrix}
P_{\tau/2}(y)\\
\s^\beta(\tau)
\end{bmatrix}_\times \lvert x'\rvert\left(\frac{1+\lvert x'\rvert^m+\lvert y \rvert^m}{1+\lvert x\rvert ^m}\right)\d y\d \tau\\
&\quad \leq C(m,g,T)\int_{\tau=0}^t \s^\beta(\tau)\d \tau\left(\frac{\lvert x'\rvert+\lvert x'\rvert^{m+1}}{1+\lvert x\rvert ^m}\right),
\end{align*}
where the second inequality uses \eqref{ratiot}, as well as the inequality $(b+c)^m\leq C(m)(b^m+c^m)$ for all $b,c\geq 0$. The last inequality proves \eqref{m:spatial}. 

To prove \eqref{m:temporal}, we consider the following difference for all $t',t\in [0,T]$, using the definition \eqref{def:mg0} of $\mathfrak m_g$:
\begin{align}
&\m_g(x,t+t')-\m_g(x,t)\notag\\
&\quad =\int_{\tau=0}^{t+t'} \int_{y}\begin{bmatrix}
P_{\tau/2}(x,y)\\
\s^\beta(\tau)
\end{bmatrix}_\times
\int_{\R^4}
\begin{bmatrix}
P_{t+t'-\tau}(y, z_2)\\
P_{t+t'-\tau}(y, z_1)
\end{bmatrix}_\times 
\begin{bmatrix}
g(z_2)\\
g(z_1)
\end{bmatrix}_\times \begin{bmatrix}
\d z_2\\
\d z_1
\end{bmatrix}_\otimes\d y\d \tau\notag\\
&\quad\quad  -\int_{\tau=0}^{t} \int_{y}\begin{bmatrix}
P_{\tau/2}(x,y)\\
\s^\beta(\tau)
\end{bmatrix}_\times
\int_{\R^4}
\begin{bmatrix}
P_{t-\tau}(y, z_2)\\
P_{t-\tau}(y, z_1)
\end{bmatrix}_\times 
\begin{bmatrix}
g(z_2)\\
g(z_1)
\end{bmatrix}_\times \begin{bmatrix}
\d z_2\\
\d z_1
\end{bmatrix}_\otimes\d y\d \tau\notag\\
&\quad =\Delta_{\ref{mgtdiff},1}+\Delta_{\ref{mgtdiff},2},\label{mgtdiff}
\end{align}
where
\begin{align*}
\Delta_{\ref{mgtdiff},1}&\,\defeq\,\int_{\tau=0}^{t+t'} \int_{y}\begin{bmatrix}
P_{\tau/2}(x,y)\\
\s^\beta(\tau)
\end{bmatrix}_\times
\int_{\R^4}
\begin{bmatrix}
P_{t+t'-\tau}(y, z_2)\\
P_{t+t'-\tau}(y, z_1)
\end{bmatrix}_\times 
\begin{bmatrix}
g(z_2)\\
g(z_1)
\end{bmatrix}_\times \begin{bmatrix}
\d z_2\\
\d z_1
\end{bmatrix}_\otimes\d z'\d \tau\\
&\quad -\int_{\tau=0}^{t} \int_{y}\begin{bmatrix}
P_{\tau/2}(x,y)\\
\s^\beta(\tau)
\end{bmatrix}_\times
\int_{\R^4}
\begin{bmatrix}
P_{t+t'-\tau}(y, z_2)\\
P_{t+t'-\tau}(y, z_1)
\end{bmatrix}_\times 
\begin{bmatrix}
g(z_2)\\
g(z_1)
\end{bmatrix}_\times \begin{bmatrix}
\d z_2\\
\d z_1
\end{bmatrix}_\otimes\d y\d \tau,\\
\Delta_{\ref{mgtdiff},2}
&\,\defeq\, \int_{\tau=0}^{t} \int_{y}\begin{bmatrix}
P_{\tau/2}(x,y)\\
\s^\beta(\tau)
\end{bmatrix}_\times
\int_{\R^4}
\begin{bmatrix}
P_{t+t'-\tau}(y, z_2)\\
P_{t+t'-\tau}(y, z_1)
\end{bmatrix}_\times 
\begin{bmatrix}
g(z_2)\\
g(z_1)
\end{bmatrix}_\times \begin{bmatrix}
\d z_2\\
\d z_1
\end{bmatrix}_\otimes\d y\d \tau\\
&\quad -\int_{\tau=0}^{t} \d \tau\int_{y}\begin{bmatrix}
P_{\tau/2}(x,y)\\
\s^\beta(\tau)
\end{bmatrix}_\times
\int_{\R^4}
\begin{bmatrix}
P_{t-\tau}(y, z_2)\\
P_{t-\tau}(y, z_1)
\end{bmatrix}_\times 
\begin{bmatrix}
g(z_2)\\
g(z_1)
\end{bmatrix}_\times \begin{bmatrix}
\d z_2\\
\d z_1
\end{bmatrix}_\otimes\d y\d \tau.
\end{align*}
The two differences in \eqref{mgtdiff} can be bounded as follows. First,
it follows from \eqref{ineq:heat2-2} that
\begin{align}\label{mgtdiff1}
\lvert\Delta_{\ref{mgtdiff},1}\rvert\leq \frac{C(m,g,T)}{1+\lvert x\rvert ^m}\int_{\tau=t}^{t+t'}\s^\beta(\tau)\d \tau,\quad \forall\;m\in \Bbb Z_+.
\end{align}
To bound $\Delta_{\ref{mgtdiff},2}$, we employ an argument similar to \eqref{4BM:GTA}, thus obtaining
\begin{align}
&\bigl|\E\bigl[G(A_0+\mc W_{t+t'-\tau})- G(A_0+\mc W_{t-\tau})\bigr]\bigr|\notag\\
&\quad \leq \E\bigl[\bigl|G(A_0+\sqrt{t+t'-\tau}\mc W_{1})- G(A_0+\sqrt{t-\tau}\mc W_{1})\bigr|\bigr]\notag\\
&\quad \leq\E\biggl[\sqrt{t'}\lvert\mc W_1\rvert\sup_{0\leq a\leq 1}\bigl|\nabla G\bigl(A_0+(a\sqrt{t+t'-\tau}+(1-a)\sqrt{t-\tau})\mc W_{1} \bigr)\bigr|\biggr] \notag\\
&\quad \leq \frac{C(m,G,T)\sqrt{t'}}{1+\lvert A_0\rvert^m},\quad \forall\; m\in \Bbb Z_+,\label{mgtdiff2}
\end{align}
where the second inequality uses the inequality $\lvert\sqrt{b}-\sqrt{c}\rvert\leq \sqrt{b-c}$ for all  $b\geq c\geq 0$. By applying \eqref{mgtdiff1}, \eqref{mgtdiff2} and \eqref{ineq:heat2-2} to \eqref{mgtdiff} and choosing $A_0\equiv (y,y)$ and $G(z_1,z_2)\equiv g(z_1)g(z_2)$ for \eqref{mgtdiff2}, we establish the following inequality, thereby proving \eqref{m:temporal}:
\begin{align*}
&\lvert\m_g(x,t+t')-\m_g(x,t)\rvert\\
&\quad \leq \frac{C(m,g,T)}{1+\lvert x\rvert ^m}\int_{\tau=t}^{t+t'}\s^\beta(\tau)\d \tau+ \frac{C(m,g,T)}{1+\lvert x\rvert ^m}\int_{\tau=0}^t \s^\beta(\tau)\d \tau\sqrt{t'}, \quad \forall\;t',t\in [0,T],\;m\in \Bbb Z_+.
\end{align*}
The proof of Lemma~\ref{lem:m} (1$\cc$) is complete.
\end{proof}

\begin{proof}[Proof of Lemma~\ref{lem:m} (2$\cc$)]
We begin by deriving a bound that features the multiplicative factor $(1+\lvert x\rvert ^m)^{-1}$, as appears on the right-hand side of \eqref{mg:abint}. According to the definition \eqref{def:mTg0} of $\m_g^T$, the left-hand side of \eqref{mg:abint} can be expressed as the following integral: 
\[
\int_{x'}\int_{t'=0}^{T-t}\lvert\mathfrak m_g(x-x',T-t-t')-\mathfrak m_g(x,T-t)\rvert P_{t'}(x')P_{t'}(x')\d t'\d x'.
\]
We can bound this integral by setting $S=T-t$. Specifically, for all $m\in \Bbb Z_+$ and $0< S\leq T$,
\begin{align*}
&\int_{x'}\int_{t'=0}^{S}\lvert\mathfrak m_g(x-x',S-t')-\mathfrak m_g(x,S)\rvert P_{t'}(x')P_{t'}(x')\d t'\d x'\\
&\quad \leq \int_{x'}\int_{t'=0}^S\biggl[C(m,g,T)\left(\frac{\lvert x'\rvert+\lvert x'\rvert^{m+1}}{1+\lvert x\rvert ^m}\right)\int_{\tau=0}^T \s^\beta(\tau)\d \tau\biggr]P_{t'}(x')P_{t'}(x')\d t'\d x'\\
&\quad \quad +\int_{x'}\int_{t'=0}^{S}\biggl[\frac{C(m,g,T)}{1+\lvert x\rvert ^m}\int_{\tau=S-t'}^{S}\s^\beta(\tau)\d \tau+ \frac{C(m,g,T)}{1+\lvert x\rvert ^m}\int_{\tau=0}^T \s^\beta(\tau)\d \tau\sqrt{t'}
\biggr]P_{t'}(x')P_{t'}(x')\d t'\d x',
\end{align*}
where the final inequality uses \eqref{m:spatial} and \eqref{m:temporal}. Therefore, by \eqref{ineq:sbeta} and \eqref{ineq:heat2-1},
\begin{align*}
&\int_{x'}\int_{t'=0}^{S}\lvert\mathfrak m_g(x-x',S-t')-\mathfrak m_g(x,S)\rvert P_{t'}(x')P_{t'}(x')\d t'\d x'\\
&\quad \leq \frac{C(m,g,T,\beta)}{1+\lvert x\rvert ^m}\left(1+\int_{x'}\int_{t'=0}^S \int_{\tau=S-t'}^S \s^\beta(\tau)\d \tau P_{t'}(x')P_{t'}(x')\d t'\d x'\right),\quad 0< S\leq T,
\end{align*}
which establishes the required bound involving $(1+\lvert x\rvert ^m)^{-1}$.

To complete the proof of \eqref{mg:abint}, it suffices to show that
\begin{align}\label{m:int}
\int_{x'}\int_{t'=0}^S \int_{\tau=S-t'}^S \s^\beta(\tau)\d \tau P_{t'}(x')P_{t'}(x')\d t'\d x'\leq C(T,\beta),\quad \forall\;0<S\leq T<\infty.
\end{align}
We consider $S\geq 1/4$ and $0<S<1/4$ separately. First, we perform the following calculation: 
\begin{align}
\forall\;S:0<S\leq T\;\&\;S\geq 1/4,&\quad \int_{x'}\int_{t'=0}^S \int_{\tau=S-t'}^S \s^\beta(\tau)\d \tau P_{t'}(x')P_{t'}(x')\d t'\d x'\notag\\
&\quad\quad \less \int_{t'=0}^S\int_{\tau=S-t'}^S \s^{\beta}(\tau)\d \tau \frac{\d t'}{t'}\notag\\
&\quad\quad \less \int_{t'=0}^{1/8}\int_{\tau=S-t'}^S\s^\beta(\tau)\d \tau\frac{\d t'}{t'}+\biggl(\int_{t'=1/8}^{S}\frac{\d t'}{t'}\biggr)\biggl(\int_{\tau=0}^S \s^\beta(\tau)\d \tau\biggr)\notag\\
&\quad\quad \less  \max_{1/8\leq \tau\leq S}\s^\beta(\tau)+\biggl(\int_{t'=1/8}^{S}\frac{\d t'}{t'}\biggr)\biggl(\int_{\tau=0}^S \s^\beta(\tau)\d \tau\biggr)\notag\\
&\quad\quad \leq C(T,\beta),\label{m:int1}
\end{align}
where the third inequality follows since $S-t'\geq 1/8$ for all $S\geq 1/4$ and
 $0\leq t'\leq 1/8$, and the last inequality holds by \eqref{ineq:sbeta}. To handle the case of $0<S<1/4$, we use
\eqref{ineq:sbeta}. Hence,
\begin{align}
\forall\;0<S<1/4,\quad
& \int_{x'}\int_{t'=0}^S \int_{\tau=S-t'}^S \s^\beta(\tau)\d \tau P_{t'}(x')P_{t'}(x')\d t'\d x'\notag\\
&\quad \leq 
C(\beta)\int_{t'=0}^S \frac{\log S-\log (S-t') }{(\log S)\log (S-t')}\frac{\d t'}{t'}\notag\\
&\quad =
C(\beta)\int_{x=0}^1 \frac{\log S-\log (S-Sx) }{(\log S)\log (S-Sx)}\frac{\d x}{x}\quad\mbox{($\because x=t'/S$)}\notag\\
&\quad \leq\frac{C(\beta)}{-\log S}\int_{x=0}^{3/4}\frac{\d x}{-\log S-\log (1-x)}+\frac{C(\beta)}{-\log S}\int_{x=3/4}^1 \frac{-\log (1-x)\d x}{-\log S-\log (1-x)}\notag\\
&\quad \leq \frac{C(\beta)}{\log^2 S}+\frac{C(\beta)}{-\log S}.\label{m:int2}
\end{align}
Here, the second inequality follows since both $-\log S$ and $-\log(1-x)$ are positive, and to obtain the first term on the right-hand side of the same inequality, we also use the fact that $\lim_{x \searrow 0} -x^{-1} \log(1-x) = 1$. Combining \eqref{m:int1} and \eqref{m:int2} proves \eqref{m:int}. The proof of (2$\cc$) is complete.
\end{proof}

\begin{proof}[Proof of Lemma~\ref{lem:m} (3$\cc$)]
By Remark~\ref{rmk:rDBG} (2$\cc$), it suffices to prove that the mapping $(x,s)\mapsto \mathfrak L\m^T_g(x,s,T)$ is continuous on $\R^2\times[0,T)$. By using \eqref{m:spatial} and \eqref{m:temporal}, this task reduces to proving that the following mapping is continuous on $ \R^2\times (0,\infty)$:
\begin{align}\label{m:contmain}
(x,S)\mapsto \int_{x'}\int_{t'=0}^{S} [\m_g(x-x',S-t')-\m_g(x,S)] P_{t'}(x')P_{t'}(x')\d t'\d x'.
\end{align}
To address this property, we consider the following equation for $\wh{x},\wt{x}\in \R^2$ and $\wh{S}, \wt{S}\in (0,\infty)$:
\begin{align}
&\int_{x'}\int_{t'=0}^{\wt{S}} [\m_g(\wt{x}-x',\wt{S}-t')-\m_g(\wt{x},\wt{S})] P_{t'}(x')P_{t'}(x')\d t'\d x'\notag\\
&\quad -\int_{x'}\int_{t'=0}^{\wh{S}} [\m_g(\wh{x}-x',\wh{S}-t')-\m_g(\wh{x},\wh{S})] P_{t'}(x')P_{t'}(x')\d t'\d x'\notag\\
&\quad =\int_{x'}\int_{t'=0}^{\wt{S}} [\m_g(\wt{x}-x',\wt{S}-t')-\m_g(\wt{x},\wt{S}-t')] P_{t'}(x')P_{t'}(x')\d t'\d x'\notag\\
&\quad \quad +\int_{x'}\int_{t'=0}^{\wt{S}} [\m_g(\wt{x},\wt{S}-t')-\m_g(\wt{x},\wt{S})] P_{t'}(x')P_{t'}(x')\d t'\d x'\notag\\
&\quad \quad -\int_{x'}\int_{t'=0}^{\wh{S}} [\m_g(\wh{x}-x',\wh{S}-t')-\m_g(\wh{x},\wh{S}-t')] P_{t'}(x')P_{t'}(x')\d t'\d x'\notag\\
&\quad \quad -\int_{x'}\int_{t'=0}^{\wh{S}} [\m_g(\wh{x},\wh{S}-t')-\m_g(\wh{x},\wh{S})] P_{t'}(x')P_{t'}(x')\d t'\d x'.
\label{Lm:cont}
\end{align}
Here, splitting the two singular integrals on the left-hand side is justified for two reasons:
(1) Lemma~\ref{lem:m} (2$\cc$) establishes their absolute convergence; (2) the first and third integrals on the right-hand side are also absolutely convergent by \eqref{m:spatial} and \eqref{ineq:heat2-1}.

Part (2) of the above justification of splitting singular integrals extends to the continuity of those integrals on the right-hand side. Specifically, by \eqref{m:spatial} with $m=0$, 
\begin{align*}
&\Bigl|\1_{\{t'\leq \wt{S}\}}[\m_g(\wt{x}-x',\wt{S}-t')-\m_g(\wt{x},\wt{S}-t')]+\1_{\{t'\leq \wh{S}\}}[\m_g(\wh{x}-x',\wh{S}-t')-\m_g(\wh{x},\wh{S}-t')]\Bigr|\\
&\quad \leq C(g,\ol{S},\beta)\lvert x'\rvert,\quad \forall\;\wh{S},\wt{S}\in (0,\ol{S}],\;0<\ol{S}<\infty,
\end{align*}
and by \eqref{ineq:heat2-1}, $\int_{x'}\int_{t'=0}^{\ol{S}} \lvert x'\rvert P_{t'}(x')P_{t'}(x')\d t'\d x'<\infty$. Hence, by dominated convergence and the continuity of $\m_g(x,t)$ in $\R^2\times \R_+$, the following limit holds for every $\ol{S}\in (0,\infty)$:
\begin{align}
&\lim_{\lvert\wt{x}-\wh{x}\rvert\to 0\atop \wt{S}-\wh{S}\searrow 0,\;\wh{S},\wt{S}\in (0,\ol{S}]}\biggl|\int_{x'}\int_{t'=0}^{\wt{S}} [\m_g(\wt{x}-x',\wt{S}-t')-\m_g(\wt{x},\wt{S}-t')] P_{t'}(x')P_{t'}(x')\d t'\d x'\notag\\
&\quad -\int_{x'}\int_{t'=0}^{\wh{S}} [\m_g(\wh{x}-x',\wh{S}-t')-\m_g(\wh{x},\wh{S}-t')] P_{t'}(x')P_{t'}(x')\d t'\d x'\biggr|=0.\label{m:lim1}
\end{align}

The remainder of the right-hand side of \eqref{Lm:cont} is given by the following difference: 
\begin{align*}
&\int_{x'}\int_{t'=0}^{\wt{S}} [\m_g(\wt{x},\wt{S}-t')-\m_g(\wt{x},\wt{S})] P_{t'}(x')P_{t'}(x')\d t'\d x'\\
&\quad -\int_{x'}\int_{t'=0}^{\wh{S}} [\m_g(\wh{x},\wh{S}-t')-\m_g(\wh{x},\wh{S})] P_{t'}(x')P_{t'}(x')\d t'\d x'.
\end{align*}
To show that this difference converges to zero as $\lvert\wt{x}-\wh{x}\rvert \to 0$ and $\wt{S} - \wh{S} \searrow 0$, we use a modified version of the derivation in \eqref{m:int2} to obtain a dominating bound. For all $0 < S_0 < 1/2 < S_1 < \infty$, with $\wh{S}, \wt{S} \in (S_0, S_1)$ and $0 < t' < S_1$, \eqref{m:temporal} with $m = 0$ yields:
\begin{align}
& \Bigl|[\m_g(\wt{x},\wt{S}-t')-\m_g(\wt{x},\wt{S})]\1_{\{t'<\wt{S}\}}-[\m_g(\wh{x},\wh{S}-t')-\m_g(\wh{x},\wh{S})]\1_{\{t'<\wh{S}\}}\Bigr|\notag\\
&\quad \leq C(g,S_1)\int_{\tau=\wt{S}-t'}^{\wt{S}}\s^\beta(\tau)\d \tau\1_{\{t'<\wt{S}\}}+C(g,S_1)\int_{\tau=\wh{S}-t'}^{\wh{S}}\s^\beta(\tau)\d \tau\1_{\{t'<\wh{S}\}}+C(g,S_1,\beta)\sqrt{t'}\notag\\
&\quad \leq C(g,S_0,S_1,\beta)\biggl[\1_{\{t'<S_0/2\}}\biggl(-\log \biggl(1-\frac{t'}{\wt{S}}\biggr)\biggr)+\1_{\{t'<S_0/2\}}\biggl(-\log \biggl(1-\frac{t'}{\wh{S}}\biggr)\biggr)\notag\\
&\quad\quad  +t'+\1_{\{t'\geq S_0/2\}}\int_{\tau=0}^{S_1}\s^\beta(\tau)\d \tau+\sqrt{t'}\biggr].\label{m:limtime0}
\end{align}
Here, the last inequality justified as follows: Let $S'\in \{\wt{S},\wh{S}\}$, and write
\begin{gather*}
\{t'<S'\}=\{t'<S_0/2,t'<S'<3/8\}\cup \{t'<S_0/2,t'<S',S'\geq 3/8\}\cup \{t'\geq S_0/2,t'<S'\}.
\end{gather*}
By splitting the next integral according to this decomposition of $\{t'<S'\}$, we find that 
\begin{align*}
&\int_{\tau=S'-t'}^{S'}\s^\beta(\tau)\d \tau\1_{\{t'<S'\}}\\
&\quad \leq \int_{\tau=S'-t'}^{S'}\s^\beta(\tau)\d \tau\1_{\{t'<S_0/2, t'<S'<3/8\}}+C(\beta',S_0,S_1)t'+ \1_{\{t'\geq S_0/2\}}\int_{\tau=0}^{S_1}\s^\beta(\tau)\d \tau\\
&\quad \leq C(\beta,S_1)\frac{\log S'-\log(S'-t')}{(\log S')\log (S'-t')}\1_{\{t'<S_0/2, t'<S'<3/8\}}+C(\beta,S_0,S_1)t'+\1_{\{t'\geq S_0/2\}}\int_{\tau=0}^{S_1}\s^\beta(\tau)\d \tau\\
&\quad \leq C(\beta,S_0,S_1)\biggl(-\log \biggl(1-\frac{t'}{S'}\biggr)\biggr)\1_{\{t'<S_0/2\}}+C(\beta,S_0,S_1)t'+\1_{\{t'\geq S_0/2\}}\int_{\tau=0}^{S_1}\s^\beta(\tau)\d \tau,
\end{align*}
where the first term on the right-hand side of the second inequality is obtained from \eqref{ineq:sbeta}. The last inequality is sufficient to justify \eqref{m:limtime0}. Furthermore, by monotonicity, \eqref{m:limtime0} gives the following bound, which is the dominating bound we need:
\begin{align}
& \Bigl|[\m_g(\wt{x},\wt{S}-t')-\m_g(\wt{x},\wt{S})]\1_{\{t'<\wt{S}\}}-[\m_g(\wh{x},\wh{S}-t')-\m_g(\wh{x},\wh{S})]\1_{\{t'<\wh{S}\}}\Bigr|\notag\\
&\quad \leq C(g,\beta,S_0,S_1)\biggl[\1_{\{t'<S_0/2\}}\left(-\log \left(1-\frac{t'}{S_0}\right)\right)\notag\\
&\quad\quad  +t'+\1_{\{t'\geq S_0/2\}}\int_{\tau=0}^{S_1}\s^\beta(\tau)\d \tau+\sqrt{t'}\biggr].\label{m:limtime}
\end{align}

To complete the proof, note that the right-hand side of \eqref{m:limtime}, viewed as a function of $t'$, lies in $L^1((0, S_1), (1/t')\d t')$ because $\lim_{x \searrow 0} -x^{-1} \log(1-x) = 1$. Therefore, by dominated convergence and the continuity of the mapping $(x, t) \mapsto \m_g(x, t)$ on $\Bbb R^2 \times \Bbb R_+$, the following limit holds for any $0 < S_0 < S_1 < \infty$:
\begin{align}
&\lim_{\lvert\wt{x}-\wh{x}\rvert\to 0\atop \wt{S}-\wh{S}\searrow 0,\;\wh{S},\wt{S}\in (S_0,S_1)}\biggl|\int_{x'}\int_{t'=0}^{\wt{S}} [\m_g(\wt{x},\wt{S}-t')-\m_g(\wt{x},\wt{S})] P_{t'}(x')P_{t'}(x')\d t'\d x'\notag\\
&\quad -\int_{x'}\int_{t'=0}^{\wh{S}} [\m_g(\wh{x},\wh{S}-t')-\m_g(\wh{x},\wh{S})] P_{t'}(x')P_{t'}(x')\d t'\d x'\biggr|=0.\label{m:lim2} 
\end{align}

Finally, applying \eqref{m:lim1} and \eqref{m:lim2} to \eqref{Lm:cont} proves the required continuity of the mapping in \eqref{m:contmain}. The proof of Lemma~\ref{lem:m} (3$\cc$) is complete.
\end{proof}

\begin{proof}[Proof of Proposition~\ref{prop:m} (2$\cc$)]
Set 
\begin{align}\label{def:dualeq}
I_{\ref{def:dualeq}}\,\defeq\,\int_{t=0}^T \int_x\int_{x'}\int_{t'=0}^{T-t}[\wt{f}(x-x',t+t')-\wt{f}(x,t)]P_{t'}(x')P_{t'}(x')\d t'\d x' \mathfrak m_g(x,t)\d x\d t.
\end{align}
By the definition \eqref{def:Lv2} of $\mathfrak L$, it suffices to show that 
\begin{align}
I_{\ref{def:dualeq}} &=\int_{t=0}^T\int_x \wt{f}(x,t)\biggl(\int_{x'}\int_{t'=0}^{t}[\mathfrak m_g(x-x',t-t')-\mathfrak m_g(x,t)]P_{t'}(x')P_{t'}(x')\d t'\d x'\notag\\
&\quad\quad +\frac{\log t-\log (T-t)}{4\pi}\mathfrak m_g(x,t)\biggr)\d x\d t.\label{dualeq}
\end{align}
More precisely, this equation suffices because, by the definition \eqref{def:mTg0} of $\mathfrak m^T_g(x,t)$, the right-hand side of \eqref{def:dualeq} equals
\begin{align*}
&\int_{t=0}^T\int_x \wt{f}(x,t)\biggl(\int_{x'}\int_{t'=0}^{T-(T-t)}[\mathfrak m^T_g(x-x',T-t+t')-\mathfrak m^T_g(x,T-t)]P_{t'}(x')P_{t'}(x')\d t'\d x'\\
&\quad +\frac{\log [T-(T-t)]}{4\pi}\mathfrak m^T_g(x,T-t)-\frac{\log (T-t)}{4\pi}\mathfrak m_g(x,t)\biggr)\d x\d t.
\end{align*}
Thus, by rearranging, \eqref{dualeq} is sufficient to recover \eqref{m3:goal}. Furthermore, by \eqref{I22:bound} and Lemma~\ref{lem:m} (2$\cc$), both sides of \eqref{m3:goal} are absolutely convergent integrals and therefore well-defined. In the argument below, we establish \eqref{dualeq} by a series of rearrangements and substitutions, using Fubini's theorem repeatedly. The need for these multiple arrangements also arises from the justification of Fubini's theorem.
\medskip 

\noindent {\bf Step 1.} We claim that 
\begin{align}
I_{\ref{def:dualeq}}&= \int_{t=0}^T\int_{t'=t}^{T}\int_{x'}\int_{x}[\wt{f}(x-x',t')-\wt{f}(x,t)]P_{t'-t}(x')P_{t'-t}(x')\mathfrak m_g(x,t)\d x\d x'\d t'\d t.\label{dual:1}
\end{align}
To justify this, observe that by \eqref{I22:bound}, Fubini's theorem allows us to interchange the order of integration on the right-hand side of \eqref{dualeq}: $\int_x \int_{x'} \int_{t'=0}^{T-t}$ becomes $\int_{t'=0}^{T-t} \int_{x'} \int_x$. Performing a substitution for $t'$ then yields \eqref{dual:1}.\medskip

\noindent {\bf Step 2.} We claim that \eqref{dual:1} implies 
\begin{align}
I_{\ref{def:dualeq}}&=\int_{t=0}^T\int_{t'=t}^{T}\int_{x}\int_{x''}P_{t'-t}(x'')P_{t'-t}(x'')[\wt{f}(x,t')\mathfrak m_g(x-x'',t)-\wt{f}(x,t)\mathfrak m_g(x,t)]\notag\\
&\quad\quad  \times \d x''\d x\d t'\d t.\label{dual:2}
\end{align}
To justify this, we use the fact that for all $0<t<t'<T$, 
\begin{align*}
&\int_{x'}\int_{x}\wt{f}(x-x',t')P_{t'-t}(x')P_{t'-t}(x')\mathfrak m_g(x,t)\d x\d x'\\
&\quad =\int_{x''}\int_{\wt{x}} \wt{f}(\wt{x},t')P_{t'-t}(x'')P_{t'-t}(x'')\mathfrak m_g(\wt{x}-x'',t)\d \wt{x}\d x''\quad (\because \wt{x}=x-x',\;x''=-x').
\end{align*}
Since the integral on the right-hand side and
$\int_{x'}\int_{x}\wt{f}(x,t)P_{t'-t}(x')P_{t'-t}(x')\mathfrak m_g(x,t)\d x\d x'$
are both absolutely convergent, \eqref{dual:2} follows upon using the right-hand side of \eqref{dual:1}. \medskip 

\noindent {\bf Step 3.} We claim that \eqref{dual:2} implies
\begin{align}
I_{\ref{def:dualeq}}&=\int_{x}\int_{x''} \biggl(\int_{t=0}^{T}\int_{t'=t}^T \wt{f}(x,t')P_{t'-t}(x'')P_{t'-t}(x'')\mathfrak m_g(x-x'',t)\d t'\d t\notag\\
&\quad  -\int_{t=0}^T\int_{t'=0}^{T-t}  \wt{f}(x,t)P_{t'}(x'')P_{t'}(x'')\mathfrak m_g(x,t)\d t'\d t \biggr)\d x''\d x.\label{dual:3}
\end{align}
To justify this, observe that the integral on the right-hand side of \eqref{dual:2} is absolutely convergent. This follows by expressing
\begin{align*}
& \wt{f}(x,t')\mathfrak m_g(x-x'',t)-\wt{f}(x,t)\mathfrak m_g(x,t)\\
&\quad =[\wt{f}(x,t')-\wt{f}(x,t)]\mathfrak m_g(x-x'',t)+\wt{f}(x,t)[\mathfrak m_g(x-x'',t)-\mathfrak m_g(x,t)].
\end{align*}
Here, the first term permits the use of \eqref{ineq:heat2-1}, since $\wt{f} \in \C_{pd}^1(\mathbb{R}^2 \times [0, T])$ by assumption. The second term similarly allows applications of both \eqref{m:spatial} and \eqref{ineq:heat2-1}. Absolute convergence of the integral on the right-hand side of \eqref{dual:2} justifies invoking Fubini's theorem to switch the order of integration: $\int_{t=0}^T \int_{t'=t}^{T} \int_x \int_{x''}$ becomes $\int_x \int_{x''} \int_{t=0}^T \int_{t'=t}^{T}$. Finally, substituting for $t'$, \eqref{dual:2} yields \eqref{dual:3}.\medskip 

\noindent {\bf Step 4.}
We perform the following calculation, where the first equality follows by applying Fubini's theorem to interchange the order of integration on the right-hand side of \eqref{dual:3} so that $\int_{t=0}^T \int_{t'=t}^T$ becomes $\int_{t'=0}^T \int_{t=0}^{t'}$:
\begin{align}
I_{\ref{def:dualeq}}&=\int_{x}\int_{x''} \biggl(\int_{t'=0}^T\int_{t=0}^{t'} \wt{f}(x,t')P_{t'-t}(x'')P_{t'-t}(x'')\mathfrak m_g(x-x'',t)\d t\d t'\notag\\
&\quad-\int_{t=0}^T\int_{t'=0}^{T-t}  \wt{f}(x,t)P_{t'}(x'')P_{t'}(x'')\mathfrak m_g(x,t)\d t'\d t \biggr)\d x''\d x\notag\\
& =\int_{x}\int_{x''} \biggl(\int_{t=0}^T\int_{t'=0}^{t} \wt{f}(x,t)P_{t'}(x'')P_{t'}(x'')\mathfrak m_g(x-x'',t-t')\d t'\d t\notag\\
&\quad -\int_{t=0}^T\int_{t'=0}^{T-t}  \wt{f}(x,t)P_{t'}(x'')P_{t'}(x'')\mathfrak m_g(x,t)\d t'\d t \biggr)\d x''\d x\notag\\
&=\int_{x}\int_{x''} \biggl(\int_{t=0}^T\int_{t'=0}^{t} \wt{f}(x,t)P_{t'}(x'')P_{t'}(x'')[\mathfrak m_g(x-x'',t-t')-\m_g(x,t)]\d t'\d t\notag\\
&\quad -\int_{t=0}^T\int_{t'=t}^{T-t}  \wt{f}(x,t)P_{t'}(x'')P_{t'}(x'')\mathfrak m_g(x,t)\d t'\d t \biggr)\d x''\d x,\label{dual:6}
\end{align}
where the last equality used the convention that $\int_{t'=t}^{T-t}=-\int_{t'=T-t}^t$ when $T-t<t$. \medskip 

\noindent {\bf Step 5.} The final step is to establish the following equivalent form of \eqref{dualeq}, which completes the present proof:
\begin{align}
I_{\ref{def:dualeq}}&=\int_{t=0}^T\int_x \wt{f}(x,t)\int_{x'}\int_{t'=0}^{t}[\mathfrak m_g(x-x',t-t')-\mathfrak m_g(x,t)] P_{t'}(x')P_{t'}(x')\d t'\d x'\d x\d t\notag\\
&\quad \quad +\int_{t=0}^T\int_x \wt{f}(x,t)\left(\frac{\log t-\log (T-t)}{4\pi }\right)\mathfrak m_g(x,t) \d x\d t.\label{dual:7}
\end{align}
To obtain this equation, we consider \eqref{dual:6} and observe that 
\begin{align}
\int_{t'=t}^{T-t}\int_{x''}P_{t'}(x'')P_{t'}(x'')\d x''\d t'
&=\frac{\log (T-t)-\log t}{4\pi}.\label{dual:7-1} 
\end{align}
Therefore, the following integral is absolutely convergent, and the subsequent equality is justified by Fubini's theorem: 
\begin{align}
&\int_{x}\int_{x''} \int_{t=0}^T\int_{t'=t}^{T-t}  \wt{f}(x,t)P_{t'}(x'')P_{t'}(x'')\mathfrak m_g(x,t)\d t'\d t\d x''\d x\notag\\
&\quad =\int_{t=0}^T\int_x \wt{f}(x,t)\left(\frac{\log (T-t)-\log t}{4\pi }\right)\mathfrak m_g(x,t) \d x\d t.\label{dual:7-2}
\end{align} 
Since the left-hand side of \eqref{dual:7-2} is absolutely convergent, the required equation \eqref{dual:7} follows upon applying \eqref{dual:7-2} to \eqref{dual:6}. 
\end{proof}

\begin{proof}[Proof of Propostion~\ref{prop:m} (4$\cc$)]
The required identity is obtained by comparing the right-hand sides of \eqref{L:aux} and \eqref{m3:goal}, with Lemma~\ref{lem:m} (3$\cc$) providing the reinforcement needed to establish the identity pointwise.
\end{proof}

\begin{proof}[Proof of Theorem~\ref{thm:qv}]
The first part of this theorem, concerning the solution to the problem in \eqref{problem:qv}, has already been established in Proposition~\ref{prop:m} (3$\cc$). Therefore, it remains to construct a sequence $\{f_n(x, t)\}_{n \in \mathbb{N}} \subset \C_{pd}^1(\mathbb{R}^2 \times [0, T])$ such that \eqref{bound:Lfn} holds pointwise for all $(x, t) \in \mathbb{R}^2 \times [0, T)$ and any $m \in \mathbb{Z}_+$. To define these functions, let $\varrho_{\rm space}(y) \in \C_c^{\infty}(\mathbb{R}^2)$ and $\varrho_{\rm time}(r) \in \C_c^{\infty}(\mathbb{R})$ be probability densities supported on $\lvert y \rvert < 1$ and $\lvert r\rvert < 1$, respectively. Then define 
\begin{align*}
\mathbf M_n(\d y,\d r)&\,\defeq\, n^{2}\varrho_{\rm space}(ny)n\varrho_{\rm time}(nr)\d y\d r, \quad (y,r)\in \R^2\times \R,\;n\in \Bbb N,\\
f_n(x,t)&\,\defeq\,\int_{r}\int_{y} \m_{g}\big(x-y,(T+r-t)^+\big)\mathbf M_n(\d y,\d r),\quad 0\leq t\leq T.
\end{align*}
It is standard to verify that $f_n(x, t) \in \C_{pd}^1(\mathbb{R}^2 \times [0, T])$, and that $f_n(x, t)$ converges to $\m_g(x, T-t)$ for all $(x, t) \in \mathbb{R}^2 \times [0, T]$. Furthermore, applying \eqref{m:temporal} with $t = 0$, along with \eqref{ineq:sbeta} and \eqref{ineq:heat2-2}, demonstrates that $f_n$ satisfies the first bound in \eqref{bound:Lfn}. The following Steps~1--4 establish the two properties of $\mathfrak{L} f_n(x, t,T)$ claimed in \eqref{bound:Lfn}, assuming $(x,t)\in \R^2\times [0,T)$. \medskip 

\noindent{\bf Step 1.} We first write out the expression for $\mathfrak{L} f_n(x, t,T)$ by applying \eqref{def:Lv2}: for every $(x,t)\in \R^2\times [0,T)$.
\begin{align}
&\mathfrak Lf_n(x,t,T)\notag\\
&\quad =-\frac{1}{2\pi}\biggl(\frac{\log (T-t)}{2}+\log 2+\lambda -\frac{\gamma_{\sf EM}}{2}-\int_{y''}\int_{y'}\phi(y'')\phi(y')\log \lvert y''-y'\rvert \d y'\d y''\biggr)
\notag\\
&\quad\quad  \times \int_{r}\int_{y}\m_{g}\big(x-y,(T+r-t)^+\big)
\mathbf M_n(\d y,\d r)\notag\\
&\quad \quad -\int_{x'}\int_{t'=0}^{T-t}
\int_{r}\int_{y}\bigl[ \m_{g}\big(x-x'-y,(T+r-t-t')^+\big)- \m_{g}\big(x-y,(T+r-t)^+\big)\bigr]\notag\\
&\quad\quad \times \mathbf M_n(\d y,\d r)
P_{t'}(x')P_{t'}(x')\d t'\d x'.\label{eq:qv-0}
\end{align}
The next two steps analyze two particular parts of the right-hand side of \eqref{eq:qv-0}. \medskip 

\noindent{\bf Step 2.} 
We start by considering the last term on the right-hand side of \eqref{eq:qv-0}, treating the cases $r < 0$ and $r \geq 0$ separately. The argument uses the fact that whenever $T + r - t > 0$ and $\lvert r\rvert < 1$, Lemma~\ref{lem:m} (2$\cc$) applies. Together with \eqref{ratiot}, this yields the following inequality:
\begin{align}
&\int_{x'}\int_{t'=0}^{T+r-t}\bigl| \m_{g}(x-x'-y,T+r-t-t') - \m_{g}(x-y,T+r-t)\bigr|P_{t'}(x')P_{t'}(x')\d t'\d x'\notag\\
&\quad \leq \frac{C(m,g,T,\beta)(1+\lvert y \rvert^m)}{1+\lvert x\rvert ^m}.\label{eq:qv-2-0}
\end{align}

We now address the last term on the right-hand side of \eqref{eq:qv-0} as follows.
Whenever $-1<r<0$ and $T+r-t>0$, 
\begin{align}
&\int_{x'}\int_{t'=0}^{T-t}\bigl[ \m_{g}\big(x-x'-y,(T+r-t-t')^+\big)- \m_{g}\big(x-y,(T+r-t)^+\big)\bigr]P_{t'}(x')P_{t'}(x')\d t'\d x'\notag\\
&\quad = \int_{x'}\int_{t'=0}^{T+r-t}\bigl[ \m_{g}\big(x-x'-y,T+r-t-t'\big)- \m_{g}\big(x-y,T+r-t\big)\bigr]P_{t'}(x')P_{t'}(x')\d t'\d x'\notag\\
&\quad \quad + \int_{x'}\int_{t'=T+r-t}^{T-t}\bigl[ \m_{g}\big(x-x'-y,0\big)- \m_{g}\big(x-y,T+r-t\big)\bigr]P_{t'}(x')P_{t'}(x')\d t'\d x'.\label{eq:qv-2}
\end{align}
The two integrals on the right-hand side can be bounded as follows:
\begin{align}
&\int_{x'}\int_{t'=0}^{T+r-t}\bigl| \m_{g}\big(x-x'-y,T+r-t-t'\big)- \m_{g}\big(x-y,T+r-t\big)\bigr|P_{t'}(x')P_{t'}(x')\d t'\d x'\notag\\
&\quad \leq \frac{C(m,g,T,\beta)(1+\lvert y \rvert^m)}{1+\lvert x\rvert ^m},\label{eq:qv-3}\\
&\int_{x'}\int_{t'=T+r-t}^{T-t}\bigl| \m_{g}\big(x-x'-y,0\big)- \m_{g}\big(x-y,T+r-t\big)\bigr|P_{t'}(x')P_{t'}(x')\d t'\d x'\notag\\
&\quad \leq C(m,g,T,\beta)\frac{\log (T-t)-\log (T+r-t)}{(1+\lvert x-y\rvert^m)\log [(T+r-t)\wedge \frac{1}{2}]^{-1}}\label{eq:qv-4}\\
&\quad \leq \frac{C(m,g,T,\beta)(1+\lvert y \rvert^m)}{1+\lvert x\rvert ^m}.\label{eq:qv-5}
\end{align}
Here, \eqref{eq:qv-3} follows from \eqref{eq:qv-2-0}, and since $\m_g(y,0)\equiv 0$, \eqref{eq:qv-4} is obtained by modifying the proof of \eqref{rDBG:bdd}, specifically by invoking \eqref{def:rDBG}, \eqref{ineq:sbeta} and \eqref{ineq:heat2-2}. Furthermore, \eqref{eq:qv-5} follows from the assumption $r < 0$ and the bound in \eqref{ratiot}. Additionally, whenever $0\leq r<1$, \eqref{eq:qv-2-0} directly gives
\begin{align}
&\int_{x'}\int_{t'=0}^{T-t}\bigl| \m_{g}\big(x-x'-y,(T+r-t-t')^+\big)- \m_{g}\big(x-y,(T+r-t)^+\big)\bigr|P_{t'}(x')P_{t'}(x')\d t'\d x'\notag\\
&\quad \leq \frac{C(m,g,T,\beta)(1+\lvert y \rvert^m)}{1+\lvert x\rvert ^m}.\label{eq:qv-1}
\end{align}
We note, for later use, that whenever $0 \leq r < 1$, the following holds:
\begin{align}
&\int_{x'}\int_{t'=T-t}^{T-t+r}\bigl| \m_{g}\big(x-x'-y,(T+r-t-t')^+\big)- \m_{g}\big(x-y,(T+r-t)^+\big)\bigr|P_{t'}(x')P_{t'}(x')\d t'\d x'\notag\\
&\quad \leq C(m,g,T,\beta)\min\left\{\frac{(1+\lvert y \rvert^m)}{1+\lvert x\rvert ^m},\frac{\log (T+r-t)-\log (T-t)}{(1+\lvert x-y\rvert^m)\log [(T+r-t)\wedge \frac{1}{2}]^{-1}}
\right\}.\label{eq:qv-11}
\end{align}
This follows from \eqref{eq:qv-2-0} and a modification of the argument used to derive \eqref{eq:qv-4}, except that now bounding both of $\m_{g}\big(x-x’-y,(T+r-t-t’)^+\big)$ and $\m_{g}\big(x-y,(T+r-t)^+\big)$.
\medskip 

\noindent{\bf Step 3.} 
We analyze the singular part of the remainder on the right-hand side of \eqref{eq:qv-0}. Its absolute value equals the left-hand side of the first inequality below, while the subsequent inequality follows by an argument similar to that used for \eqref{eq:qv-4}: for $(x,t)\in \R^2\times [0,T)$, 
\begin{align}
&\left|-\frac{1}{2\pi}\cdot \frac{\log (T-t)}{2} \int_r\int_y\m_{g}\big(x-y,(T+r-t)^+\big)
\mathbf M_n(\d y,\d r)\right|\notag\\
&\quad \leq \int_{r}\int_y \frac{C(m,g,T,\beta)\lvert\log (T-t)\rvert}{(1+\lvert x-y\rvert^m)\log [(T+r-t)\wedge \frac{1}{2}]^{-1}}\mathbf M_n(\d y,\d r)\notag\\
&\quad \leq \frac{C(m,g,T,\beta)}{1+\lvert x\rvert ^m}\int_{0}^{1}\frac{\lvert\log (T-t)\rvert}{\log[(T+r/n-t)\wedge \frac{1}{2}]^{-1}}\varrho_{\rm time}(r)\d r\notag\\
&\quad \leq  C(m,g,T,\beta) \frac{\lvert\log (T-t)\rvert}{\log[(T+1-t)\wedge \frac{1}{2}]^{-1}}\notag\\
&\quad \leq C(m,g,T,\beta).\label{eq:qv-6}
\end{align}
Here, the second inequality follows from \eqref{ineq:heat2-2} and a substitution, while the third relies on the fact that $a \mapsto \log (a\wedge \frac{1}{2})^{-1}$ is decreasing for $a>0$. \medskip 

\noindent{\bf Step 4.} In this step, we show that $\mathfrak{L} f_n(x, t)$ satisfies the limit and bound required in \eqref{bound:Lfn}. By using \eqref{m:temporal}, \eqref{eq:qv-2}, \eqref{eq:qv-3}, \eqref{eq:qv-5}, \eqref{eq:qv-1} and \eqref{eq:qv-6}, we obtain the necessary bound in \eqref{bound:Lfn}. Additionally, \eqref{eq:qv-2}, \eqref{eq:qv-3}, \eqref{eq:qv-5} and \eqref{eq:qv-1} altogether justify the application of Fubini's theorem, so that \eqref{eq:qv-0} leads to the following expression: for every $(x,t)\in \R^2\times [0,T)$, 
\begin{align*}
\mathfrak Lf_n(x,t,T)
&=-\frac{1}{2\pi}\biggl(\frac{\log (T-t)}{2}+\log 2+\lambda -\frac{\gamma_{\sf EM}}{2}-\int_{y''}\int_{y'}\phi(y'')\phi(y')\log \lvert y''-y'\rvert \d y'\d y''\biggr)
\\
&\quad  \times \int_{r}\int_{y}\m_{g}\big(x-y,(T+r-t)^+\big)
\mathbf M_n(\d y,\d r)\\
& \quad -
\int_{r}\int_{y}\int_{x'}\int_{t'=0}^{T-t}\bigl[ \m_{g}\big(x-x'-y,(T+r-t-t')^+\big)- \m_{g}\big(x-y,(T+r-t)^+\big)\bigr]\\
& \quad \times P_{t'}(x')P_{t'}(x')\d t'\d x' \mathbf M_n(\d y,\d r).
\end{align*}

By using the previous equality, we can evaluate the limit of $\mathfrak{L} f_n(x, t)$ as $n \to \infty$ as follows. The limit of the first term on the right-hand side is straightforward, while the limit of the second term can be obtained by the following calculation:
\begin{align*}
&\lim_{n\to\infty}-\int_{r}\int_{y}\int_{x'}\int_{t'=0}^{T-t}\bigl[ \m_{g}\big(x-x'-y,(T+r-t-t')^+\big)- \m_{g}\big(x-y,(T+r-t)^+\big)\bigr]\\
&\quad  \times P_{t'}(x')P_{t'}(x')\d t'\d x' \mathbf M_n(\d y,\d r)\\
&\quad=\lim_{n\to\infty}-\int_{r:T+r/n-t>0}\int_{y}\int_{x'}\int_{t'=0}^{T+r/n-t}\\
&\quad \quad \times\bigl[ \m_{g}\big(x-x'-y/n,T+r/n-t-t'\big)- \m_{g}\big(x-y/n,T+r/n-t\big)\bigr]\\
&\quad \quad \times P_{t'}(x')P_{t'}(x')\d t'\d x' \mathbf M_1(\d y,\d r)\\
&\quad=-\int_{x'}\int_{t'=0}^{T-t}[ \m_{g}(x-x'-y,T-t-t')- \m_{g}(x-y,T-t)]P_{t'}(x')P_{t'}(x')\d t'\d x' .
\end{align*}
Here, the first equality follows from \eqref{eq:qv-2}, \eqref{eq:qv-4} and \eqref{eq:qv-11}, before making substitutions for $y$ and $r$. The final limit is obtained by using the continuity property in Lemma~\ref{lem:m} (3$\cc$), together with \eqref{eq:qv-3}, \eqref{eq:qv-1} and \eqref{eq:qv-11} to justify the application of the dominated convergence theorem. We have shown that $\mathfrak{L} f_n(x, t, T)$ converges to $\mathfrak{L} \m^T_g(x, t, T)$ for every $(x,t)\in \R^2\times [0,T)$. Since Proposition~\ref{prop:m} (3$\cc$) already establishes that $\mathfrak{L} \mathbf{m}^T_g(x, t, T) = [P_{T-t} g(x)]^2$, the required limit in \eqref{bound:Lfn} for $\mathfrak{L} f_n(x, t, T)$ follows.
\end{proof}

\begin{proof}[Proof of  Lemma~\ref{lem:main3-prelim}]
For any nonnegative $g\in \C_{pd}^1(\R^2)$, combining Theorem~\ref{thm:main2} and Theorem~\ref{thm:qv} yields the following equation: 
\begin{align}
&\int_{t=0}^T\int_{x}\mathfrak L\mathfrak m_g^T(x,t,T)\la M_\infty,M_\infty\ra(\d x,\d t)\notag\\
&\quad =\int_{t=0}^T\int_{x}\mathfrak m^T_g(x,t)[P_tX_0(x)]^2\d x\d t\notag\\
&\quad \quad +2\int_{t'=0}^T\int_{z_1} \biggl( \int_{t=t'}^T \int_{x} \mathfrak m^T_g(x,t)\int_{z_2}\begin{bmatrix}
P_{t-t'}(x,z_2)\\
P_{t-t'}(x,z_1)
\end{bmatrix}_\times X_\infty(\d z_2,t')\d x\d t\biggr) \M_\infty(\d z_1,\d t').\label{mild:goal1final}
\end{align}
A more detailed justification of \eqref{mild:goal1final} is as follows. First, Theorem~\ref{thm:main1} ensures that \eqref{mild:goal1final} holds with $\m^T_g$ replaced by the sequence $f_n$ from Theorem~\ref{thm:qv}. Second, to pass to the limit as $n \to \infty$ in the Lebesgue-integral terms using $f_n$ in place of $\m^T_g$, we employ the dominated convergence theorem. This is justified by Theorem~\ref{thm:main1} (2$\cc$) and the bounds in \eqref{bound:Lfn}. By using the limits in \eqref{bound:Lfn}, we also determine the limiting values of these Lebesgue-integral terms, given by the left-hand side of \eqref{mild:goal1final} and the first term on its right-hand side. Furthermore, applying Lemma~\ref{lem:I4alt} with $n=\infty$ and the dominated convergence theorem for stochastic integrals shows that \eqref{bound:Lfn} allows the limit in probability for the stochastic integral term with $f_n$, thereby recovering the last term in \eqref{mild:goal1final}.

To finish the proof of Lemma~\ref{lem:main3-prelim}, it remains to rewrite the right-hand side of \eqref{mild:goal1final} so that it matches \eqref{mild:goal2}. This alignment follows directly from the definitions in \eqref{def:rDBG}, \eqref{def:mg0}, and \eqref{def:mTg0}, along with appropriate substitutions.  
\end{proof}

\subsection{End of the proof of Theorem~\ref{thm:main3}}\label{sec:thmmain3end}
We start by proving the following proposition on space-time approximations to the identity, which is a rigorous version of \eqref{special:ati}.  

\begin{prop}\label{prop:spdelta}
Let $\Phi\in \C_c(\R^2)$ be a probability density, and set $\Phi_\eta(y)\,\defeq\,\eta^{-2}\Phi(\eta^{-1}y)$.
For any $T>0$, $x_0\in \R^2$ and bounded $\Psi(y,t):\R^2\times [0,T]\to \R$ continuous at $(x_0,0)$, we have
\begin{align}\label{eq:spdelta}
\lim_{\eta\searrow 0}\frac{2\pi}{\log \eta^{-1}}\int_0^T\int_{\R^2}[P_{t}\Phi_\eta(y)]^2\Psi(x_0-y,t)\d y\d t=\Psi(x_0,0).
\end{align}
Moreover, if, for some $0<T_0<\infty$ and $p\in (0,\infty)$, $\sup_{0\leq t\leq T_0}\lvert\Psi(y,t)\rvert\leq C(\Psi,T_0)/(1+\lvert y \rvert^p)$ for all $y\in \R^2$, then  we have:
\begin{align}\label{bdd:spdelta}
\sup_{T\in (0,T_0]}\sup_{\eta\in (0,1/2)}\left|\frac{2\pi}{\log \eta^{-1}}\int_0^T\int_{\R^2}[P_{t}\Phi_\eta(y)]^2\Psi(x_0-y,t)\d y\d t\right|\leq \frac{C(\Psi,T_0,p,\Phi)}{1+\lvert x_0\rvert^p}.
\end{align}
\end{prop}
\begin{proof}
To prove \eqref{eq:spdelta}, we begin by noting that for any $\delta\in (0,T)$ and any $\gamma>0$,
\begin{align*}
&\lim_{\eta\searrow 0}\frac{2\pi}{\log \eta^{-1}}\int_\delta^T \int_{\R^2}[P_t\Phi_\eta(y)]^2\d y\d t=0,\\
&\lim_{\eta\searrow 0}\frac{2\pi}{\log \eta^{-1}}\int_0^\delta\int_{\lvert y \rvert\geq \gamma} [P_t\Phi_\eta(y)]^2\d y\d t =\lim_{\eta\searrow 0}\frac{2\pi }{\log \eta^{-1}}\int_0^\delta \int_{\lvert y \rvert\geq \gamma} \left(\int_{\R^2} P_t(y,\eta z)\Phi(z)\d z\right)^2\d y\d t=0.
\end{align*}
Therefore, given the continuity of $\Psi$ at $(x_0,0)$, establishing \eqref{eq:spdelta} reduces to verifying the desired limit when $\Psi \equiv 1$. Under this condition, the integrals in \eqref{eq:spdelta} take the following form:
\begin{align} 
&\frac{2\pi}{\log \eta^{-1}}\int_0^T\int_{\R^2}[P_{t}\Phi_\eta(y)]^2\d y\d t\notag\\
&\quad =\frac{2\pi}{\log \eta^{-1}}\int_0^T \iiint_{(\R^2)^3} P_t(y,\eta z_1)P_t(y,\eta z_2)\Phi(z_1)\Phi(z_2)\d z_1\d z_2\d y\d t\notag\\
&\quad =\frac{2\pi}{\log \eta^{-1}}\int_0^T \iiint_{(\R^2)^3} P_t(\eta y',\eta z_1)P_t(\eta y',\eta z_2)\Phi(z_1)\Phi(z_2)\eta^2 \d z_1\d z_2\d y'\d t\label{eq:spacetimed-1}\\
&\quad =\frac{2\pi }{\log \eta^{-1}}\int_0^{\eta^{-2}T} \iiint_{(\R^2)^3} P_{t'}( y', z_1)P_{t'}( y', z_2)\Phi(z_1)\Phi(z_2) \d z_1\d z_2\d y'\d t'\label{eq:spacetimed-2}\\
&\quad =\frac{2\pi}{\log \eta^{-1}}\int_0^{\eta^{-2}T} \iint_{(\R^2)^2} P_{2t'}( z_1,z_2)\Phi(z_1)\Phi(z_2) \d z_1\d z_2\d t', \label{eq:spacetimed-3}
\end{align}
where \eqref{eq:spacetimed-1} results from the change of variables $\eta y' = y$, \eqref{eq:spacetimed-2} from $t' = t / \eta^2$, and \eqref{eq:spacetimed-3} by applying the Chapman--Kolmogorov equation. Since $\Phi$ has compact support and is a probability density, \eqref{eq:spacetimed-3} implies
\begin{align*}
&\lim_{\eta\searrow 0}\frac{2\pi}{\log \eta^{-1}}\int_0^T\int_{\R^2}[P_{t}\Phi_\eta(y)]^2\d y\d t= \lim_{\eta\searrow 0}\frac{2\pi}{\log \eta^{-1}}\cdot \frac{\log (\eta^{-2}T)}{4\pi }=1.
\end{align*}
That is, the required limit with $\Psi\equiv 1$ holds. This suffices to prove \eqref{eq:spdelta} in general.

To obtain \eqref{bdd:spdelta}, recall the assumed decay condition on $\Psi$, and note the bound $1/(1+\lvert y-y'\rvert^p)\leq C(p)(1+\lvert y \rvert^p)/(1+\lvert y' \rvert^p)$ from \eqref{ratiot}. Suppose that $\Phi(y)$ has support contained in $\lvert y \rvert\leq M$. Then for any $T\in (0,T_0]$ and $\eta\in (0,1/2)$, we have
\begin{align*}
&\left|\frac{2\pi}{\log \eta^{-1}}\int_0^T\int_{\R^2}[P_{t}\Phi_\eta(y)]^2\Psi(x_0-y,t)\d y\d t\right|\\
&\quad \leq \frac{2\pi}{\log \eta^{-1}}\int_0^T\int_{\R^2}[P_{t}\Phi_\eta(y)]^2\frac{C(\Psi,T_0,p)(1+\lvert y \rvert^p)}{1+\lvert x_0\rvert^p}\d y\d t\\
&\quad \leq \frac{C(\Psi,T_0,p)}{1+\lvert x_0\rvert^p}\biggl(\frac{2\pi}{\log \eta^{-1}}\int_0^T\int_{\R^2}[P_{t}\Phi_\eta(y)]^2\d y\d t+\frac{2\pi}{\log \eta^{-1}}\int_0^T\int_{\lvert y \rvert\geq M}[P_{t}\Phi_\eta(y)]^2\lvert y \rvert^p\d y\d t\biggr)\\
&\quad \leq   \frac{C(\Psi,T_0,p,\Phi)}{1+\lvert x_0\rvert^p}\biggl(1+\int_0^T\int_{\lvert y \rvert\geq M}\biggl(\int_{\R^2}P_{t}(y,z)\Phi_\eta(z)\d z\biggr)^2\lvert y \rvert^p\d y\d t\biggr)\\
&\quad \leq   \frac{C(\Psi,T_0,p,\Phi)}{1+\lvert x_0\rvert^p}\biggl(1+\int_0^T\frac{1}{t}\int_{\R^2}P_{t}(y/2)\lvert y \rvert^p\d y\d t\biggr)\\
&\quad \leq \frac{C(\Psi,T_0,p,\Phi)}{1+\lvert x_0\rvert^p}.
\end{align*}
Here, the third inequality bounds $[(2\pi)/(\log \eta^{-1})]\int_0^T\int_{\R^2}[P_{t}\Phi_\eta(y)]^2\d y\d t$ by $C(T,\Phi)$ via \eqref{eq:spdelta}. The fourth inequality follows from the facts that $\int \Phi_\eta(z)\d z = 1$ and $P_t(y,z) \leq P_t(y/2)$ for any $z$ in $\supp(\Phi_\eta) \subseteq \{z' ; \lvert z'\rvert \leq M/2\}$ and any $\lvert y \rvert \geq M$. The last inequality uses the assumption $p \in (0,\infty)$ and the property $\E^B_0[\lvert B_t\rvert^p] = t^{p/2} \E^B_0[\lvert B_1\rvert^p]$ for two-dimensional Brownian motion $\{B_t\}$. The final inequality is sufficient to establish \eqref{bdd:spdelta}.
\end{proof}

\begin{proof}[Proofs of (\ref{eq:main3-1}) and (\ref{eq:main3-2}) in Theorem~\ref{thm:main3}]
The proof is organized into seven steps. Steps~1--5 focus on proving
the following result: For all $0<T<\infty$ and $f(x)\in \C^0_{pd}(\R^2)$,
\begin{align}
\begin{split}\label{main3-0}
&\int_{t=0}^{T}\int_{x}\f(x)\la M_\infty,M_\infty\ra(\d x,\d t) =\int_{t=0}^{T}\int_{x}f(x)\rDBG_{t}X^{\otimes 2}_0(x)\d x\d t\\
&\quad \quad +2\int_{t=0}^{T}\int_{t'=0}^t\int_{z_1}
\int_{x}
\f(x)\mathcal K(x,t,z_1,t')\d x\M_\infty(\d z_1,\d t')\d t,
\end{split}
\end{align}
where we use the shorthand notation
\begin{align}
\mathcal K(x,t,z_1,t')&\,\defeq\,
\int_{s=t'}^t\int_{y}\int_{z_2} 
\begin{bmatrix}
P_{(t-s)/2}(x,y)\\
\s^\beta(s-t)
\end{bmatrix}_\times  \begin{bmatrix}
P_{s-t'}(y,z_2)\\
P_{s-t'}(y,z_1)
\end{bmatrix}_\times
X_\infty(\d z_2,t')\d y\d s.\label{def:K2}
\end{align}
In Step~5, we complete the proof of \eqref{main3-0} and further justify \eqref{eq:main3-1} with $f(x,t)\equiv f(x)\in \C^0_{pd}(\R^2)$. Steps~6 and 7 then provide extensions of
\eqref{eq:main3-2} and \eqref{eq:main3-1} to general $f(x,t)\in \C^0_{pd}(\R^2)$, respectively, thus completing the proof of Theorem~\ref{thm:main3}. 
 \medskip 

\noindent {\bf Step~1.} This proof begins with an application of Lemma~\ref{lem:main3-prelim}. We consider the following equation for all $S\in (0,\infty)$, obtained by using \eqref{mild:goal2}:
\begin{align}
\begin{split}\label{mild:goal2-2}
& \int_{t=0}^{S}\int_{x}[P_{S-t}g(x)]^2\la M_\infty,M_\infty\ra(\d x,\d t)
 =  \int_{t=0}^{S} \int_{x}[P_{S-t}g(x)]^2
 \rDBG_{t}X^{\otimes 2}_0(x)\d x\d t\\
&\quad\quad +2\int_{t'=0}^S\int_{x'} \biggl(
\int_{t=t'}^{S}\int_{x}
[P_{S-t}g(x)]^2\mathcal K(x,t,z_1,t')  \d x\d t\biggr)\M_\infty(\d z_1,\d t').
\end{split}
\end{align}

To prove \eqref{main3-0}, we apply \eqref{mild:goal2-2} with the following choice of $g$.
Let $\Phi\in\C_c^\infty(\R^2)$ be a fixed probability density. For any $x_0\in \R^2$ and $\eta\in (0,1/2)$, choose
\[
g(y)\equiv \sqrt{\frac{2\pi}{\log \eta^{-1}}}\Phi_\eta(y-x_0),
\]
where $\Phi_\eta(y)\,\defeq\,\eta^{-2}\Phi(\eta^{-1}y)$. This function $g$ is nonnegative and belongs to $\C_{pd}^1(\R^2)$, and the function $P_tg(x)$ admits the following representation: 
\begin{align*}
P_tg(x)
=\int_{y} P_t(x,y+x_0)\sqrt{\frac{2\pi}{\log \eta^{-1}}}\Phi_\eta(y)\d y=\sqrt{\frac{2\pi}{\log \eta^{-1}}}P_t\Phi_\eta(x-x_0).
\end{align*}
Hence, \eqref{mild:goal2-2} applies, and we get the following equation: 
\begin{align}
&\int_{t=0}^{S}\int_{x}\frac{2\pi}{\log \eta^{-1}}[P_{S-t}\Phi_{\eta}(x-x_0)]^2\la M_\infty,M_\infty\ra(\d x,\d t)\notag\\
&\quad = \int_{t=0}^{S} \int_{x}
\frac{2\pi}{\log \eta^{-1}}[P_{S-t}\Phi_{\eta}(x-x_0)]^2 \rDBG_{t}X^{\otimes 2}_0(x)\d x\d t\notag\\
&\quad\quad+2\int_{t'=0}^S\int_{x'} \biggl(
\int_{t=t'}^{S}\int_{x}
\frac{2\pi}{\log \eta^{-1}}[P_{S-t}\Phi_{\eta}(x-x_0)]^2\mathcal K(x,t,z_1,t')  \d x\d t\biggr)\notag\\
&\quad \quad \times \M_\infty(\d z_1,\d t).\label{mild:goal2-3}
\end{align}
This equation is how our key tool of this proof, Proposition~\ref{prop:spdelta}, enters. \medskip 

\noindent {\bf Step~2.}
In this step, we initiate the proof of the following equation: For all $T\in (0,\infty)$ and $f(x)\in \C^0_{pd}(\R^2)$, 
\begin{align}
&\int_{t=0}^{T}\int_{x} \left(\int_{S=t}^{T} \int_{x_0}f(x_0)\frac{2\pi}{\log \eta^{-1}}[P_{S-t}\Phi_{\eta}(x-x_0)]^2\d x_0\d S\right)\la M_\infty,M_\infty\ra(\d x,\d t)\notag\\
&\quad =\int_{t=0}^{T}\int_{x}\left( \int_{S=t}^{T} \int_{x_0} f(x_0) 
\frac{2\pi}{\log \eta^{-1}}[P_{S-t}\Phi_{\eta}(x-x_0)]^2 \d x_0\d S\right)\rDBG_{t}X^{\otimes 2}_0(x)\d x\d t\notag\\
&\quad\quad+2\int_{t=0}^{T}\int_{t'=0}^t\int_{x'} \biggl[
\int_{x}
\left(\int_{S=t}^T\int_{x_0} f(x_0)
\frac{2\pi}{\log \eta^{-1}}[P_{S-t}\Phi_{\eta}(x-x_0)]^2\d x_0\d S \right)\mathcal K(x,t,z_1,t')\d x \biggr]\notag\\
&\quad \quad \times \M_\infty(\d z_1,\d t')\d t.\label{mild:goal2-4}
\end{align}
 
To obtain \eqref{mild:goal2-4}, first, we integrate both sides of \eqref{mild:goal2-3} against $f(x_0)\d x_0\d S$ for $0\leq S\leq T$ and $x_0\in \R^2$ to get the following equation: 
\begin{align}
&\int_{S=0}^{T}\int_{x_0}f(x_0)\int_{t=0}^{S}\int_{x}\frac{2\pi}{\log \eta^{-1}}[P_{S-t}\Phi_{\eta}(x-x_0)]^2\la M_\infty,M_\infty\ra(\d x,\d t)\d x_0\d S\notag\\
&=\int_{S=0}^{T}\int_{x_0}f(x_0)  \int_{t=0}^{S} \int_{x}
\frac{2\pi}{\log \eta^{-1}}[P_{S-t}\Phi_{\eta}(x-x_0)]^2 \rDBG_{t}X^{\otimes 2}_0(x)\d x\d t\d x_0\d S\notag\\
&\quad+2\int_{S=0}^T\int_{x_0} f(x_0)\int_{t'=0}^S\int_{x'} \biggl(
\int_{t=t'}^{S}\int_{x}
\frac{2\pi}{\log \eta^{-1}}[P_{S-t}\Phi_{\eta}(x-x_0)]^2\mathcal K(x,t,z_1,t')  \d x\d t\biggr)\notag\\
&\quad \times \M_\infty(\d z_1,\d t')\d x_0\d S.\label{mild:goal2-5}
\end{align}
By Fubini's theorem, the left-hand side and the first term on the right-hand side correspond to those in \eqref{mild:goal2-4}. Thus, it remains to show that the final terms in \eqref{mild:goal2-4} and the preceding equation are equal. This follows by justifying two applications of the stochastic Fubini theorem according to the following two changes of order of integration:
 \begin{align}\label{SFT:2ndform}
 M_\infty(\d z_1,\d t')\d x_0\d S=\d x_0\d SM_\infty(\d z_1,\d t'),\quad
 \d tM_\infty(\d z_1,\d t')=M_\infty(\d z_1,\d t')\d t,
\end{align}
noting that $T\geq S\geq t\geq t'\geq 0$ in both of those terms of \eqref{mild:goal2-4} and \eqref{mild:goal2-5}.  

To justify the two applications of the stochastic Fubini theorem mentioned above, we first note that $\lvert f(x_0)\rvert\d x_0 \d S$ on $\R^2 \times [0, T]$ and $\d t$ on $[0, T]$ are both finite measures.  Hence, by Walsh’s $L^2$-condition for the stochastic Fubini theorem (see \cite[pp.296–297]{Walsh}), it suffices to check the following two integrability conditions for any {\bf fixed} $\eta \in (0, 1/2)$:
\begin{align}
&\E\biggl[\int_{S=0}^T\int_{x_0} \lvert f(x_0)\rvert\int_{t'=0}^S\int_{x'} \biggl(
\int_{t=t'}^{S}\int_{x}
\frac{2\pi}{\log \eta^{-1}}[P_{S-t}\Phi_{\eta}(x-x_0)]^2\mathcal K(x,t,x',t')  \d x\d t\biggr)^2\notag\\
&\quad \times \la\M_\infty,M_\infty\ra(\d x',\d t')\d x_0\d S\biggr]<\infty,\label{main3-cond-1}\\
&\E\biggl[\int_{t=0}^{T}\int_{t'=0}^t\int_{x'} \biggl[
\int_{x}
\left(\int_{S=t}^T\int_{x_0} f(x_0)
\frac{2\pi}{\log \eta^{-1}}[P_{S-t}\Phi_{\eta}(x-x_0)]^2\d x_0\d S \right)\mathcal K(x,t,x',t')\d x \biggr]^2\notag\\
&\quad \times \la \M_\infty,M_\infty\ra(\d x',\d t')\d t\biggr]<\infty,\label{main3-cond-2}
\end{align}
which correspond to the first and second equalities in \eqref{SFT:2ndform}, respectively. In more detail, \eqref{main3-cond-1} and \eqref{main3-cond-2} are relevant because the time variables satisfy $T \geq S \geq t \geq t' \geq 0$, matching the structure in the last terms of \eqref{mild:goal2-4} and \eqref{mild:goal2-5}. Steps 3 and 4 below verify \eqref{main3-cond-1} and \eqref{main3-cond-2}, respectively. Our main tools are Theorem~\ref{thm:main1} (2$\cc$) and the definition of $\mathcal K$ in \eqref{def:K2}. \medskip 

\noindent {\bf Step~3.}
In this step, we verify \eqref{main3-cond-1} for any fixed $\eta \in (0, 1/2)$. To begin, we rewrite its left-hand side. Consider the following portion of \eqref{main3-cond-1}, and use \eqref{def:K2} to obtain the first equality below:
\begin{align}
&\int_{t=t'}^{S}\int_{x}
\frac{2\pi}{\log \eta^{-1}}[P_{S-t}\Phi_{\eta}(x-x_0)]^2\mathcal K(x,t,x',t')  \d x\d t\notag\\
&\quad=\int_{t=t'}^{S}\int_{x}
\frac{2\pi}{\log \eta^{-1}}[P_{S-t}\Phi_{\eta}(x-x_0)]^2\notag\\
&\quad \quad \times\int_{s=t'}^t\int_{y}\int_{z} 
\begin{bmatrix}
P_{(t-s)/2}(x,y)\\
\s^\beta(t-s)
\end{bmatrix}_\times  \begin{bmatrix}
P_{s-t'}(y,z_2)\\
P_{s-t'}(y,x')
\end{bmatrix}_\times
X_\infty(\d z_2,t')\d y\d t\d x\d s\notag\\
&\quad =\int_{s=t'}^S \int_{z}\int_{y}\Psi_\eta(S,s,x_0,y)  \begin{bmatrix}
P_{s-t'}(y,z_2)\\
P_{s-t'}(y,x')
\end{bmatrix}_\times \d yX_\infty(\d z_2,t')\d s,\label{Psieta:bdd}
\end{align}
where the second uses the rule $\int_{t=t'}^S\int_{s=t'}^t=\int_{s=t'}^S\int_{t=s}^S$, and
$\Psi_\eta(S,s,x_0,y)$ is a shorthand notation defined as follows: 
\begin{align}\label{def:Psieta}
\Psi_\eta(S,s,x_0,y)\,\defeq\,\int_{t=s}^S\int_{x}\frac{2\pi}{\log \eta^{-1}}[P_{S-t}\Phi_{\eta}(x-x_0)]^2 \begin{bmatrix}
P_{(t-s)/2}(x,y)\\
\s^\beta(t-s)
\end{bmatrix}_\times \d x \d t.
\end{align}
As a result of \eqref{Psieta:bdd}, the expectation in \eqref{main3-cond-1} can be rewritten as follows: 
\begin{align}
&\E\biggl[\int_{S=0}^T\int_{x_0} \lvert f(x_0)\rvert\int_{t'=0}^S\int_{x'} \biggl(
\int_{t=t'}^{S}\int_{x}
\frac{2\pi}{\log \eta^{-1}}[P_{S-t}\Phi_{\eta}(x-x_0)]^2\mathcal K(x,t,x',t')  \d x\d t\biggr)^2\notag\\
&\quad\quad  \times \la\M_\infty,M_\infty\ra(\d x',\d t')\d x_0\d S\biggr]\notag\\
&\quad =\E\biggl[ \int_{S=0}^T\int_{x_0}\lvert f(x_0)\rvert \int_{t'=0}^S\int_{x'}
\biggl(\int_{s=t'}^S \int_{z}\int_{y}\Psi_\eta(S,s,x_0,y)  \begin{bmatrix}
P_{s-t'}(y,z_2)\\
P_{s-t'}(y,x')
\end{bmatrix}_\times \d yX_\infty(\d z_2,t')\d s\biggr)\notag\\
&\quad\quad  \times\biggl(\int_{\wt{s}=t'}^S \int_{\wt{z}}\int_{\wt{y}}\Psi_\eta(S,\wt{s},x_0,\wt{y})  \begin{bmatrix}
P_{\wt{s}-t'}(\wt{y},\wt{z}_2)\\
P_{\wt{s}-t'}(\wt{y},x')
\end{bmatrix}_\times \d \wt{y}X_\infty(\d \wt{z}_2,t')\d \wt{s}\biggr)\notag\\
&\quad\quad  \times \la M_\infty,\M_\infty\ra(\d x',\d t')\d x_0\d S\biggr].\notag
\end{align}

We proceed to the proof of \eqref{main3-cond-1} by analyzing the right-hand side of the last equality. By integrating out $x'$ after using \eqref{X4BMbdd} with $p=1/2$, we obtain from the last equality that 
\begin{align}
&\E\biggl[\int_{S=0}^T\int_{x_0} \lvert f(x_0)\rvert\int_{t'=0}^S\int_{x'} \biggl(
\int_{t=t'}^{S}\int_{x}
\frac{2\pi}{\log \eta^{-1}}[P_{S-t}\Phi_{\eta}(x-x_0)]^2\mathcal K(x,t,x',t')  \d x\d t\biggr)^2\notag\\
&\quad\quad  \times \la\M_\infty,M_\infty\ra(\d x',\d t')\d x_0\d S\biggr]\notag\\
&\quad \leq  C(\lambda,\phi,\lVert X_0\rVert_\infty,T)\int_{S=0}^T\int_{x_0}\lvert f(x_0)\rvert\int_{t'=0}^S\int_{s=t'}^S\int_y\int_{\wt{s}=t'}^S\int_{\wt{y}}\frac{\Psi_\eta(S,s,x_0,y)
\Psi_\eta(S,\wt{s},x_0,\wt{y})}{(s+\wt{s}-2t')\min\{s-t',\wt{s}-t'\}^{1/2}}\notag\\
&\quad\quad  \times \d \wt{y}\d \wt{s}\d y\d s
\d t'\d x_0\d S\notag\\
 &\quad =  C(\lambda,\phi,\lVert X_0\rVert_\infty,T)\int_{S=0}^T\int_{x_0}\lvert f(x_0)\rvert\int_{t'=0}^S\int_{t=t'}^S\int_{s=t'}^t\int_x\int_{\wt{t}=t'}^S\int_{\wt{s}=t'}^{\wt{t}}\int_{\wt{x}}\notag\\
 &\quad \quad \times \frac{1}{(s+\wt{s}-2t')\min\{s-t',\wt{s}-t'\}^{1/2}}\notag\\
&\quad\quad \times\frac{2\pi}{\log \eta^{-1}}[P_{S-t}\Phi_{\eta}(x-x_0)]^2\s^\beta(t-s)\frac{2\pi}{\log \eta^{-1}}[P_{S-\wt{t}}\Phi_\eta(\wt{x}-x_0)]^2\s^\beta(\wt{t}-\wt{s})\notag\\
&\quad \quad \times \d\wt{x}  \d \wt{s}\d \wt{t}\d x\d s\d t
\d t'\d x_0\d S
,\label{main3-cond-1-1}
\end{align}
where the last equality is obtained by writing out $\Psi_\eta(S,s,x_0,y)
\Psi_\eta(S,\wt{s},x_0,\wt{y})$ according to the definition \eqref{def:Psieta} of $\Psi_\eta$, integrating out $y$ and $\wt{y}$, and finally, using
 the rule $\int_{s=t'}^S\int_{t=s}^S=\int_{t=t'}^S\int_{s=t'}^t$. The right-hand side of \eqref{main3-cond-1-1} can be further simplified to a multiple time integral as follows. First, we consider the bound from the following calculation: 
\begin{align}
P_t\Phi_{\eta}(x-x_0)&= \int_{y}P_t( y)\Phi_\eta(x-x_0-y)\d y\notag\\
&\leq C(\Phi,m)\int_y P_t( y)\biggl(\eta^{-2}\frac{1+\lvert \eta^{-1}x_0\rvert^m+\lvert \eta^{-1}y\rvert^m}{1+\lvert \eta^{-1}x\rvert^m}\biggr)\d y\notag\\
&\leq \frac{C(T,\eta,\Phi,m)(1+\lvert x_0\rvert^m)}{1+\lvert x\rvert ^m},
\quad \forall\;0<t\leq T,\;m\in \Bbb N,\label{main3-cond-1-1-cnst}
\end{align}
where the first inequality holds by \eqref{ratiot}. Note that the constant in \eqref{main3-cond-1-1-cnst} \emph{is allowed to depend on $\eta$} as we fix $\eta$ now. 
By \eqref{main3-cond-1-1-cnst} with $m=10$, the bound in \eqref{main3-cond-1-1} implies the following estimate in terms of a multiple time integral: with $C_{\ref{main3-cond-1-2}}=C(\lambda,\phi,\lVert X_0\rVert_\infty,T,\eta,\Phi,f)$, 
\begin{align}
&\E\biggl[\int_{S=0}^T\int_{x_0} \lvert f(x_0)\rvert\int_{t'=0}^S\int_{x'} \biggl(
\int_{t=t'}^{S}\int_{x}
\frac{2\pi}{\log \eta^{-1}}[P_{S-t}\Phi_{\eta}(x-x_0)]^2\mathcal K(x,t,x',t')  \d x\d t\biggr)^2\notag\\
&\quad  \times \la\M_\infty,M_\infty\ra(\d x',\d t')\d x_0\d S\biggr]\notag\\
&\quad\leq  C_{\ref{main3-cond-1-2}}\int_{t'=0}^T\int_{t=t'}^T\int_{s=t'}^t\int_{\wt{t}=t'}^T\int_{\wt{s}=t'}^{\wt{t}}\frac{\s^\beta(t-s)\s^\beta(\wt{t}-\wt{s}) }{(s+\wt{s}-2t')\min\{s-t',\wt{s}-t'\}^{1/2}} \d \wt{s}\d \wt{t}\d s\d t
\d t'.\label{main3-cond-1-2}
\end{align}

To bound the right-hand side of \eqref{main3-cond-1-2}, we show the following two bounds for all $0<t'<T$:
\begin{align}\label{sbeta:bdd-main}
&\int_{t=t'}^T\int_{s=t'}^t\int_{\wt{t}=t'}^T\int_{\wt{s}=t'}^{\wt{t}}\frac{\s^\beta(t-s)\s^\beta(\wt{t}-\wt{s})}{(s-t')^{3/4}(\wt{s}-t')^{3/4}}   \d \wt{s}\d \wt{t}\d s\d t\leq C(\beta,T),\\
&\int_{t=t'}^T\int_{s=t'}^t\int_{\wt{t}=t'}^T\int_{\wt{s}=t'}^{\wt{t}}\frac{\s^\beta(t-s)\s^\beta(\wt{t}-\wt{s}) }{(s+\wt{s}-2t')\min\{s-t',\wt{s}-t'\}^{1/2}} \d \wt{s}\d \wt{t}\d s\d t\leq C(\beta,T).\label{sbeta:bdd-main1}
\end{align}
To prove \eqref{sbeta:bdd-main}, we use \eqref{ineq:sbeta}. Hence, for $0<T'\leq T$ and $a\in (0,1)$, 
\begin{align}
\int_0^{T'}\frac{\s^\beta(t)\d t}{(T'-t)^{a}}&= \int_0^{T'/2}\frac{\s^\beta(t)\d t}{(T'-t)^{a}}+\int_{0}^{T'/2}\frac{\s^\beta(T'-t)\d t}{t^{a}}\notag\\
& \leq \frac{C(\beta,T)}{(T'/2)^{a}\log \{[(T'/2)\wedge \frac{1}{2}]^{-1}\}}+\frac{C(\beta,T,a)(T'/2)^{1-a}}{[(T'/2)\wedge \frac{1}{2}]\log^2[(T'/2)\wedge \frac{1}{2}]}\notag\\
&\leq  \frac{C(\beta,T,a)}{(T'/2)^{a}\log \{[(T'/2)\wedge \frac{1}{2}]^{-1}\}}.\label{sbeta:bdd:conv}
\end{align}
In more detail, the first term on the right-hand side of the second inequality comes from integrating $\s^\beta(t)$, while the second term arises from integrating $t^{-a}$. Applying \eqref{sbeta:bdd:conv} with $a=3/4$ yields the first inequality below: 
\begin{align}
&\int_{t=t'}^T\int_{s=t'}^t\int_{\wt{t}=t'}^T\int_{\wt{s}=t'}^{\wt{t}}\frac{\s^\beta(t-s)\s^\beta(\wt{t}-\wt{s})}{(s-t')^{3/4}(\wt{s}-t')^{3/4}}   \d \wt{s}\d \wt{t}\d s\d t\notag\\
&\quad =\biggl( \int_{t=t'}^T\int_{s=t'}^t\frac{\s^\beta(t-s)}{(s-t')^{3/4}}\d s\d t\biggr)^2 \leq \biggl(\int_{t=t'}^T \frac{C(\beta,T)\d t}{(\frac{t-t'}{2})^{3/4}\log [(\frac{t-t'}{2})\wedge \frac{1}{2}]^{-1}}\biggr)^2\leq C(\beta,T).\notag
\end{align}
This proves \eqref{sbeta:bdd-main}. To prove 
 \eqref{sbeta:bdd-main1}, it suffices to use 
\eqref{sbeta:bdd-main} and the following bound:
\begin{align}
\forall\;s>t',\;\wt{s}>t',\quad 
(s+\wt{s}-2t')^{-1}\min\{s-t',\wt{s}-t'\}^{-1/2}
&\leq 2 (s-t')^{-3/4}(\wt{s}-t')^{-3/4}.\notag
\end{align}
This inequality can be obtained as follows:  
\begin{align}
(s+\wt{s}-2t')^{-1}\min\{s-t',\wt{s}-t'\}^{-1/2}&\leq (s+\wt{s}-2t')^{-1}(s-t')^{-1/2}+(s+\wt{s}-2t')^{-1}(\wt{s}-t')^{-1/2}\notag\\
&\leq 2 (s-t')^{-3/4}(\wt{s}-t')^{-3/4},\label{sbeta:bdd-2}
\end{align}
where the last inequality holds because 
\[
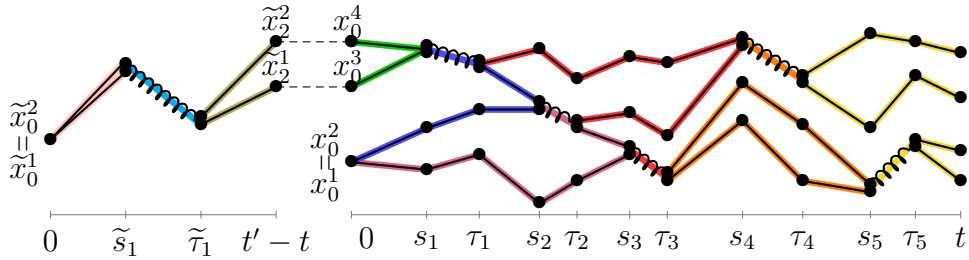
s+\wt{s}-2t'=(s+\wt{s}-2t')^{3/4}(s+\wt{s}-2t')^{1/4}\geq (\wt{s}-t')^{3/4}(s-t')^{1/4},
\]
and similarly, $s+\wt{s}-2t'\geq (s-t')^{3/4}(\wt{s}-t')^{1/4}$. We have proved both \eqref{sbeta:bdd-main} and \eqref{sbeta:bdd-main1}.

We complete the verification of \eqref{main3-cond-1} now. By using \eqref{main3-cond-1-2}  and \eqref{sbeta:bdd-main1}, we get
\begin{align*}
&\E\biggl[\int_{S=0}^T\int_{x_0} \lvert f(x_0)\rvert\int_{t'=0}^S\int_{x'} \biggl(
\int_{t=t'}^{S}\int_{x}
\frac{2\pi}{\log \eta^{-1}}[P_{S-t}\Phi_{\eta}(x-x_0)]^2\mathcal K(x,t,x',t')  \d x\d t\biggr)^2\notag\\
&\quad  \times \la\M_\infty,M_\infty\ra(\d x',\d t')\d x_0\d S\biggr]\leq C(\lambda,\phi,\lVert X_0\rVert_\infty,T,\eta,\Phi,f)<\infty.
\end{align*}
The verification of \eqref{main3-cond-1} for any fixed $\eta\in (0,1/2)$ is complete. \medskip 

\noindent {\bf Step~4.}
To verify \eqref{main3-cond-2}, we begin with the following inequality, which is implied by Proposition~\ref{prop:spdelta} with $\Psi(y,t)\equiv f(y)$: for any $\eta\in (0,1/2)$ and $m\in\Bbb N$,
\begin{align}
&\E\biggl[\int_{t=0}^{T}\int_{t'=0}^t\int_{x'} \biggl[
\int_{x}
\left(\int_{S=t}^T\int_{x_0} f(x_0)
\frac{2\pi}{\log \eta^{-1}}[P_{S-t}\Phi_{\eta}(x-x_0)]^2\d x_0\d S \right)\mathcal K(x,t,x',t')\d x \biggr]^2\notag\\
&\quad \quad \times \la \M_\infty,M_\infty\ra(\d x',\d t')\d t\biggr]\notag\\
&\quad \leq \E\biggl[\int_{t=0}^{T}\int_{t'=0}^t\int_{x'} \biggl(
\int_{x}\frac{C(f,T,\Phi,m)}{1+\lvert x\rvert ^m}\mathcal K(x,t,x',t')\d x \biggr)^2 \la \M_\infty,M_\infty\ra(\d x',\d t')\d t\biggr].
\label{main3-est-10}
\end{align}
Note that the constant $C(f, T, \Phi, m)$ used here is independent of $\eta$. Because of this, \eqref{main3-est-10} will be applied again in Step~5. 

From the definition \eqref{def:K2} of $\mathcal K(x, t, x', t')$ and the moment bound in \eqref{X4BMbdd} with $p=1/2$, it follows that for all integers $m \geq 10$,
\begin{align}
&\E\biggl[\int_{t=0}^{T}\int_{t'=0}^t\int_{x'} \biggl(
\int_{x}\frac{1}{1+\lvert x\rvert ^m}\mathcal K(x,t,x',t')\d x \biggr)^2 \la \M_\infty,M_\infty\ra(\d x',\d t')\d t\biggr]\notag\\
&\quad \leq C(\lambda,\phi,\lVert X_0\rVert_\infty,T)
\int_{t=0}^T\int_{t'=0}^t \int_{x'}\int_{x}\int_{s=t'}^t\int_{y}\int_{\wt{x}}\int_{\wt{s}=t'}^t\int_{\wt{y}}\notag\\
&\quad\quad \times \frac{1}{1+\lvert x\rvert ^m} 
\begin{bmatrix}
P_{(t-s)/2}(x,y)\\
\s^\beta(t-s)
\end{bmatrix}_\times  
P_{s-t'}(y,x')\times \frac{1}{1+\lvert\wt{x}\rvert^m} 
\begin{bmatrix}
P_{(t-\wt{s})/2}(\wt{x},\wt{y})\\
\s^\beta(t-\wt{s})
\end{bmatrix}_\times  
P_{\wt{s}-t'}(\wt{y},x')\notag\\
&\quad \quad \times \frac{\d \wt{y}\d \wt{s} \d \wt{x}\d y\d s\d x \d x'\d t'\d t}{\min\{s-t',\wt{s}-t'\}^{1/2}}\notag\\
&\quad \leq C(\lambda,\phi,\lVert X_0\rVert_\infty,T,m)
\int_{t=0}^T \int_{t'=0}^t\int_{s=t'}^t\int_{\wt{s}=t'}^t\frac{ \s^\beta(t-s)\s^\beta(t-\wt{s})\d \wt{s}\d s\d t'\d t}{{\min\{s-t',\wt{s}-t'\}^{1/2}}}<\infty.\label{main3-est-1}
\end{align}
The two inequalities in \eqref{main3-est-1} are obtained by using the following consideration.
To get the first inequality, we integrate out $\wt{y},\wt{x},y,x,x'$; this uses in particular the following identity: 
\[
\int_{\wt{x}} \int_{\wt{y}}\int_{x'}\int_y P_{(t-s)/2}(x,y)P_{s-t'}(y,x')P_{\wt{s}-t'}(\wt{y},x')P_{(t-\wt{s})/2}(\wt{x},\wt{y}) \d y\d x'\d \wt{y}\d \wt{x}=1,\quad \forall\;x\in \R^2.
\] 
Also, to justify the last inequality in \eqref{main3-est-1}, we have used \eqref{main3-main-est} which we now derive: 
\begin{align}
&\int_{s=t'}^t\int_{\wt{s}=t'}^t\frac{\s^\beta(t-s)\s^\beta(t-\wt{s})}{\min\{s-t',\wt{s}-t'\}^{1/2}}\d s\d \wt{s}\notag\\
&\quad =2\int_{s=t'}^t \s^\beta(t-s)\int_{\wt{s}=t'}^s\frac{\s^\beta(t-\wt{s})}{(\wt{s}-t')^{1/2}}\d \wt{s}\d s\notag\\
&\quad =2\int_{\widecheck{s}=0}^{t-t'} \s^\beta(\widecheck{s})\int_{\wh{s}=t-s}^{t-t'}\frac{\s^\beta(\wh{s})}{(t-t'-\wh{s})^{1/2}}\d\wh{s}\d \widecheck{s}\quad\mbox{($\because \wh{s}=t-\wt{s},\;\widecheck{s}=t-s$)}\notag\\
&\quad \leq \frac{C(\beta,T)}{[(t'-t)\wedge \frac{1}{2}]^{1/2}\log^2[(t-t')\wedge\frac{1}{2}]^{-1}}
,\label{main3-main-est}
\end{align}
where the last inequality follows by using \eqref{sbeta:bdd:conv} with $T'=t-t'$ and $a=1/2$, as well as \eqref{ineq:sbeta}. By \eqref{main3-est-10} and \eqref{main3-est-1}, we have verified \eqref{main3-cond-2}. \medskip 

\noindent {\bf Step~5.} The final step in the proof of \eqref{main3-0} is to take the limits of both sides of \eqref{mild:goal2-4} as $\eta \searrow 0$. Note that by Theorem~\ref{thm:main1} (2$\cc$),
 $\int_{t'=0}^T \int_{x'} (1+\lvert x'\rvert^{10})^{-1} \langle M_\infty, M_\infty \rangle (\d x', \d t') < \infty$ almost surely. Therefore, by dominated convergence, Proposition~\ref{prop:spdelta} implies that as $\eta \searrow 0$, the left-hand side of \eqref{mild:goal2-4} converges almost surely to the left-hand side of \eqref{main3-0}. The same proposition ensures that the first term on the right-hand side of \eqref{mild:goal2-4} converges almost surely to the first term on the right-hand side of \eqref{main3-0}. Finally, by using Proposition~\ref{prop:spdelta} and the bound in \eqref{main3-est-1}, together with It\^{o}'s isometry, the second term on the right-hand side of \eqref{mild:goal2-4} converges in $L^2(\P)$ to the second term on the right-hand side of \eqref{main3-0}. This establishes \eqref{main3-0}, thus \eqref{eq:main3-2} with $f(x,t)\equiv f(x)\in \C^0_{pd}(\R^2)$.
 
To establish \eqref{eq:main3-1} for $f(x, t) \equiv f(x) \in \C^0_{pd}(\R^2)$, we apply Walsh's $L^2$-condition from the stochastic Fubini theorem \cite[pp.296--297]{Walsh}, with justification provided by \eqref{main3-est-1}. \medskip 

\noindent {\bf Step~6.} In this step, we prove \eqref{eq:main3-2} using a standard approximation argument based on \eqref{main3-0} and \eqref{main3-est-1}. Given $f(x,t)\in \C^0_{pd}(\R^2 \times [0, T])$, we begin with the following sequence of functions, which converges uniformly to $f$: 
\[
f_N(x,t)\,\defeq\,f(x,0)\1_{\{0\}}(t)+\sum_{n=1}^{2^N} f\left(x,\frac{nT}{2^N}\right)\1_{\left(\frac{(n-1)T}{2^N},\frac{nT}{2^N}\right]}(t),\quad N\in \Bbb N.
\]
Then, by linearity and \eqref{main3-0}, equation \eqref{eq:main3-2} holds with $f$ replaced by $f_N$.

Our goal in this step is to show that, when $f$ is replaced by $f_N$, the term on the left-hand side of \eqref{eq:main3-2} and the two terms on the right-hand side of \eqref{eq:main3-1} all admit taking the limit as $N \to \infty$ under the integral signs. Specifically, we show below that the three limits
in \eqref{conv:1}, \eqref{conv:2} and \eqref{conv:3-1} below, hold:
\begin{align}
&\int_{t=0}^{T}\int_{x}f_N(x,t)\la M_\infty,M_\infty\ra(\d x,\d t) \xrightarrow[N\to\infty]{\rm a.s.}\int_{t=0}^{T}\int_{x}f(x,t)\la M_\infty,M_\infty\ra(\d x,\d t),\label{conv:1}\\
&\int_{t=0}^{T}\int_{x}f_N(x,t)\rDBG_{t}X^{\otimes 2}_0(x)\d x\d t \xrightarrow[N\to\infty]{\rm a.s.}\int_{t=0}^{T}\int_{x}f(x,t)\rDBG_{t}X^{\otimes 2}_0(x)\d x\d t,\label{conv:2}\\
&\int_{t=0}^T\int_{t'=0}^t\int_{z_1}
\int_{x}
f_N(x,t)\mathcal K(x,t,z_1,t')\d x\M_\infty(\d z_1,\d t')\d t\notag\\
&\quad =\int_{t'=0}^T\int_{z_1}
\int_{t=t'}^T\int_{x}
f_N(x,t)\mathcal K(x,t,z_1,t')\d x\d s\M_\infty(\d z_1,\d t')\label{conv:3-0}\\
&\quad \xrightarrow[N\to\infty]{L^2(\P)}\int_{t'=0}^T\int_{z_1}
\int_{t=t'}^T\int_{x}
f(x,t)\mathcal K(x,t,z_1,t')\d x\d t\M_\infty(\d z_1,\d t')\label{conv:3-1}\\
&\quad \quad \quad \quad =\int_{t=0}^T\int_{t'=0}^t\int_{z_1}
\int_{x}
f(x,t)\mathcal K(x,s,z_1,t')\d x\M_\infty(\d z_1,\d t')\d t.\label{conv:3}
\end{align}
Here, the equalities in \eqref{conv:3-0} and \eqref{conv:3} follow from applying the stochastic Fubini theorem. This application is justified, as \eqref{main3-est-1} suffices for the relevant version of Walsh's $L^2$-condition \cite[pp.296--297]{Walsh}.

The limits in \eqref{conv:1}, \eqref{conv:2}, and \eqref{conv:3-1} can be justified as follows. Both \eqref{conv:1} and \eqref{conv:2} follow from the dominated convergence theorem. In particular, the verification of \eqref{conv:1} relies on the almost sure bound $\int_{t'=0}^T \int_{x'} (1+\lvert x'\rvert^{10})^{-1} \langle M, M \rangle (\d x', \d t') < \infty$, as noted in Step~5. To prove \eqref{conv:3-1}, it suffices to observe that
\begin{align*}
&\lim_{N\to\infty}\E\biggl[\biggl(\int_{t'=0}^T\int_{z_1}
\left(\int_{t=t'}^T\int_{x}
(f-f_N)(x,t)\mathcal K(x,t,z_1,t')\d x\d s\right)\M_\infty(\d z_1,\d t')\biggr)^2\bigg]\\
&\quad =\lim_{N\to\infty}\E\biggl[\int_{t'=0}^T\int_{z_1}
\left(\int_{t=t'}^T\int_{x}
(f-f_N)(x,t)\mathcal K(x,s,z_1,t')\d x\d t\right)^2\la \M_\infty,\M_\infty\ra(\d z_1,\d t')\biggr] =0.
\end{align*}
Here, the last equality holds by dominated convergence, because for all integers $m\geq 10$, 
\begin{align}
&\E\biggl[\int_{t'=0}^T\int_{x'}
\left(\int_{t=t'}^T\int_{x}\frac{1}{1+\lvert x\rvert ^{m}}\mathcal K(x,t,x',t')\d x\d t\right)^2\la \M_\infty,\M_\infty\ra(\d x',\d t')\biggr]\notag\\
&\quad \leq \E\biggl[\int_{t=0}^{T}\int_{t'=0}^s\int_{x'} \biggl(
\int_{x}\frac{1}{1+\lvert x\rvert ^m}\mathcal K(x,t,x',t')\d x \biggr)^2T \la \M_\infty,M_\infty\ra(\d x',\d t')\d t\biggr]<\infty,\notag
\end{align}
where the first inequality follows by applying the Cauchy--Schwarz inequality to the $\d s$-integral on its left-hand side, and the last inequality follows from \eqref{main3-est-1}.

By \eqref{conv:1}, \eqref{conv:2}, and \eqref{conv:3-1}, we have extended \eqref{main3-0} such that for all $f(x,t)\in \C^0_{pd}(\R^2\times [0,T])$, \eqref{eq:main3-2} holds almost surely for all $0<T<\infty$.\medskip 

\noindent {\bf Step 7.} The proof of \eqref{eq:main3-1} has already been provided in Step~6, since the limit in \eqref{conv:3-1} directly enables an extension by invoking \eqref{eq:main3-1} with $f(x,t)\equiv f(x)\in \C^0_{pd}(\R^2)$. The proof of Theorem~\ref{thm:main3} is now complete.
\end{proof}

\begin{proof}[Proof of (\ref{eq:main3-3}) of Theorem~\ref{thm:main3}]
We consider the second term on the right-hand side of \eqref{eq:main3-1} when $f(x,t)\equiv f(x)\in \C_{pd}^0(\R^2)$. This term can be written explicitly as follows: 
\begin{align}
&2\int_{t'=0}^T\int_{z_1} \int_{t=t'}^T
 \int_{z_2}\int_x
f(x)
\rDBG_{t-t'}(x;z_1,z_2)
 X(\d z_2,t') \d x \d t M(\d z_1,\d t')\notag\\
&\quad =2 \int_{t'=0}^T\int_{z_1}\int_{t=0}^{T-t'}\int_{z_2}\int_{s=0}^{t}\int_y \begin{bmatrix}
P_{s/2}f(y)\\
\s^\beta(s)
\end{bmatrix}_\times \begin{bmatrix}
P_{t-s}(y,z_2)\\
P_{t-s}(y,z_1)
\end{bmatrix}_\times \d y\d s
X_\infty(\d z_2,t')\d  tM_\infty(\d z_1,\d t')\notag\\
&\quad =2 \int_{t'=0}^T\int_{z_1}\int_{s=0}^{T-t'}\int_y 
\begin{bmatrix}
P_{s/2}f(y)\\
\s^\beta(s)
\end{bmatrix}_\times\int_{z_2}
\int_{t=s}^{T-t'}  \begin{bmatrix}
P_{t-s}(y,z_2)\\
P_{t-s}(y,z_1)
\end{bmatrix}_\times\notag\\
&\quad\quad  \times \d t
X_\infty(\d z_2,t')\d y\d s  M_\infty(\d z_1,\d t').\label{main3-3-prelim}
\end{align}
Here, the third equality follows from interchanging the order of integration.

To establish \eqref{eq:main3-3}, it suffices to verify Walsh's $L^2$-condition for the stochastic Fubini theorem \cite[pp.296--297]{Walsh}, which, in this case, reduces to the following condition: 
\begin{align}
&\int_{s=0}^{T}\int_y  \begin{bmatrix}
P_{s/2}f(y)\\
\s^\beta(s)
\end{bmatrix}_\times
\E\biggl[ \int_{t'=0}^{T-s}\int_{x'}
\biggl(\int_{z_2}\int_{t=s}^{T-t'}  \begin{bmatrix}
P_{t-s}(y,z_2)\\
P_{t-s}(y,x')
\end{bmatrix}_\times \d t
X_\infty(\d z_2,t')\biggr)\notag\\
&\quad\times \biggl(\int_{\wt{z}_2} \int_{\wt{t}=s}^{T-t'} \begin{bmatrix}
P_{\wt{t}-s}(y,\wt{z}_2)\\
P_{\wt{t}-s}(y,x')
\end{bmatrix}_\times \d \wt{t}
X_\infty(\d \wt{z}_2,t')\biggr)\la M_\infty,M_\infty\ra(\d x',\d t')\biggr] \d y\d s<\infty\label{main3-3-SFT}
\end{align}
To bound the left-hand side, we first apply \eqref{X4BMbdd}, integrating out $z_2$ and $\wt{z}_2$, and then integrate out $x'$. Thus, for any $p \in (0,1)$, up to a multiplicative constant depending only on $(\lambda, \phi, \lVert X_0\rVert_\infty, T, p)$, the left-hand side of \eqref{main3-3-SFT} is bounded by the following integral: 
\begin{align*}
&\int_{s=0}^{T}\int_y 
\begin{bmatrix}
P_{s/2}f(y)\\
\s^\beta(s)
\end{bmatrix}_\times 
\int_{t'=0}^{T-s}
\int_{t=s}^{T-t'}
\int_{\wt{t}=s}^{T-t'}\frac{1}{\min\{t-s,\wt{t}-s\}^p}\frac{1}{t+\wt{t}-2s}\d \wt{t}\d t\d t'\d y\d s.
\end{align*}
This integral is finite because
\begin{align*}
\int_{t=0}^T\int_{\wt{t}=0}^{T}\frac{1}{\min\{t,\wt{t}\}^p(t+\wt{t})}\d \wt{t}\d t
&=2\int_{t=0}^T\int_{\wt{t}=0}^{t}\frac{1}{(\wt{t})^p(t+\wt{t})}\d \wt{t}\d t \leq C(p)\int_{t=0}^T\frac{\d t}{t^p}<\infty;
\end{align*}
see Lemma~\ref{lem:gibdd} for justification of the penultimate inequality. Therefore, \eqref{main3-3-SFT} is satisfied, allowing us to interchange $\d y \d s$ and $M_\infty(\d z_1,\d t')$ in \eqref{main3-3-prelim} to obtain \eqref{eq:main3-3}. 
\end{proof}

\section{Moments and their applications}\label{sec:mom}
In this section, we consider solutions $X_\vep$ to the approximate SHE in \eqref{def:SHE:vep}, where the initial conditions $X_0(\cdot)$ are bounded, nonnegative, and independent of $\vep$. Our aim is to prove Theorem~\ref{thm:mombdd}, Theorem~\ref{thm:tight}, and Theorem~\ref{thm:main1}~(2$\cc$). The proofs appear in Sections~\ref{sec:dual}, \ref{sec:jointC}, and \ref{sec:covaridbdd}, respectively.

\subsection{Proof of Theorem~\ref{thm:mombdd}: approximate mixed moment bounds}\label{sec:dual}
We first provide preliminary results and definitions in Sections~\ref{sec:series} and~\ref{sec:graphical}, following the methodology of \cite{C:DBG}. The proof of Theorem~\ref{thm:mombdd} is then presented in Section~\ref{sec:mombdd}. Note that our notation below differs slightly from that of \cite{C:DBG}.
 
 \subsubsection{Diagrammatic expansions}\label{sec:series}
 In the sequel, we work with the following generalization of \eqref{Mom:dual}. For any $\vep\in (0,\ovl]$, integer $N\geq 1$, and $x_0=(x_0^1,\cdots,x_0^N)$ with $x_0^k\in \R^2$ for all $1\leq k\leq N$, it holds that
 \begin{linenomath*}\begin{align}\label{eq:momdual2}
 \begin{split}
&\E\Bigg[\prod_{i=1}^N X_\vep(x_0^i,t)\Bigg]=\E^{B}_{x_0}\Bigg[\exp\Bigg\{\sum_{\bi\in \mc E_N}\lv \int_0^t \d r \varphi_\vep(B^\bi_r)\Bigg\}\prod_{i=1}^NX_0(B_t^i)\Bigg].
\end{split}
\end{align}
\end{linenomath*}
Here, we use the following notations: $\mathcal E_N\,\defeq\,\{(i\prime,i)\in \Bbb N^2;1\leq i<i\prime\leq N\}$; $\lv$ is chosen to satisfy \eqref{def:lv}; with $\phi_\vep$ defined in \eqref{def:L1},
\begin{linenomath*}\begin{align}
\varphi_\vep(x)&\,\defeq \,\phi_\vep(\two x)=\vep^{-2}\varphi(\vep^{-1}x)\quad (\varphi\equiv \varphi_1);\label{def:varphi}
\end{align}\end{linenomath*}
 the processes $B^1,\cdots,B^{N}$ are independent two-dimensional standard Brownian motions collectively written as $B=(B^1,\cdots,B^N)$ with $B^k_0=x^k_0$ under $\E^B_{x_0}$ so that
\begin{linenomath*}\begin{align}
B^{\bi}_t\;\defeq\, \frac{B^{i\prime}_t-B_t^i}{\two},\quad \bi=(i\prime,i)\in \mathcal E_N,\label{def:Bi}
\end{align}\end{linenomath*} 
define two-dimensional standard Brownian motions. We also assume that $B$ is independent of $X_\vep$.

The proof of Theorem~\ref{thm:mombdd} relies crucially on the series expansion in the proposition below, where the Dirac delta functions are used rigorously through the following definition: for any $a\in I\subset \R$ and function $F$ defined on $I$,
\[
\int_{I}\d \tau \delta_{a}(\tau) F(\tau)\,\defeq\, F(a).
\] 

\begin{prop}\label{prop:graphical}
For any given nonnegative function $f$ defined on $(\R^{2})^N$, we have
\begin{linenomath*}\begin{align}
\begin{split}
&\E^B_{x_0}\Bigg[\exp\Bigg\{ \sum_{\bi\in \mathcal E_N}\lv \int_0^t \d r\varphi_\vep(B^{\bi}_r)\Bigg\}f(B_t)\Bigg]\\
&\quad =
\E^B_{x_0}[f(B_t)] +\sum_{m=1}^\infty \sum_{\stackrel{\scriptstyle \bi_1,\cdots,\bi_m\in \mathcal E_N}{\bi_1\neq \cdots \neq \bi_m}}\sum_{\spin \in \{0,1\}^m} \int_{0< s_1< \cdots< s_m< t}\d \bs s_m Q^{\lambda,\phi;\bi_1,\cdots,\bi_m;\spin}_{\vep;s_1,\cdots,s_m,t}f( x_0),\label{P:vepseries}
\end{split}
\end{align}\end{linenomath*}
where
$\bs s_m\,\defeq\,(s_1,\cdots,s_m)$, $Q^{\lambda,\phi;\bi_1,\cdots,\bi_m;\spin}_{\vep;s_1,\cdots,s_m,t}f( x_0)$ are iterated integrals
defined as follows: 
\begin{align}
&Q^{\lambda,\phi;\bi_1,\cdots,\bi_m;\spin}_{\vep;s_1,\cdots,s_m,t}f(x_0)\notag\\
&\quad\defeq\,\int_{[s_1,s_2)}\d \tau_1\cdots \int_{[s_m,t)}  \d \tau_m \Bigg(\prod_{\scriptstyle \ell:\sigma_\ell=0 \atop \scriptstyle 1\leq \ell\leq m}\delta_{s_\ell}(\tau_\ell)\Bigg)\int \d \mathcal G^\spin_{\vep;0} \int \d \mathcal G^{\bs \sigma}_{\vep;1}  \cdots \int \d \mathcal G^\spin_{\vep;m}  f,\label{int:P:graphical0}
\end{align}
and $\int \d \mathcal G^\spin_{\vep;\ell}F$, given by iterated Lebesgue integrals depending on $(\tau_1,\cdots,\tau_m,s_1,\cdots,s_m,t)$, are the {\bf graphical integrals} defined in \cite{C:DBG}.
\end{prop}

For $N \geq 3$, Proposition~\ref{prop:graphical} follows by combining \cite[(2.12), (4.1), and (4.43)]{C:DBG}. For $N = 2$, the triple series in \eqref{P:vepseries} simplifies to two terms with $m = 1$. While \cite{C:DBG} does not directly connect this case to graphical integrals, \eqref{P:vepseries} can be derived more simply using \cite[(2.39) and Lemma~4.4]{C:DBG}, as shown in Lemma~\ref{lem:mix2mom}.

While implicit due to the graphical integrals, the key feature of the expansion in \eqref{P:vepseries} is that, at any given time, no more than one pair of particles interacts through the ``binding'' of $B^{i\prime}$ and $B^i$ with $\varphi_\vep$. This property arises from the following elementary integral expansion: for all $t\in (0,\infty)$, finite sets $\mathcal E$ such that $E=\#\mathcal E\geq 1$, and Borel measurable functions $\varphi(\bi ,r):\mathcal E\times [0,t]\to \R$ either nonnegative or bounded, we have
\begin{linenomath*}\begin{align}
\begin{split}
&\exp\Bigg\{\sum_{\bi\in \mathcal E}\int_0^t \varphi(\bi,r)\d r\Bigg\}\\
&\quad  =1+\sum_{m=1}^\infty\sum_{\stackrel{\scriptstyle \bi_1,\cdots,\bi_m\in \mathcal E}{\bi_1\neq \cdots \neq \bi_m}}\int_{0< s_1< \cdots<s_m< t}\prod_{n=1}^{m}\varphi(\bi_n,s_n)\exp\Bigg\{\int_{s_{n}}^{s_{n+1}}\varphi(\bi_n,r)\d r\Bigg\}\d \bs s_m.\label{exp:id-alt}
\end{split}
\end{align}\end{linenomath*}
Here, under the integral over $0<s_1<\cdots<s_m<t$, we write $s_{m+1}=t$ and $\bs s_m=( s_1,\cdots, s_m)$. The original proof of \eqref{exp:id-alt} for $N \geq 3$ in \cite{C:DBG} is based on Poisson random measures. This approach is motivated by the fact that Poisson point processes in $\mathcal E$ exhibit nonconcurrent ($0 < s_1 < \cdots < s_m < t$) and nonconsecutive ($\bi_1 \neq \cdots \neq \bi_m$) arrivals, as seen in an earlier resolvent expansion in \cite{GQT} when establishing limiting moments for the case that $X_0(\cdot)$ is in $L^2$. However, \eqref{exp:id-alt} can also be derived more directly by starting with the following argument and iterating indefinitely when $\varphi(\bi,r)$ is bounded: by applying the fundamental theorem of calculus to the functions $s\mapsto \mathcal E_{s,t}\,\defeq\,\exp\{\int_s^t \sum_{\bj\in \mc E}\varphi(\bj,r)\d r\}$ and $s\mapsto \exp\{\sum_{\bj:\bj\neq \bi}\int_s^t \varphi(\bj,r)\d r\}$,
\begin{align*}
\mathcal E_{0,t}
&=1+\sum_{\bi_1\in \mc E}\int_0^t \varphi(s_1,\bi_1)\mc E_{s_1,t}\d s_1\\
&=1+\sum_{\bi_1\in \mc E}\int_0^t \varphi(s_1,\bi_1)\mc E^{\bi_1}_{s_1,t}\d s_1+\sum_{\bi_1\in \mc E}\sum_{\bi_2:\bi_2\neq \bi_1}\int_0^t \varphi(s_1,\bi_1)\int_{s_1}^t\mc E^{\bi_1}_{s_1,s_2}\varphi(s_2,\bi_2)\mathcal E_{s_2,t}\d s_2 \d s_1,
\end{align*}
where $\mathcal E^\bi_{s,t}\,\defeq\, \exp\{\int_s^t \varphi(\bi,r)\d r\}$.

To facilitate the forthcoming applications of Proposition~\ref{prop:graphical},
we recall the definitions of the graphical integrals that validate \eqref{int:P:graphical0} in Section~\ref{sec:graphical} and incorporate some minor extensions as well. Also, we often adopt the following change of variables: 
\begin{align}
\begin{split}\label{def:uv123}
u_1&=s_1,\; u_2=s_2-\tau_1,\cdots,\; u_m=s_m-\tau_{m-1},\;u_{m+1}=t-\tau_m,\\
v_1&=\tau_1-s_1,\;v_2=\tau_2-s_2,\cdots, v_m=\tau_m-s_m.
\end{split}
\end{align}
This and \eqref{int:P:graphical0} yield the following alternative representations of the summands in \eqref{P:vepseries}: 
\begin{align}
&\int_{0< s_1< \cdots< s_m< t}\d \bs s_mQ^{\lambda,\phi;\bi_1,\cdots,\bi_m;\spin}_{\vep;s_1,\cdots,s_m,t}f(x_0)\notag\\
&\quad =\int_{0< s_1< \cdots< s_m< t}\d \bs s_m\int_{[s_1,s_2)}\d \tau_1\cdots \int_{[s_m,t)}  \d \tau_m \Bigg(\prod_{\scriptstyle \ell:\sigma_\ell=0 \atop \scriptstyle 1\leq \ell\leq m}\delta_{s_\ell}(\tau_\ell)\Bigg)\notag\\
&\quad \quad \times \int \d \mathcal G^\spin_{\vep;0}  \int \d \mathcal G^{\bs \sigma}_{\vep;1} \cdots \int \d \mathcal G^\spin_{\vep;m}  f\notag\\
&\quad=\int_{\Delta_m(t)}\d \bs u_m\otimes \d \bs v_m\Bigg( \prod_{\scriptstyle \ell:\sigma_\ell=0 \atop \scriptstyle 1\leq \ell\leq m}\delta_{0}(v_\ell)\Bigg)\int \d \mathcal G^\spin_{\vep;0}  \int \d \mathcal G^{\bs \sigma}_{\vep;1}
 \cdots\int \d \mathcal G^\spin_{\vep;m} f,\label{int:P:graphical}
\end{align}
where $\bs u_m\,\defeq\,(u_1,\cdots,u_m)$, $\bs v_m\,\defeq\,(v_1,\cdots,v_m)$, and
 \begin{align}\label{def:Deltam(t)}
\Delta_m(t)\,\defeq\Biggl\{(\bs u_m,\bs v_m)\in (0,t)^m\times [0,t)^m; \sum_{\ell=1}^m (u_\ell+v_\ell)<t\Biggr\}.
\end{align}

\subsubsection{Definition of the graphical integrals}\label{sec:graphical}
Throughout Section~\ref{sec:graphical}, we fix $\vep\in (0,\ovl]$, $m\in \Bbb N$, $\bi_1,\cdots,\bi_m\in \mc E_N$ with $\bi_1\neq \cdots\neq \bi_m$, and $\spin=(\sigma_1,\cdots,\sigma_m)\in \{0,1\}^m$.
Also, $s_1,\cdots,s_m,x_0$ satisfy $0<s_1<\cdots<s_m=t$ and $x_0=(x_0^1,\cdots,x_0^N)\in \R^{2N}$. With $m$ held fixed, it will be convenient to use the following notations as well:
\begin{align}\label{preset}
s_0\,\defeq\,0,\quad s_{m+1}\defeq\, t,\quad \bi_{m+1}\,\defeq\,\varnothing.
\end{align}

The graphical integrals in \eqref{int:P:graphical0} are defined in Step~1--3 below. Roughly speaking,  these steps define the following objects: (1) weighted graphs $\mathcal G^\spin_\vep$ in space and time; (2) subgraphs $\{\mathcal G^{\spin}_{\vep;0}, \ldots, \mathcal G^{\spin}_{\vep;m}\}$ decomposing $\mathcal G^\spin_\vep\setminus\{(x_0^k,s_0)\}_{k=1}^N$; and (3) the graphical integrals according to the subgraphs. Figure~\ref{fig:1} gives an example illustrating the definitions in the first two steps.\medskip 

\begin{figure}[t]
\begin{center}
 \begin{tikzpicture}[scale=1]
    \draw[gray, thin] (0,0) -- (7.6,0);
    \foreach \i in {0.0, 0.5, 2.0, 2.5, 3.2, 3.7, 4.7, 6.4, 7.6} {\draw [gray] (\i,-.05) -- (\i,.05);}
    \draw (0.0,-0.03) node[below]{$0$};
    \draw (0.5,-0.07) node[below]{$s_1$};
    \draw (2.0,-0.07) node[below]{$s_2$};
    \draw (2.5,-0.07) node[below]{$\tau_2$};
    \draw (3.2,-0.07) node[below]{$s_3$};
    \draw (3.7,-0.07) node[below]{$\tau_3$};
    \draw (4.7,-0.07) node[below]{$s_4$};
    \draw (6.4,-0.07) node[below]{$s_5$};
    \draw (7.6,-0.03) node[below]{$t$};
    \draw [line width=2.5pt, color=red!80!white] (0.0,2.30) -- (0.5,2.20);  
    \draw [line width=2.5pt, color=red!80!white] (0.5,2.20) -- (2.0,2.20);  
    \draw [line width=2.5pt, color=red!80!white] (2.0,2.20) -- (2.5,1.80);   
    \draw [line width=2.5pt, color=red!80!white] (2.5,1.80)--(3.2,2.10);   
    \draw [line width=2.5pt, color=red!80!white] (2.5,1.25)--(3.2,1.35);  
    \draw [line width=2.5pt, color=red!80!white] (3.2,1.35)--(3.7,1.05);        
    \draw [line width=2.5pt, color=red!80!white] (3.2,2.10)--(3.7,2.02);      
    \draw [line width=2.5pt, color=red!80!white] (3.7,1.05)--(4.7,2.25);    
    \draw [line width=2.5pt, color=red!80!white] (3.7,2.02)--(4.7,2.35);  
    \draw [line width=2.5pt, color=green!80!black] (0.0,0.40) -- (0.5,0.55);  
    \draw [line width=2.5pt, color=green!80!black] (0.0,1.20) -- (0.5,0.65);  
    \draw [line width=2.5pt, color=blue!70!white] (0.0,1.80) -- (0.5,1.15);    
    \draw [line width=2.5pt, color=blue!70!white] (0.5,0.65) -- (2.0,1.4);  
    \draw [line width=2.5pt, color=blue!70!white] (0.5,1.15) -- (2.0,1.5);  
    \draw [line width=2.5pt, color=purple!60!white] (0.5,0.55) -- (2.0,0.15); 
    \draw [line width=2.5pt, color=purple!60!white] (2.0,0.15) -- (2.5,0.45);  
    \draw [line width=2.5pt, color=purple!60!white] (2.0,1.46)--(2.5,1.19);  
    \draw [line width=2.5pt, color=purple!60!white] (2.5,0.45)--(3.2,0.8);  
    \draw [line width=2.5pt, color=purple!60!white] (2.5,1.15)--(3.2,0.9);  
    \draw [line width=2.5pt, color=orange] (3.7,0.45)--(4.7,1.25);  
    \draw [line width=2.5pt, color=orange] (3.7,0.55)--(4.7,1.75);  
    \draw [line width=2.5pt, color=orange] (4.7,1.25)--(6.4,0.3);  
    \draw [line width=2.5pt, color=orange] (4.7,1.75)--(6.4,0.4);  
    \draw [line width=2.5pt, color=red!80!white] (3.2,0.87)--(3.7,0.49);  
    \draw [line width=2.5pt, color=yellow!90!black] (4.7,2.35)--(6.4,2.40);   
    \draw [line width=2.5pt, color=yellow!90!black] (4.7,2.25)--(6.4,1.15);   
    \draw [line width=2.5pt, color=yellow!90!black] (6.4,0.3)--(7.6,0.45);  
    \draw [line width=2.5pt, color=yellow!90!black] (6.4,0.4)--(7.6,0.85);  
    \draw [line width=2.5pt, color=yellow!90!black] (6.4,1.15)--(7.6,1.55);  
    \draw [line width=2.5pt, color=yellow!90!black] (6.4,2.40)--(7.6,2.15);  
    \node at (0.0,0.40) {$\bullet$};
    \node at (0.0,1.20) {$\bullet$};
    \node at (0.0,1.80) {$\bullet$};
    \node at (0.0,2.30) {$\bullet$};
    \draw (0.0,0.40) node [left] {$x_0^1$};
    \draw (0.0,1.20) node [left] {$x_0^2$};
    \draw (0.0,1.80) node [left] {$x_0^3$};
    \draw (0.0,2.30) node [left] {$x_0^4$};
    \draw [thick, color=black] (0.0,0.40) -- (0.5,0.55);
    \draw [thick, color=black] (0.0,1.20) -- (0.5,0.65);
    \draw [thick, color=black] (0.0,1.80) -- (0.5,1.15);    
    \draw [thick, color=black] (0.0,2.30) -- (0.5,2.20);  
    \node at (0.5,0.55) {$\bullet$};
    \node at (0.5,0.65) {$\bullet$};
    \node at (0.5,1.15) {$\bullet$};
    \node at (0.5,2.20) {$\bullet$};
    \draw [thick, color=black] (0.5,0.65) -- (2.0,1.4);
    \draw [thick, color=black] (0.5,1.15) -- (2.0,1.5);    
    \draw [thick, color=black] (0.5,2.20) -- (2.0,2.20);    
    \draw [thick, color=black] (0.5,0.55) -- (2.0,0.15);    
    \node at (2.0,0.15) {$\bullet$};
    \node at (2.0,1.4) {$\bullet$};
    \node at (2.0,1.5) {$\bullet$};
    \node at (2.0,2.20) {$\bullet$};
    \draw [thick, color=black] (2.0,0.15) -- (2.5,0.45);
    \draw [thick, color=black, snake=coil, segment length=4pt] (2.0,1.44) -- (2.5,1.20);    
    \draw [thick, color=black] (2.0,2.20) -- (2.5,1.80);    
    \node at (2.5,0.45) {$\bullet$};
    \node at (2.5,1.25) {$\bullet$};
    \node at (2.5,1.15) {$\bullet$};
    \node at (2.5,1.80) {$\bullet$};
    \draw [thick, color=black] (2.5,0.45)--(3.2,0.8);  
    \draw [thick, color=black] (2.5,1.15)--(3.2,0.9);    
    \draw [thick, color=black] (2.5,1.25)--(3.2,1.35);  
    \draw [thick, color=black] (2.5,1.80)--(3.2,2.10);        
    \node at (3.2,0.8) {$\bullet$};
    \node at (3.2,0.9) {$\bullet$};
    \node at (3.2,1.35) {$\bullet$};
    \node at (3.2,2.10) {$\bullet$};
    \draw [thick, color=black, snake=coil, segment length=4pt] (3.2,0.85)--(3.7,0.50) ;
    \draw [thick, color=black] (3.2,1.35)--(3.7,1.05);    
    \draw [thick, color=black] (3.2,2.10)--(3.7,2.02);  
    \node at (3.7,0.45) {$\bullet$};
    \node at (3.7,0.55) {$\bullet$};
    \node at (3.7,1.05) {$\bullet$};
    \node at (3.7,2.02) {$\bullet$};
    \draw [thick, color=black] (3.7,0.45)--(4.7,1.25);
    \draw [thick, color=black] (3.7,0.55)--(4.7,1.75);    
    \draw [thick, color=black] (3.7,1.05)--(4.7,2.25);    
    \draw [thick, color=black] (3.7,2.02)--(4.7,2.35);    
    \node at (4.7,1.25) {$\bullet$};
    \node at (4.7,1.75) {$\bullet$};
    \node at (4.7,2.25) {$\bullet$};
    \node at (4.7,2.35) {$\bullet$};
    \draw [thick, color=black] (4.7,1.25)--(6.4,0.3);
    \draw [thick, color=black] (4.7,1.75)--(6.4,0.4);   
    \draw [thick, color=black] (4.7,2.25)--(6.4,1.15); 
    \draw [thick, color=black] (4.7,2.35)--(6.4,2.40); 
    \node at (6.4,0.3) {$\bullet$};
    \node at (6.4,0.4) {$\bullet$};
    \node at (6.4,1.15) {$\bullet$};
    \node at (6.4,2.40) {$\bullet$};
    \draw [thick, color=black] (6.4,0.3)--(7.6,0.45);
    \draw [thick, color=black] (6.4,0.4)--(7.6,0.85);
    \draw [thick, color=black] (6.4,1.15)--(7.6,1.55);    
    \draw [thick, color=black] (6.4,2.40)--(7.6,2.15);
    \node at (7.6,0.45) {$\bullet$};
    \node at (7.6,0.85) {$\bullet$};
    \node at (7.6,1.55) {$\bullet$};
    \node at (7.6,2.15) {$\bullet$};
\end{tikzpicture}
\end{center}
\vspace{-.5cm}
\caption{This figure illustrates the unweighted graph $\mathcal G_{\oslash}^\spin$ associated with $Q^{\lambda,\phi;\bi_1,\bi_2,\bi_3,\bi_4,\bi_5;\spin}_{\vep;s_1,s_2,s_3,s_4,s_5,t}f(x_0)$ for the case of $N=4$, where $x_0^k$’s are distinct, $\bi_1=(2,1)$, $\bi_2=(3,2)$, $\bi_3=(2,1)$, $\bi_4=(4,3)$, $\bi_5=(2,1)$, and $\spin=(0,1,1,0,0)$. A pair of joined bullet points represents a vertex of the form $(x^{i_\ell\prime}_\ell, x^{i_\ell}_\ell, s_\ell)$ or $(z'_\ell, z_\ell, \tau_\ell)$, while any other bullet point indicates a vertex of another type described in \eqref{def:vertex}. The figure highlights six subgraphs, each marked with a different color (see the electronic version of the paper). This coloring convention is used in subsequent figures as well.}
\label{fig:1}
\end{figure}

\noindent {\bf Step~1 (Weighted graphs $\mathcal G^\spin_{\vep}$).}
\begin{defi}[Vertices of $\mathcal G^\spin_\vep$]
The vertex set $\mathcal V^\spin_\vep$ of $\mathcal G^\spin_\vep$ consists of the following points:
\begin{subequations}
\label{def:vertex}
\begin{align}
&(x_0^k,s_0)=(x_0^k,0),\;\mbox{for }1\leq k\leq N;\label{def:vertex1}\\
&(x_\ell^{i_\ell\prime},  x_\ell^{i_\ell},s_\ell),\; (x_\ell^k,s_\ell),\; \mbox{for }1\leq \ell\leq m, \;k\notin \bi_\ell;\label{def:vertex2}\\
& ( z'_\ell, z_\ell,\tau_\ell),\; (z_\ell^k,\tau_\ell),\; \mbox{for $1\leq\ell\leq m$ such that }\sigma_\ell=1,\;  k\notin \bi_\ell; \label{def:vertex3}\\
&(x_{m+1}^k,s_{m+1})=(x_{m+1}^k,t),\;\mbox{for } 1\leq k\leq N,\label{def:vertex4}
\end{align}
\end{subequations}
where $x_\ell^k,z_\ell',z_\ell, z_\ell^k$, for all $1\leq \ell\leq m+1$ and $k$, are \emph{distinct} $\R^2$-variables such that they are also different from $x_0^k$, $1\leq k\leq N$, and $s_\ell,\tau_\ell$ are time variables subject to the following condition:
\begin{linenomath*}\begin{align}\label{def:stau-convention}
\begin{split}
0<s_1\leq \tau_1<s_2\leq\tau_2<\cdots<s_m\leq \tau_m<t\\
\mbox{with }s_\ell<\tau_\ell\mbox{ if $\sigma_\ell=1$ and $s_\ell=\tau_\ell$ if $\sigma_\ell=0$.}
\end{split}
\end{align}\end{linenomath*}
We allow $x_0^{k}=x_0^\ell$ for some $k\neq \ell$, in which case $(x_0^k,s_0)$ and $(x_0^\ell,s_0)$ are identified as the same vertex.  
Also, we call $x_0^k,x_\ell^{i_\ell\prime},x_\ell^{i_\ell},x_\ell^k,z_\ell',z_\ell,z_\ell^k,x_{m+1}^k$ in \eqref{def:vertex} the {\bf spatial components} of the vertices, and ${\sf Space}(\mathsf v)$ denotes the set of all the spatial components of a vertex $\mathsf v$. The {\bf temporal components} of the vertices refer to $s_0,s_\ell,\tau_\ell,t$ in \eqref{def:vertex}, and ${\sf Time}(\mathsf v)$ denotes the temporal component of $\mathsf v$. \hfill $\blacklozenge$ 
\end{defi}

To define the weights of edges, we use the following notations, where
 $\bs w=(w^1,\cdots,w^N)\in \R^{2N}$, $\bi=(i\prime,i)\in \mathcal E_N$, $W$ is a two-dimensional standard Brownian motion, and $0\leq \tau_1<\tau_2<\infty$:
\begin{align}
\widehat{\sigma}_\ell&\,\defeq\, 1-\sigma_\ell,\label{dualsigma}\\
\label{unitary}
w^{\bi\prime}&\,\defeq\,\frac{w^{i\prime }+w^i}{\two},\; 
w^{\bi}\,\defeq\,\frac{w^{i\prime }-w^i}{\two}\;\Longrightarrow\; w^{i\prime}=\frac{w^{\bi\prime }+w^{\bi}}{\two},\; 
w^{i}=\frac{w^{\bi\prime }-w^{\bi}}{\two},\\
k/\bi &\,\defeq 
\begin{cases}
i, &k=i\prime,\\
k, &k\in\{1,\cdots,N\}\setminus\{i\prime\},
\end{cases}\label{def:slash}\\
k^{/\varnothing}&\,\defeq\, k, \quad \forall\;k\in \{1,\cdots,N\},\label{def:wempty}\\
\bs w^{/\bj}&\,\defeq\, (w^{1/\bj},\cdots,w^{N/\bj})\in \R^{2N},\quad \forall\;\bj\in \{\bi,\varnothing\},\label{def:wslash}\\
P_{\tau_1,\tau_2}(x,\widetilde{y})&\,\defeq\, P_{\tau_2-\tau_1}(\widetilde{y}-x),\label{def:Pdelta}\\
\label{def:svep}
\s^{\lambda,\phi}_\vep(\tau_1,\tau_2;x,\widetilde{y})&\,\defeq\, \E_{\vep x}^{W}\left[ \lv^2\exp\left\{\lv\int_0^{\tau_2-\tau_1} \d r\varphi_\vep(W_r) \right\} \delta_{ \vep \widetilde{y}}(W_{\tau_2-\tau_1})\right].
\end{align}
The variable $x$ in \eqref{def:Pdelta} and \eqref{def:svep} will be called the {\bf initial states} of the heat kernel and $\s^{\lambda,\phi}_\vep$; the variable $\widetilde{y}$ will be called their {\bf terminal states}. Note that $\s^{\lambda,\phi}_\vep(\tau_1,\tau_2;x,\cdot)$ is understood as a density incorporating the informal effect of $\int \delta_{\vep \widetilde{y}}(W_\tau)F(\widetilde{y})\d\widetilde{y}=\vep^{-2}F(\vep^{-1}W_\tau)$, so we may assume $\s^{\lambda,\phi}_\vep>0$ pointwise. We will use this density only through the integrals.

Given $\bi_1,\ldots,\bi_m \in \mc E_N$ with $\bi_1 \neq \cdots \neq \bi_m$, $\spin = (\sigma_1, \ldots, \sigma_m) \in \{0,1\}^m$, and the vertices defined in \eqref{def:vertex}, we now specify the weights for the weighted graph in question. Roughly speaking, these weights reflect the following interactions, using the kernels in  \eqref{def:Pdelta} and \eqref{def:svep} such that $\s^{\lambda,\phi}_\vep$ governs the interactions: when $\sigma_\ell = 1$, the two particles indexed by $i_\ell'$ and $i_\ell$ are required to ``almost contact," in the sense of approaching each other within a distance of order $\vep$, at the times $s_\ell$ and $\tau_\ell$ and possibly at many other points within the interval $[s_\ell, \tau_\ell)$; when $\sigma_\ell = 0$, it is only required that such a near contact occur at time $s_\ell$. For more details on how these interactions arise in the expectation in \eqref{P:vepseries}, see \cite[Section~4, especially Section~4.1]{C:DBG}.

Specifically, the weights formalizing the above interpretation are defined below, collectively providing a dynamical picture of interacting particles that experience almost-contact interactions when viewed sequentially in increasing order of $1 \leq \ell \leq m$:
\begin{align}
{\mathsf w}_{0,1;k}&\,\defeq\,
 \left\{
 \arraycolsep=1pt\def\arraystretch{1.2}
\begin{array}{lr}\label{weight1}
 P_{0,s_1}(x_0^{i_1\prime},x_1^{i_1}+\vep\two x^{\bi_1}_1),&\hspace{4.35cm}k=i_1\prime,\\
  P_{0,s_1}(x_0^{i_1},x_1^{i_1}),&\hspace{4.35cm} k=i_1,\\
  P_{0,s_1}(x_0^{k},x_1^{k}),& 
 \hspace{4.35cm} k\in\{1,\cdots,N\}\setminus \bi_1,
\end{array}
\right.\\
\label{weight2}
 {\mathsf w}_{\ell,\ell.5;k}&\,\defeq\,
 \left\{
 \arraycolsep=1pt\def\arraystretch{1.2}
\begin{array}{lr}
\sigma_\ell 
\begin{bmatrix}
P_{s_\ell,\tau_{\ell}}(\two x_\ell^{i_\ell}+\vep x_\ell^{\bi_\ell},z'_\ell)\\
\vspace{-.3cm}\\
\s^{\lambda,\phi}_\vep (s_\ell,\tau_\ell; x_\ell^{\bi_{\ell}},z_\ell)
\end{bmatrix}_\times, &\hspace{3.3cm} k=i_\ell.5,\\
 \sigma_\ell 
 P_{s_\ell,\tau_\ell}(x_\ell^k,z_\ell^k),
 & \hspace{3.3cm} k\in \{1,\cdots,N\}\setminus \bi_\ell,
 \end{array}
 \right.\\
{\mathsf w}_{\ell.5,\ell+1;k}&\,\defeq\,
 \left\{
 \arraycolsep=1pt\def\arraystretch{1.6}
\begin{array}{lr}\label{weight3}
\displaystyle  
\sigma_\ell  P_{\tau_\ell,s_{\ell+1}}\Big(\frac{z_\ell'+\vep z_\ell}{\two},x_{\ell+1}^{i_\ell\prime/\bi_{\ell+1}}+\1_{\big\{{i_\ell\prime=i_{\ell+1}\prime,\atop 1\leq \ell\leq m-1}
\big\}}\vep\two x_{\ell+1}^{\bi_{\ell+1}}
\Big), &\hspace{-.5cm}k=i_\ell\prime,\\
\displaystyle \sigma_\ell  P_{\tau_\ell,s_{\ell+1}}\Big(\frac{z_\ell'-\vep z_\ell}{\two},x_{\ell+1}^{i_\ell/\bi_{\ell+1}}+\1_{\big\{{i_\ell=i_{\ell+1}\prime,
\atop 1\leq \ell\leq m-1
}
\big\}}\vep\two x_{\ell+1}^{\bi_{\ell+1}}
\Big),&\hspace{-.5cm}k=i_\ell,\\ 
\sigma_\ell  P_{\tau_\ell,s_{\ell+1}}\Big(z_\ell^{k},x_{\ell+1}^{k/\bi_{\ell+1}}+\1_{\big\{{k=i_{\ell+1}\prime,\atop 1\leq\ell\leq m-1}\big\}}\vep\two x_{\ell+1}^{\bi_{\ell+1}}
\Big),&\hspace{-.5cm} k\in \{1,\cdots,N\}\setminus \bi_\ell,
\end{array}
\right.\\
{\mathsf w}_{\ell,\ell+1;k}&\,\defeq\,
 \left\{
  \arraycolsep=1pt\def\arraystretch{1.6}
\begin{array}{lr}\label{weight4}
\displaystyle
\widehat{\sigma}_\ell  \lv^{1/2} P_{s_\ell,s_{\ell+1}}\Big(x_\ell^{i_\ell}+\vep\two x_\ell^{\bi_\ell},x_{\ell+1}^{i_\ell\prime/\bi_{\ell+1}}+\1_{\big\{{i_\ell\prime=i_{\ell+1}\prime,\atop 1\leq \ell\leq m-1}\big\}}\vep\two x_{\ell+1}^{\bi_{\ell+1}}\Big),&\hspace{-1cm} k=i_\ell\prime,\\
\displaystyle
\widehat{\sigma}_\ell \lv^{1/2} P_{s_\ell,s_{\ell+1}}\Big(x_\ell^{i_\ell},x_{\ell+1}^{i_\ell/\bi_{\ell+1}}+\1_{\big\{{i_\ell=i_{\ell+1}\prime,\atop 1\leq \ell\leq m-1}\big\}}\vep\two x_{\ell+1}^{\bi_{\ell+1}}\Big),&\hspace{-1cm} k=i_\ell,\\
\widehat{\sigma}_\ell  P_{s_\ell,s_{\ell+1}}\Big(x_\ell^{k},x_{\ell+1}^{k/\bi_{\ell+1}}+\1_{\big\{{k=i_{\ell+1}\prime,\atop 1\leq \ell\leq m-1}\big\}}\vep\two x_{\ell+1}^{\bi_{\ell+1}}\Big),
 &\hspace{-2cm} k\in \{1,\cdots,N\}\setminus \bi_\ell,
\end{array}
\right.
\end{align}
where $x_{m+1}^{k/\bi_{m+1}}=x_{m+1}^k$ by \eqref{preset} an \eqref{def:wempty}, and $[\cdot]_\times$ denotes a multiplication column, as defined in \eqref{def:column}. 

Let us explain the definitions of the weights in more detail. First, all the weights in \eqref{weight2} and \eqref{weight3} are nonzero if and only if $\sigma_\ell=1$ as  $\s^{\lambda,\phi}_\vep>0$, and all the weights in \eqref{weight4} are nonzero if and only if $\sigma_\ell=0$. So these two subcollections of weights complement each other. 
Second, in each of \eqref{weight3} and 
\eqref{weight4}, $x_{\ell+1}^{k/\bi_{\ell+1}}$ takes into account the possible interaction at time $s_{\ell+1}$, and one and only one indicator function takes the value $1$. Furthermore,
the heat kernels in \eqref{weight1}--\eqref{weight4} show $\mathcal O(\vep)$-perturbations in the initial and terminal states.
 Without these perturbations, some of the states become identical. For example, the heat kernels in the weights labelled by $i_\ell\prime$ and $i_\ell$ in \eqref{weight3} use the same initial condition $z_\ell'/\two$ when the perturbations are removed.  
 
 \begin{defi}
 (1$\cc$) By writing ``$\mathsf w:\mathsf v_0\edge \mathsf v_1$,'' we draw an edge for a nonzero weight $\mathsf w$ as a straight line segment to link the two vertices $\mathsf v_0$ and $\mathsf v_1$. By ``$\mathsf w:\mathsf v_0\cedge \mathsf v_1$,'' a coiled line segment is drawn instead.
\medskip 
 
\noindent (2$\cc$) Given a nonzero weight $\mathsf w$ as defined in \eqref{weight1}--\eqref{weight4}, we define ${\sf Initial}(\mathsf w)$ as the subset of $\bigcup_{\mathsf v \in \mathcal V^\spin_\vep} {\sf Space}(\mathsf v)$ that specifies the initial states for all heat kernels and any $\s^{\lambda,\phi}_{\vep}$ involved in $\mathsf w$. Similarly, ${\sf Terminal}(\mathsf w)$ denotes the corresponding subset determined by the terminal states.   \hfill $\blacklozenge$ 
 \end{defi}
 
 \begin{defi}[Weighted edges of $\mathcal G^\spin_\vep$]
For any nonzero weight $\mathsf w$ from \eqref{weight1}--\eqref{weight4},
we draw a unique edge for $\mathsf w$ according to the following rule:
\begin{align*}
\mathsf  w:
\begin{cases}
\mathsf  v_0\cedge \mathsf  v_{1}, &\mbox{if $\mathsf w= {\mathsf w}_{\ell,\ell.5;k}$ for some $1\leq \ell\leq m$, and $k=i_\ell.5$ (cf. \eqref{weight2})},\\
\mathsf  v_0\edge \mathsf  v_{1},&\mbox{otherwise,}
\end{cases}
\end{align*}
where $\mathsf v_0,\mathsf v_1$ are the unique vertices satisfying ${\sf Initial}(\mathsf w)\cap {\sf Space}(\mathsf  v_0)\neq \varnothing$ and ${\sf Terminal}(\mathsf  w)\cap {\sf Space}(\mathsf  v_{1})\neq \varnothing$.
 \hfill $\blacklozenge$ 
\end{defi}

The uniqueness of $(\mathsf v_0, \mathsf v_1)$ follows from our criterion for distinguishing vertices, as specified below \eqref{def:vertex}, and our rule of identifying $(x_0^k,s_0)$ and $(x_0^\ell,s_0)$ as the same vertex when $x_0^k=x_0^\ell$, as mentioned below \eqref{def:stau-convention}.  
 Consequently, there is at most one edge connecting any two distinct vertices with temporal components in $[s_1, t]$. In contrast, multiple edges can connect two vertices whose temporal components lie in $[0, s_1]$. For example, in the graph over $[0, t]$ in Figure~\ref{fig:2}, the two green edges illustrate these multiple edges over $[0, s_1]$. This example is essentially the only case in which multiple edges can occur, as made precise in the following lemma. This scenario will be central to our analysis of the mixed moments. 

\begin{lem}
There are multiple edges in $\mathcal G_\vep^\spin$ if and only if $x_0^{i_1\prime}=x_0^{i_1}$. In this case, the multiple edges occur in linking $(x_0^{i_1\prime},s_0)=(x_0^{i_1},s_0)$ and $(x_1^{i_1\prime},x_1^{i_1},s_1)$; there are exactly two such edges. 
\end{lem}

Up to this point, let us clarify how $\mathcal G_\vep^\spin$ depend on the underlying variables by writing
\[
\mathcal G_\vep^\spin 
=\mathcal G^{\bi_1,\cdots,\bi_m;\spin}_{\vep;x_0,\bs x_{m},x_{m+1}^1,\cdots,x_{m+1}^N,\bs z}(\bs u_{m},u_{m+1},\bs v_m).
\]
Here, $\bs x_{m}\defeq \bigcup_{\mathsf  v}{\sf Space}(\mathsf  v)$ for $\mathsf  v$ ranging over the vertices in \eqref{def:vertex2}, $\bs z\,\defeq \bigcup_{\mathsf  v}{\sf Space}(\mathsf  v)$  for $\mathsf  v$ ranging over the vertices in \eqref{def:vertex3}, and we use the notations $\bs u_m,u_{m+1},\bs v_m$ in \eqref{def:uv123} and \eqref{def:Deltam(t)}. Also, $\mathcal G_\oslash^\spin$ denotes
the unweighted graph associated with $\mathcal G_\vep^\spin$. \medskip 

\noindent {\bf Step~2 (Subgraph decomposition).}
In this step, we introduce a decomposition of $\mathcal{G}_\vep^{\spin}\setminus\{(x_0^k,0)\}_{k=1}^N$ into subgraphs denoted by
\begin{linenomath*}\begin{align}\label{def:Gell}
\mathcal G^{\spin}_{\vep;\ell}=\mathcal G^{\bi_1,\cdots,\bi_m;\spin}_{\vep;\ell;x_0,\bs x_{m},x_{m+1}^1,\cdots,x_{m+1}^N,\bs z}(\bs u_{m},u_{m+1},\bs v_m),\quad \ell\in \{0,1,\cdots,m\}.
\end{align}\end{linenomath*}
Each subgraph $\mathcal{G}^{\spin}_{\vep;\ell}$ is to be defined via its \emph{unweighted} counterpart, denoted $\mathcal G_{\oslash;\ell}^\spin$; the edge weights for $\mathcal G_{\vep;\ell}^\spin$ are then inherited from the corresponding edges in $\mathcal G_\vep^\spin$. Moreover, our construction focuses exclusively on specifying the edge sets of $\mathcal G_{\oslash;\ell}^\spin$, as the vertex set $\mathcal V^\spin_\ell$ of $\mathcal G_{\oslash;\ell}^\spin$ consists of all vertices $\mathsf v$ in $\mathcal G_\oslash^\spin$ for which there is an edge in $\mathcal G_{\oslash;\ell}^\spin$ adjacent to $\mathsf v$ from the left of $\mathsf v$. Here and throughout, we interpret \emph{``left'' and ``right''} with respect to increasing time, from left to right, as shown in Figure~\ref{fig:1}. 

To specify $\{\mathcal{G}^{\spin}_{\oslash;\ell}\}_{\ell=0}^m$, we begin by introducing the following terminology on paths and related notions that can be drawn readily from $\mathcal{G}_{\oslash}^{\spin}$. 

\begin{defi}\label{def:path}
{\rm (1$\cc$)} A {\bf  path} $\mathcal P$ in $\mathcal G_\oslash^\spin$ is a finite sequence of distinct edges such that (i) one is adjacent to the next and (ii) any given vertex of $\mathcal G_\oslash^\spin$ is adjacent to at most two edges of $\mathcal P$. Moreover, if the set of ${\sf Time}(\mathsf  v)$, defined below \eqref{def:stau-convention}, for all vertices $\mathsf v$ adjacent to edges of $\mathcal P$ can be linearly ordered and $t_{\min}$ and $t_{\max}$ denote the minimal and maximum of these ${\sf Time}(\mathsf  v)$, then $\mathcal P$ is said to begin at $t_{\min}$ and end at $t_{\max}$. Moreover, if we specify when a path ends, it implicitly implies that no two edges in the path span the same time interval. \medskip

\noindent {\rm (2$\cc$)} A path is {\bf  entanglement-free} if it contains neither (i) an edge represented by a coiled line segment nor (ii) a pair of edges adjacent to a vertex $\mathsf{v} = (x_\ell^{i_\ell'}, x_\ell^{i_\ell}, s_\ell)$, with one edge on the left and the other on the right of $\mathsf{v}$. \medskip

\noindent {\rm (3$\cc$)} A path is {\bf  left-maximal} if it cannot be extended to an entanglement-free path by adding more edges drawn as straight line segments over smaller time intervals.\medskip

\noindent {\rm (4$\cc$)} A {\bf  usual path} is a path such that (i) it begins at time $0$ and ends at time $t$ and (ii) none of its edges are adjacent to a vertex of the form $(x_\ell^{i_\ell\prime},x_\ell^{i_\ell},s_\ell)$ or to a vertex of the form $ ( z'_\ell, z_\ell,\tau_\ell)$.\medskip

\noindent {\rm (5$\cc$)} The {\bf  duration} of an edge adjacent to vertices $\mathsf  v_0$ and $\mathsf  v_1$ is $\lvert {\sf Time}(\mathsf  v_0)-{\sf Time}(\mathsf  v_1)\rvert$. 
The duration of a path $\mathcal P$, denoted by $\lvert \mathcal P\rvert$, is the sum of the durations of its edges.\medskip

\noindent {\rm (6$\cc$)} Given a subgraph $\mathcal G'$, a vertex $\mathsf  v$ is a {\bf  left end-vertex} of $\mathcal G'$ if $\mathsf  v$ is adjacent to $\mathsf  e$ for some edge $\mathsf  e$ of $\mathcal G'$, but there is no edge $\mathsf  e'$ of $\mathcal G'$ such that $\mathsf  v$ is adjacent to $\mathsf  e'$ and $\mathsf  v$ defines the right endpoint of $\mathsf  e'$. Also, an edge of $\mathcal G'$ is a {\bf  left end-edge} of $\mathcal G'$ if it is adjacent to a left end-vertex of $\mathcal G'$. The {\bf  right end-vertices} and {\bf  right end-edges} of $\mathcal G'$ are similarly defined. Moreover, since these definitions of left and right end-vertices use only edges of $\mathcal G'$, they allow straightforward extensions to any given set of edges. 
\hfill $\blacklozenge$ 
\end{defi} 

The following properties are implicit in the above definitions: First, the definition of paths in (1$\cc$) does not preclude the existence of two edges traversing the same time interval. For example, Figure~\ref{fig:1} depicts two paths, marked green and blue, that stretch backward at times $s_1$ and $s_2$, respectively. Also, Figure~\ref{fig:1} shows two edges marked yellow such that they span $[s_5,t]$ and are adjacent to the same vertex at time $s_5$. These two edges thus define a path as well. Moreover, these paths in Figure~\ref{fig:1} contrast the usual paths defined in (4$\cc$), as the latter paths now obey the rule in (1$\cc$) that two edges within a path spanning the same time interval are prohibited if we specify when the path ends. Finally, an entanglement-free path is equivalent to a path that contains no edges represented as colied line segments and ``deflects''  at any vertex in \eqref{def:vertex2}--\eqref{def:vertex3} represented by a triplet.

\begin{defi}[Subgraphs $\{\mathcal{G}^{\spin}_{\oslash;\ell}\}_{\ell=0}^m$]\label{def:subgraphedge}
{\rm (1$\cc$)} The edge set of $\mathcal G_{\oslash;0}^\spin$ consists of the two unique edges whose right end-vertices coincide at time $s_1$. \medskip 

\noindent 
{\rm (2$\cc$)} For any $\ell=1,\cdots,m-1$, the edge set of $\mathcal G^\spin_{\oslash;\ell}$ consists of the following edges: 
\begin{itemize}
\item  the edges of the two unique left-maximal entanglement-free paths, denoted by $\mathcal P'_\ell$ and $\mathcal P_{\ell}$, such that they end at time $s_{\ell+1}$ and meet at the vertex $(x_{\ell+1}^{i_{\ell+1}\prime},x_{\ell+1}^{i_{\ell+1}},s_{\ell+1})$;
\item if $\sigma_\ell=1$, an additional edge given by the unique coiled line segment over $[s_{\ell},\tau_{\ell}]$. 
\end{itemize}
Moreover, by concatenating the paths $\mathcal P'_\ell$ and $\mathcal P_\ell$, we obtain an entanglement-free path, which we call the {\bf spine} of $\mathcal G^\spin_{\oslash;\ell}$ and denote by $\mathcal P'_\ell \oplus \mathcal P_\ell$. \medskip 

\noindent 
{\rm (3$\cc$)} The edge set of $\mathcal G^\spin_{\oslash;m}$ consists of the following edges:
\begin{itemize}
\item the edges of all the left-maximal entanglement-free paths ending at time $t$;
\item if $\sigma_m=1$,  an additional edge given by the unique coiled line segment over $[s_m,\tau_m]$.\hfill $\blacklozenge$
\end{itemize} 
\end{defi}

Let us highlight several properties of these subgraphs to facilitate later applications. First, note that the edge sets for $\ell = 1, 2, \ldots, m$ satisfy the following: 
\begin{itemize}
\item For $\ell=1,2,\cdots,m-1$, the durations $\lvert \mathcal P'_\ell\rvert$ and $\lvert\mathcal P_\ell\rvert$ are bounded below by $s_{\ell+1}-\tau_\ell$; recall \eqref{def:stau-convention}. 
Lemma~\ref{lem:path} below extends this observation by using the assumption $\bi_{\ell}\neq \bi_{\ell+1}$ and will play a key role in the forthcoming proofs.

\begin{lem}[Lemma~4.9 of \cite{C:DBG}]\label{lem:path}
For any fixed $\ell=1,\cdots,m-1$, at least one of the paths $\mathcal P_\ell'$ and $\mathcal P_\ell$ begins at time $\tau_{\ell-1}$ or a smaller time (with $s_0=\tau_0=0$ for $\ell=1$). Then $\lvert\mathcal P_\ell'\oplus \mathcal P_\ell\rvert$  is bounded below by $2(s_{\ell+1}-\tau_{\ell})+(s_\ell-\tau_{\ell-1})=2u_{\ell+1}+u_\ell$, where the relationship between $s_\ell$ and $\tau_\ell$ in  \eqref{def:stau-convention} is used and
$u_\ell$ and $u_{\ell+1}$ are defined in \eqref{def:uv123}.
\end{lem}

\item For $\ell=m$, since every usual path is a left-maximal entanglement-free path ending at time $t$ by definition, all the edges of the usual paths are edges of  $\mathcal G^\spin_{\oslash;m}$. In particular, there are no usual paths if $\bi_1\cup\cdots \cup \bi_m=\{1,\cdots,N\}$. 
\end{itemize}
The next lemma concludes this step by confirming that $\{\mathcal G^\spin_{\vep;\ell}\}_{\ell=0}^m$ decomposes $\mathcal G^\spin_{\vep}$.

\begin{lem}\label{lem:edgesubG}
The edge set of $\mathcal G^\spin_{\vep}$ is the disjoint union of the edge sets of $\{\mathcal G^\spin_{\vep;\ell}\}_{\ell=0}^m$; the vertex set of $\mathcal{G}^{\spin}_{\vep}$ is the disjoint union of $\{(x_0^k,0)\}_{k=1}^N$ and the vertex sets of $\{\mathcal{G}^{\spin}_{\vep;\ell}\}_{\ell=0}^m$. 
\end{lem}
\begin{proof}
The required property about edge sets has been established in \cite[Lemma~4.10]{C:DBG}. Consequently, the property for vertex sets holds if no two subgraphs share a vertex. To see this, it suffices to note that subgraphs are connected only via the vertices in \eqref{def:vertex} represented as triplets, while each such vertex belongs to one and only one subgraph by the definition of $\mathcal V^\spin_\ell$ given below \eqref{def:Gell}. 
\end{proof}

\noindent {\bf Step~3 (Graphical integrals).}
The {\bf graphical integrals} refer to the unweighted integrals $\int \d \mathcal V_{\ell}^\spin F$ and the weighted integrals $\int \d \mathcal G_{\vep;\ell}^\spin F$, which are defined as follows: 
\begin{subequations}\label{def:graphicalintegral}
\begin{align}
\int \d \mathcal V_{\ell}^\spin F&\,\defeq \int_{\R^{2\lvert\mathcal I_\ell\rvert}}\d \mu_\ell(a_i;i\in \mathcal I_\ell) F(a_i;i\in \mathcal I_\ell),\label{def:graphicalintegral1}\\
\int \d \mathcal G_{\vep;\ell}^\spin F&
\,\defeq\int \d \mathcal V_{\ell}^\spin \Pi^{\spin}_{\vep;\ell}F.\label{def:graphicalintegral2}
\end{align}\end{subequations}
Here, $\{a_i;i\in \mathcal I_\ell\}=\bigcup\{{\sf Space}(\mathsf  v);\mathsf  v\in \mathcal V_{\ell}^\spin\}$ with $a_i\in \R^2$, 
\begin{linenomath*}\begin{align}\label{def:Pell}
\Pi^{\spin}_{\vep;\ell}=\Pi^{\lambda,\phi;\bi_1,\cdots,\bi_m;\spin}_{\vep;\ell;\bs x_{m},x_{m+1}^1,\cdots,x_{m+1}^N,\bs z}(\bs u_{m},u_{m+1},\bs v_m)\defeq \prod\{\mbox{edge weights in $\mathcal G_{\vep;\ell}^\spin$}\},
\end{align}\end{linenomath*}
and $\d\mu_\ell=\d \mu_\ell(a_i;i\in \mathcal I_\ell)$ is the product of all measures 
from the following list that use $a_i$ for some $i\in \mathcal I_\ell$ as one variable of integration:
\begin{subequations}
\label{def:productmeasure}
\begin{align}
&\eqspace\d x_\ell^{i_\ell\prime}\otimes \d  x_\ell^{i_\ell}\varphi(x_\ell^{\bi_\ell}),\; \d x_\ell^k,\; \mbox{for }1\leq \ell\leq m,\;k\notin \bi_\ell;\label{def:productmeasure1}\\
&\eqspace \d z'_\ell \otimes\d z_\ell \varphi(z_\ell),\; \d z_\ell^k,\; \mbox{for $1\leq\ell\leq m$ such that }\sigma_\ell=1,\; k\notin \bi_\ell; \label{def:productmeasure2}\\
&\eqspace\d x_{m+1}^k,\;\mbox{for } 1\leq k\leq N.\label{def:productmeasure3}
\end{align}
\end{subequations}
For example, given $0\leq \ell\leq m-1$, as $(x_{\ell+1}^{i_{\ell+1}\prime},x_{\ell+1}^{i_{\ell+1}},s_{\ell+1})$ must be a vertex of $\mathcal V_{\ell}^\spin$, the measure $\d x_{\ell+1}^{i_{\ell+1}\prime}\otimes \d  x_{\ell+1}^{i_{\ell+1}}\varphi(x_{\ell}^{\bi_{\ell+1}})$ is part of $\d\mu_\ell$. By the definition of $\mathcal V_\vep^\spin$ in \eqref{def:vertex} and Lemma~\ref{lem:edgesubG}, all of the measures in 
\eqref{def:productmeasure} are used to define the graphical integrals, but none of them are used twice to define $\int \d \mathcal V_{\ell}^\spin F$ and $\int \d \mathcal V_{\ell'}^\spin F'$ for $\ell\neq \ell'$.  
 
\subsubsection{Bounds on the mixed moments}\label{sec:mombdd}
We now begin the proof of Theorem~\ref{thm:mombdd}, focusing primarily on establishing (1$\cc$); the arguments for (2$\cc$) and (3$\cc$) will follow by appropriately modifying certain steps of this proof. To prove Theorem~\ref{thm:mombdd} (1$\cc$), we first introduce two lemmas, stated as Lemmas~\ref{lem:mix2mom} and \ref{lem:FKmom1} below. Lemma~\ref{lem:FKmom1} demonstrates series expansions for the mixed moments in Theorem~\ref{thm:mombdd} (1$\cc$), relying on the series expansions for second moments given in Lemma~\ref{lem:mix2mom} and an extension of the series in \eqref{P:vepseries}. 

For the statement of Lemma~\ref{lem:mix2mom}, recall that $[\cdot]_{\times}$ and $[\cdot]_{\otimes }$ are defined in and below \eqref{def:column}.

\begin{lem}\label{lem:mix2mom}
With $\bi_0\,\defeq\,(2,1)$, it holds that, for all $f\geq0$,  
\begin{align}
&\E^{({B}^1, {B}^2)}_{( {x}^1_0, {x}^2_0)}\Bigg[\exp\Bigg\{\lv \int_0^{ {t}}\d r \varphi_\vep( {B}^{\bi_0}_r)\Bigg\}f( {B}_{ {t}})\Bigg]\label{def:Qvep01000}\\
&\quad=\E^{( {B}^1, {B}^2)}_{( {x}^1_0, {x}^2_0)}[f( {B}_{ {t}})]+\int_0^{ {t}}\d  {s}_1\E^{( {B}^1, {B}^2)}_{( {x}^1_0, {x}^2_0)}\Bigg[\varphi_\vep( {B}^{\bi_0}_{ {s}_1}) \lv \exp\Bigg\{\lv \int_{ {s}_1}^{ {t}} \d r\varphi_\vep( {B}^{\bi_0}_r)\Bigg\} f( {B}_{ {t}}) \Bigg]\notag\\
&\quad =\E^{( {B}^1, {B}^2)}_{( {x}^1_0, {x}^2_0)}[f( {B}_{ {t}})]+\int_0^{ {t}}\d  {s}_1\notag\\
&\quad\quad \E^{( {B}^1, {B}^2)}_{( {x}^1_0, {x}^2_0)}\Bigg[\varphi_\vep( {B}^{\bi_0}_{ {s}_1})\left(\lv+ \int_{ {s}_1}^{{ {t}}}\d  {\tau}_1\lv^2 \exp\Bigg\{\lv \int_{ {s}_1}^{ {\tau}_1} \d r\varphi_\vep( {B}^{\bi_0}_r)\Bigg\}\varphi_\vep(B^{\bi_0}_{ {\tau}_1})\right) f( {B}_{ {t}})\Bigg]\label{def:Qvep0}
\\
&\quad =\E^{( {B}_1, {B}_2)}_{( {x}_0^1, {x}_0^2)}[ {f}( {B}_{ {t}})]+\sum_{ {\spin} \in \{0,1\}}\int_0^t \d s_1\int_{[s_1,t)}\d \tau_1  \Bigg(\prod_{\ell: {\sigma}_\ell=0}\delta_{s_\ell}( {\tau}_\ell)\Bigg)\int \d  {\mathcal G}^{ {\spin}}_{\vep;0}\int \d  {\mathcal G}^{ {\spin}}_{\vep;1}f.\label{graphical22}
\end{align}
\end{lem}
\begin{proof}
The first two equalities follow from the fundamental theorem of calculus, applied to functions of the form $s\mapsto \exp\{\int_s^t a(r)\,\mathrm{d}r\}$ and $\tau\mapsto \exp\{\int_s^\tau a(r)\,\mathrm{d}r\}$. 
To obtain \eqref{graphical22}, it suffices to identify the weights defined in \eqref{weight1}--\eqref{weight4} from the following formulas, which can be deduced from the Markov property of Brownian motion:
\begin{align}
&\int_0^t \d s_1\E^B_{x_0}[\varphi_\vep(B^{\bi_0}_{s_1})\lv f(B^1_t,B^2_t)]\notag\\
&\quad =\int_0^t \d s_1 \int_{\R^4} \begin{bmatrix}
\d x_1^{2}\\
 \d x_1^1
 \end{bmatrix}_{\otimes}
 \varphi(x_1^{\bi_0})\int_{\R^4} \begin{bmatrix}
\d x_2^{2}\\
 \d x_2^1
 \end{bmatrix}_{\otimes}
\begin{bmatrix}
 P_{0,s_1}(x_0^{2},x_1^{1}+\vep\two x^{\bi_0}_1)\\
  P_{0,s_1}(x_0^{1},x_1^{1})
\end{bmatrix}_\times\notag\\
&\quad\quad \times \begin{bmatrix}
  \lv^{1/2} P_{s_1,t}(x_1^{1}+\vep\two x_1^{\bi_0},x_{2}^{2})\\
\displaystyle
 \lv^{1/2} P_{s_1,t}(x_1^{1},x_{2}^{1})
\end{bmatrix}_\times f(x^1_2,x_2^2),\label{2mom:eq1}\\
&\int_0^t\d s_1\int_{s_1}^{t}\d \tau_1 \E^B_{x_0}\left[\varphi_\vep(B^{\bi_0}_s) \lv^2 \exp\Bigg\{\lv \int_{ {s}_1}^{ {\tau}_1} \d r\varphi_\vep( {B}^{\bi_0}_r)\Bigg\} \varphi_\vep(B^{\bi_0}_{\tau_1}) f(B_t)\right]\notag\\
&\quad =\int_0^t\d s_1\int_{s_1}^{t}\d \tau_1 \int_{\R^4} \begin{bmatrix}
\d x_1^{2}\\
 \d x_1^1
 \end{bmatrix}_{\otimes}\varphi(x_1^{\bi_0})\int_{\R^4} \begin{bmatrix}
\d z_1'\\
 \d z_1
 \end{bmatrix}_{\otimes} \varphi(z_1)\int_{\R^4} \begin{bmatrix}
\d x_2^{2}\\
 \d x_2^1
 \end{bmatrix}_{\otimes}
\begin{bmatrix}
 P_{0,s_1}(x_0^{2},x_1^{1}+\vep\two x^{\bi_0}_1)\\
  P_{0,s_1}(x_0^{1},x_1^{1})
\end{bmatrix}_\times\notag\\
&\quad\quad \times
\begin{bmatrix}
P_{s_1,\tau_{1}}(\two x_1^{1}+\vep x_1^{\bi_0},z'_1)\\
\vspace{-.3cm}\\
\s^{\lambda,\phi}_\vep (s_1,\tau_1; x_1^{\bi_{0}},z_1)
\end{bmatrix}_\times
\begin{bmatrix}
 P_{\tau_1,t}(\frac{z_1'+\vep z_1}{\two},x_{2}^{2}
)\\
  P_{\tau_1,t}(\frac{z_1'-\vep z_1}{\two},x_{2}^{1}
)
\end{bmatrix}_\times f(x^1_2,x_2^2).\label{2mom:eq2}
\end{align}
More details of the derivation of \eqref{2mom:eq1} and \eqref{2mom:eq2} can be found in \cite[Lemma~4.4]{C:DBG}. Equation~\eqref{graphical22} now follows by using the definitions in \eqref{def:graphicalintegral}. 
\end{proof}

The next result, Lemma~\ref{lem:FKmom1}, requires combining the definitions from Section~\ref{sec:graphical} for both $N=2$ and $N=4$, since its formulation simultaneously involves a two-particle and a four-particle system. To distinguish between these systems, we will append ``$\;\,\widetilde{\mbox{}}\,\;$'' to the notation associated with a two-particle system. For example, $B^1$ and $B^2$ in \eqref{def:Qvep01000} with $N=2$ will henceforth be rewritten as $\widetilde{B}^1$ and $\widetilde{B}^2$. Let us also adopt the following convention: for $0 < t < t' < \infty$,
\begin{align}\label{2-4}
\begin{split}
&\int \d \widetilde{\mathcal G}^{\widetilde{\spin}}_{\vep;0}\int \d \widetilde{\mathcal G}^{\widetilde{\spin}}_{\vep;1}\int \d \mathcal G^\spin_{\vep;0} \int \d \mathcal G^{\bs \sigma}_{\vep;1} 
 \cdots\int \d \mathcal G^\spin_{\vep;m} f =\int \d \widetilde{\mathcal G}^{\widetilde{\spin}}_{\vep;0}\int \d \widetilde{\mathcal G}^{\widetilde{\spin}}_{\vep;1}F
 \end{split}
\end{align}
for $F$ defined by the following function: 
\[
F:(\widetilde{x}_2^1,\widetilde{x}_2^2)\mapsto \int \d \mathcal G^\spin_{\vep;0} 
 \int \d \mathcal G^{\bs \sigma}_{\vep;1} \cdots\int \d \mathcal G^\spin_{\vep;m} f\mid_{x_0^3=\widetilde{x}^1_2,x_0^4=\widetilde{x}^2_2}
 \]
and $\int \d \widetilde{\mathcal G}^{\widetilde{\spin}}_{\vep;0}\int \d \widetilde{\mathcal G}^{\widetilde{\spin}}_{\vep;1}F$ defined by using $\widetilde{\mathcal G}^{\widetilde{\spin}}_\vep$ over $[0,t'-t]$. Here, as before, the iterated graphical integral $\int \d \mathcal G^\spin_{\vep;0} \int \d \mathcal G^{\bs \sigma}_{\vep;1} \cdots \int \d \mathcal G^\spin_{\vep;m} f$ refers to integration with respect to $\mathcal G^{\spin}_\vep$ over $[0, t]$ with $N = 4$, where the initial state at time $0$ is $(x_0^1, x_0^2, x_0^3, x_0^4)$ and the terminal state at time $t$ is $(x_{m+1}^1, x_{m+1}^2, x_{m+1}^3, x_{m+1}^4)$. In contrast, for $\widetilde{\mathcal G}^{\widetilde{\spin}}_\vep$, the particles begin at $(\wt{x}_0^1, \wt{x}_0^2) \in \R^4$ at time $0$ and evolve to $(\wt{x}_2^1, \wt{x}_2^2) \in \R^4$ at time $t' - t$. 

\begin{lem}\label{lem:FKmom1}
For any $\vep\in (0,\ovl]$, $0< t<t'<\infty$, and $x_0^1,x_0^2,\widetilde{x}_0^1,\widetilde{x}_0^2\in \R^2$, it holds that
\begin{align}
&\lv^2\E[X_\vep(x^1_0,t)X_\vep(x^2_0,t) X_\vep(\widetilde{x}_0^1,t')X_\vep(\widetilde{x}_0^2,t') ]\notag\\
\begin{split}
&\quad =\lv^2\E^{(\widetilde{B}^1,\widetilde{B}^2)}_{(\widetilde{x}_0^1,\widetilde{x}_0^2)}\Biggl[\exp\Biggl\{\lv \int_0^{t'-t} \d r \varphi_\vep(\widetilde{B}^{\bi_0}_r)\Biggr\}\\
&\quad\quad \times \E^{(B^1,B^2,B^3,B^4)}_{(x_0^1,x_0^2,\widetilde{B}^1_{t'-t},\widetilde{B}^2_{t'-t})}\Biggl[\exp\Biggl\{\sum_{\bi\in \mc E_4}\lv \int_0^{t} \d r \varphi_\vep(B^\bi_r)\Biggr\}\prod_{j=1}^4 X_0(B^j_t)\Biggr]\Biggr].\label{eq:FKmom}
\end{split}
\end{align}
Also, with the notation used in \eqref{int:P:graphical} and the function $\widetilde{f}_t(\widetilde{x}^1_2,\widetilde{x}^2_2)\,\defeq\, \E^{(B^1,B^2,B^3,B^4)}_{(x_0^1,x_0^2,\widetilde{x}^1_2,\widetilde{x}^2_2)}[f(B_t)]$
for nonnegative $f$ defined on $(\R^2)^4$, we have 
\begin{align}
&\lv^2\E^{(\widetilde{B}^1,\widetilde{B}^2)}_{(\widetilde{x}_0^1,\widetilde{x}_0^2)}\Biggl[\exp\Bigg\{\lv \int_0^{t'-t} \d r \varphi_\vep(\widetilde{B}^{\bi_0}_r)\Biggr\}\notag\\
&\times\E^{(B^1,B^2,B^3,B^4)}_{(x_0^1,x_0^2,\widetilde{B}^1_{t'-t},\widetilde{B}^2_{t'-t})}\Biggl[\exp\Biggl\{\sum_{\bi\in \mc E_4}\lv \int_0^{t} \d r \varphi_\vep(B^\bi_r)\Biggr\}f(B_t)\Biggr]\Biggr]\notag\\
&\quad =\lv^2\E^{(\widetilde{B}^1,\widetilde{B}^2)}_{(\widetilde{x}_0^1,\widetilde{x}_0^2)}\Biggl[\exp\Biggl\{\lv \int_0^{t'-t} \d r \varphi_\vep(\widetilde{B}^{\bi_0}_r)\Biggr\}\widetilde{f}_t(\wt{B}_{t'-t} )\Biggr]\notag\\
&\quad \quad +\sum_{m=1}^\infty \sum_{\stackrel{\scriptstyle \bi_1,\cdots,\bi_m\in \mathcal E_4}{\bi_1\neq \cdots \neq \bi_m}}\sum_{\spin \in \{0,1\}^m} \lv^2
\E^{(\widetilde{B}^1,\widetilde{B}^2)}_{(\widetilde{x}_0^1,\widetilde{x}_0^2)}\Biggl[\int_{\Delta_m(t)}\d \bs u_m\otimes \d \bs v_m\Biggl(\, \prod_{\scriptstyle \ell:\sigma_\ell=0 \atop \scriptstyle 1\leq \ell\leq m}\delta_{0}(v_\ell)\Biggr)\notag\\
&\quad \quad\times \int \d \mathcal G^\spin_{\vep;0} 
 \cdots\int \d \mathcal G^\spin_{\vep;m} f\mid_{(x_0^3,x_0^4)=(\widetilde{B}^1_{t'-t},\widetilde{B}^2_{t'-t}) }\Biggr]\notag\\
 &\quad \quad +\sum_{m=1}^\infty \sum_{\stackrel{\scriptstyle \bi_1,\cdots,\bi_m\in \mathcal E_4}{\bi_1\neq \cdots \neq \bi_m}}\sum_{\scriptstyle \widetilde{\spin}\in \{0,1\}\atop\scriptstyle \spin \in \{0,1\}^m}\lv^2
\int_0^{t'-t} \d \wt{s}_1\int_{[\wt{s}_1,t'-t)}\d \wt{\tau}_1  \int_{\Delta_m(t)}\d \bs u_m\otimes \d \bs v_m\notag\\
 &\quad \quad\times \Biggl(\,\prod_{\scriptstyle \ell:\widetilde{\sigma}_\ell=0\atop \scriptstyle \ell=1}\delta_{\wt{s}_\ell}(\widetilde{\tau}_\ell)\Biggr) \Biggl(\, \prod_{\scriptstyle \ell:\sigma_\ell=0 \atop \scriptstyle 1\leq \ell\leq m}\delta_{0}(v_\ell)\Biggr)\notag\\
 &\quad \quad \times \int \d \widetilde{\mathcal G}^{\widetilde{\spin}}_{\vep;0}\int \d \widetilde{\mathcal G}^{\widetilde{\spin}}_{\vep;1}\int \d \mathcal G^\spin_{\vep;0}  \int \d \mathcal G^{\bs \sigma}_{\vep;1}
 \cdots\int \d \mathcal G^\spin_{\vep;m} f.\label{mom:series}
\end{align}
\end{lem}

\begin{proof}
To prove \eqref{eq:FKmom}, 
define a Borel measure $\mathbf m^{\vep,N}_{x_0;t}(\d y_0)$ on $\R^{2N}$ by the following equation: 
\begin{align}\label{def:mmom}
\int_{\R^{2N}}\mathbf m^{\vep,N}_{x_0;t}(\d y_0) f(y_0)=\E^{B}_{x_0}\Biggl[\exp\Biggl\{\sum_{\bi\in \mc E_N}\lv \int_0^t \d r \varphi_\vep(B^\bi_r)\Biggr\}f(B_t)\Biggr].
\end{align}
Then for $0<t<t'<\infty$,
 by using the Markov property of $X_\vep$ at time $t$, \eqref{eq:momdual2} with $N=2$, and \eqref{def:mmom} with $N=2$, we get
\begin{align*}
&\lv^2\E[X_\vep(x^1_0,t)X_\vep(x^2_0,t) X_\vep(\widetilde{x}_0^1,t')X_\vep(\widetilde{x}_0^2,t') ]\\
&\quad =\lv^2\E\left[X_\vep(x^1_0,t)X_\vep(x^2_0,t)\int_{\R^4}\mathbf m^{\vep,2}_{\widetilde{x}^1_0,\widetilde{x}^2_0;t'-t}(\d y^3,\d y^4) X_\vep(y^3,t)X_\vep(y^4,t) \right]\\
&\quad =\lv^2\int_{\R^4}\mathbf m^{\vep,2}_{\widetilde{x}^1_0,\widetilde{x}^2_0;t'-t}(\d y^3,\d y^4)\E[X_\vep(x^1_0,t)X_\vep(x^2_0,t) X_\vep(y^3,t)X_\vep(y^4,t) ].
\end{align*}
To finish the proof of \eqref{eq:FKmom}, it suffices to apply \eqref{eq:momdual2} with $N=4$ and
\eqref{def:mmom} with $N=2$ in the same order to the right-hand side. 

Next, to get \eqref{mom:series}, we use \eqref{P:vepseries} and \eqref{int:P:graphical}. Note that all of the summands in \eqref{P:vepseries} are nonnegative. Hence, by the monotone convergence theorem, we can write the left-hand side of \eqref{mom:series}
 as a series of expectations: 
\begin{align}
&\lv^2\E^{(\widetilde{B}^1,\widetilde{B}^2)}_{(\widetilde{x}_0^1,\widetilde{x}_0^2)}\Biggl[\exp\Biggl\{\lv \int_0^{t'-t} \d r \varphi_\vep(\widetilde{B}^{\bi_0}_r)\Biggr\}\E^{(B^1,B^2,B^3,B^4)}_{(x_0^1,x_0^2,\widetilde{B}^1_{t'-t},\widetilde{B}^2_{t'-t})}\Biggl[\exp\Biggl\{\sum_{\bi\in \mc E_4}\lv \int_0^{t} \d r \varphi_\vep(B^\bi_r)\Biggr\}f(B_t)\Biggr]\Biggr]\notag\\
&\quad =\lv^2\E^{(\widetilde{B}^1,\widetilde{B}^2)}_{(\widetilde{x}_0^1,\widetilde{x}_0^2)}\Biggl[\exp\Biggl\{\lv \int_0^{t'-t} \d r \varphi_\vep(\widetilde{B}^{\bi_0}_r)\Biggr\}\E^{(B^1,B^2,B^3,B^4)}_{(x_0^1,x_0^2,\widetilde{B}^1_{t'-t},\widetilde{B}^2_{t'-t})}[f(B_t)]\Biggr]\notag\\
&\quad \quad+\sum_{m=1}^\infty \sum_{\stackrel{\scriptstyle \bi_1,\cdots,\bi_m\in \mathcal E_4}{\bi_1\neq \cdots \neq \bi_m}}\sum_{\spin \in \{0,1\}^m} \lv^2 \E^{(\widetilde{B}^1,\widetilde{B}^2)}_{(\widetilde{x}_0^1,\widetilde{x}_0^2)}\Biggl[\exp\Biggl\{\lv \int_0^{t'-t} \d r \varphi_\vep(\widetilde{B}^{\bi_0}_r)\Biggr\}\int_{\Delta_m(t)}\d \bs u_m\otimes \d \bs v_m\notag\\
&\quad\quad \times \Biggl(\, \prod_{\scriptstyle \ell:\sigma_\ell=0 \atop \scriptstyle 1\leq \ell\leq m}\delta_{0}(v_\ell)\Biggr)\int \d \mathcal G^\spin_{\vep;0}  \int \d \mathcal G^{\bs \sigma}_{\vep;1}
 \cdots\int \d \mathcal G^\spin_{\vep;m} f\mid_{x_0^3=\widetilde{B}^1_{t'-t},x_0^4=\widetilde{B}^2_{t'-t}}\Biggr].\label{mom:series0}
\end{align} 
With the function $\widetilde{f}_t$ defined below \eqref{eq:FKmom}, the first term on the right-hand side of \eqref{mom:series0} can be written as the first term on the right-hand side of \eqref{mom:series}.
Also, applying \eqref{graphical22} with $t$ replaced by $t'-t$ to the second term on the right-hand side of \eqref{mom:series0} yields the last two terms on the right-hand side of \eqref{mom:series}. The proof is complete.
\end{proof}

\begin{figure}
[t]
\begin{center}
 \begin{tikzpicture}[scale=1]
     \draw[gray, thin] (-4.5,0) -- (-1.5,0);
         \foreach \i in {-4.5,-3.5,-2.5,-1.5} {\draw [gray] (\i,-.05) -- (\i,.05);}
    \draw[gray, thin] (-.5,0) -- (7.6,0);
    \foreach \i in {-.5, 0.5, 1.2, 2.0, 2.5, 3.2, 3.7, 4.7, 5.5, 6.4, 7.0, 7.6} {\draw [gray] (\i,-.05) -- (\i,.05);}
    \draw (-4.5,-0.07) node[below]{$0$};
    \draw (-3.5,-0.03) node[below]{$\widetilde{s}_{1}$};
    \draw (-2.5,-0.03) node[below]{$\widetilde{\tau}_{1}$};
    \draw (-1.5,-0.03) node[below]{$t'-t$};
    \draw (-.3,-0.03) node[below]{$0$};
    \draw (0.5,-0.07) node[below]{$s_1$};
    \draw (1.2,-0.07) node[below]{$\tau_1$};
    \draw (2.0,-0.07) node[below]{$s_2$};
    \draw (2.5,-0.07) node[below]{$\tau_2$};
    \draw (3.2,-0.07) node[below]{$s_3$};
    \draw (3.7,-0.07) node[below]{$\tau_3$};
    \draw (4.7,-0.07) node[below]{$s_4$};
    \draw (5.5,-0.07) node[below]{$\tau_4$};
    \draw (6.4,-0.07) node[below]{$s_5$};
    \draw (7.0,-0.07) node[below]{$\tau_5$};
    \draw (7.6,-0.03) node[below]{$t$};
    \draw [line width=2.5pt, color=red!80!white] (-.5,2.30) -- (0.5,2.20);  
    \draw [line width=2.5pt, color=cyan] (-4.5,1.00) edge[bend right=15] (-3.5,1.90);  
    \draw [line width=2.5pt, color=cyan] (-4.5,1.00) edge[bend left=15] (-3.5,2.05);  
    \draw [line width=2.5pt, color=olive!70!white] (-3.5,2.00) -- (-2.5,1.20);  
    \draw [line width=2.5pt, color=olive!70!white] (-2.5,1.20) -- (-1.5,1.70);  
    \draw [line width=2.5pt, color=olive!70!white] (-2.5,1.3) -- (-1.5,2.30);  
    \draw [line width=2.5pt, color=red!80!white] (1.2,2.00) -- (2.0,2.20);      
    \draw [line width=2.5pt, color=red!80!white] (0.5,2.20) -- (1.2,2.00);  
    \draw [line width=2.5pt, color=red!80!white] (2.0,2.20) -- (2.5,1.80);      
    \draw [line width=2.5pt, color=red!80!white] (2.5,1.80)--(3.2,2.10);   
    \draw [line width=2.5pt, color=red!80!white] (2.5,1.25)--(3.2,1.35);    
    \draw [line width=2.5pt, color=red!80!white] (3.2,1.35)--(3.7,1.05);    
    \draw [line width=2.5pt, color=red!80!white] (3.2,2.10)--(3.7,2.02);         
    \draw [line width=2.5pt, color=red!80!white] (3.7,1.05)--(4.7,2.25);    
    \draw [line width=2.5pt, color=red!80!white] (3.7,2.02)--(4.7,2.35);   
    \draw [line width=2.5pt, color=green!80!black] (-.5,0.70) edge[bend right=15] (0.5,0.54);  
    \draw [line width=2.5pt, color=green!80!black] (-0.5,0.70) edge[bend left=15] (0.5,0.66); 
    \draw [line width=2.5pt, color=blue!70!white] (-.5,1.70) -- (0.5,1.15);      
    \draw [line width=2.5pt, color=blue!70!white] (0.5,0.60) -- (1.2,0.80);  
    \draw [line width=2.5pt, color=blue!70!white] (0.5,1.15) -- (1.2,1.40);    
    \draw [line width=2.5pt, color=blue!70!white] (1.2,0.85) -- (2.0,1.4);  
    \draw [line width=2.5pt, color=blue!70!white] (1.2,1.40) -- (2.0,1.5);  
    \draw [line width=2.5pt, color=purple!60!white] (1.2,0.75) -- (2.0,0.15);  
    \draw [line width=2.5pt, color=purple!60!white] (2.0,0.15) -- (2.5,0.45);  
    \draw [line width=2.5pt, color=purple!60!white] (2.0,1.46)--(2.5,1.18);  
    \draw [line width=2.5pt, color=purple!60!white] (2.5,0.45)--(3.2,0.8);  
    \draw [line width=2.5pt, color=purple!60!white] (2.5,1.15)--(3.2,0.9);  
    \draw [line width=2.5pt, color=orange] (3.7,0.45)--(4.7,1.25);  
    \draw [line width=2.5pt, color=orange] (3.7,0.55)--(4.7,1.75);  
    \draw [line width=2.5pt, color=orange] (4.7,1.25)--(5.5,0.45);  
    \draw [line width=2.5pt, color=orange] (4.7,1.75)--(5.5,1.20);  
    \draw [line width=2.5pt, color=orange] (5.5,0.45)--(6.4,0.3);  
    \draw [line width=2.5pt, color=orange] (5.5,1.20)--(6.4,0.4);  
    \draw [line width=2.5pt, color=red!80!white] (3.22,0.85)--(3.7,0.49);  
    \draw [line width=2.5pt, color=orange] (4.7,2.32)--(5.5,1.79);   
    \draw [line width=2.5pt, color=yellow!90!black] (5.5,1.75)--(6.4,1.15);  
    \draw [line width=2.5pt, color=yellow!90!black] (5.5,1.85)--(6.4,2.40);  
    \draw [line width=2.5pt, color=yellow!90!black] (6.4,0.34)--(7.0,0.99);  
    \draw [line width=2.5pt, color=yellow!90!black] (6.4,1.15)--(7.0,1.85);  
    \draw [line width=2.5pt, color=yellow!90!black] (6.4,2.40)--(7.0,2.30);  
    \draw [line width=2.5pt, color=yellow!90!black] (7.0,0.9)--(7.6,0.45);  
    \draw [line width=2.5pt, color=yellow!90!black] (7.0,1.0)--(7.6,0.85);  
    \draw [line width=2.5pt, color=yellow!90!black] (7.0,1.85)--(7.6,1.55);  
    \draw [line width=2.5pt, color=yellow!90!black] (7.0,2.30)--(7.6,2.15);  
    \node at (-4.5,1.00) {$\bullet$};
    \node at (-3.5,1.90) {$\bullet$};
    \node at (-3.5,2.00) {$\bullet$};
    \node at (-2.5,1.30) {$\bullet$};
    \node at (-2.5,1.20) {$\bullet$};
    \node at (-1.5,2.30) {$\bullet$};
    \node at (-1.5,1.70) {$\bullet$};
    \draw [thick, color=black, snake=coil, segment length=4pt] (-3.5,1.95)--(-2.5,1.2) ;
    \draw [thick, color=black] (-4.5,1.00) edge[bend right=15] (-3.5,1.90);    
    \draw [thick, color=black] (-4.5,1.00) edge[bend left=15] (-3.5,2.05);    
    \draw [thick, color=black] (-2.5,1.3)--(-1.5,2.3);    
    \draw [thick, color=black] (-2.5,1.2)--(-1.5,1.7);   
    \draw [densely dashed] (-1.5,2.3) -- (-.5,2.3); 
    \draw [densely dashed] (-1.5,1.7) -- (-.5,1.7);          
    \node at (-.5,2.30) {$\bullet$};
    \node at (-.5,2.30) {$\bullet$};
    \node at (-.5,1.70) {$\bullet$};
    \node at (-.5,0.70) {$\bullet$};
    \draw [thick, color=black] (-.5,2.30)--(0.5,2.20);
    \draw [thick, color=black] (-.5,0.70) edge[bend right=15] (0.5,0.54);    
    \draw [thick, color=black] (-.5,0.70) edge[bend left=15] (0.5,0.66);    
    \draw [thick, color=black] (-.5,1.7)--(0.5,1.15);    
    \draw (-4.5,.6) node [left] {$\wt{x}_0^1$};
    \draw (-4.64,.96) node [left] {$\shortparallel $};
    \draw (-4.5,1.3) node [left] {$\wt{x}^2_0$};
    \draw (-.5,0.4) node [left] {$x_0^1$};
    \draw (-.62,0.72) node [left] {$\shortparallel $};
    \draw (-.5,1.) node [left] {$x_0^2$};
    \draw (-.2,1.97) node [left] {$x_0^3$};
    \draw (-.2,2.57) node [left] {$x_0^4$};
    \draw (-1.15,1.97) node [left] {$\wt{x}_2^1$};
    \draw (-1.15,2.57) node [left] {$\wt{x}_2^2$};
    \node at (0.5,0.55) {$\bullet$};
    \node at (0.5,0.65) {$\bullet$};
    \node at (0.5,1.15) {$\bullet$};
    \node at (0.5,2.20) {$\bullet$};
    \draw [thick, color=black, snake=coil, segment length=4pt] (0.5,0.60) -- (1.2,0.80);
    \draw [thick, color=black] (0.5,1.15) -- (1.2,1.40);    
    \draw [thick, color=black] (0.5,2.20) -- (1.2,2.00);    
    \node at (1.2,0.75) {$\bullet$};
    \node at (1.2,0.85) {$\bullet$};
    \node at (1.2,1.40) {$\bullet$};
    \node at (1.2,2.00) {$\bullet$};
    \draw [thick, color=black] (1.2,0.75) -- (2.0,0.15);
    \draw [thick, color=black] (1.2,0.85) -- (2.0,1.4);
    \draw [thick, color=black] (1.2,1.40) -- (2.0,1.5);    
    \draw [thick, color=black] (1.2,2.00) -- (2.0,2.20);  
    \node at (2.0,0.15) {$\bullet$};
    \node at (2.0,1.4) {$\bullet$};
    \node at (2.0,1.5) {$\bullet$};
    \node at (2.0,2.20) {$\bullet$};
    \draw [thick, color=black] (2.0,0.15) -- (2.5,0.45);
    \draw [thick, color=black, snake=coil, segment length=4pt] (2.0,1.44) -- (2.52,1.185);    
    \draw [thick, color=black] (2.0,2.20) -- (2.5,1.80);    
    \node at (2.5,0.45) {$\bullet$};
    \node at (2.5,1.25) {$\bullet$};
    \node at (2.5,1.15) {$\bullet$};
    \node at (2.5,1.80) {$\bullet$};
    \draw [thick, color=black] (2.5,0.45)--(3.2,0.8);  
    \draw [thick, color=black] (2.5,1.15)--(3.2,0.9);    
    \draw [thick, color=black] (2.5,1.25)--(3.2,1.35);  
    \draw [thick, color=black] (2.5,1.80)--(3.2,2.10);        
    \node at (3.2,0.80) {$\bullet$};
    \node at (3.2,0.90) {$\bullet$};
    \node at (3.2,1.35) {$\bullet$};
    \node at (3.2,2.10) {$\bullet$};
    \draw [thick, color=black, snake=coil, segment length=4pt] (3.2,0.85)--(3.734,0.45) ;
    \draw [thick, color=black] (3.2,1.35)--(3.7,1.05);    
    \draw [thick, color=black] (3.2,2.10)--(3.7,2.02);  
    \node at (3.7,0.45) {$\bullet$};
    \node at (3.7,0.55) {$\bullet$};
    \node at (3.7,1.05) {$\bullet$};
    \node at (3.7,2.02) {$\bullet$};
    \draw [thick, color=black] (3.7,0.45)--(4.7,1.25);
    \draw [thick, color=black] (3.7,0.55)--(4.7,1.75);    
    \draw [thick, color=black] (3.7,1.05)--(4.7,2.25);    
    \draw [thick, color=black] (3.7,2.02)--(4.7,2.35);    
    \node at (4.7,1.25) {$\bullet$};
    \node at (4.7,1.75) {$\bullet$};
    \node at (4.7,2.25) {$\bullet$};
    \node at (4.7,2.35) {$\bullet$};
    \draw [thick, color=black] (4.7,1.25)--(5.5,0.45);
    \draw [thick, color=black] (4.7,1.75)--(5.5,1.20);   
    \draw [thick, color=black, snake=coil, segment length=4pt] (4.7,2.30)--(5.5,1.80); 
    \node at (5.5,0.45) {$\bullet$};
    \node at (5.5,1.20) {$\bullet$};
    \node at (5.5,1.75) {$\bullet$};
    \node at (5.5,1.85) {$\bullet$};
    \draw [thick, color=black] (5.5,0.45)--(6.4,0.3);
    \draw [thick, color=black] (5.5,1.20)--(6.4,0.4); 
    \draw [thick, color=black] (5.5,1.75)--(6.4,1.15);    
    \draw [thick, color=black] (5.5,1.85)--(6.4,2.40); 
    \node at (6.4,0.3) {$\bullet$};
    \node at (6.4,0.4) {$\bullet$};
    \node at (6.4,1.15) {$\bullet$};
    \node at (6.4,2.40) {$\bullet$};
    \draw [thick, color=black, snake=coil, segment length=4pt] (6.4,0.35)--(7.0,0.95);
    \draw [thick, color=black] (6.4,1.15)--(7.0,1.85);    
    \draw [thick, color=black] (6.4,2.40)--(7.0,2.30);
    \node at (7.0,0.9) {$\bullet$};
    \node at (7.0,1.0) {$\bullet$};
    \node at (7.0,1.85) {$\bullet$};
    \node at (7.0,2.30) {$\bullet$};
    \draw [thick, color=black] (7.0,0.9)--(7.6,0.45);
    \draw [thick, color=black] (7.0,1.0)--(7.6,0.85);    
    \draw [thick, color=black] (7.0,1.85)--(7.6,1.55); 
    \draw [thick, color=black] (7.0,2.30)--(7.6,2.15); 
    \node at (7.6,0.45) {$\bullet$};
    \node at (7.6,0.85) {$\bullet$};
    \node at (7.6,1.55) {$\bullet$};
    \node at (7.6,2.15) {$\bullet$};
\end{tikzpicture}
\end{center}
\vspace{-.5cm}
\caption{This figure illustrates an example of $(\widetilde{\mathcal G}^{\widetilde{\spin}}_\vep, \mathcal G^{\spin}_\vep)$ for $\bi_1 = (2, 1)$. On the left, the graph $\widetilde{\mathcal G}^{\widetilde{\spin}}_\vep$ spans the interval $[0, t'-t]$ and depicts two particles undergoing attractive interactions, which reach the states $\wt{x}_2^1$ and $\wt{x}_2^2$ at time $t'-t$. On the right, the graph $\mathcal G^{\spin}_\vep$ covers $[0, t]$ and represents four particles with initial states $x^1_0$, $x^2_0$, $\wt{x}_2^1 = x^3_0$, and $\wt{x}_2^2 = x^4_0$, where $x^2_0$, $\wt{x}_2^1$, and $\wt{x}_2^2$ are distinct, and $x^1_0=x_0^2$. The identification $\wt{x}_2^1 = x^3_0$ and $\wt{x}_2^2 = x^4_0$ ``glues'' the left graph to the right graph at the associated vertices. Note that this figure is in the ``extremal case'' of the initial states such that $\wt{x}_0^1 = \wt{x}_0^2$ and $x_0^1 = x_0^2$.}\label{fig:2}
\begin{center}
 \begin{tikzpicture}[scale=1]
     \draw[gray, thin] (-4.5,0) -- (-1.5,0);
         \foreach \i in {-4.5,-3.5,-2.5,-1.5} {\draw [gray] (\i,-.05) -- (\i,.05);}
    \draw[gray, thin] (-.5,0) -- (7.6,0);
    \foreach \i in {-.5, 0.5, 1.2, 2.0, 2.5, 3.2, 3.7, 4.7, 5.5, 6.4, 7.0, 7.6} {\draw [gray] (\i,-.05) -- (\i,.05);}
    \draw (-4.5,-0.07) node[below]{$0$};
    \draw (-3.5,-0.03) node[below]{$\widetilde{s}_{1}$};
    \draw (-2.5,-0.03) node[below]{$\widetilde{\tau}_{1}$};
    \draw (-1.5,-0.03) node[below]{$t'-t$};
    \draw (-.3,-0.03) node[below]{$0$};
    \draw (0.5,-0.07) node[below]{$s_1$};
    \draw (1.2,-0.07) node[below]{$\tau_1$};
    \draw (2.0,-0.07) node[below]{$s_2$};
    \draw (2.5,-0.07) node[below]{$\tau_2$};
    \draw (3.2,-0.07) node[below]{$s_3$};
    \draw (3.7,-0.07) node[below]{$\tau_3$};
    \draw (4.7,-0.07) node[below]{$s_4$};
    \draw (5.5,-0.07) node[below]{$\tau_4$};
    \draw (6.4,-0.07) node[below]{$s_5$};
    \draw (7.0,-0.07) node[below]{$\tau_5$};
    \draw (7.6,-0.03) node[below]{$t$};
    \draw [line width=2.5pt, color=green!80!black] (-.5,2.30) -- (0.5,2.20);  
    \draw [line width=2.5pt, color=cyan] (-4.5,1.00) edge[bend right=15] (-3.5,1.90);  
    \draw [line width=2.5pt, color=cyan] (-4.5,1.00) edge[bend left=15] (-3.5,2.05);  
    \draw [line width=2.5pt, color=olive!70!white] (-3.5,2.00) -- (-2.5,1.20);  
    \draw [line width=2.5pt, color=olive!70!white] (-2.5,1.20) -- (-1.5,1.70);  
    \draw [line width=2.5pt, color=olive!70!white] (-2.5,1.3) -- (-1.5,2.30);  
    \draw [line width=2.5pt, color=red!80!white] (1.2,2.00) -- (2.0,2.20);  
    \draw [line width=2.5pt, color=blue!70!white] (0.5,2.20) -- (1.2,2.00);  
    \draw [line width=2.5pt, color=red!80!white] (2.0,2.20) -- (2.5,1.80);       
    \draw [line width=2.5pt, color=red!80!white] (2.5,1.80)--(3.2,2.10);     
    \draw [line width=2.5pt, color=red!80!white] (2.5,1.25)--(3.2,1.35);  
    \draw [line width=2.5pt, color=red!80!white] (3.2,1.35)--(3.7,1.05);     
    \draw [line width=2.5pt, color=red!80!white] (3.2,2.10)--(3.7,2.02);     
    \draw [line width=2.5pt, color=red!80!white] (3.7,1.05)--(4.7,2.25);       
    \draw [line width=2.5pt, color=red!80!white] (3.7,2.02)--(4.7,2.35);  
    \draw [line width=2.5pt, color=purple!60!white] (-.5,0.70) -- (0.5,0.6);  
    \draw [line width=2.5pt, color=blue!70!white] (-0.5,0.70) -- (0.5,1.15);  
    \draw [line width=2.5pt, color=green!80!black] (-.5,1.70) -- (0.5,2.20);   
    \draw [line width=2.5pt, color=purple!60!white] (0.5,0.60) -- (1.2,0.80);  
    \draw [line width=2.5pt, color=blue!70!white] (0.5,1.15) -- (1.2,1.40);  
    \draw [line width=2.5pt, color=blue!70!white] (1.2,1.4) -- (2.0,1.4);  
    \draw [line width=2.5pt, color=blue!70!white] (1.2,2.00) -- (2.0,1.5);  
    \draw [line width=2.5pt, color=purple!60!white] (1.2,0.8) -- (2.0,0.15);  
    \draw [line width=2.5pt, color=purple!60!white] (2.0,0.15) -- (2.5,0.45);  
    \draw [line width=2.5pt, color=purple!60!white] (2.0,1.46)--(2.5,1.18);  
    \draw [line width=2.5pt, color=purple!60!white] (2.5,0.45)--(3.2,0.8);  
    \draw [line width=2.5pt, color=purple!60!white] (2.5,1.15)--(3.2,0.9);  
    \draw [line width=2.5pt, color=orange] (3.7,0.45)--(4.7,1.25); 
    \draw [line width=2.5pt, color=orange] (3.7,0.55)--(4.7,1.75);  
    \draw [line width=2.5pt, color=orange] (4.7,1.25)--(5.5,0.45);  
    \draw [line width=2.5pt, color=orange] (4.7,1.75)--(5.5,1.20); 
    \draw [line width=2.5pt, color=orange] (5.5,0.45)--(6.4,0.3);  
    \draw [line width=2.5pt, color=orange] (5.5,1.20)--(6.4,0.4);  
    \draw [line width=2.5pt, color=red!80!white] (3.22,0.85)--(3.7,0.49);  
    \draw [line width=2.5pt, color=orange] (4.7,2.32)--(5.5,1.79);   
    \draw [line width=2.5pt, color=yellow!90!black] (5.5,1.75)--(6.4,1.15);  
    \draw [line width=2.5pt, color=yellow!90!black] (5.5,1.85)--(6.4,2.40);  
    \draw [line width=2.5pt, color=yellow!90!black] (6.4,0.34)--(7.0,0.99);  
    \draw [line width=2.5pt, color=yellow!90!black] (6.4,1.15)--(7.0,1.85);  
    \draw [line width=2.5pt, color=yellow!90!black] (6.4,2.40)--(7.0,2.30);  
    \draw [line width=2.5pt, color=yellow!90!black] (7.0,0.9)--(7.6,0.45);  
    \draw [line width=2.5pt, color=yellow!90!black] (7.0,1.0)--(7.6,0.85); 
    \draw [line width=2.5pt, color=yellow!90!black] (7.0,1.85)--(7.6,1.55);  
    \draw [line width=2.5pt, color=yellow!90!black] (7.0,2.30)--(7.6,2.15);  
    \node at (-4.5,1.00) {$\bullet$};
    \node at (-3.5,1.90) {$\bullet$};
    \node at (-3.5,2.00) {$\bullet$};
    \node at (-2.5,1.30) {$\bullet$};
    \node at (-2.5,1.20) {$\bullet$};
    \node at (-1.5,2.30) {$\bullet$};
    \node at (-1.5,1.70) {$\bullet$};
    \draw [thick, color=black, snake=coil, segment length=4pt] (-3.5,1.95)--(-2.5,1.2) ;
    \draw [thick, color=black] (-4.5,1.00) edge[bend right=15] (-3.5,1.90);    
    \draw [thick, color=black] (-4.5,1.00) edge[bend left=15] (-3.5,2.05);    
    \draw [thick, color=black] (-2.5,1.3)--(-1.5,2.3);    
    \draw [thick, color=black] (-2.5,1.2)--(-1.5,1.7);   
    \draw [densely dashed] (-1.5,2.3) -- (-.5,2.3); 
    \draw [densely dashed] (-1.5,1.7) -- (-.5,1.7);           
    \node at (-.5,2.30) {$\bullet$};
    \node at (-.5,2.30) {$\bullet$};
    \node at (-.5,1.70) {$\bullet$};
    \node at (-.5,0.70) {$\bullet$};
    \draw [thick, color=black] (-.5,2.30)--(0.5,2.20);
    \draw [thick, color=black] (-.5,0.70)--(0.5,0.6);    
    \draw [thick, color=black] (-.5,0.70)--(0.5,1.15);    
    \draw [thick, color=black] (-.5,1.7)--(0.5,2.20);    
    \draw (-4.5,.6) node [left] {$\wt{x}_0^1$};
    \draw (-4.64,.96) node [left] {$\shortparallel $};
    \draw (-4.5,1.3) node [left] {$\wt{x}^2_0$};
    \draw (-.5,0.4) node [left] {$x_0^1$};
    \draw (-.62,0.72) node [left] {$\shortparallel $};
    \draw (-.5,1.) node [left] {$x_0^2$};
    \draw (-.2,1.97) node [left] {$x_0^3$};
    \draw (-.2,2.57) node [left] {$x_0^4$};
    \draw (-1.15,1.97) node [left] {$\wt{x}_2^1$};
    \draw (-1.15,2.57) node [left] {$\wt{x}_2^2$};
    \node at (0.5,0.60) {$\bullet$};
    \node at (0.5,1.15) {$\bullet$};
    \node at (0.5,2.15) {$\bullet$};
    \node at (0.5,2.25) {$\bullet$};
    \draw [thick, color=black] (0.5,0.60) -- (1.2,0.80);
    \draw [thick, color=black] (0.5,1.15) -- (1.2,1.40);    
    \draw [thick, color=black, , snake=coil, segment length=4pt] (0.5,2.20) -- (1.2,2.00);    
    \node at (1.2,0.80) {$\bullet$};
    \node at (1.2,1.40) {$\bullet$};
    \node at (1.2,2.05) {$\bullet$};
    \node at (1.2,1.95) {$\bullet$};
    \draw [thick, color=black] (1.2,0.8) -- (2.0,0.15);
    \draw [thick, color=black] (1.2,1.4) -- (2.0,1.4);
    \draw [thick, color=black] (1.2,2.00) -- (2.0,1.5);    
    \draw [thick, color=black] (1.2,2.00) -- (2.0,2.20);  
    \node at (2.0,0.15) {$\bullet$};
    \node at (2.0,1.4) {$\bullet$};
    \node at (2.0,1.5) {$\bullet$};
    \node at (2.0,2.20) {$\bullet$};
    \draw [thick, color=black] (2.0,0.15) -- (2.5,0.45);
    \draw [thick, color=black, snake=coil, segment length=4pt] (2.0,1.44) -- (2.52,1.185);    
    \draw [thick, color=black] (2.0,2.20) -- (2.5,1.80);    
    \node at (2.5,0.45) {$\bullet$};
    \node at (2.5,1.25) {$\bullet$};
    \node at (2.5,1.15) {$\bullet$};
    \node at (2.5,1.80) {$\bullet$};
    \draw [thick, color=black] (2.5,0.45)--(3.2,0.8);  
    \draw [thick, color=black] (2.5,1.15)--(3.2,0.9);    
    \draw [thick, color=black] (2.5,1.25)--(3.2,1.35);  
    \draw [thick, color=black] (2.5,1.80)--(3.2,2.10);        
    \node at (3.2,0.80) {$\bullet$};
    \node at (3.2,0.90) {$\bullet$};
    \node at (3.2,1.35) {$\bullet$};
    \node at (3.2,2.10) {$\bullet$};
    \draw [thick, color=black, snake=coil, segment length=4pt] (3.2,0.85)--(3.734,0.45) ;
    \draw [thick, color=black] (3.2,1.35)--(3.7,1.05);    
    \draw [thick, color=black] (3.2,2.10)--(3.7,2.02);  
    \node at (3.7,0.45) {$\bullet$};
    \node at (3.7,0.55) {$\bullet$};
    \node at (3.7,1.05) {$\bullet$};
    \node at (3.7,2.02) {$\bullet$};
    \draw [thick, color=black] (3.7,0.45)--(4.7,1.25);
    \draw [thick, color=black] (3.7,0.55)--(4.7,1.75);    
    \draw [thick, color=black] (3.7,1.05)--(4.7,2.25);    
    \draw [thick, color=black] (3.7,2.02)--(4.7,2.35);    
    \node at (4.7,1.25) {$\bullet$};
    \node at (4.7,1.75) {$\bullet$};
    \node at (4.7,2.25) {$\bullet$};
    \node at (4.7,2.35) {$\bullet$};
    \draw [thick, color=black] (4.7,1.25)--(5.5,0.45);
    \draw [thick, color=black] (4.7,1.75)--(5.5,1.20);   
    \draw [thick, color=black, snake=coil, segment length=4pt] (4.7,2.30)--(5.5,1.80); 
    \node at (5.5,0.45) {$\bullet$};
    \node at (5.5,1.20) {$\bullet$};
    \node at (5.5,1.75) {$\bullet$};
    \node at (5.5,1.85) {$\bullet$};
    \draw [thick, color=black] (5.5,0.45)--(6.4,0.3);
    \draw [thick, color=black] (5.5,1.20)--(6.4,0.4); 
    \draw [thick, color=black] (5.5,1.75)--(6.4,1.15);    
    \draw [thick, color=black] (5.5,1.85)--(6.4,2.40); 
    \node at (6.4,0.3) {$\bullet$};
    \node at (6.4,0.4) {$\bullet$};
    \node at (6.4,1.15) {$\bullet$};
    \node at (6.4,2.40) {$\bullet$};
    \draw [thick, color=black, snake=coil, segment length=4pt] (6.4,0.35)--(7.0,0.95);
    \draw [thick, color=black] (6.4,1.15)--(7.0,1.85);    
    \draw [thick, color=black] (6.4,2.40)--(7.0,2.30);
    \node at (7.0,0.9) {$\bullet$};
    \node at (7.0,1.0) {$\bullet$};
    \node at (7.0,1.85) {$\bullet$};
    \node at (7.0,2.30) {$\bullet$};
    \draw [thick, color=black] (7.0,0.9)--(7.6,0.45);
    \draw [thick, color=black] (7.0,1.0)--(7.6,0.85);    
    \draw [thick, color=black] (7.0,1.85)--(7.6,1.55); 
    \draw [thick, color=black] (7.0,2.30)--(7.6,2.15); 
    \node at (7.6,0.45) {$\bullet$};
    \node at (7.6,0.85) {$\bullet$};
    \node at (7.6,1.55) {$\bullet$};
    \node at (7.6,2.15) {$\bullet$};
\end{tikzpicture}
\end{center}
\vspace{-.5cm}
\caption{$\bi_1=(4,3)$}
\label{fig:3}
\end{figure}

\begin{figure}[!htb]
\begin{center}
 \begin{tikzpicture}[scale=1]
     \draw[gray, thin] (-4.5,0) -- (-1.5,0);
         \foreach \i in {-4.5,-3.5,-2.5,-1.5} {\draw [gray] (\i,-.05) -- (\i,.05);}
    \draw[gray, thin] (-.5,0) -- (7.6,0);
    \foreach \i in {-.5, 0.5, 1.2, 2.0, 2.5, 3.2, 3.7, 4.7, 5.5, 6.4, 7.0, 7.6} {\draw [gray] (\i,-.05) -- (\i,.05);}
    \draw (-4.5,-0.07) node[below]{$0$};
    \draw (-3.5,-0.03) node[below]{$\widetilde{s}_{1}$};
    \draw (-2.5,-0.03) node[below]{$\widetilde{\tau}_{1}$};
    \draw (-1.5,-0.03) node[below]{$t'-t$};
    \draw (-.3,-0.03) node[below]{$0$};
    \draw (0.5,-0.07) node[below]{$s_1$};
    \draw (1.2,-0.07) node[below]{$\tau_1$};
    \draw (2.0,-0.07) node[below]{$s_2$};
    \draw (2.5,-0.07) node[below]{$\tau_2$};
    \draw (3.2,-0.07) node[below]{$s_3$};
    \draw (3.7,-0.07) node[below]{$\tau_3$};
    \draw (4.7,-0.07) node[below]{$s_4$};
    \draw (5.5,-0.07) node[below]{$\tau_4$};
    \draw (6.4,-0.07) node[below]{$s_5$};
    \draw (7.0,-0.07) node[below]{$\tau_5$};
    \draw (7.6,-0.03) node[below]{$t$};
    \draw [line width=2.5pt, color=blue!70!white] (-.5,2.30) -- (0.5,2.20);  
    \draw [line width=2.5pt, color=cyan] (-4.5,1.00) edge[bend right=15] (-3.5,1.90);  
    \draw [line width=2.5pt, color=cyan] (-4.5,1.00) edge[bend left=15] (-3.5,2.05);  
    \draw [line width=2.5pt, color=olive!70!white] (-3.5,2.00) -- (-2.5,1.20); 
    \draw [line width=2.5pt, color=olive!70!white] (-2.5,1.20) -- (-1.5,1.70);  
    \draw [line width=2.5pt, color=olive!70!white] (-2.5,1.3) -- (-1.5,2.30);  
    \draw [line width=2.5pt, color=blue!70!white] (1.2,2.00) -- (2.0,2.20);      
    \draw [line width=2.5pt, color=blue!70!white] (0.5,2.20) -- (1.2,2.00);  
    \draw [line width=2.5pt, color=purple!60!white] (2.0,2.20) -- (2.5,1.80);      
    \draw [line width=2.5pt, color=red!80!white] (2.5,1.80)--(3.2,2.10);      
    \draw [line width=2.5pt, color=red!80!white] (2.5,1.80)--(3.2,1.35);  
    \draw [line width=2.5pt, color=red!80!white] (3.2,1.35)--(3.7,1.05);       
    \draw [line width=2.5pt, color=red!80!white] (3.2,2.10)--(3.7,2.02);         
    \draw [line width=2.5pt, color=red!80!white] (3.7,1.05)--(4.7,2.25);       
    \draw [line width=2.5pt, color=red!80!white] (3.7,2.02)--(4.7,2.35);    
    \draw [line width=2.5pt, color=purple!60!white] (-.5,0.70) -- (0.5,0.6);  
    \draw [line width=2.5pt, color=green!80!black] (-0.5,0.70) -- (0.5,1.15);  
    \draw [line width=2.5pt, color=green!80!black] (-.5,1.70) -- (0.5,1.15);     
    \draw [line width=2.5pt, color=purple!60!white] (0.5,0.60) -- (1.2,0.80);  
    \draw [line width=2.5pt, color=blue!70!white] (0.5,1.15) -- (1.2,1.40);    
    \draw [line width=2.5pt, color=purple!60!white] (1.2,1.4) -- (2.0,1.45);  
    \draw [line width=2.5pt, color=blue!70!white] (1.2,1.40) -- (2.0,2.20);  
    \draw [line width=2.5pt, color=purple!60!white] (1.2,0.80) -- (2.0,0.15);  
    \draw [line width=2.5pt, color=purple!60!white] (2.0,0.15) -- (2.5,0.45);  
    \draw [line width=2.5pt, color=purple!60!white] (2.0,1.45)--(2.5,1.15);  
    \draw [line width=2.5pt, color=purple!60!white] (2.5,0.45)--(3.2,0.8);  
    \draw [line width=2.5pt, color=purple!60!white] (2.5,1.15)--(3.2,0.9);  
    \draw [line width=2.5pt, color=orange] (3.7,0.45)--(4.7,1.25);  
    \draw [line width=2.5pt, color=orange] (3.7,0.55)--(4.7,1.75);  
    \draw [line width=2.5pt, color=orange] (4.7,1.25)--(5.5,0.45); 
    \draw [line width=2.5pt, color=orange] (4.7,1.75)--(5.5,1.20);  
    \draw [line width=2.5pt, color=orange] (5.5,0.45)--(6.4,0.3);  
    \draw [line width=2.5pt, color=orange] (5.5,1.20)--(6.4,0.4);  
    \draw [line width=2.5pt, color=red!80!white] (3.22,0.85)--(3.7,0.49);  
    \draw [line width=2.5pt, color=orange] (4.7,2.32)--(5.5,1.79);   
    \draw [line width=2.5pt, color=yellow!90!black] (5.5,1.75)--(6.4,1.15);  
    \draw [line width=2.5pt, color=yellow!90!black] (5.5,1.85)--(6.4,2.40);  
    \draw [line width=2.5pt, color=yellow!90!black] (6.4,0.34)--(7.0,0.99); 
    \draw [line width=2.5pt, color=yellow!90!black] (6.4,1.15)--(7.0,1.85);  
    \draw [line width=2.5pt, color=yellow!90!black] (6.4,2.40)--(7.0,2.30);  
    \draw [line width=2.5pt, color=yellow!90!black] (7.0,0.9)--(7.6,0.45);  
    \draw [line width=2.5pt, color=yellow!90!black] (7.0,1.0)--(7.6,0.85);  
    \draw [line width=2.5pt, color=yellow!90!black] (7.0,1.85)--(7.6,1.55);  
    \draw [line width=2.5pt, color=yellow!90!black] (7.0,2.30)--(7.6,2.15);  
    \node at (-4.5,1.00) {$\bullet$};
    \node at (-3.5,1.90) {$\bullet$};
    \node at (-3.5,2.00) {$\bullet$};
    \node at (-2.5,1.30) {$\bullet$};
    \node at (-2.5,1.20) {$\bullet$};
    \node at (-1.5,2.30) {$\bullet$};
    \node at (-1.5,1.70) {$\bullet$};
    \draw [thick, color=black, snake=coil, segment length=4pt] (-3.5,1.95)--(-2.5,1.2) ;
    \draw [thick, color=black] (-4.5,1.00) edge[bend right=15] (-3.5,1.90);    
    \draw [thick, color=black] (-4.5,1.00) edge[bend left=15] (-3.5,2.05);    
    \draw [thick, color=black] (-2.5,1.3)--(-1.5,2.3);    
    \draw [thick, color=black] (-2.5,1.2)--(-1.5,1.7);   
    \draw [densely dashed] (-1.5,2.3) -- (-.5,2.3); 
    \draw [densely dashed] (-1.5,1.7) -- (-.5,1.7);        
    \node at (-.5,2.30) {$\bullet$};
    \node at (-.5,2.30) {$\bullet$};
    \node at (-.5,1.70) {$\bullet$};
    \node at (-.5,0.70) {$\bullet$};
    \draw [thick, color=black] (-.5,2.30)--(0.5,2.20);
    \draw [thick, color=black] (-.5,0.70)--(0.5,1.15);    
    \draw [thick, color=black] (-.5,0.70)--(0.5,0.60);    
    \draw [thick, color=black] (-.5,1.7)--(0.5,1.15);    
    \draw (-4.5,.6) node [left] {$\wt{x}_0^1$};
    \draw (-4.64,.96) node [left] {$\shortparallel $};
    \draw (-4.5,1.3) node [left] {$\wt{x}^2_0$};
    \draw (-.5,0.4) node [left] {$x_0^1$};
    \draw (-.62,0.72) node [left] {$\shortparallel $};
    \draw (-.5,1.) node [left] {$x_0^2$};
    \draw (-.2,1.97) node [left] {$x_0^3$};
    \draw (-.2,2.57) node [left] {$x_0^4$};
    \draw (-1.15,1.97) node [left] {$\wt{x}_2^1$};
    \draw (-1.15,2.57) node [left] {$\wt{x}_2^2$};
    \node at (0.5,0.60) {$\bullet$};
    \node at (0.5,1.10) {$\bullet$};
    \node at (0.5,1.20) {$\bullet$};
    \node at (0.5,2.20) {$\bullet$};
    \draw [thick, color=black] (0.5,0.60) -- (1.2,0.80);
    \draw [thick, color=black,  snake=coil, segment length=4pt] (0.5,1.15) -- (1.2,1.40);    
    \draw [thick, color=black] (0.5,2.20) -- (1.2,2.00);    
    \node at (1.2,0.80) {$\bullet$};
    \node at (1.2,1.35) {$\bullet$};
    \node at (1.2,1.45) {$\bullet$};
    \node at (1.2,2.00) {$\bullet$};
    \draw [thick, color=black] (1.2,0.80) -- (2.0,0.15);
    \draw [thick, color=black] (1.2,1.4) -- (2.0,1.45);
    \draw [thick, color=black] (1.2,1.40) -- (2.0,2.20);    
    \draw [thick, color=black] (1.2,2.00) -- (2.0,2.20);  
    \node at (2.0,0.15) {$\bullet$};
    \node at (2.0,1.45) {$\bullet$};
    \node at (2.0,2.25) {$\bullet$};
    \node at (2.0,2.15) {$\bullet$};
    \draw [thick, color=black] (2.0,0.15) -- (2.5,0.45);
    \draw [thick, color=black] (2.0,1.44) -- (2.5,1.15);    
    \draw [thick, color=black, snake=coil, segment length=4pt] (2.0,2.20) -- (2.5,1.80);    
    \node at (2.5,0.45) {$\bullet$};
    \node at (2.5,1.15) {$\bullet$};
    \node at (2.5,1.85) {$\bullet$};
    \node at (2.5,1.75) {$\bullet$};
    \draw [thick, color=black] (2.5,0.45)--(3.2,0.8);  
    \draw [thick, color=black] (2.5,1.15)--(3.2,0.9);    
    \draw [thick, color=black] (2.5,1.80)--(3.2,1.35);  
    \draw [thick, color=black] (2.5,1.80)--(3.2,2.10);        
    \node at (3.2,0.80) {$\bullet$};
    \node at (3.2,0.90) {$\bullet$};
    \node at (3.2,1.35) {$\bullet$};
    \node at (3.2,2.10) {$\bullet$};
    \draw [thick, color=black, snake=coil, segment length=4pt] (3.2,0.85)--(3.734,0.45) ;
    \draw [thick, color=black] (3.2,1.35)--(3.7,1.05);    
    \draw [thick, color=black] (3.2,2.10)--(3.7,2.02);  
    \node at (3.7,0.45) {$\bullet$};
    \node at (3.7,0.55) {$\bullet$};
    \node at (3.7,1.05) {$\bullet$};
    \node at (3.7,2.02) {$\bullet$};
    \draw [thick, color=black] (3.7,0.45)--(4.7,1.25);
    \draw [thick, color=black] (3.7,0.55)--(4.7,1.75);    
    \draw [thick, color=black] (3.7,1.05)--(4.7,2.25);    
    \draw [thick, color=black] (3.7,2.02)--(4.7,2.35);    
    \node at (4.7,1.25) {$\bullet$};
    \node at (4.7,1.75) {$\bullet$};
    \node at (4.7,2.25) {$\bullet$};
    \node at (4.7,2.35) {$\bullet$};
    \draw [thick, color=black] (4.7,1.25)--(5.5,0.45);
    \draw [thick, color=black] (4.7,1.75)--(5.5,1.20);   
    \draw [thick, color=black, snake=coil, segment length=4pt] (4.7,2.30)--(5.5,1.80); 
    \node at (5.5,0.45) {$\bullet$};
    \node at (5.5,1.20) {$\bullet$};
    \node at (5.5,1.75) {$\bullet$};
    \node at (5.5,1.85) {$\bullet$};
    \draw [thick, color=black] (5.5,0.45)--(6.4,0.3);
    \draw [thick, color=black] (5.5,1.20)--(6.4,0.4); 
    \draw [thick, color=black] (5.5,1.75)--(6.4,1.15);    
    \draw [thick, color=black] (5.5,1.85)--(6.4,2.40); 
    \node at (6.4,0.3) {$\bullet$};
    \node at (6.4,0.4) {$\bullet$};
    \node at (6.4,1.15) {$\bullet$};
    \node at (6.4,2.40) {$\bullet$};
    \draw [thick, color=black, snake=coil, segment length=4pt] (6.4,0.35)--(7.0,0.95);
    \draw [thick, color=black] (6.4,1.15)--(7.0,1.85);    
    \draw [thick, color=black] (6.4,2.40)--(7.0,2.30);
    \node at (7.0,0.9) {$\bullet$};
    \node at (7.0,1.0) {$\bullet$};
    \node at (7.0,1.85) {$\bullet$};
    \node at (7.0,2.30) {$\bullet$};
    \draw [thick, color=black] (7.0,0.9)--(7.6,0.45);
    \draw [thick, color=black] (7.0,1.0)--(7.6,0.85);    
    \draw [thick, color=black] (7.0,1.85)--(7.6,1.55); 
    \draw [thick, color=black] (7.0,2.30)--(7.6,2.15); 
    \node at (7.6,0.45) {$\bullet$};
    \node at (7.6,0.85) {$\bullet$};
    \node at (7.6,1.55) {$\bullet$};
    \node at (7.6,2.15) {$\bullet$};
\end{tikzpicture}
\end{center}
\begin{center}
 \begin{tikzpicture}[scale=1]
     \draw[gray, thin] (-4.5,0) -- (-1.5,0);
         \foreach \i in {-4.5,-3.5,-2.5,-1.5} {\draw [gray] (\i,-.05) -- (\i,.05);}
    \draw[gray, thin] (-.5,0) -- (7.6,0);
    \foreach \i in {-.5, 0.5, 1.2, 2.0, 2.5, 3.2, 3.7, 4.7, 5.5, 6.4, 7.0, 7.6} {\draw [gray] (\i,-.05) -- (\i,.05);}
    \draw (-4.5,-0.07) node[below]{$0$};
    \draw (-3.5,-0.03) node[below]{$\widetilde{s}_{1}$};
    \draw (-2.5,-0.03) node[below]{$\widetilde{\tau}_{1}$};
    \draw (-1.5,-0.03) node[below]{$t'-t$};
    \draw (-.3,-0.03) node[below]{$0$};
    \draw (0.5,-0.07) node[below]{$s_1$};
    \draw (1.2,-0.07) node[below]{$\tau_1$};
    \draw (2.0,-0.07) node[below]{$s_2$};
    \draw (2.5,-0.07) node[below]{$\tau_2$};
    \draw (3.2,-0.07) node[below]{$s_3$};
    \draw (3.7,-0.07) node[below]{$\tau_3$};
    \draw (4.7,-0.07) node[below]{$s_4$};
    \draw (5.5,-0.07) node[below]{$\tau_4$};
    \draw (6.4,-0.07) node[below]{$s_5$};
    \draw (7.0,-0.07) node[below]{$\tau_5$};
    \draw (7.6,-0.03) node[below]{$t$};
    \draw [line width=2.5pt, color=red!80!white] (-.5,2.30) -- (0.5,2.20);  
    \draw [line width=2.5pt, color=cyan] (-4.5,1.00) edge[bend right=15] (-3.5,1.90);  
    \draw [line width=2.5pt, color=cyan] (-4.5,1.00) edge[bend left=15] (-3.5,2.05);  
    \draw [line width=2.5pt, color=olive!70!white] (-3.5,2.00) -- (-2.5,1.20);  
    \draw [line width=2.5pt, color=olive!70!white] (-2.5,1.20) -- (-1.5,1.70);  
    \draw [line width=2.5pt, color=olive!70!white] (-2.5,1.3) -- (-1.5,2.30);  
    \draw [line width=2.5pt, color=red!80!white] (1.2,2.00) -- (2.0,2.20);   
    \draw [line width=2.5pt, color=red!80!white] (0.5,2.20) -- (1.2,2.00);  
    \draw [line width=2.5pt, color=red!80!white] (2.0,2.20) -- (2.5,1.80);   
    \draw [line width=2.5pt, color=red!80!white] (2.5,1.80)--(3.2,2.10);     
    \draw [line width=2.5pt, color=red!80!white] (3.2,1.35)--(3.7,1.05);      
    \draw [line width=2.5pt, color=red!80!white] (3.2,2.10)--(3.7,2.02);      
    \draw [line width=2.5pt, color=red!80!white] (3.7,1.05)--(4.7,2.25);    
    \draw [line width=2.5pt, color=red!80!white] (3.7,2.02)--(4.7,2.35);     
    \draw [line width=2.5pt, color=blue!70!white] (-.5,0.70) -- (0.5,0.6);  
    \draw [line width=2.5pt, color=green!80!black] (-0.5,0.70) -- (0.5,1.15);  
    \draw [line width=2.5pt, color=green!80!black] (-.5,1.70) -- (0.5,1.15);      
    \draw [line width=2.5pt, color=blue!70!white] (0.5,0.60) -- (1.2,0.80);  
    \draw [line width=2.5pt, color=blue!70!white] (0.5,1.15) -- (1.2,1.40);  
    \draw [line width=2.5pt, color=purple!60!white ] (1.2,1.4) -- (2.0,1.45);  
    \draw [line width=2.5pt, color=blue!70!white] (1.2,1.40) -- (2.0,0.15);  
    \draw [line width=2.5pt, color=blue!70!white] (1.2,0.80) -- (2.0,0.15);  
    \draw [line width=2.5pt, color=purple!60!white] (2.0,0.15) -- (2.5,0.45); 
    \draw [line width=2.5pt, color=purple!60!white] (2.0,1.45)--(2.5,1.15);  
    \draw [line width=2.5pt, color=orange] (2.5,0.45)--(3.2,0.85); 
    \draw [line width=2.5pt, color=purple!60!white] (2.5,1.15)--(3.2,1.35); 
    \draw [line width=2.5pt, color=purple!60!white] (2.5,0.45)--(3.2,1.35);  
    \draw [line width=2.5pt, color=orange] (3.7,0.5)--(4.7,1.25); 
    \draw [line width=2.5pt, color=orange] (3.7,1.05)--(4.7,1.75);  
    \draw [line width=2.5pt, color=orange] (4.7,1.25)--(5.5,0.45);  
    \draw [line width=2.5pt, color=orange] (4.7,1.75)--(5.5,1.20);  
    \draw [line width=2.5pt, color=orange] (5.5,0.45)--(6.4,0.3);  
    \draw [line width=2.5pt, color=orange] (5.5,1.20)--(6.4,0.4);  
    \draw [line width=2.5pt, color=orange] (3.2,0.85)--(3.7,0.5);  
    \draw [line width=2.5pt, color=orange] (4.7,2.32)--(5.5,1.79);   
    \draw [line width=2.5pt, color=yellow!90!black] (5.5,1.75)--(6.4,1.15);  
    \draw [line width=2.5pt, color=yellow!90!black] (5.5,1.85)--(6.4,2.40);  
    \draw [line width=2.5pt, color=yellow!90!black] (6.4,0.34)--(7.0,0.99); 
    \draw [line width=2.5pt, color=yellow!90!black] (6.4,1.15)--(7.0,1.85);  
    \draw [line width=2.5pt, color=yellow!90!black] (6.4,2.40)--(7.0,2.30);  
    \draw [line width=2.5pt, color=yellow!90!black] (7.0,0.9)--(7.6,0.45);  
    \draw [line width=2.5pt, color=yellow!90!black] (7.0,1.0)--(7.6,0.85);  
    \draw [line width=2.5pt, color=yellow!90!black] (7.0,1.85)--(7.6,1.55);  
    \draw [line width=2.5pt, color=yellow!90!black] (7.0,2.30)--(7.6,2.15);  
    \node at (-4.5,1.00) {$\bullet$};
    \node at (-3.5,1.90) {$\bullet$};
    \node at (-3.5,2.00) {$\bullet$};
    \node at (-2.5,1.30) {$\bullet$};
    \node at (-2.5,1.20) {$\bullet$};
    \node at (-1.5,2.30) {$\bullet$};
    \node at (-1.5,1.70) {$\bullet$};
    \draw [thick, color=black, snake=coil, segment length=4pt] (-3.5,1.95)--(-2.5,1.2) ;
    \draw [thick, color=black] (-4.5,1.00) edge[bend right=15] (-3.5,1.90);    
    \draw [thick, color=black] (-4.5,1.00) edge[bend left=15] (-3.5,2.05);    
    \draw [thick, color=black] (-2.5,1.3)--(-1.5,2.3);    
    \draw [thick, color=black] (-2.5,1.2)--(-1.5,1.7);   
    \draw [densely dashed] (-1.5,2.3) -- (-.5,2.3); 
    \draw [densely dashed] (-1.5,1.7) -- (-.5,1.7);        
    \node at (-.5,2.30) {$\bullet$};
    \node at (-.5,2.30) {$\bullet$};
    \node at (-.5,1.70) {$\bullet$};
    \node at (-.5,0.70) {$\bullet$};
    \draw [thick, color=black] (-.5,2.30)--(0.5,2.20);
    \draw [thick, color=black] (-.5,0.70)--(0.5,1.15);    
    \draw [thick, color=black] (-.5,0.70)--(0.5,0.60);    
    \draw [thick, color=black] (-.5,1.7)--(0.5,1.15);    
    \draw (-4.5,.6) node [left] {$\wt{x}_0^1$};
    \draw (-4.64,.96) node [left] {$\shortparallel $};
    \draw (-4.5,1.3) node [left] {$\wt{x}^2_0$};
    \draw (-.5,0.4) node [left] {$x_0^1$};
    \draw (-.62,0.72) node [left] {$\shortparallel $};
    \draw (-.5,1.) node [left] {$x_0^2$};
    \draw (-.2,1.97) node [left] {$x_0^3$};
    \draw (-.2,2.57) node [left] {$x_0^4$};
    \draw (-1.15,1.97) node [left] {$\wt{x}_2^1$};
    \draw (-1.15,2.57) node [left] {$\wt{x}_2^2$};
    \node at (0.5,0.60) {$\bullet$};
    \node at (0.5,1.10) {$\bullet$};
    \node at (0.5,1.20) {$\bullet$};
    \node at (0.5,2.20) {$\bullet$};
    \draw [thick, color=black] (0.5,0.60) -- (1.2,0.80);
    \draw [thick, color=black,  snake=coil, segment length=4pt] (0.5,1.15) -- (1.2,1.40);    
    \draw [thick, color=black] (0.5,2.20) -- (1.2,2.00);    
    \node at (1.2,0.80) {$\bullet$};
    \node at (1.2,1.35) {$\bullet$};
    \node at (1.2,1.45) {$\bullet$};
    \node at (1.2,2.00) {$\bullet$};
    \draw [thick, color=black] (1.2,0.80) -- (2.0,0.15);
    \draw [thick, color=black] (1.2,1.4) -- (2.0,1.45);
    \draw [thick, color=black] (1.2,1.40) -- (2.0,0.15);    
    \draw [thick, color=black] (1.2,2.00) -- (2.0,2.20);  
    \node at (2.0,0.20) {$\bullet$};
    \node at (2.0,0.10) {$\bullet$};
    \node at (2.0,1.45) {$\bullet$};
    \node at (2.0,2.20) {$\bullet$};
    \draw [thick, color=black,snake=coil, segment length=4pt] (2.0,0.15) -- (2.5,0.45);
    \draw [thick, color=black] (2.0,1.45) -- (2.5,1.15);    
    \draw [thick, color=black] (2.0,2.20) -- (2.5,1.80);    
    \node at (2.5,0.50) {$\bullet$};
    \node at (2.5,0.40) {$\bullet$};
    \node at (2.5,1.15) {$\bullet$};
    \node at (2.5,1.80) {$\bullet$};
    \draw [thick, color=black] (2.5,0.45)--(3.2,1.35);  
    \draw [thick, color=black] (2.5,0.45)--(3.2,0.85);  
    \draw [thick, color=black] (2.5,1.15)--(3.2,1.35);    
    \draw [thick, color=black] (2.5,1.80)--(3.2,2.10);        
    \node at (3.2,0.85) {$\bullet$};
    \node at (3.2,1.40) {$\bullet$};
    \node at (3.2,1.30) {$\bullet$};
    \node at (3.2,2.10) {$\bullet$};
    \draw [thick, color=black] (3.2,0.85)--(3.7,0.5) ;
    \draw [thick, color=black,snake=coil, segment length=4pt] (3.2,1.35)--(3.7,1.05);    
    \draw [thick, color=black] (3.2,2.10)--(3.7,2.02);  
    \node at (3.7,0.5) {$\bullet$};
    \node at (3.7,1.10) {$\bullet$};
    \node at (3.7,1.00) {$\bullet$};
    \node at (3.7,2.02) {$\bullet$};
    \draw [thick, color=black] (3.7,0.5)--(4.7,1.25);
    \draw [thick, color=black] (3.7,1.05)--(4.7,1.75);    
    \draw [thick, color=black] (3.7,1.05)--(4.7,2.25);    
    \draw [thick, color=black] (3.7,2.02)--(4.7,2.35);    
    \node at (4.7,1.25) {$\bullet$};
    \node at (4.7,1.75) {$\bullet$};
    \node at (4.7,2.25) {$\bullet$};
    \node at (4.7,2.35) {$\bullet$};
    \draw [thick, color=black] (4.7,1.25)--(5.5,0.45);
    \draw [thick, color=black] (4.7,1.75)--(5.5,1.20);   
    \draw [thick, color=black, snake=coil, segment length=4pt] (4.7,2.30)--(5.5,1.80); 
    \node at (5.5,0.45) {$\bullet$};
    \node at (5.5,1.20) {$\bullet$};
    \node at (5.5,1.75) {$\bullet$};
    \node at (5.5,1.85) {$\bullet$};
    \draw [thick, color=black] (5.5,0.45)--(6.4,0.3);
    \draw [thick, color=black] (5.5,1.20)--(6.4,0.4); 
    \draw [thick, color=black] (5.5,1.75)--(6.4,1.15);    
    \draw [thick, color=black] (5.5,1.85)--(6.4,2.40); 
    \node at (6.4,0.3) {$\bullet$};
    \node at (6.4,0.4) {$\bullet$};
    \node at (6.4,1.15) {$\bullet$};
    \node at (6.4,2.40) {$\bullet$};
    \draw [thick, color=black, snake=coil, segment length=4pt] (6.4,0.35)--(7.0,0.95);
    \draw [thick, color=black] (6.4,1.15)--(7.0,1.85);    
    \draw [thick, color=black] (6.4,2.40)--(7.0,2.30);
    \node at (7.0,0.9) {$\bullet$};
    \node at (7.0,1.0) {$\bullet$};
    \node at (7.0,1.85) {$\bullet$};
    \node at (7.0,2.30) {$\bullet$};
    \draw [thick, color=black] (7.0,0.9)--(7.6,0.45);
    \draw [thick, color=black] (7.0,1.0)--(7.6,0.85);    
    \draw [thick, color=black] (7.0,1.85)--(7.6,1.55); 
    \draw [thick, color=black] (7.0,2.30)--(7.6,2.15); 
    \node at (7.6,0.45) {$\bullet$};
    \node at (7.6,0.85) {$\bullet$};
    \node at (7.6,1.55) {$\bullet$};
    \node at (7.6,2.15) {$\bullet$};
\end{tikzpicture}
\end{center}
\vspace{-.5cm}
\caption{$\bi_1=(3,2)$}
\label{fig:4}
\end{figure}

The following definition will facilitate subsequent applications of \eqref{2-4} and Lemma~\ref{lem:FKmom1}. 

\begin{defi}\label{def:glue}
Let $\widetilde{\mathcal G}^{\wt{\spin}}_\vep\oplus \mathcal G^\spin_\vep$ denote the graph obtained by glueing $\widetilde{\mathcal G}_\vep^{\wt{\spin}}$ over $[0,t'-t]$ and $\mathcal G_\vep^{\spin}$ over $[0,t]$ such that the following identification of vertices holds:
$(\wt{x}_2^1,t'-t)$ is identified with $(x_0^3,0)$ and $(\wt{x}_2^2,t'-t)$ is identified with $(x_0^4,0)$. \hfill $\blacklozenge$ 
\end{defi}

In the remainder of Section~\ref{sec:mombdd}, our objective is to establish bounds for the right-hand side of \eqref{mom:series}, and thereby, end the proof of Theorem~\ref{thm:mombdd}. To obtain bounds for the right-hand side of \eqref{mom:series}, we will show that the first term there can be uniformly bounded with respect to $(\widetilde{x}^1_0, \widetilde{x}^2_0) \in \R^4$ and $t \in [0, T]$ for any $0 < T < \infty$. The primary focus, however, is on the second and third terms, which require a more intricate analysis that simultaneously encompasses both two- and four-particle systems. 

To address the second and third terms, each represented by a triple series, the key observation is that the bounds from \cite[Proposition~4.14]{C:DBG} for other moments of $X_\vep$ can be extended when we interpret the summands in those triple series in terms of the combined graph $\widetilde{\mathcal G}^{\widetilde{\spin}}_\vep \oplus \mathcal G^{\spin}_\vep$, rather than considering $\widetilde{\mathcal G}^{\widetilde{\spin}}_\vep$ and $\mathcal G^{\spin}_\vep$ separately. This extension relies on the fact that by using edges of $\widetilde{\mathcal G}^{\widetilde{\spin}}_\vep$, some of the paths in $\mathcal G^{\spin}_\vep$ can be extended in any case of $\bi_1$. Morover, to identify these extended paths, it suffices to distinguish $\bi_1$ in three scenarios: (i) $\bi_1=(2,1)$; (ii) $\bi_1=(4,3)$; and (iii) $\bi_1=(i,j)$ for $i \in \{4,3\}$ and $j \in \{2,1\}$. [This classification is due to symmetry; recall \eqref{eq:FKmom} and \eqref{mom:series}. More details appear in the proof of Theorem~\ref{thm:mombdd}.] Thus, we focus on the following three scenarios: 
\begin{itemize}
\item [(i)] $\bi_1 = (2, 1)$: the duration of the spine in $\mathcal G^{\spin}_{1;\vep}$ can be increased by $t' - t - \widetilde{\tau}_1$ (see Figure~\ref{fig:2}).
\item [(ii)] $\bi_1 = (4, 3)$: the durations of both edges in $\mathcal G^{\spin}_{\vep;0}$ can be extended by $t' - t - \widetilde{\tau}_1$ (see Figure~\ref{fig:3}).
\item [(iii)] $\bi_1 = (3, 2)$: 
the duration of an edge in $\mathcal G^{\spin}_{\vep;0}$ can be extended by $t' - t - \widetilde{\tau}_1$ (see Figure~\ref{fig:4}).  
\end{itemize}

We now introduce a collection of tools and establish two principal sets of bounds, formulated as Propositions~\ref{prop:mombdd1} and \ref{prop:mombdd2}.  First, we use the following notation extending \eqref{def:svep}:
\begin{align}
\begin{split}
\label{def:os}
\overline{\s}^{\lambda,\phi}_\vep(\tau)&\,\defeq \int_{\R^4}\d x\otimes \d \tilde{y}\phi(x)\s^{\lambda,\phi}_\vep(0,\tau; x/\two,\tilde{y})\varphi(\tilde{y})\\
&=\int_{\R^4}\d x^{i_\ell\prime}_\ell \otimes\d z_\ell \varphi(x_\ell^{\bi_\ell})\s^{\lambda,\phi}_\vep(0,\tau;  x^{\bi_\ell}_\ell,z_\ell)\varphi(z_\ell)\\
&=\int_{\R^2}\d x\phi(x)
\E_{\vep  x}^{\two W}\left[ \lv^2\exp\left\{\lv\int_0^\tau \d r \varphi_\vep(W_r) \right\}\varphi_\vep(W_\tau)\right],
\end{split}\\
\mathfrak S^{\lambda,\phi}_\vep(q)&\,\defeq \int_0^\infty \d \tau\e^{-q\tau}\overline{\s}^{\lambda,\phi}_\vep(\tau).\label{def:Svep}
\end{align}
The following lemma collects some key properties of $\overline{\s}^{\lambda,\phi}_\vep(\tau)$ and $\mathfrak S^{\lambda,\phi}_\vep(q)$. 

\begin{lem}\label{lem:Sbounds}
Recall that $\lv$ is defined by \eqref{def:lv}, using a given constant $\lambda\in \R$. \medskip

\noindent {\rm (1$\cc$)} 
 There exists $q_{\ref{def:Sconv}}=q_{\ref{def:Sconv}}(\lambda,\lVert\varphi\rVert_\infty)\in (0,\infty)$ such that 
\begin{align}\label{def:Sconv}
\lim_{\vep\to 0}\mathfrak S^{\lambda,\phi}_\vep(q)=\frac{4\pi}{\log (q/\beta)},\quad \forall\;q\in (q_{\ref{def:Sconv}},\infty),
\end{align}
where $\beta>0$ is defined as in \eqref{def:betaintro}. \medskip 

\noindent {\rm (2$\cc$)} 
For all $t\in (0,\infty)$ and $\vep\in (0,\ovl]$,
\begin{linenomath*}\begin{align}
 \sup_{y\in \R^2}\int_0^t \d \tau\int_{\R^2}\d \tilde{y}\s^{\lambda,\phi}_\vep(0,\tau; y,\tilde{y})\varphi(\tilde{y})
&= \sup_{y\in \R^2}\int_0^{t}\d \tau\lv^2\E^W_{\vep y}\left[\exp\left\{\lv\int_0^\tau \d r \varphi_\vep(W_r) \right\}\varphi_\vep(W_\tau)\right]\notag\\
&\leq C(\lambda,\lVert\varphi\rVert_\infty)\lvert\lv\rvert [1+\log^+ (\vep^{-2}t)]\e^{q_{\ref{ineq:apriori}}(\lambda,\lVert\varphi\rVert_\infty)t}.\label{ineq:apriori}
\end{align}\end{linenomath*}
\end{lem}
\begin{proof}
The proof of \eqref{def:Sconv} and  the proof of \eqref{ineq:apriori} whenever $\vep\in (0,\vep(\lambda,\lVert\varphi\rVert_\infty))$ are already obtained in 
\cite[Proposition~5.2, p.178]{C:DBG} and \cite[Lemma~5.4, p.180]{C:DBG}, respectively, provided that $\vep(\lambda,\lVert\varphi\rVert_\infty)>0$ is a small enough constant. 
Moreover, by enlarging $q_{\ref{ineq:apriori}}(\lambda,\lVert\varphi\rVert_\infty)$ if necessary, it is plain that \eqref{ineq:apriori} also holds for all $\vep\in [\vep(\lambda,\lVert\varphi\rVert_\infty),\ovl]$.   
\end{proof}

We now proceed to the proofs of two principal sets of bounds. The following proposition establishes the first set of results, concerning the iterated graphical integrals in \eqref{mom:series} for $f\equiv \1$. For the statement, recall the change of variables in \eqref{def:uv123}, and we use the following notations: $\prod_{\ell=j}^{k}=\prod_{\ell\in \varnothing}\equiv 1$ for all $k<j$, and
\begin{align*}
\widehat{\spin}=(\wh{\sigma}_1,\cdots,\wh{\sigma}_m)\,\defeq\,(1-\sigma_1,\cdots,1-\sigma_m),\quad  \lVert\spin\rVert\,\defeq\,\sum_{\ell=1}^m \sigma_\ell,\quad \forall\;\spin\in \{0,1\}^m.
\end{align*}

\begin{prop}\label{prop:mombdd1}
Let $\vep\in (0,\ovl]$, $0< t<t'<\infty$, $x_0^1,x_0^2,\widetilde{x}_0^1,\widetilde{x}_0^2\in \R^2$, $m\in \Bbb N$, $\bi_1\neq \cdots \neq \bi_m$ with $\bi_1,\cdots,\bi_m\in \mc E_4$, $\widetilde{\spin}\in \{0,1\}$, and $\spin\in \{0,1\}^m$. Consider the corresponding iterated graphical integrals in \eqref{mom:series}.
\begin{itemize}
\item [\rm (1$\cc$)] For $\bi_1=(2,1)$,
\begin{align}
&\int \d \widetilde{\mathcal G}^{\widetilde{\bs \sigma}}_{\vep;0}\int \d \widetilde{\mathcal G}^{\widetilde{\bs \sigma}}_{\vep;1} \int \d \mathcal G^{\bs \sigma}_{\vep;0}\int \d \mathcal G^{\bs \sigma}_{\vep;1}\cdots \int \d \mathcal G^{\bs \sigma}_{\vep;m}\1\notag\\
&\quad \leq \int_{\R^2}\d \wt{\overline{z}}
 P_{2\wt{s}_1}(\wt{x}_0^{2}-\vep\wt{\overline{z}},\wt{x}_0^{1})\phi(\wt{\overline{z}}) \Biggl(\int_{\R^2} \d \wt{z}_1\s^{\lambda,\phi}_\vep(\wt{s}_1,\wt{\tau}_1;\wt{\overline{z}}/\two,\wt{z}_1) \varphi(\wt{z}_1)\wt{\sigma}_1+  \widehat{\wt{\sigma}}_1\Biggr)\notag\\
&\quad\quad \times \int_{\R^2}
\d \overline{z}
  P_{2u_1}(x_0^{2}-\vep \overline{z},x_0^{1})  \phi(\overline{z})
  \Biggl(\int_{\R^2}\d z_1\s^{\lambda,\phi}_\vep(0,v_1;\overline{z}/\two,z_1)\varphi(z_1) \sigma_1+ \widehat{\sigma}_1\Biggr) \notag\\
&\quad\quad  \times\Biggl(\prod_{\scriptstyle \ell:\sigma_\ell=1\atop \scriptstyle 2\leq \ell\leq m}\overline{\s}^{\lambda,\phi}_\vep(v_\ell)\Biggr)\lv^{\lVert\wh{\wt{\spin}}\rVert+\lVert\wh{\spin}\rVert}\notag\\
&\quad\quad \times \Biggl(\frac{\1_{m\geq 2}}{t'-t-\widetilde{\tau}_1+2u_2+u_1}+\1_{m=1}\Biggr)
\left(\prod_{\ell=2}^{m-1}\frac{1}{2u_{\ell+1}+u_\ell}\right),\label{mom:bdd1}\\
&\E^{(\widetilde{B}^1,\widetilde{B}^2)}_{(\widetilde{x}_0^1,\widetilde{x}_0^2)}\Bigg[\int \d \mathcal G^\spin_{\vep;0} \int \d \mathcal G^{\bs \sigma}_{\vep;1}
 \cdots\int \d \mathcal G^\spin_{\vep;m} \1\mid_{(x_0^3,x_0^4)=(\widetilde{B}^1_{t'-t},\widetilde{B}^2_{t'-t}) }\Bigg]\notag\\
 &\quad \leq   \int_{\R^2}
\d \overline{z}
  P_{2u_1}(x_0^{2}-\vep \overline{z},x_0^{1})  \phi(\overline{z})
  \Biggl(\int_{\R^2}\d z_1\s^{\lambda,\phi}_\vep(0,v_1;\overline{z}/\two,z_1)\varphi(z_1) \sigma_1+ \widehat{\sigma}_1\Biggr)\notag\\
&\quad\quad \times\Biggl(\prod_{\scriptstyle \ell:\sigma_\ell=1\atop \scriptstyle 2\leq \ell\leq m}\overline{\s}^{\lambda,\phi}_\vep(v_\ell)\Biggr)\lv^{\lVert\wh{\spin}\rVert}\Biggl(\frac{\1_{m\geq 2}}{t'-t+2u_2+u_1}+\1_{m=1}\Biggr)
\left(\prod_{\ell=2}^{m-1}\frac{1}{2u_{\ell+1}+u_\ell}\right) .\label{mom:bdd1-1}
 \end{align}
 \item [\rm (2$\cc$)] For $\bi_1=(4,3)$ or $(3,2)$,   
\begin{align}
&\int \d \widetilde{\mathcal G}^{\widetilde{\bs \sigma}}_{\vep;0}\int \d \widetilde{\mathcal G}^{\widetilde{\bs \sigma}}_{\vep;1} \int \d \mathcal G^{\bs \sigma}_{\vep;0}\int \d \mathcal G^{\bs \sigma}_{\vep;1}\cdots \int \d \mathcal G^{\bs \sigma}_{\vep;m}\1\notag\\
&\quad \leq \int_{\R^2}\d \wt{\overline{z}}
 P_{2\wt{s}_1}(\wt{x}_0^{2}-\vep \wt{\overline{z}},\wt{x}_0^{1})\phi( \wt{\overline{z}})\left(\int_{\R^2} 
 \d \wt{z}_1
\s^{\lambda,\phi}_\vep (\wt{s}_1,\wt{\tau}_1; \wt{\overline{z}}/\two,\wt{z}_1)\varphi(\wt{z}_1)\wt{\sigma}_1 + \widehat{\wt{\sigma}}_1\right)  \notag\\
  &\quad \quad \times  \Biggl(\prod_{\scriptstyle \ell:\sigma_\ell=1\atop\scriptstyle 1\leq \ell\leq m}\overline{\s}^{\lambda,\phi}_\vep(v_\ell)\Biggr)\lv^{\lVert\wh{\wt{\spin}}\rVert+\lVert\wh{\spin}\rVert}
\Biggl(\frac{1}{t'-t-\widetilde{\tau}_1+u_1}\Biggr)\left(\prod_{\ell=1}^{m-1}\frac{1}{2u_{\ell+1}+u_\ell}\right),
\label{mom:bdd2}\\
&\E^{(\widetilde{B}^1,\widetilde{B}^2)}_{(\widetilde{x}_0^1,\widetilde{x}_0^2)}\Bigg[\int \d \mathcal G^\spin_{\vep;0} \int \d \mathcal G^{\bs \sigma}_{\vep;1}
 \cdots\int \d \mathcal G^\spin_{\vep;m} \1\mid_{(x_0^3,x_0^4)=(\widetilde{B}^1_{t'-t},\widetilde{B}^2_{t'-t}) }\Bigg]\notag\\
 &\quad \leq  \Biggl(\prod_{\scriptstyle \ell:\sigma_\ell=1\atop\scriptstyle 1\leq \ell\leq m}\overline{\s}^{\lambda,\phi}_\vep(v_\ell)\Biggr)\lv^{\lVert\wh{\spin}\rVert}
\Biggl(\frac{1}{t'-t+u_1}\Biggr)\left(\prod_{\ell=1}^{m-1}\frac{1}{2u_{\ell+1}+u_\ell}\right).\label{mom:bdd2-1}
 \end{align}
\end{itemize}
\end{prop}

\begin{proof}
This proof is divided into five steps. Steps~1 proves an explicit formula for $\int \d \mathcal G^\spin_{\vep;m}\1$. Step~2 proves a bound for $\int \d \mathcal G^{\bs \sigma}_{\vep;2}\cdots \int \d \mathcal G^{\bs \sigma}_{\vep;m}\1$ when $m\geq 2$. By using the results of Steps~1 and 2,
Step~3 proves both bounds in (1$\cc$), Step~4 proves \eqref{mom:bdd2} and \eqref{mom:bdd2-1} for $\bi=(4,3)$, and Step~5 proves \eqref{mom:bdd2} and \eqref{mom:bdd2-1} for $\bi=(3,2)$. In particular, by following these steps, the reader may regard this proof as an extension of \cite[Proof of (4.53) for Proposition~4.14 (2$\cc$), pp.167--170]{C:DBG}, with a modified method for bounding the graphical integrals $\int \d \mathcal G^\spin_{\vep;0}\1$ and $\int \d \mathcal G^\spin_{\vep;1}\1$. See Steps~3--5 for these modifications. Also, we note the following simple identity that will be used repeatedly in what follows: for all $x',x\in \R^2$ and $s,t,\vep\in (0,\infty)$ and all nonnegative $F$, 
\begin{align}\label{ineq:gausstime}
\int_{\R^4}
\begin{bmatrix}
\d z'\\
\d z
\end{bmatrix}_\otimes
\begin{bmatrix}
P_{t}(x',z+\vep(z'-z))\\
P_s(x,z)
\end{bmatrix}_{\times}
F(z'-z)&=\int_{\R^2} \d \bar{z} P_{s+t}(x'-\vep \bar{z},x) F(\bar{z}),
\end{align}
which holds by the Chapman--Kolmogorov equation.\medskip 

\noindent {\bf Step 1.}
In this step, we prove the following identity:
\begin{linenomath*}\begin{align}\label{Gm:bdd}
\int \d \mathcal G^\spin_{\vep;m}\1 =\left(
 \int_{\R^2}\d z_m \s^{\lambda,\phi}_\vep (0,v_m; x_m^{\bi_{m}},z_m)\varphi(z_m) \sigma_m+
\widehat{\sigma}_m\right)  (\lv^{1/2})^{\#_m},
\end{align}\end{linenomath*}
where $v_m$ is defined in \eqref{def:uv123}, and $\#_m$ denotes the number of left-maximal entanglement-free paths ending at time $t$
such that the left end-vertices are given by $(x_{\ell_0}^{i_{\ell_0}\prime},x_{\ell_0}^{i_{\ell_0}},s_{\ell_0})$ for some $1\leq \ell_0\leq m$ satisfying $\sigma_{\ell_0}=0$. 

Identity \eqref{Gm:bdd} holds by considering the following. Recall that $\int \d \mathcal G^\spin_{\vep;m}$ is defined with the measures from \eqref{def:productmeasure} such that the variables of integration are given by the spatial components of the vertices of $ \mathcal G^\spin_{\vep;m}$. By this definition and the Chapman--Kolmogorov equation, the right-hand side of \eqref{Gm:bdd} for $\sigma_m=1$ follows by integrating out the spatial components of the right end-vertices of $\mathcal G^{\spin}_{\vep;m}$. In particular, in doing the integration, we need to take into account the possibility that some left-maximal entanglement-free paths ending at time $t$ have left end-vertices given by $(x_{\ell_0}^{i_{\ell_0}\prime},x_{\ell_0}^{i_{\ell_0}},s_{\ell_0})$ for some $1\leq \ell_0<m$ and $\sigma_{\ell_0}=0$ [e.g. the two paths marked yellow and begun at time $s_4$ in Figure~\ref{fig:1}]. For such a left-maximal entanglement-free path, one of the first two weights in \eqref{weight4} has to come into play, and it brings in a multiplicative factor of $\lv^{1/2}$. This is how $(\lv^{1/2})^{\#_m}$ in \eqref{Gm:bdd} arises from $\int \d \mathcal G^\spin_{\vep;m}\1$. Also, recall that when $\sigma_m=1$, the weight of the coiled line segment in 
$\mathcal G^\spin_{\vep;m}$ is given by the multiplication column in \eqref{weight2} for $\ell=m$, and $v_m=\tau_m-s_m$ according to \eqref{def:uv123}. Hence, when $\sigma_m=1$, we get
\begin{linenomath*}\begin{align*}
\int \d \mathcal G^\spin_{\vep;m} \1
&=\int_{\R^4}
\begin{bmatrix}
\d z_m'\\
 \d z_m
 \end{bmatrix}_\otimes 
 \varphi(z_m)\begin{bmatrix}
P_{v_m}(\two x_m^{i_m}+\vep x_m^{\bi_m},z'_m)\\
\vspace{-.3cm}\\
\s^{\lambda,\phi}_\vep (0,v_m; x_m^{\bi_{m}},z_m)
\end{bmatrix}_\times(\lv^{1/2})^{\#_m}\\
&=\int_{\R^2}\d z_m \s^{\lambda,\phi}_\vep (0,v_m; x_m^{\bi_{m}},z_m)\varphi(z_m)(\lv^{1/2})^{\#_m},
\end{align*}\end{linenomath*}
which proves \eqref{Gm:bdd} for $\sigma_m=1$. Identity \eqref{Gm:bdd} for $\sigma_m=0$ follows from a simpler argument. \medskip 

\noindent {\bf Step 2.}
In this step, we prove that for all $m\geq 2$ and $\ell_0\in \{1,2\}$,
\begin{align}
\int \d \mathcal G^{\bs \sigma}_{\vep;\ell_0}\cdots \int \d \mathcal G^{\bs \sigma}_{\vep;m}\1
&\leq 
\Biggl(\int_{\R^2}\d z_{\ell_0}\s^{\lambda,\phi}_\vep(0,v_{\ell_0};x_{\ell_0}^{\bi_{\ell_0}},z_{\ell_0})\varphi(z_{\ell_0})\sigma_{\ell_0}+ \widehat{\sigma}_{\ell_0}\Biggr)\notag\\
&\quad \times\Biggl(\prod_{\scriptstyle \ell:\sigma_\ell=1\atop\scriptstyle \ell_0+1\leq \ell\leq m}\overline{\s}^{\lambda,\phi}_\vep(v_\ell)\Biggr)(\lv^{1/2})^{\sum_{\ell=\ell_0}^m\#_\ell}\left(\prod_{\ell=\ell_0}^{m-1}\frac{1}{2u_{\ell+1}+u_\ell}\right),\label{mom:bdd-intermediate}
\end{align}
where $\overline{\s}^{\lambda,\phi}_\vep(\cdot)$ is defined in \eqref{def:os}, $u_\ell$'s are defined in \eqref{def:uv123}, and $\#_\ell$ for $1\leq \ell\leq m-1$ are analogues of $\#_m$, defined as follows: 
\begin{align}
\#_\ell\,\defeq\,\#\left\{\mathcal P\in \{\mathcal P_\ell',\mathcal P_\ell\}; 
\begin{array}{ll}
\mbox{the left end-vertex of $\mathcal P$ is $(x_{\ell_0}^{i_{\ell_0\prime}},x_{\ell_0}^{i_{\ell_0}},s_{\ell_0})$}\\
\mbox{for some $1\leq \ell_{0}\leq \ell$ and $\sigma_{\ell_0}=0$}
\end{array}
\right\}\in \{0,1,2\}.\label{def:number1}
\end{align}
[Recall that $\mathcal P_\ell'$ and $\mathcal P_\ell$ are paths defining the spine of $\mathcal G^\spin_{\oslash;\ell}$, specified in Definition~\ref{def:subgraphedge} (2$\cc$).]
The numbers $\#_\ell$ defined in \eqref{def:number1} serve a purpose similar to $\#_m$: we use this $\#_\ell$ to keep track of the number of weights in  $\mathcal G^{\bs \sigma}_{\vep;\ell}$ such that they are equal to one of the first two weights in \eqref{weight4}. The application of $\#_1$ is deferred until Steps~3--5, though. 

By \eqref{Gm:bdd}, it suffices to prove \eqref{mom:bdd-intermediate} for $m\geq 3$. In this case,
we begin by expressing $\int \d \mathcal G^\spin_{\vep;m-1} \int \d \mathcal G_{\vep;m}^\spin\1$ more explicitly;
see the orange subgraphs in Figures~\ref{fig:1} and~\ref{fig:2} for examples of $ \mathcal G^{\spin}_{\vep;m-1}$.
By \eqref{Gm:bdd}, we have
\begin{linenomath*}\begin{align}
&\int \d \mathcal G^\spin_{\vep;m-1} \int \d \mathcal G_{\vep;m}^\spin\1\notag\\
&\quad =\int \d \mathcal G^\spin_{\vep;m-1}\left( \int_{\R^2}\d z_m \s^{\lambda,\phi}_\vep (0,v_m; x_m^{\bi_{m}},z_m)\varphi(z_m)\sigma_m+\widehat{\sigma}_m\right) (\lv^{1/2})^{\#_m} \label{m-1-1100-1}\\
&\quad =\int \d \mu \int_{\R^4}
\begin{bmatrix}
\d x^{i_{m}\prime}_{m}\\
 \d x^{i_{m}}_{m}
 \end{bmatrix}_\otimes
 \varphi(x^{\bi_m}_m) (\lv^{1/2})^{\#_{m-1}}\begin{bmatrix}
P_{\lvert\mathcal P'_{m-1}\rvert}(A',x_m^{i_m}+\vep\two x_m^{\bi_m})\\
P_{\lvert\mathcal P_{m-1}\rvert}(A,x_m^{i_m})
\end{bmatrix}_\times\notag\\
&\quad\quad \times\left(\int_{\R^2}\d z_m \s^{\lambda,\phi}_\vep (0,v_m; x_m^{\bi_{m}},z_m)\varphi(z_m)\sigma_m+\widehat{\sigma}_m\right)(\lv^{1/2})^{\#_m},
\label{m-1-1100}
\end{align}\end{linenomath*}
where $\lvert\mathcal P\rvert$ is defined to be the duration of a path $\mathcal P$ [Definition~\ref{def:path} (5$\cc$)]. Let us provide further details on why \eqref{m-1-1100} holds and explain the new notation introduced therein. First, the second multiplication column in \eqref{m-1-1100} follows by using the Chapman--Kolmogorov equation and
the spine $\mathcal P'_{m-1}\oplus \mathcal P_{m-1}$ of $\mathcal G^\spin_{\vep;m-1}$ so that  
$x^{i_m}_{m}+\vep \two x_m^{\bi_m}$ and $x^{i_m}_m$ arise as the terminal states of some of the heat kernels in \eqref{weight3} and \eqref{weight4} for $\ell=m-1$. Also, \eqref{m-1-1100} uses a measure $\d \mu$ and linear combinations $A',A$ of spatial components of vertices, chosen as follows:  
\begin{itemize}
\item When $\sigma_{m-1}=1$, $\d \mu$ is defined as follows, using the measure $\d z'_{m-1} \otimes\d z_{m-1} \varphi(z_{m-1})$ from \eqref{def:productmeasure2} and
a weight from \eqref{weight2} with $\ell=m-1$:
\[
\d \mu= \begin{bmatrix}
\d z_{m-1}'\\
\d z_{m-1}
\end{bmatrix}_{\otimes}\varphi(z_{m-1}) \begin{bmatrix}
P_{s_{m-1},\tau_{m-1}}(\two x_{m-1}^{i_{m-1}}+\vep x_{m-1}^{\bi_{m-1}},z'_{m-1})\\
\vspace{-.3cm}\\
\s^{\lambda,\phi}_\vep (s_{m-1},\tau_{m-1}; x_{m-1}^{\bi_{m-1}},z_{m-1})
\end{bmatrix}_\times .
\]
\item When $\sigma_{m-1}=0$, $\int \d\mu$ only acts as a ``placeholder,'' so we only require that $\d\mu$ satisfy $\int \d \mu\1=1$. Specifically, $\mu$ can be an arbitrary measure satisfying $\int \d \mu\1=1$, but not using any spatial components of vertices for its variables of integration. 
\item The spatial components defining $A'$ and $A$ in \eqref{m-1-1100} are chosen from the left end-vertices of $\mathcal P'_{m-1}$ and $\mathcal P_{m-1}$, and $A'$ and $A$ may depend on some of the variables of integration in $\d \mu$ in the case of $\sigma_{m-1}=1$. In particular, $A'$ and $A$ are independent of $x_m^{i_m\prime}$ and $x_m^{i_m}$. 
\end{itemize}

Next, we bound $\int \d \mathcal G^\spin_{\vep;m-1} \int \d \mathcal G_{\vep;m}^\spin\1$. By \eqref{def:varphi}, \eqref{ineq:gausstime} and \eqref{m-1-1100}, we get:
\begin{linenomath*}\begin{align}
\begin{split}
\int \d \mathcal G^\spin_{\vep;m-1} \int \d \mathcal G_{\vep;m}^\spin\1&=\int \d \mu \int_{\R^2}\d \bar{z}\phi(\bar{z})(\lv^{1/2})^{\#_{m-1}}
P_{\lvert\mathcal P'_{m-1}\oplus \mathcal P_{m-1}\rvert}(A'-\vep \bar{z},A) \\
&\quad \times\left(\int_{\R^2}\d z_m \s^{\lambda,\phi}_\vep (0,v_m; \overline{z}/\two,z_m)\varphi(z_m)\sigma_m+\widehat{\sigma}_m\right)(\lv^{1/2})^{\#_m}\notag
\end{split}\\
&\leq  \left(\int \d \mu\1\right)
\frac{(\lv^{1/2})^{\#_{m-1}}}{2\pi \lvert\mathcal P'_{m-1}\oplus \mathcal P_{m-1}\rvert }\int_{\R^2}\d \bar{z} \phi(\bar{z})\notag\\
&\quad\times\left(\int_{\R^2}\d z_m \s^{\lambda,\phi}_\vep (0,v_m; \overline{z}/\two,z_m)\varphi(z_m)\sigma_m+\widehat{\sigma}_m\right)(\lv^{1/2})^{\#_m}\notag\\
&\leq \left(\int \d \mu\1\right)(\lv^{1/2})^{\#_{m-1}}\left(\frac{1}{2u_{m}+u_{m-1}}\right)\left(\overline{\s}^{\lambda,\phi}_\vep (v_m)\sigma_m+\widehat{\sigma}_m\right)(\lv^{1/2})^{\#_m}\notag\\
&=  \biggl(\int_{\R^2}\d z_{m-1} \s^{\lambda,\phi}_\vep (0,v_{m-1};x_{m-1}^{\bi_{m-1}},z_{m-1})\varphi(z_{m-1})\sigma_{m-1}+ \widehat{\sigma}_{m-1}\biggr)\notag\\
 &\quad\times\Biggl(\prod_{\scriptstyle \ell:\sigma_\ell=1\atop\scriptstyle m\leq \ell\leq m}\overline{\s}^{\lambda,\phi}_\vep(v_\ell)\Biggr)(\lv^{1/2})^{\sum_{\ell=m-1}^m\#_\ell}\Bigg(\prod_{\ell={m-1}}^{m-1} \frac{1}{2u_{\ell+1}+u_\ell}\Bigg),\label{bdd:chain}
\end{align}\end{linenomath*}
where the first inequality uses the bound $P_r(x)\leq 1/(2\pi r)$,
 the second inequality uses the bound in Lemma~\ref{lem:path}, the definition \eqref{def:os} of $\overline{\s}^{\lambda,\phi}_\vep$, and the assumption that $\phi$ is a probability density, and
the last equality uses the choice of $\d \mu$ specified below \eqref{m-1-1100}.

Finally, by recalling how \eqref{Gm:bdd} is used in \eqref{m-1-1100-1}, it should become apparent that the derivation of \eqref{bdd:chain} naturally extends to provide similar bounds for $\int \d \mathcal G_{\vep;\ell}^\spin \cdots \int \d \mathcal G^\spin_{\vep;m} \1$, proceeding inductively in the order $\ell = m-2, m-3, \ldots, \ell_0$, for all $\ell_0\in \{1,2\}$. Consequently, \eqref{mom:bdd-intermediate} also holds when $\ell = \ell_0$.\medskip 

\noindent {\bf Step 3.}
We prove the bounds in (1$\cc$) in this step, starting with the proof of \eqref{mom:bdd1}, which uses $\bi_1=(2,1)$. See Figure~\ref{fig:2} for an example of the underlying graphs. 

For $m\geq 2$, one of $\mathcal P'_1,\mathcal P_1$ defining the spine of $\mathcal G^{\spin}_{\vep;1}$ in $\mathcal G^{\spin}_\vep$ can be extended to a left-maximal entanglement-free path in $\widetilde{\mathcal G}_\vep^{\widetilde{\spin}}\oplus \mathcal G_\vep^{\spin}$ (in the sense analogous to the one in Definition~\ref{def:path}) such that the duration is extended by $t'-t-\widetilde{\tau}_1$, where $\widetilde{\mathcal G}^{\widetilde{\spin}}_{\vep}\oplus\mathcal G^{\spin}_{\vep}$ is defined in Definition~\ref{def:glue} and, as in \eqref{def:stau-convention}, $\wt{\tau}_1=\wt{s}_1$ if $\wt{\sigma}_1=0$. More precisely, the extension must be via the vertex $(\wt{x}_2^1,t'-t)=(x_0^3,0)$ or the vertex $(\wt{x}_2^2,t'-t)=(x_0^4,0)$, as illustrated in Figure~\ref{fig:2}.
 Write $\mathcal P_1^{\prime\rm e},\mathcal P_1^{\rm e}$ for $\mathcal P_1',\mathcal P_1$ after this extension. Then with notation analogous to $\lvert\mathcal P^{\prime}_1\oplus \mathcal P_1\rvert$,
\begin{align}\label{pathbdd:var}
\lvert\mathcal P^{\prime\e}_1\oplus \mathcal P^\e_1\rvert\geq t'-t-\widetilde{\tau}_1+ 2u_2+u_1,
\end{align}
where the contribution $2u_2+u_1$ follows from Lemma~\ref{lem:path} for $\ell=1$. 

With \eqref{pathbdd:var}, we can bound the left-hand side of our main objective of Step~3, namely \eqref{mom:bdd1}, by applying the derivation of  \eqref{mom:bdd-intermediate} with $\ell_0=2$ after a slight modification. Specifically, we have the following inequalities, starting by applying  \eqref{mom:bdd-intermediate} with $\ell_0=2$ to get the first inequality:  
\begin{align}
&  \int \d \widetilde{\mathcal G}^{\widetilde{\spin}}_{\vep;0}\int \d \widetilde{\mathcal G}^{\widetilde{\spin}}_{\vep;1}\int \d \mathcal G^{\bs \sigma}_{\vep;0}
\int \d \mathcal G^{\bs \sigma}_{\vep;1}
\cdots \int \d \mathcal G^{\bs \sigma}_{\vep;m}\1\notag\\
&\quad \leq
 \int_{\R^4} \begin{bmatrix}
\d \widetilde{x}_1^{2}\\
 \d\wt{x}_1^1
 \end{bmatrix}_{\otimes}\varphi(\wt{x}_1^{\bi_0})\begin{bmatrix}
 P_{\wt{s}_1}(\wt{x}_0^{2},\wt{x}_1^{1}+\vep\two \wt{x}^{\bi_0}_1)\\
  P_{\wt{s}_1}(\wt{x}_0^{1},\wt{x}_1^{1})
\end{bmatrix}_\times\notag\\
&\quad\quad \times\Bigg(\int_{\R^4} \begin{bmatrix}
\d \wt{z}_1'\\
 \d \wt{z}_1
 \end{bmatrix}_{\otimes} \varphi(\wt{z}_1)
\begin{bmatrix}
P_{\wt{v}_{1}}(\two \wt{x}_1^{1}+\vep \wt{x}_1^{\bi_0},\wt{z}'_1)\\
\vspace{-.3cm}\\
\s^{\lambda,\phi}_\vep (\wt{s}_1,\wt{\tau}_1; \wt{x}_1^{\bi_{0}},\wt{z}_1)
\end{bmatrix}_\times\wt{\sigma}_1+ \lv\widehat{\wt{\sigma}}_1\Bigg)\notag\\
&\quad\quad \times \int_{\R^4}
\begin{bmatrix}
\d x^{2}_1\\
\d x^{1}_1
\end{bmatrix}_\otimes\varphi(x^{\bi_1}_1)
 \begin{bmatrix}
  P_{u_1}(x_0^{2},x_1^{1}+\vep\two x^{\bi_1}_1)\\
  P_{u_1}(x_0^{1},x_1^{1})
  \end{bmatrix}_\times  \notag\\
&\quad\quad \times  \Biggl(\int_{\R^2}\d z_1\s^{\lambda,\phi}_\vep(0,v_1;x_1^{\bi_1},z_1)\varphi(z_1) \sigma_1+\widehat{\sigma}_1\Biggr)\notag\\
&\quad\quad  \times\Biggl(\prod_{\scriptstyle \ell:\sigma_\ell=1\atop\scriptstyle 2\leq \ell\leq m}\overline{\s}^{\lambda,\phi}_\vep(v_\ell)\Biggr)(\lv^{1/2})^{\sum_{\ell=1}^m\#_\ell} \Biggl(\frac{1}{t'-t-\widetilde{\tau}_1+2u_2+u_1}\Biggr)
\left(\prod_{\ell=2}^{m-1}\frac{1}{2u_{\ell+1}+u_\ell}\right)\label{mom:bdd1:aux-m=1}\\
&\quad \leq   \int_{\R^2}\d \wt{\overline{z}}
 P_{2\wt{s}_1}(\wt{x}_0^{2}-\vep\wt{\overline{z}},\wt{x}_0^{1})\phi(\wt{\overline{z}}) \Biggl(\int_{\R^2} \d \wt{z}_1\s^{\lambda,\phi}_\vep(\wt{s}_1,\wt{\tau}_1;\wt{\overline{z}}/\two,\wt{z}_1) \varphi(\wt{z}_1)\wt{\sigma}_1+ \widehat{\wt{\sigma}}_1\Biggr)\notag\\
&\quad\quad \times \int_{\R^2}
\d \overline{z}
  P_{2u_1}(x_0^{2}-\vep \overline{z},x_0^{1})  \phi(\overline{z})
  \Biggl(\int_{\R^2}\d z_1\s^{\lambda,\phi}_\vep(0,v_1;\overline{z}/\two,z_1)\varphi(z_1) \sigma_1+ \widehat{\sigma}_1\Biggr)\notag\\
&\quad\quad  \times\Biggl(\prod_{\scriptstyle \ell:\sigma_\ell=1\atop \scriptstyle 2\leq \ell\leq m}\overline{\s}^{\lambda,\phi}_\vep(v_\ell)\Biggr)\lv^{\lVert\wh{\wt{\spin}}\rVert+\lVert\widehat{\spin}\rVert} \Biggl(\frac{1}{t'-t-\widetilde{\tau}_1+2u_2+u_1}\Biggr)
\left(\prod_{\ell=2}^{m-1}\frac{1}{2u_{\ell+1}+u_\ell}\right).\label{mom:bdd1:aux}
 \end{align}
In more detail, the integrators $\d \wt{x}^2_1\otimes \d \wt{x}^1_1\varphi(\wt{x}_1^{\bi_0})$ and $\d x^2_1\otimes \d x^1_1\varphi(x^{\bi_1}_1)$
on the right-hand side of \eqref{mom:bdd1:aux-m=1} are the measures, from \eqref{def:productmeasure1}, defining $\int \d \wt{\mathcal G}^{\wt{\spin}}_{\vep;0}$ and $\int \d \mathcal G^{\spin}_{\vep;0}$, respectively. Also, to get \eqref{mom:bdd1:aux}, we have used \eqref{ineq:gausstime} twice after integrating out $ \wt{z}_1'$, and 
\begin{align}\label{sigmasum}
\sum_{\ell=1}^m \#_\ell=2\sum_{\ell=1}^m\widehat{\sigma}_\ell=2\lVert\widehat{\spin}\rVert
\end{align}
because the total number of weights that appear in $\mathcal G^{\spin}_\vep$ as one of the first two in \eqref{weight4} is $2\sum_{\ell=1}^m\widehat{\sigma}_\ell$. The last inequality is the required inequality in \eqref{mom:bdd1} for $m\geq 2$. The case of $m=1$ is easier. We can just use \eqref{Gm:bdd}, rather than 
\eqref{mom:bdd-intermediate} with $\ell_0=2$, to obtain a simpler analogue of \eqref{mom:bdd1:aux-m=1} as follows and simplify in the same way as before: 
\begin{align*}
  \int \d \widetilde{\mathcal G}^{\widetilde{\spin}}_{\vep;0}\int \d \widetilde{\mathcal G}^{\widetilde{\spin}}_{\vep;1}\int \d \mathcal G^{\bs \sigma}_{\vep;0}
\int \d \mathcal G^{\bs \sigma}_{\vep;1}\1
& =
 \int_{\R^4} \begin{bmatrix}
\d \widetilde{x}_1^{2}\\
 \d\wt{x}_1^1
 \end{bmatrix}_{\otimes}\varphi(\wt{x}_1^{\bi_0})\begin{bmatrix}
 P_{\wt{s}_1}(\wt{x}_0^{2},\wt{x}_1^{1}+\vep\two \wt{x}^{\bi_0}_1)\\
  P_{\wt{s}_1}(\wt{x}_0^{1},\wt{x}_1^{1})
\end{bmatrix}_\times\\
&\quad \times\Bigg(\int_{\R^4} \begin{bmatrix}
\d \wt{z}_1'\\
 \d \wt{z}_1
 \end{bmatrix}_{\otimes} \varphi(\wt{z}_1)
\begin{bmatrix}
P_{\wt{v}_{1}}(\two \wt{x}_1^{1}+\vep \wt{x}_1^{\bi_0},\wt{z}'_1)\\
\vspace{-.3cm}\\
\s^{\lambda,\phi}_\vep (\wt{s}_1,\wt{\tau}_1; \wt{x}_1^{\bi_{0}},\wt{z}_1)
\end{bmatrix}_\times\wt{\sigma}_1+ \lv\widehat{\wt{\sigma}}_1\Bigg)\\
&\quad \times \int_{\R^4}
\begin{bmatrix}
\d x^{2}_1\\
\d x^{1}_1
\end{bmatrix}_\otimes\varphi(x^{\bi_1}_1)
 \begin{bmatrix}
  P_{u_1}(x_0^{2},x_1^{1}+\vep\two x^{\bi_1}_1)\\
  P_{u_1}(x_0^{1},x_1^{1})
  \end{bmatrix}_\times  \notag\\
&\quad \times \Biggl(
 \int_{\R^2}\d z_1 \s^{\lambda,\phi}_\vep (0,v_1; x_1^{\bi_{1}},z_1)\varphi(z_1) \sigma_1+
\widehat{\sigma}_1\Biggr)  (\lv^{1/2})^{\#_1}.
\end{align*}

The proof of  \eqref{mom:bdd1-1} is similar. For $m\geq 2$, an analogue of the bound in \eqref{pathbdd:var}
applies, now improved to
$\lvert\mathcal P^{\prime\e}_1\oplus \mathcal P^\e_1\rvert\geq t'-t+ 2u_2+u_1$. 
An inspection of the derivation of \eqref{mom:bdd1:aux} then leads to \eqref{mom:bdd1-1} for $m\geq 2$. The case of $m=1$ is again simpler, as the required bound can be deduced from the following identity, which is implied by \eqref{Gm:bdd}:
\begin{align*}
\int \d \mathcal G^{\bs \sigma}_{\vep;0}
\int \d \mathcal G^{\bs \sigma}_{\vep;1}\1
& = \int_{\R^4}
\begin{bmatrix}
\d x^{2}_1\\
\d x^{1}_1
\end{bmatrix}_\otimes\varphi(x^{\bi_1}_1)
 \begin{bmatrix}
  P_{u_1}(x_0^{2},x_1^{1}+\vep\two x^{\bi_1}_1)\\
  P_{u_1}(x_0^{1},x_1^{1})
  \end{bmatrix}_\times \\
  &\quad  \times \Biggl(
 \int_{\R^2}\d z_1 \s^{\lambda,\phi}_\vep (0,v_1; x_1^{\bi_{1}},z_1)\varphi(z_1) \sigma_1+
\widehat{\sigma}_1\Biggr)  (\lv^{1/2})^{\#_1}.
\end{align*}
 
 \noindent {\bf Step 4.} We now begin the proof of (2$\cc$), establishing \eqref{mom:bdd2} and \eqref{mom:bdd2-1} for $\bi_1=(4,3)$ in this step. The proofs of these inequalities for $\bi_1=(3,2)$ will be presented in Step~5.
 
To obtain \eqref{mom:bdd2} for $\bi_1=(4,3)$ and all $m\geq 1$,  it is necessary to obtain a sufficiently sharp bound for the iterated graphical integral $ \int \d \widetilde{\mathcal G}^{\widetilde{\spin}}_{\vep;0}\int \d \widetilde{\mathcal G}^{\widetilde{\spin}}_{\vep;1}\int \d \mathcal G^{\bs \sigma}_{\vep;0}
\int \d \mathcal G^{\bs \sigma}_{\vep;1}
\cdots \int \d \mathcal G^{\bs \sigma}_{\vep;m}\1$ as a starting point. The crucial observation for this purpose is the following: the two edges defining $\mathcal G^{\spin}_{\vep;0}$ can be extended to left-maximal entanglement-free paths within $\widetilde{\mathcal G}_\vep^{\widetilde{\spin}}\oplus \mathcal G_\vep^{\spin}$. More precisely, the duration of each edge is extended from $u_1=s_1$ to $t'-t-\wt{\tau}_1+u_1$: one extension passes through the vertex $(\wt{x}_2^1,t'-t)=(x_0^3,0)$, and the other through $(\wt{x}_2^2,t'-t)=(x_0^4,0)$. Figure~\ref{fig:3} illustrates $\widetilde{\mathcal G}_\vep^{\widetilde{\spin}}\oplus \mathcal G_\vep^{\spin}$ in this scenario. Moreover, this extension of paths can be performed using the Chapman--Kolmogorov equation when writing out $ \int \d \widetilde{\mathcal G}^{\widetilde{\spin}}_{\vep;0}\int \d \widetilde{\mathcal G}^{\widetilde{\spin}}_{\vep;1}\int \d \mathcal G^{\bs \sigma}_{\vep;0}
\int \d \mathcal G^{\bs \sigma}_{\vep;1}
\cdots \int \d \mathcal G^{\bs \sigma}_{\vep;m}\1$. The remainder of this iterated graphical integral is then handled as follows: 
\begin{itemize}
\item Expand the remainder of $\int \d \widetilde{\mathcal G}^{\widetilde{\spin}}_{\vep;0} \int \d \widetilde{\mathcal G}^{\widetilde{\spin}}_{\vep;1}$ as in \eqref{2mom:eq1} and \eqref{2mom:eq2}.
\item For $m=1$, apply \eqref{Gm:bdd}. For $m\geq 2$, utilize \eqref{mom:bdd-intermediate} with $\ell_0=1$ to bound $\int \d \mathcal G^{\bs \sigma}_{\vep;1}
\cdots \int \d \mathcal G^{\bs \sigma}_{\vep;m}\1$. 
\end{itemize}
Hence, we obtain the inequality below, which also handles the cases of $\wt{\sigma}_1=1$ and $\wt{\sigma}_1=0$ simultaneously: 
  \begin{align}
& \int \d \widetilde{\mathcal G}^{\widetilde{\spin}}_{\vep;0}\int \d \widetilde{\mathcal G}^{\widetilde{\spin}}_{\vep;1}\int \d \mathcal G^{\bs \sigma}_{\vep;0}
\int \d \mathcal G^{\bs \sigma}_{\vep;1}
\cdots \int \d \mathcal G^{\bs \sigma}_{\vep;m}\1\notag\\
&\quad \leq \wt{\sigma}_1
 \int_{\R^4} \begin{bmatrix}
\d \widetilde{x}_1^{2}\\
 \d\wt{x}_1^1
 \end{bmatrix}_{\otimes}\varphi(\wt{x}_1^{\bi_0})\begin{bmatrix}
 P_{\wt{s}_1}(\wt{x}_0^{2},\wt{x}_1^{1}+\vep\two \wt{x}^{\bi_0}_1)\\
  P_{\wt{s}_1}(\wt{x}_0^{1},\wt{x}_1^{1})
\end{bmatrix}_\times\int_{\R^4} \begin{bmatrix}
\d \wt{z}_1'\\
 \d \wt{z}_1
 \end{bmatrix}_{\otimes} \varphi(\wt{z}_1)
\begin{bmatrix}
P_{\wt{s}_1,\wt{\tau}_1}(\two \wt{x}_1^{1}+\vep \wt{x}_1^{\bi_0},\wt{z}'_1)\\
\vspace{-.3cm}\\
\s^{\lambda,\phi}_\vep (\wt{s}_1,\wt{\tau}_1; \wt{x}_1^{\bi_{0}},\wt{z}_1)
\end{bmatrix}_\times\notag\\
&\quad\quad \times \int_{\R^4}
\begin{bmatrix}
\d x^{4}_1\\
\d x^{3}_1
\end{bmatrix}_\otimes  \varphi(x^{\bi_1}_1)
 \begin{bmatrix}
  P_{t'-t-\wt{\tau}_1+u_1}\left(\frac{\wt{z}_1'+\vep \wt{z}_1}{\two},x_1^{3}+\vep\two x^{\bi_1}_1\right)\\
  P_{t'-t-\wt{\tau}_1+u_1}\left(\frac{\wt{z}_1'-\vep \wt{z}_1}{\two},x_1^{3}\right)
  \end{bmatrix}_\times\notag\\
  &\quad\quad  \times 
  \Biggl(\int_{\R^2}\d z_1\s^{\lambda,\phi}_\vep(0,v_1;x_1^{\bi_1},z_1)\varphi(z_1) \sigma_1+ \widehat{\sigma}_1\Biggr)\notag\\
  &\quad\quad  \times\Biggl(\prod_{\scriptstyle \ell:\sigma_\ell=1\atop\scriptstyle 2\leq \ell\leq m}\overline{\s}^{\lambda,\phi}_\vep(v_\ell)\Biggr)(\lv^{1/2})^{\sum_{\ell=1}^m\#_\ell} 
\left(\prod_{\ell=1}^{m-1}\frac{1}{2u_{\ell+1}+u_\ell}\right)\notag\\
&\quad\quad  +\widehat{ \wt{\sigma}}_1\int_{\R^4} \begin{bmatrix}
\d \widetilde{x}_1^{2}\\
 \d\wt{x}_1^1
 \end{bmatrix}_{\otimes}\varphi(\wt{x}_1^{\bi_0})\begin{bmatrix}
 P_{\wt{s}_1}(\wt{x}_0^{2},\wt{x}_1^{1}+\vep\two \wt{x}^{\bi_0}_1)\\
  P_{\wt{s}_1}(\wt{x}_0^{1},\wt{x}_1^{1})
\end{bmatrix}_\times\notag\\
&\quad\quad \times \int_{\R^4}
\begin{bmatrix}
\d x^{4}_1\\
\d x^{3}_1
\end{bmatrix}_\otimes\varphi(x^{\bi_1}_1)
 \begin{bmatrix}
\lv^{1/2}  P_{t'-t-\wt{\tau}_1+u_1}(\wt{x}_1^1+\vep\two \wt{x}_1^{\bi_0} ,x_1^{3}+\vep\two x^{\bi_1}_1)\\
\lv^{1/2}  P_{t'-t-\wt{\tau}_1+u_1}(\wt{x}_1^1,x_1^{3})
  \end{bmatrix}_\times\notag\\
  &\quad \quad \times  
  \Biggl(\int_{\R^2}\d z_1\s^{\lambda,\phi}_\vep(0,v_1;x_1^{\bi_1},z_1)\varphi(z_1) \sigma_1+ \widehat{\sigma}_1\Biggr)\notag\\
&\quad \quad \times\Biggl(\prod_{\scriptstyle \ell:\sigma_\ell=1\atop\scriptstyle 2\leq \ell\leq m}\overline{\s}^{\lambda,\phi}_\vep(v_\ell)\Biggr)(\lv^{1/2})^{\sum_{\ell=1}^m\#_\ell} 
\left(\prod_{\ell=1}^{m-1}\frac{1}{2u_{\ell+1}+u_\ell}\right).\label{mom:bdd2:aux11}
\end{align}

To simplify the right-hand side of \eqref{mom:bdd2:aux11}, we apply \eqref{ineq:gausstime} multiple times, starting with the integrals with respect to $\d x_1^4\otimes \d x_1^3$, and use the identity $\sum_{\ell=1}^m \#_\ell=2\lVert\widehat{\spin}\rVert$ from \eqref{sigmasum}. Hence, 
\begin{align}
& \int \d \widetilde{\mathcal G}^{\widetilde{\spin}}_{\vep;0}\int \d \widetilde{\mathcal G}^{\widetilde{\spin}}_{\vep;1}\int \d \mathcal G^{\bs \sigma}_{\vep;0}
\int \d \mathcal G^{\bs \sigma}_{\vep;1}
\cdots \int \d \mathcal G^{\bs \sigma}_{\vep;m}\1\notag\\
&\quad \leq\wt{\sigma}_1 \int_{\R^2}\d \wt{\overline{z}}\phi( \wt{\overline{z}})
 P_{2\wt{s}_1}(\wt{x}_0^{2}-\vep \wt{\overline{z}},\wt{x}_0^{1})\int_{\R^2} 
 \d \wt{z}_1\varphi(\wt{z}_1)
\s^{\lambda,\phi}_\vep (\wt{s}_1,\wt{\tau}_1; \wt{\overline{z}}/\two,\wt{z}_1)\notag\\
&\quad \quad \times \int_{\R^2}\d \overline{z} \phi(\overline{z})
  P_{2(t'-t-\wt{\tau}_1+u_1)}(\vep \two \wt{z}_1-\vep\overline{z})
  \Biggl(\int_{\R^2}\d z_1\s^{\lambda,\phi}_\vep(0,v_1;\overline{z}/\two,z_1)\varphi(z_1) \sigma_1+ \widehat{\sigma}_1\Biggr)\notag\\
  &\quad\quad  \times\Biggl(\prod_{\scriptstyle \ell:\sigma_\ell=1\atop\scriptstyle 2\leq \ell\leq m}\overline{\s}^{\lambda,\phi}_\vep(v_\ell)\Biggr)\lv^{\lVert\wh{\spin}\rVert}
\left(\prod_{\ell=1}^{m-1}\frac{1}{2u_{\ell+1}+u_\ell}\right)\notag\\
&\quad\quad  +\widehat{\wt{\sigma}}_1 \int_{\R^2}
\d \wt{\overline{z}} \phi( \wt{\overline{z}})P_{2\wt{s}_1}(\wt{x}_0^{2}-\vep \wt{\overline{z}},\wt{x}_0^{1}) 
\notag\\
  &\quad\quad  \times  \int_{\R^2}\d \overline{z}\phi(\overline{z})
 P_{2(t'-t-\wt{\tau}_1+u_1)}(\vep \wt{\overline{z}}-\vep\overline{z})
  \Biggl(\int_{\R^2}\d z_1\s^{\lambda,\phi}_\vep(0,v_1;\overline{z}/\two,z_1)\varphi(z_1) \sigma_1+ \widehat{\sigma}_1\Biggr)\notag\\
  &\quad\quad  \times\Biggl(\prod_{\scriptstyle \ell:\sigma_\ell=1\atop\scriptstyle 2\leq \ell\leq m}\overline{\s}^{\lambda,\phi}_\vep(v_\ell)\Biggr)\lv^{1+\lVert\wh{\spin}\rVert} 
\left(\prod_{\ell=1}^{m-1}\frac{1}{2u_{\ell+1}+u_\ell}\right).\label{mom:bdd2:aux1}
 \end{align}
We now obtain \eqref{mom:bdd2} for $\bi_1=(4,3)$ and $m\geq 1$ upon applying the bound 
\[
\max\{P_{2(t'-t-\wt{\tau}_1+u_1)}(\vep \two \wt{z}_1-\vep\overline{z}),
 P_{2(t'-t-\wt{\tau}_1+u_1)}(\vep \wt{\overline{z}}-\vep\overline{z})\}\leq (t'-t-\wt{\tau}_1+u_1)^{-1} 
\]
and the definition \eqref{def:os} of $\overline{\s}^{\lambda,\phi}_\vep(\cdot)$ 
to the right-hand side of \eqref{mom:bdd2:aux1}. 

To establish \eqref{mom:bdd2-1} for $\bi=(4,3)$ and $m\geq 1$, we use the following inequality, which results from a modification of the derivation of \eqref{mom:bdd2:aux11}:
  \begin{align}
&\E^{(\widetilde{B}^1,\widetilde{B}^2)}_{(\widetilde{x}_0^1,\widetilde{x}_0^2)}\Bigg[\int \d \mathcal G^\spin_{\vep;0}  \int \d \mathcal G^{\bs \sigma}_{\vep;1}
 \cdots\int \d \mathcal G^\spin_{\vep;m} \1\mid_{(x_0^3,x_0^4)=(\widetilde{B}^1_{t'-t},\widetilde{B}^2_{t'-t}) }\Bigg]\notag\\
 &\quad \leq  \int_{\R^4}
\begin{bmatrix}
\d x^{4}_1\\
\d x^{3}_1
\end{bmatrix}_\otimes  \varphi(x^{\bi_1}_1)
 \begin{bmatrix}
  P_{t'-t+u_1}(\wt{x}_0^2,x_1^{3}+\vep\two x^{\bi_1}_1)\\
  P_{t'-t+u_1}(\wt{x}_0^1,x_1^{3})
  \end{bmatrix}_\times\notag\\
  &\quad \quad \times 
  \Biggl(\int_{\R^2}\d z_1\s^{\lambda,\phi}_\vep(0,v_1;x_1^{\bi_1},z_1)\varphi(z_1) \sigma_1+ \widehat{\sigma}_1\Biggr)\Biggl(\prod_{\scriptstyle \ell:\sigma_\ell=1\atop\scriptstyle 2\leq \ell\leq m}\overline{\s}^{\lambda,\phi}_\vep(v_\ell)\Biggr)\lv^{\lVert\wh{\spin}\rVert} 
\left(\prod_{\ell=1}^{m-1}\frac{1}{2u_{\ell+1}+u_\ell}\right)\notag\\
&\quad =\int_{\R^2}\d \ol{z}\phi(\ol{z})  P_{2(t'-t+u_1)}(\wt{x}_0^2-\vep \ol{z},\wt{x}_0^1)\Biggl(\int_{\R^2}\d z_1\s^{\lambda,\phi}_\vep(0,v_1;\ol{z}/\two,z_1)\varphi(z_1) \sigma_1+ \widehat{\sigma}_1\Biggr)\notag\\
&\quad\quad  \times\Biggl(\prod_{\scriptstyle \ell:\sigma_\ell=1\atop\scriptstyle 2\leq \ell\leq m}\overline{\s}^{\lambda,\phi}_\vep(v_\ell)\Biggr)\lv^{\lVert\wh{\spin}\rVert}
\left(\prod_{\ell=1}^{m-1}\frac{1}{2u_{\ell+1}+u_\ell}\right).\notag
 \end{align}
Accordingly, we can apply a similar argument to that used below \eqref{mom:bdd2:aux1} to derive \eqref{mom:bdd2-1} for $\bi=(4,3)$ and $m \geq 1$. Thus, all the required inequalities for $\bi=(4,3)$ have been established.\medskip 

\noindent {\bf Step 5.} To prove \eqref{mom:bdd2} for $\bi_1=(3,2)$ and $m\geq 1$, we use the observation that exactly one of the two edges defining $\mathcal G^{\spin}_{\vep;0}$ can be extended to a left-maximal entanglement-free path in $\widetilde{\mathcal G}_\vep^{\widetilde{\spin}}\oplus \mathcal G_\vep^{\spin}$, and the duration of this edge increases from $u_1=s_1$ to $t'-t-\wt{\tau}_1+u_1$. Moreover, that edge in $\mathcal G^{\spin}_{\vep;0}$ must be the one whose left-end vertex is $(\wt{x}^1_2,t'-t)=(x_0^3,0)$, as illustrated in the two examples of Figure~\ref{fig:4}. We will not use the other extension via $(\wt{x}^2_2,t'-t)=(x_0^4,0)$, though. 

By using the above observation, the first bound for $\int \d \widetilde{\mathcal G}^{\widetilde{\spin}}_{\vep;0}\int \d \widetilde{\mathcal G}^{\widetilde{\spin}}_{\vep;1}\int \d \mathcal G^{\bs \sigma}_{\vep;0}
\int \d \mathcal G^{\bs \sigma}_{\vep;1}
\cdots \int \d \mathcal G^{\bs \sigma}_{\vep;m}\1$ can be obtained similarly as in \eqref{mom:bdd2:aux11}. Thus, we have the following: 
\begin{align}
& \int \d \widetilde{\mathcal G}^{\widetilde{\spin}}_{\vep;0}\int \d \widetilde{\mathcal G}^{\widetilde{\spin}}_{\vep;1}\int \d \mathcal G^{\bs \sigma}_{\vep;0}
\int \d \mathcal G^{\bs \sigma}_{\vep;1}
\cdots \int \d \mathcal G^{\bs \sigma}_{\vep;m}\1\notag\\
&\quad \leq \wt{\sigma}_1
 \int_{\R^4} \begin{bmatrix}
\d \widetilde{x}_1^{2}\\
 \d\wt{x}_1^1
 \end{bmatrix}_{\otimes}\varphi(\wt{x}_1^{\bi_0})\begin{bmatrix}
 P_{\wt{s}_1}(\wt{x}_0^{2},\wt{x}_1^{1}+\vep\two \wt{x}^{\bi_0}_1)\\
  P_{\wt{s}_1}(\wt{x}_0^{1},\wt{x}_1^{1})
\end{bmatrix}_\times\int_{\R^4} \begin{bmatrix}
\d \wt{z}_1'\\
 \d \wt{z}_1
 \end{bmatrix}_{\otimes} \varphi(\wt{z}_1)
\begin{bmatrix}
P_{\wt{s}_1,\wt{\tau}_1}(\two \wt{x}_1^{1}+\vep \wt{x}_1^{\bi_0},\wt{z}'_1)\\
\vspace{-.3cm}\\
\s^{\lambda,\phi}_\vep (\wt{s}_1,\wt{\tau}_1; \wt{x}_1^{\bi_{0}},\wt{z}_1)
\end{bmatrix}_\times\notag\\
&\quad\quad \times \int_{\R^4}
\begin{bmatrix}
\d x^{3}_1\\
\d x^{2}_1
\end{bmatrix}_\otimes \varphi(x^{\bi_1}_1)
 \begin{bmatrix}
  P_{t'-t-\wt{\tau}_1+u_1}\left(\frac{\wt{z}_1'-\vep \wt{z}_1}{\two},x_1^{2}+\vep \two x_1^{\bi_1}\right)\\
  P_{u_1}(x_0^2,x_1^2)
  \end{bmatrix}_\times\notag\\
&\quad\quad  \times  \Biggl(\int_{\R^2}\d z_1\s^{\lambda,\phi}_\vep(0,v_1;x_1^{\bi_1},z_1)\varphi(z_1) \sigma_1+ \widehat{\sigma}_1\Biggr)\notag\\
  &\quad \quad \times\Biggl(\prod_{\scriptstyle \ell:\sigma_\ell=1\atop\scriptstyle 2\leq \ell\leq m}\overline{\s}^{\lambda,\phi}_\vep(v_\ell)\Biggr)(\lv^{1/2})^{\sum_{\ell=1}^m\#_\ell} 
\left(\prod_{\ell=1}^{m-1}\frac{1}{2u_{\ell+1}+u_\ell}\right)\notag\\
&\quad\quad  +\widehat{\wt{\sigma}}_1 \int_{\R^4} \begin{bmatrix}
\d \widetilde{x}_1^{2}\\
 \d\wt{x}_1^1
 \end{bmatrix}_{\otimes}\varphi(\wt{x}_1^{\bi_0})\begin{bmatrix}
 P_{\wt{s}_1}(\wt{x}_0^{2},\wt{x}_1^{1}+\vep\two \wt{x}^{\bi_0}_1)\\
  P_{\wt{s}_1}(\wt{x}_0^{1},\wt{x}_1^{1})
\end{bmatrix}_\times\notag\\
&\quad\quad \times \int_{\R^4}
\begin{bmatrix}
\d x^{3}_1\\
\d x^{2}_1
\end{bmatrix}_\otimes\varphi(x^{\bi_1}_1)
 \begin{bmatrix}
\lv  P_{t'-t-\wt{\tau}_1+u_1}(\wt{x}_1^1,x_1^{2}+\vep \two x_1^{\bi_1})\\
  P_{u_1}(x_0^{2},x_1^{2})
  \end{bmatrix}_\times\notag
\\
&\quad\quad \times  \Biggl(\int_{\R^2}\d z_1\s^{\lambda,\phi}_\vep(0,v_1;x_1^{\bi_1},z_1)\varphi(z_1) \sigma_1+ \widehat{\sigma}_1\Biggr)\notag\\
&\quad\quad  \times\Biggl(\prod_{\scriptstyle \ell:\sigma_\ell=1\atop\scriptstyle 2\leq \ell\leq m}\overline{\s}^{\lambda,\phi}_\vep(v_\ell)\Biggr)(\lv^{1/2})^{\sum_{\ell=1}^m\#_\ell} 
\left(\prod_{\ell=1}^{m-1}\frac{1}{2u_{\ell+1}+u_\ell}\right).\label{mom:bdd2:aux2}
\end{align}
Here, when $\widehat{\wt{\sigma}}_1=1$, the prefactor $\lv $ of $\lv P_{t'-t-\wt{\tau}_1+u_1}(\wt{x}_1^1,x_1^{2}+\vep \two x_1^{\bi_1})$ emerges from the edge weights of $ \widetilde{\mathcal G}^{\widetilde{\spin}}_{\vep;1}$. 
To analyze the right-hand side of \eqref{mom:bdd2:aux2}, we apply \eqref{def:os} and \eqref{ineq:gausstime} to obtain the following two bounds:
\begin{align}
&\int_{\R^4}
\begin{bmatrix}
\d x^{3}_1\\
\d x^{2}_1
\end{bmatrix}_\otimes \varphi(x^{\bi_1}_1)
 \begin{bmatrix}
  P_{t'-t-\wt{\tau}_1+u_1}\left(\frac{\wt{z}_1'-\vep \wt{z}_1}{\two},x_1^{2}+\vep \two x_1^{\bi_1}\right)\\
  P_{u_1}(x_0^2,x_1^2)
  \end{bmatrix}_\times\notag\\
&\times  \Biggl(\int_{\R^2}\d z_1\s^{\lambda,\phi}_\vep(0,v_1;x_1^{\bi_1},z_1)\varphi(z_1) \sigma_1+ \widehat{\sigma}_1\Biggr)\notag\\
  &\quad =\int_{\R^2}
\d \overline{z} \phi(\overline{z})
  P_{t'-t-\wt{\tau}_1+2u_1}\left(\frac{\wt{z}_1'-\vep \wt{z}_1}{\two}-\vep \overline{z},x_0^2\right)
  \Biggl(\int_{\R^2}\d z_1\s^{\lambda,\phi}_\vep(0,v_1;\overline{z}/\two,z_1)\varphi(z_1) \sigma_1+ \widehat{\sigma}_1\Biggr)\notag\\
  &\quad \leq \int_{\R^2}\d \overline{z}
  \frac{\phi(\overline{z})}{t'-t-\wt{\tau}_1+2u_1}
  \Biggl(\int_{\R^2}\d z_1\s^{\lambda,\phi}_\vep(0,v_1;\overline{z}/\two,z_1)\varphi(z_1) \sigma_1+ \widehat{\sigma}_1\Biggr)\notag\\
  &\quad =\frac{\overline{\s}^{\lambda,\phi}_\vep(v_1)\sigma_1+\widehat{\sigma}_1}{t'-t-\wt{\tau}_1+2u_1},\label{mom:bdd2:aux2-1}
  \end{align}
 and
 \begin{align}
&\int_{\R^4}
\begin{bmatrix}
\d x^{3}_1\\
\d x^{2}_1
\end{bmatrix}_\otimes\varphi(x^{\bi_1}_1)
 \begin{bmatrix}
\lv  P_{t'-t-\wt{\tau}_1+u_1}(\wt{x}_1^1,x_1^{2}+\vep \two x_1^{\bi_1})\\
  P_{u_1}(x_0^{2},x_1^{2})
  \end{bmatrix}_\times  \Biggl(\int_{\R^2}\d z_1\s^{\lambda,\phi}_\vep(0,v_1;x_1^{\bi_1},z_1)\varphi(z_1) \sigma_1+ \widehat{\sigma}_1\Biggr)\notag\\
&\quad  \leq\lv \frac{\overline{\s}^{\lambda,\phi}_\vep(v_1)\sigma_1+\widehat{\sigma}_1}{t'-t-\wt{\tau}_1+2u_1}.\label{mom:bdd2:aux2-2}
\end{align}  
Additionally, the remainder of the right-hand side of \eqref{mom:bdd2:aux2} can be simplified by 
using \eqref{ineq:gausstime} and \eqref{sigmasum}. Hence, by \eqref{mom:bdd2:aux2}--\eqref{mom:bdd2:aux2-2}, we obtain the following bound:
\begin{align*}
&  \int \d \widetilde{\mathcal G}^{\widetilde{\spin}}_{\vep;0}\int \d \widetilde{\mathcal G}^{\widetilde{\spin}}_{\vep;1}\int \d \mathcal G^{\bs \sigma}_{\vep;0}
\int \d \mathcal G^{\bs \sigma}_{\vep;1}
\cdots \int \d \mathcal G^{\bs \sigma}_{\vep;m}\1\\
&\quad \leq \wt{\sigma}_1\int_{\R^2}\d \wt{\overline{z}}
 P_{2\wt{s}_1}(\wt{x}_0^{2}-\vep \wt{\overline{z}},\wt{x}_0^{1})\phi( \wt{\overline{z}})\int_{\R^2} 
 \d \wt{z}_1
\s^{\lambda,\phi}_\vep (\wt{s}_1,\wt{\tau}_1; \wt{\overline{z}}/\two,\wt{z}_1) \varphi(\wt{z}_1)
  \\
  &\quad\quad  \times  \Biggl(\prod_{\scriptstyle \ell:\sigma_\ell=1\atop\scriptstyle 1\leq \ell\leq m}\overline{\s}^{\lambda,\phi}_\vep(v_\ell)\Biggr)\lv^{\lVert\wh{\spin}\rVert}
\left(\frac{1}{t'-t-\wt{\tau}_1+2u_1}\right)\left(\prod_{\ell=1}^{m-1}\frac{1}{2u_{\ell+1}+u_\ell}\right)\notag\\
&\quad \quad +\widehat{\wt{\sigma}}_1\lv \int_{\R^2}
\d \wt{\overline{z}} P_{2\wt{s}_1}(\wt{x}_0^{2}-\vep \wt{\overline{z}},\wt{x}_0^{1})\phi( \wt{\overline{z}})\\
&\quad \quad \times\Biggl(\prod_{\scriptstyle \ell:\sigma_\ell=1\atop\scriptstyle 1\leq \ell\leq m}\overline{\s}^{\lambda,\phi}_\vep(v_\ell)\Biggr)\lv^{\lVert\wh{\spin}\rVert}
\left(\frac{1}{t'-t-\wt{\tau}_1+2u_1}\right)\left(\prod_{\ell=1}^{m-1}\frac{1}{2u_{\ell+1}+u_\ell}\right),
 \end{align*}
which is enough to get \eqref{mom:bdd2} in the case of $\bi_1=(3,2)$ and $m\geq 1$.

To derive \eqref{mom:bdd2-1} for $\bi=(3,2)$ and $m \geq 1$, observe that the bound in \eqref{mom:bdd2:aux2} can be adjusted to yield the following result:
 \begin{align}
&\E^{(\widetilde{B}^1,\widetilde{B}^2)}_{(\widetilde{x}_0^1,\widetilde{x}_0^2)}\Bigg[\int \d \mathcal G^\spin_{\vep;0}  \int \d \mathcal G^{\bs \sigma}_{\vep;1}
 \cdots\int \d \mathcal G^\spin_{\vep;m} \1\mid_{(x_0^3,x_0^4)=(\widetilde{B}^1_{t'-t},\widetilde{B}^2_{t'-t}) }\Bigg]\notag\\
&\quad \leq  \int_{\R^4}
\begin{bmatrix}
\d x^{3}_1\\
\d x^{2}_1
\end{bmatrix}_\otimes \varphi(x^{\bi_1}_1)
 \begin{bmatrix}
  P_{t'-t+u_1}(\wt{x}_0^1,x_1^{2}+\vep \two x_1^{\bi_1})\\
  P_{u_1}(x_0^2,x_1^2)
  \end{bmatrix}_\times\notag\\
&\quad\quad  \times  \Biggl(\int_{\R^2}\d z_1\s^{\lambda,\phi}_\vep(0,v_1;x_1^{\bi_1},z_1)\varphi(z_1) \sigma_1+ \widehat{\sigma}_1\Biggr)\Biggl(\prod_{\scriptstyle \ell:\sigma_\ell=1\atop\scriptstyle 2\leq \ell\leq m}\overline{\s}^{\lambda,\phi}_\vep(v_\ell)\Biggr)\lv^{\lVert\wh{\spin}\rVert}
\left(\prod_{\ell=1}^{m-1}\frac{1}{2u_{\ell+1}+u_\ell}\right)\notag\\
&\quad=\int_{\R^2}\d \overline{z}
  P_{t'-t+2u_1}\left(\wt{x}_0^1-\vep \overline{z},x_0^2\right) \phi(\overline{z})\notag\\
&\quad\quad  \times  \Biggl(\int_{\R^2}\d z_1\s^{\lambda,\phi}_\vep(0,v_1;\overline{z}/\two,z_1)\varphi(z_1) \sigma_1+ \widehat{\sigma}_1\Biggr)\Biggl(\prod_{\scriptstyle \ell:\sigma_\ell=1\atop\scriptstyle 2\leq \ell\leq m}\overline{\s}^{\lambda,\phi}_\vep(v_\ell)\Biggr)\lv^{\lVert\wh{\spin}\rVert}
\left(\prod_{\ell=1}^{m-1}\frac{1}{2u_{\ell+1}+u_\ell}\right).\notag
\end{align}
We have proved all of the bounds stated in Proposition~\ref{prop:mombdd1}. The proof is complete.
\end{proof}

Lemma~\ref{lem:gibdd} studies some elementary $k$-fold iterated integrals and obtains bounds, given by \eqref{CSZ:bdd}, which are enough for the forthcoming applications of Proposition~\ref{prop:mombdd1}. More specifically, we will use these bounds simultaneously for all $k\in \Bbb N$. Note that \eqref{CSZ:bdd} with $T=1$ and $p=1/2$ slightly improves the second inequality in \cite[Lemma~5.4]{CSZ:19a} by using smaller constants and a more straightforward argument. On the other hand, the method of that proof in \cite{CSZ:19a} may be of independent interest for its random-walk counting argument. 

\begin{lem}\label{lem:gibdd}
For all $p\in (0,1)$ and $a\in (0,\infty)$,
\begin{align}\label{int:sabddp}
\int_0^{1}\frac{\d s}{s^p(s+a)}\leq \frac{C_{\ref{int:sabddp}}(p)}{a^{p}},
\end{align}
and for all integers $k\geq 1$, $s_0,T\in (0,\infty)$ and $p\in (0,1)$, 
\begin{align}\label{CSZ:bdd}
\int_{(0,T)^k}\d (s_1,\cdots,s_k)\prod_{j=1}^k\frac{1}{2s_j+s_{j-1}}
\leq C_{\ref{int:sabddp}}(p)^k\left(\frac{T}{s_0}\right)^p.
\end{align}
In particular, we can take $C_{\ref{int:sabddp}}(\tfrac{1}{2})$ to be $\pi$.
\end{lem}

\begin{rmk}\label{rmk:gibdd}
The fact that the bounds in \eqref{CSZ:bdd} allow arbitrary $p\in (0,1)$ may suggest an improvement to the logarithmic order. Nevertheless, although the integral in \eqref{CSZ:bdd} with $k=1$ allows a sharp bound of the order $\log s_0^{-1}$ as $s_0\searrow 0$, the integral with $k=2$ is of the order $\log^2 s_0^{-1}$ as $s_0\searrow 0$. See also the first inequality in  \cite[Lemma~5.4]{CSZ:19a}. \hfill $\blacklozenge$
\end{rmk}

\begin{proof}[Proof of Lemma~\ref{lem:gibdd}]
To prove \eqref{int:sabddp} for general $p\in (0,1)$, first, we use the change of variables $as'=s$ to get
\[
\int_0^{1}\frac{\d s}{s^p(s+a)}=\frac{1}{a^{p}}\int_0^{1/a}\frac{\d s'}{(s')^p(s'+1)}.
\]
To bound the right-hand side, note that
\begin{align*}
\forall\; a\in [1,\infty),\quad 
\frac{1}{a^p}\int_0^{1/a}\frac{\d s'}{(s')^p(s'+1)}&\asymp \frac{1}{a^p}\int_0^{1/a}\frac{\d s'}{(s')^p}\asymp  \frac{C(p)}{a}\leq \frac{C(p)}{a^p};\\
\forall\; a\in (0,1),\quad 
\frac{1}{a^p}\int_1^{1/a}\frac{\d s'}{(s')^p(s'+1)}&\asymp \frac{1}{a^p}\int_1^{1/a}\frac{\d s'}{(s')^{p+1}}\asymp  C(p).
\end{align*}
For the case of $a\in (0,1)$, we still have $\int_0^{1}\d s'/[(s')^p(s'+1)]$ to take into account. Hence, \eqref{int:sabddp} holds and is sharp as $a\searrow 0$. Also, to see that we can take $C_{\ref{int:sabddp}}(\tfrac{1}{2})$ to be $\pi$, it is enough to note the following identity and $\lvert\tan^{-1}(\cdot)\rvert \leq \pi/2$:
\begin{align}\label{int:sa}
\int\frac{\d s}{\sqrt{s}(s+a)}=\frac{2}{\sqrt{a}}\tan^{-1}\left(\frac{\sqrt{s}}{\sqrt{a}}\right)+C,\quad s,a\in (0,\infty).
\end{align}

We prove \eqref{CSZ:bdd} now. By the change of variables $T\,\widetilde{s}_j=s_j$ for all $j$, it suffices to prove \eqref{CSZ:bdd} for $T=1$. In this case, the required inequality with $k=1$ holds because
\begin{align}
 \int_0^1\frac{\d s_1}{2s_{1}+s_{0}}
\leq \int_0^1\frac{\d s_1}{s_1^p(s_1+s_{0})}\leq \frac{C_{\ref{int:sabddp}}(p)}{s_0^{p}},\quad \forall\;s_0\in (0,\infty),\label{int:sk}
\end{align}
where the first inequality uses the bound $x+y\geq x^py$, which is equivalent to $x+(1-x^p)y\geq 0$, for all $x\in [0,1]$ and $y\geq 0$ by taking $x=s_1$ and $y=s_1+s_0$, and the last equality uses \eqref{int:sabddp}. To obtain  \eqref{CSZ:bdd} for $T=1$ and all $k\geq 2$ and $s_0\in (0,\infty)$, note that \eqref{int:sk} and \eqref{int:sabddp} imply the first and second inequalities below, respectively, where $k'\geq 2$: 
\begin{align*}
\int_{(0,1)^{k}}\d (s_1,\cdots,s_{k})\prod_{j=1}^{k}\frac{1}{2s_j+s_{j-1}}
&\leq \int_{(0,1)^{k-1}}\d (s_1,\cdots,s_{k-1})\frac{C_{\ref{int:sabddp}}(p)}{s_{k-1}^p} \prod_{j=1}^{k-1}\frac{1}{2s_j+s_{j-1}},\\
 \int_{(0,1)^{k'}}\d (s_1,\cdots,s_{k'})\frac{1}{s_{k'}^p} \prod_{j=1}^{k'}\frac{1}{2s_j+s_{j-1}}
 &\leq  \int_{(0,1)^{k'-1}}\d (s_1,\cdots,s_{k'-1})\frac{C_{\ref{int:sabddp}}(p)}{s_{k'-1}^p} \prod_{j=1}^{k'-1}\frac{1}{2s_j+s_{j-1}}.
 \end{align*}
By iterating the second inequality, we see that the first inequality gives
\[
\int_{(0,1)^{k}}\d (s_1,\cdots,s_{k})\prod_{j=1}^{k}\frac{1}{2s_j+s_{j-1}}\leq\int_{0}^1\frac{\d s_1 C_{\ref{int:sabddp}}(p)^{k-1} }{s_{1}^p(2s_1+s_{0})}\leq \frac{ C_{\ref{int:sabddp}}(p)^k}{s_0^p},
\]
where the last inequality uses \eqref{int:sabddp}. We have proved \eqref{CSZ:bdd} for $T=1$ and all $k\geq 1$ and $s_0\in (0,\infty)$, and hence, \eqref{CSZ:bdd} for general $T$ as well. 
\end{proof}

We now turn to establishing the main bounds needed for the proof of Theorem~\ref{thm:mombdd}, with Lemma~\ref{lem:FKmom1} again connecting these bounds to the theorem. Unlike the earlier moment bounds (e.g., \cite{C:DBG}), \eqref{mom:finalbdd-1-1} and \eqref{mom:finalbdd-2-1}  stated below involve integrated forms with $\int_{T_1}^{T_2}\d t\int_{t}^{T_2}\d t'$. For further details on their application in \eqref{mom:finalbdd-1-1}, see Remark~\ref{rmk:useint} following \eqref{mom:finalbdd-1-1:2}.

\begin{prop}\label{prop:mombdd2}
Let $\vep\in (0,\ovl]$, $0< t<t'<\infty$, $x_0^1,x_0^2,\widetilde{x}_0^1,\widetilde{x}_0^2\in \R^2$ with $x_0^2=x_0^1+\vep x_0^0$ and $\wt{x}_0^2=\wt{x}_0^1+\vep \wt{x}_0^0$ for $x_0^0,\wt{x}_0^0\in \R^2$,  $m\in \Bbb N$, $\widetilde{\spin}\in \{0,1\},\spin\in \{0,1\}^m$, and $\bi_1\neq \cdots \neq \bi_m$ with $\bi_1,\cdots,\bi_m\in \mc E_4$. Consider the corresponding iterated graphical integrals in \eqref{mom:series}. Additionally, let $q\geq 0$, $0<T<\infty$, and $p\in (0,1)$, and recall $\mathfrak S^{\lambda,\phi}_\vep(q)$ defined in \eqref{def:Svep}. 
\begin{itemize}
\item [\rm (1$\cc$)] If $\bi_1=(2,1)$, then for all $0\leq T_1< T_2\leq T$,
\begin{align}
&\lv^2\int_{T_1}^{T_2}\d t\int_{t}^{T_2}\d t'\int_0^{t'-t}\d \wt{s}_1\int_{[\wt{s}_1,t'-t)}\d \wt{\tau_1}\int_{\Delta_m(t)}\d \bs u_m\otimes \d \bs v_m\notag\\
& \times\Bigg(\prod_{\ell:\wt{\sigma}_\ell=0}\delta_{\wt{s}_\ell}(\wt{\tau}_\ell)\Bigg)\Bigg(\prod_{\scriptstyle \ell:\sigma_\ell=0 \atop \scriptstyle 1\leq \ell\leq m}\delta_0(v_\ell)\Bigg)\int \d \widetilde{\mathcal G}^{\widetilde{\bs \sigma}}_{\vep;0}\int \d \widetilde{\mathcal G}^{\widetilde{\bs \sigma}}_{\vep;1} \int \d \mathcal G^{\bs \sigma}_{\vep;0} \int \d \mathcal G^{\bs \sigma}_{\vep;1}\cdots \int \d \mathcal G^{\bs \sigma}_{\vep;m}\1\notag\\
&\quad \leq C(\lambda,\lVert\varphi\rVert_\infty,q,T,p)
(T_2-T_1)^{2-p}\notag\\
&\quad \quad\times \left(\int_0^T\d \wt{s}_1\int_{\R^2}\d \wt{\ol{z}}P_{2\wt{s}_1}(\vep \wt{x}_0^0-\vep \wt{\ol{z}})\phi(\wt{\ol{z}})
\right) \left(\int_0^T \d u_1\int_{\R^2}\d \ol{z}P_{2u_1}(\vep x_0^0-\vep\ol{z})\phi(\ol{z})\right)\notag\\
&\quad \quad \times C(T,p)^m\mathfrak S^{\lambda,\phi}_\vep(q)^{\lVert\spin\rVert-\sigma_1} (\wt{\sigma}_1+\lv\wh{\wt{\sigma}}_1)\lv^{2+\lVert\wh{\spin}\rVert},\label{mom:finalbdd-1-1}
\end{align}
and for all $0<t<t'\leq T$, 
\begin{align}
&\lv^2\E^{(\widetilde{B}^1,\widetilde{B}^2)}_{(\widetilde{x}_0^1,\widetilde{x}_0^2)}\Bigg[
\int_{\Delta_m(t)}\d \bs u_m\otimes \d \bs v_m\Bigg(\prod_{\scriptstyle \ell:\sigma_\ell=0 \atop \scriptstyle 1\leq \ell\leq m}\delta_0(v_\ell)\Bigg)\notag\\
&\quad \times \int \d \mathcal G^{\bs \sigma}_{\vep;0} \int \d \mathcal G^{\bs \sigma}_{\vep;1}\cdots \int \d \mathcal G^{\bs \sigma}_{\vep;m}\1\mid_{(x_0^3,x_0^4)=(\widetilde{B}^1_{t'-t},\widetilde{B}^2_{t'-t}) }\Bigg]\notag\\
&\quad \leq  \frac{C(\lambda,\lVert\varphi\rVert_\infty,q,T,p)}{(t'-t)^p}
 \left(\int_0^T \d u_1\int_{\R^2}\d \ol{z}P_{2u_1}(\vep x_0^0-\vep\ol{z})\phi(\ol{z})\right)\notag\\
 &\quad\quad \times C(T,p)^m\mathfrak S^{\lambda,\phi}_\vep(q)^{\lVert\spin\rVert-\sigma_1} \lv^{2+\lVert\wh{\spin}\rVert}.
\label{mom:finalbdd-1-2}
\end{align}

\item [\rm (2$\cc$)] If $\bi_1=(4,3)$ or $\bi_1=(3,2)$, then for all $0\leq T_1< T_2\leq T$,
\begin{align}
&\lv^2\int_{T_1}^{T_2}\d t\int_{t}^{T_2}\d t'\int_0^{t'-t}\d \wt{s}_1\int_{[\wt{s}_1,t'-t)}\d \wt{\tau_1}\int_{\Delta_m(t)}\d \bs u_m\otimes \d \bs v_m\notag\\
&\times\Bigg(\prod_{\ell:\wt{\sigma}_\ell=0}\delta_{\wt{s}_\ell}(\wt{\tau}_\ell)\Bigg)\Bigg(\prod_{\scriptstyle \ell:\sigma_\ell=0 \atop \scriptstyle 1\leq \ell\leq m}\delta_0(v_\ell)\Bigg)\int \d \widetilde{\mathcal G}^{\widetilde{\bs \sigma}}_{\vep;0}\int \d \widetilde{\mathcal G}^{\widetilde{\bs \sigma}}_{\vep;1} \int \d \mathcal G^{\bs \sigma}_{\vep;0} \int \d \mathcal G^{\bs \sigma}_{\vep;1}\cdots \int \d \mathcal G^{\bs \sigma}_{\vep;m}\1\notag\\
&\quad \leq  C(\lambda,\lVert\varphi\rVert_\infty,q,T,p)
(T_2-T_1)^{2-p}\Bigg(\int_0^T\d \wt{s}_1\int_{\R^2}\d \wt{\ol{z}}P_{2\wt{s}_1}(\vep \wt{x}_0^0-\vep \wt{\ol{z}})\phi(\wt{\ol{z}})
\Bigg)\notag\\
&\quad\quad \times C(T,p)^m\mathfrak S^{\lambda,\phi}_\vep(q)^{\lVert\spin\rVert} (\wt{\sigma}_1+\lv\wh{\wt{\sigma}}_1)\lv^{2+\lVert\wh{\spin}\rVert}, \label{mom:finalbdd-2-1}
\end{align} 
and for all $0<t<t'\leq T$, 
\begin{align}
&\lv^2\E^{(\widetilde{B}^1,\widetilde{B}^2)}_{(\widetilde{x}_0^1,\widetilde{x}_0^2)}\Bigg[
\int_{\Delta_m(t)}\d \bs u_m\otimes \d \bs v_m\Bigg(\prod_{\scriptstyle \ell:\sigma_\ell=0 \atop \scriptstyle 1\leq \ell\leq m}\delta_0(v_\ell)\Bigg)\notag \\
&\quad \times \int \d \mathcal G^{\bs \sigma}_{\vep;0} \int \d \mathcal G^{\bs \sigma}_{\vep;1}\cdots \int \d \mathcal G^{\bs \sigma}_{\vep;m}\1\mid_{(x_0^3,x_0^4)=(\widetilde{B}^1_{t'-t},\widetilde{B}^2_{t'-t}) }\Bigg]\notag\\
&\quad \leq  \frac{C(\lambda,\lVert\varphi\rVert_\infty,q,T,p)}{(t'-t)^p}C(T,p)^m
\mathfrak S^{\lambda,\phi}_\vep(q)^{\lVert\spin\rVert} \lv^{2+\lVert\wh{\spin}\rVert}.
\label{mom:finalbdd-2-2}
\end{align}
\end{itemize}
\end{prop}

\begin{proof}[Proof of Proposition~\ref{prop:mombdd2} (1$\cc$)]
By the present choice of $(x_0^1,x_0^2,\widetilde{x}_0^1,\widetilde{x}_0^2)$, we can write
\begin{align}\label{diff:0}
\wt{x}_0^2-\vep\wt{\overline{z}}-\wt{x}_0^1=\vep\wt{x}^0_0-\vep \wt{\overline{z}},\quad x_0^2-\vep\overline{z}-x_0^1=\vep x^0_0-\vep \overline{z}.
\end{align}
Hence, by \eqref{mom:bdd1}, 
\begin{align}
&\lv^2\int_0^{t'-t}\d \wt{s}_1\int_{[\wt{s}_1,t'-t)}\d \wt{\tau_1}\int_{\Delta_m(t)}\d \bs u_m\otimes \d \bs v_m\Bigg(\prod_{\ell:\wt{\sigma}_\ell=0}\delta_{\wt{s}_\ell}(\wt{\tau}_\ell)\Bigg)\Bigg(\prod_{\scriptstyle \ell:\sigma_\ell=0 \atop \scriptstyle 1\leq \ell\leq m}\delta_0(v_\ell)\Bigg)\notag\\
&\times\int \d \widetilde{\mathcal G}^{\widetilde{\bs \sigma}}_{\vep;0}\int \d \widetilde{\mathcal G}^{\widetilde{\bs \sigma}}_{\vep;1} \int \d \mathcal G^{\bs \sigma}_{\vep;0} \int \d \mathcal G^{\bs \sigma}_{\vep;1}\cdots \int \d \mathcal G^{\bs \sigma}_{\vep;m}\1\notag\\
&\quad \leq \lv^2\int_0^{t'-t}\d \wt{s}_1\int_{[\wt{s}_1,t'-t)}\d \wt{\tau_1}\int_{\Delta_m(t)}\d \bs u_m\otimes \d \bs v_m\Bigg(\prod_{\ell:\wt{\sigma}_\ell=0}\delta_{\wt{s}_\ell}(\wt{\tau}_\ell)\Bigg)\Bigg(\prod_{\scriptstyle \ell:\sigma_\ell=0 \atop \scriptstyle 1\leq \ell\leq m}\delta_0(v_\ell)\Bigg)\notag\\
&\quad\quad \times \int_{\R^2}\d \wt{\overline{z}}
 P_{2\wt{s}_1}(\vep\wt{x}^0_0-\vep \wt{\overline{z}})\phi(\wt{\overline{z}}) \left(\int_{\R^2} \d \wt{z}_1\s^{\lambda,\phi}_\vep(\wt{s}_1,\wt{\tau}_1;\wt{\overline{z}}/\two,\wt{z}_1) \varphi(\wt{z}_1)\wt{\sigma}_1+  \widehat{\wt{\sigma}}_1\right)\notag\\
&\quad\quad \times \int_{\R^2}
\d \overline{z}
  P_{2u_1}(\vep x^0_0-\vep \overline{z})  \phi(\overline{z})
  \Biggl(\int_{\R^2}\d z_1\s^{\lambda,\phi}_\vep(0,v_1;\overline{z}/\two,z_1)\varphi(z_1) \sigma_1+ \widehat{\sigma}_1\Biggr)\notag\\
&\quad\quad  \times\Biggl(\prod_{\scriptstyle \ell:\sigma_\ell=1\atop \scriptstyle 2\leq \ell\leq m}\overline{\s}^{\lambda,\phi}_\vep(v_\ell)\Biggr)\lv^{\lVert\wh{\wt{\spin}}\rVert+\lVert\wh{\spin}\rVert} \notag\\
&\quad \quad \times \left(\frac{\1_{m\geq 2}}{t'-t-\widetilde{\tau}_1+2u_2+u_1}+\1_{m=1}\right)
\left(\prod_{\ell=2}^{m-1}\frac{1}{2u_{\ell+1}+u_\ell}\right).\label{mom:bdd1-1:1}
\end{align}

To simplify the right-hand side of \eqref{mom:bdd1-1:1}, we analyze the integrals with respect to $\bs v_m$ and $\bs u_m$ separately as follows. With $\Delta_m(t)_{\bs u_m}$ denoting the section of the set $\Delta_m(t)$, defined in \eqref{def:Deltam(t)},
for fixed $\bs u_m$, we obtain from the simple inequality $1\leq \e^{qt}\e^{-q\sum_{\ell=1}^m v_\ell}$ for any $q\geq 0$ that 
\begin{align}
&\int_{\Delta_m(t)_{\bs u_m}}\d \bs v_m\Bigg(\prod_{\scriptstyle \ell:\sigma_\ell=0 \atop \scriptstyle 1\leq \ell\leq m}\delta_0(v_\ell)\Bigg) \Biggl(\int_{\R^2}\d z_1\s^{\lambda,\phi}_\vep(0,v_1;\overline{z}/\two,z_1)\varphi(z_1) \sigma_1+ \widehat{\sigma}_1\Biggr)\Biggl(\prod_{\scriptstyle \ell:\sigma_\ell=1\atop \scriptstyle 2\leq \ell\leq m}\overline{\s}^{\lambda,\phi}_\vep(v_\ell)\Biggr)\notag \\
&\quad \leq  \e^{qt}\Biggl(\int_0^t \d v_1\int_{\R^2}\d z_1\s^{\lambda,\phi}_\vep(0,v_1;\overline{z}/\two,z_1)\varphi(z_1) \sigma_1+ \widehat{\sigma}_1\Biggr) \Biggl(\prod_{\scriptstyle \ell:\sigma_\ell=1\atop \scriptstyle 2\leq \ell\leq m}\int_0^t \d v_\ell\e^{-qv_\ell}\overline{\s}^{\lambda,\phi}_\vep(v_\ell)\Biggr)\notag\\
&\quad \leq C(\lambda,\lVert\varphi\rVert_\infty,q,T)
\mathfrak S^{\lambda,\phi}_\vep(q)^{\lVert\spin\rVert-\sigma_1},\label{mom:bdd1-1:1-100000}
\end{align}
where the last inequality uses \eqref{ineq:apriori} and the definition \eqref{def:Svep} of $\mathfrak S^{\lambda,\phi}_\vep(q)$. Note that the right-hand side of the foregoing inequality is independent of $\bs u_m$. Also, 
\begin{align}
&\int_{(0,t)^m}\d \bs u_m \int_{\R^2}
\d \overline{z} P_{2u_1}(\vep x^0_0-\vep \overline{z}) \phi(\overline{z})  \left(\frac{\1_{m\geq 2}}{t'-t-\widetilde{\tau}_1+2u_2+u_1}+\1_{m=1}\right)\left(\prod_{\ell=2}^{m-1}\frac{1}{2u_{\ell+1}+u_\ell}\right)\notag\\
&\quad \leq \int_0^{t}\d u_1 \int_{\R^2}
\d \overline{z} P_{2u_1}( \vep x_0^0-\vep \overline{z}) \phi(\overline{z})\Bigg[\int_{(0,t)^{m-1}}\d (u_2,\cdots,u_m)\left(\frac{\1_{m\geq 2}}{t'-t-\widetilde{\tau}_1+2u_2}\right)\notag\\
&\quad\quad  \times\left(\prod_{\ell=2}^{m-1}\frac{1}{2u_{\ell+1}+u_\ell}\right) +\1_{m=1}\Bigg]\notag\\
&\quad\leq C(T,p) \int_0^{t}\d u_1 \int_{\R^2}
\d \overline{z}P_{2u_1}( \vep x_0^0-\vep \overline{z}) \phi(\overline{z})\left(\frac{C(T,p)^{m}}{(t'-t-\wt{\tau}_1)^p}\1_{m\geq 2}+\1_{m=1}\right)\notag\\
&\quad \leq C(T,p)\int_0^{t}\d u_1\int_{\R^2}
\d \overline{z} P_{2u_1}( \vep x_0^0-\vep \overline{z})  \phi(\overline{z})\frac{C(T,p)^{m}}{(t'-t-\wt{\tau}_1)^p},\label{mom:bdd1-1:1-100}
\end{align} 
where the next to the last inequality follows from \eqref{CSZ:bdd} with the choice of
\[
k\equiv m-1,\;s_0=t'-t-\wt{\tau}_1,\;s_1\equiv u_2,\;s_2\equiv u_3,\cdots,\;s_{m-1}\equiv u_m.
\]
Hence, by considering separately $\wt{\sigma}_1=1$ and $\wt{\sigma}_1=0$, the above two displays applied to 
\eqref{mom:bdd1-1:1} give the following bound:
\begin{align}
&\lv^2\int_0^{t'-t}\d \wt{s}_1\int_{[\wt{s}_1,t'-t)}\d \wt{\tau_1}\int_{\Delta_m(t)}\d \bs u_m\otimes \d \bs v_m\Bigg(\prod_{\ell:\wt{\sigma}_\ell=0}\delta_{\wt{s}_\ell}(\wt{\tau}_\ell)\Bigg)\Bigg(\prod_{\scriptstyle \ell:\sigma_\ell=0 \atop \scriptstyle 1\leq \ell\leq m}\delta_0(v_\ell)\Bigg)\notag\\
&\times\int \d \widetilde{\mathcal G}^{\widetilde{\bs \sigma}}_{\vep;0}\int \d \widetilde{\mathcal G}^{\widetilde{\bs \sigma}}_{\vep;1} \int \d \mathcal G^{\bs \sigma}_{\vep;0} \int \d \mathcal G^{\bs \sigma}_{\vep;1}\cdots \int \d \mathcal G^{\bs \sigma}_{\vep;m}\1\notag\\
&\quad \leq\wt{\sigma}_1C(\lambda,\lVert\varphi\rVert_\infty,q,T,p)
\int_0^{t'-t}\d \wt{s}_1\int_{\wt{s}_1}^{t'-t}\d \wt{\tau_1}\notag\\
&\quad\quad \times \int_{\R^2}\d \wt{\overline{z}}
 P_{2\wt{s}_1}(\vep\wt{x}^0_0-\vep \wt{\overline{z}})\phi(\wt{\overline{z}}) \int_{\R^2} \d \wt{z}_1\frac{\s^{\lambda,\phi}_\vep(\wt{s}_1,\wt{\tau}_1;\wt{\overline{z}}/\two,\wt{z}_1) \varphi(\wt{z}_1)}{(t'-t-\wt{\tau}_1)^p}\notag\\
&\quad\quad \times  \int_0^{t}\d u_1\int_{\R^2}\d \ol{z}P_{2u_1}( \vep x_0^0-\vep \overline{z})\phi(\ol{z})C(T,p)^m\mathfrak S^{\lambda,\phi}_\vep(q)^{\lVert\spin\rVert-\sigma_1} \lv^{2+\lVert\wh{\spin}\rVert} \notag\\
&\quad\quad +\wh{\wt{\sigma}}_1C(\lambda,\lVert\varphi\rVert_\infty,q,T,p) 
\int_0^{t'-t}\d \wt{s}_1 \int_{\R^2}\d \wt{\overline{z}}
 P_{2\wt{s}_1}(\vep\wt{x}^0_0-\vep \wt{\overline{z}})\phi(\wt{\overline{z}})\frac{1}{(t'-t-\wt{s}_1)^p}\notag\\
&\quad\quad \times  \int_0^{t}\d u_1\int_{\R^2}\d \ol{z}P_{2u_1}( \vep x_0^0-\vep \overline{z})\phi(\ol{z})C(T,p)^m\mathfrak S^{\lambda,\phi}_\vep(q)^{\lVert\spin\rVert-\sigma_1}\lv^{3+\lVert\wh{\spin}\rVert} .\label{mom:finalbdd-1-1:2}
\end{align}

\begin{rmk}\label{rmk:useint}
To bound the right-hand side of \eqref{mom:finalbdd-1-1:2}, our main focus is on the first term, especially the integral with respect to $\d \wt{s}_1\d \wt{\tau}_1$. At this stage, it is necessary to introduce $\int_{T_1}^{T_2}\d t\int_t^{T_2}\d t'$ to bypass singularities and deal with a complex convolution specified below.\hfill $\blacklozenge$
\end{rmk}

To proceed, we show the following bound: for $0\leq T_1<T_2\leq T$,
\begin{align}
&\int_{T_1}^{T_2}\d t\int_t^{T_2}\d t'\int_0^{t'-t}\d \wt{s}_1\int_{\wt{s}_1}^{t'-t}\d \wt{\tau_1}  \int_{\R^2}\d \wt{\overline{z}}
 P_{2\wt{s}_1}(\vep\wt{x}^0_0-\vep \wt{\overline{z}})\phi(\wt{\overline{z}})
\int_{\R^2} \d \wt{z}_1\frac{\s^{\lambda,\phi}_\vep(\wt{s}_1,\wt{\tau}_1;\wt{\overline{z}}/\two,\wt{z}_1) \varphi(\wt{z}_1)}{(t'-t-\wt{\tau}_1)^p}\notag\\
&\quad \leq C(\lambda,\lVert\varphi\rVert_\infty,T,p)(T_2-T_1)^{2-p} \left(\int_0^T\d \wt{s}_1\int_{\R^2}\d \wt{\ol{z}}P_{2\wt{s}_1}(\vep \wt{x}_0^0-\vep \wt{\ol{z}})\phi(\wt{\ol{z}})
\right).\label{T1T2:intfinal}
\end{align}
First, write
\begin{align}
&\int_0^{t'-t}\d \wt{s}_1\int_{\wt{s}_1}^{t'-t}\d \wt{\tau_1}  
 P_{2\wt{s}_1}(\vep\wt{x}^0_0-\vep \wt{\overline{z}})\phi(\wt{\overline{z}})
\frac{\s^{\lambda,\phi}_\vep(\wt{s}_1,\wt{\tau}_1;\wt{\overline{z}}/\two,\wt{z}_1) \varphi(\wt{z}_1)}{(t'-t-\wt{\tau}_1)^p}\notag\\
&\quad =\int_0^{t'-t}\d \wt{s}_1\int_{0}^{t'-t-\wt{s}_1}\d \wt{\tau_1}  
 P_{2\wt{s}_1}(\vep\wt{x}^0_0-\vep \wt{\overline{z}})\phi(\wt{\overline{z}})
\frac{\s^{\lambda,\phi}_\vep(0,\wt{\tau}_1;\wt{\overline{z}}/\two,\wt{z}_1) \varphi(\wt{z}_1)}{(t'-t-\wt{s}_1-\wt{\tau}_1)^p},\notag
\end{align}
which follows by changing variables. We apply $\int_{T_1}^{T_2}\d t\int_t^{T_2}\d t'$ to both sides of the foregoing equality. Hence,
\begin{align}
&\int_{T_1}^{T_2}\d t\int_t^{T_2}\d t'\int_0^{t'-t}\d \wt{s}_1\int_{\wt{s}_1}^{t'-t}\d \wt{\tau_1}  \int_{\R^2}\d \wt{\overline{z}}
 P_{2\wt{s}_1}(\vep\wt{x}^0_0-\vep \wt{\overline{z}})\phi(\wt{\overline{z}})
\int_{\R^2} \d \wt{z}_1\frac{\s^{\lambda,\phi}_\vep(\wt{s}_1,\wt{\tau}_1;\wt{\overline{z}}/\two,\wt{z}_1) \varphi(\wt{z}_1)}{(t'-t-\wt{\tau}_1)^p}\notag\\
&\quad =\int_{T_1}^{T_2}\d t\int_t^{T_2}\d t'\int_0^{t'-t}\d \wt{s}_1\int_{0}^{t'-t-\wt{s}_1}\d \wt{\tau_1}  
\int_{\R^2}\d \wt{\overline{z}} P_{2\wt{s}_1}(\vep\wt{x}^0_0-\vep \wt{\overline{z}})\phi(\wt{\overline{z}})\notag\\
&\quad \quad \times \int_{\R^2} \d \wt{z}_1
\frac{\s^{\lambda,\phi}_\vep(0,\wt{\tau}_1;\wt{\overline{z}}/\two,\wt{z}_1) \varphi(\wt{z}_1)}{(t'-t-\wt{s}_1-\wt{\tau}_1)^p}\notag\\
&\quad =\int_{\R^2}\d \wt{\overline{z}}\int_{\R^2} \d \wt{z}_1\int_{T_1}^{T_2}\d t\int_t^{T_2}\d t'\int_0^{t'-t}\d \wt{s}_1\int_{0}^{t'-t-\wt{s}_1}\d \wt{\tau_1}  
P_{2\wt{s}_1}(\vep\wt{x}^0_0-\vep \wt{\overline{z}})\phi(\wt{\overline{z}})\notag\\
&\quad \quad \times
\frac{\s^{\lambda,\phi}_\vep(0,\wt{\tau}_1;\wt{\overline{z}}/\two,\wt{z}_1) \varphi(\wt{z}_1)}{(t'-t-\wt{s}_1-\wt{\tau}_1)^p},\label{mom:finalbdd-1-1:3}
\end{align}
and the right-hand side 
takes the form of the left-hand side of the following inequality:  
\begin{align}
&\int_{\R^2}\d \wt{\overline{z}}\int_{\R^2} \d \wt{z}_1\int_{T_1}^{T_2} \d t\int_{t}^{T_2} \d t'(f\star g\star h)(t'-t)\notag\\
&\quad \leq \int_{\R^2}\d \wt{\overline{z}}\int_{\R^2} \d \wt{z}_1\int_{T_1}^{T_2}\d t\left(\int_0^{T_2-t}\d t'f(t')\right)\left(\int_0^{T_2-t}\d t'g(t')\right)\left(\int_0^{T_2-t}\d t'h(t')\right).\label{mom:finalbdd-1-1:3-1}
\end{align}
Specifically, we take $f(t')\equiv P_{2t'}(\vep\wt{x}^0_0-\vep \wt{\overline{z}})\phi(\wt{\overline{z}})$, $g(t')\equiv \s^{\lambda,\phi}_\vep(0,t';\wt{\overline{z}}/\two,\wt{z}_1)\varphi(\wt{z}_1)$ and $h(t')\equiv 1/(t')^p$. 
Moreover,
\begin{gather}
\int_0^{T_2-t}\d t'\int_{\R^2} \d \wt{z}_1\s^{\lambda,\phi}_\vep(0,t';\wt{\overline{z}}/\two,\wt{z}_1) \varphi(\wt{z}_1)\leq C(\lambda,\lVert\varphi\rVert_\infty,T),\label{T1T2:int0}\\
\int_{T_1}^{T_2}\d t\int_0^{T_2-t}\frac{\d t'}{(t')^p}\leq C(p) \int_{T_1}^{T_2} \d t(T_2-t)^{1-p}\leq C(p) (T_2-T_1)^{2-p},\label{T1T2:int}
\end{gather}
where \eqref{T1T2:int0} uses \eqref{ineq:apriori}. Combining \eqref{mom:finalbdd-1-1:3}--\eqref{T1T2:int} then proves \eqref{T1T2:intfinal}.

We are ready to finish the proof of \eqref{mom:finalbdd-1-1}. 
By \eqref{T1T2:intfinal} and 
a similar, but simpler, bound for the second term in \eqref{mom:finalbdd-1-1:2}, we obtain from \eqref{mom:finalbdd-1-1:2} that
\begin{align*}
&\lv^2\int_{T_1}^{T_2}\d t\int_{t}^{T_2}\d t'\int_0^{t'-t}\d \wt{s}_1\int_{[\wt{s}_1,t'-t)}\d \wt{\tau_1}\int_{\Delta_m(t)}\d \bs u_m\otimes \d \bs v_m\notag\\
& \times\Bigg(\prod_{\ell:\wt{\sigma}_\ell=0}\delta_{\wt{s}_\ell}(\wt{\tau}_\ell)\Bigg)\Bigg(\prod_{\scriptstyle \ell:\sigma_\ell=0 \atop \scriptstyle 1\leq \ell\leq m}\delta_0(v_\ell)\Bigg)\int \d \widetilde{\mathcal G}^{\widetilde{\bs \sigma}}_{\vep;0}\int \d \widetilde{\mathcal G}^{\widetilde{\bs \sigma}}_{\vep;1} \int \d \mathcal G^{\bs \sigma}_{\vep;0} \int \d \mathcal G^{\bs \sigma}_{\vep;1}\cdots \int \d \mathcal G^{\bs \sigma}_{\vep;m}\1\notag\\
&\quad \leq C(\lambda,\lVert\varphi\rVert_\infty,q,T,p)
(T_2-T_1)^{2-p}\left(\int_0^T\d \wt{s}_1\int_{\R^2}\d \wt{\ol{z}}P_{2\wt{s}_1}(\vep \wt{x}_0^0-\vep \wt{\ol{z}})\phi(\wt{\ol{z}})
\right)\\
&\quad \quad \times \left(\int_0^T \d u_1\int_{\R^2}\d \ol{z}P_{2u_1}(\vep x_0^0-\vep\ol{z})\phi(\ol{z})\right)C(T,p)^m\mathfrak S^{\lambda,\phi}_\vep(q)^{\lVert\spin\rVert-\sigma_1} (\wt{\sigma}_1+\lv\wh{\wt{\sigma}}_1)\lv^{2+\lVert\wh{\spin}\rVert},
\end{align*}
as required in \eqref{mom:finalbdd-1-1}. 

The proof of \eqref{mom:finalbdd-1-2} is similar. For $0<t<t'\leq T$ and $q\geq 0$, we obtain the first inequality below by using \eqref{mom:bdd1-1} and \eqref{diff:0}:
\begin{align}
&\lv^2\E^{(\widetilde{B}^1,\widetilde{B}^2)}_{(\widetilde{x}_0^1,\widetilde{x}_0^2)}\Bigg[
\int_{\Delta_m(t)}\d \bs u_m\otimes \d \bs v_m\Bigg(\prod_{\scriptstyle \ell:\sigma_\ell=0 \atop \scriptstyle 1\leq \ell\leq m}\delta_0(v_\ell)\Bigg)\notag\\
&\quad \times \int \d \mathcal G^{\bs \sigma}_{\vep;0} \int \d \mathcal G^{\bs \sigma}_{\vep;1}\cdots \int \d \mathcal G^{\bs \sigma}_{\vep;m}\1\Big|_{(x_0^3,x_0^4)=(\widetilde{B}^1_{t'-t},\widetilde{B}^2_{t'-t}) }\Bigg]\notag\\
&\quad \leq\lv^2 \int_{\Delta_m(t)}\d \bs u_m\otimes \d \bs v_m\Bigg(\prod_{\scriptstyle \ell:\sigma_\ell=0 \atop \scriptstyle 1\leq \ell\leq m}\delta_0(v_\ell)\Bigg) 
\int_{\R^2}
\d \overline{z}
  P_{2u_1}(\vep x_0^{0}-\vep \overline{z})  \phi(\overline{z})\notag\\
&\quad\quad \times  \Biggl(\int_{\R^2}\d z_1\s^{\lambda,\phi}_\vep(0,v_1;\overline{z}/\two,z_1)\varphi(z_1) \sigma_1+ \widehat{\sigma}_1\Biggr)\Biggl(\prod_{\scriptstyle \ell:\sigma_\ell=1\atop \scriptstyle 2\leq \ell\leq m}\overline{\s}^{\lambda,\phi}_\vep(v_\ell)\Biggr)\lv^{\lVert\wh{\spin}\rVert}\notag\\
&\quad\quad  \times\left(\frac{\1_{m\geq 2}}{t'-t+2u_2+u_1}+\1_{m=1}\right)
\left(\prod_{\ell=2}^{m-1}\frac{1}{2u_{\ell+1}+u_\ell}\right)\notag\\
&\quad \leq \frac{C(\lambda,\lVert\varphi\rVert_\infty,q,T,p)}{(t'-t)^p}
 \left(\int_0^T \d u_1\int_{\R^2}\d \ol{z}P_{2u_1}(\vep x_0^0-\vep\ol{z})\phi(\ol{z})\right) \notag\\
 &\quad \quad \times C(T,p)^m\mathfrak S^{\lambda,\phi}_\vep(q)^{\lVert\spin\rVert-\sigma_1} \lv^{2+\lVert\wh{\spin}\rVert},\notag
\end{align}
where the last inequality uses \eqref{mom:bdd1-1:1-100000} and adjusts the derivation of \eqref{mom:bdd1-1:1-100} by 
using \eqref{CSZ:bdd} with the choice of
\[
k\equiv m-1,\;s_0=t'-t,\;s_1\equiv u_2,\;s_2\equiv u_3,\cdots,\;s_{m-1}\equiv u_m.
\]
The last inequality proves \eqref{mom:finalbdd-1-2}. 
\end{proof}

\begin{proof}[Proof of Proposition~\ref{prop:mombdd2} (2$\cc$)]
We first prove \eqref{mom:finalbdd-2-1}. By \eqref{mom:bdd2} and  \eqref{diff:0}, for $0\leq T_1< T_2\leq T$ and $q\geq 0$, 
\begin{align}
&\lv^2\int_{T_1}^{T_2}\d t\int_{t}^{T_2}\d t'\int_0^{t'-t}\d \wt{s}_1\int_{[\wt{s}_1,t'-t)}\d \wt{\tau_1}\int_{\Delta_m(t)}\d \bs u_m\otimes \d \bs v_m\Bigg(\prod_{\ell:\wt{\sigma}_\ell=0}\delta_{\wt{s}_\ell}(\wt{\tau}_\ell)\Bigg)\Bigg(\prod_{\scriptstyle \ell:\sigma_\ell=0 \atop \scriptstyle 1\leq \ell\leq m}\delta_0(v_\ell)\Bigg)\notag\\
&\times\int \d \widetilde{\mathcal G}^{\widetilde{\bs \sigma}}_{\vep;0}\int \d \widetilde{\mathcal G}^{\widetilde{\bs \sigma}}_{\vep;1} \int \d \mathcal G^{\bs \sigma}_{\vep;0} \int \d \mathcal G^{\bs \sigma}_{\vep;1}\cdots \int \d \mathcal G^{\bs \sigma}_{\vep;m}\1\notag\\
&\quad \leq  \lv^2\int_{T_1}^{T_2}\d t\int_{t}^{T_2}\d t'\int_0^{t'-t}\d \wt{s}_1\int_{[\wt{s}_1,t'-t)}\d \wt{\tau_1}\int_{\Delta_m(t)}\d \bs u_m\otimes \d \bs v_m\notag\\
&\quad \quad \times \Bigg(\prod_{\ell:\wt{\sigma}_\ell=0}\delta_{\wt{s}_\ell}(\wt{\tau}_\ell)\Bigg)\Bigg(\prod_{\scriptstyle \ell:\sigma_\ell=0 \atop \scriptstyle 1\leq \ell\leq m}\delta_0(v_\ell)\Bigg)\notag\\
&\quad\quad  \times \int_{\R^2}\d \wt{\overline{z}}
 P_{2\wt{s}_1}(\vep \wt{x}_0^0-\vep \wt{\ol{z}})\phi( \wt{\overline{z}})\left(\int_{\R^2} 
 \d \wt{z}_1
\s^{\lambda,\phi}_\vep (\wt{s}_1,\wt{\tau}_1; \wt{\overline{z}}/\two,\wt{z}_1)\varphi(\wt{z}_1)\wt{\sigma}_1 + \widehat{\wt{\sigma}}_1\right)  \notag\\
  &\quad\quad  \times  \Biggl(\prod_{\scriptstyle \ell:\sigma_\ell=1\atop\scriptstyle 1\leq \ell\leq m}\overline{\s}^{\lambda,\phi}_\vep(v_\ell)\Biggr)\lv^{\lVert\wh{\wt{\spin}}\rVert+\lVert\wh{\spin}\rVert}
\left(\frac{1}{t'-t-\widetilde{\tau}_1+u_1}\right)\left(\prod_{\ell=1}^{m-1}\frac{1}{2u_{\ell+1}+u_\ell}\right)\notag\\
&\quad \leq  C(\lambda,\lVert\varphi\rVert_\infty,q,T,p)\int_{T_1}^{T_2}\d t\int_{t}^{T_2}\d t'\int_0^{t'-t}\d \wt{s}_1\int_{[\wt{s}_1,t'-t)}\d \wt{\tau_1}\Bigg(\prod_{\ell:\wt{\sigma}_\ell=0}\delta_{\wt{s}_\ell}(\wt{\tau}_\ell)\Bigg)\notag\\
&\quad\quad  \times \int_{\R^2}\d \wt{\overline{z}}
 P_{2\wt{s}_1}(\vep \wt{x}_0^0-\vep \wt{\ol{z}})\phi( \wt{\overline{z}})\left(\int_{\R^2} 
 \d \wt{z}_1
\s^{\lambda,\phi}_\vep (\wt{s}_1,\wt{\tau}_1; \wt{\overline{z}}/\two,\wt{z}_1)\varphi(\wt{z}_1)\wt{\sigma}_1 + \widehat{\wt{\sigma}}_1\right)\notag\\
&\quad \quad \times\left(\frac{1}{(t'-t-\widetilde{\tau}_1)^p}\right)  C(T,p)^m\mathfrak S^{\lambda,\phi}_\vep(q)^{\lVert\spin\rVert} (\wt{\sigma}_1+\lv\wh{\wt{\sigma}}_1)\lv^{2+\lVert\wh{\spin}\rVert}\notag\\
&\quad \leq C(\lambda,\lVert\varphi\rVert_\infty,q,T,p)
(T_2-T_1)^{2-p}\left(\int_0^T\d \wt{s}_1\int_{\R^2}\d \wt{\ol{z}}P_{2\wt{s}_1}(\vep \wt{x}_0^0-\vep \wt{\ol{z}})\phi(\wt{\ol{z}})
\right)\notag\\
&\quad\quad  \times C(T,p)^m\mathfrak S^{\lambda,\phi}_\vep(q)^{\lVert\spin\rVert} (\wt{\sigma}_1+\lv\wh{\wt{\sigma}}_1)\lv^{2+\lVert\wh{\spin}\rVert}.\notag
\end{align} 
where the second and third inequalities follow by modifying the argument after \eqref{mom:bdd1-1:1} until \eqref{T1T2:int}. The last inequality proves \eqref{mom:finalbdd-2-1}.

Finally, we show \eqref{mom:finalbdd-2-2}. By \eqref{mom:bdd2-1},  for all $0<t<t'\leq T$, 
\begin{align}
&\lv^2\E^{(\widetilde{B}^1,\widetilde{B}^2)}_{(\widetilde{x}_0^1,\widetilde{x}_0^2)}\Bigg[
\int_{\Delta_m(t)}\d \bs u_m\otimes \d \bs v_m\Bigg(\prod_{\scriptstyle \ell:\sigma_\ell=0 \atop \scriptstyle 1\leq \ell\leq m}\delta_0(v_\ell)\Bigg) \int \d \mathcal G^{\bs \sigma}_{\vep;0}\cdots \int \d \mathcal G^{\bs \sigma}_{\vep;m}\1\Big|_{(x_0^3,x_0^4)=(\widetilde{B}^1_{t'-t},\widetilde{B}^2_{t'-t}) }\Bigg]\notag\\
&\quad \leq\lv^2 \int_{\Delta_m(t)}\d \bs u_m\otimes \d \bs v_m\Bigg(\prod_{\scriptstyle \ell:\sigma_\ell=0 \atop \scriptstyle 1\leq \ell\leq m}\delta_0(v_\ell)\Bigg) \Biggl(\prod_{\scriptstyle \ell:\sigma_\ell=1\atop\scriptstyle 1\leq \ell\leq m}\overline{\s}^{\lambda,\phi}_\vep(v_\ell)\Biggr)\lv^{\lVert\wh{\spin}\rVert}\notag\\
&\quad \quad \times
\left(\frac{1}{t'-t+u_1}\right)\left(\prod_{\ell=1}^{m-1}\frac{1}{2u_{\ell+1}+u_\ell}\right)\notag\\
&\quad \leq  \frac{C(\lambda,\lVert\varphi\rVert_\infty,q,T,p)}{(t'-t)^p}C(T,p)^m
\mathfrak S^{\lambda,\phi}_\vep(q)^{\lVert\spin\rVert}\lv^{2+\lVert\wh{\spin}\rVert},\notag
\end{align}
where is obtained by adjusting the proof of \eqref{mom:finalbdd-1-2}.
We have proved \eqref{mom:finalbdd-2-2}.
\end{proof}

We now prove Theorem~\ref{thm:mombdd}. Since $X_0(\cdot)$ is assumed bounded and nonnegative, \eqref{eq:FKmom} reduces the argument to the case $X_0 \equiv \1$. Thus, for parts (1$\cc$) and (2$\cc$) of Theorem~\ref{thm:mombdd}, it suffices to analyze \eqref{eq:FKmom} accordingly and \eqref{mom:series} with $f \equiv \1$.\medskip 

\begin{proof}[Proof of Theorem~\ref{thm:mombdd} (1$\cc$)]
We first bound the first term on the right-hand side of \eqref{mom:series}. By 
\eqref{graphical22}, \eqref{2mom:eq1}, and \eqref{2mom:eq2}, 
\begin{align}
&\lv^2\E^{(\widetilde{B}^1,\widetilde{B}^2)}_{(\widetilde{x}_0^1,\wt{x}_0^1+\vep\wt{x}_0^0)}\Bigg[\exp\Bigg\{\lv \int_0^{t'-t} \d r \varphi_\vep(\widetilde{B}^{\bi_0}_r)\Bigg\}\widetilde{\1}_t(\wt{B}_{t'-t})\Bigg]\notag\\
&\quad \leq \lv^2+\wh{\wt{\sigma}}_1\lv^3\int_0^{t'-t} \d \wt{s}_1 \int_{\R^4} \begin{bmatrix}
\d \wt{x}_1^{2}\\
 \d \wt{x}_1^1
 \end{bmatrix}_{\otimes}
 \varphi(\wt{x}_1^{\bi_0})
\begin{bmatrix}
 P_{0,\wt{s}_1}(\wt{x}_0^1+\vep\wt{x}_0^0,\wt{x}_1^{1}+\vep\two \wt{x}^{\bi_0}_1)\\
  P_{0,\wt{s}_1}(\wt{x}_0^{1},\wt{x}_1^{1})
\end{bmatrix}_\times\notag\\
&\quad\quad  +\wt{\sigma}_1\lv^2\int_0^{t'-t}\d \wt{s}_1\int_{\wt{s}_1}^{t'-t}\d \wt{\tau}_1 \int_{\R^4} \begin{bmatrix}
\d \wt{x}_1^{2}\\
 \d \wt{x}_1^1
 \end{bmatrix}_{\otimes}\varphi(\wt{x}_1^{\bi_0})\int_{\R^2} \d \wt{z}_1 \varphi(\wt{z}_1)\notag\\
&\quad\quad  \times\begin{bmatrix}
 P_{0,\wt{s}_1}(\wt{x}_0^1+\vep\wt{x}_0^0,\wt{x}_1^{1}+\vep\two \wt{x}^{\bi_0}_1)\\
  P_{0,\wt{s}_1}(\wt{x}_0^{1},\wt{x}_1^{1})
\end{bmatrix}_\times
\s^{\lambda,\phi}_\vep (\wt{s}_1,\wt{\tau}_1; \wt{x}_1^{\bi_{0}},\wt{z}_1)
\notag\\
&\quad \leq \lv^2+\wh{\wt{\sigma}}_1\lv^3\int_0^{t'-t} \d \wt{s}_1 \int_{\R^2}\d \wt{\ol{z}}\phi(\wt{\ol{z}})P_{2\wt{s}_1}(\vep \wt{x}_0^0-\vep \wt{\ol{z}})\notag\\
&\quad \quad +\wt{\sigma}_1\lv^2C(\lambda,\lVert\varphi\rVert_\infty,T)\int_0^{t'-t}\d \wt{s}_1 \int_{\R^4} \begin{bmatrix}
\d \wt{x}_1^{2}\\
 \d \wt{x}_1^1
 \end{bmatrix}_{\otimes}\varphi(\wt{x}_1^{\bi_0})\begin{bmatrix}
 P_{0,\wt{s}_1}(\wt{x}_0^1+\vep\wt{x}_0^0,\wt{x}_1^{1}+\vep\two \wt{x}^{\bi_0}_1)\\
  P_{0,\wt{s}_1}(\wt{x}_0^{1},\wt{x}_1^{1})
\end{bmatrix}_\times\notag\\
&\quad \leq \lv^2+\lv^2C(\lambda,\lVert\varphi\rVert_\infty,T)\int_0^{t'-t} \d \wt{s}_1 \int_{\R^2}\d \wt{\ol{z}}\phi(\wt{\ol{z}})P_{2\wt{s}_1}(\vep\wt{x}_0^0-\vep \wt{\ol{z}}),\label{heatgrowth}
\end{align}
where the second equality uses \eqref{ineq:apriori}  and \eqref{ineq:gausstime}, and the last inequality uses \eqref{ineq:gausstime} again. Hence, for all $\vep\in (0,\ovl]$ and $0\leq T_1\leq T_2\leq T$,  
\begin{align}
&\int_{T_1}^{T_2}\d t\int_{t}^{T_2}\d t'  \lv^2\E^{(\widetilde{B}^1,\widetilde{B}^2)}_{(\widetilde{x}_0^1,\wt{x}_0^1+\vep\wt{x}_0^0)}\Bigg[\exp\Bigg\{\lv \int_0^{t'-t} \d r \varphi_\vep(\widetilde{B}^{\bi_0}_r)\Bigg\}\widetilde{\1}_t(B_{t'-t})\Bigg]\notag\\
&\quad \leq C(\lambda,\phi,T)\left(\lv^2\int_0^{T} \d \wt{s}_1 \int_{\R^2}\d \wt{\ol{z}}\phi(\wt{\ol{z}})P_{2\wt{s}_1}(\vep\wt{x}_0^0-\vep \wt{\ol{z}})+\lv^2\right)(T_2-T_1)^2.\label{mod:final1aux}
\end{align}

Next, we analyze the $\int_{T_1}^{T_2}\d t\int_{t}^{T_2}\d t' $-integral of the remaining terms on the right-hand side of \eqref{mom:series}. To use the bounds obtained in Proposition~\ref{prop:mombdd2}, we distinguish $\bi_1$ in three scenarios: (i) $\bi_1=(2,1)$; (ii) $\bi_2=(4,3)$; and (iii) $\bi_1=(i,j)$ for $i\in \{4,3\}$ and $j\in \{2,1\}$. The first two scenarios allow direct application of the bounds in Proposition~\ref{prop:mombdd2}. By relabelling the initial conditions of the four particles, all cases in the third scenario can use the bounds established for $\bi_1=(3,2)$ because $f\equiv \1$ is used. Specifically, for $\bi=(3,1)$, the bounds in \eqref{mom:finalbdd-2-1} and \eqref{mom:finalbdd-2-2} are still valid because mapping $(3,1)$ to $(3,2)$ does not lead to swapping $\wt{x}_0^1$ and $\wt{x}_0^2$, and the right-hand sides of those bounds depend on $(x_0^1,x_0^2,\wt{x}_0^1,\wt{x}_0^2)$ only through $(\wt{x}_0^1,\wt{x}_0^2)$. For $\bi\in \{(4,1),(4,2)\}$, a similar argument shows that the bound in \eqref{mom:finalbdd-2-2} still applies, but since the roles of $\wt{x}_0^1$ and $\wt{x}_0^2$ are exchanged, the bound in \eqref{mom:finalbdd-2-1} must be modified so that $\vep\wt{x}_0^0=\wt{x}_0^2-\wt{x}_0^1$ in \eqref{mom:finalbdd-2-1} is replaced by $-\vep\wt{x}_0^0=\wt{x}_0^1-\wt{x}_0^2$. The corresponding bound is given below:
\begin{align}
&C(\lambda,\lVert\varphi\rVert_\infty,q,T,p)
(T_2-T_1)^{2-p}\Bigg(\int_0^T\d \wt{s}_1\int_{\R^2}\d \wt{\ol{z}}P_{2\wt{s}_1}(-\vep \wt{x}_0^0-\vep \wt{\ol{z}})\phi(\wt{\ol{z}})
\Bigg)\notag\\
&\quad \times C(T,p)^m\mathfrak S^{\lambda,\phi}_\vep(q)^{\lVert\spin\rVert} (\wt{\sigma}_1+\lv\wh{\wt{\sigma}}_1)\lv^{2+\lVert\wh{\spin}\rVert}\notag\\
&\quad =C(\lambda,\lVert\varphi\rVert_\infty,q,T,p)
(T_2-T_1)^{2-p}\Bigg(\int_0^T\d \wt{s}_1\int_{\R^2}\d \wt{\ol{z}}P_{2\wt{s}_1}(\vep \wt{x}_0^0-\vep \wt{\ol{z}})\phi(\wt{\ol{z}})
\Bigg)\notag\\
&\quad\quad  \times C(T,p)^m\mathfrak S^{\lambda,\phi}_\vep(q)^{\lVert\spin\rVert} (\wt{\sigma}_1+\lv\wh{\wt{\sigma}}_1)\lv^{2+\lVert\wh{\spin}\rVert},\label{new:mod:final1aux}
\end{align}
where we first change $\phi(\wt{\ol{z}})$ to $\phi(-\wt{\ol{z}})$ by changing variables and then use the even parity of $\phi$ from Proposition~\ref{prop:XM} (1$\cc$). 

We proceed to an explicit bound of the $\int_{T_1}^{T_2}\d t\int_{t}^{T_2}\d t' $-integral of the remaining terms on the right-hand side of \eqref{mom:series} now. First, 
for any $q\geq 0$ and $0\leq T_1\leq T_2\leq T$, we obtain from \eqref{mom:finalbdd-1-2}, \eqref{mom:finalbdd-2-2}, the aforementioend extension of \eqref{mom:finalbdd-2-2}, \eqref{T1T2:int}, and the rule
$(a+b)^m=\sum_{\spin \{0,1\}^m}a^{\lVert\spin\rVert}b^{\lVert\wh{\spin}\rVert}$ for $a,b\in \R$ that
\begin{align*}
&\sum_{m=1}^\infty \sum_{\stackrel{\scriptstyle \bi_1,\cdots,\bi_m\in \mathcal E_4}{\bi_1\neq \cdots \neq \bi_m}}\sum_{\spin \in \{0,1\}^m} 
\int_{T_1}^{T_2}\d t\int_{t}^{T_2}\d t'\lv^2\E^{(\widetilde{B}^1,\widetilde{B}^2)}_{(\widetilde{x}_0^1,\widetilde{x}_0^1+\vep \wt{x}_0^0)}\Bigg[\int_{\Delta_m(t)}\d \bs u_m\otimes \d \bs v_m\Biggl( \prod_{\scriptstyle \ell:\sigma_\ell=0 \atop \scriptstyle 1\leq \ell\leq m}\delta_{0}(v_\ell)\Biggr)\\
&\quad \times \int \d \mathcal G^\spin_{\vep;0}  \int \d \mathcal G^{\bs \sigma}_{\vep;1}
 \cdots\int \d \mathcal G^\spin_{\vep;m} \1\mid_{(x_0^3,x_0^4)=(\widetilde{B}^1_{t'-t},\widetilde{B}^2_{t'-t}) }\Bigg]\\
 &\quad \leq \sum_{m=1}^\infty \sum_{\stackrel{\scriptstyle \bi_1,\cdots,\bi_m\in \mathcal E_4}{\bi_1\neq \cdots \neq \bi_m}}
 C(\lambda,\lVert\varphi\rVert_\infty,q,T,p)
(T_2-T_1)^{2-p}\\
&\quad\quad  \times \left(\lv^2\int_0^T \d u_1\int_{\R^2}\d \ol{z}P_{2u_1}(\vep x_0^0-\vep\ol{z})\phi(\ol{z})+\lv^2\right)C(T,p)^m\left(\mathfrak S^{\lambda,\phi}_\vep(q) +\lv\right)^m \left(\mathfrak S^{\lambda,\phi}_\vep(q)^{-1}\vee 1\right)\\
&\quad \leq \sum_{m=1}^\infty C(\lambda,\lVert\varphi\rVert_\infty,q,T,p)\left(\lv^2\int_0^T \d u_1\int_{\R^2}\d \ol{z}P_{2u_1}(\vep x_0^0-\vep\ol{z})\phi(\ol{z})+\lv^2\right)\\
&\quad\quad  \times {4\choose 2}^m C(T,p)^m \left(\mathfrak S^{\lambda,\phi}_\vep(q) +\lv\right)^m \left(\mathfrak S^{\lambda,\phi}_\vep(q)^{-1}\vee 1\right)
(T_2-T_1)^{2-p}.
 \end{align*}
To make the last series converge, we use \eqref{def:Sconv} to choose $q_{\ref{q:modulus}}=q_{\ref{q:modulus}}(\lambda,\varphi)$ and $\vep(\lambda,\phi,q_{\ref{q:modulus}})\in (0,\ovl)$ such that for all  $\vep\in (0,\vep(\lambda,\phi,q_{\ref{q:modulus}}))$, 
\begin{align}
{4\choose 2}\cdot C(T,p)\cdot (\mathfrak S^{\lambda,\phi}_\vep(q_{\ref{q:modulus}}) +\lv)
<\frac{1}{2}\quad\mbox{and}\quad  \mathfrak S^{\lambda,\phi}_\vep(q_{\ref{q:modulus}})\geq \frac{1}{10}\cdot \frac{4\pi}{\log (q_{\ref{q:modulus}} /\beta)}.\label{q:modulus}
\end{align}
Then, we get, for all $\vep\in (0,\vep(\lambda,\phi,q_{\ref{q:modulus}}))$ and all $0\leq T_1\leq T_2\leq T$, 
\begin{align}
&\sum_{m=1}^\infty \sum_{\stackrel{\scriptstyle \bi_1,\cdots,\bi_m\in \mathcal E_4}{\bi_1\neq \cdots \neq \bi_m}}\sum_{\spin \in \{0,1\}^m} 
\int_{T_1}^{T_2}\d t\int_{t}^{T_2}\d t' \lv^2\E^{(\widetilde{B}^1,\widetilde{B}^2)}_{(\widetilde{x}_0^1,\widetilde{x}_0^1+\vep \wt{x}_0^0)}\Bigg[\int_{\Delta_m(t)}\d \bs u_m\otimes \d \bs v_m\notag\\
&\quad \times \Biggl( \prod_{\scriptstyle \ell:\sigma_\ell=0 \atop \scriptstyle 1\leq \ell\leq m}\delta_{0}(v_\ell)\Biggr)\int \d \mathcal G^\spin_{\vep;0}  \int \d \mathcal G^{\bs \sigma}_{\vep;1}
 \cdots\int \d \mathcal G^\spin_{\vep;m} \1\Big|_{(x_0^3,x_0^4)=(\widetilde{B}^1_{t'-t},\widetilde{B}^2_{t'-t}) }\Bigg]\notag\\
 &\quad \leq C(\lambda,\phi,T,p)\left(\lv^2\int_0^T \d u_1\int_{\R^2}\d \ol{z}P_{2u_1}(\vep x_0^0-\vep\ol{z})\phi(\ol{z})+\lv^2\right)  (T_2-T_1)^{2-p}. \label{mod:final2}
\end{align}
Moreover, a similar argument shows the following bound for the $\int_{T_1}^{T_2}\d t\int_{t}^{T_2}\d t' $-integral of the third term on the right-hand side of \eqref{mom:series}: by \eqref{mom:finalbdd-1-1}, \eqref{mom:finalbdd-2-1}, and the aforementioned extension of \eqref{mom:finalbdd-2-1} (using the bound in \eqref{new:mod:final1aux} when needed): for all $\vep\in (0,\vep(\lambda,\phi,q_{\ref{q:modulus}}))$, 
\begin{align}
&\sum_{m=1}^\infty \sum_{\stackrel{\scriptstyle \bi_1,\cdots,\bi_m\in \mathcal E_4}{\bi_1\neq \cdots \neq \bi_m}}\sum_{\scriptstyle \widetilde{\spin}\in \{0,1\}\atop\scriptstyle \spin \in \{0,1\}^m}\lv^2\int_{T_1}^{T_2}\d t\int_{t}^{T_2} \d t'
\int_0^{t'-t} \d \wt{s}_1\int_{[\wt{s}_1,t'-t)}\d \wt{\tau}_1  \int_{\Delta_m(t)}\d \bs u_m\otimes \d \bs v_m\notag\\
 & \times\Biggl(\prod_{\scriptstyle \ell:\widetilde{\sigma}_\ell=0\atop \scriptstyle \ell=1}\delta_{\wt{s}_\ell}(\widetilde{\tau}_\ell)\Biggr) \Biggl( \prod_{\scriptstyle \ell:\sigma_\ell=0 \atop \scriptstyle 1\leq \ell\leq m}\delta_{0}(v_\ell)\Biggr)\int \d \widetilde{\mathcal G}^{\widetilde{\spin}}_{\vep;0}\int \d \widetilde{\mathcal G}^{\widetilde{\spin}}_{\vep;1}\int \d \mathcal G^\spin_{\vep;0} 
 \cdots\int \d \mathcal G^\spin_{\vep;m} \1\notag\\
 &\quad \leq C(\lambda,\phi,T,p) \biggl[\biggl(\lv\int_0^T \d \wt{s}_1\int_{\R^2}\d \wt{\ol{z}}P_{2\wt{s}_1}(\vep \wt{x}_0^0-\vep\wt{\ol{z}})\phi(\wt{\ol{z}})\biggr)\notag\\
 &\quad \quad \times \biggl(\lv\int_0^T \d u_1\int_{\R^2}\d \ol{z}P_{2u_1}(\vep x_0^0-\vep\ol{z})\phi(\ol{z})\biggr)+\lv^2\biggr](T_2-T_1)^{2-p}. \label{mod:final3}
 \end{align}

In summary, by \eqref{mom:series}, we can apply \eqref{mod:final1aux}, \eqref{mod:final2} and \eqref{mod:final3} to get \eqref{modulus:bdd4} for all $\vep\in (0,\vep(\lambda,\phi,q_{\ref{q:modulus}}))$. Moreover, since the bound in \eqref{modulus:bdd4} holds for $\vep\in [\vep(\lambda,\phi,q_{\ref{q:modulus}}),\ovl]$ by \eqref{eq:momdual2}, it follows that \eqref{modulus:bdd4} holds for all $\vep\in (0,\ovl]$. This completes the proof of (1$\cc$). 
\end{proof}

\begin{proof}[Proof of Theorem~\ref{thm:mombdd} (2$\cc$)]
By \eqref{eq:momdual2}, \eqref{P:vepseries} and \eqref{int:P:graphical}, it holds that, for all $t,a,a'\in (0,T]$, 
\begin{align*}
&\lv\E_{(\wt{x}_0^1,\wt{x}_0^1+\vep \wt{x}_0^0)}^{(\wt{B}^1,\wt{B}^2)}\E\left[X_\vep(x_0^1,t)X_\vep(x_0^1+\vep x_0^0,t)X_\vep (B^1_{a},t)X_\vep(B^2_{a'},t)\right]\\
&\quad \leq\lv\lVert X_0\rVert_\infty \\
&\quad\quad + \lVert X_0\rVert_\infty \sum_{m=1}^\infty \sum_{\stackrel{\scriptstyle \bi_1,\cdots,\bi_m\in \mathcal E_4}{\bi_1\neq \cdots \neq \bi_m}}\sum_{\spin \in \{0,1\}^m} \lv\E^{(\widetilde{B}^1,\widetilde{B}^2)}_{(\widetilde{x}_0^1,\wt{x}_0^1+\vep \wt{x}_0^0)}\Bigg[
\int_{\Delta_m(t)}\d \bs u_m\otimes \d \bs v_m\Bigg(\prod_{\scriptstyle \ell:\sigma_\ell=0 \atop \scriptstyle 1\leq \ell\leq m}\delta_0(v_\ell)\Bigg) \\
&\quad \quad \times \int \d \mathcal G^{\bs \sigma}_{\vep;0} \int \d \mathcal G^{\bs \sigma}_{\vep;1}\cdots \int \d \mathcal G^{\bs \sigma}_{\vep;m}\1\Big|_{(x_0^3,x_0^4)=(\widetilde{B}^1_{a},\widetilde{B}^2_{a'}) }\Bigg],
\end{align*}
where the graphical integrals being taken as expectations are subject to the initial condition of $(x_0^1,x_0^2,x_0^3,x_0^4)$ with $x_0^2=x_0^1+\vep x_0^0$. 

When $a'=a=a_0$, we apply \eqref{mom:finalbdd-1-2}, \eqref{mom:finalbdd-2-2}, and the extension of \eqref{mom:finalbdd-2-2} mentioned in the proof of (1$\cc$), dividing each by $\lv$, and selecting $t'-t=a_0$. With the same choice of $q_{\ref{q:modulus}}=q_{\ref{q:modulus}}(\lambda,\varphi)$ and $\vep(\lambda,\phi,q_{\ref{q:modulus}})\in (0,1)$ as in \eqref{q:modulus}, we obtain \eqref{modulus:bdd2} for all $\vep\in (0,\vep(\lambda,\phi,q_{\ref{q:modulus}}))$. By \eqref{eq:momdual2}, this establishes \eqref{modulus:bdd2} for all $\vep\in (0,\ovl]$ when $a'=a=a_0$. Moreover, the right-hand side of \eqref{modulus:bdd2} in this case is independent of $(\wt{x}_0^1, \wt{x}_0^1 + \vep \wt{x}_0^0)$. Therefore, for $a < a'$, we can apply the Markov property of $\wt{B}^2$ at time $a' - a$ and use the result established for $a=a_0$ to obtain \eqref{modulus:bdd2}. This argument extends \eqref{modulus:bdd2} to all $a', a \in (0, T]$.
\end{proof}

\begin{proof}[Proofs of Theorem~\ref{thm:mombdd} (3$\cc$) and (\ref{thmmombdd:aux})]
For $N=2$, the infinite series in \eqref{P:vepseries} reduces to just one term, corresponding to $m=1$; see Lemma~\ref{lem:mix2mom}.  Therefore, the proof of Theorem~\ref{thm:mombdd}~(3$\cc$) follows immediately by applying a single step from the argument for Theorem~\ref{thm:mombdd}~(1$\cc$); see \eqref{heatgrowth}. Additionally, \eqref{thmmombdd:aux} follows from the simple calculation below:
\begin{align}
\lv\int_0^{T} \d s_1 \int_{\R^2}\d\ol{z}\phi(\ol{z})P_{2s_1}(\vep x_0^0-\vep \ol{z})&= \lv\int_{\R^2}\d \ol{z}\phi(\ol{z})
\int_0^{\vep^{-2}\lvert x_0^0-\ol{z}\rvert^{-2}T}\frac{\d u}{4\pi u}\exp\left(-\frac{1}{4u}\right)
\notag\\
&\leq C(\lambda,T)\left(\int_{\R^2}\d \ol{z}\phi(\ol{z})\lvert\log \lvert x_0^0-\ol{z}\rvert \rvert+1\right),
\end{align}
where the last inequality uses the definition \eqref{def:lv} of $\lv$. 
\end{proof}

\subsection{Proof of Theorem~\ref{thm:tight}: joint tightness}\label{sec:jointC}
Recall the definitions of $X^\sharp_\vep$, $\nu^\sharp_\vep$, $\mu^\sharp_\vep$, and $\mathring{\mu}^\sharp_\vep$ in \eqref{def:sharp}, and observe the following properties:
\begin{itemize}
\item [\bf (P1)] By \cite[2.4 Proposition, p.107]{EK:MP}, the family of laws of $\{(X^\sharp_\vep,\nu^\sharp_\vep,\mu^\sharp_\vep,\mathring{\mu}^\sharp_\vep)\}_{\vep\in (0,\ovl]}$ is tight as probability measures on the product Polish space $C_{\mathcal M_f(\R^2)}[0,\infty)\times  C_{\mathcal M_f(\R^2)}[0,\infty)\times C_{\mathcal M_f(\R^4)}[0,\infty)\times C_{\mathcal M_f(\R^4)}[0,\infty)$ if and only if the marginals are tight. Here, $\mathcal{M}_f(\R^d)$ is Polish for each $d \in \mathbb{N}$, as established by normalizing measures and applying \cite[1.7 Theorem, p.101]{EK:MP}, and \emph{marginal tightness} means that the family of laws from each set below is tight as probability measures on $C_{\mathcal M_f(\R^{r_j})}[0,\infty)$, with $r_1=r_2=2$ for the first two sets and $r_3=r_4=4$ for the last two sets: 
\begin{align}\label{thm:tight-1}
\{X^\sharp_\vep\}_{\vep\in (0,\ovl]}, \; \{\nu^\sharp_\vep\}_{\vep\in (0,\ovl]}, \; \{\mu^\sharp_\vep\}_{\vep\in (0,\ovl]}, \; \{\mathring{\mu}_\vep^\sharp\}_{\vep\in (0,\ovl]}.
\end{align} 

\item [\bf (P2)] The set $\C^2_b(\R^d)$, consisting of functions on $\R^d$ with bounded continuous partial derivatives up to second order, forms a separating class in $\C_b(\R^d)$ and includes the constant function $1$. A separating class is a collection of functions such that if two finite measures agree on all functions in the class, the measures are identical (see the definition preceding \cite[Theorem~II.4.1]{Perkins:DW}).
\end{itemize}
By (P1), it suffices to establish marginal tightness. Let $D_{\mathcal{M}_f(\R^{r_j})}[0, \infty)$ denote the space of $\mathcal{M}_f(\R^{r_j})$-valued c\`adl\`ag functions on $[0,\infty)$ with Skorokhod's $J_1$-topology.  Because all processes in \eqref{thm:tight-1} are continuous, \cite[Problem 25, p.153]{EK:MP} further reduces this to verifying that, for each set, the corresponding family of laws is tight as probability measures on $D_{\mathcal{M}_f(\R^{r_j})}[0, \infty)$ and that the family of laws is $C$-tight: every subsequential limit law corresponds to an almost surely continuous process. To establish this, we apply \cite[Theorem~II.4.1]{Perkins:DW} with $(E, D_0)$ taken as $(\R^2, \C^2_b(\R^2))$ or $(\R^4, \C^2_b(\R^4))$, as justified by (P2). It remains to verify conditions (C1)--(C4) below:

\begin{itemize}
\item [\bf (C1)]For any $\eta>0$ and $0<T<\infty$, there exist $0\leq f_{0}\in \C^2_b(\R^2)$ and $K\in(0,\infty)$ with $f_{0}(x)=1$ for all $\lvert x\rvert >K$ such that 
\[
\sup_{\vep\in (0,\ovl]}\P\Biggl(\sup_{0\leq t\leq T}\rho_\vep(f_{0},t)>\eta\Biggr)<\eta,\quad \forall\;\rho\in \{X^\sharp,\nu^\sharp\}.
\]   
\item [\bf (C2)] For any $f\in \C^2_b(\R^2)$ and $\rho\in \{X^\sharp,\nu^\sharp\}$, the family of laws of $\{\rho_\vep(f,t)\}_{\vep\in (0,\ovl]}$ is $C$-tight as probability measures on $D_\R[0,\infty)$.
\item [\bf (C3)]For any $\eta>0$ and $0<T<\infty$, there exist $0\leq g_{0}\in \C^2_b(\R^4)$ and $L\in (0,\infty)$ with $g_{0}(x,y)=1$ for all $\lvert (x,y)\rvert>L$, $(x,y)\in \R^4$, such that 
\[
\sup_{\vep\in (0,\ovl]}\P\Biggl(\sup_{0\leq t\leq T}\rho_\vep(g_{0},t)>\eta\Biggr)<\eta,\quad \forall\;\rho\in \{\mu^\sharp,\mathring{\mu}^\sharp\}.
\]   
\item [\bf (C4)] For any $g\in \C^2_b(\R^4)$ and $\rho\in \{\mu^\sharp,\mathring{\mu}^\sharp\}$, the family of laws of $\{\rho_\vep(g,t)\}_{\vep\in (0,\ovl]}$ is $C$-tight as probability measures on $D_\R[0,\infty)$.
\end{itemize}

The verifications of (C1) and (C2) are provided in Steps~1 and 2 below. Because the definitions of $\nu_\vep$, $\mu_\vep$, and $\mathring{\mu}_\vep$ differ only slightly, we address (C3) and (C4) in Step~3 by highlighting the key inequalities; the remaining arguments closely follow those for (C1) and (C2) with $\rho \equiv \nu^\sharp$. \medskip 

\noindent {\bf Step 1.} To verify (C1), we consider the case of  $\rho\equiv X^\sharp$ first. The key bound for this purpose is given by the second inequality below, while the first inequality is deduced from the weak formulation in \eqref{def:weak}: for any $0\leq f\in \C_b^2(\R^2)$,
\begin{align*}
&\E\biggl[\sup_{0\leq t\leq T}X^\sharp_\vep(f,t)\biggr]\\
&\quad \leq X_\vep (f\psi,0)+\E\left[\int_{t=0}^TX_\vep\left(\left|\frac{\Delta (f\psi)}{2}\right|,t\right)\d t\right]+\E\biggl[\sup_{0\leq t\leq T}\lvert M_\vep(f\psi,t)\rvert \biggr]\\
&\quad \less \int_x(f\psi)(x)X_0(x)\d x+\E\left[\int_{t=0}^TX_\vep\left(\left|\frac{\Delta (f\psi)}{2}\right|,t\right)\d t\right]\\
&\quad\quad +\E\biggl[\int_{t=0}^T\int_{x_0^1} \int_{x_0^0}(f\psi)(x_0^1)(f\psi)(x_0^1+\vep x_0^0)\phi(x_0^0)\lv X_\vep(x_0^1,t)X_\vep(x_0^1+\vep x_0^0,t)\d x_0^0  \d x_0^1\d t\biggr]^{1/2}.
\end{align*}
In more detail, $\psi(x)\,\defeq\,(1+\lvert x\rvert ^{10})^{-1}$ by the definition of $X^\sharp$ in \eqref{def:sharp1}, and the last term follows by applying the Burkholder--Davis--Gundy inequality \cite[(4.1) Theorem, p.160]{RY}, \eqref{eq:covarM-2}, and Jensen's inequality. For the right-hand side of the last inequality, we can bound the first term by using the assumption that $X_0(\cdot)$ is bounded, the second term by using the mild form \eqref{def:mild:vep}, and the third term by using Theorem~\ref{thm:mombdd} (3$\cc$). Hence, by dominated convergence, for any $\eta>0$, there exists $(f_0,K)$ such that $0\leq f_0\in \C^2_b(\R^2)$ and $K\in (0,\infty)$ with $f_0(x)=1$ for all $\lvert x\rvert >K$, and 
\[
\sup_{\vep\in (0,\ovl]}\E\biggl[\sup_{0\leq t\leq T}X^\sharp_\vep(f_0,t)\biggr]<\eta^2.
\]
By Markov's inequality, the foregoing inequality is enough to verify (C1) for $\rho\equiv X^\sharp$.

The verification of (C1) for $\rho\equiv \nu^\sharp$ is similar. The main difference is that we proceed with the following bound now instead of the key bound mentioned above for $\rho\equiv X^\sharp$:
\begin{align}
\forall\;0\leq f\in \C_b^2(\R^2),&\quad \;\E\biggl[\sup_{0\leq t\leq T}\nu^\sharp_\vep(f,t)\biggr]
\leq \E\biggl[\int_{t=0}^T\int_{x}(f\psi)(x)\lv X_\vep(x,t)^2\d x \d t\biggr],\label{thm:tight-2}
\end{align} 
which follows immediately from the definition of $\nu^\sharp_\vep$ in \eqref{def:sharp1}. By Theorem~\ref{thm:mombdd}~(3$\cc$), the preceding argument verifying (C1) for $\rho \equiv X^\sharp$ also applies here, yielding the required condition in (C1) with $\rho \equiv \nu^\sharp$. This completes the verification of (C1). \medskip 

\noindent {\bf Step 2.} We verify (C2) in this step by showing that the processes under consideration satisfy
Kolmogorov's criterion for weak compactness \cite[(1.8) Theorem, pp.517--518]{RY}.

For $\rho\equiv X^\sharp$, we work with $\widetilde{X}_\vep$ as the solution to \eqref{def:SHE:vep} when the initial condition is the constant function $\1(\cdot)$. Then by the weak formulation \eqref{def:weak} and the inequality $(a+b)^4\less a^4+b^4$ for all $a,b\geq 0$, we get, for any $0\leq T_1\leq T_2\leq T$, 
\begin{align}
&\E\left[\lvert X_\vep^\sharp(f,T_2)-X^\sharp_\vep(f,T_1)\rvert^4\right]\notag\\
&\quad\less \E\biggl[\biggl|\int_{t=T_1}^{T_2}X_\vep\left(\frac{\Delta (f\psi)}{2},t\right)\d t\biggr|^4\biggr]+\E\biggl[\biggl|\int_{t=T_1}^{T_2}\int_{x}(f\psi)(x) \M_\vep(\d x,\d t)\biggr|^4\biggr]\notag\\
&\quad\leq C(\lVert X_0\rVert_\infty) \E\biggl[\wt{X}_\vep\left(\biggl|\frac{\Delta (f\psi)}{2}\biggr|,T\right)^4\biggr](T_2-T_1)^4
\notag\\
&\quad \quad+\E\biggl[\biggl(\int_{t=T_1}^{T_2}\int_{x_0^1}\int_{x_0^0} (f\psi)(x_0^1)(f\psi)(x_0^1+\vep x_0^0)\phi(x_0^0)\lv X_\vep(x_0^1,t)X_\vep(x_0^1+\vep x_0^0,t)\d x_0^0\d x_0^1  \d t\biggr)^2\biggr].\notag
\end{align}
Here, the first term on the right-hand side follows by applying to the first term on the left-hand side Jensen's inequality and then the increasing monotonicity of $t\mapsto \E[\prod_{j=1}^4\widetilde{X}_\vep(x_0^j,t)]$ for all $x_0^1,x_0^2,x_0^3,x_0^4\in \R^2$, which can be seen by using \eqref{eq:momdual2}, and the last term in the foregoing inequality follows by using the Burkholder--Davis--Gundy inequality~\cite[(4.1) Theorem, p.160]{RY} and \eqref{eq:covarM-2}. To use the last inequality, note that
the expectation in the first term admits a bound uniformly in $\vep\in (0,\ovl]$, which is a particular consequence of \cite[Theorem~2.3]{C:DBG}, and
we can bound the last term by Theorem~\ref{thm:mombdd} (1$\cc$) with $p=1/2$. Hence,
\[
\sup_{\vep\in (0,\ovl]}\E\left[\lvert X_\vep^\sharp(f,T_2)-X^\sharp_\vep(f,T_1)\rvert ^4\right]\leq C(\lambda,\phi,\lVert X_0\rVert_\infty, T,f)(T_2-T_1)^{3/2},\;\forall\;0\leq T_1\leq T_2\leq T.
\]
We have verified (C2) for $\rho\equiv X^\sharp$ since Kolmogorov's criterion for weak compactness applies. 

To verify (C2) for $\rho\equiv \nu^\sharp$, Theorem~\ref{thm:mombdd} (1$\cc$) with $p=1/2$ is applied again to get the second inequality below:
\begin{align}
&\sup_{\vep\in (0,\ovl]}\E\left[\lvert \nu_\vep^\sharp(f,T_2)-\nu^\sharp_\vep(f,T_1)\rvert^2\right]\notag\\
&\quad \leq
\sup_{\vep\in (0,\ovl]}\E\biggl[\biggl(\int_{t=T_1}^{T_2}\int_{x}(f\psi)(x)\lv X_\vep(x,t)^2 \d x\d t\biggr)^2\biggr]\label{thm:tight-3}\\
&\quad \leq  C(\lambda,\phi,\lVert X_0\rVert_\infty, T,f)(T_2-T_1)^{3/2},\;\forall\; 0\leq T_1\leq T_2\leq T.\notag
\end{align}
The required $C$-tightness now follows by using Kolmogorov's criterion for weak compactness. \medskip

\noindent {\bf Step~3.} To verify (C3) and (C4), note that the definitions of $\mu^\sharp_\vep$ and $\mathring{\mu}^\sharp_\vep$ in \eqref{def:sharp2} yield the following bounds, where $0\leq g\in \C^2_b(\R^4)$:
\begin{align*}
&\E\biggl[\sup_{0\leq t\leq T}\mu^\sharp_\vep(g,t)\biggr]\\
&\quad \;\leq \E\left[\int_{t=0}^T\int_{x}\int_{y}g(\vep y+x,x)\psi(\vep y+x)\psi(x)\lv X_\vep(\vep y+x,t)X_\vep(x,t)\phi(y)\d y\d x \d t\right],\\
&\E\biggl[\sup_{0\leq t\leq T}\mathring{\mu}^\sharp_\vep(g,t)\biggr]\\
&\quad\;\leq \E\biggl[\int_{t=0}^T\int_{x}\int_{y}g( x,y)\psi(x)\psi(y)\lv X_\vep(\vep y+x,t)X_\vep(x,t)\phi(y)\d y\d x \d t\biggr],\\
&\E\biggl[\lvert \mu_\vep^\sharp(g,T_2)-\mu^\sharp_\vep(g,T_1)\rvert^2\biggr]\\
&\quad \;\leq
\E\biggl[\left(\int_{t=T_1}^{T_2}\int_{x}\int_{y}g(\vep y+x,x)\psi(\vep y+x)\psi(x)\lv X_\vep(\vep y+x,t)X_\vep(x,t)\phi(y)\d y\d x\d t\right)^2\biggr],\\
&\E\left[\lvert \mathring{\mu}_\vep^\sharp(g,T_2)-\mathring{\mu}^\sharp_\vep(g,T_1)\rvert^2\right]\\
&\quad\;\leq
\E\biggl[\left(\int_{t=T_1}^{T_2}\int_{x}\int_{y}g( x,y)\psi(x)\psi(y)\lv X_\vep(\vep y+x,t)X_\vep(x,t)\phi(y)\d y\d x\d t\right)^2\biggr].
\end{align*} 
Here, the first two inequalities are comparable to \eqref{thm:tight-2}, and the last two inequalities are comparable to \eqref{thm:tight-3}. We omit the remaining details for the reason mentioned before Step~1. 

\subsection{Proof of Theorem~\ref{thm:main1} (2$\cc$): covariation identity and bounds}\label{sec:covaridbdd}

\begin{proof}[Proof of (\ref{m1:apriori}) in Theorem~\ref{thm:main1} (2$\cc$)]
We begin by proving the following identity: 
\begin{align}
\forall\;h\in \C_c(\R^2),\quad 
\E\left[\int_{t=0}^T\int_xh(x)\nu_\infty(\d x,\d t)\right]
&=
\lim_{n\to\infty}\int_{t=0}^T \int_x h(x)\E[\Lambda_{n}X_n(x,t)^2]\d x\d t\label{mg:density}.
\end{align}
To justify \eqref{mg:density}, we use the following direct consequence of Corollary~\ref{cor:muinfty}: 
\[
\int_{t=0}^T \int_x h(x)\nu_\infty(\d x,\d t)=\lim_{n\to\infty}\int_{t=0}^T \int_x h(x)\nu_n(\d x,\d t)=\lim_{n\to\infty}\int_{t=0}^T \int_x h(x)\Lambda_nX_n(x,t)^2\d x\d t.
\]
Theorem~\ref{thm:mombdd} (1$\cc$) ensures that the sequence of random variables in the rightmost limit is $L^2(\P)$-bounded, and thus forms a uniformly integrable family. Consequently, by a standard result on uniform integrability (see \cite[6.5.2 Theorem, p.263]{Ash}), \eqref{mg:density} follows.

The rest of this proof is devoted to justifying the following identity:
\begin{align}
\lim_{n\to\infty}\int_{t=0}^T \int_x h(x)\E[\Lambda_{n}X_n(x,t)^2]\d x\d t
&=\int_{t=0}^T\int_{x}h(x)\rDBG_tX_0^{\otimes 2}(x)\d x\d t,\label{formula:mg}
\end{align}
which we address through several steps. Step~1 demonstrates that the identity follows from the right-hand side of \eqref{mg:density}. This step adapts the method used in \cite[Proposition~2.6]{C:DBG} and assumes that
\begin{align}
&\lim_{n\to \infty}\Lambda_n\E^{B^{\bi_0}}_0\left[\exp\left(\Lambda_n\int_{r=0}^t \varphi_{\vep_n}(B^{\bi_0}_r)\d r\right)F_{x,T_0}(B^{\bi_0}_t)\right]\notag\\
&\quad =\int_{\tau=0}^t\s^\beta(\tau)P_{t-\tau}F_{x,T_0}(0)\d \tau ,\quad \forall\; x\in \R^2, \;0\leq t<\infty,\;0<T_0<\infty.\label{2mom:goal}
\end{align}
Here, we use the notation below \eqref{eq:momdual2} with $\bi_0=(2,1)$, and define
\begin{gather}
B^{\bi_0\prime}_t\,\defeq\,\frac{B^2_t+B^1_t}{\two},\quad 
F_{x,T_0}(y)\;\defeq\; \E^{B^{\bi_0\prime}}_{\two x}\biggl[X_0\biggl(\frac{B^{\bi_0\prime}_{T_0}+y}{\two}\biggr)X_0\biggl(\frac{B^{\bi_0\prime}_{T_0}-y}{\two}\biggr)\biggr]. \notag
\end{gather}
Step~2 turns to the proof of \eqref{2mom:goal}, where we modify the argument from \cite[Proof of (5.13) for Proposition~5.2, p.183]{C:DBG}. \medskip 

\noindent {\bf Step~1.} To prove \eqref{formula:mg}, we first utilize the moment duality in  \eqref{Mom:dual}. Hence, we have
\begin{align*}
\Lambda_n\E[X_n(x,t)^2]&=\Lambda_{\vep}\E^{(B^1,B^2)}_{(x,x)}\biggl[\exp\left(\Lambda_n\int_{r=0}^t \varphi_{\vep}(B^{\bi_0}_r)\d r\right)\prod_{i=1}^2X_0(B^i_t)\biggr]\\
&=\Lambda_{\vep}\E^{(B^1,B^2)}_{(x,x)}\biggl[\exp\biggl(\Lambda_n\int_{r=0}^t \varphi_{\vep}(B^{\bi_0}_r)\d r\biggr)X_0\biggl(\frac{B^{\bi_0\prime}_t+B^{\bi_0}_t}{\two}\biggr)X_0\biggl(\frac{B^{\bi_0\prime}_t-B^{\bi_0}_t}{\two}\biggr)\biggr].
\end{align*}
Observe that $B^{\bi_0\prime}$ and $B^{\bi_0}$ are independent two-dimensional standard Brownian motions. By integrating out $B^{\bi_0\prime}$, we obtain
\begin{align*}
\Lambda_n\E[X_n(x,t)^2]=\Lambda_n\E^{B^{\bi_0}}_0\left[\exp\left(\Lambda_n\int_{r=0}^t \varphi_{\vep_n}(B^{\bi_0}_r)\d r\right)F_{x,t}(B^{\bi_0}_t)\right].
\end{align*}
Then by \eqref{2mom:goal} with $T_0=t$, we get 
\begin{align}
\lim_{n\to \infty}\Lambda_n \E[X_n(x,t)^2]&=\int_{\tau=0}^t \s^\beta(\tau)\E^{B^{\bi_0}}_0[F_{x,t}(B^{\bi_0}_{t-\tau})]\d \tau\notag\\
&=\int_{\tau=0}^t \s^\beta(\tau)\E^{B^{\bi_0\prime}}_{\two x}\biggl[\E^{(B^{\bi_0},B^{\bi_0\prime})}_{(0,B^{\bi_0\prime}_\tau)}\E\biggl[X_0\biggl(\frac{B^{\bi_0\prime}_{t-\tau}+B^{\bi_0}_{t-\tau}}{\two}\biggr)X_0\biggl(\frac{B^{\bi_0\prime}_{t-\tau}-B^{\bi_0}_{t-\tau}}{\two}\biggr)\biggr]\biggr]\d \tau\notag\\
&=\int_{\tau=0}^t\s^\beta(\tau)\E^{B^{\bi_0\prime}}_{\two x}\biggl[\E^{(B^{1},B^{2})}_{(B^{\bi_0\prime}_\tau/\two,B^{\bi_0\prime}_\tau/\two)}\E\biggl[\prod_{i=1}^2X_0(B^i_{t-\tau})\biggr]\biggr] \d \tau\notag\\
&=\int_{\tau=0}^t \int_{z'}\begin{bmatrix}
P_{\tau/2}(x,z')\\
\s^\beta(\tau)
\end{bmatrix}_\times
\int_{\R^4}
\begin{bmatrix}
P_{t-\tau}(z', x_2^2)\\
P_{t-\tau}(z', x_2^1)
\end{bmatrix}_\times 
\begin{bmatrix}
X_0(x_2^2)\\
X_0(x_2^1)
\end{bmatrix}_\times \begin{bmatrix}
\d x_2^2\\
\d x_2^1
\end{bmatrix}_\otimes\d z'\d \tau\notag\\
&=\rDBG_t X_0^{\otimes 2}(x),\label{mg:limeqn}
\end{align}
where $\rDBG_t F(x)$ is as defined in \eqref{def:rDBG}. Moreover, by Theorem~\ref{thm:mombdd} (3$\cc$), we can apply \eqref{mg:limeqn} to obtain \eqref{formula:mg} via dominated convergence, since $h(x) \in \C_c(\R^2)$ by assumption. \medskip 

\noindent {\bf Step~2.} In this step, we prove \eqref{2mom:goal}, which will be accomplished in Step~2-3 below after demonstrating two limits in Steps~2-1 and 2-2.\medskip 

\noindent {\bf Step~2-1.} The first limit is as follows: 
\begin{align}
&\lim_{n\to \infty}\int_{t=0}^\infty\e^{-qt} \Lambda_n^2\E^{B^{\bi_0}}_0\left[\exp\left(\Lambda_n\int_{r=0}^t \varphi_{\vep_n}(B^{\bi_0}_r)\d r\right)\varphi_{\vep_n}(B^{\bi_0}_t)\right]\d t\notag\\
&\quad -\int_{t=0}^\infty\e^{-qt} \Lambda_n^2\int_y \phi(y)\E^{B^{\bi_0}}_{{\vep_n} y/\two}\left[\exp\left(\Lambda_n\int_{r=0}^t \varphi_{\vep_n}(B^{\bi_0}_r)\d r\right)\varphi_{\vep_n}(B^{\bi_0}_t)\right]\d y\d t=0.\label{2mom:0}
\end{align}

To prove \eqref{2mom:0}, we use the expansion $\e^{\int_0^t h(s)\d s}=1+\int_0^t h(s)\e^{\int_s^th(r)\d r}\d s$. Hence, for all $x\in\R^2$, we have
\begin{align}
&\int_0^\infty \e^{-qt}\Lambda_n^2\E^{B^{\bi_0}}_{{\vep_n} x}\left[\exp\left(\Lambda_n\int_{r=0}^t \varphi_{{\vep_n}}(B^{\bi_0}_r)\d r\right)\varphi_{\vep_n}(B^{\bi_0}_t)\right]\d t\notag\\
&\quad =\int_0^\infty \e^{-qt}\Lambda_n^2\int_{y'} P_t({\vep_n} x,{\vep_n} y')\varphi(y')\d y'\d t+\int_0^\infty \e^{-qt}\Lambda_n^3\int_{s=0}^t  \int_{y'} P_s({\vep_n} x,{\vep_n} y')\varphi(y')\notag\\
&\quad\quad  \times \E^{B^{\bi_0}}_{{\vep_n} y'}\left[\exp\left(\Lambda_n\int_{r=0}^{t-s} \varphi_{\vep_n}(B^{\bi_0}_r)\d r\right)\varphi_{\vep_n}(B^{\bi_0}_{t-s})\right]\d y'\d s\d t.\label{fotc:expansion}
\end{align}
We use this identity for $x = 0$ and $x = y/\two$, integrate out $y$ against $\phi(y)\d y$, and subtract the resulting equations. We claim that in the limit, the two resulting terms on the right-hand side vanish. For the first term, we perform the following calculation, utilizing the fact that $\phi$ is a probability density, as stated in 
Proposition~\ref{prop:XM} (1$\cc$):
\begin{align}
&\biggl|\int_0^\infty \e^{-qt}\Lambda_n^2\int_{y'} P_t(0,{\vep_n} y')\varphi(y')\d y'\d t
-\int_0^\infty \e^{-qt}\Lambda_n^2\int_y\int_{y'} P_t({\vep_n} y/\two,{\vep_n} y')\varphi(y')\d y'\phi(y)\d y\d t\biggr|
\notag\\
&\quad \leq\int_{t=0}^\infty \e^{-q{\vep}_n^2 t}\Lambda_n^2 \int_y\int_{y'} \lvert P_t(y')-P_t(y/\two,y')\rvert \varphi(y')\d y'\phi(y)\d y\d t\notag\\
&\quad \leq \Lambda_{n}^2C(\varphi)  \int_{t=0}^\infty\frac{1}{(t\vee 1)^2} \d t\xrightarrow[n\to \infty]{} 0,\label{2mom:1}
\end{align}
where the last inequality holds because $1-\e^{-a}\leq a$ for all $a\geq 0$, and $\lvert\lvert y' \rvert^2 - \lvert y/\two - y'\rvert ^2\rvert \less \lvert y \rvert^2 + \lvert y' \rvert^2$. To show that the other limit vanishes, we use the following two observations:
\begin{gather*}
\lim_{n\to \infty}\sup_{y\in \supp(\phi)}\Lambda_n\int_{t=0}^\infty \e^{-qt }\int_{y'} \lvert P_t({\vep_n} y')-P_t(\vep y/\two,{\vep_n} y')\rvert \varphi(y')\d y'\d t=0,\\
\limsup_{n\to \infty}\int_{t=0}^\infty q\e^{-qt}\sup_{y''\in \R^2}\Lambda_n^2\int_{s=0}^t \E^{B^{\bi_0}}_{{\vep_n} y''}\left[\exp\left(\Lambda_n\int_{r=0}^{s} \varphi_{{\vep_n} }(B^{\bi_0}_r)\d r\right)\varphi_{\vep_n}(B^{\bi_0}_{s})\right]\d s\d t<\infty.
\end{gather*}
Here, the equality in the first line is justified in the same way as in the proof of \eqref{2mom:1}, and the inequality in the second line holds for all sufficiently large $q$ as a direct consequence of \cite[Lemma~5.4]{C:DBG}. For any such $q$, we perform the following calculation to justify the limit: 
\begin{align}
&\biggl|\int_0^\infty \e^{-qt}\Lambda_n^3\int_{s=0}^t  \int_{y'} P_s(0,{\vep_n} y')\varphi(y')\E^{B^{\bi_0}}_{{\vep_n} y'}\left[\exp\left(\Lambda_n\int_{r=0}^{t-s} \varphi_{\vep_n}(B^{\bi_0}_r)\d r\right)\varphi_{\vep_n}(B^{\bi_0}_{t-s})\right]\d y'\d s\d t\notag\\
&\quad -\int_0^\infty \e^{-qt}\Lambda_n^3\int_{s=0}^t\int_y  \int_{y'} P_s({\vep_n}y/\two,{\vep_n} y')\varphi(y')\notag\\
&\quad \times \E^{B^{\bi_0}}_{{\vep_n} y'}\left[\exp\left(\Lambda_n\int_{r=0}^{t-s} \varphi_{\vep_n}(B^{\bi_0}_r)\d r\right)\varphi_{\vep_n}(B^{\bi_0}_{t-s})\right]\d y'\phi(y)\d y\d s\d t\biggr|\notag\\
&\quad \leq \int_{t=0}^\infty \e^{-qt }\Lambda_{n}^3\int_{s=0}^t \int_y\int_{y'} \lvert P_s({\vep_n} y')-P_s({\vep_n} y/\two,{\vep_n} y')\rvert \varphi(y')\notag\\
&\quad\quad  \times \E^{B^{\bi_0}}_{{\vep_n} y'}\left[\exp\left(\Lambda_n\int_{r=0}^{t-s} \varphi_{{\vep_n}}(B^{\bi_0}_r)\d r\right)\varphi_{\vep_n}(B^{\bi_0}_{t-s})\right]\d y'\phi(y)\d y\d s\d t\notag\\
&\quad \leq \Lambda_{n}\int_y \int_{y'}\biggl(\int_{t=0}^\infty \e^{-qt } \lvert P_t({\vep_n} y')-P_t({\vep_n} y/\two,{\vep_n} y')\rvert \varphi(y')\d t\biggr)\notag\\
&\quad\quad  \times \biggl(\int_{t=0}^\infty \e^{-qt}\Lambda_n^2\E^{B^{\bi_0}}_{{\vep_n} y'}\left[\exp\left(\Lambda_n\int_{r=0}^{t} \varphi_{{\vep_n}}(B^{\bi_0}_r)\d r\right)\varphi_{\vep_n}(B^{\bi_0}_{t})\right]\d t\biggr)\d y'\phi(y)\d y\notag\\
&\quad = \Lambda_n\int_y \int_{y'}\biggl(\int_{t=0}^\infty \e^{-qt } \lvert P_t({\vep_n} y')-P_t({\vep_n} y/\two,{\vep_n} y')\rvert \varphi(y')\d t\biggr)\notag \\
&\quad\quad \times\biggl(\int_{t=0}^\infty q\e^{-qt}\Lambda_{n }^2 \int_{s=0}^t \E^{B^{\bi_0}}_{{\vep_n} y'}\left[\exp\left(\Lambda_{n}\int_{r=0}^{s} \varphi_{n}(B^{\bi_0}_r)\d r\right)\varphi_{\vep_n}(B^{\bi_0}_{s})\right]\d s\d t\biggr)\d y' \phi(y)\d y \notag\\
&\quad  \xrightarrow[n\to\infty]{}0,\label{2mom:2}
\end{align}
where the second inequality and the last equality both use elementary properties of Laplace transforms, and the limit can be justified by using the two observations mentioned above. 

In summary, \eqref{2mom:0} follows from \eqref{2mom:1} and \eqref{2mom:2} via \eqref{fotc:expansion}.\medskip 

\noindent {\bf Step 2-2.} The second limit states that for all large enough $q>0$,
\begin{align}
\lim_{n\to \infty}\int_{t=0}^\infty\e^{-qt} \Lambda_n^2\E^{B^{\bi_0}}_0\left[\exp\left(\Lambda_n\int_{r=0}^t \varphi_{\vep_n}(B^{\bi_0}_r)\d r\right)\varphi_{\vep_n}(B^{\bi_0}_t)\right]\d t =\int_{t=0}^\infty \e^{-qt}\s^\beta(t)\d t.\label{2mom:goal1}
\end{align}

To see \eqref{2mom:goal1}, we use \eqref{2mom:0}, thereby obtaining
\begin{align*}
&\lim_{n\to \infty}\int_{t=0}^\infty\e^{-qt} \Lambda_n^2\E^{B^{\bi_0}}_0\left[\exp\left(\Lambda_n\int_{r=0}^t \varphi_{\vep_n}(B^{\bi_0}_r)\d r\right)\varphi_{\vep_n}(B^{\bi_0}_t)\right]\d t\\
&\quad =\lim_{n\to \infty}\int_{t=0}^\infty\e^{-qt} \Lambda_n^2\int_y \phi(y)\E^{B^{\bi_0}}_{{\vep_n} y/\two}\left[\exp\left(\Lambda_n\int_{r=0}^t \varphi_{{\vep_n}}(B^{\bi_0}_r)\d r\right)\varphi_{\vep_n}(B^{\bi_0}_t)\right]\d y\d t\\
&\quad =\frac{4\pi}{\log (q/\beta)}=\int_{t=0}^\infty \e^{-qt}\s^\beta(t)\d t,
\end{align*}
where the second equality follows from \cite[(5.12)]{C:DBG} (see \cite[Section~4.4]{C:Additive} for an alternative proof of this equality), and the last equality uses \cite[(5.2)]{C:DBG}. \medskip 

\noindent {\bf Step 2-3.} We now proceed to prove \eqref{2mom:goal}. By using L\'evy's continuity argument (see \cite[Lemma~4.18, p.165]{C:DBG}), \eqref{2mom:goal1} yields
\[
\lim_{n\to\infty}\int_{\tau=0}^t \Lambda_n^2\E^{B^{\bi_0}}_0\left[\exp\left(\Lambda_n\int_{r=0}^\tau \varphi_{\vep_n}(B^{\bi_0}_r)\d r\right)\varphi_{\vep_n}(B^{\bi_0}_\tau)\right]\d \tau=\int_{\tau=0}^t \s^\beta(\tau)\d \tau,\quad \forall\;t\geq 0.
\]
Hence, by the bounded continuity of $\tau\mapsto P_{t-\tau}F_{x,T_0}(0)$, $0\leq \tau\leq t$, we get
\begin{align*}
&\int_{\tau=0}^t \s^\beta(\tau)P_{t-\tau}F_{x,T_0}(0)\d \tau\\
&\quad =\lim_{n\to \infty}\int_{\tau=0}^t\Lambda_n^2\E^{B^{\bi_0}}_0\left[\exp\left(\Lambda_n\int_{r=0}^\tau \varphi_{\vep_n}(B^{\bi_0}_r)\d r\right)\varphi_{\vep_n}(B^{\bi_0}_\tau)\right]P_{t-\tau}F_{x,T_0}(0)\d \tau\\
&\quad =\lim_{n\to \infty}\Lambda_n\E^{B^{\bi_0}}_0\left[\exp\left(\Lambda_n\int_{r=0}^t \varphi_{\vep_n}(B^{\bi_0}_r)\d r\right)F_{x,T_0}(B^{\bi_0}_t)\right],\notag
\end{align*}
where the last equality holds since $\e^{\int_0^t h(s)\d s}=1+\int_0^t h(s)\e^{\int_s^th(r)\d r}\d s$. This proves \eqref{2mom:goal}. 
\end{proof}

\begin{proof}[Proof of (\ref{MMbdd}) in Theorem~\ref{thm:main1} (2$\cc$)]
For any nonnegative $g_1,g_2\in \C^0_{pd}(\R^2\times \R^2)$, it holds almost surely that
\begin{align*}
&\int_{s_1=T_1}^{T_2}\int_{x_1}\int_{s_2=s_1}^{T_2}
\int_{x_2}\prod_{j=1}^2g_j(x_j,x_j)\la M_\infty,M_\infty\ra(\d x_2,\d s_2)
\la M_\infty,M_\infty\ra(\d x_1,\d s_1)\\
&\quad =\lim_{n\to\infty}\int_{s_1=T_1}^{T_2}\int_{x_1}\int_{x_1'}\int_{s_2=s_1}^{T_2}
\int_{x_2}\int_{x_2'}\prod_{j=1}^2g_j(x_j',x_j)
 \la M_n,M_n\ra(\d x_2',\d x_2,\d s_2) \\
 &\quad \quad \times \la M_n,M_n\ra(\d x_1',\d x_1,\d s_1)\\
 &\quad =\lim_{n\to\infty}\int_{s_1=T_1}^{T_2}\int_{x_1}\int_{\wt{x}_1}\int_{s_2=s_1}^{T_2}
\int_{x_2}\int_{\wt{x}_2}
\prod_{j=1}^2
g_j(\vep_n\wt{x}_j+x_j,x_j)\Lambda_n
X_n(\vep_n \wt{x}_j+x_j,s_j)X_n(x_j,s_j)\phi(\wt{x}_j)\\
&\quad \quad \times \d \wt{x}_2\d x_2\d s_2\d \wt{x}_1\d x_1\d s_1,
\end{align*}
where the first equality follows from Corollary~\ref{cor:muinfty} and \eqref{Mncovar:derivation:eq}, and the last equality applies \eqref{eq:covarM-2}. Hence, applying Theorem~\ref{thm:mombdd} (1$\cc$) via Fatou's lemma yields that, for any $p\in (0,1)$, 
\begin{align*}
&\E\biggl[\int_{s_1=T_1}^{T_2}\int_{x_1}\int_{s_2=s_1}^{T_2}
\int_{x_2}\prod_{j=1}^2g_j(x_j,x_j)\la M_\infty,M_\infty\ra(\d x_2,\d s_2)
\la M_\infty,M_\infty\ra(\d x_1,\d s_1)\biggr]\\
&\quad \leq  \liminf_{n\to\infty}\int_{s_1=T_1}^{T_2}\int_{x_1}\int_{\wt{x}_1}\int_{s_2=s_1}^{T_2}
\int_{x_2}\int_{\wt{x}_2}
\biggl(\prod_{j=1}^2
g_j(\vep_n\wt{x}_j+x_j,x_j) \phi(\wt{x}_j)\biggr)\\
&\quad \quad \times \Lambda_n^2\E[
X_n(\vep_n \wt{x}_1+x_1,s_1)X_n(x_1,s_1)X_n(\vep_n \wt{x}_2+x_2,s_2)X_n(x_2,s_2)] \d \wt{x}_2\d x_2\d s_2\d \wt{x}_1\d x_1\d s_1\\
&\quad \leq C(\lambda,\phi,\lVert X_0\rVert_\infty,T,p)\biggl(\prod_{j=1}^2\int_{x_j}g_j(x_j,x_j)\d x_j\biggr)(T_2-T_1)^{2-p}.
\end{align*}
Since $g_1,g_2$ are arbitrary, the last inequality is enough to prove \eqref{MMbdd}.
\end{proof}

\begin{proof}[Proof of (\ref{X4BMbdd}) in Theorem~\ref{thm:main1} (2$\cc$)]
The proof proceeds similarly to that of \eqref{MMbdd}, except we now use Theorem~\ref{thm:mombdd} (2$\cc$) in place of (1$\cc$), together with the following almost sure convergence:
\begin{align*}
&\lim_{n\to\infty}\int_{s=0}^t\int_{x}\int_{x'}f^{1/2}(x')f^{1/2}(x)
\begin{bmatrix}
\int_{z_2}P_{a_2}(y_2,z_2)X_n(\d z_2,s)\\
\int_{z_1} P_{a_1}(y_1,z_1)X_n(\d z_1,s)
\end{bmatrix}_\times \la M_n,M_n\ra(\d x',\d x,\d s)\\
&\quad =\int_{s=0}^t\int_{x}f(x)\begin{bmatrix}
\int_{z_2}P_{a_2}(y_2,z_2)X_\infty(\d z_2,s)\\
\int_{z_1} P_{a_1}(y_1,z_1)X_\infty(\d z_1,s)
\end{bmatrix}_\times
\la M_\infty,M_\infty\ra(\d x,\d s).
\end{align*}
More precisely, besides \eqref{eq:Sko}, this convergence depends on the uniform convergence of the integrands as functions of $s \in [0, T]$ as $n \to \infty$ (see Corollary~\ref{cor:muinfty}).
 \end{proof}

\section{Heat kernel expansions and integral bounds}\label{sec:heat}
In this section, we collect two sets of properties of the two-dimensional heat kernels  
$P_t(x,y)=P_t(x-y)$ defined in \eqref{def:Pt}. The properties of the first set are gathered in Lemma~\ref{lem:heat1} below. It shows expansions that refine or restate some of the results established in \cite{C:Whittaker,C:DBG}. 

\begin{lem}\label{lem:heat1}
{\rm (1$\cc$)} For all $y\in \R^2\setminus\{0\}$ and $T\in (0,\infty)$, it holds that
\begin{linenomath*}\begin{align}\label{limT:glue}
\begin{split}
\int_0^{T}P_{2r}(y)\d r&=\frac{1}{4\pi}\log \left(\frac{4T}{\lvert y \rvert^2}\right) -\frac{\gamma_{\sf EM}}{4\pi}+ \vep_{\ref{limT:glue}}\left(\frac{4T}{\lvert y \rvert^2}\right).
\end{split} 
\end{align}\end{linenomath*}
Here, $\vep_{\ref{limT:glue}}(a)$ is a nonnegative function defined by
\begin{linenomath*}\begin{align*}
\vep_{\ref{limT:glue}}(a)&\;\defeq \,\int_{0}^{\tfrac{1}{a}}\frac{(\1_{(0,1]}(t)-\e^{-t})}{4\pi t}\d t+\frac{\log^- a}{4\pi},\quad a\in (0,\infty),
\end{align*}\end{linenomath*}
where $\log^- x\,\defeq-\log(x\wedge 1)$ is the negative part  of $\log x$.
In particular,
\begin{align}\label{keyvep:bdd}
0\leq \vep_{\ref{limT:glue}}(a)\less a^{-1}\1(a\geq 1)+(1+\log^- a)\1 (a<1).
\end{align}

\noindent {\rm (2$\cc$)} For all $\vep,\delta \in (0,\infty)$ with 
$\vep/\delta^{1/2}\leq 1$
and all $M\in (0,\infty)$,
\begin{linenomath*}\begin{align}\label{gauss:vep}
\begin{split}
&\sup_{\lvert x\rvert \leq M}\lvert P_t(\vep x+y)-P_t(y)\rvert \leq C(M)(\vep/\delta^{1/2})P_{2t}(y),\; \forall\;y\in \R^2,\;t\geq \delta.
\end{split}
\end{align}\end{linenomath*}
Also, for all $\delta\in(0,\infty)$, 
\begin{align}\label{kernel:tmod}
\begin{split}
\lvert P_{t_1}(x)-P_{t_0}(x)\rvert &\leq C(\delta)\lvert t_1-t_0\rvert ,\quad\forall\;x\in \R^2,\;t_1,t_0\geq \delta,\\
\lvert P_{t_1}(x)-P_{t_0}(x)\rvert &\leq \lvert t_1-t_0\rvert \e^{-\lvert x\rvert ^2/\delta}\left(\frac{\lvert x\rvert ^2+2\max\{t_0,t_1\}}{4\pi \delta^3}\right),\quad \forall\;x\in \R^2,\;t_1,t_0\geq \delta.
\end{split}
\end{align}
\end{lem}

\begin{rmk}\label{rmk:Green}
The following asymptotic expansion  as $\vep\to 0$, a consequence of an integral representation of the Macdonald function $K_0$  \cite[(5.10.25), p.119]{Lebedev}
and a series expansion of $K_0$ \cite[(5.7.11), p.110]{Lebedev}, can also be deduced from \eqref{limT:glue}: for all $q\in (0,\infty)$ and $x\in \R^2\setminus\{0\}$, 
\begin{align}\label{rmk:EM}
\int_0^\infty\frac{\e^{-q t}}{4\pi t}\exp\left(-\frac{\vep^2\lvert x\rvert ^2}{4t}\right)\d t=\frac{\log \vep^{-1}}{2\pi}+\frac{1}{2\pi }\left(\frac{\log 4-\log q}{2}-\log \lvert x\rvert -\EM\right)+o(1).
\end{align}
More specifically, to obtain \eqref{rmk:EM} from \eqref{limT:glue}, note $\int_0^\infty \e^{-qt}f(t)\d t=\int_0^\infty q\e^{-qT}\int_0^T f(t)\d t\d T$ and the following representation of $\EM$: $\EM=-\int_0^\infty \e^{-T'}\log T'\d T'$
(cf. \cite[(1.3.19), p.8]{Lebedev}). \hfill $\blacklozenge$
\end{rmk}

\begin{details}
\begin{proof}[Details of Remark~\ref{rmk:Green}]
By the known expansion, we have
\begin{align*}
\int_0^\infty\frac{ \e^{-qt}}{4\pi t}\exp\left(-\frac{\vep^2\lvert x\rvert ^2}{4t}\right)\d t&=\frac{1}{2}\int_0^\infty\frac{ \e^{-t}}{2\pi t}\exp\left(-\frac{\vep^2\lvert x\rvert ^2q/2}{2t}\right)\d t\\
&=\frac{1}{2\pi}\log \frac{\sqrt{2}}{\vep \lvert x\rvert \sqrt{q}}+\frac{1}{2\pi}\left(\frac{\log 2-\log 1}{2}-\EM\right)\\
&\quad +\mathcal O\left(\frac{\vep^2\lvert x\rvert ^2q}{2}\log \frac{1}{\sqrt{\vep^2\lvert x\rvert ^2q/2}}\right)\\
&=\frac{\log \vep^{-1}}{2\pi}+\frac{1}{2\pi }\left(\frac{\log 4-\log q}{2}-\log \lvert x\rvert -\EM\right)+o(1).
\end{align*}
On the other hand, to recover this expansion, note that \eqref{limT:glue} implies
\begin{align*}
&\int_0^\infty q\e^{-qT}\left[\frac{1}{4\pi}\log \left(\frac{4T}{\vep^2\lvert x\rvert ^2}\right)-\frac{\EM}{4\pi}+\vep_{\ref{limT:glue}}\left(\frac{4T}{\vep^2\lvert x\rvert ^2}\right)\right]\d T\\
&\quad =\frac{\log (4q^{-1}\vep^{-2}\lvert x\rvert ^{-2})}{4\pi}-\frac{\EM}{4\pi}+\frac{1}{4\pi}\int_0^\infty \e^{-T'}\log T'\d T'\\
&\quad \quad +\mathcal O\left(\int_{\vep^2\lvert x\rvert ^2/4}^\infty q\e^{-qT}\frac{\vep^2\lvert x\rvert ^2}{4T}\d T\right) +\mathcal O\left(\int_0^{\vep^2\lvert x\rvert ^2/4} q\e^{-qT}\left(1+\log^-\frac{4T}{\vep^2\lvert x\rvert ^2}\right)\d T\right).
\end{align*}
\end{proof}
\end{details}

\begin{proof}[Proof of Lemma~\ref{lem:heat1}]
(1$\cc$) The expansion in \eqref{limT:glue} follows from a slight modification of the proof of 
\cite[Theorem~4.18, pp.2946--2947]{C:Whittaker}. By the change of variables $v=r/\lvert y \rvert^2$, we can write 
\begin{align}
\int_0^TP_{2r}(y)\d r&=\int_0^{\frac{T}{\lvert y \rvert^2}}\frac{(\e^{-\frac{1}{4v}}-\1_{[1,\infty)}(4v))}{4\pi v}\d v+\1_{[1,\infty)}\left(\frac{4T}{\lvert y \rvert^2}\right)\cdot \frac{1}{4\pi}\log \left(\frac{4T}{\lvert y \rvert^2}\right)\notag\\
&= \frac{1}{4\pi}\log \left(\frac{4T}{\lvert y \rvert^2}\right)-\int_0^{\infty}\frac{(\1_{[1,\infty)}(4v)-\e^{-\frac{1}{4v}})}{4\pi v}\d v\notag\\
&\quad +\int_{(\frac{4T}{\lvert y \rvert^2})/4}^{\infty}\frac{(\1_{[1,\infty)}(4v)-\e^{-\frac{1}{4v}})}{4\pi v}\d v-\1_{ (0,1)
}\left(\frac{4T}{\lvert y \rvert^2}\right)\cdot \frac{1}{4\pi}\log \left(\frac{4T}{\lvert y \rvert^2}\right).\label{EM:proof}
\end{align}
Note that the second and third terms on the right-hand side of \eqref{EM:proof} are equal to the next integral with $a=0$ and with $a=4T/\lvert y \rvert^2$, respectively: 
\begin{align}\label{EM:partial}
\int_{a/4}^\infty\frac{(\1_{[1,\infty)}(4v)-\e^{-\frac{1}{4v}})}{ 4\pi v}\d v=\int_0^{\frac{1}{a}}\frac{(\1_{(0,1]}(t)-\e^{-t})}{4\pi t}\d t,
\end{align}
where the equality follows by using the change of variables $t=1/(4v)$. In particular, 
when $a=0$, the right-hand side of \eqref{EM:partial} equals $\EM/(4\pi)$ by a standard integral representation of $\EM$ \cite[the display next to (1.3.19), p.8]{Lebedev}. Hence, \eqref{EM:proof} is enough to get \eqref{limT:glue}. 

It remains to show \eqref{keyvep:bdd}. To get the nonnegativity of $\vep_{\ref{limT:glue}}(a)$, it is enough to
note that the integral on the right-hand side of \eqref{EM:partial} is positive when $a\geq 1$ and decreases to $\EM/(4\pi)>0$ over $0<a<1$ as $a\searrow 0$.
Also, the second inequality in \eqref{keyvep:bdd} can be seen by applying the bound $1-\e^{-t}\leq t$, valid for all $t\geq 0$, to the integral in $\vep_{\ref{limT:glue}}(a)$. We have proved all of the properties stated in (1$\cc$). 
\medskip

\noindent {(2$\cc$)} The estimates in \eqref{gauss:vep} and \eqref{kernel:tmod} restate \cite[Lemma~4.16 (2$\cc$) and (3$\cc$), pp.161--162]{C:DBG}. 
\end{proof}

The next lemma gives the second set of properties. It shows bounds for various integrals used in the earlier sections; the simple bound in \eqref{ratiot} is the central tool. 

\begin{lem}\label{lem:heat2}
For all $T_0,T_1\in (0,\infty)$, $a\in (0,\infty)$, 
 $p,p'\in \R_+$, $p_1,p_2\in (0,\infty)$ and $y,y'\in \R^2$, the following holds:
\begin{align}
\frac{1}{1+\lvert y-y'\rvert^p}&\leq C(p)\frac{1+\lvert y' \rvert^p}{1+\lvert y \rvert^p},\label{ratiot}\\
 \int_{x'}\int_{t=0}^{T_0} \frac{(\lvert x'\rvert^{p_1}+t^{p_2})P_{at}(0)P_{t}(x')}{1+\lvert y-x'\rvert^p}\d t\d x'   &\leq  C(a,p,p_1,p_2)  \frac{\sum_{j=1}^2(T_0^{jp_j/2}+T_0^{(jp_j+p)/2})}{1+\lvert y \rvert^p},\label{ineq:heat2-1}\\
\int_{x'}\frac{\lvert x'\rvert^{p'}P_{T_0}(x')}{1+\lvert y-x'\rvert^p}\d x' &\leq C(p,p')\frac{T_0^{p'/2}+T_0^{(p'+p)/2}}{1+\lvert y \rvert^p},\label{ineq:heat2-2}\\
\int_{x'}\frac{P_{T_0}(y,x')P_{T_1}(x',y')}{1+\lvert x'\rvert^p}\d x'  &\leq C(p)\frac{1+T_1^{p/2}}{1+\lvert y' \rvert^p}P_{T_0+2T_1}(y-y').\label{ineq:heat2-4}
\end{align}
Here, we adopt the continuity convention that $\lvert y'' \rvert^0=1$ when $0=y''\in \R^2$.
\end{lem}
\begin{proof} 
We prove \eqref{ratiot}--\eqref{ineq:heat2-4} in Steps~1--4 below, respectively. \medskip

\noindent {\bf Step 1.} The required inequality in \eqref{ratiot} follows upon noting that
\begin{align*}
\frac{1}{1+\lvert y-y'\rvert^p}&=
\frac{1}{1+\lvert y \rvert^p}\frac{1+\lvert y \rvert^p}{1+\lvert y-y'\rvert^p}
\leq \frac{C(p)}{1+\lvert y \rvert^p} \frac{1+\lvert y-y'\rvert^p+\lvert y' \rvert^p}{1+\lvert y-y'\rvert^p}
\leq C(p)\frac{1+\lvert y' \rvert^p}{1+\lvert y \rvert^p},
\end{align*}
where the first equality uses the bound
 $\lvert y \rvert^p\leq (\lvert y-y'\rvert+\lvert y' \rvert)^p\leq C(p)(\lvert y-y'\rvert^p+\lvert y' \rvert^p)$.
\medskip  
 
\noindent {\bf Step~2.} By \eqref{ratiot}, we get
\begin{align}
&\int_{x'}\int_{t=0}^{T_0} \frac{(\lvert x'\rvert^{p_1}+t^{p_2})P_{at}(0)P_{t}(x')}{1+\lvert y-x'\rvert^p}\d t\d x'\notag\\
&\quad \leq \frac{C(p)}{1+\lvert y \rvert^p} \int_{x'}(\lvert x'\rvert^{p_1}+\lvert x'\rvert^{p_1+p})\int_{t=0}^{T_0} \frac{1}{4\pi^2at^2}\exp\left(-\frac{\lvert x'\rvert^2}{2t}\right)\d t\d x'\notag\\
&\quad\quad +\frac{C(p)}{1+\lvert y \rvert^p}\int_{t=0}^{T_0}t^{-1+p_2}\int_{x'}(1+\lvert x'\rvert^{p})
\frac{1}{4\pi^2at}\exp\left(-\frac{\lvert x'\rvert^2}{2t}\right)\d x'\d t\notag\\
&\quad \leq\frac{C(a,p)}{1+\lvert y \rvert^p}\left[ \int_{x'}\left(\lvert x'\rvert^{p_1-2}+\lvert x'\rvert^{p_1+p-2}\right)\exp\left(-\frac{\lvert x'\rvert^2}{2T_0}\right)\d x'+\int_{t=0}^{T_0}(t^{-1+p_2}+t^{-1+p_2+p/2})\d t\right]\notag\\
&\quad \leq C(a,p,p_1,p_2)\frac{T_0^{p_1/2}+T_0^{(p_1+p)/2}+T_0^{ p_2}+T_0^{ p_2+p/2}}{1+\lvert y \rvert^p},\label{ineq:heat2-1-final}
\end{align}
where the second inequality applies the formula $\int t^{-2}\e^{-b/t}\d t=\e^{-b/t}/b+C$, for $b\neq 0$ and $t>0$, and the following scaling property: 
\begin{align}
\forall\;\wt{p}\in (-2,\infty),\;t\in(0,\infty),\quad 
\int_{x'}\lvert x'\rvert^{\wt{p}}\exp\left(-\frac{\lvert x'\rvert^2}{2t}\right)\d x'
=C(\wt{p})t^{(\wt{p}+2)/2},\label{int:scaling}
\end{align}
and the inequality in \eqref{ineq:heat2-1-final} is also obtained by using \eqref{int:scaling}. We obtain  \eqref{ineq:heat2-1} from \eqref{ineq:heat2-1-final}. 
\medskip  
 
\noindent {\bf Step~3.}
The proof of \eqref{ineq:heat2-2} is similar. Applying \eqref{ratiot} and \eqref{int:scaling} in the same order, we get
\begin{align*}
\int_{x'}\frac{\lvert x'\rvert^{p'}P_{T_0}(x')}{1+\lvert y-x'\rvert^p}\d x'&\leq C(p)\int_{x'}\frac{(\lvert x'\rvert^{p'}+\lvert x'\rvert^{p'+p})P_{T_0}(x')}{1+\lvert y \rvert^p}\d x'\leq C(p,p')\frac{T_0^{p'/2}+T_0^{(p'+p)/2}}{1+\lvert y \rvert^p}.
\end{align*}
\mbox{}

\noindent {\bf Step~4.}
To prove \eqref{ineq:heat2-4}, we first change variables to get the following equality:
\begin{align*}
\int_{x'}\frac{P_{T_0}(y,x')P_{T_1}(x',y')}{1+\lvert x'\rvert^p}\d x'&=\int_{x''}\frac{P_{T_0}(y-y',x'')P_{T_1}(x'')}{1+\lvert x''+y'\rvert^p}\d x''\\
&\leq C(p) \int_{x''}\frac{P_{T_0}(y-y',x'')P_{T_1}(x'')(1+\lvert x''\rvert^p)}{1+\lvert y' \rvert^p}\d x'',
\end{align*}
where the last inequality uses \eqref{ratiot}. To continue, we rewrite the right-hand side so that 
\begin{align*}
\int_{x'}\frac{P_{T_0}(y,x')P_{T_1}(x',y')}{1+\lvert x'\rvert^p}\d x'
&\leq \frac{C(p)}{1+\lvert y' \rvert^p}\int_{x''}P_{T_0}(y-y',x'')P_{T_1}(x'')\d x''\\
&\quad+\frac{C(p)T_1^{p/2}}{1+\lvert y' \rvert^p}\int_{x''}P_{T_0}(y-y',x'')P_{T_1}(x'')\left|\frac{x''}{\sqrt{2T_1}}\right|^p\d x''\\
&\leq \frac{C(p)}{1+\lvert y' \rvert^p}P_{T_0+T_1}(y,y') +\frac{C(p)T_1^{p/2}}{1+\lvert y' \rvert^p}\int_{x''}P_{T_0}(y-y',x'')P_{2T_1}(x'')\d x''\\
&\leq C(p)\frac{1+T_1^{p/2}}{1+\lvert y' \rvert^p}P_{T_0+2T_1}(y-y').
\end{align*}
Here, we have applied the Chapman--Kolmogorov equation to derive the two terms on the right-hand side of the second inequality; the derivation of the second term also uses the bound $\e^{-a}a^p\leq C(p) \e^{-a/2}$ for all $a\geq 0$. The last inequality proves \eqref{ineq:heat2-4}. The proof of Lemma~\ref{lem:heat2} is complete.
\end{proof}

\section{AI disclosure} 
The paper has used AI for language editing. All ideas and proofs in this work were developed by the author.

\end{document}